\documentclass[10pt]{amsart}
\usepackage{amssymb,amsmath,amsfonts,amsthm,graphics,mathrsfs,amscd,hyperref}
\usepackage[makeroom]{cancel}
\usepackage[hmargin=1in,vmargin=1in]{geometry}
\usepackage[all,cmtip]{xy}

\theoremstyle{plain}

\newtheorem{definition}[equation]{Definition}

\newtheorem{lemma}[equation]{Lemma}
\newtheorem{proposition}[equation]{Proposition}
\newtheorem{theorem}[equation]{Theorem}

\newtheorem{remark}[equation]{Remark}

\numberwithin{equation}{subsection}

\title[Functors of Artin rings]{On deformations of foliated complex analytic structures II : Functors of Artin rings}
\author{Chunghoon Kim}

\begin{document}

\maketitle
\begin{abstract}
We study deformations of foliated nonsingular varieties defined by locally free subsheaves of tangent sheaves and locally free subshseaves of cotangent sheaves in the language of functors of Artin rings. We identify first-order deformations and obstructions.
\end{abstract}

\tableofcontents

 Our presentation is based on \cite{Ser06}. We will always denote by $k=\mathbb{C}$ the field of complex numbers. All schemes will be assumed to be defined over $k$, locally-noetherian and separated. If $s\in S$, we denote by $k(s)=\mathcal{O}_{S,s}/\mathfrak{m}_s$ the residue field of $S$ at $s$. We denote by $\bold{Art}$ the category of local artinian $k$-algebras with residue field $k$, and by $\bold{\widehat{Art}}$ the category of complete local noetherian $k$-algebras with residue field $k$, and by $\bold{Art}^*$ the category of local noetherian $k$-algebras with residue field $k$. 
 
\begin{definition}[small extension]
We say that for $\left( \tilde{A}, \tilde{\mathfrak{m}} \right), \left(A, \mathfrak{m} \right)\in \bold{Art}$, an exact sequence of the form $0\to (t)\to \tilde{A} \to A \to 0$ is a small extension if $t\in \tilde{\mathfrak{m}}$ is annihilated by $\tilde{\mathfrak{m}}$ $($i.e. $t\cdot \tilde{\mathfrak{m}}=0)$, so that $(t)$ is an one dimensional $k$-vector space.
\end{definition}

\section{Generalities on foliated deformations in terms of tangent sheaves}

We extend the presentation \cite{Ser06} p.20-p.21 in the context of foliated deformations in terms of tangent sheaves. Let $X$ be a nonsingular variety over $k$ and $T_{\mathcal{F}_0}$ be a singular foliation on $X$. In other words,  we have an exact sequence
\begin{align*}
0\to T_{\mathcal{F}_0}\to T_X \to T_X/T_{\mathcal{F}_0} \to 0
\end{align*}
where $T_{\mathcal{F}_0}$ is a coherent subsheaf of the tangent sheaf $T_X:=\mathscr{H}om_{\mathcal{O}_X}\left(\Omega_{X/k}^1 , \mathcal{O}_X \right)$ such that $\left[ T_{\mathcal{F}_0}, T_{\mathcal{F}_0} \right]\subset T_{\mathcal{F}_0}$ under the Lie bracket and   $T_X/T_{\mathcal{F}_0}$ is a torsion-free $\mathcal{O}_X$-module, so that $T_{\mathcal{F}_0}$ is reflexive.  A cartesian diagram of morphism of algebraic schemes
\begin{center}
$\eta:\begin{CD}
\left(X, T_{\mathcal{F}_0} \right)@>>> \left(\mathcal{X}, T_\mathcal{F} \right) \\
@VVV @VV\pi V\\
\textnormal{Spec}(k)@>s>> S
\end{CD}$
\end{center}
where $\pi$ is flat, smooth and surjective, and $S$ is a connected algebraic scheme, is called a family of (singularly) foliated deformations of $\left(X, T_{\mathcal{F}_0} \right)$ in terms of tangent sheaves parametrized by $S$ if there exists an exact sequence
\begin{align}\label{a1}
0\to  T_\mathcal{F} \to T_{\mathcal{X}/S}:= \mathscr{H}om_{\mathcal{O}_\mathcal{X}}\left(\Omega_{\mathcal{X}/S}^1,\mathcal{O}_\mathcal{X} \right)\to T_{\mathcal{X}/S}/T_\mathcal{F}\to 0
\end{align}
where $\left[T_{\mathcal{F}}, T_{\mathcal{F}}\right] \subset T_{\mathcal{F}}$ and $T_{\mathcal{X}/S}/T_\mathcal{F}$ is flat over $S$, and for each $k$-rational point $t\in S$ the exact sequence $(\ref{a1})$ induces an exact sequence on the fiber $X_t$
\begin{align*}
0\to  T_{\mathcal{F}_t}\to T_{X_t}:=\mathscr{H}om_{\mathcal{O}_{X_t}}\left(\Omega_{X_t/k}^1, \mathcal{O}_{X_t} \right) \to T_{X_t}/T_{\mathcal{F}_t}\to 0
\end{align*}
such that $T_{X_t}/T_{\mathcal{F}_t}$ is a torsion-free $\mathcal{O}_{X_t}$-module. When $S=\textnormal{Spec}(A)$ with $A\in \textnormal{ob}\left(\bold{Art}^*\right)$ and $s\in S$ is the closed point we have a local family of foliated deformations of $\left(X, T_{\mathcal{F}}\right)$ in terms of tangent sheaves over $A$. The foliated deformation $\eta$ in terms of tangent sheaves will be denoted by $\left( S, \eta\right)$ or $\left( A, \eta\right)$ when $S=\textnormal{Spec}(A)$. The local foliated deformation $\left( A, \eta\right)$ is infinitesimal (resp. first-order) if $A\in \bold{Art}$ (resp. $A=k[\epsilon]$ the ring of dual numbers). Given another foliated deformation of $\left( X, T_{\mathcal{F}_0}\right)$ over $S$ in terms of tangent sheaves
\begin{center}
$\xi: \begin{CD}
\left( X, T_{\mathcal{F}_0}\right) @>>>  \left( \mathcal{Y}, T_\mathcal{G}\right)\\
@VVV @VVV\\
\textnormal{Spec}(k) @>>> S
\end{CD}$
\end{center}
An isomorphism between foliated deformations $\eta$ and $\xi$ is an foliated $S$-isomorphism $\phi:\left(\mathcal{X}, T_{\mathcal{F}} \right)\to \left(  \mathcal{Y}, T_{\mathcal{G}} \right)$ with $T_\mathcal{F}=\phi^* T_{\mathcal{G}}$ as a subsheaf of $T_{\mathcal{X}/S}$ inducing the identity on $\left( X, T_{\mathcal{F}_0}\right)$, or equivalently the commutative diagram
\begin{center}
$\begin{CD}
0 @>>> T_\mathcal{F} @>>> T_{\mathcal{X}/S} \\
@. @V\cong VV @| \\
0@>>> \phi^* T_\mathcal{G} @>>> \phi^*T_{\mathcal{Y}/S}=T_{\mathcal{X}/S}
\end{CD}$
\end{center}

By a pointed scheme, we will mean a pair $(S,s)$ where $S$ is a scheme and $s\in S$. We call $\left(S,s\right)$ a $k$-pointed scheme if $k=k(s)$.

\begin{definition}[Trivial foliated deformations in terms of tangent sheaves]
Let $\left( X, T_{\mathcal{F}_0} \right)$ be a foliated nonsingular variety over $k$, and $(S,s)$ be a $k$-pointed scheme. We define the trivial foliated family induced by $\left( X, T_{\mathcal{F}_0} \right)$ and $(S,s)$ to be the following foliated deformation of $\left( X, T_{\mathcal{F}_0}\right)$ over $S$,
\begin{center}
$\begin{CD}
\left( X, T_{\mathcal{F}_0}\right) @>>> \left(  X\times_k S    , \pi_1^* T_{\mathcal{F}_0} \right)\\
@VVV @VV\pi_2 V\\
\textnormal{Spec}(k) @>s>> S
\end{CD}$
\end{center}
where $\pi_1:X\times S\to X$ and $\pi_2:X\times S \to S$ are the projections. A flat foliated deformation of $\left( X, T_{\mathcal{F}_0}\right)$ over $S$ is called trivial if it is isomorphic to the trivial foliated family over $S$.
\end{definition}

\begin{definition}[Rigid foliated deformations in terms of tangent sheaves]\label{a5}
Let $\left( X, T_{\mathcal{F}_0}\right)$ be a foliated nonsingular variety. $\left(X, T_{\mathcal{F}_0}\right)$ is called rigid in foliated deformations in terms of tangent sheaves if every infinitesimal foliated deformation of $\left( X, T_{\mathcal{F}_0}\right)$ over $A$ is trivial for every $A \in \bold{Art}$.

\end{definition}

\section{Foliated deformations of $\left(X, T_{\mathcal{F}_0} \right)$ with $T_{\mathcal{F}_0}$ locally free such that the underlying deformations of $(X, T_{\mathcal{F}_0})$ are trivial}

Let $\left(X, T_{\mathcal{F}_0}\right)$ be a foliated nonsingular variety such that $T_{\mathcal{F}_0}$ is a locally free subsheaf of $T_X$, so that we have an exact sequence $0\to T_{\mathcal{F}_0}\to T_X \to T_X/ T_{\mathcal{F}_0} \to 0$ where $T_X/T_{\mathcal{F}_0} $ is a torsion-free $\mathcal{O}_X$-module. We start from the most simple foliated deformations in terms of tangent sheaves. In this section, we study flat foliated deformations of $\left(X, T_{\mathcal{F}_0} \right)$ such that the underlying deformations of a pair of a nonsingular variety $X$ and a locally free sheaf $T_{\mathcal{F}_0}$ are trivial.

\begin{definition}
Let $(X, T_{\mathcal{F}_0})$ be a foliated nonsingular variety as above. An infinitesimal foliated deformations over $A\in \bold{Art}$ of $\left(X, T_{\mathcal{F}_0}\right)$  with $T_{\mathcal{F}_0}$ and $X$ fixed is a locally free subsheaf $T_{\mathcal{F}}$ of $T_{\frac{X\times \textnormal{Spec}(A)}{\textnormal{Spec}(A)}}=T_X\otimes A$ such that $T_{\mathcal{F}} \cong T_{\mathcal{F}_0}\otimes A$ and $\left[ T_{\mathcal{F}}, T_{\mathcal{F}} \right] \subset T_{\mathcal{F}}$ and $\frac{T_{X\times \textnormal{Spec}(A)}}{T_\mathcal{F}}$ is flat over $A$. Equivalently an infinitesimal foliated deformations over $A$ of $\left(X, T_{\mathcal{F}_0}\right)$  with $T_{\mathcal{F}_0}$ and $X$ fixed is an exact sequence on $\mathcal{X}= X\times_k \textnormal{Spec}(A)$
\begin{align}
0\to T_{\mathcal{F}_0}\otimes A \xrightarrow{l} T_{\mathcal{X}/\textnormal{Spec}(A)}=T_X\otimes A \to \frac{T_X\otimes A}{ l\left(T_{\mathcal{F}_0}\otimes A\right)}\to 0
\end{align}
which induces $0\to T_{\mathcal{F}_0} \to T_X \to \frac{T_X}{T_{\mathcal{F}_0}}\to 0$ over $\textnormal{Spec}(k)$ such that $\left[ l\left(T_{\mathcal{F}_0}\otimes A\right) , l\left(T_{\mathcal{F}_0}\otimes A\right)\right]\subset l\left( T_{\mathcal{F}_0}\otimes A\right)$ under the Lie bracket, and $T_{\mathcal{X}/A}/l\left(T_{\mathcal{F}_0} \otimes A \right)$ is flat over $\textnormal{Spec}(A)$. Then we have a functor of Artin rings
\begin{align*}
\textnormal{FDef}_{\left(X, T_{\mathcal{F}_0}\right)_{\textnormal{fixed}}}:\bold{Art}&\to \textnormal{Sets}\\
                                                                              A &\mapsto        \{\textnormal{infinitesimal foliated deformations of $(X, T_{\mathcal{F}_0})$ over $A$  with $T_{\mathcal{F}_0}$ and $X$ fixed} \}
\end{align*}
\end{definition}

Then we will show that 
{\small{\begin{align*}
\mathscr{H}om_{\mathcal{O}_X}\left(T_{\mathcal{F}_0}, T_X \right)^\bullet: \mathscr{H}om_{\mathcal{O}_X}\left(T_{\mathcal{F}_0}, T_X \right) \xrightarrow{D_1'}\mathscr{H}om_{\mathcal{O}_X}\left(\bigwedge^2 T_{\mathcal{F}_0}, T_X/T_{\mathcal{F}_0} \right) \xrightarrow{D_2} \mathscr{H}om_{\mathcal{O}_X}\left(\bigwedge^3 T_{\mathcal{F}_0}, T_X/T_{\mathcal{F}_0}\right) \xrightarrow{D_3} \cdots
\end{align*}}}
controls infinitesimal foliated deformations of  $\left( X, T_{\mathcal{F}_0} \right)$ with $X$ and $T_{\mathcal{F}_0}$ fixed. Here $D_1'$ is defined by the composition of the natural map $\mathscr{H}om_{\mathcal{O}_X}\left( T_{\mathcal{F}_0}, T_X  \right)\to \mathscr{H}om_{\mathcal{O}_M}\left( T_{\mathcal{F}_0}, \frac{T_X}{T_{\mathcal{F}_0}}\right)$ induced by the quotient $T_X\to \frac{T_X}{T_{\mathcal{F}_0}}$ with $D_1$ (see Part I). We will denote the $i$-th cohomology group of $\mathscr{H}om_{\mathcal{O}_X}\left(T_{\mathcal{F}_0}, T_X \right)^\bullet$ by $\mathbb{H}^i\left( X, \mathscr{H}om_{\mathcal{O}_X}\left(T_{\mathcal{F}_0}, T_X \right)^\bullet\right)$.

\begin{proposition}\label{at3}
Let $\left( X, T_{\mathcal{F}_0} \right)$ be a foliated nonsingular variety over $k$ with $T_{\mathcal{F}_0}$ locally free. Then
\begin{enumerate}
\item there is a $1-1$ correspondence 
\begin{align*}
\textnormal{FDef}_{\left( X, T_{\mathcal{F}_0}\right)_{\textnormal{fixed}}} \left(k[\epsilon] \right) \cong \mathbb{H}^0\left( X, \mathscr{H}om_{\mathcal{O}_X}\left(T_{\mathcal{F}_0}, T_X \right)^\bullet \right)
\end{align*}
\item Let $A\in \bold{Art}$ and $\eta=\left( \mathcal{X}, T_\mathcal{F} \right)$ be an infinitesimal foliated deformations over $A$ of $\left( X, T_{\mathcal{F}_0} \right)$ with $X$ and $T_{\mathcal{F}_0}$ fixed. To every small extension $0\to (t)\to \tilde{A} \to A \to 0$, we can associate an element $o_\eta(e)\in \mathbb{H}^1\left( X, \mathscr{H}om_{\mathcal{O}_X}\left( T_{\mathcal{F}_0}, T_X \right)^\bullet \right)$ called the obstruction lifting of $\eta$ to $\tilde{A}$, which is $0$ if and only if a lifting of $\eta$ to $\tilde{A}$ exists.
\end{enumerate}
\end{proposition}

\begin{proof}

Given a first-order foliated deformation $\left(X\times \textnormal{Spec}\left(k[\epsilon]\right), T_{\mathcal{F}} \right)$ of $\left(X, T_{\mathcal{F}_0} \right)$ such that the underlying deformation of $\left(X, T_{\mathcal{F}_0}\right)$ is trivial, we choose an affine open covering $\mathcal{U}=\{U_i\}$ of $X$ such that $\Gamma\left( U_i, T_{\mathcal{F}_0}\right)$ is generated by $T_{0i}^1,..., T_{0i}^p \in \Gamma\left( U_i, T_X \right)$ with $T_{0i}^\alpha =\sum_{\beta=1}^p r_{0ij}^{\alpha\beta}T_{0j}^\beta$ for some $r_{0ij}^{\alpha\beta}\in \Gamma\left( U_{ij}, \mathcal{O}_X\right),\alpha=1,...,p$ and $\left[ T_{0i}^\alpha, T_{0i}^\beta \right] = \sum_{\gamma=1}^p g_{0i\alpha\beta}^\gamma T_{0i}^\gamma$ for some $g_{0i\alpha\beta}^\gamma \in \Gamma\left(U_i, \mathcal{O}_X \right), \alpha,\beta=1,...,p$, and $\Gamma\left( U_i, T_{\mathcal{F}} \right)$ is generated by $T_{0i}^\alpha+ \epsilon T_{i}'^\alpha$ for some $T_{i}'^\alpha\in \Gamma\left( U_i, T_X\right)$ with $T_{0i}^\alpha+ \epsilon T_{i}'^\alpha= \sum_{\beta=1}^p r_{0ij}^{\alpha\beta}\left(T_{0j}^\beta + \epsilon T_{j}'^\beta \right)$ and $\left[ T_{0i}^\alpha + \epsilon T_{i}'^\alpha, T_{0i}^\beta + \epsilon  T_{i}'^\beta \right] = \sum_{\gamma=1}^p \left( g_{0i\alpha\beta}^\gamma + \epsilon c_{i\alpha\beta}^\gamma \right)\left(T_{0i}^\gamma + \epsilon T_i'^\gamma \right)$ for some $c_{i\alpha\beta}^\gamma \in \Gamma\left( U_i, \mathcal{O}_X \right)$. Then we have
\begin{align}
&T_i'^\alpha = \sum_{\beta=1}^p r_{0ij}^{\alpha\beta} T_j'^\beta \label{at1}\\
\left[ T_{0i}^\alpha, T_i'^\beta\right]& + \left[ T_i'^\alpha, T_{0i}^\beta\right] = \sum_{\gamma=1}^p g_{0i\alpha\beta}^\gamma T_i'^\gamma + \sum_{\gamma=1}^p c_{i\alpha\beta}^\gamma T_{0i}^\gamma \label{at2}
\end{align}
We define $T_i'\in \Gamma\left( U_i, \mathscr{H}om_{\mathcal{O}_X}\left( T_{\mathcal{F}_0}, T_X \right) \right)$ in the following way:
\begin{align*}
T_i':\Gamma\left( U_i, T_{\mathcal{F}_0} \right) &\to \Gamma\left( X, T_X\right)\\
 T_{0i}^\alpha &\mapsto  T_{i}'^\alpha
\end{align*}
and linearly extends to $\Gamma\left(U_i, T_{\mathcal{F}_0} \right)$. Then $(\ref{at1})$ implies that $T':=\{T_i'\} \in H^0\left(X, \mathscr{H}om_{\mathcal{O}_X}\left(T_{\mathcal{F}_0},T_X  \right) \right)$ and $(\ref{at2})$ implies that $D_1'(T_i')=0$, so that $T'\in \mathbb{H}^0\left(X, \mathscr{H}om_{\mathcal{O}_X}\left(T_{\mathcal{F}_0}, T_X \right)^\bullet \right)$. 

We define an inverse map. Given a $T'\in \mathbb{H}^0\left(X, \mathscr{H}om_{\mathcal{O}_X}\left(T_{\mathcal{F}_0}, T_X \right)^\bullet \right)$, we have $T_i':\Gamma\left(U_i, T_{\mathcal{F}_0}\right)\to \Gamma\left( U_i, T_X\right), T_{0i}^\alpha \mapsto T_i'^\alpha$ and $\left[ T_{0i}^\alpha, T_i'^\beta \right] - \left[ T_{0i}^\beta, T_i'^\alpha \right]- \sum_{\gamma=1}^p g_{0i\alpha\beta}^\gamma T_i'^\gamma = \sum_{\eta=1}^p c_{i\alpha\beta}^\gamma T_{0i}^\gamma \in \Gamma\left( U_i, T_{\mathcal{F}_0} \right)$ for some $c_{i\alpha\beta}^\gamma T_{0i}^\gamma$. By reversing the above process, we obtain a first-order foliated deformation of $\left(X, T_{\mathcal{F}_0}\right)$ with $X$ and $T_{\mathcal{F}_0}$ fixed.

Next we identify obstructions. Given a small extension $0\to (t)\to \tilde{A} \to A \to 0$ in $\bold{Art}$ and an infinitesimal foliated deformation $\eta=\left(\mathcal{X}, T_{\mathcal{F}} \right)$ over $A$ of $\left( X, T_{\mathcal{F}_0}\right)$ with $X$ and $T_{\mathcal{F}_0}$ fixed, we choose an affine open covering $\mathcal{U}=\{U_i\}$ of $X$ such that $\Gamma\left(U_i, T_\mathcal{F} \right)$ is generated by $T_i^1,..., T_i^p \in \Gamma\left( U_i, T_X \right)\otimes A$ which induces $T_{0i}^\alpha$ by modulo $\mathfrak{m}_A$ and $T_i^\alpha = \sum_{\beta=1}^p r_{0ij}^{\alpha\beta} T_j^\beta$ and $\left[ T_i^\alpha, T_i^\beta \right] = \sum_{\gamma=1}^p g_{i\alpha\beta}^\gamma T_i^\gamma$ for some $g_{i\alpha\beta}^\gamma \in \Gamma\left(  U_i, \mathcal{O}_X \right)\otimes A$ which induces $g_{0i\alpha\beta}^\gamma$ by modulo $\mathfrak{m}_A$. In order to see if there exists a lifting of $\left(\mathcal{X}, T_\mathcal{F} \right)$ to $\tilde{A}$, we choose an arbitrary lifting $\tilde{T}_i^\alpha\in \Gamma\left(U_i, T_X \right)\otimes \tilde{A}$ of $T_i^\alpha$ and arbitrary lifting $\tilde{g}_{i\alpha\beta}^\gamma \in \Gamma\left( U_i, \mathcal{O}_X \right) \otimes \tilde{A}$ of $g_{i\alpha\beta}^\gamma$ for $\alpha,\beta=1,...,p$. Then we have
\begin{align*}
\tilde{T}_i^\alpha-\sum_{\beta=1}^p r_{0ij}^{\alpha\beta} \tilde{T}_j^\beta= t W_{ij}^{\alpha}\,\,\,\,\,\,\,\,\,\,\text{for some}\,\,\,  W_{ij}^{\alpha}\in \Gamma(U_{ij}, T_X)\\
\left[\tilde{T}_i^\alpha, \tilde{T}_i^\beta\right]-\sum_{\eta=1}^p \tilde{g}_{i\alpha\beta}^{\eta} \tilde{T}_i^\eta = t P_i^{\alpha\beta}\,\,\,\,\,\,\,\,\textnormal{for some}\,\,\,\,\,P_i^{\alpha\beta}\in \Gamma(U_i, T_X)
\end{align*}
Then we define $W_{ij}\in \Gamma\left( U_{ij}, \mathscr{H}om_{\mathcal{O}_X}\left( T_{\mathcal{F}_0}, T_X \right) \right)$ in the following way:
\begin{align*}
W_{ij}:\Gamma\left( U_{ij}, T_{\mathcal{F}_0} \right) &\to \Gamma\left( U_i, T_X \right) \\
 T_{0i}^\alpha &\mapsto  W_{ij}^\alpha
\end{align*}
and linearly extends to $\Gamma\left( U_{ij}, T_{\mathcal{F}_0}\right)$, and we define $P_i\in \Gamma\left( U_i, \mathscr{H}om_{\mathcal{O}_X}\left( \bigwedge^2 T_{\mathcal{F}_0}, \frac{T_X}{T_{\mathcal{F}_0}} \right) \right)$ in the following way:
\begin{align*}
P_i:\Gamma\left( U_i, \bigwedge^2 T_{\mathcal{F}_0} \right) &\to \Gamma\left( U_i, \frac{T_X}{T_{\mathcal{F}_0}} \right) \\
T_{0i}^\alpha \wedge T_{0i}^\beta &\mapsto \overline{P_i^{\alpha\beta}}
\end{align*}
and linearly extends to $\Gamma\left(U_i, \bigwedge^2 T_{\mathcal{F}_0} \right)$. Here $\overline{P_i^{\alpha\beta}}$ is the image of $P_i^{\alpha\beta}\in \Gamma\left( U_i, T_X \right)$ in $\Gamma\left( U_i, \frac{T_X}{T_{\mathcal{F}_0}} \right)$. Then we claim that $\left( \left\{ P_i\right\}, \left\{ -W_{ij} \right\} \right)\in C^0\left(\mathcal{U}, \mathscr{H}om_{\mathcal{O}_X}\left( \bigwedge^2 T_{\mathcal{F}_0}, \frac{ T_X}{T_{\mathcal{F}_0}} \right) \right) \bigoplus C^1\left(\mathcal{U}, \mathscr{H}om_{\mathcal{O}_X}\left( T_{\mathcal{F}_0}, T_X \right) \right)$ defines a $1$-cocycle in the following \v Cech resolution of $\mathscr{H}om_{\mathcal{O}_X}\left( T_{\mathcal{F}_0}, T_X \right)^\bullet$:
\begin{center}
$\begin{CD}
\cdots \\
@AD_2AA\\
C^0\left( \mathcal{U}, \mathscr{H}om_{\mathcal{O}_X}\left( \bigwedge^2 T_{\mathcal{F}_0}, \frac{T_X}{T_{\mathcal{F}_0}} \right) \right) @>-\delta>> \cdots \\
@AD_1'AA @AD_1' AA \\
C^0\left(\mathcal{U}, \mathscr{H}om_{\mathcal{O}_X}\left(T_{\mathcal{F}_0}, T_X \right) \right) @>\delta>> C^1\left( \mathcal{U}, \mathscr{H}om_{\mathcal{O}_X}\left(  T_{\mathcal{F}_0}, T_X \right) \right) @>-\delta>> 
\end{CD}$
\end{center}

We note that
\begin{align*}
tW_{ij}^\alpha-tW_{ik}^\alpha+t\sum_{\beta=1}^p r_{0ij}^{\alpha\beta} W_{jk}^\beta&=\tilde{T}_i^\alpha-\sum_{\beta=1}^p r_{0ij}^{\alpha\beta} \tilde{T}_j^\beta -\tilde{T}_i^\alpha+\sum_{\beta=1}^p r_{0ik}^{\alpha\beta} \tilde{T}_k^\beta+\sum_{\beta=1}^p r_{0ij}^{\alpha\beta} \tilde{T}_j^\beta-\sum_{\beta,\gamma=1}^p r_{0ij}^{\alpha\beta} r_{0jk}^{\beta\gamma} \tilde{T}_k^\gamma=0
\end{align*}
This implies that $-\delta\left(\left\{ -W_{ij} \right\} \right)=0$. We note that
\begin{align*}
0&=\left[T_i^\alpha, T_i^\beta \right]-\sum_{\eta=1}^p g_{i\alpha\beta}^\eta T_i^\eta=\sum_{\gamma,\eta=1}^p \left[ r_{0ij}^{\alpha\gamma} T_j^\gamma, r_{0ij}^{\beta\eta} T_j^\eta \right]-\sum_{\eta,\gamma=1}^p g_{i\alpha\beta}^\eta r_{0ij}^{\eta\gamma} T_j^\gamma\\
&= \sum_{\gamma,\eta=1}^p r_{0ij}^{\alpha\gamma} r_{0ij}^{\beta\eta}\left[T_j^\gamma, T_j^\eta \right] +\sum_{\gamma,\eta=1}^p r_{0ij}^{\alpha\gamma} T_j^\gamma\left(r_{0ij}^{\beta\eta}\right) T_j^\eta -\sum_{\gamma,\eta=1}^p r_{0ij}^{\beta\eta} T_j^\eta\left(r_{0ij}^{\alpha\gamma}\right) T_j^\gamma-\sum_{\eta,\gamma=1}^p g_{i\alpha\beta}^\eta r_{0ij}^{\eta\gamma} T_j^\gamma\\
&=\sum_{\xi=1}^p\left(\sum_{\gamma,\eta=1}^p r_{0ij}^{\alpha\gamma}r_{0ij}^{\beta\eta} g_{j\gamma\eta}^\xi+\sum_{\gamma=1}^p r_{0ij}^{\alpha\gamma} T_j^\gamma\left(r_{0ij}^{\beta\xi}\right)-\sum_{\eta=1}^p r_{0ij}^{\beta\eta} T_j^\eta\left(r_{0ij}^{\alpha\xi}\right)-\sum_{\eta=1}^p g_{i\alpha\beta}^\eta r_{0ij}^{\eta\xi} \right) T_j^\xi
\end{align*}
so that we have
\begin{align}
\sum_{\gamma,\eta=1}^p r_{0ij}^{\alpha\gamma}r_{0ij}^{\beta\eta} g_{j\gamma\eta}^\xi+\sum_{\gamma=1}^p r_{0ij}^{\alpha\gamma} T_j^\gamma\left(r_{0ij}^{\beta\xi}\right)-\sum_{\eta=1}^p r_{0ij}^{\beta\eta} T_j^\eta \left(r_{0ij}^{\alpha\xi}\right)-\sum_{\eta=1}^p g_{i\alpha\beta}^\eta r_{0ij}^{\eta\xi}=0 \label{at23}
\end{align}
Then we have
\begin{align}\label{at5}
\sum_{\gamma,\eta=1}^p r_{0ij}^{\alpha\gamma} r_{0ij}^{\beta\eta} \tilde{g}_{j\gamma\eta}^\xi+\sum_{\gamma=1}^p r_{0ij}^{\alpha\gamma} \tilde{T}_j^\gamma\left(r_{0ij}^{\beta\xi}\right)-\sum_{\eta=1}^p r_{0ij}^{\beta\eta} \tilde{T}_j^\eta\left(r_{0ij}^{\alpha\xi}\right)-\sum_{\eta=1}^p \tilde{g}_{i\alpha\beta}^\eta r_{0ij}^{\eta\xi} = t e_{ij}^{\alpha\beta\xi}\,\,\,\,\,\textnormal{for some}\,\,\, e_{ij}^{\alpha\beta \xi}\in \Gamma\left( U_i, \mathcal{O}_X \right)
\end{align}
Then from $(\ref{at5})$, we have
{\small{\begin{align*}
tP_i^{\alpha\beta}&=\left[\tilde{T}_i^\alpha, \tilde{T}_i^\beta \right]-\sum_{\eta=1}^p \tilde{g}_{i\alpha\beta}^\eta \tilde{T}_i^\eta=\left[tW_{ij}^\alpha+\sum_{\gamma=1}^p  r_{0ij}^{\alpha\gamma} \tilde{T}_j^\gamma, tW_{ij}^\beta+\sum_{\eta=1}^p r_{0ij}^{\beta\eta} \tilde{T}_j^\eta \right] -\sum_{\eta=1}^p\tilde{g}_{i\alpha\beta}^\eta\left(tW_{ij}^\eta+\sum_{\gamma=1}^p r_{0ij}^{\eta\gamma}\tilde{T}_j^\gamma \right)\\
&=t\left[W_{ij}^\alpha, T_{0i}^\beta \right]+t\left[T_{0i}^\alpha, W_{ij}^\beta \right]-t\sum_{\eta=1}^p g_{0i\alpha\beta}^\eta W_{ij}^\eta +\sum_{\gamma,\eta=1}^p r_{0ij}^{\alpha\gamma} r_{0ij}^{\beta\eta}\left(\left[\tilde{T}_j^\gamma, \tilde{T}_j^\eta \right] -\sum_{\xi=1}^p \tilde{g}_{j\gamma\eta}^\xi \tilde{T}_j^\xi \right) +\sum_{\gamma,\eta,\xi=1}^p r_{0ij}^{\alpha\gamma}r_{0ij}^{\beta\eta} \tilde{g}_{j\gamma\eta}^\xi \tilde{T}_j^\xi\\
&+\sum_{\gamma,\eta=1}^p r_{0ij}^{\alpha\gamma}\tilde{T}_j^\gamma \left( r_{0ij}^{\beta\eta}\right)\tilde{T}_j^\eta - \sum_{\gamma,\eta=1}^pr_{0ij}^{\beta\eta} \tilde{T}_j^\eta\left( r_{0ij}^{\alpha\gamma} \right) \tilde{T}_j^\gamma -\sum_{\eta,\gamma=1}^p \tilde{g}_{i\alpha\beta}^\eta r_{0ij}^{\eta\gamma} \tilde{T}_j^\gamma\\
&=t\left(\left[W_{ij}^\alpha, T_{0i}^\beta \right]+\left[ T_{0i}^\alpha, W_{ij}^\beta \right]-\sum_{\eta=1}^p g_{0i\alpha\beta}^\eta W_{ij}^\eta +\sum_{\gamma,\eta=1}^p r_{0ij}^{\alpha\gamma}r_{0ij}^{\beta\eta} P_j^{\gamma\eta}\right)\\
&+\sum_{\xi=1}^p\left(\sum_{\gamma,\eta=1}^p r_{0ij}^{\alpha\gamma} r_{0ij}^{\beta\eta} \tilde{g}_{j\gamma\eta}^\xi+\sum_{\gamma=1}^p r_{0ij}^{\alpha\gamma} \tilde{T}_j^\gamma\left( r_{0ij}^{\beta\xi} \right) - \sum_{\eta=1}^p r_{0ij}^{\beta\eta} \tilde{T}_j^\eta\left(r_{0ij}^{\alpha\xi} \right) -\sum_{\eta=1}^p \tilde{g}_{i\alpha\beta}^\eta r_{0ij}^{\eta\xi}\right) \tilde{T}_j^\xi\\
&=t\left(\left[W_{ij}^\alpha, T_{0i}^\beta \right]+\left[ T_{0i}^\alpha, W_{ij}^\beta \right]-\sum_{\eta=1}^p g_{0i\alpha\beta}^\eta W_{ij}^\eta +\sum_{\gamma,\eta=1}^p r_{0ij}^{\alpha\gamma}r_{0ij}^{\beta\eta} P_j^{\gamma\eta}\right)+t\sum_{\xi=1}^p e_{ij}^{\alpha\beta\xi} T_{0j}^\xi
\end{align*}}}
This implies that $-\delta(\{P_i\})+D_1(\{ -W_{ij}\})=0$. On the other hand, we note that
{\small{\begin{align*}
0&=\left[T_i^\gamma, \left[T_i^\alpha, T_i^\beta \right]\right]- \sum_{\eta=1}^p \left[T_i^\gamma, g_{i\alpha\beta}^\eta T_i^\eta\right]=\left[\left[T_i^\gamma, T_i^\alpha\right], T_i^\beta \right]- \left[\left[T_i^\gamma, T_i^\beta\right], T_i^\alpha\right]-\sum_{\eta=1}^p T_i^\gamma\left(g_{i\alpha\beta}^\eta\right)T_i^\eta -\sum_{\eta=1}^p g_{i\alpha\beta}^\eta \left[T_i^\gamma, T_i^\eta\right]\\
&=\left[\sum_{\eta=1}^p g_{i\gamma\alpha}^\eta T_i^\eta, T_i^\beta \right] -\left[ \sum_{\eta=1}^p g_{i\gamma\beta}^\eta T_i^\eta, T_i^\alpha \right] -\sum_{\eta=1}^p T_i^\gamma\left(g_{i\alpha\beta}^\eta \right)T_i^\eta -\sum_{\eta,\xi=1}^p g_{i\alpha\beta}^\eta g_{i\gamma\eta}^\xi T_i^\xi\\
&= \sum_{\eta,\xi=1}^p g_{i\gamma\alpha}^\eta g_{i\eta\beta}^\xi T_i^\xi -\sum_{\eta=1}^p T_i^\beta \left(g_{i\gamma\alpha}^\eta\right) T_i^\eta-\sum_{\eta,\xi=1}^p g_{i\gamma\beta}^\eta g_{i\eta\alpha}^\xi T_i^\xi +\sum_{\eta=1}^p T_i^\alpha\left(g_{i\gamma\beta}^\eta\right)T_i^\eta - \sum_{\eta=1}^p T_i^\gamma\left(g_{i\alpha\beta}^\eta\right)T_i^\eta-\sum_{\eta,\xi=1}^p g_{i\alpha\beta}^\eta g_{i\gamma\eta}^\xi T_i^\xi\\
&= \sum_{\xi=1}^p\left(\sum_{\eta=1}^p g_{i\gamma\alpha}^\eta g_{i\eta\beta}^\xi- T_i^\beta\left(g_{i\gamma\alpha}^\xi\right)-\sum_{\eta=1}^p g_{i\gamma\beta}^\eta g_{i\eta\alpha}^\xi+ T_i^\alpha\left(g_{i\gamma\beta}^\xi\right) - T_i^\gamma\left(g_{i\alpha\beta}^\xi\right)-\sum_{\eta=1}^p g_{i\alpha\beta}^\eta g_{i\gamma\eta}^\xi \right) T_i^\xi
\end{align*}}}
Then we have
\begin{align*}
\sum_{\eta=1}^p g_{i\gamma\alpha}^\eta g_{i\eta\beta}^\xi- T_i^\beta\left(g_{i\gamma\alpha}^\xi \right)-\sum_{\eta=1}^p g_{i\gamma\beta}^\eta g_{i\eta\alpha}^\xi+ T_i^\alpha\left(g_{i\gamma\beta}^\xi \right) - T_i^\gamma\left(g_{i\alpha\beta}^\xi\right)-\sum_{\eta=1}^p g_{i\alpha\beta}^\eta g_{i\gamma\eta}^\xi =0,
\end{align*}
so that we have
\begin{align}\label{at4}
\sum_{\eta=1}^p \tilde{g}_{i\gamma\alpha}^\eta \tilde{g}_{i\eta\beta}^\xi- \tilde{T}_i^\beta\left(\tilde{g}_{i\gamma\alpha}^\xi\right)-\sum_{\eta=1}^p \tilde{g}_{i\gamma\beta}^\eta \tilde{g}_{i\eta\alpha}^\xi+ \tilde{T}_i^\alpha\left(\tilde{g}_{i\gamma\beta}^\xi\right) - \tilde{T}_i^\gamma\left(\tilde{g}_{i\alpha\beta}^\xi\right)-\sum_{\eta=1}^p \tilde{g}_{i\alpha\beta}^\eta \tilde{g}_{i\gamma\eta}^\xi =tA_{i\alpha\beta\gamma}^\xi\,\,\,\textnormal{for some}\,\,\,\, A_{i\alpha\beta\gamma}^\xi \in \Gamma\left( U_i, \mathcal{O}_X  \right)
\end{align}

Then from $(\ref{at4})$, we have
{\small{\begin{align*}
&t\left[T_{0i}^\gamma, P_i^{\alpha\beta}\right]=\left[\tilde{T}_i^\gamma, \left[\tilde{T}_i^\alpha, \tilde{T}_i^\beta \right] \right]-\sum_{\eta=1}^p \left[\tilde{T}_i^\gamma, \tilde{g}_{i\alpha\beta}^\eta \tilde{T}_i^\eta \right]= \left[\left[\tilde{T}_i^\gamma, \tilde{T}_i^\alpha\right], \tilde{T}_i^\beta\right]-\left[\left[\tilde{T}_i^\gamma, \tilde{T}_i^\beta\right],\tilde{T}_i^\alpha\right]-\sum_{\eta=1}^p \tilde{T}_i^\gamma\left(\tilde{g}_{i\alpha\beta}^\eta\right)\tilde{T}_i^\eta-\sum_{\eta=1}^p \tilde{g}_{i\alpha\beta}^\eta\left[\tilde{T}_i^\gamma, \tilde{T}_i^\eta\right]\\
 &=\left[tP_i^{\gamma\alpha}+\sum_{\eta=1}^p \tilde{g}_{i\gamma\alpha}^\eta \tilde{T}_i^\eta, \tilde{T}_i^\beta \right] - \left[tP_i^{\gamma\beta}+\sum_{\eta=1}^p \tilde{g}_{i\gamma\beta}^\eta \tilde{T}_i^\eta, \tilde{T}_i^\alpha \right] - \sum_{\eta=1}^p \tilde{T}_i^\gamma \left(\tilde{g}_{i\alpha\beta}^\eta \right)\tilde{T}_i^\eta-\sum_{\eta=1}^p \tilde{g}_{i\alpha\beta}^\eta\left(tP_i^{\gamma\eta}+\sum_{\xi=1}^p \tilde{g}_{i\gamma\eta}^\xi \tilde{T}_i^\xi \right)\\
 &=t\left[P_i^{\gamma\alpha} , T_{0i}^\beta \right] +\sum_{\eta=1}^p \tilde{g}_{i\gamma\alpha}^\eta\left(tP_i^{\eta\beta}+\sum_{\xi=1}^p\tilde{g}_{i\eta\beta}^\xi \tilde{T}_i^\xi \right)-\sum_{\eta=1}^p \tilde{T}_i^\beta\left( \tilde{g}_{i\gamma\alpha}^\eta\right) \tilde{T}_i^\eta-t\left[P_i^{\gamma\beta}, T_{0i}^\alpha \right] - \sum_{\eta=1}^p \tilde{g}_{i\gamma\beta}^\eta\left( tP_i^{\eta\alpha}+\sum_{\xi=1}^p \tilde{g}_{i\eta\alpha}^\xi \tilde{T}_i^\xi\right)+\sum_{\eta=1}^p \tilde{T}_i^\alpha\left(\tilde{g}_{i\gamma\beta}^\eta \right) \tilde{T}_i^\eta\\
 &-\sum_{\eta=1}^p \tilde{T}_i^\gamma\left( \tilde{g}_{i\alpha\beta}^\eta \right) \tilde{T}_i^\eta - t\sum_{\eta=1}^p g_{0i\alpha\beta}^\eta P_i^{\gamma\eta} -\sum_{\eta,\xi=1}^p \tilde{g}_{i\alpha\beta}^\eta \tilde{g}_{i\gamma\eta}^\xi \tilde{T}_i^\xi\\
 &=t\left(\left[ P_i^{\gamma\alpha}, T_{0i}^\beta\right] +\sum_{\eta=1}^p g_{0i\gamma\alpha}^\eta P_i^{\eta\beta}-\left[ P_i^{\gamma\beta}, T_{0i}^\alpha\right] -\sum_{\eta=1}^p g_{0i\gamma\beta}^\eta P_i^{\eta\alpha}-\sum_{\eta=1}^p g_{0i\alpha\beta}^\eta P_i^{\gamma\eta}\right) \\
 &+\sum_{\eta,\xi=1}^p \tilde{g}_{i\gamma\alpha}^\eta \tilde{g}_{i\eta\beta}^\xi T_i^\xi - \sum_{\eta=1}^p \tilde{T}_i^\beta\left( \tilde{g}_{i\gamma\alpha}^\eta \right) \tilde{T}_i^\eta -\sum_{\eta,\xi=1}^p \tilde{g}_{i\gamma\beta}^\eta \tilde{g}_{i\eta\alpha}^\xi \tilde{T}_i^\xi+\sum_{\eta=1}^p \tilde{T}_i^\alpha\left( \tilde{g}_{i\gamma\beta}^\eta \right) \tilde{T}_i^\eta - \sum_{\eta=1}^p \tilde{T}_i^\gamma\left( \tilde{g}_{i\alpha\beta}^\eta \right) \tilde{T}_i^\eta  - \sum_{\eta,\xi=1}^p \tilde{g}_{i\alpha\beta}^\eta \tilde{g}_{i\gamma\eta}^{\xi} \tilde{T}_i^{\xi}\\
  &=t\left(\left[ P_i^{\gamma\alpha}, T_{0i}^\beta\right] +\sum_{\eta=1}^p g_{0i\gamma\alpha}^\eta P_i^{\eta\beta}-\left[ P_i^{\gamma\beta}, T_{0i}^\alpha\right] -\sum_{\eta=1}^p g_{0i\gamma\beta}^\eta P_i^{\eta\alpha}-\sum_{\eta=1}^p g_{0i\alpha\beta}^\eta P_i^{\gamma\eta}\right) \\
 &+\sum_{\xi=1}^p \left(\sum_{\eta=1}^p \tilde{g}_{i\gamma\alpha}^\eta \tilde{g}_{i\eta\beta}^\xi -\tilde{T}_i^\beta\left( \tilde{g}_{i\gamma\alpha}^\xi \right)  -\sum_{\eta=1}^p \tilde{g}_{i\gamma\beta}^\eta \tilde{g}_{i\eta\alpha}^\xi +\tilde{T}_i^\alpha\left( \tilde{g}_{i\gamma\beta}^\xi \right) - \tilde{T}_i^\gamma\left( \tilde{g}_{i\alpha\beta}^\xi \right)   - \sum_{\eta,\xi=1}^p \tilde{g}_{i\alpha\beta}^\eta \tilde{g}_{i\gamma\eta}^{\xi}\right) \tilde{T}_i^\xi\\
   &=t\left(\left[ P_i^{\gamma\alpha}, T_{0i}^\beta\right] +\sum_{\eta=1}^p g_{0i\gamma\alpha}^\eta P_i^{\eta\beta}-\left[ P_i^{\gamma\beta}, T_{0i}^\alpha\right] -\sum_{\eta=1}^p g_{0i\gamma\beta}^\eta P_i^{\eta\alpha}-\sum_{\eta=1}^p g_{0i\alpha\beta}^\eta P_i^{\gamma\eta}\right) + t\sum_{\xi=1}^p A_{i\alpha\beta\gamma}^\xi T_{0i}^\xi
\end{align*}}}
This implies that $D_2\left( \left\{ P_i\right\} \right)=0$. This proves the claim. We show that for another lifting $\tilde{T}_i'^\alpha$ of $T_i^\alpha,\alpha=1,...p$, and $\tilde{g}_{i\alpha\beta}^\gamma$ of $g_{i\alpha\beta}^\gamma, \alpha,\beta=1,...p$, the associated $1$-cocycle $\beta:=\left( \left\{ P_i'\right\}, \left\{ - W_{ij}' \right\} \right)$ is cohomologous to $\alpha:=\left( \left\{ P_i \right\}, \left\{ - W_{ij} \right\} \right)$ as above. Then $\tilde{T}_i'^\alpha= \tilde{T}_i^\alpha + t Y_i^\alpha $ for some $Y_i^\alpha \in \Gamma\left(U_i, T_X \right)$ and $\tilde{g}_{i\alpha\beta}'^\gamma= \tilde{g}_{i\alpha\beta}^\gamma + t e_{i\alpha\beta}^\gamma$ for some $e_{i\alpha\beta}^\gamma \in \Gamma\left(U_i, \mathcal{O}_X \right),\alpha,\beta=1,...,p$. We define $Y_i\in \Gamma\left(U_i, \mathscr{H}om_{\mathcal{O}_X}\left(T_{\mathcal{F}_0}, T_X \right) \right)$ by
\begin{align}\label{at10}
Y_i:\Gamma\left( U_i, T_{\mathcal{F}_0} \right) &\to \Gamma\left( U_i, T_X\right)\\
 T_{0i}^\alpha &\mapsto Y_i^\alpha \notag
\end{align}
Then we see that
{\small{\begin{align}\label{at11}
t\left(P_i'^{\alpha\beta}-P_i^{\alpha\beta} \right) &= \left[ \tilde{T}_i^\alpha+t Y_i^\alpha, \tilde{T}_i^\beta  + t Y_i^\beta \right] - \sum_{\eta=1}^p \left(\tilde{g}_{i\alpha\beta}^\eta  + t e_{i\alpha\beta}^\gamma \right)\left(\tilde{T}_i^\eta + t Y_i^\eta \right) - \left[ \tilde{T}_i^\alpha, \tilde{T}_i^\beta \right] + \sum_{\eta=1}^p \tilde{g}_{i\alpha\beta}^\eta \tilde{T}_i^\eta\\
&= t\left[ Y_i^\alpha, T_{0i}^\beta \right] + t \left[ T_{0i}^\alpha, Y_i^\beta \right] - t\sum_{\eta=1}^p g_{0i\alpha\beta}^\eta Y_i^\eta - t\sum_{\eta=1}^p e_{i\alpha\beta}^\gamma T_{0i}^\eta \notag
\end{align}}}
This implies that $D_1'\left(\left\{ Y_i \right\} \right)= \left\{ P_i'-P_i \right\}$. On the other hand, we see that
\begin{align*}
t\left(- W_{ij}'^\alpha + W_{ij}^\alpha \right) = -\tilde{T}_i^\alpha - t Y_i^\alpha + \sum_{\beta=1}^p r_{0ij}^{\alpha\beta} \left(\tilde{T}_j^\beta + t Y_j^\beta \right) + \tilde{T}_i^\alpha - \sum_{\beta=1}^p r_{0ij}^{\alpha\beta} \tilde{T}_j^\beta =  t\left( \sum_{\beta=1}^p r_{0ij}^{\alpha\beta} Y_j^\beta- Y_i^\alpha \right)
\end{align*}
This implies that $\delta\left(\left\{ Y_i\right\} \right)=\left\{ -W_{ij}'-\left(- W_{ij}\right) \right\}$. This proves that $\alpha$ is cohomologous to $\beta$. So given a small extension $e:0\to (t)\to \tilde{A} \to A \to 0$ and an infinitesimal foliated deformation $\eta$ of $\left( X, T_{\mathcal{F}_0}\right)$ with $X$ and $T_{\mathcal{F}_0}$ fixed, we can associate an element $o_\eta(e):=$ the cohomology class of $\alpha \in \mathbb{H}^1\left( X, \mathscr{H}om_{\mathcal{O}_X}\left( T_{\mathcal{F}_0}, T_X \right)^\bullet \right)$. We note that $o_\eta(e)=0$ if and only if there exists $\left\{ X_i \right\} \in C^0\left( \mathcal{U}, \mathscr{H}om_{\mathcal{O}_X}\left( T_{\mathcal{F}_0}, T_X \right)\right)$ where $X_i:\Gamma\left(U_i, T_{\mathcal{F}_0} \right)\to \Gamma\left( U_i, T_X \right), T_{0i}^\alpha \mapsto X_i^\alpha$ such that $X_j- X_i= -W_{ij}$ and $D_1'\left(X_i \right)=P_i$, so that
{\small{\begin{align}
&t\left(\sum_{\beta=1}^p r_{0ij}^{\alpha\beta} X_j^\beta - X_i^\alpha\right)= - \tilde{T}_i^\alpha + \sum_{\beta=1}^p r_{0ij}^{\alpha\beta} \tilde{T}_j^\beta \Longrightarrow \left(\tilde{T}_i^\alpha - tX_i^\alpha\right) = \sum_{\beta=1}^p r_{0ij}^{\alpha\beta}\left(\tilde{T}_j^\beta - t X_j^\beta \right)\notag\\
&t\left(\left[ T_{0i}^\alpha, X_i^\beta \right] - \left[ T_{0i}^\beta, X_i^\alpha \right] - \sum_{\gamma=1}^p g_{0i\alpha\beta}^\gamma X_i^\gamma + \sum_{\gamma=1}^p l_{i\alpha\beta}^\gamma  T_{0i}^\gamma \right) =\left[ \tilde{T}_i^\alpha, \tilde{T}_i^\beta \right] - \sum_{\gamma=1}^p \tilde{g}_{i\alpha\beta}^\gamma \tilde{T}_i^\gamma,\,\,\,\,\,\,\,\textnormal{for some}\,\,\, l_{i\alpha\beta}^\gamma \in \Gamma\left( U_i, \mathcal{O}_X \right)\notag\\
&\Longrightarrow \left[ \tilde{T}_i^\alpha- t X_i^\alpha, \tilde{T}_i^\beta - t X_i^\beta \right] - \sum_{\gamma=1}^p\left( \tilde{g}_{i\alpha\beta}^\gamma + t l_{i\alpha\beta}^\gamma  \right)\left(\tilde{T}_i^\gamma - t X_i^\gamma \right)\label{at12}
\end{align}}}
We note that $\left\{\tilde{T}_i^\alpha - t X_i^\alpha\right\},\alpha=1,...,p$ is linearly independent. In fact, if $\sum_{\alpha=1}^p \tilde{a}_i^\alpha \left( \tilde{T}_i^\alpha - t X_i^\alpha \right)=0$ for some $\tilde{a}_i^\alpha \in \Gamma\left( U_i, \mathcal{O}_X \right) \otimes \tilde{A}$, then $\sum_{\alpha=1}^p a_i^\alpha T_i^\alpha=0$ by modulo $(t)$ where $a_i^\alpha$ is the image of $\tilde{a}_i^\alpha$ in $\Gamma\left( U_i, \mathcal{O}_X \right)\otimes A$, so that $a_i^\alpha=0$ in $A$. This implies that $\tilde{a}_i^\alpha= t b_i^\alpha$ for some $b_i^\alpha\in \Gamma\left( U_i, \mathcal{O}_X \right)$. Then $\sum_{\alpha=1}^p tb_i^\alpha\left(\tilde{T}_i^\alpha- t X_i^\alpha \right)= t\sum_{\alpha=1}^p t b_i^\alpha T_{0i}^\alpha=0$, so that $b_i^\alpha=0$. Hence $\tilde{a}_i^\alpha=0$. Hence a lifting $\tilde{\eta}$ of $\eta$ exits. This completes the proof of Proposition \ref{at3}.
\end{proof}

\section{Foliated deformations of $\left(X, T_{\mathcal{F}_0}\right)$ with $T_{\mathcal{F}_0}$ locally free under the trivial deformation of $X$}

Let $\left(X, T_{\mathcal{F}_0} \right)$ be a foliated nonsingular variety over $k$. In the previous section we study flat foliated deformations of $\left(X, T_{\mathcal{F}_0} \right)$ such that the underling flat deformations of $\left(X,T_{\mathcal{F}_0} \right)$ is trivial. In this section, in addition, we vary $T_{\mathcal{F}_0}$ with $X$ fixed.

\begin{definition}
Let $\left( X, T_{\mathcal{F}_0} \right)$ be a foliated nonsingular variety over $k$ with $T_{\mathcal{F}_0}$ locally free, so that we have an exact sequence $0\to T_{\mathcal{F}_0}\to T_X \to \frac{T_X}{T_{\mathcal{F}_0}}\to 0$.  An infinitesimal foliated deformations of $\left(X, T_{\mathcal{F}_0} \right)$ under
 the trivial deformation of $X$ over $A\in \bold{Art}$ is a  locally free subsheaf  $T_{\mathcal{F}}$ of $T_{\mathcal{X}/A}=T_{X\times \textnormal{Spec}(A)}$ such that $\left[ T_{\mathcal{F}}, T_{\mathcal{F}} \right] \subset T_{\mathcal{F}}$ and $T_{\mathcal{X}/A}/T_{\mathcal{F}}$ is flat over $\textnormal{Spec}(A)$ and the exact sequence on $\mathcal{X}:=X\times_{\textnormal{Spec}(k)} \textnormal{Spec}(A)$ 
\begin{align}\label{at2}
0\to T_{\mathcal{F}} \to T_{\mathcal{X}/A}=T_X\otimes A \to \frac{T_X\otimes A}{T_{\mathcal{F}}}\to 0
\end{align}
induces $0\to T_{\mathcal{F}_0}\to T_X \to \frac{T_X}{T_{\mathcal{F}_0}}\to 0$ over $\textnormal{Spec}(k)$. Then we have a functor of Artin rings
\begin{align*}
\textnormal{FDef}_{T_{\mathcal{F}_0}}: \bold{Art} &\to \textnormal{(Sets)} \\
                             A &\to \left\{ \textnormal{Infinitesimal foliated deformation of $\left( X, T_{\mathcal{F}_0}\right)$ with $X$ fixed}   \right\}
\end{align*}
\end{definition}

We will show that the following complex of sheaves truncated from the leaf complex (see Part I)
{\small{\begin{align*}
\mathscr{H}om_{\mathcal{O}_X}\left(T_{\mathcal{F}_0}, T_X/T_{\mathcal{F}_0} \right)^\bullet :\mathscr{H}om_{\mathcal{O}_X}\left(T_{\mathcal{F}_0}, T_X/T_{\mathcal{F}_0} \right) \xrightarrow{D_1} \mathscr{H}om_{\mathcal{O}_X}\left(\bigwedge^2 T_{\mathcal{F}_0}, T_X/T_{\mathcal{F}_0} \right) \xrightarrow{D_2}  \mathscr{H}om_{\mathcal{O}_X}\left( \bigwedge^3 T_{\mathcal{F}_0}, T_X/T_{\mathcal{F}_0}\right) \to \cdots
\end{align*}}}
controls foliated deformations of $\left(X, T_{\mathcal{F}_0} \right)$ under the trivial deformation of $X$. We will denote the $i$-th cohomology group of $\mathscr{H}om_{\mathcal{O}_X}\left(  T_{\mathcal{F}_0}, \frac{T_X}{T_{\mathcal{F}_0}} \right)^\bullet$ by $\mathbb{H}^i\left( X, \mathscr{H}om_{\mathcal{O}_X}\left(  T_{\mathcal{F}_0}, \frac{T_X}{T_{\mathcal{F}_0}} \right)^\bullet\right)$.

\begin{proposition}\label{at8}
Let $\left(X, T_{\mathcal{F}_0} \right)$ be a foliated nonsingular variety with $T_{\mathcal{F}_0}$ locally free.
\begin{enumerate}
\item there is a $1-1$ correspondence
\begin{align*}
\textnormal{FDef}_{T_{\mathcal{F}_0}}\left(k[\epsilon] \right)\cong \mathbb{H}^0\left(X, \mathscr{H}om_{\mathcal{O}_X}\left(  T_{\mathcal{F}_0}, \frac{T_X}{T_{\mathcal{F}_0}} \right)^\bullet \right)
\end{align*}
\item Let $A\in \bold{Art}$ and $\eta=\left( \mathcal{X}, T_\mathcal{F} \right)$ be an infinitesimal foliated deformations over $A$ of $\left( X, T_{\mathcal{F}_0} \right)$ with $X$ fixed. To every small extension $0\to (t)\to \tilde{A} \to A \to 0$, we can associate an element $o_\eta(e)\in \mathbb{H}^1\left( X, \mathscr{H}om_{\mathcal{O}_X}\left( T_{\mathcal{F}_0}, \frac{T_X}{T_{\mathcal{F}_0}} \right)^\bullet \right)$ called the obstruction lifting of $\eta$ to $\tilde{A}$, which is $0$ if and only if a lifting of $\eta$ to $\tilde{A}$ exists.
\end{enumerate}
\end{proposition}

\begin{proof}
Given a first-order foliated deformation $\left(X\times \textnormal{Spec}(k[\epsilon]), T_{\mathcal{F}} \right)$ of $\left(X, T_{\mathcal{F}_0} \right)$ under the trivial deformation of $X$, we choose an affine open covering $\mathcal{U}=\{U_i\}$ of $X$ such that $\Gamma\left( U_i, T_{\mathcal{F}_0}\right)$ is generated by $T_{0i}^1,..., T_{0i}^p \in \Gamma\left( U_i, T_X \right)$ with $T_{0i}^\alpha =\sum_{\beta=1}^p r_{0ij}^{\alpha\beta}T_{0j}^\beta$ for some $r_{0ij}^{\alpha\beta}\in \Gamma\left( U_{ij}, \mathcal{O}_X\right),\alpha=1,...,p$ and $\left[ T_{0i}^\alpha, T_{0i}^\beta \right] = \sum_{\gamma=1}^p g_{0i\alpha\beta}^\gamma T_{0i}^\gamma$ for some $g_{0i\alpha\beta}^\gamma \in \Gamma\left(U_i, \mathcal{O}_X \right), \alpha,\beta=1,...,p$, and $\Gamma\left( U_i, T_{\mathcal{F}} \right)$ is generated by $T_{0i}^\alpha + \epsilon T_i'^\alpha$ for some $T_i'^\alpha\in \Gamma\left(U_i, T_X \right)$ with $T_{0i}^\alpha + \epsilon T_i'^\alpha= \sum_{\beta=1}^p \left(r_{0ij}^{\alpha\beta} + \epsilon b_{ij}^{\alpha\beta} \right) \left(T_{0j}^\beta + \epsilon T_j'^\beta \right)$ for some $b_{ij}^{\alpha\beta}\in \Gamma\left(U_{ij}, \mathcal{O}_X \right)\alpha=1,...,p$ and $\left[ T_{0i}^\alpha + \epsilon T_{i}'^\alpha, T_{0i}^\beta + \epsilon  T_{i}'^\beta \right] = \sum_{\gamma=1}^p \left( g_{0i\alpha\beta}^\gamma + \epsilon c_{i\alpha\beta}^\gamma \right)\left(T_{0i}^\gamma + \epsilon T_i'^\gamma \right)$ for some $c_{i\alpha\beta}^\gamma \in \Gamma\left( U_i, \mathcal{O}_X \right)$. Then we have
\begin{align}
&T_i'^\alpha = \sum_{\beta=1}^p r_{0ij}^{\alpha\beta} T_j'^\beta + \sum_{\beta=1}^p b_{ij}^{\alpha\beta}T_{0j}^\beta \label{at6}\\
\left[ T_{0i}^\alpha, T_i'^\beta\right]& + \left[ T_i'^\alpha, T_{0i}^\beta\right] = \sum_{\gamma=1}^p g_{0i\alpha\beta}^\gamma T_i'^\gamma + \sum_{\gamma=1}^p c_{i\alpha\beta}^\gamma T_{0i}^\gamma \label{at7}
\end{align}
We define $\overline{T_i'}\in \Gamma\left( U_i, \mathscr{H}om_{\mathcal{O}_X}\left( T_{\mathcal{F}_0}, \frac{T_X}{T_{\mathcal{F}_0}} \right) \right)$ in the following way:
\begin{align}\label{at20}
\overline{T_i'}:\Gamma\left( U_i, T_{\mathcal{F}_0} \right) &\to \Gamma\left( X, \frac{T_X}{T_{\mathcal{F}_0}}\right)\\
 T_{0i}^\alpha &\mapsto  \overline{T_{i}'^\alpha}\notag
\end{align}
and linearly extends to $\Gamma\left(U_i, T_{\mathcal{F}_0} \right)$. Here $\overline{T_i'^{\alpha}}$ is the image of $T_i^{\alpha}\in \Gamma\left( U_i, T_X \right)$ in $\Gamma\left( U_i, \frac{T_X}{T_{\mathcal{F}_0}} \right)$. Then $(\ref{at6})$ implies that $T':=\{T_i'\} \in H^0\left(X, \mathscr{H}om_{\mathcal{O}_X}\left(T_{\mathcal{F}_0},\frac{T_X}{T_{\mathcal{F}_0}}  \right) \right)$ and $(\ref{at7})$ implies that $D_1(T_i')=0$, so that $T'\in \mathbb{H}^0\left(X, \mathscr{H}om_{\mathcal{O}_X}\left(T_{\mathcal{F}_0}, \frac{T_X}{T_{\mathcal{F}_0}} \right)^\bullet \right)$. 

We define an inverse map. Given a $T'\in \mathbb{H}^0\left(X, \mathscr{H}om_{\mathcal{O}_X}\left(T_{\mathcal{F}_0}, \frac{T_X}{T_{\mathcal{F}_0}} \right)^\bullet \right)$, we have $T_i':\Gamma\left(U_i, T_{\mathcal{F}_0}\right)\to \Gamma\left( U_i, \frac{T_X}{T_{\mathcal{F}_0}}\right), T_{0i}^\alpha \mapsto \overline{T_i'^\alpha}$ for some $T_i'^\alpha\in \Gamma\left(U_i, T_X\right)$ and $\left[ T_{0i}^\alpha, T_i'^\beta \right] - \left[ T_{0i}^\beta, T_i'^\alpha \right]- \sum_{\gamma=1}^p g_{0i\alpha\beta}^\gamma T_i'^\gamma = \sum_{\eta=1}^p c_{i\alpha\beta}^\gamma T_{0i}^\gamma \in \Gamma\left( U_i, T_{\mathcal{F}_0} \right)$ for some $c_{i\alpha\beta}^\gamma T_{0i}^\gamma$. Then by reversing the above process, we obtain a first-order foliated deformation $\left( \mathcal{X}, T_{\mathcal{F}}\right)$ with $X$ fixed.

Next we identify obstructions. Given a small extension $0\to (t)\to \tilde{A} \to A \to 0$ in $\bold{Art}$ and a foliated deformation $\left(\mathcal{X}, T_{\mathcal{F}} \right)$ over $A$ of $\left( X, T_{\mathcal{F}_0}\right)$ with $X$ fixed, we choose an affine open covering $\mathcal{U}=\{U_i\}$ of $X$ such that $\Gamma\left(U_i, T_\mathcal{F} \right)$ is generated by $T_i^1,..., T_i^p \in \Gamma\left( U_i, T_X \right)\otimes A$ which induces $T_{0i}^\alpha$ by modulo $\mathfrak{m}_A$ and $T_i^\alpha = \sum_{\beta=1}^p r_{ij}^{\alpha\beta} T_j^\beta$ for some $r_{ij}^{\alpha\beta}\in \Gamma\left( U_{ij}, \mathcal{O}_X \right)$ which induces $r_{0ij}^{\alpha\beta}$ modulo $\mathfrak{m}_A$ and $\left[ T_i^\alpha, T_i^\beta \right] = \sum_{\gamma=1}^p g_{i\alpha\beta}^\gamma T_i^\gamma$ for some $g_{i\alpha\beta}^\gamma \in \Gamma\left(  U_i, \mathcal{O}_X \right)\otimes A$ which induces $g_{0i\alpha\beta}^\gamma$ by modulo $\mathfrak{m}_A$. In order to see if there exists a lifting of $\left(\mathcal{X}, T_\mathcal{F} \right)$ to $\tilde{A}$, we choose an arbitrary lifting $\tilde{T}_i^\alpha\in \Gamma\left(U_i, T_X \right)\otimes \tilde{A}$ of $T_i^\alpha$ and arbitrary lifting $\tilde{r}_{ij}^{\alpha\beta}\in \Gamma\left(U_{ij}, \mathcal{O}_X \right)\otimes \tilde{A}$ of $r_{ij}^{\alpha\beta}$ and arbitrary lifting $\tilde{g}_{i\alpha\beta}^\gamma \in \Gamma\left( U_i, \mathcal{O}_X \right) \otimes \tilde{A}$ of $g_{i\alpha\beta}^\gamma$ for $\alpha,\beta=1,...,p$. Then we have
\begin{align*}
\tilde{T}_i^\alpha-\sum_{\beta=1}^p \tilde{r}_{ij}^{\alpha\beta} \tilde{T}_j^\beta= t W_{ij}^{\alpha}\,\,\,\,\,\,\,\,\,\,&\text{for some}\,\,\,  W_{ij}^{\alpha}\in \Gamma(U_{ij}, T_X)\\
\left[\tilde{T}_i^\alpha, \tilde{T}_i^\beta\right]-\sum_{\eta=1}^p \tilde{g}_{i\alpha\beta}^{\eta} \tilde{T}_i^\eta = t P_i^{\alpha\beta}\,\,\,\,\,\,\,\,&\textnormal{for some}\,\,\,\,\,P_i^{\alpha\beta}\in \Gamma(U_i, T_X)\\
\tilde{r}_{ik}^{\alpha\gamma}- \sum_{\beta=1}^p \tilde{r}_{ij}^{\alpha \beta}\tilde{r}_{jk}^{\beta \gamma}= t C_{ijk}^{\alpha \gamma}\,\,\,\,\,\,&\textnormal{for some}\,\,\,\,\, C_{ijk}^{\alpha\gamma}\in \Gamma\left(U_{ijk}, \mathcal{O}_X \right)
\end{align*}
Then we define $W_{ij}\in \Gamma\left( U_{ij}, \mathscr{H}om_{\mathcal{O}_X}\left( T_{\mathcal{F}_0}, \frac{T_X}{T_{\mathcal{F}_0}} \right) \right)$ in the following way:
\begin{align} \label{at21}
W_{ij}:\Gamma\left( U_{ij}, T_{\mathcal{F}_0} \right) &\to \Gamma\left( U_i, \frac{T_X}{T_{\mathcal{F}_0}} \right) \\
 T_{0i}^\alpha &\mapsto  \overline{W_{ij}^\alpha} \notag
\end{align}
and linearly extends to $\Gamma\left( U_{ij}, T_{\mathcal{F}_0}\right)$. Here $\overline{W_{ij}^\alpha}$ is the image of $W_{ij}^\alpha$ in $\Gamma\left( U_{ij}, \frac{T_X}{T_{\mathcal{F}_0}}\right)$, and we define $P_i\in \Gamma\left( U_i, \mathscr{H}om_{\mathcal{O}_X}\left( \bigwedge^2 T_{\mathcal{F}_0}, \frac{T_X}{T_{\mathcal{F}_0}} \right) \right)$ in the following way:
\begin{align}\label{at22}
P_i:\Gamma\left( U_i, \bigwedge^2 T_{\mathcal{F}_0} \right) &\to \Gamma\left( U_i, \frac{T_X}{T_{\mathcal{F}_0}} \right) \\
T_{0i}^\alpha \wedge T_{0i}^\beta &\mapsto \overline{P_i^{\alpha\beta}} \notag
\end{align}
and linearly extends to $\Gamma\left(U_i, \bigwedge^2 T_{\mathcal{F}_0} \right)$. Here $\overline{P_i^{\alpha\beta}}$ is the image of $P_i^{\alpha\beta}\in \Gamma\left( U_i, T_X \right)$ in $\Gamma\left( U_i, \frac{T_X}{T_{\mathcal{F}_0}} \right)$. Then we claim that $\left( \left\{ P_i\right\}, \left\{ - W_{ij} \right\} \right)\in C^0\left(\mathcal{U}, \mathscr{H}om_{\mathcal{O}_X}\left( \bigwedge^2 T_{\mathcal{F}_0}, \frac{ T_X}{T_{\mathcal{F}_0}} \right) \right) \bigoplus C^1\left(\mathcal{U}, \mathscr{H}om_{\mathcal{O}_X}\left( T_{\mathcal{F}_0}, \frac{T_X}{T_{\mathcal{F}_0}} \right) \right)$ defines a $1$-cocycle in the following \v Cech resolution of $\mathscr{H}om_{\mathcal{O}_X}\left( T_{\mathcal{F}_0}, \frac{T_X}{T_{\mathcal{F}_0}} \right)^\bullet$:
\begin{center}
$\begin{CD}
\cdots \\
@AD_2AA\\
C^0\left( \mathcal{U}, \mathscr{H}om_{\mathcal{O}_X}\left( \bigwedge^2 T_{\mathcal{F}_0}, \frac{T_X}{T_{\mathcal{F}_0}} \right) \right) @>-\delta>> \cdots \\
@AD_1'AA @AD_1' AA \\
C^0\left(\mathcal{U}, \mathscr{H}om_{\mathcal{O}_X}\left(T_{\mathcal{F}_0}, \frac{T_X}{T_{\mathcal{F}_0}} \right) \right) @>\delta>> C^1\left( \mathcal{U}, \mathscr{H}om_{\mathcal{O}_X}\left(  T_{\mathcal{F}_0}, \frac{T_X}{T_{\mathcal{F}_0}} \right) \right) @>-\delta>> 
\end{CD}$
\end{center}
We note that
\begin{align}\label{at9}
tW_{ij}^\alpha-tW_{ik}^\alpha+t\sum_{\beta=1}^p r_{0ij}^{\alpha\beta} W_{jk}^\beta&=\tilde{T}_i^\alpha-\sum_{\beta=1}^p \tilde{r}_{ij}^{\alpha\beta} \tilde{T}_j^\beta -\tilde{T}_i^\alpha+\sum_{\beta=1}^p \tilde{r}_{ik}^{\alpha\beta} \tilde{T}_k^\beta+\sum_{\beta=1}^p \tilde{r}_{ij}^{\alpha\beta} \tilde{T}_j^\beta-\sum_{\beta,\gamma=1}^p \tilde{r}_{ij}^{\alpha\beta} \tilde{r}_{jk}^{\beta\gamma} \tilde{T}_k^\gamma\\
&=\sum_{\gamma=1}^p \left(\tilde{r}_{ik}^{\alpha\gamma}-\sum_{\beta=1}^p \tilde{r}_{ij}^{\alpha\beta}\tilde{r}_{jk}^{\beta\gamma} \right) \tilde{T}_k^\gamma=t\sum_{\gamma=1}^p t C_{ijk}^{\alpha\gamma} T_{0k}^\gamma \notag
\end{align}
This implies that $-\delta\left(\left\{- W_{ij} \right\} \right)=0$. As in the proof of Proposition \ref{at3}, $-\delta\left( \left\{ P_i \right\} \right) + D_1'\left( \left\{ -W_{ij} \right\} \right)=0$ and $D_2\left(\left\{ P_i \right\} \right)=0$. We show that for another lifting $\tilde{T}_i'^\alpha$ of $T_i^\alpha$, and another lifting $\tilde{r}_{ij}'^{\alpha\beta}$ of $r_{ij}^{\alpha\beta}$ and another lifting of $\tilde{g}_{i\alpha\beta}'^\gamma$ of $g_{i\alpha\beta}^\gamma,\alpha,\beta=1,...,p$, the associated $1$-cocycle $\beta:=\left(\left\{ P_i' \right\}, \left\{ - W_{ij}' \right\} \right)$ is cohomologous to $\alpha:=\left( \left\{ P_i \right\}, \left\{- W_{ij} \right\} \right)$. Then $\tilde{T}_i'^\alpha= \tilde{T}_i^\alpha + t Y_i^\alpha$ for some $Y_i^\alpha$ and $\tilde{r}_{ij}'^{\alpha\beta}= \tilde{r}_{ij}^{\alpha\beta} + t s_{ij}^{\alpha\beta}$ for some $s_{ij}^{\alpha\beta}\in \Gamma\left( U_{ij}, \mathcal{O}_X \right)$ and $\tilde{g}_{i\alpha\beta}'^\gamma = \tilde{g}_{i\alpha\beta}^\gamma + t e_{i\alpha\beta}^\gamma$ for some $e_{i\alpha\beta}^\gamma \in \Gamma\left( U_i, \mathcal{O}_X \right)$. We define $Y_i \in \Gamma\left( U_i, \mathscr{H}om_{\mathcal{O}_X}\left( T_{\mathcal{F}_0}, \frac{T_X}{T_{\mathcal{F}_0}} \right) \right)$ by $Y_i:\Gamma\left( U_i, T_{\mathcal{F}_0} \right)\to \Gamma\left( U_i, \frac{T_X}{T_{\mathcal{F}_0}} \right), T_{0i}^\alpha \mapsto \overline{Y_i^\alpha}$. Then as in $(\ref{at11})$, we have $D_1\left( \left\{ Y_i \right\} \right)= P_i'- P_i$. On the other hand, we see that
\begin{align}\label{at13}
t\left(- W_{ij}'^\alpha + W_{ij}^\alpha \right) &= -\tilde{T}_i^\alpha - t Y_i^\alpha + \sum_{\beta=1}^p \left(\tilde{r}_{ij}^{\alpha\beta}+ t s_{ij}^{\alpha\beta} \right) \left(\tilde{T}_j^\beta + t Y_j^\beta \right) + \tilde{T}_i^\alpha - \sum_{\beta=1}^p \tilde{r}_{ij}^{\alpha\beta} \tilde{T}_j^\beta  \\
&=  t\left( \sum_{\beta=1}^p r_{0ij}^{\alpha\beta} Y_j^\beta- Y_i^\alpha \right) + t\sum_{\beta=1}^p  s_{ij}^{\alpha\beta} T_{0j}^\beta \notag
\end{align}
This implies that $\delta\left(\left\{ Y_i\right\} \right)=\left\{ -W_{ij}'-\left(- W_{ij}\right) \right\}$. This proves that $\alpha$ is cohomologous to $\beta$. So given a small extension $e:0\to (t)\to \tilde{A} \to A \to 0$ and an infinitesimal foliated deformation $\eta$ of $\left( X, T_{\mathcal{F}_0}\right)$ with $X$ fixed, we can associate an element $o_\eta(e):=$ the cohomology class of $\alpha \in \mathbb{H}^1\left( X, \mathscr{H}om_{\mathcal{O}_X}\left( T_{\mathcal{F}_0}, \frac{T_X}{T_{\mathcal{F}_0}} \right)^\bullet \right)$. We note that $o_\eta(e)=0$ if and only if there exists $\left\{ X_i \right\} \in C^0\left( \mathcal{U}, \mathscr{H}om_{\mathcal{O}_X}\left( T_{\mathcal{F}_0}, \frac{T_X}{T_{\mathcal{F}_0}} \right)\right)$ where $X_i:\Gamma\left(U_i, T_{\mathcal{F}_0} \right)\to \Gamma\left( U_i, T_X \right), T_{0i}^\alpha \mapsto \overline{X_i^\alpha}$ with $X_i^\alpha \in \Gamma\left( U_i, T_X \right)$ such that $X_j- X_i= -W_{ij}$ and $D_1\left(X_i \right)=P_i$, so that
{\small{\begin{align}
&t\left(\sum_{\beta=1}^p r_{0ij}^{\alpha\beta} X_j^\beta - X_i^\alpha +  \sum_{\beta=1}^p v_{ij}^{\alpha\beta} T_{0j}^\beta\right)= - \tilde{T}_i^\alpha + \sum_{\beta=1}^p \tilde{r}_{ij}^{\alpha\beta} \tilde{T}_j^\beta \notag \\
& \Longrightarrow \tilde{T}_i^\alpha - tX_i^\alpha  = \sum_{\beta=1}^p \left(\tilde{r}_{ij}^{\alpha\beta}- tv_{ij}^{\alpha\beta} \right)\left(\tilde{T}_j^\beta - t X_j^\beta \right)\,\,\,\textnormal{for some} \,\,\, v_{ij}^{\alpha\beta}\in \Gamma\left( U_{ij}, \mathcal{O}_X \right) \label{at14}
\end{align}}}
and as in $(\ref{at12})$, we have  $\left[ \tilde{T}_i^\alpha- t X_i^\alpha, \tilde{T}_i^\beta - t X_i^\beta \right] - \sum_{\gamma=1}^p\left( \tilde{g}_{i\alpha\beta}^\gamma + t l_{i\alpha\beta}^\gamma  \right)\left(\tilde{T}_i^\gamma - t X_i^\gamma \right)$ for some $l_{i\alpha\beta}^\gamma \in \Gamma\left(U_i, \mathcal{O}_X \right)$. Hence a lifting $\tilde{\eta}$ of $\eta$ exits. This completes the proof of Proposition \ref{at8}.
\end{proof}

\subsection{Foliated deformations of $T_{\mathcal{F}_0}$ with $X$ fixed in terms of vector bundles} \

Let $\left( X, T_{\mathcal{F}_0} \right)$ be a foliated nonsingular variety with $T_{\mathcal{F}_0}$ locally free. The underlying deformations of locally free sheaf $T_{\mathcal{F}_0}$ is implicit in the complex $\mathscr{H}om_{\mathcal{O}_X}\left(T_{\mathcal{F}_0}, \frac{T_X}{T_{\mathcal{F}_0}} \right)^\bullet$. In order to make the underlying deformations of $T_{\mathcal{F}_0}$ explicit, we define another complex of sheaves also controlling flat foliated deformations of $\left(X, T_{\mathcal{F}_0} \right)$ under the trivial deformation of $X$ and present foliated deformations of $\left(X, T_{\mathcal{F}_0}\right)$ with $X$ fixed by using the complex. As is known, deformations of a locally free sheaf $T_{\mathcal{F}_0}$ on $X$ is controlled by $\mathscr{H}om_{\mathcal{O}_X}\left( T_{\mathcal{F}_0}, T_{\mathcal{F}_0} \right)$ (see \cite{Har10}). We note that we have an injection $j:\mathscr{H}om_{\mathcal{O}_X}\left( T_{\mathcal{F}_0}, T_{\mathcal{F}_0} \right) \to \mathscr{H}om_{\mathcal{O}_X}\left( T_{\mathcal{F}_0}, T_X \right)$ induced by $T_{\mathcal{F}_0} \hookrightarrow T_X$. Then let us consider an exact sequence of complexes of sheaves
{\small{\begin{center}
$\begin{CD}
@. \cdots @. \cdots @. \cdots \\
@. @AAA @AD_3AA @AD_3AA \\
0@>>>0 @>>> \mathscr{H}om_{\mathcal{O}_X}\left(\bigwedge^3 T_{\mathcal{F}_0}, \frac{T_X}{T_{\mathcal{F}_0}} \right) @>>> \mathscr{H}om_{\mathcal{O}_X}\left( \bigwedge^3 T_{\mathcal{F}_0}, \frac{T_X}{T_{\mathcal{F}_0}}\right)@>>> 0 \\
@.@AAA @AD_2AA @AD_2AA \\
0@>>>0 @>>> \mathscr{H}om_{\mathcal{O}_X}\left(\bigwedge^2 T_{\mathcal{F}_0}, \frac{T_X}{T_{\mathcal{F}_0}} \right) @>>> \mathscr{H}om_{\mathcal{O}_X}\left( \bigwedge^2 T_{\mathcal{F}_0}, \frac{T_X}{T_{\mathcal{F}_0}}\right) @>>> 0\\
@. @AAA @AD_1'AA @AD_1AA \\
0@>>> \mathscr{H}om_{\mathcal{O}_X}\left(T_{\mathcal{F}_0}, T_{\mathcal{F}_0} \right) @>>> \mathscr{H}om_{\mathcal{O}_X}\left(T_{\mathcal{F}_0},T_X \right) @>>> \mathscr{H}om_{\mathcal{O}_X}\left(T_{\mathcal{F}_0},\frac{T_X}{T_{\mathcal{F}_0}} \right) @>>> 0 \\
@. @A\textnormal{Id}AA @AjAA @AAA \\
0@>>> \mathscr{H}om_{\mathcal{O}_X}\left(T_{\mathcal{F}_0}, T_{\mathcal{F}_0} \right) @>>> \mathscr{H}om_{\mathcal{O}_X}\left( T_{\mathcal{F}_0}, T_{\mathcal{F}_0} \right) @>>> 0@>>> 0
\end{CD}$
\end{center}}}
Let us denote the middle complex of sheaves by $\mathscr{H}om_{\mathcal{O}_X}\left(T_{\mathcal{F}_0}, T_{\mathcal{F}_0} \right)^\bullet$:
{\Small{\begin{align}\label{tf1}
\mathscr{H}om_{\mathcal{O}_X}\left(T_{\mathcal{F}_0}, T_{\mathcal{F}_0} \right)^\bullet: \mathscr{H}om_{\mathcal{O}_X}\left(T_{\mathcal{F}_0}, T_{\mathcal{F}_0}\right) \to  \mathscr{H}om_{\mathcal{O}_X}\left(T_{\mathcal{F}_0}, T_X \right) \to \mathscr{H}om_{\mathcal{O}_X}\left(\bigwedge^2 T_{\mathcal{F}_0}, \frac{T_X}{T_{\mathcal{F}_0}} \right)\to \mathscr{H}om_{\mathcal{O}_X}\left( \bigwedge^3 T_{\mathcal{F}_0}, \frac{T_X}{T_{\mathcal{F}_0}}\right)\to \cdots
\end{align}}}
We will denote the $i$-th cohomology group by $\mathbb{H}^i\left( X, \mathscr{H}om_{\mathcal{O}_X}\left(T_{\mathcal{F}_0}, T_{\mathcal{F}_0} \right)^\bullet \right)$. From the exact sequence we see that
\begin{align*}
\mathbb{H}^i\left(X, \mathscr{H}om_{\mathcal{O}_X}\left(T_{\mathcal{F}_0}, T_{\mathcal{F}_0} \right)^\bullet \right)\cong \mathbb{H}^{i-1}\left( X, \mathscr{H}om_{\mathcal{O}_X}\left( T_{\mathcal{F}_0}, T_X/ T_{\mathcal{F}_0}\right)^\bullet\right),\,\,\,\,\,i\geq 0 . 
\end{align*}
This shows that $\mathscr{H}om_{\mathcal{O}_X}\left( T_{\mathcal{F}_0}, T_{\mathcal{F}_0}\right)^\bullet $ also controls flat foliated deformations of $\left( X, T_{\mathcal{F}_0} \right)$ with $X$ fixed. Moreover the complex $(\ref{tf1})$ explicitly contains $\mathscr{H}om_{\mathcal{O}_X}\left( T_{\mathcal{F}_0}, T_{\mathcal{F}_0} \right)$ as its component. We will reinterpret flat foliated deformations of $\left(X, T_{\mathcal{F}_0} \right)$ with $X$ fixed by using the complex $(\ref{tf1})$.

\begin{remark}
To emphasize deformations of the underlying locally free sheaf $T_{\mathcal{F}_0}$, we can define infinitesimal deformations of $\left(X, T_{\mathcal{F}_0}\right)$ with $X$ fixed equivalently in the following way: an infinitesimal foliated deformation $\left(\mathcal{X}=X\times \textnormal{Spec}(A), T_{\mathcal{F}}\right)$ of $\left(X, T_{\mathcal{F}_0} \right)$ under the trivial deformation of $X$ over $A\in \bold{Art}$ is a locally free sheaf $T_{\mathcal{F}}$ and an exact sequence on $\mathcal{X}=X\times \textnormal{Spec}(A)$
\begin{align*}
0 \to T_{\mathcal{F}}\xrightarrow{l} T_{\mathcal{X}/A}=T_X\otimes A \to \frac{T_X\otimes A}{l\left(T_{\mathcal{F}} \right)}\to 0
\end{align*}
such that $\left[ l\left( T_{\mathcal{F}} \right), l\left(T_{\mathcal{F}} \right)\right]\subset l\left(T_{\mathcal{F}} \right)$ and $\frac{T_{\mathcal{X}/A} }{l\left(T_{\mathcal{F}}\right)}$ is flat over $A$, and induces $0\to T_{\mathcal{F}_0} \to T_X \to \frac{T_X}{T_{\mathcal{F}_0}}\to 0$ over $\textnormal{Spec}(k)$. Two infinitesimal deformations $\left(\mathcal{X}, T_{\mathcal{F}}\right)$ and $\left(\mathcal{X}, T_{\mathcal{F}'}\right)$ of $\left(X, T_{\mathcal{F}_0} \right)$ with $X$ fixed are equivalent  if there is an isomorphism $\alpha: T_{\mathcal{F}}\to T_{\mathcal{F}'}$ over $\mathcal{X}$ such that $l= l'\circ \alpha$ where $0 \to T_{\mathcal{F}'}\xrightarrow{l'} T_{\mathcal{X}/A} \to \frac{T_{\mathcal{X}/A}}{T_{\mathcal{F}'}}\to 0$. Then we have a functor of Artin rings
\begin{align*}
\textnormal{FDef}_{T_{\mathcal{F}_0}}: \bold{Art} &\to \textnormal{(Sets)} \\
                             A &\to \left\{ \textnormal{Infinitesimal foliated deformation of $\left( X, T_{\mathcal{F}_0}\right)$ with $X$ fixed}   \right\}/\textnormal{isomorphisms}
\end{align*}
\end{remark}

\begin{proposition}\label{at15}
Let $\left(X, T_{\mathcal{F}_0} \right)$ be a foliated nonsingular variety with $T_{\mathcal{F}_0}$ locally free.
\begin{enumerate}
\item there is a $1-1$ correspondence
\begin{align*}
\textnormal{FDef}_{T_{\mathcal{F}_0}}\left(k[\epsilon] \right)\cong \mathbb{H}^1\left(X, \mathscr{H}om_{\mathcal{O}_X}\left(  T_{\mathcal{F}_0}, T_{\mathcal{F}_0} \right)^\bullet \right)
\end{align*}
\item Let $A\in \bold{Art}$ and $\eta=\left( \mathcal{X}, T_\mathcal{F} \right)$ be an infinitesimal foliated deformations over $A$ of $\left( X, T_{\mathcal{F}_0} \right)$ with $X$ fixed. To every small extension $0\to (t)\to \tilde{A} \to A \to 0$, we can associate an element $o_\eta(e)\in \mathbb{H}^2\left( X, \mathscr{H}om_{\mathcal{O}_X}\left( T_{\mathcal{F}_0}, T_{\mathcal{F}_0} \right)^\bullet \right)$ called the obstruction lifting of $\eta$ to $\tilde{A}$, which is $0$ if and only if a lifting of $\eta$ to $\tilde{A}$ exists.
\end{enumerate}
\end{proposition}

\begin{proof}
We keep the notations in the proof of Proposition \ref{at8}. Given a first-order foliated deformation $\left(\mathcal{X}, T_{\mathcal{F}} \right)$ of $\left( X, T_{\mathcal{F}_0}\right)$ with $X$ fixed, let us define $T_i' \in \Gamma\left( U_i, \mathscr{H}om_{\mathcal{O}_X}\left(T_{\mathcal{F}_0}, T_X \right)\right)$ by $T_i':\Gamma\left(U_i, T_{\mathcal{F}_0} \right)\to \Gamma\left(U_i, T_X \right), T_{0i}^\alpha \mapsto T_i'^\alpha$ and define $B_{ij}\in \Gamma\left(U_{ij}, \mathscr{H}om_{\mathcal{O}_X}\left( T_{\mathcal{F}_0}, T_{\mathcal{F}_0} \right) \right)$ by
\begin{align*}
B_{ij}: \Gamma\left( U_{ij}, T_{\mathcal{F}_0}\right) &\to \Gamma\left( U_{ij}, T_{\mathcal{F}_0} \right) \\
 T_{0i}^\alpha &\mapsto \sum_{\beta=1}^p b_{ij}^{\alpha\beta} T_{0j}^\beta
\end{align*}
Then we claim that $\left( \left\{ T_i' \right\} , \left\{ B_{ij}  \right\} \right) \in C^0\left( \mathcal{U}, \mathscr{H}om_{\mathcal{O}_X}\left( T_{\mathcal{F}_0}, T_X \right) \right) \bigoplus C^1\left(\mathcal{U}, \mathscr{H}om_{\mathcal{O}_X}\left(  T_{\mathcal{F}_0}, T_{\mathcal{F}_0} \right) \right)$ is a $1$-cocycle in the following \v Cech resolution of $\mathscr{H}om_{\mathcal{O}_X}\left( T_{\mathcal{F}_0}, T_{\mathcal{F}_0} \right)^\bullet$:{\small{\begin{center}
$\begin{CD}
\cdots \\
@AD_2AA \\
C^0\left(\mathcal{U}, \mathscr{H}om_{\mathcal{O}_X}\left( \bigwedge^2 T_{\mathcal{F}_0}, \frac{T_X}{T_{\mathcal{F}_0}} \right) \right) @>-\delta>> \cdots \\
@AD_1'AA @AD_1'AA \\
C^0\left( \mathcal{U}, \mathscr{H}om_{\mathcal{O}_X}\left(T_{\mathcal{F}_0}, T_X \right) \right) @>\delta>> C^1\left( \mathcal{U}, \mathscr{H}om_{\mathcal{O}_X}\left( T_{\mathcal{F}_0}, T_X \right) \right) @>-\delta>> \cdots \\
@AjAA @AjAA @AjAA \\
C^0\left( \mathcal{U}, \mathscr{H}om_{\mathcal{O}_X}\left( T_{\mathcal{F}_0}, T_{\mathcal{F}_0}\right) \right) @>-\delta>> C^1\left( \mathcal{U}, \mathscr{H}om_{\mathcal{O}_X}\left(T_{\mathcal{F}_0}, T_{\mathcal{F}_0} \right) \right) @>\delta>> C^2\left( \mathcal{U}, \mathscr{H}om_{\mathcal{O}_X}\left( T_{\mathcal{F}_0}, T_{\mathcal{F}_0} \right) \right) @>>> \cdots
\end{CD}$
\end{center}}}
As in the proof of Proposition  \ref{at8}, we have $D_1'\left(\left\{  T_i' \right\} \right)=0$ and $(\ref{at6})$ implies that $\delta\left( \left\{ T_i' \right\} \right)+ j\left(\left\{ B_{ij} \right\} \right)=0$ and we note that
\begin{align*}
&r_{0ik}^{\alpha \gamma} + \epsilon b_{ik}^{\alpha \gamma} = \sum_{\beta=1}^p \left( r_{0ij}^{\alpha \beta} + \epsilon b_{ij}^{\alpha\beta} \right) \left( r_{0jk}^{\beta\gamma} + \epsilon b_{jk}^{\beta \gamma} \right) \Longrightarrow b_{ik}^{\alpha \gamma}- \sum_{\beta=1}^p r_{0ij}^{\alpha \beta} b_{jk}^{\beta \gamma} - \sum_{\beta=1}^p b_{ij}^{\alpha\beta} r_{0jk}^{\beta \gamma} =0 \\
&\Longrightarrow \left(B_{ij} - B_{ik} + B_{jk} \right)\left(T_{0i}^\alpha \right) = \sum_{\beta=1}^p b_{ij}^{\alpha\beta} T_{0j}^\beta - \sum_{\beta=1}^p b_{ik}^{\alpha\beta} T_{0k}^\beta + \sum_{\beta,\gamma=1}^p r_{0ij}^{\alpha\beta} b_{jk}^{\beta \gamma} T_{0k}^\gamma = 0
\end{align*}
This implies that $\delta\left(\left\{ B_{ij} \right\} \right)=0$. This proves the claim. For another equivalent foliated deformations $\beta:=\left( \left\{ T_i''\right\}, \left\{ B_{ij}'\right\} \right)$ of $\left(X, T_{\mathcal{F}_0} \right)$ with $X$ fixed, we show that $\beta$ is cohomologous to $\alpha:=\left( \left\{ T_i'\right\}, \left\{ B_{ij}\right\} \right)$. We note that we have $T_{0i}^\alpha+ \epsilon T_i''^\alpha= \sum_{\beta=1}^p \left( \delta_\alpha^\beta + \epsilon a_{i}^{\alpha\beta} \right)\left(T_{0i}^\beta + \epsilon T_i'^\beta \right)$ for some $a_i^{\alpha\beta} \in \Gamma\left( U_i, \mathcal{O}_X \right)$, so that we have $T_i''^\alpha= \sum_{\beta=1}^p a_i^{\alpha\beta} T_{0i}^\beta + T_i'^\alpha$.  Then we define $A_i\in \Gamma\left( U_i, \mathscr{H}om_{\mathcal{O}_X}\left( T_{\mathcal{F}_0}, T_{\mathcal{F}_0} \right) \right)$ by 
\begin{align*}
A_i:\Gamma\left( U_i, T_{\mathcal{F}_0}\right) \to \Gamma\left( U_i, T_{\mathcal{F}_0}\right)\\
 T_{0i}^\alpha \mapsto \sum_{\beta=1}^p a_i^{\alpha\beta} T_{0i}^\beta
\end{align*}
Then $j(A_i)(T_{0i}^\alpha)= T_i''^\alpha- T_i'^\alpha$ and since we have from $(\ref{at6})$
{\small{\begin{align*}
\left(B_{ij}'- B_{ij}\right)\left(T_{0i}^\alpha \right)= \sum_{\beta=1}^p b_{ij}'^{\alpha\beta} T_{0j}^\beta- \sum_{\beta=1}^p b_{ij}^{\alpha\beta}T_{0j}^\beta = T_i''^\alpha - T_i'^\alpha - \sum_{\beta=1}^p r_{0ij}^{\alpha\beta}\left( T_j''^\beta- T_j'^\beta\right)= \sum_{\beta=1}^p a_i^{\alpha\beta} T_{0i}^\beta - \sum_{\beta,\gamma=1}^p r_{0ij}^{\alpha\beta} a_j^{\beta \gamma} T_{0j}^\gamma
\end{align*}}}
we see that $-\delta\left(\left\{ A_i \right\} \right)= B_{ij}'-B_{ij}$. Hence $\alpha$ is cohomologus to $\beta$. As in Proposition $\ref{at8}$ we define an inverse map.

Next we identify obstructions. We define $W_{ij}\in \Gamma\left( U_{ij}, \mathscr{H}om_{\mathcal{O}_X}\left( T_{\mathcal{F}_0}, T_X \right) \right)$ by $W_{ij}:\Gamma\left(U_{ij}, T_{\mathcal{F}_0} \right) \to \Gamma\left( U_{ij}, T_X \right), T_{0i}^\alpha \mapsto W_{ij}^\alpha$ and $P_i\in \Gamma\left( U_i, \mathscr{H}om_{\mathcal{O}_X}\left( \bigwedge^2 T_{\mathcal{F}_0} , \frac{T_X}{T_{\mathcal{F}_0}} \right) \right)$ as in the proof of Proposition \ref{at8}. We define $C_{ijk}\in \Gamma\left( U_{ijk}, \mathscr{H}om_{\mathcal{O}_X}\left( T_{\mathcal{F}_0}, T_{\mathcal{F}_0} \right) \right)$ by
\begin{align*}
C_{ijk}:\Gamma\left(U_{ijk}, T_{\mathcal{F}_0} \right)&\to \Gamma\left( U_{ijk}, T_{\mathcal{F}_0} \right) \\
 T_{0i}^\alpha &\mapsto \sum_{\gamma=1}^p C_{ijk}^{\alpha \gamma} T_{0k}^\gamma
\end{align*}
Then $\left\{ C_{ijk}   \right\}\in C^2\left(\mathcal{U}, \mathscr{H}om_{\mathcal{O}_X}\left( T_{\mathcal{F}_0}, T_{\mathcal{F}_0} \right) \right)$. We note that
\begin{align}\label{ot46}
&tC_{ijl}^{\alpha\gamma}= \tilde{r}_{il}^{\alpha \gamma}- \sum_{\beta=1}^p\tilde{r}_{ij}^{\alpha \beta} \tilde{r}_{jl}^{\beta\gamma} = \sum_{\beta=1}^p \tilde{r}_{ik}^{\alpha\beta} \tilde{r}_{kl}^{\beta\gamma}+ t C_{ikl}^{\alpha\gamma} -\sum_{\beta=1}^p\tilde{r}_{ij}^{\alpha\beta}\left(\sum_{\delta=1}^p \tilde{r}_{jk}^{\beta\delta}\tilde{r}_{kl}^{\delta\gamma} + tC_{jkl}^{\beta\gamma} \right)\\
&\iff \sum_{\beta=1}^p \tilde{r}_{ik}^{\alpha\beta} \tilde{r}_{kl}^{\beta\gamma}+t C_{ikl}^{\alpha\gamma}-\sum_{\delta=1}^p\left( \tilde{r}_{ik}^{\alpha\delta}-tC_{ijk}^{\alpha\delta}\right)\tilde{r}_{kl}^{\delta\gamma}-t\sum_{\beta=1}^p r_{0ij}^{\alpha\beta} C_{jkl}^{\beta\gamma}-tC_{ijl}^{\alpha\gamma}=0 \notag \\
&\iff C_{ikl}^{\alpha\gamma}+\sum_{\delta=1}^p C_{ijk}^{\alpha\delta} r_{0kl}^{\delta\gamma}-\sum_{\beta=1}^p r_{0ij}^{\alpha\beta} C_{jkl}^{\beta\gamma}- C_{ijl}^{\alpha\gamma}=0 \notag\\
&\Longrightarrow \left(C_{jkl} - C_{ikl} + C_{ijl} - C_{ijk} \right)\left(T_{0i}^\alpha \right) = \sum_{\beta,\gamma=1}^p r_{0ij}^{\alpha\beta} C_{jkl}^{\beta \gamma}T_{0l}^\gamma - \sum_{\gamma=1}^p C_{ikl}^{\alpha \gamma} T_{0l}^\gamma + \sum_{\gamma=1}^p C_{ijl}^{\alpha \gamma} T_{0l}^\gamma - \sum_{\beta=1}^p C_{ijk}^{\alpha \beta} T_{0k}^\beta=0\notag
\end{align}
This implies that $-\delta\left(\left\{ -C_{ijk} \right\} \right)=0$. From $(\ref{at9})$, we have $-\delta\left( \left\{ -W_{ij} \right\} \right) + j\left(  \left\{ - C_{ijk} \right\} \right)=0$. As in the proof of Proposition $(\ref{at3})$, $-\delta\left( \left\{P_i \right\} \right)+ D_1'\left( \left\{- W_{ij}\right\} \right)=0$ and $D_2\left( \left\{ P_i \right\} \right)=0$. Hence we see that
{\small{\begin{align*}
\left( \left\{ P_i \right\}, \left\{ -W_{ij} \right\}, \left\{ - C_{ijk} \right\} \right) \in C^0\left( \mathcal{U}, \mathscr{H}om_{\mathcal{O}_X}\left( \bigwedge^2 T_{\mathcal{F}_0}, \frac{T_X}{T_{\mathcal{F}_0}} \right) \right) \bigoplus C^1\left( \mathcal{U}, \mathscr{H}om_{\mathcal{O}_X}\left( T_{\mathcal{F}_0}, T_X \right) \right) \bigoplus C^2\left( \mathcal{U}, \mathscr{H}om_{\mathcal{O}_X}\left( T_{\mathcal{F}_0} , T_{\mathcal{F}_0} \right) \right)
\end{align*}}}
defines a $2$-cocycle in the above \v Cech resolution of $\mathscr{H}om_{\mathcal{O}_X}\left( T_{\mathcal{F}_0}, T_{\mathcal{F}_0} \right)^\bullet$. We show that for another lifting $\tilde{T}_i'^\alpha$ of $T_i^\alpha$, and another lifting $\tilde{r}_{ij}'^{\alpha\beta}$ of $r_{ij}^{\alpha\beta}$ and another lifting $\tilde{g}_{i\alpha\beta}'^\gamma$ of $g_{i\alpha\beta}^\gamma,\alpha,\beta=1,...,p$, the associated $2$-cocycle $\beta:=\left(\left\{ P_i' \right\}, \left\{ - W_{ij}' \right\} , \left\{ - C_{ijk}' \right\} \right)$ is cohomologous to $\alpha:=\left( \left\{ P_i \right\}, \left\{ -W_{ij} \right\} , \left\{ -C_{ijk} \right\}\right)$. Then $\tilde{T}_i'^\alpha= \tilde{T}_i^\alpha + t Y_i^\alpha$ for some $Y_i^\alpha$ and $\tilde{r}_{ij}'^{\alpha\beta}= \tilde{r}_{ij}^{\alpha\beta} + t s_{ij}^{\alpha\beta}$ for some $s_{ij}^{\alpha\beta}\in \Gamma\left( U_{ij}, \mathcal{O}_X \right)$ and $\tilde{g}_{i\alpha\beta}'^\gamma = \tilde{g}_{i\alpha\beta}^\gamma + t e_{i\alpha\beta}^\gamma$ for some $e_{i\alpha\beta}^\gamma \in \Gamma\left( U_i, \mathcal{O}_X \right)$. We define $Y_i \in \Gamma\left( U_i, \mathscr{H}om_{\mathcal{O}_X}\left( T_{\mathcal{F}_0}, T_X \right) \right)$ by $Y_i:\Gamma\left( U_i, T_{\mathcal{F}_0} \right)\to \Gamma\left( U_i, T_X \right), T_{0i}^\alpha \mapsto Y_i^\alpha$ and define $C_{ij}\in \Gamma\left( U_{ij}, \mathscr{H}om_{\mathcal{O}_X}\left( T_{\mathcal{F}_0}, T_{\mathcal{F}_0} \right) \right)$ by
\begin{align*}\label{ot47}
C_{ij}:\Gamma\left( U_{ij}, T_{\mathcal{F}_0}  \right) &\to \Gamma\left( U_{ij}, T_{\mathcal{F}_0} \right)\\
  T_{0i}^\alpha &\mapsto  \sum_{\beta=1}^p s_{ij}^{\alpha\beta} T_{0j}^\beta           
\end{align*}
Then as in $(\ref{at11})$, we have $D_1'\left( \left\{ Y_i \right\} \right)= P_i'- P_i$ and $(\ref{at13})$ implies that $\delta\left(\left\{ Y_i \right\} \right) + j\left(\left\{ C_{ij} \right\} \right)= - W_{ij}'+ W_{ij}$. We note that
\begin{align}
t\left(- C_{ijk}' + C_{ijk} \right)\left(T_{0i}^\alpha \right)& = - \sum_{\gamma=1}^p \left(\tilde{r}_{ik}^{\alpha \gamma}+ t s_{ik}^{\alpha\gamma} - \sum_{\beta=1}^p \left( \tilde{r}_{ij}^{\alpha\beta}+ t s_{ij}^{\alpha\beta} \right) \left(\tilde{r}_{jk}^{\beta \gamma} + t s_{jk}^{\beta \gamma} \right) - \tilde{r}_{ik}^{\alpha \gamma} + \sum_{\beta=1}^p  \tilde{r}_{ij}^{\alpha \beta} \tilde{r}_{jk}^{\beta \gamma}    \right) \tilde{T}_{k}^\gamma \\
&= - t \sum_{\gamma=1}^p \left( s_{ik}^{\alpha \gamma}- \sum_{\beta=1}^p r_{0ij}^{\alpha\beta} s_{jk}^{\beta \gamma} - \sum_{\beta=1}^p s_{ij}^{\alpha\beta} r_{0jk}^{\beta \gamma}  \right) T_{0k}^\gamma = t \left( - C_{ik}\left(T_{0i}^\alpha \right)  + C_{jk}\left( T_{0i}^\alpha \right)  + C_{ij}\left(T_{0i}^\alpha \right)  \right) \notag
\end{align}
This implies that $\delta\left( \left\{ C_{ij} \right\} \right) = \left\{ - C_{ijk}'  +  C_{ijk} \right\}$. This proves that $\alpha$ is cohomologous to $\beta$. So given a small extension $e:0\to (t)\to \tilde{A} \to A \to 0$ and an infinitesimal foliated deformation $\eta$ of $\left( X, T_{\mathcal{F}_0}\right)$ with $X$ fixed, we can associate an element $o_\eta(e):=$ the cohomology class of $\alpha \in \mathbb{H}^2\left( X, \mathscr{H}om_{\mathcal{O}_X}\left( T_{\mathcal{F}_0}, T_{\mathcal{F}_0} \right)^\bullet \right)$. We note that $o_\eta(e)=0$ if and only if there exists $\left\{ X_i \right\} \in C^0\left( \mathcal{U}, \mathscr{H}om_{\mathcal{O}_X}\left( T_{\mathcal{F}_0}, T_X \right)\right)$ where $X_i:\Gamma\left(U_i, T_{\mathcal{F}_0} \right)\to \Gamma\left( U_i, T_X \right), T_{0i}^\alpha \mapsto X_i^\alpha$ and $V_{ij} \in \Gamma\left( U_{ij}, \mathscr{H}om_{\mathcal{O}_X}\left( T_{\mathcal{F}_0}, T_{\mathcal{F}_0} \right) \right)$ where $V_{ij}:\Gamma\left(U_{ij}, T_{\mathcal{F}_0} \right) \to \Gamma\left( U_{ij}, T_{\mathcal{F}_0} \right), T_{0i}^\alpha \mapsto \sum_{\beta=1}^p v_{ij}^{\alpha\beta} T_{0j}^\beta$ such that $X_j- X_i + j\left(V_{ij}\right)= - W_{ij}$ and $D_1'\left( X_i\right)=P_i$ and $V_{ij}- V_{ik} + V_{jk}= - C_{ijk}$. Then as in $(\ref{at14})$ we have $\tilde{T}_i^\alpha - tX_i^\alpha  = \sum_{\beta=1}^p \left(\tilde{r}_{ij}^{\alpha\beta}- tv_{ij}^{\alpha\beta} \right)\left(\tilde{T}_j^\beta - t X_j^\beta \right)$ and as in $(\ref{at12})$, we have  $\left[ \tilde{T}_i^\alpha- t X_i^\alpha, \tilde{T}_i^\beta - t X_i^\beta \right] - \sum_{\gamma=1}^p\left( \tilde{g}_{i\alpha\beta}^\gamma + t l_{i\alpha\beta}^\gamma  \right)\left(\tilde{T}_i^\gamma - t X_i^\gamma \right)$ for some $l_{i\alpha\beta}^\gamma \in \Gamma\left(U_i, \mathcal{O}_X \right)$ and
\begin{align}\label{ot32}
&t\left(V_{ij}- V_{ik} + V_{jk} \right)\left(T_{0i}^\alpha \right)= - \sum_{\gamma=1}^p \left( \tilde{r}_{ik}^{\alpha \gamma} - \sum_{\beta=1}^p \tilde{r}_{ij}^{\alpha \beta } \tilde{r}_{jk}^{\beta \gamma}  \right) {T}_{0k}^\gamma \\
&\Longrightarrow t\left( \sum_{\beta=1}^p v_{ij}^{\alpha\beta} T_{0j}^\beta - \sum_{\gamma=1}^p v_{ik}^{\alpha \gamma} T_{0k}^\gamma + \sum_{\beta,\gamma=1}^p r_{0ij}^{\alpha\beta} v_{jk}^{\beta \gamma} T_{0k}^\gamma    \right)=  - \sum_{\gamma=1}^p \left( \tilde{r}_{ik}^{\alpha \gamma} - \sum_{\beta=1}^p \tilde{r}_{ij}^{\alpha \beta } \tilde{r}_{jk}^{\beta \gamma}  \right) {T}_{0k}^\gamma \notag \\
&\Longrightarrow  t\left(\sum_{\beta,\gamma=1}^p v_{ij}^{\alpha\beta} r_{0jk}^{\beta\gamma} -  v_{ik}^{\alpha \gamma} + \sum_{\beta, \gamma=1}^p r_{0ij}^{\alpha\beta} v_{jk}^{\beta \gamma} \right)= - \tilde{r}_{ik}^{\alpha \gamma} + \sum_{\beta=1}^p \tilde{r}_{ij}^{\alpha \beta } \tilde{r}_{jk}^{\beta \gamma}  \Longrightarrow \tilde{r}_{ik}^{\alpha \gamma}- tv_{ik}^{\alpha \gamma} = \sum_{\beta=1}^p\left(\tilde{r}_{ij}^{\alpha \beta} - t v_{ij}^{\alpha \beta} \right)\left( \tilde{r}_{jk}^{\beta\gamma} - t v_{jk}^{\beta \gamma} \right) \notag
\end{align}
Hence a lifting $\tilde{\eta}$ of $\eta$ exits. This completes the proof of Proposition \ref{at15}.
\end{proof}

\subsection{Relationship between foliated deformations of $\left(X, T_{\mathcal{F}_0} \right)$ with $X$ fixed and the underlying deformations of $T_{\mathcal{F}_0}$} \

In the previous section we have shown that $\mathscr{H}om_{\mathcal{O}_X}\left( T_{\mathcal{F}_0}, T_{\mathcal{F}_0} \right)^\bullet$ controls flat foliated deformations of $\left(X, T_{\mathcal{F}_0} \right)$ with $X$ fixed. As is known, the underlying flat deformations of the locally free sheaf $T_{\mathcal{F}_0}$ with $X$ fixed is controlled by $\mathscr{H}om_{\mathcal{O}_X}\left(T_{\mathcal{F}_0}, T_{\mathcal{F}_0} \right)$. Let us denote by $\textnormal{Def}_{T_{\mathcal{F}_0}}$ the functor of Artin rings defined by
\begin{align*}
\textnormal{Def}_{T_{\mathcal{F}_0}}:\bold{Art} &\to \left(\textnormal{Sets} \right) \\
 A &\mapsto \left\{ \textnormal{locally free sheaves on $X\times \textnormal{Spec}(A)$ inducing  $T_{\mathcal{F}_0}$}   \right\}/\textnormal{isomorphisms}
\end{align*}
Then we have a morphism of deformation functors
\begin{align*}
\textnormal{FDef}_{T_{\mathcal{F}_0}} \to \textnormal{Def}_{T_{\mathcal{F}_0}}
\end{align*}
by forgetting foliated structures. We also have a canonical morphism of deformation functors 
\begin{align*}
\textnormal{FDef}_{\left( X, T_{\mathcal{F}_0} \right)_{\textnormal{fixed}}}\to \textnormal{FDef}_{T_{\mathcal{F}_0}}
\end{align*}
Then the relationship of those three functors $\textnormal{FDef}_{\left(X, T_{\mathcal{F}_0} \right)}, \textnormal{FDef}_{T_{\mathcal{F}_0}}$ and $\textnormal{Def}_{T_{\mathcal{F}_0}}$ is expressed by the following exact sequence of complexes of sheaves
{\small{\begin{center}
$\begin{CD}
@. \cdots @. \cdots @. \cdots \\
@. @AD_3AA @AD_3AA @AAA \\
0 @>>> \mathscr{H}om_{\mathcal{O}_X}\left( \bigwedge^3 T_{\mathcal{F}_0}, \frac{T_X}{T_{\mathcal{F}_0}} \right) @>>> \mathscr{H}om_{\mathcal{O}_X}\left( \bigwedge^3 T_{\mathcal{F}_0}, \frac{T_X}{T_{\mathcal{F}_0}} \right) @>>> 0 @>>> 0\\
@. @AD_2AA @AD_2AA @AAA \\
0 @>>> \mathscr{H}om_{\mathcal{O}_X}\left( \bigwedge^2 T_{\mathcal{F}_0}, \frac{T_X}{T_{\mathcal{F}_0}} \right) @>>> \mathscr{H}om_{\mathcal{O}_X}\left( \bigwedge^2 T_{\mathcal{F}_0}, \frac{T_X}{T_{\mathcal{F}_0}} \right) @>>> 0 @>>> 0\\
@. @AD_1'AA @AD_1'AA @AAA \\
0 @>>> \mathscr{H}om_{\mathcal{O}_X}\left(T_{\mathcal{F}_0}, T_X  \right) @>>>  \mathscr{H}om_{\mathcal{O}_X}\left(T_{\mathcal{F}_0}, T_X \right) @>>> 0 @>>> 0\\
@. @AAA @A jAA @AAA \\
0 @>>> 0 @>>> \mathscr{H}om_{\mathcal{O}_X}\left( T_{\mathcal{F}_0}, T_{\mathcal{F}_0} \right) @>>> \mathscr{H}om_{\mathcal{O}_X}\left( T_{\mathcal{F}_0},  T_{\mathcal{F}_0} \right) @>>> 0
\end{CD}$
\end{center}}}
Then from the long exact sequence of cohomology groups, we have
\begin{align*}
 \mathbb{H}^0\left(X, \mathscr{H}om_{\mathcal{O}_X}\left(T_{\mathcal{F}_0}, T_X \right)^\bullet \right) \to \mathbb{H}^1\left( X, \mathscr{H}om_{\mathcal{O}_X}\left(T_{\mathcal{F}_0} , T_{\mathcal{F}_0} \right)^\bullet \right) \to H^1\left(X, \mathscr{H}om_{\mathcal{O}_X}\left(T_{\mathcal{F}_0}, T_{\mathcal{F}_0} \right) \right)
\end{align*}
which gives a morphism of first-order deformations 
\begin{align*}
\textnormal{FDef}_{\left( X, T_{\mathcal{F}_0} \right)_{\textnormal{fixed}}}\left(k[\epsilon] \right)\to \textnormal{FDef}_{T_{\mathcal{F}_0}} \left( k[\epsilon] \right) \to \textnormal{Def}_{T_{\mathcal{F}_0}} \left( k[\epsilon ] \right)
\end{align*}

and the long exact sequence also gives an exact sequence 
\begin{align*}
\mathbb{H}^1\left(X, \mathscr{H}om_{\mathcal{O}_X}\left(T_{\mathcal{F}_0}, T_X \right)^\bullet \right) \to \mathbb{H}^2\left( X, \mathscr{H}om_{\mathcal{O}_X}\left(T_{\mathcal{F}_0} , T_{\mathcal{F}_0} \right)^\bullet \right) \to H^2\left(X, \mathscr{H}om_{\mathcal{O}_X}\left(T_{\mathcal{F}_0}, T_{\mathcal{F}_0} \right) \right)
\end{align*}
which gives obstruction maps of morphisms of deformation functors $\textnormal{FDef}_{\left( X, T_{\mathcal{F}_0} \right)_{\textnormal{fixed}}} \to \textnormal{FDef}_{T_{\mathcal{F}_0}} \to \textnormal{Def}_{T_{\mathcal{F}_0}}$.

\subsection{Relationship between foliated deformations of $\left(X, T_{\mathcal{F}_0} \right)$ with $X$ fixed and the underlying flat deformations of $\frac{T_X}{T_{\mathcal{F}_0}}$}\

By definition of foliated deformations of $\left(X, T_{\mathcal{F}_0} \right)$ with $T_{\mathcal{F}_0}$ locally free under the trivial deformation of $X$, it also involves flat deformations of the underlying coherent sheaf $\frac{T_X}{T_{\mathcal{F}_0}}$. We denote by $\textnormal{Def}_{\frac{T_X}{T_{\mathcal{F}_0}}}$ the functor of Artin rings of flat deformations of $T_X/T_{\mathcal{F}_0}$ on the fixed underlying nonsingular variety $X$. Then we have a natural morphism of deformation functor 
\begin{align*}
\textnormal{FDef}_{T_{\mathcal{F}_0}} \to \textnormal{Def}_{\frac{T_X}{T_{\mathcal{F}_0}}}
\end{align*}
We will identify the relationships between first-order deformations and obstructions of two deformation functors. This subsection is a direct consequence of \cite{Har10} Theorem 7.2 and Theorem 7.3 p.54-p.56.

We note that we have an exact sequence of complex of sheaves
{\small{\begin{center}
$\begin{CD}
@. \cdots @. \cdots @. \cdots \\
@. @AAA @AAA @AAA \\
0@>>> \mathscr{H}om_{\mathcal{O}_X}\left(\bigwedge^3 T_{\mathcal{F}_0}, \frac{T_X}{T_{\mathcal{F}_0}} \right) @>>> \mathscr{H}om_{\mathcal{O}_X}\left( \bigwedge^3 T_{\mathcal{F}_0}, \frac{T_X}{T_{\mathcal{F}_0}}\right) @>>> 0 @>>> 0 \\
@. @AAA @AAA @AAA \\
0@>>> \mathscr{H}om_{\mathcal{O}_X}\left(\bigwedge^2 T_{\mathcal{F}_0}, \frac{T_X}{T_{\mathcal{F}_0}} \right) @>>> \mathscr{H}om_{\mathcal{O}_X}\left( \bigwedge^2 T_{\mathcal{F}_0}, \frac{T_X}{T_{\mathcal{F}_0}}\right) @>>> 0 @>>> 0 \\
@. @AAA @AAA @AAA \\
0@>>> 0 @>>> \mathscr{H}om_{\mathcal{O}_X}\left(T_{\mathcal{F}_0}, \frac{T_X}{T_{\mathcal{F}_0}} \right) @>>> \mathscr{H}om_{\mathcal{O}_X}\left(T_{\mathcal{F}_0},  \frac{T_X}{T_{\mathcal{F}_0}}\right) @>>> 0
\end{CD}$
\end{center}}}  

Then we have the natural indued morphisms
\begin{align}
\mathbb{H}^0\left(X, \mathscr{H}om_{\mathcal{O}_X}\left(  T_{\mathcal{F}_0}, T_X / T_{\mathcal{F}_0} \right)^\bullet \right) \to H^0\left(X, \mathscr{H}om_{\mathcal{O}_X}\left( T_{\mathcal{F}_0}, T_X / T_{\mathcal{F}_0} \right) \right)\cong \textnormal{Hom}\left( T_{\mathcal{F}_0}, T_X/ T_{\mathcal{F}_0} \right) \label{tr1} \\
\mathbb{H}^1\left(X, \mathscr{H}om_{\mathcal{O}_X}\left( T_{\mathcal{F}_0} , T_X / T_{\mathcal{F}_0} \right)^\bullet \right) \to H^1\left(X, \mathscr{H}om_{\mathcal{O}_X}\left( T_{\mathcal{F}_0}, T_X / T_{\mathcal{F}_0} \right) \right)\cong \textnormal{Ext}^1\left( T_{\mathcal{F}_0} , T_X / T_{\mathcal{F}_0} \right) \notag
\end{align}
Here we note that since $T_{\mathcal{F}_0}$ is locally free, we have
\begin{align*}
\textnormal{Ext}^i\left( T_{\mathcal{F}_0} , T_X / T_{\mathcal{F}_0} \right)\cong \textnormal{Ext}^i\left(\mathcal{O}_X, \left( T_{\mathcal{F}_0}  \right)^*\otimes T_X /T_{\mathcal{F}_0} \right) \cong H^i\left(X, \mathscr{H}om_{\mathcal{O}_X}\left( T_{\mathcal{F}_0}, T_X/ T_{\mathcal{F}_0} \right) \right)
\end{align*}

On the other hand, from the exact sequence $0\to T_{\mathcal{F}_0} \to  T_X \to  T_X /T_{\mathcal{F}_0} \to 0$, we have the exact sequence
\begin{align}
0&\to \textnormal{Hom}\left( T_X  / T_{\mathcal{F}_0} , T_X /  T_{\mathcal{F}_0}  \right)\to  \textnormal{Hom}\left( T_X ,  T_X / T_{\mathcal{F}_0}  \right)\to  \textnormal{Hom}\left(T_{\mathcal{F}_0} , T_X / T_{\mathcal{F}_0} \right)\notag\\
 &\to \textnormal{Ext}^1\left( T_X / T_{\mathcal{F}_0} , T_X /T_{\mathcal{F}_0} \right) \to \textnormal{Ext}^1\left(T_X, T_X / T_{\mathcal{F}_0} \right) \to \textnormal{Ext}^1\left( T_{\mathcal{F}_0} ,T_X/ T_{\mathcal{F}_0}\right) \label{tr2}\\
 &\to \textnormal{Ext}^2\left(T_X / T_{\mathcal{F}_0} , T_X / T_{\mathcal{F}_0} \right) \to \textnormal{Ext}^2\left(T_X ,T_X / T_{\mathcal{F}_0} \right) \to \textnormal{Ext}^2\left( T_{\mathcal{F}_0} , T_X / T_{\mathcal{F}_0} \right)\to \cdots \notag
\end{align}
Now we assume that $X$ is projective. Then we know that $\textnormal{Def}_{T_X/T_{\mathcal{F}_0}}\left(k[\epsilon]\right)\cong \textnormal{Ext}^1\left(T_X / T_{\mathcal{F}_0} , T_X / T_{\mathcal{F}_0} \right)$ and the obstruction space for $\textnormal{Def}_{T_X / T_{\mathcal{F}_0} }$ is  $\textnormal{Ext}^2\left(T_X / T_{\mathcal{F}_0} , T_X /T_{\mathcal{F}_0} \right)$. By combining $(\ref{tr1})$ and $(\ref{tr2})$, we have
\begin{align*}
\mathbb{H}^0\left(X, \mathscr{H}om_{\mathcal{O}_X}\left( T_{\mathcal{F}_0}, T_X/ T_{\mathcal{F}_0} \right)^\bullet \right) \to \textnormal{Hom}\left( T_{\mathcal{F}_0} , T_X /T_{\mathcal{F}_0} \right) \to \textnormal{Ext}^1\left( T_X /T_{\mathcal{F}_0}, T_X / T_{\mathcal{F}_0} \right)
\end{align*}
which gives a morphism of first order deformations between two deformation functors
\begin{align*}
\textnormal{FDef}_{T_{\mathcal{F}_0}}\left(k[\epsilon]\right)\to \textnormal{Def}_{T_X /T_{\mathcal{F}_0} }\left(k[\epsilon] \right)
\end{align*}

On the other hand, we have
\begin{align*}
\mathbb{H}^1\left(X, \mathscr{H}om_{\mathcal{O}_X}\left( T_{\mathcal{F}_0} , T_X /T_{\mathcal{F}_0} \right)^\bullet \right)\to \textnormal{Ext}^1\left( T_{\mathcal{F}_0} , T_X /T_{\mathcal{F}_0}  \right) \to \textnormal{Ext}^2\left(T_X / T_{\mathcal{F}_0} , T_X /T_{\mathcal{F}_0} \right)
\end{align*}
which gives the obstruction map between two deformation functors.

\section{Flat foliated deformations of $\left(X, T_{\mathcal{F}_0} \right)$ with $T_{\mathcal{F}_0}$ locally free}

Let $\left( X, T_{\mathcal{F}_0}\right)$ be a foliated nonsingular variety over $k$ with $T_{\mathcal{F}_0}$ locally free. In this section we now vary the underlying nonsingular variety $X$. We recall that

\begin{definition}[Flat foliated deformations in terms of tangent sheaves]
Let $\left(X, T_{\mathcal{F}_0} \right)$ be a foliated nonsingular variety with $T_{\mathcal{F}_0}$ locally free. An infinitesimal foliated deformation $\eta$ of $\left( X, T_{\mathcal{F}_0} \right)$ over $A\in \bold{Art}$ is a cartesian diagram of morphisms of schemes
\begin{center}
$\eta:\begin{CD}
\left( X, T_{\mathcal{F}_0} \right) @>>> \left( \mathcal{X}, T_\mathcal{F} \right) \\
@VVV @VV\pi V \\
\textnormal{Spec}(k) @>>> \textnormal{Spec}(A)
\end{CD}$
\end{center}
where $\pi$ is flat, $T_{\mathcal{F}}$ is a locally free subsheaf of $T_{\mathcal{X}/\textnormal{Spec}(A)}$ such that $\left[ T_\mathcal{F}, T_\mathcal{F} \right] \subset T_\mathcal{F}$, and $\left(\mathcal{X}, T_\mathcal{F} \right)$ induces $\left( X, T_{\mathcal{F}_0}\right)$. Two infinitesimal foliated deformations $\left( \mathcal{X}, T_\mathcal{F} \right)$ and $\left(\mathcal{X}', T_{\mathcal{F}'} \right)$ of $\left( X, T_{\mathcal{F}_0} \right)$ over $A$ are isomorphic if there is an isomorphism $\phi:\mathcal{X}\to \mathcal{X}'$ over $A$ inducing the identity on $\left(X, T_{\mathcal{F}_0} \right)$ such that $\phi^* T_{\mathcal{F}'}= T_\mathcal{F}$. Then we have a functor of Artin rings
\begin{align*}
\textnormal{FDef}_{\left( X, T_{\mathcal{F}_0}  \right)}: \bold{Art} &\to \left(\textnormal{Sets} \right)\\
 A &\to \left\{ \textnormal{Infinitesimal foliated deformations of $\left(X, T_{\mathcal{F}_0} \right)$ over $A$}  \right\}/\textnormal{isomorphisms}
\end{align*}
\end{definition}

We will show that the leaf complex (see Part I)
\begin{align}\label{at16}
T_{\mathcal{F}_0}^\bullet: T_X \xrightarrow{D_0} \mathscr{H}om_{\mathcal{O}_X}\left( T_{\mathcal{F}_0}, \frac{T_X}{T_{\mathcal{F}_0}} \right) \xrightarrow{D_1} \mathscr{H}om_{\mathcal{O}_X}\left( \bigwedge^2 T_{\mathcal{F}_0}, \frac{T_X}{T_{\mathcal{F}_0}} \right) \xrightarrow{D_2}  \mathscr{H}om_{\mathcal{O}_X}\left( \bigwedge^3 T_{\mathcal{F}_0}, \frac{T_X}{T_{\mathcal{F}_0}} \right)\xrightarrow{D_3} \cdots
\end{align}
controls foliated deformations of $\left( X, T_{\mathcal{F}_0} \right)$.

\begin{lemma}[compare \cite{Ser06} Lemma 1.2.6 p.26] \label{at17}
Let $B_0$ be a $k$-algebra with an element $T_0\in \textnormal{Hom}_{B_0}\left(\Omega_{B_0/k}^1 , B_0 \right)$, and $e:0\to (t)\to \tilde{A}\to A \to 0$ be a small extension in $\bold{Art}$. Let $T\in \textnormal{Hom}_{B_0}\left(\Omega_{B_0/k}^1, B_0 \right)\otimes A$ which induces $T_0$, and $T_1,T_2\in \textnormal{Hom}_{B_0}\left( \Omega_{B_0/k}^1, B_0 \right)\otimes \tilde{A}$ which induces $T$, so that there exists $T'\in \textnormal{Hom}_{B_0}\left(\Omega_{B_0/k}^1, B_0 \right)$ such that $T_1-T_2= tT'$. Then there is one to one correspondence
\begin{align*}
\left\{ \textnormal{isomorphisms on $\left( B_0\otimes_k \tilde{A}, T_1 \right)$ and $\left( B_0\otimes_k \tilde{A}, T_2\right)$ inducing the identity on $\left(B_0\otimes_k  A,  T \right)$} \right\}\\
\to \left\{ P\in \textnormal{Der}_k\left(B_0, B_0 \right)=\textnormal{Hom}_{B_0}\left(\Omega_{B_0/k}^1, B_0 \right)| T'+\left[ P, T_0\right]=0 \right\}
\end{align*}
In particular, when $T_1=T_2$, there is a canonical isomorphism of groups
\begin{align*}
\left\{ \textnormal{automorphisms on $\left(B_0\otimes_k \tilde{A}, T_1 \right)$ inducing the identity on $\left(B_0\otimes_k A, T \right)$} \right\} \to \left\{P\in \textnormal{Der}_k\left(B_0, B_9 \right) |\left[ P, T_0\right]=0 \right\}
\end{align*}
\end{lemma}

\begin{proof}
Let $\theta:\left( B_0\otimes_k \tilde{A}, T_1 \right)\to \left( B_0\otimes_k \tilde{A}, T_2 \right)$ be an isomorphism sending $T_1$ to $T_2$ which induces the identity on $\left( B_0\otimes_k A, T\right)$. Then we have $\theta(x)= x+ tP\left(dx\right)$ for some $P\in \textnormal{Der}_k\left( B_0, B_0\right)= \textnormal{Hom}_{B_0}\left( \Omega_{B_0/k}^1, B_0 \right)$. We define the correspondence by $\theta\mapsto P$. On the other hand, since $\theta$ is compatible with $T_1,T_2$, we see that for $x \in B_0\otimes_k \tilde{A}$, we have
\begin{align*}
&\theta\left( T_1(dx) \right) = T_2 \left( d\theta(x)\right) \iff T_1\left(d x \right) + t P\left( d\left(  T_1(dx)\right) \right) = T_2\left( d x+ t d\left( P(dx) \right)   \right) \\
&\iff \left(T_1- T_2 \right)(dx) = t \left[ T_0, P \right](dx) \iff tT'(dx)=t \left[ T_0,P \right](dx) \iff T'+\left[ P, T_0 \right]=0
\end{align*}
Since $\theta$ is determined by $P$, the correspondence is one to one. Now assume that $T_1=T_2$. Then $\theta$ corresponds to $P$ with $\left[ P, T_0 \right]=0$ and the correspondence $\theta\to P$ is a group isomorphism.
\end{proof}

\begin{proposition}\label{at25}
Let $\left( X, T_{\mathcal{F}_0} \right)$ be a foliated nonsingular variety with $T_{\mathcal{F}_0}$ locally free.
\begin{enumerate}
\item there is a $1-1$ correspondence
\begin{align*}
\textnormal{FDef}_{\left(X, T_{\mathcal{F}_0}\right)}\left(k[\epsilon] \right)\cong \mathbb{H}^1\left(X, T_{\mathcal{F}_0}^\bullet \right)
\end{align*}
\item Let $A\in \bold{Art}$ and $\eta=\left( \mathcal{X}, T_\mathcal{F} \right)$ be an infinitesimal foliated deformations of $\left( X, T_{\mathcal{F}_0} \right)$ over $A$. To every small extension $0\to (t)\to \tilde{A} \to A \to 0$, we can associate an element $o_\eta(e)\in \mathbb{H}^2\left( X,  T_{\mathcal{F}_0}^\bullet \right)$ called the obstruction lifting of $\eta$ to $\tilde{A}$, which is $0$ if and only if a lifting of $\eta$ to $\tilde{A}$ exists.
\end{enumerate}
\end{proposition}

\begin{proof}
Given a first-order foliated deformation $\left(\mathcal{X}, T_\mathcal{F} \right)$ of $\left(X, T_{\mathcal{F}_0} \right)$ over $k[\epsilon]$, we choose an affine open covering $\mathcal{U}=\{U_i\}$ of $X$ such that $\Gamma\left( U_i, T_{\mathcal{F}_0}\right)$ is generated by $T_{0i}^1,..., T_{0i}^p \in \Gamma\left( U_i, T_X \right)$ with $T_{0i}^\alpha =\sum_{\beta=1}^p r_{0ij}^{\alpha\beta}T_{0j}^\beta$ for some $r_{0ij}^{\alpha\beta}\in \Gamma\left( U_{ij}, \mathcal{O}_X\right),\alpha=1,...,p$ and $\left[ T_{0i}^\alpha, T_{0i}^\beta \right] = \sum_{\gamma=1}^p g_{0i\alpha\beta}^\gamma T_{0i}^\gamma$ for some $g_{0i\alpha\beta}^\gamma \in \Gamma\left(U_i, \mathcal{O}_X \right), \alpha,\beta=1,...,p$, and we have $\mathcal{X}|_{U_i}\cong U_i\times \textnormal{Spec}(k[\epsilon])$ is trivial for all $i$ with the induced foliated structure $T_{\mathcal{F}_i}$ on $U_i\times \textnormal{Spec}(k[\epsilon])$, and $\Gamma(U_i, T_{\mathcal{F}_i})$ is generated by $T_{0i}^\alpha+\epsilon T_i'^\alpha,\alpha=1,...,p$, where $T_i'^\alpha\in \Gamma(U_i, T_X)$, and $\left[ T_{0i}^\alpha + \epsilon T_{i}'^\alpha, T_{0i}^\beta + \epsilon  T_{i}'^\beta \right] = \sum_{\gamma=1}^p \left( g_{0i\alpha\beta}^\gamma + \epsilon c_{i\alpha\beta}^\gamma \right)\left(T_{0i}^\gamma + \epsilon T_i'^\gamma \right)$ for some $c_{i\alpha\beta}^\gamma \in \Gamma\left( U_i, \mathcal{O}_X \right)$. For all $i$, we have a foliated isomorphism
\begin{align*}
\theta_i: \left(U_i\times \textnormal{Spec}\left(k[\epsilon]\right), T_{\mathcal{F}_i}\right)\to \left(\mathcal{X}|_{U_i}, T_\mathcal{F}|_{U_i} \right)
\end{align*}
Then for each $i,j,\theta_{ij}:=\theta_j^{-1}\theta_i:\left(U_{ij}\times \textnormal{Spec}(k[\epsilon]), T_{\mathcal{F}_i}\right)\to \left(U_{ij}\times \textnormal{Spec}\left(k[\epsilon]\right), T_{\mathcal{F}_j}\right)$ is a foliated isomorphism inducing the identity on $\left(U_{ij}, T_{\mathcal{F}_0}\right)$ by modulo $\epsilon$. Then by Lemma \ref{at17}, $\theta_{ij}$ correspondes to $\textnormal{Id}+\epsilon p_{ij}:\left(\Gamma\left(U_{ij}, \mathcal{O}_X \right)\otimes k[\epsilon], T_{\mathcal{F}_j} \right)\to \left(\Gamma\left( U_{ij}, \mathcal{O}_X \right)\otimes k[\epsilon], T_{\mathcal{F}_i} \right)$ for some $p_{ij}\in \Gamma\left( U_{ij}, T_X\right)$ such that
\begin{align}\label{at18}
\theta_{ji}^*\left( T_{0i}^\alpha + \epsilon T_i'^\alpha  \right)=T_{0i}^\alpha + \epsilon T_{i}'^\alpha - \epsilon \left[ p_{ij},  T_{0i}^\alpha \right] = \sum_{\beta=1}^p \left( r_{0ij}^{\alpha\beta} + \epsilon b_{ij}^{\alpha\beta}  \right)\left( T_{0j}^\beta + \epsilon T_j'^\beta \right)
\end{align}
for some $b_{ij}^{\alpha\beta}\in \Gamma\left( U_{ij}, \mathcal{O}_X \right)$. Then $(\ref{at18})$ implies that 
\begin{align}\label{at19}
T_i'^\alpha -\left[ p_{ij}, T_{0i}^\alpha \right]=\sum_{\beta=1}^p r_{0ij}^{\alpha\beta} T_j'^\beta + \sum_{\beta=1}^p b_{ij}^{\alpha\beta} T_{0j}^\beta
\end{align}
We define an element $\overline{T_i'} \in \Gamma\left( U_i, \mathscr{H}om_{\mathcal{O}_X}\left( T_{\mathcal{F}_0}, \frac{T_X}{T_{\mathcal{F}_0}} \right) \right)$ as in $(\ref{at20})$. Then $(\ref{at7})$ and $(\ref{at19})$ implies that $\left( \left\{ \overline{T_i'} \right\}, \left\{ p_{ij} \right\} \right) \in C^0\left( \mathcal{U}, \mathscr{H}om_{\mathcal{O}_X}\left( T_{\mathcal{F}_0}, \frac{T_X}{T_{\mathcal{F}_0}} \right) \right) \bigoplus C^1 \left( \mathcal{U}, T_X \right)$ defines a $1$-cocycle in the following \v Cech resolution of $T_{\mathcal{F}_0}^\bullet$:
{\Small{\begin{center}
$\begin{CD}
\cdots\\
@AD_3AA \\
C^0\left(\mathcal{U}, \mathscr{H}om_{\mathcal{O}_X}\left( \bigwedge^3 T_{\mathcal{F}_0}, \frac{T_X}{T_{\mathcal{F}_0}}\right) \right)@>\delta>> \cdots\\
@AD_2AA @AD_2AA \\
C^0\left(\mathcal{U}, \mathscr{H}om_{\mathcal{O}_X}\left(\bigwedge^2 T_{\mathcal{F}_0}, \frac{T_X}{T_{\mathcal{F}_0}} \right) \right) @>-\delta>> C^1\left(\mathcal{U}, \mathscr{H}om_{\mathcal{O}_X}\left(  \bigwedge^2 T_{\mathcal{F}_0}, \frac{T_X}{T_{\mathcal{F}_0}} \right) \right) @>\delta>> \cdots\\ 
@AD_1AA @AD_1AA @AD_1AA \\
C^0\left(\mathcal{U}, \mathscr{H}om_{\mathcal{O}_X}\left( T_{\mathcal{F}_0}, \frac{T_X}{T_{\mathcal{F}_0}}\right) \right)@>\delta>> C^1\left(\mathcal{U}, \mathscr{H}om_{\mathcal{O}_X}\left( T_{\mathcal{F}_0}, \frac{T_X}{T_{\mathcal{F}_0}}\right) \right) @>-\delta>> C^2\left(\mathcal{U}, \mathscr{H}om_{\mathcal{O}_X}\left( T_{\mathcal{F}_0}, \frac{T_X}{T_{\mathcal{F}_0}}\right) \right) \\
@AD_0AA @AD_0AA @AD_0AA \\
C^0\left(\mathcal{U},T_X\right)@>-\delta>> C^1\left(\mathcal{U},T_X \right) @>\delta>> C^2 \left(\mathcal{U},T_X\right) @>-\delta>> \cdots
\end{CD}$
\end{center}}}
We show that for two equivalent first-order foliated deformations of $\left( X, T_{\mathcal{F}_0}\right)$, the associated $1$-cocycles are equivalent. Let $\left( \mathcal{X}', T_{\mathcal{F}'} \right)$ be another foliated deformations of $\left( X, T_{\mathcal{F}_0} \right)$ over $k[\epsilon]$ which induces $\left( \left\{ \overline{T_i''} \right\}, \left\{ p_{ij}'\right\}  \right)$ and $\Phi:\left( \mathcal{X}, T_\mathcal{F} \right) \to \left( \mathcal{X}', T_{\mathcal{F}'} \right)$ is a foliated isomorphism inducing the identity on $\left(X, T_{\mathcal{F}_0} \right)$. Then there is a foliated isomorphism $\alpha_i$ on $U_i\times \textnormal{Spec}\left(k[\epsilon ] \right)$ given by 
\begin{align*}
\alpha_i:\left(U_i\times \textnormal{Spec}\left(k[\epsilon]\right) , T_{\mathcal{F}_i} \right) \xrightarrow{\theta_i} \left( \mathcal{X}|_{U_i}, T_{\mathcal{F}}|_{U_i} \right) \xrightarrow{\Phi |_{U_i}} \left( \mathcal{X}'|_{U_i}, T_{\mathcal{F}'}|_{U_i} \right) \xrightarrow{\theta_i'^{-1}} \left( U_i \times \textnormal{Spec}\left(k[\epsilon] \right) , T_{\mathcal{F}_i'} \right)
\end{align*}
which corresponds to $1+\epsilon a_i:\left(\Gamma\left( U_i, \mathcal{O}_X  \right)\otimes k[\epsilon], T_{\mathcal{F}_i'} \right)\to \left( \Gamma\left(U_i, \mathcal{O}_X \right) \otimes k[\epsilon] , T_{\mathcal{F}_i} \right)$ for some $a_i\in \Gamma\left( U_i, T_X \right)$ such that $T_{0i}^\alpha + \epsilon T_i'^\alpha - \epsilon \left[ a_i, T_{0i}^\alpha \right] = \sum_{\beta=1}^p \left(\delta_\alpha^\beta + \epsilon b_i^{\alpha\beta} \right) \left( T_{0i}^\beta + \epsilon T_i''^\beta \right)$ for some $b_i^{\alpha\beta} \in \Gamma\left( U_i, \mathcal{O}_X \right)$ by Lemma \ref{at17}. Then we have $\left[ a_i, T_{0i}^\alpha\right]= T_i'^\alpha- T_i''^\alpha - \sum_{\beta=1}^p  b_i^{\alpha\beta} T_{0i}^\beta $. On the other hand, since $\theta_i'\alpha_i= \Phi |_{U_i} \theta_i$, we have $\alpha_j^{-1} \theta_{ij}' \alpha_i=\left(\theta_j' \alpha_j \right)^{-1}\left( \theta_i' \alpha_j \right) = \theta_j^{-1} \Phi |_{U_{ij}}^{-1} \Phi|_{U_{ij}} \theta_i= \theta_{ij}$ which corresponds to $\left(\textnormal{Id}+\epsilon a_i \right)\left(\textnormal{Id}+ \epsilon p_{ij}' \right)\left( \textnormal{Id} - \epsilon a_j \right)= \textnormal{Id}+t p_{ij}$, so that $a_i- a_j=  p_{ij}- p_{ij}'$. This implies that $\left( \left\{ \overline{T_i'} \right\}, \left\{ p_{ij} \right\} \right)$ and $\left( \left\{ \overline{T_i''} \right\}, \left\{  p_{ij}' \right\} \right)$ are cohomologous.

We define an inverse map. Given an element in $\mathbb{H}^1\left(X, T_{\mathcal{F}_0}^\bullet \right)$, we represent it by a \v Cech 1-cocycle $\left(\left\{ \overline{T_i'} \right\}, \left\{ p_{ij} \right\} \right)$ for an affine open covering $\mathcal{U}=\{U_i\}$ of $X$ satisfying $(\ref{at7})$ and $(\ref{at19})$ and $p_{ij}- p_{ik}+ p_{jk}=0$. By reversing the above process, we obtain a first-order foliated deformation of $\left( X, T_{\mathcal{F}_0} \right)$.

Next we identify obstructions. Let $e:0\to (t)\to \tilde{A}\to A\to 0$ be a small extension in $\bold{Art}$ and $\eta:=(\mathcal{X}, T_\mathcal{F})$ be a foliated deformation of $\left(X, T_{\mathcal{F}_0}\right)$ over $A$. Let $\mathcal{U}=\{U_i\}$ be an affine covering of $X$ such that we have foliated isomorphisms $\theta_i: (U_i\times \textnormal{Spec}(A), T_{\mathcal{F}_i})\to (\mathcal{X}|_{U_i}, T_\mathcal{F}|_{U_i} )$, where $T_{\mathcal{F}_i}\subset T_X|_{U_i}\otimes A$ is the induced foliation on $U_i\times \textnormal{Spec}(A)$ from $T_\mathcal{F}$ via $\theta_i$. Let $\Gamma(U_i, T_{\mathcal{F}_i})\subset \Gamma(U_i, T_X)\otimes A$ be generated by $T_i^\alpha,\alpha=1,...,p$ which induces $T_{0i}^\alpha$ modulo $\mathfrak{m}_A$ with $\left[T_i^\alpha, T_i^\beta\right]=\sum_{\gamma=1}^p g_{i\alpha\beta}^\gamma T_i^\gamma$ for some $g_{i\alpha\beta}^\gamma \in\Gamma\left( U_i, \mathcal{O}_X \right)\otimes A$ which induces $\left[ T_{0i}^\alpha, T_{0i}^\beta \right]= \sum_{\gamma=1}^p g_{0i\alpha\beta}^\gamma T_{0i}^\gamma$ modulo $\mathfrak{m}_A$, and $\theta_{ij}:=\theta_j^{-1}\theta_i: \left(U_{ij}\times \textnormal{Spec}(A), T_{\mathcal{F}_i}\right)\to \left(U_{ij}\times \textnormal{Spec}(A), T_{\mathcal{F}_j}\right)$ is a foliated isomorphism with $\theta_{jk}\theta_{ij}=\theta_{ik}$ on $U_{ijk}\times \textnormal{Spec}(A)$ and $\theta_{ij}^*\left(T_{\mathcal{F}_j}\right)=T_{\mathcal{F}_i}$ and $\theta_{ji}^* \left(T_i^\alpha \right)=\sum_{\beta=1}^p r_{ij}^{\alpha\beta} T_j^\beta$ for some $r_{ij}^{\alpha\beta}\in \Gamma\left(U_{ij},\mathcal{O}_X\right)\otimes A$, which induces $T_{0i}^\alpha= \sum_{\beta=1}^p r_{0ij}^{\alpha\beta} T_{0j}^\beta$ modulo $\mathfrak{m}_A$. We let $f_{ij}:\left(\Gamma\left( U_{ij}, \mathcal{O}_X \right)\otimes A, T_{\mathcal{F}_j} \right) \to \left( \Gamma\left( U_{ij}, \mathcal{O}_X \right)\otimes A,  T_{\mathcal{F}_i} \right)$ be the ring homomorphism corresponding to $\theta_{ij}$. Then we have $f_{ij}f_{jk}=f_{ik}$. To give a lifting $\tilde{\eta}$ of $\eta$ to $\tilde{A}$ is equivalent to give a collection $\left\{T_{\tilde{\mathcal{F}}_i} \right\}$ where $\Gamma\left(U_i, T_{\tilde{\mathcal{F}}_i} \right)\subset \Gamma\left( U_i, T_X  \right)\otimes \tilde{A}$ is generated by $\tilde{T}_i^\alpha\in \Gamma\left(U_i, T_X\right)\otimes \tilde{A},\alpha= 1,...,p$ which induces $T_i^\alpha$ such that $\left[\tilde{T}_i^\alpha, \tilde{T}_i^\beta\right]=\sum_{\gamma=1}^p \tilde{g}_{i\alpha\beta}^\gamma \tilde{T}_i^\gamma$ for some $\tilde{g}_{i\alpha\beta}^\gamma\in \Gamma\left(U_i,\mathcal{O}_X\right)\otimes \tilde{A}$ which is a lifting of $g_{i\alpha\beta}^\gamma$ and a collection of foliated isomorphisms $\left\{ \tilde{\theta}_{ij} \right\}$ which induces $\theta_{ij}$, where $\tilde{\theta}_{ij}: \left(U_{ij}\times \textnormal{Spec}\left(\tilde{A}\right), T_{\tilde{\mathcal{F}}_i}\right)\to \left(U_{ij}\times \textnormal{Spec}\left(\tilde{A}\right), T_{\tilde{\mathcal{F}}_j}\right)$ such that $\tilde{\theta}_{ji}^*\left(\tilde{T}_i^\alpha\right)=\sum_{\beta=1}^p \tilde{r}_{ij}^{\alpha\beta} \tilde{T}_j^\beta$ for some $\tilde{r}_{ij}^{\alpha\beta}\in \Gamma(U_{ij}, \mathcal{O}_X)\otimes \tilde{A}$ which induces $r_{ij}^{\alpha\beta}$, and $\tilde{\theta}_{jk}\tilde{\theta}_{ij}=\tilde{\theta}_{ik}$ (equivalently $\tilde{f}_{ij}\tilde{f}_{jk}=\tilde{f}_{ik})$. From such data, we can glue together $\left(U_i\times \textnormal{Spec}(\tilde{A}), T_{\tilde{\mathcal{F}}_i}\right)$ to make a foliated deformation $\left(\tilde{X}, T_{\tilde{\mathcal{F}}}\right)$ inducing $\left(\mathcal{X}, T_\mathcal{F}\right)$. Now given a foliated deformation $\eta=\left(\mathcal{X}, T_\mathcal{F}\right)$ over $A$ as above, we associate an element $o_\xi(e)\in \mathbb{H}^2 \left(X, T_{\mathcal{F}_0}^\bullet \right)$. In order to see if there exists a lifting of $\left(\mathcal{X}, T_\mathcal{F} \right)$ to $\tilde{A}$, we choose arbitrary automorphism $\left\{\tilde{\theta}_{ij}\right\}$ which restricts to $\theta_{ij}$ on $U_{ij}\times \textnormal{Spec}\left(A\right)$ and arbitrary lifting $\tilde{T}_i^\alpha\in \Gamma\left(U_i, T_X\right)\otimes \tilde{A}$ which restricts to $T_i^\alpha$ and arbitrary lifting $\tilde{r}_{ij}^{\alpha\beta} \in \Gamma\left( U_{ij}, \mathcal{O}_X \right)\otimes \tilde{A}$ of $r_{ij}^{\alpha\beta}$ and arbitrary lifting $\tilde{g}_{i\alpha\beta}^\gamma \in \Gamma\left( U_i, \mathcal{O}_X \right) \otimes \tilde{A}$ of $g_{i\alpha\beta}^\gamma$ for $\alpha,\beta=1,...,p$. Since $\tilde{\theta}_{ijk}=\tilde{\theta}_{ki}\tilde{\theta}_{jk}\tilde{\theta}_{ij}$ (equivalently $\tilde{f}_{ijk}=\tilde{f}_{ij}\tilde{f}_{jk}\tilde{f}_{ki}$) is an automorphism on $U_{ijk}\times \textnormal{Spec}\left(\tilde{A}\right)$ inducing the identity on $U_{ijk}\times \textnormal{Spec}(A)$, $\tilde{\theta}_{ijk}$ corresponds to $\tilde{d}_{ijk}\in \Gamma\left(U_{ijk}, T_X\right)$ with $\tilde{d}_{jkl}-\tilde{d}_{ikl}+\tilde{d}_{ijl}-\tilde{d}_{jkl}=0$, so that $\delta(\{\tilde{d}_{ijk}\})=0$. On the other hand, we have
\begin{align}
\tilde{\theta}_{ji}^* \left(\tilde{T}_i^\alpha \right) -\sum_{\beta=1}^p \tilde{r}_{ij}^{\alpha\beta} \tilde{T}_j^\beta= t W_{ij}^{\alpha}\,\,\,\,\,\,\,\,\,\,&\text{for some}\,\,\,  W_{ij}^{\alpha}\in \Gamma(U_{ij}, T_X) \label{ts5}\\
\left[\tilde{T}_i^\alpha, \tilde{T}_i^\beta\right]-\sum_{\eta=1}^p \tilde{g}_{i\alpha\beta}^{\eta} \tilde{T}_i^\eta = t P_i^{\alpha\beta}\,\,\,\,\,\,\,\,&\textnormal{for some}\,\,\,\,\,P_i^{\alpha\beta}\in \Gamma(U_i, T_X) \notag
\end{align}
We define $W_{ij}\in \Gamma\left( U_{ij}, \mathscr{H}om_{\mathcal{O}_X}\left(  T_{\mathcal{F}_0}, \frac{T_X}{T_{\mathcal{F}_0}}\right) \right)$ in a similar way to $(\ref{at21})$ and $P_i\in \Gamma\left( U_i, \mathscr{H}om_{\mathcal{O}_X}\left( \bigwedge^2 T_{\mathcal{F}_0}, \frac{T_X}{T_{\mathcal{F}_0} } \right) \right)$ in a similar way to $(\ref{at22})$. We note that by applying $\theta_{ij}^*$ to $\theta_{ji}^*\left( T_i^\alpha \right)= \sum_{\beta=1}^p r_{ij}^{\alpha\beta} T_j^\beta$, we have $T_i^\alpha= \sum_{\beta=1}^p f_{ij}\left(r_{ij}^{\alpha\beta} \right)\left(\sum_{\gamma=1}^p r_{ji}^{\beta \gamma} T_i^\gamma \right)$, so that $\sum_{\beta=1}^p f_{ij}\left( r_{ij}^{\alpha\beta} \right)r_{ji}^{\beta \gamma}= \delta_\alpha^\gamma$. Then we have $\sum_{\beta=1}^p \tilde{f}_{ij}\left( \tilde{r}_{ij}^{\alpha\beta} \right) \tilde{r}_{ji}^{\beta \gamma}= \delta_\alpha^\gamma + t b_{ij\alpha}^\gamma$ for some $b_{ij\alpha}^\gamma \in \Gamma\left( U_{ij}, \mathcal{O}_X \right)$. Then we have
{\small{\begin{align}\label{atb1}
tW_{ij}^\alpha&= \tilde{\theta}_{ij}^*\left(t W_{ij}^\alpha \right)=\tilde{\theta}_{ij}^*\tilde{\theta}_{ji}^*\left(\tilde{T}_i^\alpha\right)-\sum_{\beta=1}^p \tilde{f}_{ij}\left(\tilde{r}_{ij}^{\alpha\beta}\right)\tilde{\theta}_{ij}^*\left(\tilde{T}_j^\beta\right)=\tilde{T}_i^\alpha-\sum_{\beta=1}^p \tilde{f}_{ij} \left( \tilde{r}_{ij}^{\alpha\beta} \right) \left(tW_{ji}^\beta+\sum_{\gamma=1}^p \tilde{r}_{ji}^{\beta\gamma} \tilde{T}_i^\gamma \right) \\
   &=\tilde{T}_i^\alpha-t\sum_{\beta=1}^p r_{0ij}^{\alpha\beta} W_{ji}^\beta -\sum_{\beta,\gamma=1}^p  \tilde{f}_{ij}\left(\tilde{r}_{ij}^{\alpha\beta} \right) \tilde{r}_{ji}^{\beta\gamma} \tilde{T}_i^\gamma \notag \\
                    &=\tilde{T}_i^\alpha-t\sum_{\beta=1}^p r_{0ij}^{\alpha\beta} W_{ji}^\beta-\sum_{\gamma=1}^p\left(\delta_{\alpha}^\gamma+t b_{ij}^\gamma \right)\tilde{T}_i^\gamma= -t\sum_{\beta=1}^p r_{0ij}^{\alpha\beta} W_{ji}^\beta -t \sum_{\gamma=1}^p b_{ij\alpha}^\gamma T_{0i}^\gamma\Longrightarrow \overline{W_{ij}^\alpha}=-\sum_{\beta=1}^p r_{0ij}^{\alpha\beta}\overline{ W_{ji}^\beta} \notag
\end{align}}}
This implies that $\left\{ W_{ij} \right\} \in C^1\left( \mathcal{U}, \mathscr{H}om_{\mathcal{O}_X}\left( T_{\mathcal{F}_0}, \frac{T_X}{T_{\mathcal{F}_0}} \right) \right)$. Then we claim that $\left( \left\{  P_i \right\}, \left\{ - W_{ij} \right\} , \left\{ d_{ijk} \right\}\right)\in C^0\left( \mathcal{U}, \mathscr{H}om_{\mathcal{O}_X}\left( \bigwedge^2 T_{\mathcal{F}_0}, \frac{T_X}{T_{\mathcal{F}_0}} \right) \right) $ $ \bigoplus C^1\left( \mathcal{U}, \mathscr{H}om_{\mathcal{O}_X}\left( T_{\mathcal{F}_0}, \frac{T_X}{T_{\mathcal{F}_0}} \right) \right) \bigoplus C^2\left(\mathcal{U}, T_X \right)$ defines a $2$-cocycle in the above \v Cech resolution of $T_{\mathcal{F}_0}^\bullet$. We have $-\delta\left( \left\{ d_{ijk} \right\} \right)=0$ as above and we have $D_2\left( P_i \right)=0$ as in the proof of Proposition \ref{at3}. We show that $-\delta\left( \left\{P_i \right\} \right) + D_1\left( \left\{ - W_{ij} \right\} \right)=0$. In fact, first we note that
\begin{align*}
0&=\theta_{ji}^* \left(\left[T_i^\alpha, T_i^\beta \right]-\sum_{\eta=1}^p g_{i\alpha\beta}^\eta T_i^\eta \right)=\sum_{\gamma,\eta=1}^p \left[ r_{ij}^{\alpha\gamma} T_j^\gamma, r_{ij}^{\beta\eta} T_j^\eta \right]-\sum_{\eta,\gamma=1}^p f_{ji}\left(g_{i\alpha\beta}^\eta\right) r_{ij}^{\eta\gamma} T_j^\gamma\\
&= \sum_{\gamma,\eta=1}^p r_{ij}^{\alpha\gamma} r_{ij}^{\beta\eta}\left[T_j^\gamma, T_j^\eta \right] +\sum_{\gamma,\eta=1}^p r_{ij}^{\alpha\gamma} T_j^\gamma\left(r_{ij}^{\beta\eta}\right) T_j^\eta -\sum_{\gamma,\eta=1}^p r_{ij}^{\beta\eta} T_j^\eta\left(r_{ij}^{\alpha\gamma}\right) T_j^\gamma-\sum_{\eta,\gamma=1}^p f_{ji}\left(g_{i\alpha\beta}^\eta \right) r_{ij}^{\eta\gamma} T_j^\gamma\\
&=\sum_{\xi=1}^p\left(\sum_{\gamma,\eta=1}^p r_{ij}^{\alpha\gamma}r_{ij}^{\beta\eta} g_{j\gamma\eta}^\xi+\sum_{\gamma=1}^p r_{ij}^{\alpha\gamma} T_j^\gamma\left(r_{ij}^{\beta\xi}\right)-\sum_{\eta=1}^p r_{ij}^{\beta\eta} T_j^\eta\left(r_{ij}^{\alpha\xi}\right)-\sum_{\eta=1}^p f_{ji}\left(g_{i\alpha\beta}^\eta \right) r_{ij}^{\eta\xi} \right) T_j^\xi,
\end{align*}
so that we have
\begin{align}
\sum_{\gamma,\eta=1}^p r_{ij}^{\alpha\gamma}r_{ij}^{\beta\eta} g_{j\gamma\eta}^\xi+\sum_{\gamma=1}^p r_{ij}^{\alpha\gamma} T_j^\gamma\left(r_{ij}^{\beta\xi}\right)-\sum_{\eta=1}^p r_{ij}^{\beta\eta} T_j^\eta \left(r_{0ij}^{\alpha\xi}\right)-\sum_{\eta=1}^p f_{ji}\left( g_{i\alpha\beta}^\eta \right) r_{ij}^{\eta\xi}=0 
\end{align}
Then we have
\begin{align}\label{at24}
\sum_{\gamma,\eta=1}^p \tilde{r}_{ij}^{\alpha\gamma} \tilde{r}_{ij}^{\beta\eta} \tilde{g}_{j\gamma\eta}^\xi+\sum_{\gamma=1}^p \tilde{r}_{ij}^{\alpha\gamma} \tilde{T}_j^\gamma\left( \tilde{r}_{ij}^{\beta\xi}\right)-\sum_{\eta=1}^p \tilde{r}_{ij}^{\beta\eta} \tilde{T}_j^\eta\left( \tilde{r}_{ij}^{\alpha\xi}\right)-\sum_{\eta=1}^p \tilde{f}_{ji}\left( \tilde{g}_{i\alpha\beta}^\eta\right) \tilde{r}_{ij}^{\eta\xi} = t e_{ij}^{\alpha\beta\xi}\,\,\,\,\,\textnormal{for some}\,\,\, e_{ij}^{\alpha\beta \xi}\in \Gamma\left( U_i, \mathcal{O}_X \right)
\end{align}
Then from $(\ref{at24})$, we have
{\small{\begin{align*}
&tP_i^{\alpha\beta}= \tilde{\theta}_{ji}^*\left( t P_{ij}^{\alpha\beta} \right)=\tilde{\theta}_{ji}^*\left( \left[ \tilde{T}_i^\alpha, \tilde{T}_i^\beta \right] - \sum_{\eta=1}^p \tilde{g}_{i\alpha\beta}^\eta \tilde{T}_i^\eta \right) =\left[ \tilde{\theta}_{ji}^*\left(\tilde{T}_i^\alpha\right), \tilde{\theta}_{ji}^* \left(\tilde{T}_i^\beta\right) \right]-\sum_{\eta=1}^p \tilde{f}_{ji}\left(\tilde{g}_{i\alpha\beta}^\eta\right) \tilde{\theta}_{ji}^* \left(\tilde{T}_i^\eta\right)\\
&=\left[tW_{ij}^\alpha+\sum_{\gamma=1}^p  \tilde{r}_{ij}^{\alpha\gamma} \tilde{T}_j^\gamma, tW_{ij}^\beta+\sum_{\eta=1}^p \tilde{r}_{ij}^{\beta\eta} \tilde{T}_j^\eta \right] -\sum_{\eta=1}^p\tilde{f}_{ji}\left(\tilde{g}_{i\alpha\beta}^\eta\right)\left(tW_{ij}^\eta+\sum_{\gamma=1}^p \tilde{r}_{ij}^{\eta\gamma}\tilde{T}_j^\gamma \right)\\
&=t\left[W_{ij}^\alpha, T_{0i}^\beta \right]+t\left[T_{0i}^\alpha, W_{ij}^\beta \right]-t\sum_{\eta=1}^p g_{0i\alpha\beta}^\eta W_{ij}^\eta +\sum_{\gamma,\eta=1}^p \tilde{r}_{ij}^{\alpha\gamma} \tilde{r}_{ij}^{\beta\eta}\left(\left[\tilde{T}_j^\gamma, \tilde{T}_j^\eta \right] -\sum_{\xi=1}^p \tilde{g}_{j\gamma\eta}^\xi \tilde{T}_j^\xi \right) +\sum_{\gamma,\eta,\xi=1}^p \tilde{r}_{ij}^{\alpha\gamma} \tilde{r}_{ij}^{\beta\eta} \tilde{g}_{j\gamma\eta}^\xi \tilde{T}_j^\xi\\
&+\sum_{\gamma,\eta=1}^p \tilde{r}_{ij}^{\alpha\gamma}\tilde{T}_j^\gamma \left( \tilde{r}_{ij}^{\beta\eta}\right)\tilde{T}_j^\eta - \sum_{\gamma,\eta=1}^p \tilde{r}_{ij}^{\beta\eta} \tilde{T}_j^\eta\left( \tilde{r}_{ij}^{\alpha\gamma} \right) \tilde{T}_j^\gamma -\sum_{\eta,\gamma=1}^p \tilde{f}_{ji}\left(\tilde{g}_{i\alpha\beta}^\eta \right) \tilde{r}_{ij}^{\eta\gamma} \tilde{T}_j^\gamma\\
&=t\left(\left[W_{ij}^\alpha, T_{0i}^\beta \right]+\left[ T_{0i}^\alpha, W_{ij}^\beta \right]-\sum_{\eta=1}^p g_{0i\alpha\beta}^\eta W_{ij}^\eta +\sum_{\gamma,\eta=1}^p r_{0ij}^{\alpha\gamma}r_{0ij}^{\beta\eta} P_j^{\gamma\eta}\right)\\
&+\sum_{\xi=1}^p\left(\sum_{\gamma,\eta=1}^p \tilde{r}_{ij}^{\alpha\gamma} \tilde{r}_{ij}^{\beta\eta} \tilde{g}_{j\gamma\eta}^\xi+\sum_{\gamma=1}^p \tilde{r}_{ij}^{\alpha\gamma} \tilde{T}_j^\gamma\left( \tilde{r}_{ij}^{\beta\xi} \right) - \sum_{\eta=1}^p \tilde{r}_{ij}^{\beta\eta} \tilde{T}_j^\eta\left( \tilde{r}_{ij}^{\alpha\xi} \right) -\sum_{\eta=1}^p \tilde{f}_{ji}\left(\tilde{g}_{i\alpha\beta}^\eta\right) \tilde{r}_{ij}^{\eta\xi}\right) \tilde{T}_j^\xi\\
&=t\left(\left[W_{ij}^\alpha, T_{0i}^\beta \right]+\left[ T_{0i}^\alpha, W_{ij}^\beta \right]-\sum_{\eta=1}^p g_{0i\alpha\beta}^\eta W_{ij}^\eta +\sum_{\gamma,\eta=1}^p r_{0ij}^{\alpha\gamma}r_{0ij}^{\beta\eta} P_j^{\gamma\eta}\right)+t\sum_{\xi=1}^p e_{ij}^{\alpha\beta\xi} T_{0j}^\xi
\end{align*}}}
This implies that $-\delta(\{P_i\})+D_1(\{ -W_{ij}\})=0$. On the other hand, we show that $-\delta\left( \left\{ - W_{ij}\right\}\right) + D_0\left(\left\{ d_{ijk} \right\} \right)=0$. In fact, first we note that 
\begin{align*}
T_i^\alpha =\theta_{ij}^* \theta_{jk}^* &\theta_{ki}^*\left( T_i^\alpha \right) = \theta_{ij}^*\theta_{jk}^*\left(\sum_{\beta=1}^p r_{ik}^{\alpha\beta} T_k^\beta \right)= \theta_{ij}^*\left( \sum_{\beta, \gamma=1}^p f_{jk}\left(r_{ik}^{\alpha\beta} \right)  r_{kj}^{\beta \gamma} T_j^\gamma \right) = \sum_{\beta, \gamma, \eta=1}^p \left(f_{ij} f_{jk} \right)\left(r_{ik}^{\alpha\beta} \right)f_{ij}\left( r_{kj}^{\beta \gamma} \right) r_{ji}^{\gamma \eta} T_i^\eta \\
&\Longrightarrow \sum_{\beta, \gamma=1}^p  \left(\tilde{f}_{ij} \tilde{f}_{jk} \right)\left(\tilde{r}_{ik}^{\alpha\beta} \right) \tilde{f}_{ij}\left( \tilde{r}_{kj}^{\beta \gamma} \right) \tilde{r}_{ji}^{\gamma \eta} = \delta_\alpha^\eta + t g_{ijk}^{\alpha \eta}\,\,\,\,\,\,\,\textnormal{for some}\,\,\,  g_{ijk}^{\alpha \eta} \in \Gamma\left( U_{ijk}, \mathcal{O}_X \right)
\end{align*}
Then we note that we have an isomorphism $\textnormal{Id}+ td_{ijk}:\Gamma\left( U_{ijk}, \mathcal{O}_X \right)\otimes \tilde{A}\to \Gamma\left(U_{ijk}, \mathcal{O}_X \right)\otimes\tilde{A}$ corresponding to $\tilde{\theta}_{ki}\tilde{\theta}_{jk} \tilde{\theta}_{ij}$ and so $\tilde{f}_{ij}\tilde{f}_{jk}\tilde{f}_{ki}=\textnormal{Id}+t d_{ijk}$. Then by Lemma \ref{at17}, 
\begin{align*}
\tilde{T}_{i}^\alpha+t\left[ d_{ijk},T_{0i}^\alpha \right]&=\tilde{\theta}_{ij}^*\tilde{\theta}_{jk}^*\tilde{\theta}_{ki}^* \left(\tilde{T}_i^\alpha\right)=\tilde{\theta}_{ij}^*\tilde{\theta}_{jk}^*\left(tW_{ik}^\alpha +\sum_{\beta=1}^p \tilde{r}_{ik}^{\alpha\beta} \tilde{T}_k^\beta \right)=\tilde{\theta}_{ij}^*\left(tW_{ik}^\alpha+\sum_{\beta=1}^p\tilde{f}_{jk}\left(\tilde{r}_{ik}^{\alpha\beta} \right)\left(tW_{kj}^\beta+\sum_{\gamma=1}^p \tilde{r}_{kj}^{\beta\gamma} \tilde{T}_j^\gamma \right) \right)\\
&=tW_{ik}^\alpha+t\sum_{\beta=1}^p r_{0ik}^{\alpha\beta} W_{kj}^\beta +\sum_{\beta,\gamma=1}^p \left(\tilde{f}_{ij}\tilde{f}_{jk}\right)\left( \tilde{r}_{ik}^{\alpha\beta}\right) \tilde{f}_{ij}\left( \tilde{r}_{kj}^{\beta\gamma}\right)\left(tW_{ji}^\gamma +\sum_{\eta=1}^p \tilde{r}_{ji}^{\gamma\eta} \tilde{T}_i^\eta \right)\\
&=t\left(W_{ik}^\alpha+ \sum_{\beta=1}^p r_{0ik}^{\alpha\beta} W_{kj}^\beta+\sum_{\gamma=1}^p r_{0ij}^{\alpha\gamma} W_{ji}^\gamma \right)+ \tilde{T}_i^\alpha+ t\sum_{\gamma=1}^p g_{ijk}^{\alpha\eta} T_{0i}^\eta
\end{align*}
This proves our claim. We show that for another lifting $\tilde{\theta}_{ij}'$ of $\theta_{ij}$, and another lifting $\tilde{T}_i'^\alpha$ of $T_i^\alpha$, and another lifting of $\tilde{r}_{ij}'^{\alpha\beta}$ of $r_{ij}^{\alpha\beta}$, and another lifting $\tilde{g}_{i\alpha\beta}'^\gamma$ of $g_{i\alpha\beta}^\gamma$, the associated $2$-cocycle $\beta:=\left(\left\{ P_i' \right\}, \left\{ - W_{ij}' \right\}, \left\{d_{ijk}' \right\} \right)$ is cohomologus to $\alpha:=\left( \left\{ P_i \right\}, \left\{ - W_{ij}\right\}, \left\{ d_{ijk} \right\} \right)$. Then $\tilde{f}_{ij}'=\tilde{f}_{ij} + t p_{ij}$ for some $p_{ij}\in \Gamma\left( U_{ij}, T_X \right)$ and $\tilde{T}_i'^\alpha= \tilde{T}_i^\alpha + t Y_i^\alpha$ for some $Y_i^\alpha \in \Gamma\left( U_i, T_X \right)$, and $\tilde{r}_{ij}'^{\alpha\beta}= \tilde{r}_{ij}^{\alpha\beta} + t s_{ij}^{\alpha\beta}$ for some $s_{ij}^{\alpha\beta}\in \Gamma\left( U_{ij}, \mathcal{O}_X \right)$, and $\tilde{g}_{i\alpha\beta}'^\gamma=\tilde{g}_{i\alpha\beta}^\gamma + t e_{i\alpha\beta}^\gamma$ for some $e_{i\alpha\beta}^\gamma\in \Gamma\left( U_i, \mathcal{O}_X\right)$. We define $Y_i\in \Gamma\left(U_i, \mathscr{H}om_{\mathcal{O}_X}\left( T_{\mathcal{F}_0}, \frac{T_X}{T_{\mathcal{F}_0}}\right) \right), T_{0i}^\alpha \mapsto \overline{Y_i^\alpha}$.  Then as in $(\ref{at11})$, we have $D_1\left( Y_i \right)=P_i'- P_i$. We note that
\begin{align}\label{ts10}
&t\left(- W_{ij}'^\alpha + W_{ij}^\alpha \right) = - \tilde{\theta}_{ji}'^*\left( \tilde{T}_i'^\alpha \right) + \sum_{\beta=1}^p \tilde{r}_{ij}'^{\alpha\beta} \tilde{T}_j'^\beta + \tilde{\theta}_{ji}^*\left( \tilde{T}_i^\alpha \right) - \sum_{\beta=1}^p \tilde{r}_{ij}^{\alpha\beta} \tilde{T}_j^\beta \\
&=\left(- \tilde{\theta}_{ji}'^* +\tilde{\theta}_{ji}^*\right) \left( \tilde{T}_i^\alpha \right)- t Y_{i}^\alpha + \sum_{\beta=1}^p\left( \tilde{r}_{ij}^{\alpha\beta} +t s_{ij}^{\alpha\beta} \right)\left( \tilde{T}_j^\beta + t Y_j^\beta \right) - \sum_{\beta=1}^p \tilde{r}_{ij}^{\alpha\beta} \tilde{T}_j^\beta \notag \\
&= - t\left[ p_{ji}, T_{0i}^\alpha \right]+ t  \left(  \sum_{\beta=1}^p r_{0ij}^{\alpha\beta} Y_j^\beta - Y_i^\alpha \right) + t \sum_{\beta=1}^p s_{ij}^{\alpha\beta} T_{0j}^\beta = t\left(\left[ p_{ij}, T_{0i}^\alpha \right] + \sum_{\beta=1}^p r_{0ij}^{\alpha\beta} Y_j^\beta - Y_i^\alpha + \sum_{\beta=1}^p s_{ij}^{\alpha\beta} T_{0j}^\beta \right) \notag
\end{align}
This implies that $\delta\left(\left\{ Y_i\right\} \right) + D_0\left( \left\{ p_{ij} \right\}\right)= \left\{ -W_{ij}' + W_{ij} \right\}$ and $\delta\left(\left\{ p_{ij} \right\} \right)= \left\{ d_{ijk}' - d_{ijk}\right\}$ (see \cite{Ser06}). This proves that $\alpha$ is cohomologous to $\beta$. So given a small extension $e:0\to (t)\to \tilde{A} \to A \to 0$ and an infinitesimal foliated deformation $\eta$ of $\left( X, T_{\mathcal{F}_0} \right)$ over $A$, we can associate an element $o_\eta(e):=$ the cohomology class of $\alpha \in \mathbb{H}^2\left( X, T_{\mathcal{F}_0}^\bullet \right)$. We note that $o_\eta(e)=0$ if and only if there exists $\left(\left\{ X_i \right\}, \left\{ d_{ij} \right\} \right) \in C^0\left( \mathcal{U}, \mathscr{H}om_{\mathcal{O}_X}\left( T_{\mathcal{F}_0}, \frac{T_X}{T_{\mathcal{F}_0}} \right) \right) \bigoplus C^1\left( \mathcal{U}, T_X \right)$ where $X_i:\Gamma\left( U_i, T_{\mathcal{F}_0} \right) \to \Gamma\left( U_i, \frac{T_X}{T_{\mathcal{F}_0}}  \right), T_{0i}^\alpha \mapsto \overline{X_i^\alpha}$ such that $D_1(X_i)= P_i$ and $X_j- X_i + D_0\left( d_{ij} \right)= - W_{ij}$ and $d_{ij}- d_{ik} + d_{jk}= d_{ijk}$, so that
\begin{align}\label{atb2}
t\left(d_{ij}- d_{ik} + d_{jk} \right) = \tilde{f}_{ij}\tilde{f}_{jk}- \tilde{f}_{ik}  \Longrightarrow \left(\tilde{f}_{ij}- t d_{ij} \right)\left(\tilde{f}_{jk}- t d_{jk} \right)= \tilde{f}_{ik} - td_{ik}
\end{align}
Let $\Phi_{ij}: U_{ij}\times \textnormal{Spec}\left(\tilde{A}\right) \to U_{ij} \times \textnormal{Spec}\left( \tilde{A}\right)$ be an isomorphism corresponding to $\tilde{f}_{ij} - t d_{ij}$. Then $\Phi_{jk}\Phi_{ij} = \Phi_{ik}$. We also note that
{\small{\begin{align}
&t\left(\left[ d_{ij}, T_{0i}^\alpha \right] +\sum_{\beta=1}^p r_{0ij}^{\alpha\beta} X_j^\beta - X_i^\alpha +  \sum_{\beta=1}^p v_{ij}^{\alpha\beta} T_{0j}^\beta\right)= - \tilde{\theta}_{ji}^* \left(\tilde{T}_i^\alpha \right) + \sum_{\beta=1}^p \tilde{r}_{ij}^{\alpha\beta} \tilde{T}_j^\beta \notag \\
& \Longrightarrow \tilde{\theta}_{ji}^*\left( \tilde{T}_i^\alpha \right) - tX_i^\alpha + t \left[ d_{ij}, T_{0i}^\alpha \right]  = \sum_{\beta=1}^p \left(\tilde{r}_{ij}^{\alpha\beta}- tv_{ij}^{\alpha\beta} \right)\left(\tilde{T}_j^\beta - t X_j^\beta \right)\,\,\,\textnormal{for some} \,\,\, v_{ij}^{\alpha\beta}\in \Gamma\left( U_{ij}, \mathcal{O}_X \right) \notag\\
&\Longrightarrow \Phi_{ji}^*\left( \tilde{T}_i^\alpha- tX_i^\alpha \right) = \sum_{\beta=1}^p \left(\tilde{r}_{ij}^{\alpha\beta}- tv_{ij}^{\alpha\beta} \right)\left(\tilde{T}_j^\beta - t X_j^\beta \right)
\end{align}}}
and as in $(\ref{at12})$, we have  $\left[ \tilde{T}_i^\alpha- t X_i^\alpha, \tilde{T}_i^\beta - t X_i^\beta \right] - \sum_{\gamma=1}^p\left( \tilde{g}_{i\alpha\beta}^\gamma + t l_{i\alpha\beta}^\gamma  \right)\left(\tilde{T}_i^\gamma - t X_i^\gamma \right)$ for some $l_{i\alpha\beta}^\gamma \in \Gamma\left(U_i, \mathcal{O}_X \right)$. Hence a lifting $\tilde{\eta}$ of $\eta$ exists. This completes the proof of Proposition \ref{at25}.

\end{proof}

\subsection{Relationship between foliated deformations of $\left(X, T_{\mathcal{F}_0} \right)$ with $X$ fixed and flat foliated deformation of $\left(X, T_{\mathcal{F}_0} \right)$ and underlying flat deformations of $X$}\

Let $\left( X, T_{\mathcal{F}_0} \right)$ be a foliated nonsingular variety with $T_{\mathcal{F}_0}$ locally free. Let $\textnormal{Def}_X$ be the functor of Artin rings of flat deformations of the underlying nonsingular variety $X$ (see \cite{Ser06}). As is known, flat deformations of the underlying nonsingular variety $X$ is controlled by the tangent sheaf $T_X$. We have a canonical morphism of deformation functors
\begin{align*}
\textnormal{FDef}_{T_{\mathcal{F}_0}}\to \textnormal{FDef}_{\left(X, T_{\mathcal{F}_0}\right)}
\end{align*}
and we have a morphism of deformation functors
\begin{align}\label{e1}
\textnormal{FDef}_{\left(X, T_{\mathcal{F}_0}\right)}\to \textnormal{Def}_X
\end{align}
by forgetting foliated structures. Then the relationship of three deformation functors $\textnormal{FDef}_{T_{\mathcal{F}_0}}, \textnormal{FDef}_{\left( X, T_{\mathcal{F}_0} \right)}$ and $\textnormal{Def}_X$ is expressed by the following exact sequence of complexes of sheaves
{\Small{\begin{equation}\label{e2}
\begin{CD}
@. \cdots @. \cdots @. \cdots \\
@. @AD_3AA @AD_3AA @AAA \\
0@>>> \mathscr{H}om_{\mathcal{O}_X}\left( \bigwedge^3 T_{\mathcal{F}_0}, \frac{T_X}{T_{\mathcal{F}_0}}\right) @>>> \mathscr{H}om_{\mathcal{O}_X}\left( \bigwedge^3 T_{\mathcal{F}_0}, \frac{T_X}{T_{\mathcal{F}_0}} \right) @>>> 0@>>> 0 \\
@. @AD_2AA @AD_2AA @AAA \\
0@>>> \mathscr{H}om_{\mathcal{O}_X}\left( \bigwedge^2 T_{\mathcal{F}_0}, \frac{T_X}{T_{\mathcal{F}_0}}\right) @>>> \mathscr{H}om_{\mathcal{O}_X}\left( \bigwedge^2 T_{\mathcal{F}_0}, \frac{T_X}{T_{\mathcal{F}_0}} \right) @>>> 0@>>> 0 \\
@. @AD_1AA @AD_1AA @AAA \\
0@>>> \mathscr{H}om_{\mathcal{O}_X}\left( T_{\mathcal{F}_0}, \frac{T_X}{T_{\mathcal{F}_0}}\right) @>>> \mathscr{H}om_{\mathcal{O}_X}\left( T_{\mathcal{F}_0}, \frac{T_X}{T_{\mathcal{F}_0}} \right) @>>> 0@>>> 0 \\
@. @AAA @AD_0AA @AAA \\
0 @>>> 0 @>>> T_X @>>> T_X @>>> 0
\end{CD}
\end{equation}}}

Then the long exact sequence induces maps of first-order deformations
{\small{\begin{align*}
\textnormal{FDef}_{T_{\mathcal{F}_0}}\left(k[\epsilon]\right)\cong \mathbb{H}^0\left(X,\mathscr{H}om_{\mathcal{O}_X}\left(T_{\mathcal{F}_0}, \frac{T_X}{T_{\mathcal{F}_0}} \right)^\bullet  \right) \to \textnormal{FDef}_{ \left(X, T_{\mathcal{F}_0} \right)}\left(k[\epsilon]\right)\cong  \mathbb{H}^1\left(X, T_{\mathcal{F}_0}^\bullet \right)\to \textnormal{Def}_X\left(k[\epsilon]\right)\cong H^1\left(X,T_X\right)
\end{align*}}}
and obstruction maps
\begin{align*}
\mathbb{H}^1\left( X, \mathscr{H}om_{\mathcal{O}_X}\left(T_{\mathcal{F}_0}, T_X/T_{\mathcal{F}_0} \right)^\bullet \right) \to \mathbb{H}^2(X, T_{\mathcal{F}_0}^\bullet) \to H^2(X, T_X)
\end{align*}
of morphisms of deformation functors $\textnormal{FDef}_{T_{\mathcal{F}_0}}\to \textnormal{FDef}_{\left(X, T_{\mathcal{F}_0} \right)} \to \textnormal{Def}_X$.

\subsection{Flat foliated deformations of $\left(X, T_{\mathcal{F}_0} \right)$ with $T_{\mathcal{F}_0}$ locally free in terms of vector bundles}\

The underlying deformation of locally free sheaf $T_{\mathcal{F}_0}$ is implicit in the leaf complex $(\ref{at16})$. In order to make the underlying deformations of $T_{\mathcal{F}_0}$ explicit, we define a complex of sheaves also controlling flat foliated deformations of $\left(X, T_{\mathcal{F}_0} \right)$ and present foliated deformations of $\left(X, T_{\mathcal{F}_0}\right)$ by using the complex. As is known, deformations of a pair of a nonsingular variety $X$ and a locally free sheaf $T_{\mathcal{F}_0}$ is controlled by the Atiyah extension (see \cite{Hua95}) which is the exact sequence
\begin{align}\label{si1}
0\to \mathscr{H}om_{\mathcal{O}_X}\left( T_{\mathcal{F}_0} ,  T_{\mathcal{F}_0} \right)  \to  \mathcal{E}_{T_{\mathcal{F}_0}} \to T_X \to 0
\end{align}
where $\mathcal{E}_{T_{\mathcal{F}_0}}$ is defined in the following way: let $\left\{U_i \right\}$ be an affine open covering of $X$ such that $\Gamma\left(U_i, T_{\mathcal{F}_0}\right)\cong \bigoplus^p \Gamma\left(U_i, \mathcal{O}_X\right)$ is trivialized in terms of the basis $\left\{T_{0i}^\alpha \right\}_{\alpha=1,...,p}$ with  the relation $T_{0i}^\alpha=\sum_{\beta=1}^p r_{0ij}^{\alpha\beta} T_{0j}^\beta,\alpha=1,...,p$. Then we define $\mathcal{E}_{T_{\mathcal{F}_0}}$ locally on $U_i$
\begin{align*}
\mathcal{E}_{T_{\mathcal{F}_0}} |_{U_i}\cong \mathscr{H}om_{\mathcal{O}_X}\left(T_{\mathcal{F}_0}, T_{\mathcal{F}_0} \right)  |_{U_i} \bigoplus T_X |_{U_i}
\end{align*}
such that $(\phi_i,X)\in \mathscr{H}om_{\mathcal{O}_X}\left( T_{\mathcal{F}_0}, T_{\mathcal{F}_0} \right)|_{U_i} \bigoplus T_X |_{U_i}$ and $(\phi_j, X)\in \mathscr{H}om_{\mathcal{O}_X}\left( T_{\mathcal{F}_0}, T_{\mathcal{F}_0} \right) \bigoplus T_X$ are identified on $U_{ij}$ if
\begin{align*}
\phi_i\left(T_{0i}^\alpha\right)=\phi_j\left(T_{0i}^\alpha\right)+\sum_{\beta=1}^p \left[X, r_{0ij}^{\alpha\beta}\right]T_{0j}^\beta
\end{align*}
Then $\mathcal{E}_{T_{\mathcal{F}_0}}$ is well-defined. In other words, if $\left(\phi_i,X \right)= \left(\phi_j,X \right)$ and $\left(\phi_j,X\right)=\left(\phi_k,X\right)$, then $\left(\phi_i,X\right)=\left(\phi_k,X\right)$. In fact, since $\phi_j\left(T_{0j}^\beta\right)=\phi_k\left(T_{0j}^\beta\right)+\sum_{\gamma=1}^p \left[X, r_{0jk}^{\beta\gamma}\right]T_{0k}^\gamma$, we have
\begin{align}
\phi_i \left(T_{0i}^\alpha \right)&=\phi_j \left(T_{0i}^\alpha \right)+\sum_{\beta=1}^p \left[X, r_{0ij}^{\alpha\beta}\right]T_{0j}^\beta=\phi_j\left(\sum_{\beta=1}^p r_{0ij}^{\alpha\beta}T_{0j}^\beta \right)  +\sum_{\beta,\gamma=1}^p \left[X, r_{0ij}^{\alpha\beta}\right]r_{0jk}^{\beta\gamma} T_{0k}^\gamma \label{si3} \\
 &=\sum_{\beta=1}^p r_{0ij}^{\alpha\beta}\left( \phi_k\left(T_{0j}^\beta\right)+\sum_{\gamma=1}^p \left[X, r_{0jk}^{\beta\gamma}\right]T_{0k}^\gamma \right) +\sum_{\beta,\gamma=1}^p \left[X, r_{0ij}^{\alpha\beta}\right]r_{0jk}^{\beta\gamma}T_{0k}^\gamma=\phi_k\left(T_{0i}^\alpha\right) +\sum_{\gamma=1}^p \left[X, r_{0ik}^{\alpha\gamma}\right]T_{0k}^\gamma \notag
\end{align}

Then we shall construct a complex of sheaves
{\small{\begin{align}\label{tv1}
\mathcal{E}_{T_{\mathcal{F}_0}}^\bullet: \mathcal{E}_{T_{\mathcal{F}_0}}\xrightarrow{D_0'} \mathscr{H}om_{\mathcal{O}_X}\left( T_{\mathcal{F}_0}, T_X \right) \xrightarrow{D_1'} \mathscr{H}om_{\mathcal{O}_X}\left( \bigwedge^2 T_{\mathcal{F}_0}, \frac{T_X}{T_{\mathcal{F}_0}}\right) \xrightarrow{D_2} \mathscr{H}om_{\mathcal{O}_X}\left( \bigwedge^3 T_{\mathcal{F}_0}, \frac{T_X}{T_{\mathcal{F}_0}} \right) \xrightarrow{D_3} \cdots
\end{align}}}
which also controls foliated deformations of $\left(X, T_{\mathcal{F}_0} \right)$. We will denote the $i$-th cohomology group of $\mathcal{E}_{T_{\mathcal{F}_0}}^\bullet$ by $\mathbb{H}^i\left( X,  \mathcal{E}_{T_{\mathcal{F}_0}}^\bullet\right)$. We define $D_0':\mathcal{E}_{T_{\mathcal{F}_0}}\to \mathscr{H}om_{\mathcal{O}_X}\left(T_{\mathcal{F}_0}, T_X \right)$ locally in the following way: 
\begin{align*}
D_0':\mathcal{E}_{T_{\mathcal{F}_0}}|_{U_i}\cong \mathscr{H}om_{\mathcal{O}_X}\left( T_{\mathcal{F}_0}, T_{\mathcal{F}_0} \right)\bigoplus T_X|_{U_i} &\to \mathscr{H}om_{\mathcal{O}_X}\left(T_{\mathcal{F}_0}, T_X \right)
\end{align*}
On $U_i$ we define $D_0'((\phi_i,X))$ in terms of basis $\{T_{0i}^\alpha\}$
\begin{align}\label{si7}
D_0' \left(\left(\phi_i,X\right)\right)\left(T_{0i}^\alpha\right)= \left(\mathcal{L}_X-\phi_i\right)\left(T_{0i}^\alpha\right)
\end{align}
Then we linearly extend to $T_{\mathcal{F}_0}|_{U_i}$. We remark that  $D_0'\left(\left(\phi_i,X\right)\right)\left(T\right)$ may be different from $\left(L_X-\phi\right)\left(T\right)$ for arbitrary $T \in T_{\mathcal{F}_0}|_{U_i}$. Similarly we define $D_0'\left(\left(\phi_j,X\right)\right)\left(T_{0j}^\beta\right)=\left(\mathcal{L}_X-\phi_j\right)\left(T_{0j}^\beta\right)$ on $U_j$ in terms of basis $\left\{T_{0j}^\beta\right\}$ and then linearly extend to $T_{\mathcal{F}_0}|_{U_j}$. We claim that $D_0'$ is well-defined. In fact, if $\left(\phi_i,X\right)=\left(\phi_j,X\right)$ on $U_{ij}$, then we have 
\begin{align}
&\left(\mathcal{L}_X-\phi_i\right)\left(T_{0i}^\alpha\right)-\sum_{\beta=1}^p r_{0ij}^{\alpha\beta}\left(\mathcal{L}_X-\phi_j\right)\left(T_{0j}^\beta\right)= \mathcal{L}_X\left(T_{0i}^\alpha\right)-\phi_i\left(T_{0i}^\alpha\right) -\sum_{\beta=1}^p r_{0ij}^{\alpha\beta}\mathcal{L}_X\left(T_{0j}^\beta\right)+\sum_{\beta=1}^p r_{0ij}^{\alpha\beta} \phi_j\left(T_{0j}^\beta\right) \label{si9} \\
& = \mathcal{L}_X\left( \sum_{\beta=1}^p r_{0ij}^{\alpha\beta} T_{0j}^\beta \right) -\phi_i \left(T_{0i}^\alpha\right)-\sum_{\beta=1}^p r_{0ij}^{\alpha\beta} \mathcal{L}_X\left(T_{0j}^\beta\right)+\phi_j\left(T_{0i}^\alpha\right)=\sum_{\beta=1}^p \left[X, r_{0ij}^{\alpha\beta}\right] T_{0j}^\beta -\phi_i \left(T_{0i}^\alpha\right)+\phi_j \left(T_{0i}^\alpha\right) = 0. \notag
\end{align}
Next we define $D_1':\mathscr{H}om_{\mathcal{O}_X}\left( T_{\mathcal{F}_0}, T_X \right) \to \mathscr{H}om_{\mathcal{O}_X}\left( \bigwedge^2 T_{\mathcal{F}_0}, \frac{T_X}{T_{\mathcal{F}_0}} \right)$ in $(\ref{tv1})$ to be the composition of the natural quotient map $T_X\to T_X/T_{\mathcal{F}_0}$ and $D_1$. Then we have an exact sequence of complexes of sheaves:
{\Small{\begin{center}
$\begin{CD}
@. \cdots @. \cdots @. \cdots \\
@. @AAA @AD_3AA @AD_3AA \\
0 @>>>0@>>>  \mathscr{H}om_{\mathcal{O}_X}\left( \bigwedge^3 T_{\mathcal{F}_0}, \frac{T_X}{T_{\mathcal{F}_0}} \right) @>>> \mathscr{H}om_{\mathcal{O}_X}\left( \bigwedge^3 T_{\mathcal{F}_0}, \frac{T_X}{T_{\mathcal{F}_0}} \right) @>>> 0 \\
@.@AAA @AD_2AA @AD_2AA \\
0@>>> 0@>>> \mathscr{H}om_{\mathcal{O}_X}\left(\bigwedge^2 T_{\mathcal{F}_0}, \frac{T_X}{T_{\mathcal{F}_0}} \right) @>>> \mathscr{H}om_{\mathcal{O}_X}\left(\bigwedge^2 T_{\mathcal{F}_0}, \frac{T_X}{T_{\mathcal{F}_0}} \right) @>>> 0 \\
@.@AAA @AD_1' AA @AD_1AA \\
0@>>> \mathscr{H}om_{\mathcal{O}_X}\left( T_{\mathcal{F}_0}, T_{\mathcal{F}_0} \right) @>>> \mathscr{H}om_{\mathcal{O}_X}\left(T_{\mathcal{F}_0}, T_X \right) @>>> \mathscr{H}om_{\mathcal{O}_X}\left(T_{\mathcal{F}_0}, \frac{T_X}{T_{\mathcal{F}_0}} \right) @>>> 0\\
@. @A-\textnormal{Id}AA @A D_0' AA @AD_0 AA \\
0@>>> \mathscr{H}om_{\mathcal{O}_X}\left( T_{\mathcal{F}_0}, T_{\mathcal{F}_0} \right) @>>> \mathcal{E}_{T_{\mathcal{F}_0}}@>>> T_X @>>> 0
\end{CD}$
\end{center}}}
This implies that
\begin{align*}
\mathbb{H}^i \left( X, \mathcal{E}_{T_{\mathcal{F}_0}}^\bullet \right) \cong \mathbb{H}^i \left(X, T_{\mathcal{F}_0}^\bullet\right)
\end{align*}

We reinterpret Proposition \ref{at25} in terms of $\mathcal{E}_{T_{\mathcal{F}_0}}^\bullet$ instead of $T_{\mathcal{F}_0}^\bullet$.

\begin{proposition}[flat foliated deformations in terms of vector bundles] \label{at27}
Let $\left( X, T_{\mathcal{F}_0} \right)$ be a foliated nonsingular variety with $T_{\mathcal{F}_0}$ locally free.
\begin{enumerate}
\item there is a $1-1$ correspondence
\begin{align*}
\textnormal{FDef}_{\left(X, T_{\mathcal{F}_0}\right)}\left(k[\epsilon] \right)\cong \mathbb{H}^1\left(X, \mathcal{E}_{T_{\mathcal{F}_0}}^\bullet \right)
\end{align*}
\item Let $A\in \bold{Art}$ and $\eta=\left( \mathcal{X}, T_\mathcal{F} \right)$ be an infinitesimal foliated deformations of $\left( X, T_{\mathcal{F}_0} \right)$ over $A$. To every small extension $0\to (t)\to \tilde{A} \to A \to 0$, we can associate an element $o_\eta(e)\in \mathbb{H}^2\left( X,  \mathcal{E}_{T_{\mathcal{F}_0}}^\bullet \right)$ called the obstruction lifting of $\eta$ to $\tilde{A}$, which is $0$ if and only if a lifting of $\eta$ to $\tilde{A}$ exists.
\end{enumerate}
\end{proposition}

\begin{proof}
Let $\left( \mathcal{X}, T_\mathcal{F} \right)$ be a first-order foliated deformation of $\left( X, T_{\mathcal{F}_0} \right)$. We keep the notations in the proof $(1)$ of Proposition \ref{at25}. We define an element $T_i'\in \Gamma\left( U_i, \mathscr{H}om_{\mathcal{O}_X}\left( T_{\mathcal{F}_0}, T_X \right)\right)$ by $T_i':\Gamma\left( U_i, T_{\mathcal{F}_0} \right)\to \Gamma\left( U_i, T_X\right), T_{0i}^\alpha \mapsto T_i'^\alpha$. We define an element $b_{ij}\in \Gamma\left( U_{ij}, \mathscr{H}om_{\mathcal{O}_X}\left( T_{\mathcal{F}_0}, T_{\mathcal{F}_0} \right) \right)$ by
\begin{align*}
b_{ij}:\Gamma\left( U_{ij}, T_{\mathcal{F}_0} \right) &\to \Gamma\left(U_{ij}, T_{\mathcal{F}_0} \right)\\
 T_{0i}^\alpha &\mapsto \sum_{\beta=1}^p b_{ij}^{\alpha\beta} T_{0j}^\beta
\end{align*}
We claim that $\left( p_{ij}, - b _{ij} \right)\in T_X|_{U_{ij}} \bigoplus \mathscr{H}om_{\mathcal{O}_X}\left( T_{\mathcal{F}_0}, T_{\mathcal{F}_0} \right)|_{U_{ij}} \subset T_X|_{U_i} \bigoplus \mathscr{H}om_{\mathcal{O}_X}\left( T_{\mathcal{F}_0}, T_{\mathcal{F}_0} \right)|_{U_i} $ is identified with $\left( -p_{ji}, b_{ji} \right)\in T_X|_{U_{ij}}\bigoplus \mathscr{H}om_{\mathcal{O}_X}\left( T_{\mathcal{F}_0}, T_{\mathcal{F}_0} \right)|_{U_{ij}}\subset T_X|_{U_j} \bigoplus \mathscr{H}om_{\mathcal{O}_X}\left( T_{\mathcal{F}_0}, T_{\mathcal{F}_0} \right)|_{U_j}$ as an element of $\Gamma\left(U_{ij}, \mathcal{E}_{T_{\mathcal{F}_0}} \right)$. In fact,
\begin{align}\label{tn1}
&T_{0i}^\alpha + \epsilon T_i'^\alpha=\theta_{ij}^* \theta_{ji}^*\left( T_{0i}^\alpha+ \epsilon T_i'^\alpha\right) = \theta_{ij}^*\left(\sum_{\beta=1}^p \left( r_{0ij}^{\alpha\beta} + \epsilon b_{ij}^{\alpha\beta} \right) \left( T_{0j}^\beta + \epsilon T_j'^\beta \right) \right) \\
&= \sum_{\beta , \gamma =1}^p\left( r_{0ij}^{\alpha\beta} + \epsilon b_{ij}^{\alpha\beta} +\epsilon \left[ p_{ij}, r_{0ij}^{\alpha\beta}\right]     \right)\left( r_{0ji}^{\beta \gamma}+\epsilon b_{ji}^{\beta \gamma} \right)\left( T_{0i}^\gamma + \epsilon T_i'^\gamma \right) \notag \\
&\Longrightarrow \sum_{\beta=1}^p \left( r_{0ij}^{\alpha\beta} + \epsilon b_{ij}^{\alpha\beta} +\epsilon \left[ p_{ij}, r_{0ij}^{\alpha\beta}\right]     \right)\left( r_{0ji}^{\beta \gamma}+\epsilon b_{ji}^{\beta \gamma} \right) = \delta_\alpha^\gamma \Longrightarrow \sum_{\beta=1}^p b_{ij}^{\alpha\beta} r_{0ji}^{\beta \gamma} + \sum_{\beta=1}^p \left[ p_{ij}, r_{0ij}^{\alpha\beta} \right] r_{0ji}^{\beta \gamma} + \sum_{\beta=1}^p r_{0ij}^{\alpha\beta} b_{ji}^{\beta \gamma}=0 \notag \\
&\Longrightarrow \sum_{\beta=1}^p b_{ij}^{\alpha\beta} T_{0j}^\beta + \sum_{\beta=1}^p \left[ p_{ij}, r_{0ij}^{\alpha\beta} \right] T_{0j}^\beta + \sum_{\beta, \gamma=1}^p r_{0ij}^{\alpha\beta} b_{ji}^{\beta \gamma} T_{0i}^\gamma =0 \Longrightarrow - b_{ij}\left( T_{0i}^\alpha \right) = \sum_{\beta=1}^p \left[ p_{ij}, r_{0ij}^{\alpha\beta} \right] T_{0j}^\beta + \sum_{\beta=1}^p r_{0ij}^{\alpha\beta} b_{ji}\left( T_{0j}^\beta \right) \notag \\
&\Longrightarrow - b_{ij}\left(  T_{0i}^\alpha \right)= b_{ji}\left( T_{0i}^\alpha \right) + \sum_{\beta=1}^p \left[p_{ij}, r_{0ij}^{\alpha\beta} \right] T_{0j}^\beta \notag
\end{align}
This implies that $\left\{ \left( p_{ij}, - b_{ij} \right)\right\} \in C^1\left( \mathcal{U}, \mathcal{E}_{T_{\mathcal{F}_0}} \right)$. Then we claim that $\left( \left\{ T_i' \right\}, \left\{\left( p_{ij}, - b_{ij} \right) \right\} \right) \in C^0\left( \mathcal{U}, \mathscr{H}om_{\mathcal{O}_X}\left( T_{\mathcal{F}_0}, T_X \right) \right)$ $ \bigoplus C^1\left( \mathcal{U}, \mathcal{E}_{T_{\mathcal{F}_0}} \right)$ defines a $1$-cocycle in the following \v Cech resolution of $\mathcal{E}_{T_{\mathcal{F}_0}}^\bullet$:
{\Small{\begin{center}
$\begin{CD}
\cdots\\
@AD_3AA \\
C^0\left(\mathcal{U}, \mathscr{H}om_{\mathcal{O}_X}\left( \bigwedge^3 T_{\mathcal{F}_0}, \frac{T_X}{T_{\mathcal{F}_0}}\right) \right)@>\delta>> \cdots\\
@AD_2AA @AD_2AA \\
C^0\left(\mathcal{U}, \mathscr{H}om_{\mathcal{O}_X}\left(\bigwedge^2 T_{\mathcal{F}_0}, \frac{T_X}{T_{\mathcal{F}_0}} \right) \right) @>-\delta>> C^1\left(\mathcal{U}, \mathscr{H}om_{\mathcal{O}_X} \left( \bigwedge^2 T_{\mathcal{F}_0}, \frac{T_X}{T_{\mathcal{F}_0}} \right) \right) @>\delta>> \cdots\\ 
@AD_1' AA @AD_1' AA @AD_1'AA \\
C^0\left(\mathcal{U}, \mathscr{H}om_{\mathcal{O}_X}\left( T_{\mathcal{F}_0}, T_X \right) \right)@>\delta>> C^1\left(\mathcal{U}, \mathscr{H}om_{\mathcal{O}_X}\left( T_{\mathcal{F}_0}, T_X \right) \right) @>-\delta>> C^2\left(\mathcal{U}, \mathscr{H}om_{\mathcal{O}_X}\left( T_{\mathcal{F}_0}, T_X \right) \right) \\
@AD_0' AA @AD_0' AA @AD_0' AA \\
C^0\left(\mathcal{U}, \mathcal{E}_{T_{\mathcal{F}_0}}\right)@>-\delta>> C^1\left(\mathcal{U}, \mathcal{E}_{T_{\mathcal{F}_0}}\right) @>\delta>> C^2 \left(\mathcal{U}, \mathcal{E}_{T_{\mathcal{F}_0}}\right) @>-\delta>> \cdots
\end{CD}$
\end{center}}}
It is sufficient to show that $\delta\left( \left\{ \left( p_{ij}, - b_{ij} \right) \right\} \right)=0$. In fact, we note that
{\small{\begin{align*}
\theta_{kj}^*\theta_{ji}^*\left( T_{0i}^\alpha + \epsilon T_i'^\alpha \right)&= \theta_{kj}^*\left(\sum_{\beta=1}^p \left( r_{0ij}^{\alpha\beta} + \epsilon b_{ij}^{\alpha\beta} \right)\left( T_{0j}^\beta + \epsilon T_j'^\beta \right) \right) = \sum_{\beta=1}^p \left(\textnormal{Id}+ \epsilon p_{kj} \right)\left( r_{0ij}^{\alpha\beta} + \epsilon  b_{ij}^{\alpha\beta} \right)\theta_{kj}^*\left( T_{0j}^\beta + \epsilon T_j'^\beta \right)\\
&= \sum_{\beta, \gamma=1}^p \left( r_{0ij}^{\alpha\beta} + \epsilon b_{ij}^{\alpha\beta} + \epsilon \left[ p_{kj}, r_{0ij}^{\alpha\beta}\right] \right) \left( r_{0jk}^{\beta \gamma} + \epsilon b_{jk}^{\beta \gamma}   \right) \left( T_{0k}^\gamma + \epsilon T_k'^\gamma   \right)\\
&=\sum_{\beta, \gamma=1}^p \left( r_{0ij}^{\alpha\beta} r_{0jk}^{\beta \gamma} + \epsilon \left( r_{0ij}^{\alpha\beta} b_{jk}^{\beta \gamma} + b_{ij}^{\alpha\beta} r_{0jk}^{\beta \gamma} + \left[p_{kj}, r_{0ij}^{\alpha\beta} \right] r_{0jk}^{\beta \gamma}    \right) \right)  \left( T_{0k}^\gamma + \epsilon T_k'^\gamma  \right)
\end{align*}}}
On the other hand, we have
\begin{align*}
\theta_{ki}^*\left( T_{0i}^\alpha + \epsilon T_i'^\alpha   \right) = \sum_{\gamma=1}^p \left( r_{0ik}^{\alpha\gamma} + \epsilon b_{ik}^{\alpha \gamma} \right)\left( T_{0k}^\gamma + \epsilon T_k'^\gamma  \right)
\end{align*}
Then since $\theta_{kj}^*\theta_{ij}^*\left( T_{0i}^\alpha + \epsilon T_i'^\alpha \right)=\theta_{ki}^*\left( T_{0i}^\alpha + \epsilon T_i'^\alpha \right)$, by considering the coefficients of $\epsilon$ we see that
\begin{align}
& \sum_{\beta=1}^p b_{ij}^{\alpha\beta} r_{0jk}^{\beta \gamma} - b_{ik}^{\alpha \gamma} + \sum_{\beta=1}^p r_{0ij}^{\alpha\beta} b_{jk}^{\beta \gamma} - \sum_{\beta=1}^p \left[ p_{jk}, r_{0ij}^{\alpha\beta}\right] r_{0jk}^{\beta \gamma} = 0 \notag \\
& \Longrightarrow \sum_{\beta=1}^p b_{ij}^{\alpha\beta} T_{0j}^\beta - \sum_{\gamma=1}^p b_{ik}^{\alpha \gamma} T_{0k}^\gamma + \sum_{\beta, \gamma=1}^p r_{0ij}^{\alpha\beta} b_{jk}^{\beta \gamma} T_{0k}^\gamma - \sum_{\beta=1}^p \left[ p_{jk}, r_{0ij}^{\alpha\beta} \right] T_{0j}^\beta = 0 \label{at28}
\end{align}
We note that $\left(p_{jk}, - b_{jk}\right)$ on $U_j$ is identified with $\left(p_{jk}, b_{jk}^i \right)$ on $U_i$, where $b_{jk}^i\left(T_{0i}^\alpha\right)= - b_{jk}\left( T_{0i}^\alpha\right) + \sum_{\beta=1}^p \left[ p_{jk}, r_{0ij}^{\alpha\beta} \right] T_{0j}^\beta$.
Then we see that from $(\ref{at28})$
\begin{align}\label{tn2}
\left(- b_{ij}+ b_{ik} + b_{jk}^i \right)\left( T_{0i}^\alpha\right) = - \sum_{\beta=1}^p b_{ij}^{\alpha\beta} T_{0j}^\beta + \sum_{\gamma=1}^p b_{ik}^{\alpha \gamma} T_{0k}^\gamma - b_{jk}\left(\sum_{\beta=1}^p r_{0ij}^{\alpha\beta} T_{0j}^\beta \right) + \sum_{\beta=1}^p \left[ p_{jk},  r_{0ij}^{\alpha\beta} \right] T_{0j}^\beta = 0.
\end{align}
This implies that $\delta\left(\left\{ p_{ij}, - b_{ij}  \right\} \right)=0$. We show that for two equivalent first-order foliated deformations of $\left( X, T_{\mathcal{F}_0}\right)$, the associated $1$-cocycles are equivalent. Let $\left( \mathcal{X}', T_{\mathcal{F}'} \right)$ be another foliated deformations of $\left( X, T_{\mathcal{F}_0} \right)$ over $k[\epsilon]$ which induces $\left( \left\{ T_i'' \right\}, \left\{ \left( p_{ij}' , - b_{ij}'\right)\right\}  \right)$ and $\Phi:\left( \mathcal{X}, T_\mathcal{F} \right) \to \left( \mathcal{X}', T_{\mathcal{F}'} \right)$ is a foliated isomorphism inducing the identity on $\left(X, T_{\mathcal{F}_0} \right)$. Then we define $b_i\in \Gamma\left( U_i, \mathscr{H}om_{\mathcal{O}_X}\left( T_{\mathcal{F}_0}, T_{\mathcal{F}_0} \right) \right)$ by $b_i:\Gamma\left( U_i, T_{\mathcal{F}_0}\right) \to \Gamma\left(U_i, T_{\mathcal{F}_0} \right), T_{0i}^\alpha \mapsto \sum_{\beta=1}^p  b_i^{\alpha\beta} T_{0i}^\beta$. We consider $\left\{\left\{ \left( a_i, -b_i \right)\right\} \right\}\in C^0\left( \mathcal{U}, \mathcal{E}_{T_{\mathcal{F}_0}}\right)$. It is sufficient to show that $- b_i+ b_j^i= - b_{ij} + b_{ij}'$, i.e. $-b_i\left( T_{0i}^\alpha \right) + b_j\left( T_{0i}^\alpha \right) - \sum_{\beta=1}^p\left[ a_j, r_{0ij}^{\alpha\beta} \right] T_{0j}^\beta = - b_{ij}\left( T_{0i}^\alpha\right) + b_{ij}'\left( T_{0i}^\alpha \right)$. In fact, since $ \theta_{ji} \alpha_j^{-1}= \alpha_i^{-1} \theta_{ji}'$, we have
\begin{align}
 \left(\alpha_j^{-1}\right)^* \theta_{ji}^*\left( T_{0i}^\alpha + \epsilon T_i'^\alpha \right)&= \left( \alpha_j^{-1} \right)^* \left(\sum_{\beta=1}^p \left( r_{0ij}^{\alpha\beta} + \epsilon b_{ij}^{\alpha\beta} \right)\left( T_{0j}^\beta + \epsilon T_j'^\beta \right)    \right) \label{tn3} \\
 &  = \sum_{\beta, \gamma=1}^p \left( r_{0ij}^{\alpha\beta} + \epsilon b_{ij}^{\alpha\beta} - \epsilon \left[ a_j, r_{0ij}^{\alpha\beta} \right]     \right)\left( \delta_\beta^\gamma + \epsilon b_j^{\beta \gamma} \right)\left( T_{0j}^\gamma + \epsilon T_{j}''^\gamma \right) \notag \\
 \theta_{ji}'^* \left(\alpha_i^{-1} \right)^* \left( T_{0i}^\alpha+ \epsilon T_i'^\alpha \right) &= \theta_{ji}'^* \left( \sum_{\beta=1}^p  \left( \delta_\alpha^\beta + \epsilon b_i^{\alpha\beta} \right)\left( T_{0i}^\beta + \epsilon T_i''^\beta \right) \right) \label{tn4}\\
 &= \sum_{\beta=1}^p \left( \delta_\alpha^\beta + \epsilon b_i^{\alpha\beta} \right)\left(  r_{0ij}^{\beta\gamma} + \epsilon b_{ij}'^{\beta\gamma} \right)\left( T_{0j}^\gamma + \epsilon T_j''^\gamma \right) \notag
\end{align}
Then we have $\sum_{\beta=1}^p b_{ij}^{\alpha\beta} T_{0j}^\beta + \sum_{\beta, \gamma=1}^p r_{0ij}^{\alpha\beta} b_j^{\beta \gamma} T_{0j}^\gamma - \sum_{\beta=1}^p \left[ a_j, r_{0ij}^{\alpha\beta} \right] T_{0j}^\beta= \sum_{\beta=1}^p  b_i^{\alpha\beta} T_{0j}^\beta + \sum_{\beta=1}^p b_{ij}'^{\alpha \beta} T_{0j}^\beta$.

Next we identity obstructions. Let $e:0\to (t) \to \tilde{A} \to A \to 0$ be a small extension in $\bold{Art}$ and $\eta:=\left( \mathcal{X}, T_\mathcal{F} \right)$ be a foliated deformation of $\left(X , T_{\mathcal{F}_0} \right)$ over $A$. We keep the notations in the proof of $(2)$ in Proposition \ref{at25}. We recall that $\sum_{\beta=1}^p \tilde{f}_{ij}\left(\tilde{r}_{ij}^{\alpha\beta} \right)\tilde{r}_{ji}^{\beta \gamma}= \delta_\alpha^\gamma + t b_{ij\alpha}^\gamma$ for some $b_{ij\alpha}^\gamma \in \Gamma\left( U_{ij}, \mathcal{O}_X \right)$. By multiplying $\sum_{\gamma=1}^q \left( \delta_\gamma^\eta- t b_{ij\gamma}^{ \eta} \right)$, we have
\begin{align*}
&\sum_{\beta, \gamma=1}^p \tilde{f}_{ij}\left( \tilde{r}_{ij}^{\alpha\beta} \right) \tilde{r}_{ji}^{\beta \gamma}\left(\delta_\gamma^\eta - t b_{ij\gamma}^{ \eta} \right) = \sum_{\gamma=1}^p \left( \delta_\alpha^\gamma +t b_{ij\alpha}^{ \gamma} \right) \left( \delta_\gamma^\eta - t b_{ij \gamma}^{ \eta} \right)  = \delta_\alpha^\eta + t b_{ij\alpha}^{ \eta}- t b_{ij\alpha}^{ \eta} = \delta_\alpha^\eta\\
&\Longrightarrow \sum_{\beta=1}^p \tilde{f}_{ij}\left( \tilde{r}_{ij}^{\alpha\beta} \right)\left(\tilde{r}_{ji}^{\beta \eta} - t \sum_{\gamma=1}^p r_{0ji}^{\beta \gamma} b_{ij\gamma}^{ \eta} \right) = \delta_\alpha^\eta
\end{align*}
Then by replacing $\tilde{r}_{ji}^{\beta \eta}$ by $\tilde{r}_{ji}^{\beta \eta} - t \sum_{\gamma=1}^p r_{0ji}^{\beta \gamma} b_{ij\gamma}^{ \eta}$ which still induces $r_{ji}^{\beta \eta}$, so that we may assume that
\begin{align}\label{ts8}
\sum_{\beta=1}^p \tilde{f}_{ij}\left(\tilde{r}_{ij}^{\alpha\beta} \right)\tilde{r}_{ji}^{\beta \gamma} = \delta_\alpha^\gamma
\end{align}
We define an element $W_{ij}\in \Gamma\left( U_{ij}, \mathscr{H}om_{\mathcal{O}_X}\left( T_{\mathcal{F}_0 } , T_X  \right) \right)$ by $W_{ij}:\Gamma\left( U_{ij}, T_{\mathcal{F}_0} \right)\to \Gamma\left( U_{ij}, T_X \right), T_{0i}^\alpha \mapsto W_{ij}^\alpha$. Then we have
\begin{align}
tW_{ij}^\alpha&= \tilde{\theta}_{ij}^*\left(t W_{ij}^\alpha \right)=\tilde{\theta}_{ij}^*\tilde{\theta}_{ji}^*\left(\tilde{T}_i^\alpha\right)-\sum_{\beta=1}^p \tilde{f}_{ij}\left(\tilde{r}_{ij}^{\alpha\beta}\right)\tilde{\theta}_{ij}^*\left(\tilde{T}_j^\beta\right)=\tilde{T}_i^\alpha-\sum_{\beta=1}^p \tilde{f}_{ij} \left( \tilde{r}_{ij}^{\alpha\beta} \right) \left(tW_{ji}^\beta+\sum_{\gamma=1}^p \tilde{r}_{ji}^{\beta\gamma} \tilde{T}_i^\gamma \right) \label{tn5}\\
   &=\tilde{T}_i^\alpha-t\sum_{\beta=1}^p r_{0ij}^{\alpha\beta} W_{ji}^\beta -\sum_{\beta,\gamma=1}^p  \tilde{f}_{ij}\left(\tilde{r}_{ij}^{\alpha\beta} \right) \tilde{r}_{ji}^{\beta\gamma} \tilde{T}_i^\gamma = - t\sum_{\beta=1}^p r_{0ij}^{\alpha\beta} W_{ji}^\beta \notag
   \end{align}
This implies that $\left\{ W_{ij} \right\}\in C^1\left(\mathcal{U}, \mathscr{H}om_{\mathcal{O}_X}\left( T_{\mathcal{F}_0}, T_X \right) \right)$. 

On the other hand, we also recall that $\sum_{\beta, \gamma=1}^p  \left(\tilde{f}_{ij} \tilde{f}_{jk} \right)\left(\tilde{r}_{ik}^{\alpha\beta} \right) \tilde{f}_{ij}\left( \tilde{r}_{kj}^{\beta \gamma} \right) \tilde{r}_{ji}^{\gamma \eta} = \delta_\alpha^\eta + t g_{ijk}^{\alpha \eta}$ for some $g_{ijk}^{\alpha \eta} \in \Gamma\left( U_{ijk}, \mathcal{O}_X \right)$. Then we have
\begin{align*}
&\sum_{\beta, \gamma=1}^p \left( \tilde{f}_{ij} \tilde{f}_{jk} \right)\left(  \tilde{r}_{ik}^{\alpha\beta} \right) \tilde{f}_{ij}\left( \tilde{r}_{kj}^{\beta \gamma} \right) \tilde{r}_{ji}^{\gamma \eta} = \delta_\alpha^\eta + t g_{ijk}^{\alpha \eta} \\
&\iff \sum_{\beta, \gamma=1}^p \tilde{f}_{jk}\left(\tilde{r}_{ik}^{\alpha\beta} \right) \cdot \tilde{r}_{kj}^{\beta \gamma} \cdot \tilde{f}_{ji}\left(\tilde{r}_{ji}^{\gamma \eta} \right) = \delta_\alpha^\eta + t g_{ijk}^{\alpha \eta} \\
& \iff \sum_{\beta, \gamma=1}^p \tilde{r}_{ik}^{\alpha\beta} \cdot \tilde{f}_{kj}\left(\tilde{r}_{kj}^{\beta \gamma} \right) \tilde{f}_{kj} \tilde{f}_{ji}\left( \tilde{r}_{ji}^{\gamma \eta} \right) = \delta_\alpha^\eta + t g_{ijk}^{\alpha \eta}
\end{align*}
We define an element $g_{ijk}\in \Gamma\left( U_{ijk}, \mathscr{H}om_{\mathcal{O}_X}\left( T_{\mathcal{F}_0}, T_{\mathcal{F}_0} \right) \right)$ by
\begin{align*}
g_{ijk}:\Gamma\left(U_{ijk}, T_{\mathcal{F}_0}  \right) &\to \Gamma\left( U_{ijk} , T_{\mathcal{F}_0} \right) \\
 T_{0i}^\alpha &\mapsto \sum_{\eta=1}^p g_{ijk}^{\alpha \eta} T_{0i}^\eta
\end{align*}
and linearly extends to $\Gamma\left( U_{ijk}, T_{\mathcal{F}_0} \right)$ in terms of the basis $\left\{ T_{0i}^\alpha \right\}$. Then we have
\begin{align*}
g_{jkl}\left( T_{0i}^\alpha \right)=\sum_{\beta=1}^p r_{0ij}^{\alpha\beta} g_{jkl}\left(T_{0j}^\beta \right) = \sum_{\beta, \gamma=1}^p r_{0ij}^{\alpha\beta} g_{jkl}^{\beta \gamma} T_{0j}^\gamma = \sum_{\beta, \gamma, \eta=1}^p r_{0ij}^{\alpha\beta} g_{jkl}^{\beta \gamma} r_{0ji}^{\gamma \eta} T_{0i}^\eta
\end{align*}
We consider $\left(\left\{d_{ijk} \right\}, \left\{ g_{ijk} \right\} \right)\in T_X|_{U_{ijk}}\bigoplus \mathscr{H}om_{\mathcal{O}_X}\left( T_{\mathcal{F}_0}, T_{\mathcal{F}_0} \right)|_{U_{ijk}} \subset T_X|_{U_i} \bigoplus \mathscr{H}om_{\mathcal{O}_X}\left( T_{\mathcal{F}_0}, T_{\mathcal{F}_0} \right)|_{U_i}$ as an element of $\Gamma\left( U_{ijk}, \mathcal{E}_{T_{\mathcal{F}_0}}\right)$. Then we claim that 
\begin{align*}
 \left( \left\{ P_i \right\}, \left\{ - W_{ij} \right\}, \left\{ \left(d_{ijk}, g_{ijk} \right) \right\} \right) \in C^0\left( \mathcal{U}, \mathscr{H}om_{\mathcal{O}_X}\left( \bigwedge^2 T_{\mathcal{F}_0}, \frac{T_X}{T_{\mathcal{F}_0}} \right)\right) \bigoplus C^1\left( \mathcal{U}, \mathscr{H}om_{\mathcal{O}_X}\left( T_{\mathcal{F}_0}, T_X \right) \right) \bigoplus C^2\left( \mathcal{U}, \mathcal{E}_{T_{\mathcal{F}_0}} \right)
\end{align*}
defines a $2$-cocycle in the above \v Cech resolution of $\mathcal{E}_{T_{\mathcal{F}_0}}^\bullet$. It is sufficient to show that $-\delta\left( \left\{ \left( d_{ijk}, g_{ijk} \right) \right\} \right)=0$. In fact, if we set matrices $G_{ijk}:=\left( g_{ijk}^{\alpha \beta} \right)$, $\tilde{r}_{ij}= \left( \tilde{r}_{ij}^{\alpha\beta} \right)$, and $r_{0ij}= \left( r_{0ij}^{\alpha\beta} \right)$, and $\textnormal{Id}$ is the identity matrix, we have (in the following, for example,  $\tilde{f}_{ij}\left( \tilde{r}_{lj} \right)$ means the matrix $\left( \tilde{f}_{ij}\left( r_{lj}^{\alpha\beta} \right) \right)$, and $\left[ {d}_{ijk},  r_{0kj} \right]$ means the matrix $\left( \left[ {d}_{ijk}, r_{0kj}^{\alpha\beta} \right] \right)$):

\begin{align*}
&\left(\textnormal{Id}+ t G_{ilk} \right)\left( \textnormal{Id} + t G_{ijl} \right)\left(\textnormal{Id} + t r_{0ij} G_{jkl} r_{0ji}\right) \left( \textnormal{Id} + t G_{ikj}  \right)\\
&=\left(\textnormal{Id}+ t G_{ilk} \right)\left( \textnormal{Id} + t G_{ijl} \right)\tilde{f}_{lj}\left(\tilde{r}_{ij} \right) \left(\textnormal{Id} + t G_{jkl} \right)  \tilde{f}_{lj}\tilde{f}_{ji}\left(\tilde{r}_{ji} \right)\tilde{f}_{lj}\left( \textnormal{Id} + t  G_{ikj}  \right)\\
&= \tilde{f}_{lk}\left(\tilde{r}_{ik} \right)\cdot \tilde{r}_{kl}  \cdot \tilde{f}_{li}\left( \tilde{r}_{li}\right)  \cdot \tilde{r}_{il} \cdot \tilde{f}_{lj}\left(\tilde{r}_{lj} \right) \cdot \tilde{f}_{lj} \tilde{f}_{ji}\left( \tilde{r}_{ji}\right) \cdot \tilde{f}_{lj}\left( \tilde{r}_{ij}\right) \cdot \tilde{r}_{jl} \cdot \tilde{f}_{lk}\left( \tilde{r}_{lk} \right) \tilde{f}_{lk} \tilde{f}_{kj}\left( \tilde{r}_{kj} \right)\\
& \cdot \tilde{f}_{lj}\tilde{f}_{ji}\left(\tilde{r}_{ji} \right) \cdot \tilde{f}_{lj}\left( \tilde{r}_{ij}\right) \cdot \tilde{f}_{lj}\tilde{f}_{jk}\left( \tilde{r}_{jk} \right) \cdot \tilde{f}_{lj} \tilde{f}_{jk} \tilde{f}_{ki}\left(\tilde{r}_{ki}\right)\\
&= \tilde{f}_{lk}\left(\tilde{r}_{ik} \right)\cdot \tilde{r}_{kl} \cdot \tilde{f}_{lj}\left(\tilde{r}_{lj} \right)  \cdot \tilde{r}_{jl} \cdot  \tilde{f}_{lk}\left( \tilde{r}_{lk} \right) \tilde{f}_{lk} \tilde{f}_{kj}\left( \tilde{r}_{kj} \right) \cdot \tilde{f}_{lj}\tilde{f}_{jk}\left( \tilde{r}_{jk} \right) \cdot \tilde{f}_{lj} \tilde{f}_{jk} \tilde{f}_{ki}\left(\tilde{r}_{ki}\right)\\
&= \tilde{f}_{lk}\left(\tilde{r}_{ik} \right)\cdot \left(\tilde{f}_{lj} + t d_{lkj} \right)\left( \tilde{r}_{kj} \right) \cdot \tilde{f}_{lj}\tilde{f}_{jk}\left( \tilde{r}_{jk} \right) \cdot \left( \tilde{f}_{lk} + t d_{ljk}\right) \tilde{f}_{ki}\left(\tilde{r}_{ki}\right)\\
&= \textnormal{Id} + tr_{0ik}\left[ d_{lkj},  r_{0kj} \right] r_{0ji} + t r_{0ik} \left[ d_{ljk}, r_{0ki}\right]= \textnormal{Id}  - t r_{0ik}\left[ d_{jkl} ,  r_{0kj} \right] r_{0ji} + t r_{0ik}\left[ d_{jkl}, r_{0ki} \right]\\
&= \textnormal{Id} + t r_{0ij}\left[ d_{jkl}, r_{0ji} \right] = \textnormal{Id} - \left[ d_{jkl},  r_{0ij} \right] r_{0ji}
\end{align*}

Then we have
\begin{align}
& G_{ilk}+ G_{ijl} + r_{0ij} G_{jkl} r_{0ji} + \left[ d_{jkl},  r_{0ij} \right] r_{0ji} + G_{ikj} = 0 \notag \\
&\iff  - G_{ikl}+  G_{ijl} + r_{0ij} G_{jkl} r_{0ji} + \left[ d_{jkl},  r_{0ij} \right] r_{0ji} - G_{ijk} = 0 \label{at26}
\end{align}

We note that $\left(d_{jk}, g_{jkl} \right)$ on $U_j$ is identified with $\left(d_{jkl}, g_{jkl}^i\right)$ on $U_i$ by
\begin{align*}
g_{jkl}^i\left( T_{0i}^\alpha\right) & =  g_{jkl}\left( T_{0i}^\alpha\right) + \sum_{\beta, \gamma=1}^p \left[ d_{jkl}, r_{0ij}^{\alpha\beta} \right] r_{0ji}^{\beta \gamma} T_{0i}^\gamma = \sum_{\beta=1}^p r_{0ij}^{\alpha\beta} g_{jkl}\left(T_{0j}^\beta\right) +  \sum_{\beta, \eta=1}^p \left[ d_{jkl}, r_{0ij}^{\alpha\beta} \right] r_{0ji}^{\beta \eta} T_{0i}^\eta \\
 &= \sum_{\beta, \gamma, \eta=1}^p r_{0ij}^{\alpha\beta} g_{jkl}^{\beta \gamma} r_{0ji}^{\gamma \eta} T_{0i}^\eta +  \sum_{\beta, \eta=1}^p \left[ d_{jkl}, r_{0ij}^{\alpha\beta} \right] r_{0ji}^{\beta \eta} T_{0i}^\eta = \sum_{\eta=1}^p \left( \sum_{\beta, \gamma=1}^p r_{0ij}^{\alpha\beta} g_{jkl}^{\beta \gamma} r_{0ji}^{\gamma \eta} + \sum_{\beta=1}^p \left[d_{jkl}, r_{0ij}^{\alpha\beta} \right] r_{0ji}^{\beta \eta} \right) T_{0i}^\eta
\end{align*}

Hence we have from $(\ref{at26})$
\begin{align*}
- g_{ikl}+ g_{ijl}  + g_{jkl}^i - g_{ijk} = 0
\end{align*}
This implies that $-\delta \left( \left\{d_{ijk}, g_{ijk}  \right\} \right)=0$. This prove the claim. We show that for another lifting $\tilde{\theta}_{ij}'$ of $\theta_{ij}$, and another lifting $\tilde{T}_i'^\alpha$ of $T_i^\alpha$, and another lifting of $\tilde{r}_{ij}'^{\alpha\beta}$ of $r_{ij}^{\alpha\beta}$, and another lifting $\tilde{g}_{i\alpha\beta}'^\gamma$ of $g_{i\alpha\beta}^\gamma$, the associated $2$-cocycle $\beta:=\left(\left\{ P_i' \right\}, \left\{ - W_{ij}' \right\}, \left\{ \left( d_{ijk}' , g_{ijk}'  \right) \right\} \right)$ is cohomologus to $\alpha:=\left( \left\{ P_i \right\}, \left\{ - W_{ij}\right\}, \left\{ \left( d_{ijk}, g_{ijk} \right) \right\} \right)$. We define $s_{ij}\in \Gamma\left( U_{ij}, \mathscr{H}om_{\mathcal{O}_X}\left( T_{\mathcal{F}_0}, T_{\mathcal{F}_0} \right) \right), T_{0i}^\alpha \mapsto \sum_{\beta=1}^p s_{ij}^{\alpha\beta} T_{0j}^\beta$. Then we consider $\left( p_{ij}, - s_{ij} \right)\in T_X|_{U_{ij}}$ $\bigoplus \mathscr{H}om_{\mathcal{O}_X}\left(T_{\mathcal{F}_0}, T_{\mathcal{F}_0} \right)|_{U_{ij}}\subset T_X|_{U_i} \bigoplus \mathscr{H}om_{\mathcal{O}_X}\left( T_{\mathcal{F}_0} , T_{\mathcal{F}_0} \right)|_{U_i}$ as an element of $\Gamma\left( U_{ij}, \mathcal{E}_{T_{\mathcal{F}_0}} \right)$.  Then we have $D_1\left( Y_i \right)= P_i'- P_i$ and $\delta\left( \left\{ Y_i\right\}\right)+ D_0'\left( \left\{ \left( p_{ij}, - s_{ij} \right) \right\} \right)=\left\{- W_{ij}'+ W_{ij} \right\}$ and we claim that $\delta\left(\left\{\left( p_{ij}, - s_{ij} \right) \right\} \right) =\left\{\left( d_{ijk}'- d_{ijk},  g_{ijk}' - g_{ijk} \right)    \right\}$. In fact, we have $\delta\left( \left\{p_{ij}\right\}\right)= \left\{ d_{ijk}'- d_{ijk}\right\}$. On the other hand, we note that
\begin{align*}
\textnormal{Id}+ tg_{ijk}'^{\alpha \eta} & = \sum_{\beta, \gamma=1}^p \left( \tilde{f}_{ij}' \tilde{f}_{jk}' \right)\left(  \tilde{r}_{ik}'^{\alpha\beta} \right) \tilde{f}_{ij}'\left( \tilde{r}_{kj}'^{\beta \gamma} \right) \tilde{r}_{ji}'^{\gamma \eta}\\
& = \sum_{\beta,\gamma=1}^p \left( \tilde{f}_{ij}+ t p_{ij} \right)\left( \tilde{f}_{jk} + t p_{jk}\right)\left( \tilde{r}_{ik}^{\alpha\beta} + t s_{ik}^{\alpha \beta}\right) \cdot \left( \tilde{f}_{ij} + t p_{ij} \right)\left( \tilde{r}_{kj}^{\beta \gamma} + t s_{kj}^{\beta \gamma} \right) \cdot \left( \tilde{r}_{ji}^{\gamma \eta} + t s_{ji}^{\gamma \eta}\right) \\
&= \textnormal{Id} + t g_{ijk}^{\alpha \eta} + t \sum_{\beta=1}^p \left[ p_{ij}, r_{0ik}^{\alpha\beta}  \right] r_{0ki}^{\beta \eta} + t\sum_{\beta=1}^p \left[ p_{jk}, r_{0ik}^{\alpha \beta} \right] r_{0ki}^{\beta \eta} + t \sum_{\beta, \gamma=1}^p r_{0ik}^{\alpha \beta} \left[ p_{ij}, r_{0kj}^{\beta \gamma} \right] r_{0ji}^{\beta \eta} \\
&+ t \sum_{\beta=1}^p s_{ik}^{\alpha \beta} r_{0ki}^{\beta \eta} + t\sum_{\beta=1}^p r_{0ik}^{\alpha\beta} s_{kj}^{\beta \gamma} r_{0ji}^{\gamma \eta} +t \sum_{\beta=1}^p r_{0ij}^{\alpha \beta} s_{ji}^{\beta \eta}
\end{align*}
Then we have
\begin{align*}
g_{ijk}'\left( T_{0i}^\alpha \right) &= g_{ijk}\left( T_{0i}^\alpha \right) + \sum_{\beta=1}^p \left[ p_{ij}, r_{0ij}^{\alpha \beta} \right] T_{0j}^\beta + \sum_{\beta=1}^p \left[ p_{jk}, r_{0ik}^{\alpha\beta} \right] T_{0k}^\beta + s_{ik}\left( T_{0i}^\alpha \right) + \sum_{\beta=1}^p r_{0ik}^{\alpha \beta} s_{kj}\left( T_{0k}^\beta \right) + \sum_{\beta=1}^p r_{0ij}^{\alpha\beta} s_{ji}\left( T_{0j}^\beta \right)\\
&= g_{ijk}\left( T_{0i}^\alpha \right)  - \sum_{\beta=1}^p \left[ p_{ji}, r_{0ij}^{\alpha \beta} \right] T_{0j}^\beta - \sum_{\beta=1}^p \left[ p_{kj}, r_{0ik}^{\alpha\beta} \right] T_{0k}^\beta  + s_{ik}\left( T_{0i}^\alpha \right) + s_{kj}\left(  T_{0i}^\alpha \right) + s_{ji}\left( T_{0i}^\alpha \right)\\
&= g_{ijk}\left( T_{0i}^\alpha \right) - s_{ij}\left( T_{0i}^\alpha \right) + s_{ik}\left( T_{0i}^\alpha \right) -\left( - s_{kj}\left(T_{0i}^\alpha \right) +\sum_{\beta=1}^p \left[ p_{kj}, r_{0ik}^{\alpha\beta} \right] T_{0k}^\beta\right)
\end{align*}
This proves that $\alpha$ is cohomologous to $\beta$. So given a small extension $e:0\to (t)\to \tilde{A} \to A \to 0$ and an infinitesimal foliated deformation $\eta$ of $\left( X, T_{\mathcal{F}_0} \right)$ over $A$, we can associate an element $o_\eta(e):=$ the cohomology class of $\alpha \in \mathbb{H}^2\left( X, T_{\mathcal{F}_0}^\bullet \right)$. We note that $o_\eta(e)=0$ if and only if there exists $\left(\left\{ X_i \right\}, \left\{ \left(  d_{ij} , -v_{ij} \right) \right\} \right) \in C^0\left( \mathcal{U}, \mathscr{H}om_{\mathcal{O}_X}\left( T_{\mathcal{F}_0}, T_X \right) \right) \bigoplus C^1\left( \mathcal{U}, \mathcal{E}_{T_{\mathcal{F}_0}} \right)$ where $X_i:\Gamma\left( U_i, T_{\mathcal{F}_0} \right) \to \Gamma\left( U_i, T_X \right), T_{0i}^\alpha \mapsto X_i^\alpha $  and $v_{ij}:\Gamma\left( U_{ij}, T_{\mathcal{F}_0} \right) $ $\to \Gamma\left( U_{ij} , T_{ \mathcal{F}_0} \right), T_{0i}^\alpha \mapsto \sum_{\beta=1}^p v_{ij}^{\alpha\beta} T_{0j}^\beta$ such that $D_1'\left(X_i\right)= P_i$ and $X_j- X_i + D_0'\left( \left(d_{ij}, -v_{ij}\right) \right)= - W_{ij}$ and $d_{ij}- d_{ik} + d_{jk}= d_{ijk}$ and $-v_{ij}+ v_{ik} - v_{jk}^i= g_{ijk}$. This implies a lifting $\tilde{\eta}$ of $\eta$ exists. This completes the proof of Proposition \ref{at27}.

\end{proof}

\subsection{Relationships between foliated deformations of $\left(X, T_{\mathcal{F}_0} \right)$ and the underlying flat deformations of a pair $\left(X, T_{\mathcal{F}_0} \right)$ of a nonsingular variety $X$ and a locally free sheaf $\mathcal{N}_{\mathcal{F}_0}^*$ }\

Let $\left( X, T_{\mathcal{F}_0} \right)$ be a foliated nonsingular variety with $T_{\mathcal{F}_0}$ locally free. Let $\textnormal{Def}_{\left(X, T_{\mathcal{F}_0} \right)}$ be the functor of Artin rings of flat deformations of the pair of the underlying nonsingular variety $X$ and the underlying locally free sheaf $T_{\mathcal{F}_0}$. As mentioned above, $\textnormal{Def}_{\left( X, T_{\mathcal{F}_0} \right)}$ is controlled by the Atiyah extension $\mathcal{E}_{T_{\mathcal{F}_0}}$. We have a morphism of deformation functors
\begin{align*}
\textnormal{FDef}_{\left( X, T_{\mathcal{F}_0} \right)}\to \textnormal{Def}_{\left( X, T_{\mathcal{F}_0} \right)}
\end{align*}
by forgetting foliated structures and we have a canonical morphism of deformation functors
\begin{align*}
\textnormal{FDef}_{\left( X, T_{\mathcal{F}_0} \right)_{\textnormal{fixed}}} \to \textnormal{FDef}_{\left( X, T_{\mathcal{F}_0} \right)}
\end{align*}
Then the relationship of three deformation functors $\textnormal{FDef}_{\left( X, T_{\mathcal{F}_0}\right)_{\textnormal{fixed}}}, \textnormal{FDef}_{\left( X, T_{\mathcal{F}_0} \right)}$ and $\textnormal{Def}_{\left(X, T_{\mathcal{F}_0}\right)}$ is expressed by the following exact sequence of complexes of sheaves
{\Small{\begin{center}
$\begin{CD}
@. \cdots @. \cdots @. \cdots \\
@. @AD_3AA @AD_3AA @AAA \\
0@>>> \mathscr{H}om_{\mathcal{O}_X}\left(\bigwedge^3 T_{\mathcal{F}_0}, \frac{T_X}{T_{\mathcal{F}_0}} \right) @>>> \mathscr{H}om_{\mathcal{O}_X}\left(\bigwedge^3 T_{\mathcal{F}_0}, \frac{T_X}{T_{\mathcal{F}_0}} \right) @>>> 0 @>>> 0 \\
@. @AD_2AA @AD_2AA @AAA \\
0@>>> \mathscr{H}om_{\mathcal{O}_X}\left(\bigwedge^2 T_{\mathcal{F}_0}, \frac{T_X}{T_{\mathcal{F}_0}} \right) @>>> \mathscr{H}om_{\mathcal{O}_X}\left( \bigwedge^2 T_{\mathcal{F}_0}, \frac{T_X}{T_{\mathcal{F}_0}}\right) @>>> 0 @>>> 0 \\
@. @AD_1'AA @AD_1AA @AAA \\
0@>>> \mathscr{H}om_{\mathcal{O}_X}\left(T_{\mathcal{F}_0}, T_X \right) @>>> \mathscr{H}om_{\mathcal{O}_X}\left( T_{\mathcal{F}_0}, T_X \right) @>>> 0 @>>> 0  \\
@. @AAA @AD_0'AA @AAA \\
0 @>>> 0 @>>> \mathcal{E}_{T_{\mathcal{F}_0}} @>>> \mathcal{E}_{T_{\mathcal{F}_0}} @>>> 0
\end{CD}$
\end{center}}}
Then the long exact sequence induces maps of first-order deformations
{\small{\begin{align*}
\textnormal{FDef}_{\left( X, T_{\mathcal{F}_0} \right)_{\textnormal{fixed}}}\left(k[\epsilon] \right)\cong \mathbb{H}^0\left( X, \mathscr{H}om_{\mathcal{O}_X}\left( T_{\mathcal{F}_0}, T_X \right)^\bullet \right) \to \textnormal{FDef}_{\left( X, T_{\mathcal{F}_0} \right)} \left( k[\epsilon] \right) \cong \mathbb{H}^1\left( X, \mathcal{E}_{T_{\mathcal{F}_0}}^\bullet  \right) \to \textnormal{Def}_{\left( X, T_{\mathcal{F}_0} \right)}\left(k[\epsilon] \right)\cong H^1\left( X, \mathcal{E}_{ T_{\mathcal{F}_0}} \right)
\end{align*}}}
and obstruction maps
\begin{align*}
\mathbb{H}^1\left( X, \mathscr{H}om_{\mathcal{O}_X}\left( T_{\mathcal{F}_0}, T_X \right)^\bullet \right) \to \mathbb{H}^2\left( X, \mathcal{E}_{ T_{\mathcal{F}_0}}^\bullet \right) \to H^2\left( X, \mathcal{E}_{ T_{\mathcal{F}_0}} \right)
\end{align*}
of morphisms of deformation functors $\textnormal{FDef}_{\left( X, T_{\mathcal{F}_0} \right)_{\textnormal{fixed}}} \to \textnormal{FDef}_{\left( X, T_{\mathcal{F}_0} \right)} \to \textnormal{Def}_{\left( X , T_{\mathcal{F}_0} \right) }$.

\section{Flat foliated deformations of $\left(X, T_{\mathcal{F}_0} \right)$ with $T_{\mathcal{F}_0}$ locally free in terms of deformations of sections} \label{nt1}

Let $\left(X, T_{\mathcal{F}_0} \right)$ be a foliated nonsingular variety with $T_{\mathcal{F}_0}$ locally free of rank $p$. The exact sequence $0\to T_{\mathcal{F}_0}\to T_X \to \frac{T_X}{T_{\mathcal{F}_0}}\to 0$ induces a section $T_0 \in H^0\left( X, \bigwedge^p T_X \otimes \mathcal{Q}_0 \right)$ where $Q_0:= \left(\bigwedge^p T_{\mathcal{F}_0}\right)^*$. An infinitesimal flat foliated deformation $\left(\mathcal{X}, T_\mathcal{F}\right)$ of $\left( X, T_{\mathcal{F}_0} \right)$ over $A\in \bold{Art}$ defined by the exact sequence $0\to T_{\mathcal{F}} \to T_{\mathcal{X}/A} \to \frac{T_{\mathcal{X}/A}}{T_{\mathcal{F}}}\to 0$ induces a section $T\in H^0\left(\mathcal{X}, \bigwedge^p T_{\mathcal{X}/A}\otimes \mathcal{Q} \right)$ where $\mathcal{Q}:= \left(\bigwedge^p T_{\mathcal{F}}\right)^*$ which is an infinitesimal deformation of the section $T_0$ and a pair of the nonsingular variety $X$ and a line bundle $Q_0$.  We will describe infinitesimal foliated deformations of $\left(X, T_{\mathcal{F}_0} \right)$ in terms deformations of the section $T_0$ and the pair of $\left(X, \mathcal{Q}_0 \right)$. Let $\mathcal{U}= \left\{U_i\right\}$ be an affine open covering of $X$ such that $\Gamma\left(U_i, T_{\mathcal{F}_0} \right)$ is generated by $T_{0i}^1,..., T_{0i}^p\in \Gamma\left(U_i, T_X \right)$ with $T_{0i}^\alpha = \sum_{\beta=1}^p r_{0ij}^{\alpha\beta}T_{0j}^\beta$ for some $r_{0ij}^{\alpha\beta}\in \Gamma\left( U_{ij}, \mathcal{O}_X \right)$. Let $\det\left(r_{0ij}^{\alpha\beta}\right)= r_{0ij}$.  Then we have
\begin{align*}
T_{0i}^1\wedge \cdots \wedge T_{0i}^p = r_{0ij} T_{0j}^1 \wedge \cdots \wedge T_{0j}^p
\end{align*}
Then $\mathcal{Q}_0=\left(\bigwedge^p T_{\mathcal{F}_0} \right)^*$ is defined by the transition functions $\{r_{0ij}\}$ with $r_{0ij}\in \Gamma(U_{ij}, \mathcal{O}_X^*)$ and $T_0$ is defined locally on $U_i$ by $T_0:=\left\{ T_{0i}^1\wedge \cdots \wedge T_{0i}^p\right\}\in H^0\left(X, \bigwedge^p T_X \otimes \mathcal{Q}_0 \right)$. As is known, deformations of a pair of a nonsingular variety $X$ and a line bundle $\mathcal{Q}_0$ is controlled by the Atiyah extension (see \cite{Ser06}) which is the exact sequence
\begin{align}\label{ts2}
0\to \mathcal{O}_X \to \mathcal{E}_{\mathcal{Q}_0}\to T_X \to 0
\end{align}
where $\mathcal{E}_{\mathcal{Q}_0}$ is defined locally on $U_i$
\begin{align*}
\mathcal{E}_{\mathcal{Q}_0}|_{U_i}:= \mathcal{O}_X|_{U_i} \bigoplus T_X|_{U_i}
\end{align*}
and $\left(a_i, b_i \right)\in \Gamma\left(U_i, \mathcal{O}_X\right)\bigoplus \Gamma\left(U_i, T_X\right)$ is identified with $\left(a_j, b_j\right)\in \Gamma\left(U_j, \mathcal{O}_X\right)\bigoplus \Gamma\left(U_j, T_X\right)$ if and only if $b_i=b_j$ and $a_j=a_i+ \left[b_i ,  \log r_{0ij} \right] $. In the following we will construct a complex of sheaves also controlling foliated deformations of $\left(X, T_{\mathcal{F}_0} \right)$ which explicitly shows the underlying deformations of the section $T_0$ and the line bundle $\mathcal{Q}_0$:
\begin{align}\label{ts1}
\mathcal{Q}_0^\bullet : \mathcal{E}_{\mathcal{Q}_0} \xrightarrow{\tilde{D}_0} \mathcal{R} \xrightarrow{\tilde{D}_1} \mathscr{H}om_{\mathcal{O}_X}\left(\bigwedge^2 T_{\mathcal{F}_0}, \frac{T_X}{T_{\mathcal{F}_0}} \right) \xrightarrow{D_2} \mathscr{H}om_{\mathcal{O}_X}\left( \bigwedge^3 T_{\mathcal{F}_0}, \frac{T_X}{T_{\mathcal{F}_0}} \right) \xrightarrow{D_3} \cdots
\end{align}
We will denote the $i$-th cohomology group of $\mathcal{E}_{\mathcal{Q}_0}^\bullet$ by $\mathbb{H}^i\left( X, \mathcal{E}_{\mathcal{Q}_0}^\bullet \right)$. We define a subsheaf $\mathcal{R} \subset  \bigwedge^p T_X \otimes \mathcal{Q}_0$ which is defined in the following way: $\Gamma(U_i, \mathcal{R})\subset \Gamma\left(U_i, \bigwedge^p T_X \otimes \mathcal{Q}_0^*\right)$ is generated by 
\begin{align*}
T_{0i}^1 \wedge \cdots \wedge \widehat{T_{0i}^\alpha} \wedge \cdots \wedge T_{0i}^p,\,\,\,\,\,\,\,\,\,\,\alpha=1,...,p.
\end{align*}
This implies that any element in $\Gamma\left(\mathcal{U}_i, \mathcal{R}\right)$ is of the form
\begin{align*}
a_i=\sum_{\alpha=1}^p T_{0i}^1 \wedge \cdots \wedge a_i^\alpha \wedge \cdots \wedge T_{0i}^p,\,\,\,\,\,\,\,\,a_i^\alpha \in \Gamma(U_i, T_X)
\end{align*}
and $a_i$ is transformed into $\Gamma\left(U_j, \mathcal{R}\right)$ in the following way:
\begin{align*}
a_j:=\frac{1}{r_{0ij}}a_i=  \frac{1}{r_{0ij}} \sum_{\alpha=1}^p T_{0i}^1 \wedge \cdots \wedge a_i^\alpha \wedge \cdots \wedge T_{0i}^p= \sum_{\beta=1}^p T_{0j}^1 \wedge \cdots \wedge a_j^\beta \wedge \cdots T_{0j}^p,\,\,\,\,\,\,\,\,\,\textnormal{for some}\,\,\,a_j^\beta\in \Gamma(U_i\cap U_j, \mathcal{O}_X)
\end{align*}
so that $\mathcal{R}$ is well-defined. On the other hand, we define a map
\begin{align*}
\beta:\mathcal{R} \to \mathscr{H}om_{\mathcal{O}_X}\left( T_{\mathcal{F}_0}, \frac{ T_X }{ T_{\mathcal{F}_0} } \right)
\end{align*}
which is defined locally in the following way: for $a_i\in \Gamma(U_i, \mathcal{R})$ as above, we define
\begin{align*}
\beta:\Gamma(U_i, \mathcal{R})&\to \Gamma\left( U_i, \mathscr{H}om_{\mathcal{O}_X}\left( T_{\mathcal{F}_0}, \frac{T_X}{T_{\mathcal{F}_0} } \right) \right) \\
a_i&\mapsto \left( T_{0i}^\alpha \mapsto \overline{a_i^\alpha} \right)
\end{align*}
and linearly extends to $\Gamma\left( U_i, T_{\mathcal{F}_0} \right)$. Then $\beta$ is well-defined: In fact,  if we have
\begin{align*}
a_i = \sum_{\alpha=1}^p  T_{0i}^1 \wedge \cdots \wedge a_i^\alpha \wedge \cdots \wedge T_{0i}^p = \sum_{\alpha=1}^p T_{0i}^1 \wedge \cdots \wedge a_i'^\alpha \wedge \cdots \wedge T_{0i}^q,
\end{align*}
then we have 
\begin{align*}
\left( a_i^\alpha - a_i'^\alpha \right) \wedge T_{0i}^1 \wedge \cdots \wedge T_{0i}^p =0
\end{align*}
Then for $x\in U_i-S$ and a neighborhood $U_{i_x}\subset  U_i-S$, we can write 
\begin{align*}
a_i^\alpha - a_i'^\alpha = \sum_{\alpha=1}^n A_{i_x}^\alpha T_{0i}^\alpha \in \Gamma \left( U_{i_x} , T_{\mathcal{F}_0}\right)
\end{align*}
where $A_{i_x}^\alpha \in \Gamma\left(U_{i_x}, \mathcal{O}_X\right)$, so that $a_i^\alpha - a_i'^\alpha \in \Gamma\left(U_i-S, T_{\mathcal{F}_0}\right)$. Since $T_{\mathcal{F}_0}$ is reflexive, we have $a_i^\alpha - a_i'^\alpha \in \Gamma\left(U_i, T_{\mathcal{F}_0}\right)$. On the other hand, $a_j$ on $U_j$ is transformed into $a_i=r_{0ij} a_j$ on $U_i$. We note that
\begin{align*}
&-\frac{1}{r_{0ij}} T_{0i}^1\wedge \cdots \wedge T_{0i}^p \wedge a_i^\alpha= \frac{1}{r_{ij}} a_i\wedge T_{0i}^\alpha=a_j\wedge T_{0i}^\alpha\\
&= \sum_{\beta=1}^p T_{0j}^1 \wedge \cdots \wedge a_j^\beta \wedge \cdots \wedge T_{0j}^p \wedge \left(\sum_{\gamma=1}^p r_{0ij}^{\alpha\gamma} T_{0j}^\gamma \right) = \sum_{\beta=1}^p T_{0j}^1 \wedge \cdots \wedge a_j^\beta \wedge \cdots \wedge T_{0j}^p \wedge \left( r_{0ij}^{\alpha\beta} T_{0j}^\beta \right)\\
                                     &= - T_{0j}^1 \wedge \cdots \wedge T_{0j}^p \wedge \left( \sum_{\beta=1}^p r_{0ij}^{\alpha\beta} a_j^\beta \right) =-\frac{1}{r_{0ij}} T_{0i}^1\wedge \cdots \wedge T_{0i}^p \wedge\left( \sum_{\beta=1}^p r_{0ij}^{\alpha\beta} a_j^\beta \right)
\end{align*}
This implies that $a_i^\alpha - \sum_{\beta=1}^p r_{0ij}^{\alpha\beta} a_j^\beta \in \Gamma \left( U_i\cap U_j, T_{\mathcal{F}_0} \right)$. Hence $\mathcal{R} \xrightarrow{\beta} \mathscr{H}om_{\mathcal{O}_X}\left( T_{\mathcal{F}_0},  \frac{T_X}{T_{\mathcal{F}_0}} \right)$ is well-defined. Then we define $\tilde{D}_1$ in $(\ref{ts1})$ by the composition of $\beta$ and $D_1: \mathscr{H}om_{\mathcal{O}_X}\left( T_{\mathcal{F}_0}, \frac{T_X}{T_{\mathcal{F}_0}} \right)\to\mathscr{H}om_{\mathcal{O}_X}\left( \bigwedge^2 T_{\mathcal{F}_0}, \frac{T_X}{T_{\mathcal{F}_0}} \right)$. We note that we have an exact sequence
\begin{align*}
0 \to \mathcal{O}_X &\to \mathcal{R} \xrightarrow{\beta} \mathscr{H}om_{\mathcal{O}_X}\left( T_{\mathcal{F}_0}, \frac{T_X}{ T_{\mathcal{F}_0}} \right) \to 0\\
                    a  &\mapsto  a T
\end{align*}
where $T_0 =\left\{T_{0i}^1\wedge \cdots \wedge T_{0i}^p \right\} \in H^0\left(X, \bigwedge^p T_{\mathcal{F}_0}\otimes \mathcal{Q}_0\right)$. In fact, It is clear that $\beta$ is surjective. Assume that $\beta(a_i)=\beta\left(\sum_{\alpha=1}^p T_{0i}^1 \wedge \cdots \wedge a_i^\alpha \wedge \cdots \wedge T_{0i}^p\right)=0$. This implies that $a_i^\alpha \in \Gamma(U_i, T_{\mathcal{F}_0})$. Then $a_i= A_i T_{0i}^1 \wedge \cdots \wedge T_{0i}^p, A_i \in  \Gamma(U_i, \mathcal{O}_X)$.

Next let us define $\tilde{D}_0$ in $(\ref{ts1})$. We define $\tilde{D}_0: \mathcal{E}_\mathcal{Q} \to \mathcal{R}$ locally on $U_i$ by
\begin{align*}
\mathcal{O}_X|_{U_i} \bigoplus T_X|_{U_i} \to & \mathcal{R}|_{U_i}\\
                 (a_i,b_i) \mapsto & \left[b_i, T_{0i}^1\wedge \cdots \wedge T_{0i}^p \right] + a_i T_{0i}^1\wedge \cdots \wedge T_{0i}^p \\
                                             &= \sum_{\alpha=1}^p T_{0i}^1 \wedge \cdots \wedge \left[b_i, T_{0i}^\alpha \right] \wedge \cdots \wedge T_{0i}^p + a_i T_{0i}^1\wedge \cdots \wedge T_{0i}^p
\end{align*}
It is clear that the image is in $\Gamma\left(U_i, \mathcal{R}\right)$. We show that it is well-defined. In fact, $(a_i, b_i)$ on $U_i$ is translated into $\left( a_i + \left[ b_i, \log r_{0ij}\right] , b_i\right)$ on $U_j$, and we have
{\small{\begin{align*}
&\frac{1}{r_{0ij}}\left(\left[ b_i, T_{0i}^1\wedge \cdots \wedge T_{0i}^p\right] + a_i T_{0i}^1\wedge \cdots \wedge T_{0i}^p \right)=\frac{1}{r_{0ij}}\left(\left[ b_i, r_{0ij}T_{0j}^1\wedge \cdots \wedge T_{0j}^p\right]+ a_ir_{0ij} T_{0j}^1\wedge \cdots \wedge T_{0j}^p \right)\\
&=\frac{1}{r_{0ij}}\left( r_{0ij}\left[ b_i, T_{0j}^1 \wedge \cdots \wedge T_{0j}^p   \right]  +\left(\left[ b_i, r_{0ij}\right] + a_i r_{0ij} \right) T_{0j}^1 \wedge \cdots \wedge T_{0j}^p   \right) = \left[b_i, T_{0j}^1\wedge \cdots \wedge T_{0j}^p \right] + \left( \left[b_i, \log r_{0ij} \right]+ a_i \right) T_{0j}^1\wedge \cdots \wedge T_{0j}^p.
\end{align*}}}
We claim that we have the commutative diagram
\begin{center}
$\begin{CD}
\mathcal{R} @>\beta>> \mathscr{H}om_{\mathcal{O}_X}\left( T_{\mathcal{F}_0}, \frac{ T_X}{T_{\mathcal{F}_0}} \right) \\
@A\tilde{D}_0 AA @AA D_0 A \\
\mathcal{E}_{\mathcal{Q}_0} @>>> T_X
\end{CD}$
\end{center}
where $\mathcal{E}_{\mathcal{Q}_0} \to T_X$ is the map in the Atiyah extension $(\ref{ts2})$. In fact, for $\left(a_i, b_i\right)\in \Gamma\left(U_i, \mathcal{E}_{\mathcal{Q}_0}\right)$, we have
\begin{align*}
\beta\left(\tilde{D}_0\left(a_i, b_i\right)\right) \left(T_{0i}^\alpha\right)= \beta\left(\sum_{\alpha=1}^p T_{0i}^1 \wedge \cdots \wedge \left[b_i, T_{0i}^\alpha\right]\wedge \cdots \wedge T_{0i}^p + a_i T_{0i}^1\wedge \cdots \wedge T_{0i}^p \right)(T_{0i}^\alpha)= \overline{\left[b_i, T_{0i}^\alpha\right]}= D_0(b_i)\left(T_{0i}^\alpha\right)
\end{align*}
This shows that $\mathcal{E}_\mathcal{Q}^\bullet$ in $(\ref{ts1})$ is a complex. Moreover we have an exact sequence of complexes of sheaves
\begin{center}
$\begin{CD}
@. 0 @. 0 @. 0 \\
@. @AAA @AD_3AA @AD_3AA \\
0 @>>> 0 @>>> \mathscr{H}om_{\mathcal{O}_X}\left( \bigwedge^3 T_{\mathcal{F}_0}, T_X/T_{\mathcal{F}_0} \right) @>>> \mathscr{H}om_{\mathcal{O}_X}\left( \bigwedge^3 T_{\mathcal{F}_0}, T_X/T_{\mathcal{F}_0} \right) @>>> 0 \\
@. @AAA @AD_2AA @AD_2AA \\
0@>>> 0 @>>> \mathscr{H}om_{\mathcal{O}_X}\left( \bigwedge^2 T_{\mathcal{F}_0}, T_X/T_{\mathcal{F}_0} \right) @>>> \mathscr{H}om_{\mathcal{O}_X}\left( \bigwedge^2 T_{\mathcal{F}_0}, \frac{ T_X}{T_{\mathcal{F}_0}} \right) @>>> 0 \\
@. @AAA @A\tilde{D}_1 AA @AD_1AA \\
0 @>>> \mathcal{O}_X @>>> \mathcal{R}@>\beta>> \mathscr{H}om_{\mathcal{O}_X}\left( T_{\mathcal{F}_0}, \frac{ T_X}{T_{\mathcal{F}_0} } \right) @>>> 0 \\
@. @A\textnormal{Id} AA @A\tilde{D}_0AA @AD_0AA \\
0 @>>> \mathcal{O}_X @>>> \mathcal{E}_{\mathcal{Q}_0} @>>> T_X @>>> 0
\end{CD}$
\end{center}
where $0\to \mathcal{O}_X\to \mathcal{E}_{\mathcal{Q}_0}\to T_X$ is the Atiyah extension. This implies that we have
\begin{align*}
\mathbb{H}^i \left(X, \mathcal{E}_{\mathcal{Q}_0}^\bullet \right) \cong \mathbb{H}^i\left(X, T_{\mathcal{F}_0}^\bullet\right)
\end{align*}

We reinterpret Proposition \ref{at25} in terms of $\mathcal{E}_{\mathcal{Q}_0}^\bullet$ instead of $T_{\mathcal{F}_0}^\bullet$.

\begin{proposition}\label{ts11}
Let $\left( X, T_{\mathcal{F}_0} \right)$ be a foliated nonsingular variety with $T_{\mathcal{F}_0}$ locally free.
\begin{enumerate}
\item there is a $1-1$ correspondence
\begin{align*}
\textnormal{FDef}_{\left(X, T_{\mathcal{F}_0}\right)}\left(k[\epsilon] \right)\cong \mathbb{H}^1\left(X, \mathcal{E}_{Q_0}^\bullet \right)
\end{align*}
\item Let $A\in \bold{Art}$ and $\eta=\left( \mathcal{X}, T_\mathcal{F} \right)$ be an infinitesimal foliated deformations of $\left( X, T_{\mathcal{F}_0} \right)$ over $A$. To every small extension $0\to (t)\to \tilde{A} \to A \to 0$, we can associate an element $o_\eta(e)\in \mathbb{H}^2\left( X,  \mathcal{E}_{\mathcal{Q}_0}^\bullet \right)$ called the obstruction lifting of $\eta$ to $\tilde{A}$, which is $0$ if and only if a lifting of $\eta$ to $\tilde{A}$ exists.
\end{enumerate}
\end{proposition}

\begin{proof}
Let $\left( \mathcal{X}, T_{\mathcal{F}} \right)$ be a first-order foliated deformation of $\left(  X, T_{\mathcal{F}_0} \right)$. We keep the notations in the proof $(1)$ of Proposition \ref{at25}. We define an element $T_i'^* :  = \sum_{\alpha=1}^p T_{0i}^1\wedge \cdots \wedge \overbrace{T_i'^\alpha}^{\alpha-\textnormal{th}} \wedge \cdots \wedge T_{0i}^q \in \Gamma\left( U_i, \mathcal{R} \right)$. By taking the wedge product of $(\ref{at19})$ for $\alpha=1,...,p$, we have
\begin{align*}
&\theta_{ji}^*\left(\left( T_{0i}^1  +\epsilon T_i'^1\right) \wedge \cdots \wedge \left( T_{0i}^p + \epsilon T_i'^p\right)     \right) \\
&=\left( T_{0i}^1+ \epsilon T_i'^1-\epsilon \left[ p_{ij}, T_{0i}^1\right] \right) \wedge \cdots \wedge \left( T_{0i}^p + \epsilon T_i'^p - \epsilon \left[ p_{ij}, T_{0i}^p \right] \right)=\det\left( r_{0ij}^{\alpha\beta}+ \epsilon b_{ij}^{\alpha\beta} \right) \left( T_{0j}^1 + \epsilon T_{j}'^1 \right) \wedge \cdots \wedge \left( T_{0j}^p + \epsilon T_j'^p \right) \\
&\Longrightarrow  \sum_{\alpha=1}^p T_{0i}^1\wedge \cdots \wedge \left( T_i'^\alpha - \left[ p_{ij}, T_{0i}^\alpha \right] \right) \wedge \cdots \wedge T_{0i}^p = r_{0ij} \left(\sum_{\alpha=1}^p T_{0j}^1 \wedge \cdots \wedge T_j'^\alpha \wedge \cdots \wedge T_{0j}^p \right) + B_{ij} T_{0j}^1 \wedge \cdots \wedge T_{0j}^p \\
&\Longrightarrow T_i'^* - \left[ p_{ij},  T_{0i}^1 \wedge \cdots \wedge T_{0i}^p \right] = r_{0ij} T_j'^* +  \frac{B_{ij}}{r_{0ij}} T_{0i}^1 \wedge \cdots \wedge T_{0i}^p
\end{align*}
Here we write $\det\left(r_{0ij}^{\alpha\beta} + \epsilon b_{ij}^{\alpha\beta} \right):= r_{0ij}+ \epsilon B_{ij}$ for some $B_{ij}\in \Gamma\left(U_{ij}, \mathcal{O}_X\right)$. We also note that 
\begin{align}\label{si31}
&\left( T_{0i}^1  +\epsilon T_i'^1\right) \wedge \cdots \wedge \left( T_{0i}^p + \epsilon T_i'^p\right)  =\theta_{ij}^* \theta_{ji}^* \left(\left( T_{0i}^1  +\epsilon T_i'^1\right) \wedge \cdots \wedge \left( T_{0i}^p + \epsilon T_i'^p \right)  \right) \\
&= \left(\textnormal{Id}+ \epsilon p_{ij} \right)\left(  r_{0ij}+ \epsilon B_{ij} \right)\cdot \left( r_{0ji} + \epsilon B_{ji} \right)\left( T_{0i}^1  +\epsilon T_i'^1\right) \wedge \cdots \wedge \left( T_{0i}^p + \epsilon T_i'^p\right) \notag\\
&\Longrightarrow \left( r_{0ij} + \epsilon B_{ij} + \epsilon \left[ p_{ij}, B_{ij} \right] \right) \cdot \left( r_{0ji} + \epsilon B_{ji}\right) = 1 \Longrightarrow B_{ij} r_{0ji}+ \left[ p_{ij}, r_{0ij} \right] r_{0ji} + r_{0ij} B_{ji}= 0 \notag \\
&\Longrightarrow -\frac{B_{ji}}{r_{0ji}}= \frac{B_{ij}}{r_{0ij}} + \left[ p_{ij}, \log r_{0ij} \right] \notag
\end{align}
This implies that $\left\{\left( p_{ij}, \frac{B_{ij}}{r_{0ij}} \right) \right\}\in C^1\left( \mathcal{U}, \mathcal{E}_{\mathcal{Q}_0} \right)$
On the other hand,
\begin{align}\label{ts3}
&\theta_{kj}^*\theta_{ji}^*\left( \left( T_{0i}^1+ \epsilon T_i'^1 \right)\wedge \cdots \wedge \left( T_{0i}^q + \epsilon T_i'^q \right) \right)= \theta_{kj}^*\left( \left( r_{0ij} + \epsilon B_{ij} \right)\left(T_{0j}^1 + \epsilon T_j'^1 \right)\wedge \cdots \wedge \left( T_{0j}^q+ \epsilon T_j'^q \right) \right) \\
&=\left(\textnormal{Id}+ \epsilon p_{kj} \right)\left( r_{0ij}+ \epsilon B_{ij} \right) \cdot \left( r_{0jk} + \epsilon B_{jk} \right)\left(T_{0k}^1 + \epsilon T_k'^1 \right) \wedge \cdots \wedge \left( T_{0k}^p + \epsilon T_k'^p \right) \notag \\
&=\left( r_{0ij} + \epsilon B_{ij} + \epsilon \left[ p_{kj}, r_{0ij} \right] \right)\cdot \left(r_{0jk}+ \epsilon B_{jk} \right)\left(T_{0k}^1 + \epsilon T_k'^1 \right) \wedge \cdots \wedge \left( T_{0k}^p + \epsilon T_k'^p \right) \notag
\end{align}
and we have
\begin{align}\label{ts4}
\theta_{ki}^*\left(\left( T_{0i}^1+ \epsilon T_i'^1 \right)\wedge \cdots \wedge \left( T_{0i}^q + \epsilon T_i'^q \right)  \right) = \left(r_{0jk} + \epsilon B_{ik} \right) \left(T_{0k}^1 + \epsilon T_k'^1 \right) \wedge \cdots \wedge \left( T_{0k}^p + \epsilon T_k'^p \right)
\end{align}
From $(\ref{ts3})$ and $(\ref{ts4})$ we have
\begin{align}\label{si32}
B_{ij} r_{0jk}-\left[ p_{jk}, r_{0ij} \right] r_{0jk}+ r_{0ij} B_{jk}= B_{ik} \Longrightarrow \frac{B_{ij}}{r_{0ij}} - \frac{B_{ik}}{r_{0ik}} + \frac{B_{jk}}{r_{0jk}} - \left[ p_{jk}, \log r_{0ij} \right]=0,
\end{align}
so that $\delta\left( \left\{ \left( p_{ij}, \frac{B_{ij}}{r_{0ij}} \right) \right\} \right)=0$. This implies that $\left( \left\{ T_i'^* \right\}, \left\{ \left( p_{ij}, \frac{B_{ij}}{r_{0ij}} \right)  \right\} \right)\in C^0\left( \mathcal{U}, \mathcal{R} \right) \bigoplus C^1\left( \mathcal{U}, \mathcal{E}_{\mathcal{Q}_0} \right)$ defines a $1$-cocycle in the following \v Cech resolution of $\mathcal{E}_{\mathcal{Q}_0}^\bullet$:
{\Small{\begin{center}
$\begin{CD}
\cdots\\
@AD_3AA \\
C^0\left(\mathcal{U}, \mathscr{H}om_{\mathcal{O}_X}\left( \bigwedge^3 T_{\mathcal{F}_0}, \frac{T_X}{T_{\mathcal{F}_0}}\right) \right)@>\delta>> \cdots\\
@AD_2AA @AD_2AA \\
C^0\left(\mathcal{U}, \mathscr{H}om_{\mathcal{O}_X}\left(\bigwedge^2 T_{\mathcal{F}_0}, \frac{T_X}{T_{\mathcal{F}_0}} \right) \right) @>-\delta>> C^1\left(\mathcal{U}, \mathscr{H}om_{\mathcal{O}_X}\left(  \bigwedge^2 T_{\mathcal{F}_0}, \frac{T_X}{T_{\mathcal{F}_0}} \right) \right) @>\delta>> \cdots\\ 
@A\tilde{D}_1AA @A\tilde{D}_1AA @A\tilde{D}_1AA \\
C^0\left(\mathcal{U}, \mathcal{R} \right)@>\delta>> C^1\left(\mathcal{U}, \mathcal{R} \right) @>-\delta>> C^2\left(\mathcal{U},  \mathcal{R} \right) \\
@A \tilde{D}_0AA @A\tilde{D}_0AA @A\tilde{D}_0AA \\
C^0\left(\mathcal{U}, \mathcal{E}_{\mathcal{Q}_0}\right)@>-\delta>> C^1\left(\mathcal{U}, \mathcal{E}_{\mathcal{Q}_0} \right) @>\delta>> C^2 \left(\mathcal{U}, \mathcal{E}_{\mathcal{Q}_0} \right) @>-\delta>> \cdots
\end{CD}$
\end{center}}}
We show that for two equivalent  first-order foliated deformations of $\left(X, T_{\mathcal{F}_0} \right)$, the associated $1$-cocycles are equivalent. Let $\left(\mathcal{X}', T_{\mathcal{F}'} \right)$ be another first-order foliated deformation of $\left(X, T_{\mathcal{F}_0} \right)$ over $k[\epsilon]$ which induces $\left(\left\{ T_i''^*\right\}, \left\{ \left( p_{ij}' , \frac{B_{ij}'}{r_{0ij}} \right) \right\} \right)$. As in the proof of Proposition \ref{at25}, $\left[ a_i, T_{0i}^\alpha \right] = T_i'^\alpha- T_i''^\alpha- \sum_{\beta=1}^p b_i^{\alpha\beta} T_{0i}^\beta$. Then we have
\begin{align*}
&\sum_{\alpha=1}^p T_{0i}^1 \wedge \cdots \wedge \left[ a_i, T_{0i}^\alpha \right]\wedge \cdots \wedge T_{0i}^p = \sum_{\alpha=1}^p T_{0i}\wedge \cdots\wedge \left( T_i'^\alpha- T_i''^\alpha - \sum_{\beta=1}^p b_i^{\alpha\beta} T_{0i}^\beta \right) \wedge \cdots \wedge T_{0i}^p\\
& \Longrightarrow \left[ a_i , T_{0i}^1 \wedge \cdots \wedge T_{0i}^p \right]+ \left( b_i^{11}+ \cdots + b_i^{pp}\right) T_{0i}^1 \wedge \cdots \wedge T_{0i}^p= T_i'^* - T_i''^*
\end{align*}
and since $\theta_{ji}\alpha_j^{-1}= \alpha_i^{-1} \theta_{ji}'$, we have
\begin{align*}
\left(\alpha_j^{-1}\right)^*\theta_{ji}^*\left(  \left(T_{0i}^1+ \epsilon T_i'^1\right) \wedge \cdots \wedge \left(T_{0i}^p + \epsilon T_i'^p \right) \right) =\left( \alpha_j^{-1}\right)^*\left( \left(r_{0ij} + \epsilon B_{ij} \right)\left(T_{0j}^1+ \epsilon T_j'^1 \right)\wedge \cdots  \wedge \left( T_{0j}^p + \epsilon T_j'^p \right) \right)\\
=\left( r_{0ij}+ \epsilon B_{ij}- \epsilon \left[ a_j, r_{0ij} \right] \right)\cdot \det\left(\delta_\beta^\gamma + \epsilon b_j^{\beta \gamma} \right)\left(T_{0j}^1+ \epsilon T_j''^1 \right)\wedge \cdots  \wedge \left( T_{0j}^p + \epsilon T_j''^p \right) \\
\theta_{ji}'^* \left( \alpha_i^{-1} \right)^* \left(  \left(T_{0i}^1+ \epsilon T_i'^1\right) \wedge \cdots \wedge \left(T_{0i}^p + \epsilon T_i'^p \right)  \right) = \theta_{ji}'^*\left(\det\left(\delta_\alpha^\beta+ \epsilon b_i^{\alpha\beta} \right) \left(T_{0i}^1+ \epsilon T_i''^1  \right)  \wedge \cdots \wedge \left( T_{0i}^p + \epsilon T_i''^p \right) \right)\\
=\det\left(\delta_\alpha^\beta + \epsilon b_i^{\alpha\beta} \right)\left(r_{0ij}+ \epsilon B_{ij}'  \right)  \left(T_{0j}^1+ \epsilon T_j''^1 \right)\wedge \cdots  \wedge \left( T_{0j}^p + \epsilon T_j''^p \right)
\end{align*}
Then we have
\begin{align*}
&B_{ij}- \left[ a_j, r_{0ij} \right] + r_{0ij}\left( b_j^{11}+ \cdots + b_j^{pp} \right)= B_{ij}' + \left(b_i^{11}+ \cdots + b_i^{pp} \right) r_{0ij} \\
&\Longrightarrow \frac{B_{ij}}{r_{0ij}}- \frac{B_{ij}'}{r_{0ij}} = \left( b_i^{11}+ \cdots + b_i^{pp} \right) - \left( \left( b_j^{11}+ \cdots + b_j^{pp} \right)  - \left[ a_j, \log r_{0ij}\right] \right)
\end{align*}
This shows that $\left( \left\{ T_i'^*\right\}, \left\{ \left( p_{ij}, \frac{B_{ij}}{ r_{0ij}} \right) \right\} \right)$ and $\left( \left\{ T_i''^*\right\}, \left\{ \left( p_{ij}', \frac{B_{ij}'}{ r_{0ij}} \right) \right\} \right)$ are cohomologous.

Next we identify obstructions. Let $e: 0 \to (t) \to \tilde{A} \to A \to 0$ be a small extension in $\bold{Art}$ and $\eta:=\left( \mathcal{X}, T_\mathcal{F} \right)$ be a foliated deformation of $\left( X, T_{\mathcal{F}_0} \right)$ over $A$. We keep the notations in the proof of $(2)$ in Proposition \ref{at25}. We set $r_{ij} := \det\left(r_{ij}^{\alpha\beta} \right)$. Since $\theta_{ji}^*\left( T_i^1\wedge \cdots T_i^p \right) = r_{ij} T_j^1 \wedge \cdots \wedge T_j^p$, we see that
\begin{align*}
T_i^1 \wedge \cdots \wedge T_i^p & = \theta_{ij}^* \theta_{jk}^* \theta_{ki}^*\left(T_i^1 \wedge \cdots \wedge T_i^p \right) = \theta_{ij}^*\theta_{jk}^*\left( r_{ik} T_k^1 \wedge \cdots \wedge T_k^p \right) = \theta_{ij}^*\left(f_{jk}\left( r_{ik} \right) r_{kj} T_j^1 \wedge \cdots \wedge T_j^p \right) \\
&= \left( f_{ij} f_{jk} \right)\left( r_{ik} \right) f_{ij}\left(r_{kj} \right) r_{ji} T_j^1 \wedge \cdots \wedge T_j^p
\end{align*}
Then we have
\begin{align}\label{si33}
\left(f_{ij}f_{jk}\right)\left(r_{ik}\right) f_{ij}\left(r_{kj}\right) r_{ji}= 1 &\Longrightarrow \left(\tilde{f}_{ij} \tilde{f}_{jk}\right)\left(\tilde{r}_{ik} \right)\tilde{f}_{ij}\left( \tilde{r}_{kj} \right)\tilde{r}_{ji}= 1 + t g_{ijk}^* \,\,\,\,\,\,\text{for some $g_{ijk}^*\in \Gamma \left(U_{ijk}, \mathcal{O}_X\right)$}
\end{align}
On the other hand, we set $\tilde{r}_{ij} :=\det\left(\tilde{r}_{ij}^{\alpha\beta} \right)$. As in $(\ref{ts8})$, we may assume that $\tilde{f}_{ij}\left(\tilde{r}_{ij} \right)\tilde{r}_{ji}=1$. From $(\ref{ts5})$, we have
\begin{align*}
\tilde{\theta}_{ji}^*\left( \tilde{T}_i^1 \wedge \cdots \wedge \tilde{T}_i^p     \right) &= \left( \sum_{\beta=1}^p \tilde{r}_{ij}^{1\beta} \tilde{T}_j^\beta + t W_{ij}^1 \right) \wedge \cdots \wedge \left( \sum_{\beta=1}^p \tilde{r}_{ij}^{p\beta} \tilde{T}_j^\beta + t W_{ij}^p  \right) \\
&= \tilde{r}_{ij} \tilde{T}_j^1 \wedge \cdots \tilde{T}_j^p + t \sum_{\alpha=1}^p T_{0i}^1 \wedge \cdots \wedge W_{ij}^\alpha \wedge \cdots \wedge T_{0i}^p
\end{align*}
We define $W_{ij}^*\in \Gamma\left( U_{ij}, \mathcal{R} \right)$ by
\begin{align*}
W_{ij}^*:= \sum_{\alpha=1}^p T_{0i}^1 \wedge \cdots W_{ij}^\alpha \wedge \cdots \wedge T_{0i}^p \in \Gamma\left( U_{ij}, \mathcal{R} \right) \subset \Gamma\left( U_{ij} , \bigwedge^p T_X \otimes \mathcal{Q}_0 \right)
\end{align*}
Then from $(\ref{atb1})$, we have
\begin{align}\label{si35}
W_{ij}^*= \sum_{\alpha=1}^p T_{0i}^1 \wedge \cdots W_{ij}^\alpha \wedge \cdots \wedge T_{0i}^p = - r_{0ij}\sum_{\alpha=1}^p T_{0j}^1 \wedge \cdots \wedge W_{ji}^\alpha \wedge \cdots \wedge T_{0j}^p = - r_{0ij}W_{ji}^*,
\end{align}
so that $\left\{ W_{ij}^*\right\}\in C^1\left( \mathcal{U}, \mathcal{R} \right)$. Then we claim that $\left( \left\{ P_i \right\}, \left\{ - W_{ij}^*\right\}, \left\{ \left( d_{ijk}, - g_{ijk}^* \right) \right\} \right)\in C^0\left(\mathcal{U}, \mathscr{H}om_{\mathcal{O}_X}\left(\bigwedge^2 T_{\mathcal{F}_0}, \frac{T_X}{ T_{\mathcal{F}_0}} \right) \right) $ $ \bigoplus C^1\left( \mathcal{U}, \mathcal{R}\right) \bigoplus C^2\left( \mathcal{U}, \mathcal{E}_{\mathcal{Q}_0} \right)$ defines a $2$-cocycle in the above \v Cech resolution of $\mathcal{E}_{\mathcal{Q}_0}^\bullet$. Here we consider $\left( d_{ijk}, - g_{ijk}^*\right)\in T_X|_{U_{ijk}} \bigoplus \mathcal{O}_X|_{U_{ijk}}\subset T_X|_{U_i} \bigoplus \mathcal{O}_X|_{U_i}$ as an element of $\Gamma\left( U_{ijk}, \mathcal{E}_{\mathcal{Q}_0} \right)$. It is sufficient to show that $-\delta\left( \left\{ - W_{ij}^* \right\} \right) +\tilde{D}_0\left( \left\{ \left( d_{ijk} , - g_{ijk}^*  \right) \right\} \right)=0 $ and $-\delta\left( \left\{ \left( d_{ijk} , - g_{ijk}^* \right) \right\} \right)=0$. In fact, we note that
\begin{align*}
\tilde{\theta}_{ij}^*\tilde{\theta}_{jk} \tilde{\theta}_{ki}^*\left(  \tilde{T}_i^1 \wedge \cdots \wedge \tilde{T}_i^p \right) & = \left( \tilde{T}_i^1 + t \left[ d_{ijk} , T_{0i}^1\right] \right) \wedge \cdots \wedge \left( \tilde{T}_i^p + t \left[ d_{ijk} , T_{0i}^p \right] \right) \\
&= \tilde{T}_i^1\wedge \cdots \wedge \tilde{T}_i^p + t\sum_{\alpha=1}^p T_{0i}^1 \wedge \cdots \wedge \left[ d_{ijk} , T_{0i}^\alpha \right] \wedge \cdots \wedge T_{0i}^p
\end{align*}
On the other hand,
\begin{align*}
\tilde{\theta}_{ij}^*\tilde{\theta}_{jk} \tilde{\theta}_{ki}^* \left( \tilde{T}_i^1 \wedge \cdots \wedge \tilde{T}_i^p \right) & = \tilde{\theta}_{ij}^*\tilde{\theta}_{jk}^* \left( \tilde{r}_{ik} \tilde{T}_k^1 \wedge \cdots \wedge \tilde{T}_k^p + t W_{ik}^* \right) = \tilde{\theta}_{ij}^*\left( \tilde{f}_{jk}\left(\tilde{r}_{ik} \right) \tilde{r}_{kj} \tilde{T}_j^1 \wedge \cdots \wedge \tilde{T}_j^p  + t r_{0ik} W_{kj}^* + t W_{ik}^* \right) \\
&= \tilde{f}_{ij}\tilde{f}_{jk}\left( \tilde{r}_{ik}\right)\tilde{r}_{kj} \tilde{r}_{ji}\tilde{T}_i^1 \wedge \cdots \tilde{T}_i^p + t r_{0ij} W_{ji}^* +t r_{0ik} W_{kj}^* + t W_{ik}^* \\
&=\left( 1+ t g_{ijk}^* \right)\tilde{T}_i^1\wedge \cdots \wedge \tilde{T}_i^p + t r_{0ij} W_{ji}^* +t r_{0ik} W_{kj}^* + t W_{ik}^*
\end{align*}
Then we have
\begin{align*}
\left[ d_{ijk}, T_{0i}^1\wedge \cdots \wedge T_{0i}^p \right]- g_{ijk}^* T_{0i}^1 \wedge \cdots \wedge T_{0i}^p = - W_{ij}^* - r_{0ij} W_{jk}^* + W_{ik}^*
\end{align*}
This implies that $-\delta\left( \left\{ - W_{ij}^* \right\} \right) +\tilde{D}_0\left( \left\{ \left( d_{ijk} , - g_{ijk}^*  \right) \right\} \right)=0$. Next we note that
\begin{align}\label{ts7}
&\left( \tilde{f}_{ij}\tilde{f}_{jk} \right)\left( \tilde{r}_{ik} \right)\tilde{f}_{ij}\left(\tilde{r}_{kj} \right)\tilde{r}_{ji}= 1+ t g_{ijk}^*  \\
&\iff \tilde{f}_{jk}\left( \tilde{r}_{ik} \right) \tilde{r}_{kj} \cdot \tilde{f}_{ji}\left( \tilde{r}_{ji} \right) = 1 + t g_{ijk}^*  \notag \\
& \iff \tilde{r}_{ik} \tilde{f}_{kj}\left( \tilde{r}_{kj} \right) \tilde{f}_{kj} \tilde{f}_{ji} \left( \tilde{r}_{ji} \right) = 1 + t g_{ijk}^* \notag
\end{align}
Then we have from $(\ref{ts7})$
\begin{align*}
&\left(1+ tg_{ilk}^* \right)\left(1+ t g_{ijl}^*\right)\left(1+ g_{jkl}^*\right)\left( 1+ g_{ikj}^*\right) \\
&= \tilde{f}_{lk}\left(\tilde{r}_{ik} \right)\cdot \tilde{r}_{kl}  \cdot \tilde{f}_{li}\left( \tilde{r}_{li}\right)  \cdot \tilde{r}_{il} \cdot \tilde{f}_{lj}\left(\tilde{r}_{lj} \right) \cdot \tilde{f}_{lj} \tilde{f}_{ji}\left( \tilde{r}_{ji}\right) \cdot \tilde{r}_{jl} \cdot \tilde{f}_{lk}\left( \tilde{r}_{lk} \right) \tilde{f}_{lk} \tilde{f}_{kj}\left( \tilde{r}_{kj} \right) \cdot \tilde{f}_{lj}\left( \tilde{r}_{ij}\right) \cdot \tilde{f}_{lj}\tilde{f}_{jk}\left( \tilde{r}_{jk} \right) \cdot \tilde{f}_{lj} \tilde{f}_{jk} \tilde{f}_{ki}\left(\tilde{r}_{ki}\right)\\
&= \tilde{f}_{lk}\left( \tilde{r}_{ik} \right) \tilde{f}_{lk} \tilde{f}_{kj} \left( \tilde{r}_{kj} \right) \tilde{f}_{lj} \tilde{f}_{jk}\left( \tilde{r}_{jk} \right) \cdot \tilde{f}_{lj} \tilde{f}_{jk} \tilde{f}_{ki}\left( \tilde{r}_{ki} \right)\\
&= \tilde{f}_{lk}\left( \tilde{r}_{ik} \right) \left( \tilde{f}_{lj} + t d_{lkj}  \right) \left( \tilde{r}_{kj}  \right) \tilde{f}_{lj} \tilde{f}_{jk} \left( \tilde{r}_{jk} \right) \cdot \left( \tilde{f}_{lk} + t \tilde{d}_{ljk} \right) \tilde{f}_{ki} \left(  \tilde{r}_{ki} \right) \\
&= 1 + t \left[d_{lkj},  \log r_{0kj} \right] + t \left[ d_{ljk},  \log r_{0ki}\right] = 1 + t \left[ d_{jkl}, \log r_{0jk}+ \log r_{0ki} \right]= 1 - t\left[ d_{jkl}, \log r_{0ij} \right]
\end{align*}
Then we have
\begin{align*}
&g_{ilk}^*+ g_{ijl}^* + g_{jkl}^* + g_{ikj}^* + \left[ d_{jkl}, \log r_{0ij} \right]=0\\
&\iff -g_{ikl}^*+ g_{ijl}^* + g_{jkl}^* - g_{ijk}^* + \left[ d_{jkl}, \log r_{0ij} \right]=0\\
&\iff g_{ikl}^* - g_{ijl}^* - g_{jkl}^* + g_{ijk}^* - \left[ d_{jkl}, \log r_{0ij} \right]=0
\end{align*}
This implies that $-\delta\left( \left\{ \left( d_{ijk} , - g_{ijk}^* \right) \right\} \right)=0$. This proves the claim. We show that for another lifting $\tilde{\theta}_{ij}'$ of $\theta_{ij}$, and another lifting $\tilde{T}_i'^\alpha$ of $T_i^\alpha$, and another lifting $\tilde{r}_{ij}'^{\alpha\beta}$ of $r_{ij}^{\alpha\beta}$, and another lifting of $\tilde{g}_{i\alpha\beta}'^\gamma$ of $g_{i\alpha\beta}^\gamma$, the associated $2$-cocycle $\beta:=\left( \left\{P_i' \right\},  \left\{ - W_{ij}'^* \right\}, \left\{ \left( d_{ijk}', - g_{ijk}'^*\right)\right\} \right)$ is cohomologous to $\alpha:=\left(\left\{P_i \right\},  \left\{ - W_{ij}^* \right\}, \left\{ \left( d_{ijk}, - g_{ijk}^*\right)\right\} \right)$. We define $Y_i^*:=\sum_{\alpha=1}^p T_{0i}^1 \wedge \cdots \wedge Y_i^\alpha\wedge \cdots \wedge T_{0i}^p \in \Gamma\left( U_i, \mathcal{R} \right)$. We note that 
\begin{align}\label{ts9}
\tilde{r}_{ij}'^{\alpha\beta}= \tilde{r}_{ij}^{\alpha\beta} + t s_{ij}^{\alpha\beta} &\Longrightarrow \sum_{\beta=1}^p \tilde{r}_{ij}'^{\alpha\beta} \tilde{f}_{ji}\left(\tilde{r}_{ji}^{\beta \gamma} \right) = \delta_\alpha^\gamma + t \sum_{\beta=1}^p s_{ij}^{\alpha\beta} r_{0ji}^{\beta \gamma} \Longrightarrow \tilde{r}_{ij}' \tilde{f}_{ji}\left(\tilde{r}_{ji} \right)= 1 + t \sum_{\alpha, \beta=1}^p s_{ij}^{\alpha\beta} r_{0ji}^{\beta \alpha} \notag \\
&\Longrightarrow \tilde{r}_{ij}'= \tilde{r}_{ij}+ t r_{0ij}\sum_{\alpha, \beta=1}^p s_{ij}^{\alpha\beta} r_{0ij}^{\beta \alpha}:= \tilde{r}_{ij} + t r_{0ij} A_{ij} ,\,\,\,\,\,\,\ \textnormal{where}\,\,\,A_{ij}:=\sum_{\alpha,\beta=1}^p  s_{ij}^{\alpha\beta} r_{0ji}^{\beta \alpha}
\end{align}
and we also note that
\begin{align}\label{si36}
\tilde{r}_{ji}+ t r_{0ji} A_{ji}= \tilde{r}_{ji}' & \Longrightarrow \tilde{r}_{ji} \tilde{f}_{ij}'\left(\tilde{r}_{ij}' \right) + t A_{ji}= 1 \Longrightarrow  \tilde{r}_{ji} \left(\tilde{f}_{ij} + t p_{ij} \right)\left( \tilde{r}_{ij} + t r_{0ij} A_{ij} \right) + t A_{ji} =1 \\
& \Longrightarrow  1 +  r_{0ji}\left[ p_{ij}, r_{0ij}\right] + t A_{ij} + tA_{ji}=1 \Longrightarrow - A_{ji}= A_{ij}-\left[ p_{ij}, \log r_{0ji}\right] \notag
\end{align}
This implies that $\left( \left\{ p_{ij}, A_{ij} \right\} \right)\in C^1\left( \mathcal{U}, \mathcal{E}_{\mathcal{Q}_0} \right)$. Then from $(\ref{ts9})$ we have
\begin{align} \label{si37}
1+ t g_{ijk}'^*&= \tilde{f}_{ij}' \tilde{f}_{jk}'\left( \tilde{r}_{ik}' \right)\tilde{f}_{ij}'\left( \tilde{r}_{kj}' \right) \tilde{r}_{ji}'  \\
&=\left(\tilde{f}_{ij} + t p_{ij} \right)\left(\tilde{f}_{jk} + tp_{jk} \right)\left( \tilde{r}_{ik}+ t r_{0ik} A_{ik} \right) \cdot \left( \tilde{f}_{ij} + t p_{ij} \right)\left( \tilde{r}_{kj} + t r_{0kj} A_{kj} \right) \cdot \left( \tilde{r}_{ji} + t r_{0ji} A_{ji} \right) \notag \\
&= 1+ t g_{ijk}^* + t \left[ p_{ij}, r_{0ik}\right] r_{0ki} + t \left[ p_{jk},  r_{0ik}\right] r_{0ki} + \left[ p_{ij}, r_{0kj}\right] r_{0jk} + tA_{ik}+ t A_{kj}+ tA_{ji} \notag \\
&=1 +t g_{ijk}^* - t \left[ p_{ji}, \log r_{0ij}\right] - t \left[ p_{kj}, \log r_{0ik}\right] +t A_{ik}+ t A_{kj} + t A_{ji} \notag \\
&= 1+  tg_{ijk}^* - tA_{ij}+ tA_{ik}- tA_{jk}^i \notag
\end{align} 
This implies that $\delta\left(\left\{ p_{ij}, A_{ij} \right\} \right) =\left( \left\{\left(d_{ijk}'-d_{ijk}, - g_{ijk}'^* - \left( - g_{ijk}^* \right) \right) \right\} \right)$. On the other hand, from $(\ref{ts10})$, we have
\begin{align*}
- W_{ij}'^* + W_{ij}^* &= \left[ p_{ij}, T_{0i}^1\wedge \cdots \wedge T_{0i}^p \right] + r_{0ij} Y_j^* - Y_i^* + \sum_{\alpha=1}^p T_{0i}^1 \wedge \cdots \wedge \left(\sum_{\beta,\gamma=1}^p s_{ij}^{\alpha\beta} r_{0ji}^{\beta \gamma} T_{0i}^\gamma \right) \wedge \cdots \wedge T_{0i}^p\\
&=  \left[ p_{ij}, T_{0i}^1\wedge \cdots \wedge T_{0i}^p \right] + r_{0ij} Y_j^* - Y_i^*  + A_{ij} T_{0i}^1 \wedge \cdots \wedge T_{0i}^p
\end{align*}
This implies that $\delta\left(\left\{ Y_i^*\right\} \right) + \tilde{D}_0\left( \left\{ \left( p_{ij} , A_{ij} \right) \right\} \right)=\left\{ -W_{ij}'^*+ W_{ij}^*\right\}$. This proves that $\alpha$ is cohomologous to $\beta$. So given a small extension $e: 0 \to (t) \to \tilde{A} \to A \to 0$ and an infinitesimal foliated deformation $\eta$ of $\left( X, T_{\mathcal{F}_0} \right)$ over $A$,  we can associate an element $o_\eta(e):=$ the cohomology class of $\alpha\in \mathbb{H}^2\left( X, \mathcal{E}_{\mathcal{Q}_0}^\bullet \right)$. We note that $o_\eta(e)=0$ if and only if there exists $\left( \left\{ X_i^*\right\}, \left\{ d_{ij}, V_{ij} \right\} \right)\in C^0\left( \mathcal{U}, \mathcal{R} \right)\bigoplus C^1\left( \mathcal{U}, \mathcal{E}_{\mathcal{Q}_0} \right)$ where $X_i^*=\sum_{\alpha=1}^p T_{0i}^1 \wedge \cdots \wedge X_i^\alpha \wedge \cdots \wedge T_{0i}^p$ such that $\tilde{D}_1\left( X_i^* \right)= P_i$ and $X_j^* - X_i^*+ \tilde{D}_0\left( \left\{ d_{ij} , V_{ij} \right\} \right)= - W_{ij}^*$ and $d_{ij}- d_{ik} + d_{jk}=d_{ijk}$ and $V_{ij}- V_{ik}+ V_{jk}^i= - g_{ijk}^*$. Then
\begin{align*}
&t \sum_{\alpha=1}^p T_{0i}^1 \wedge \cdots \wedge \left(  \left[ d_{ij}, T_{0i}^\alpha \right] + \sum_{\beta=1}^p r_{0ij}^{\alpha\beta} X_j^\beta - X_i^\alpha          \right) \wedge \cdots \wedge T_{0i}^p + t V_{ij} T_{0i}^1 \wedge \cdots \wedge T_{0i}^p  = - \tilde{\theta}_{ji}^*\left( \tilde{T}_i^1 \wedge \cdots \wedge \tilde{T}_i^p \right) + \tilde{r}_{ij} \tilde{T}_j^1 \wedge \cdots \wedge \tilde{T}_j^p \\
&\Longrightarrow  \Phi_{ji}^* \left(\left( \tilde{T}_i^1 - t X_i^1 \right) \wedge \cdots \wedge \left( \tilde{T}_i^p - t X_i^p \right) \right)      = \left( \tilde{r}_{ij} - t r_{0ij}V_{ij} \right)  \left( \tilde{T}_j^1 - t X_j^1 \right) \wedge \cdots \wedge \left( \tilde{T}_j^p - t X_j^p \right)
\end{align*}
and
\begin{align}
&t V_{ij}- t V_{ik} + t V_{jk} - t\left[ d_{jk}, \log r_{0ij} \right] = -t g_{ijk}^*= 1 - \left(\tilde{f}_{ij} \tilde{f}_{jk} \right)\left( \tilde{r}_{ik} \right) \tilde{f}_{ij}\left( \tilde{r}_{kj} \right)\tilde{r}_{ji} \notag \\
&\Longrightarrow -t V_{ji} + \left[ d_{ij}, \log r_{0ji}\right] - t V_{ik} -tV_{kj} + t\left[ d_{jk}, \log r_{0kj}\right]- t \left[ d_{jk}, \log r_{0ij}\right] = 1-  \left(\tilde{f}_{ij} \tilde{f}_{jk} \right)\left( \tilde{r}_{ik} \right) \tilde{f}_{ij}\left( \tilde{r}_{kj} \right)\tilde{r}_{ji} \notag \\
&\Longrightarrow \left( \tilde{f}_{ij} - t d_{ij}\right)\left( \tilde{f}_{jk} - t d_{jk}\right)\left(\tilde{r}_{ik}- t r_{0ik} V_{ik} \right) \cdot \left( \tilde{f}_{ij} - t d_{ij} \right)\left(\tilde{r}_{kj} -t r_{0kj} V_{kj} \right)\cdot \left(\tilde{r}_{ji} - t r_{0ji} V_{ji} \right) = 1 \label{si40}
\end{align}
Hence a lifting $\tilde{\eta}$ of $\eta$ exists. This completes the proof of Proposition \ref{ts11}.

\end{proof}

\subsection{Relationship between foliated deformations of $\left(X, T_{\mathcal{F}_0} \right)$ in terms of vector bundles and in terms of deformations of sections (in terms of line bundles)} \

Let $\left( X, T_{\mathcal{F}_0} \right)$ be a foliated nonsingular variety with $T_{\mathcal{F}_0}$ locally free. We will show that the following diagram commutes and is exact.

\begin{center}
$\begin{CD}
@. \cdots @. \cdots @. \cdots \\
@. @AAA @AD_3AA @AD_3AA \\
0@>>>0 @>>> \mathscr{H}om_{\mathcal{O}_X}\left( \bigwedge^3 T_{\mathcal{F}_0}, \frac{T_X}{T_{\mathcal{F}_0}}\right) @>>> \mathscr{H}om_{\mathcal{O}_X}\left( \bigwedge^3 T_{\mathcal{F}_0}, \frac{T_X}{T_{\mathcal{F}_0}} \right) @>>> 0 \\
@. @AAA @AD_2AA @AD_2AA \\
0@>>> 0 @>>> \mathscr{H}om_{\mathcal{O}_X}\left( \bigwedge^2 T_{\mathcal{F}_0}, \frac{T_X}{T_{\mathcal{F}_0}} \right) @>>> \mathscr{H}om_{\mathcal{O}_X}\left( \bigwedge^2 T_{\mathcal{F}_0}, \frac{T_X}{T_{\mathcal{F}_0}} \right) @>>> 0  \\
@. @AAA @AD_1' AA @A\tilde{D}_1AA \\
0@>>> \mathscr{H}om_{\mathcal{O}_X}\left(T_{\mathcal{F}_0}, T_{\mathcal{F}_0} \right)|_{\textnormal{Tr}=0} @>>>  \mathscr{H}om_{\mathcal{O}_X}\left( T_{\mathcal{F}_0}, T_X \right) @>\beta_1>> \mathcal{R}\subset \bigwedge^p T_X\otimes \mathcal{Q}_0 @>>> 0 \\
@. @A-\textnormal{Id}AA @AD_0'AA @A\tilde{D}_0AA \\
0 @>>> \mathscr{H}om_{\mathcal{O}_X}\left( T_{\mathcal{F}_0}, T_{\mathcal{F}_0} \right)|_{\textnormal{Tr}=0} @>>> \mathcal{E}_{T_{\mathcal{F}_0}} @>\beta_0>> \mathcal{E}_{\mathcal{Q}_0} @>>> 0
\end{CD}$
\end{center}

We define the map $\beta_1$ in the following way:
\begin{align*}
\beta_1:\mathscr{H}om_{\mathcal{O}_X}\left( T_{\mathcal{F}_0}, T_X \right) &\to \mathcal{R}\subset \bigwedge^p T_X\otimes \mathcal{Q}_0 \\
\phi &\mapsto \phi_i:=\sum_{\alpha=1}^p T_{0i}^1 \wedge \cdots \wedge \phi \left(T_{0i}^\alpha \right) \wedge \cdots \wedge T_{0i}^p\,\,\,\,\,\,\,\,\,\,\textnormal{on}\,\,\,U_i
\end{align*}
This is well-define since
\begin{align*}
\phi_i=\sum_{\alpha=1}^p T_{0i}^1\wedge \cdots \wedge \phi\left(T_{0i}^\alpha\right)\wedge \cdots \wedge T_{0i}^p = \det\left(r_{0ij}^{\alpha\beta}\right) \sum_{\alpha=1}^p T_{0j}^1\wedge \cdots \wedge \phi \left(T_{0j}^\alpha\right)\wedge \cdots \wedge T_{0j}^p= r_{0ij}  \phi_j
\end{align*}
This is clearly surjective. Let us find the kernel. Assume that $\beta_1\left(\phi\right)=0$. Then 
\begin{align*}
\sum_{\alpha=1}^p T_{0i}^1\wedge \cdots \wedge \phi\left(T_{0i}^\alpha\right) \wedge \cdots \wedge T_{0i}^p =0
\end{align*}
By multiplying $T_{0i}^\alpha$, we have
\begin{align*}
T_{0i}^1\wedge \cdots \wedge T_{0i}^p \wedge \phi\left(T_{0i}^\alpha\right)=0
\end{align*}
Then $\phi \left(T_{0i}^\alpha\right) \in \Gamma(U_i-S, T_{\mathcal{F}_0})$, so that by the reflexiveness of $T_{\mathcal{F}_0}$, we have $\phi\left(T_{0i}^\alpha\right)\in \Gamma\left(U_i, T_{\mathcal{F}_0}\right)$ and so $\phi \left(T_{0i}^\alpha\right)=\sum_{\beta=1}^p f_i^{\alpha\beta} T_{0i}^\beta$ for some $f_i^{\alpha\beta}\in \Gamma(U_i, \mathcal{O}_X)$. Then
\begin{align*}
0= \sum_{\alpha=1}^p T_{0i}^1\wedge \cdots \wedge \phi\left(T_{0i}^\alpha\right)\wedge \cdots \wedge T_{0i}^p = \sum_{\alpha=1}^p T_{0i}^1\wedge \cdots \wedge f_i^{\alpha\alpha} T_{0i}^\alpha \wedge \cdots \wedge T_{0i}^p = \left( f_i^{11}+ \cdots + f_i^{pp}\right) T_{0i}^1\wedge \cdots \wedge T_{0i}^p
\end{align*}
so that $\sum_{\alpha=1}^p f_i^{\alpha\alpha}=0$. This implies that $\phi\in \mathscr{H}om_{\mathcal{O}_X}\left( T_{\mathcal{F}_0}, T_{\mathcal{F}_0} \right)|_{\textnormal{Tr}=0}$.

On the other hand, we define $\beta_0:\mathcal{E}_{T_{\mathcal{F}_0}}\to \mathcal{E}_{\mathcal{Q}_0}$ locally in the following way:
\begin{align*}
\beta_0: \mathcal{E}_{T_{\mathcal{F}_0}}|_{U_i}\cong \mathscr{H}om_{\mathcal{O}_X}\left( T_{\mathcal{F}_0}, T_{\mathcal{F}_0}  \right)|_{U_i} \bigoplus T_X |_{U_i} &\to \mathcal{E}_{\mathcal{Q}_0} |_{U_i} \cong \mathcal{O}_X|_{U_i}\bigoplus T_X|_{U_i}\\
                                \left(\phi_i, X \right) &\mapsto  \left(-\textnormal{Tr}(\phi_i), X \right)
\end{align*}
Then the kernel of $\beta_0$ on $U_i$ is $\mathscr{H}om_{\mathcal{O}_X}\left( T_{\mathcal{F}_0}, T_{\mathcal{F}_0} \right)|_{\textnormal{Tr}=0}|_{U_i}$. Next we show that $\beta_0$ is well-defined, i.e. the following diagram commutes
\begin{center}
$\begin{CD}
\left(\phi_j, X\right)\in \mathscr{H}om_{\mathcal{O}_X}\left(T_{\mathcal{F}_0}, T_{\mathcal{F}_0} \right)|_{U_j} \bigoplus T_X|_{U_j} @>\beta_0>> \left(-\textnormal{Tr}\left(\phi_j \right)= - \textnormal{Tr}\left(\phi_i\right)+\left[ X, \log r_{0ij}\right], X \right) \in  \mathcal{O}_X|_{U_j} \bigoplus T_X|_{U_j} \\
@AAA @AAA \\
\left(\phi_i, X \right) \in \mathscr{H}om_{\mathcal{O}_X}\left(T_{\mathcal{F}_0}, T_{\mathcal{F}_0} \right)|_{U_i} \bigoplus T_X|_{U_i} @>\beta_0>> \left(- \textnormal{Tr}\left(\phi_i \right), X \right) \in \mathcal{O}_X|_{U_i} \bigoplus T_X|_{U_i}
\end{CD}$
\end{center}
where $\phi_i\left(T_{0i}^\alpha\right)=\phi_j\left(T_{0i}^\alpha\right)+\sum_{\beta=1}^p \left[X, r_{0ij}^{\alpha\beta}\right] T_{0j}^\beta=\phi_j\left(T_{0i}^\alpha\right)+ \left[ X, T_{0i}^\alpha \right]-\sum_{\beta=1}^p r_{0ij}^{\alpha\beta} \left[X, T_{0j}^\beta \right] $. We have to show that $-\textnormal{Tr}\left(\phi_j\right)= - \textnormal{Tr}\left(\phi_i \right)+\left[X, \log r_{0ij} \right]$.  In fact, we note that
{\small{\begin{align*}
&\textnormal{Tr}\left(\phi_i\right)T_{0i}^1\wedge \cdots \wedge T_{0i}^p=\sum_{\alpha=1}^pT_{0i}^1\wedge \cdots \wedge \phi_i\left(T_{0i}^\alpha\right)\wedge \cdots \wedge T_{0i}^p =\sum_{\alpha=1}^p T_{0i}^1\wedge \cdots \wedge \left( \phi_j \left(T_{0i}^\alpha \right) +\left[X, T_{0i}^\alpha \right]  - \sum_{\beta=1}^p r_{0ij}^{\alpha\beta} \left[ X, T_{0j}^\beta \right]  \right) \wedge \cdots \wedge T_{0i}^p\\
&= r_{0ij} \textnormal{Tr} \left(\phi_j \right) T_{0j}^1\wedge \cdots \wedge T_{0j}^p +\left[ X, T_{0i}^1\wedge \cdots \wedge T_{0i}^p   \right] - r_{0ij} \left[ X, T_{0j}^1\wedge \cdots \wedge T_{0j}^p \right]= \textnormal{Tr} \left(\phi_j \right) T_{0i}^1\wedge \cdots \wedge T_{0i}^p + \left[X, r_{0ij}\right] T_{0j}^1\wedge \cdots \wedge T_{0j}^p
\end{align*}}}
Then we see that $\textnormal{Tr}\left(\phi_i \right)= \textnormal{Tr}\left(\phi_j \right) + \left[ X, \log r_{0ij} \right]$. Next we show that the following diagram commutes
\begin{center}
$\begin{CD}
\mathscr{H}om_{\mathcal{O}_X}\left( T_{\mathcal{F}_0}, T_X \right) @>\beta_1>> \mathcal{R} \\
@AD_0' AA @AA\tilde{D}_0A \\
\mathcal{E}_{T_{\mathcal{F}_0}} @>\beta_0>> \mathcal{E}_{\mathcal{Q}_0} 
\end{CD}$
\end{center}
In fact, on $U_i$, $D_0' \left(\phi_i, X\right) \left(T_{0i}^\alpha \right)= \left(L_X-\phi_i \right)\left(T_{0i}^\alpha \right)= \left[X, T_{0i}^\alpha \right]- \phi_i\left(T_{0i}^\alpha\right)$ for $\left(\phi_i, X \right)\in \Gamma\left( U_i, \mathcal{E}_{T_{\mathcal{F}_0}}\right)$. Then we have
\begin{align*}
\beta_1D_0' \left(\phi_i, X \right)&=\sum_{\alpha=1}^p T_{0i}^1\wedge \cdots  \wedge \left( \left[X, T_{0i}^\alpha\right] -\phi_i \left(T_{0i}^\alpha \right) \right) \wedge \cdots \wedge T_{0i}^p= \left[ X, T_{0i}^1\wedge \cdots \wedge T_{0i}^p\right] - \textnormal{Tr}\left(\phi_i\right) T_{0i}^1\wedge \cdots \wedge T_{0i}^p\\
 &= \tilde{D}_0\left( - \textnormal{Tr}\left(\phi_i \right), X \right) = \tilde{D}_0 \beta_0\left(\phi_i, X \right)
\end{align*}
and by definition $D_1'= \tilde{D}_1\circ \beta_1$, and for $\phi \in \mathscr{H}om_{\mathcal{O}_X}\left(T_{\mathcal{F}_0}, T_{\mathcal{F}_0} \right)|_{\textnormal{Tr}=0}$, we have $D_0' \left(\phi,  0 \right)=-\phi$. Hence the above diagram is exact. Then we have 
\begin{align*}
\mathbb{H}^i\left( X, \mathcal{E}_{T_{\mathcal{F}_0}}^\bullet \right) \cong \mathbb{H}^i\left( X, \mathcal{E}_{\mathcal{Q}_0}^\bullet \right)
\end{align*}

Then the map $\mathcal{E}_{T_{\mathcal{F}_0}}^\bullet \xrightarrow{\beta} \mathcal{E}_{\mathcal{Q}_0}^\bullet$ above gives an isomorphism of first-order deformations and obstructions of foliated deformations of $(X, T_{\mathcal{F}_0})$ in terms of tangent sheaf between interpretations in terms of the vector bundle and in terms of the line bundle. In fact, we keep the notations in the proof of Proposition \ref{at27} and Proposition \ref{ts11}. Then we have to show that
\begin{align*}
\mathbb{H}^1\left( X, \mathcal{E}_{T_{\mathcal{F}_0}}^\bullet \right) &\xrightarrow{\beta} \mathbb{H}^1\left( X, \mathcal{E}_{\mathcal{Q}_0}^\bullet \right) \\
\left( \left\{ T_i' \right\}, \left\{ \left( p_{ij}, - b_{ij} \right) \right\} \right) &\mapsto \left( \left\{ T_i'^* \right\}, \left\{ \left( p_{ij},  \frac{B_{ij}}{r_{0ij}} \right) \right\} \right)
\end{align*}
In fact, it is clear that $T_i'\mapsto T_i'^*$ under $\beta_1$. We note that if we set two matrices $R_{0ij} :=\left( r_{0ij}^{\alpha\beta} \right)$ and $g_{ij} := \left(  b_{ij}^{\alpha\beta} \right)$, then
\begin{align}\label{si43}
r_{0ij} + \epsilon B_{ij} & = \det\left( r_{0ij}^{\alpha\beta} +\epsilon b_{ij}^{\alpha\beta} \right)  = \det \left( R_{0ij} \right) \det\left(\textnormal{Id} + \epsilon \left(R_{0ij}^{-1} b_{ij}\right) \right)  = r_{0ij} \left(1 + \epsilon \textnormal{Tr}\left( R_{0ij}^{-1} g_{ij} \right) \right) \\
& \Longrightarrow \frac{B_{ij}}{r_{0ij}}=  \textnormal{Tr}\left( g_{ij} R_{0ij}^{-1} \right) = \sum_{\alpha,\beta=1}^p  b_{ij}^{\alpha \beta} r_{0ji}^{\beta \alpha }  = \textnormal{Tr}\left( b_{ij} \right) =  - \textnormal{Tr} \left(   - b_{ij}  \right) \notag
\end{align}

Next we show that $\beta:\mathcal{E}_{T_{\mathcal{F}_0}}^\bullet \xrightarrow{\beta} \mathcal{E}_{\mathcal{Q}_0}^\bullet$ gives  the obstruction map for a small extension $e: 0\to (t)\to \tilde{A} \to A \to 0$. We have to show that 
\begin{align*}
\mathbb{H}^2\left( X,  \mathcal{E}_{T_{\mathcal{F}_0}}^\bullet \right) &\xrightarrow{\beta} \mathbb{H}^2\left( X, \mathcal{E}_{\mathcal{Q}_0}^\bullet  \right) \\
\left( \left\{ P_i\right\}, \left\{ - W_{ij} \right\}, \left\{ \left( d_{ijk}, g_{ijk} \right) \right\} \right) &\mapsto \left( \left\{ P_i \right\}, \left\{ - W_{ij}^* \right\}, \left\{ \left( d_{ijk},  - g_{ijk}^* \right) \right\} \right)
\end{align*}
In fact, it is clear that $- W_{ij} \mapsto - W_{ij}^*$ under $\beta_1$. We recall that $\sum_{\beta, \gamma=1}^p \left( \tilde{f}_{ij} \tilde{f}_{jk}  \right)\left( \tilde{r}_{ik}^{\alpha\beta} \right) \tilde{f}_{ij}\left( \tilde{r}_{kj}^{\beta \gamma} \right) \tilde{r}_{ji}^{\gamma \eta}  = \delta_\alpha^\eta +t g_{ijk}^{\alpha\eta}$. Let $G_{ijk}$ be the matrix of $\left( g_{ijk}^{\alpha\eta} \right)$. Then
\begin{align}\label{si45}
1+ t \textnormal{Tr}\left( G_{ijk} \right)=\det\left(\textnormal{Id} + t G_{ijk} \right) = \left( \tilde{f}_{ij} \tilde{f}_{jk} \right)\left( \tilde{r}_{ik} \right) \tilde{f}_{ij}\left(\tilde{r}_{kj} \right) \tilde{r}_{ji}= 1 + t g_{ijk}^* \Longrightarrow  - \textnormal{Tr}\left( g_{ijk} \right) = - g_{ijk}^*
\end{align}

We also remark  the case of foliated deformations of $(X, T_{\mathcal{F}_0})$ with $X$ fixed (i.e. without deformations of complex structures of $X$). We note that we have shown that $\textnormal{FDef}_{T_{\mathcal{F}_0}}$ is controlled by $\mathscr{H}om_{\mathcal{O}_X} \left( T_{\mathcal{F}_0},  T_{\mathcal{F}_0} \right)^\bullet$ from $(\ref{tf1})$. We consider foliated deformations of $T_{\mathcal{F}_0}$ with $X$ fixed in terms of deformations of section or in terms of the line bundle $\mathcal{Q}_0$. We have the following exact sequence of complexes of sheaves

{\small{\begin{center}
$\begin{CD}
@. \cdots @. \cdots @. \cdots \\
@. @AAA @AAA @AAA \\
0 @>>> 0 @>>> \mathscr{H}om_{\mathcal{O}_X}\left( \bigwedge^3 T_{\mathcal{F}_0}, \frac{T_X}{T_{\mathcal{F}_0}} \right) @>>> \mathscr{H}om_{\mathcal{O}_X}\left( \bigwedge^3 T_{\mathcal{F}_0}, \frac{T_X}{T_{\mathcal{F}_0}} \right) @>>> 0 \\
@. @AAA @AAA @AAA \\
0 @>>> 0 @>>> \mathscr{H}om_{\mathcal{O}_X}\left( \bigwedge^2 T_{\mathcal{F}_0}, \frac{T_X}{T_{\mathcal{F}_0}} \right) @>>> \mathscr{H}om_{\mathcal{O}_X}\left( \bigwedge^2 T_{\mathcal{F}_0}, \frac{T_X}{T_{\mathcal{F}_0}} \right) @>>> 0 \\
@. @AAA @AAA @AAA \\
0@>>> \mathscr{H}om_{\mathcal{O}_X}\left( T_{\mathcal{F}_0}, T_{\mathcal{F}_0} \right)|_{\textnormal{Tr}=0} @>>> \mathscr{H}om_{\mathcal{O}_X}\left( T_{\mathcal{F}_0}, T_X \right) @>>>  \mathcal{R} @>>> 0 \\
@. @A-\textnormal{Id}AA @AAA @AAA \\
0 @>>> \mathscr{H}om_{\mathcal{O}_X}\left( T_{\mathcal{F}_0}, T_{\mathcal{F}_0} \right)|_{\textnormal{Tr}=0}  @>>> \mathscr{H}om_{\mathcal{O}_X}\left(T_{\mathcal{F}_0}, T_{\mathcal{F}_0 } \right) @>-\textnormal{Tr}>> \mathcal{O}_X @>>> 0
\end{CD}$
\end{center}}}
We denote the complex of sheaves $\mathcal{O}_X\to \mathcal{R} \to \mathscr{H}om_{\mathcal{O}_X}\left(\bigwedge^2 T_{\mathcal{F}_0}, \frac{T_X}{T_{\mathcal{F}_0}} \right)\to \mathscr{H}om_{\mathcal{O}_X} \left( \bigwedge^3 T_{\mathcal{F}_0}, \frac{T_X}{T_{\mathcal{F}_0}} \right) \to \cdots$ by $\mathcal{Q}_0^{\bullet}$. Then the above exact sequence implies that
\begin{align*}
\mathbb{H}^i\left(X, \mathscr{H}om_{\mathcal{O}_X}\left(T_{\mathcal{F}_0},T_{\mathcal{F}_0} \right)^\bullet \right)\cong \mathbb{H}^i  \left( X, \mathcal{Q}_0^{\bullet}   \right)
\end{align*}
which gives the isomorphisms of first-oder foliated deformations and obstructions for $i=1,2$ respectively.

\section{Rigid foliated deformations, Schlessingers's criterion for $\textnormal{FDef}_{\left(X, T_{\mathcal{F}_0} \right)}$, stability and formal foliated deformations}

\subsection{Rigid foliated deformations}\

We recall Definition \ref{a5} that a foliated nonsingular variety $\left( X, T_{\mathcal{F}_0}\right)$ is called rigid in foliated deformations if every infinitesimal foliated deformation $\left(\mathcal{X}, T_\mathcal{F} \right)$ of $\left(X, T_{\mathcal{F}_0}\right)$ over $A$ is trivial for every $A$ in $\bold{Art}$, in other words, $\left( \mathcal{X},T_{\mathcal{F}} \right)\cong \left( X\times \textnormal{Spec}\left( A \right), T_{\mathcal{F}_0}\right)$.

\begin{proposition}\label{a7}
A foliated nonsingular variety $\left( X, T_{\mathcal{F}_0} \right)$ with $T_{\mathcal{F}_0}$ locally free is rigid in foliated deformations if and only if $\mathbb{H}^1\left( X, T_{\mathcal{F}_0}^\bullet \right)=0$.
\end{proposition}

\begin{proof}
Assume that $\left(X, T_{\mathcal{F}_0} \right)$ is rigid in foliated deformations. Since any infinitesimal foliated deformation (in particular, any first-order foliated deformation) is trivial, we have $\mathbb{H}^1\left( X , T_{\mathcal{F}_0}^\bullet \right)=0$ by Proposition \ref{at25}. Assume that $\mathbb{H}^1\left( X, T_{\mathcal{F}_0}^\bullet \right)=0$. First we claim that given an infinitesimal foliated deformation $\eta$ of $\left( X, T_{\mathcal{F}_0} \right)$ over $A\in \bold{Art}$ and a small extension $0\to (t)\to \tilde{A} \to A \to 0$, any two liftings $\xi, \xi'$ of $\eta$ to $\tilde{A}$ are equivalent. Let $\mathcal{U}=\left\{ U_i \right\}$ be an affine open covering of $X$ such that $\xi=\left( \mathcal{X} , T_\mathcal{F} \right)$ is represented by the following data: we have foliated isomorphisms $\theta_i:\left( U_i \times \textnormal{Spec}\left(\tilde{A} \right), T_{\mathcal{F}_i}  \right)\to \left( \mathcal{X}|_{ U_i}, T_{\mathcal{F}}|_{U_i} \right)$ where $T_{\mathcal{F}_i}$ is the foliated structure on $U_i\times \textnormal{Spec}\left( \tilde{A} \right)$ induced from $T_{\mathcal{F}}$ via $\theta_i$ and $\Gamma\left( U_i, T_{\mathcal{F}_i} \right)$ is generated by $T_i^1,..., T_i^p\in \Gamma\left( U_i, T_X\right)\otimes \tilde{A}$ with $\left[ T_i^\alpha, T_i^\beta \right] = \sum_{\gamma=1}^p g_{i\alpha\beta}^\gamma T_i^\gamma$ for $g_{i\alpha\beta}^\gamma \in \Gamma\left( U_i, \mathcal{O}_X \right)\otimes \tilde{A}$, and letting $\theta_{ij}=\theta_j^{-1}\theta_i$, we have foliated isomorphisms $\theta_{ij}:\left(U_{ij}\times \textnormal{Spec}\left( \tilde{A} \right), T_{\mathcal{F}_i} \right)\to \left( U_{ij}\times \textnormal{Spec}\left( \tilde{A} \right), T_{\mathcal{F}_j} \right)$ such that $\theta_{ji}^* T_i^\alpha= \sum_{\beta=1}^p r_{ij}^{\alpha\beta} T_j^\beta$ for $r_{ij}^{\alpha\beta}\in \Gamma\left( U_{ij}, \mathcal{O}_X \right)\otimes \tilde{A}$, and on the other hand $\xi'=\left( \mathcal{X}', T_{\mathcal{F}'} \right)$ is represented by the following data: similarly $\theta_i', \theta_{ij}'=\theta_j'^{-1}\theta_i'$ and $\left\{ T_{\mathcal{F}_i'}\right\}, T_i'^\alpha, g_{i\alpha\beta}'^\gamma ,  r_{ij}'^{\alpha\beta}$. We may assume that $T_i^\alpha$ and $T_i'^\alpha$ induces same element in $\Gamma\left( U_i, T_X\right)\otimes A$. Let $f_{ij}: \left(\Gamma\left( U_{ij}, \mathcal{O}_X \right) \otimes \tilde{A}, T_{\mathcal{F}_j} \right) \to \left( \Gamma\left( U_{ij}, \mathcal{O}_X \right) \otimes \tilde{A}, T_{\mathcal{F}_i}\right)$ be the ring homomorphism corresponding to $\theta_{ij}$ and similarly $f_{ij}'$ corresponding to $\theta_{ij}'$. Since $\xi, \xi'$ induces $\eta$ over $A$, we have $f_{ij}'= f_{ij} + t p_{ij}$ for some $p_{ij} \in \Gamma\left( U_{ij}, T_X \right)$ and $T_i'^\alpha=  T_i^\alpha  + t V_i^\alpha$ for some $V_i^\alpha \in \Gamma\left(U_i, T_X \right), \alpha=1,...,p$, and $r_{ij}'^{\alpha\beta}= r_{ij}^{\alpha\beta} +t a_{ij}^{\alpha\beta}$ for some $a_{ij}^{\alpha\beta}\in \Gamma\left(U_{ij}, \mathcal{O}_X \right)$, and $g_{i\alpha\beta}'^\gamma= g_{i\alpha\beta}'^\gamma + t b_{i\alpha\beta}^\gamma$ for some $b_{i\alpha\beta}^\gamma \in \Gamma\left( U_i, \mathcal{O}_X \right)$. We define $\overline{V_i}\in \Gamma\left( U_i, \mathscr{H}om_{\mathcal{O}_X}\left( T_{\mathcal{F}_0}, \frac{T_X}{T_{\mathcal{F}_0}} \right) \right)$ by $\overline{V_i}: \Gamma\left( U_i, T_{\mathcal{F}_0} \right) \to \Gamma\left( U_i, \frac{T_X}{T_{\mathcal{F}_0}} \right), T_{0i}^\alpha \mapsto \overline{ V_i^\alpha}$. Then we claim that $\left(\left\{ \overline{V_i} , \left\{ p_{ij} \right\}  \right\} \right) \in C^0\left( \mathcal{U}, \mathscr{H}om_{\mathcal{O}_X}\left( T_{\mathcal{F}_0}, \frac{T_X}{T_{\mathcal{F}_0}} \right) \right) \bigoplus C^1\left(\mathcal{U}, T_X \right)$ defines a $1$-cocycle of the \v Cech resolution of $T_{\mathcal{F}_0}^\bullet$. In fact, we have $p_{ij}- p_{ik} + p_{jk}=0$, and 
\begin{align*}
\left[ T_i'^\alpha, T_i'^\beta \right] &= \left[T_i^\alpha + t V_i^\alpha, T_i^\beta + tV_i^\beta \right] = \sum_{\gamma=1}^p g_{i\alpha\beta}^\gamma T_i^\gamma + t \left[ T_{0i}^\alpha, V_i^\beta \right] - t \left[ T_{0i}^\beta, V_i^\beta \right]\\
\left[ T_i'^\alpha, T_i'^\beta \right]&= \sum_{\gamma=1}^p g_{i\alpha\beta}'^\gamma T_i'^\gamma= \sum_{\gamma=1}^p \left(g_{i\alpha\beta}^\gamma + t b_{i\alpha\beta}^\gamma\right)\left( T_i^\gamma + t V_i^\gamma \right)= \sum_{\gamma=1}^p g_{i\alpha\beta}^\gamma T_i^\gamma  + t \sum_{\gamma=1}^p g_{0i\alpha\beta}^\gamma V_i^\gamma + t \sum_{\gamma=1}^p b_{i\alpha\beta}^\gamma T_{0i}^\gamma
\end{align*}
This implies that $D_1\left( \overline{V_i} \right)=0$. On the other hand, we note that 
\begin{align*}
\theta_{ji}'^*\left(  T_i'^\alpha \right) &= \theta_{ji}'^*\left( T_i^\alpha + t V_i^\alpha \right) = \theta_{ji}^* T_i^\alpha - t \left[ p_{ij}, T_{0i}^\alpha \right]+ t V_i^\alpha\\
\theta_{ji}'^* \left( T_i'^\alpha \right) &= \sum_{\beta=1}^p r_{ij}'^{\alpha\beta} T_j'^\beta = \sum_{\beta=1}^p \left( r_{ij}^{\alpha\beta} + t a_{ij}^{\alpha\beta} \right)\left( T_j^\beta + t V_j^\beta \right) = \sum_{\beta=1}^p r_{ij}^{\alpha\beta} T_j^\beta + t \sum_{\beta=1}^p r_{0ij}^{\alpha\beta} V_j^\beta + t \sum_{\beta=1}^p  a_{ij}^{\alpha\beta} T_{0j}^\beta
\end{align*}
This implies that $\delta\left( \left\{ \overline{V_i^\alpha }\right\} \right) + D_0\left( \left\{ p_{ij} \right\} \right)=0$. This proves the claim. Since $\mathbb{H}^1\left( X, T_{\mathcal{F}_0}^\bullet \right)=0$, there exists $\left\{ a_i \right\}\in C^0\left( \mathcal{U}, T_X \right)$ such that $a_i-a_j= p_{ij}$ and $\left[ a_i, T_{0i}^\alpha \right]- V_i^\alpha = \sum_{\beta=1}^p c_i^{\alpha\beta} T_{0i}^\beta $ for some $c_i^{\alpha\beta}\in \Gamma\left(  U_i, \mathcal{O}_X \right)$. Now we explicitly construct a foliated isomorphism $\left( \mathcal{X}, T_\mathcal{F} \right) \cong \left( \mathcal{X}', T_{\mathcal{F}'} \right)$. We claim that $\textnormal{Id}- t a_i:\left(\Gamma\left( U_i, \mathcal{O}_X\right)\otimes \tilde{A} , T_{\mathcal{F}_i'} \right) \to \left( \Gamma\left( U_i, \mathcal{O}_X \right)\otimes \tilde{A}, T_{\mathcal{F}_i} \right)$ is a foliated isomorphism. In fact, $T_i'^\alpha- t\left[ a_i, T_{0i}^\alpha\right]= T_i^\alpha + tV_i^\alpha - t \left[ a_i, T_{0i}^\alpha \right]= T_i^\alpha -t \sum_{\beta=1}^p c_i^{\alpha\beta} T_{0i}^\beta = \sum_{\beta=1}^p \left(\delta_\alpha^\beta - t c_i^{\alpha\beta} \right) T_i^\beta$. On the other hand, since $p_{ij}= a_i- a_j$, we have $\left( \textnormal{Id} - ta_i \right)\circ f_{ij}'= f_{ij}\circ \left(\textnormal{Id}- t a_j \right)$. Hence $\left\{ \textnormal{Id} - ta_i \right\}$ glues together to give a foliated isomorphism $\left( \mathcal{X}, T_{\mathcal{F}} \right)\cong \left( \mathcal{X}', T_{\mathcal{F}'} \right)$.

Now we prove if $\mathbb{H}^1\left( X, T_{\mathcal{F}_0}^\bullet \right)=0$, then $\left( X, T_{\mathcal{F}_0} \right)$ is rigid in foliated deformations. We prove by induction on the dimension of $A\in \bold{Art}$. Assume that any infinitesimal foliated deformation of $\left( X, T_{\mathcal{F}_0} \right)$ over $A$ with $\dim_k A \leq n-1$ is trivial. Let $\xi$ be an infinitesimal foliated deformation of $\left(X, T_{\mathcal{F}_0} \right)$ over $A$ with $\dim_k A =n$ such that $\mathfrak{m}^{p-1}\ne 0$ and $\mathfrak{m}^p=0$. Choose an element $t\ne 0 \in \mathfrak{m}^{p-1}$. Then $0\to (t)\to A \to A/(t)\to 0$ is a small extension and $\dim_k A/(t)\leq n-1$. Hence the induced foliated deformation $\overline{\xi}$ over $A/(t)$ from $\xi$ is trivial by induction hypothesis. Since $\xi$ is a lifting of $\overline{\xi}$ and the trivial foliated deformation over $A$ is also a lifting of $\overline{\xi}$, $\xi$ is trivial since we have at most one lifting of $\overline{\xi}$. This completes the proof of Proposition \ref{a7}.
\end{proof}

\subsection{Schlessinger's criterion for $\textnormal{FDef}_{\left( X, T_{\mathcal{F}_0} \right) }$} \

\begin{lemma}\label{a9}
Let $\left(X, T_{\mathcal{F}_0}\right)$ be a foliated nonsingular projective variety with $T_{\mathcal{F}_0}$ locally free. Assume that $\mathbb{H}^0\left(X, T_{\mathcal{F}_0}^\bullet \right)=0$. Then for any infinitesimal foliated deformation $(\mathcal{X}, T_\mathcal{F})$ of $\left(X, T_{\mathcal{F}_0}\right)$ over $A$ for an $A\in \bold{Art}$.
\begin{align*}
\textnormal{Aut}\left(\frac{\left( \mathcal{X}, T_\mathcal{F}\right) }{\left(X, T_{\mathcal{F}_0}\right)} \right)=\textnormal{Id}
\end{align*}
where $\textnormal{Aut}\left(\frac{\left( \mathcal{X}, T_\mathcal{F}\right) }{\left(X, T_{\mathcal{F}_0}\right)} \right):=$ the set of foliated automorphisms of $\left(\mathcal{X}, T_\mathcal{F} \right)$ restricting to the identity on $\left( X, T_{\mathcal{F}_0} \right)$.
\end{lemma}

\begin{proof}
We prove by the induction on the dimension of $A$. Let $\dim_k A=1$. Then $A=k$. So we have nothing to prove. Let's assume that Lemma \ref{a9} holds for $A$ with $\dim_k A\leq n-1$. Let $\dim_k A=n$ and $\left(\mathcal{X}, T_{\mathcal{F}} \right)$ be an infinitesimal foliated deformation of $\left( X, T_{\mathcal{F}_0} \right)$ over $A$. Assume that the maximal ideal $\mathfrak{m}$ of $A$ satisfies $\mathfrak{m}^{p-1}\ne 0$ and $\mathfrak{m}^p=0$. Choose $t\ne 0 \in \mathfrak{m}^{p-1}$. Then $A/(t)\in \bold{Art}$ with $\dim_k A/(t)\leq n-1$ and $0\to (t)\to A \to A/(t)\to 0$ is a small extension. Now let $g:\left(\mathcal{X}, T_\mathcal{F} \right)\to \left( \mathcal{X}, T_\mathcal{F} \right)$ over $A$ be a foliated automorphism with $g^* T_\mathcal{F}=T_{\mathcal{F}}$ restricting to the identity on $\left( X, T_{\mathcal{F}_0} \right)$. Let $\mathcal{U}=\{U_i\}$ be an affine cover of $X$ such that $\Gamma\left( U_i, T_{\mathcal{F}_0}\right)$ is generated by $T_{0i}^1,..., T_{0i}^p\in \Gamma\left(U_i, T_X \right)$ with $T_{0i}^\alpha= \sum_{\beta=1}^p r_{0ij}^{\alpha\beta} T_{0j}^\beta$ and $\left[T_{0i}^\alpha, T_{0i}^\beta \right]= \sum_{\gamma=1}^p g_{0i\alpha\beta}^\gamma T_{0i}^\gamma$, and $\theta_i: \left(U_i\times \textnormal{Spec}(A), T_{\mathcal{F}_i}\right)\cong \left(\mathcal{X}|_{U_i}, T_{\mathcal{F}}|_{U_i} \right)$ is an foliated isomorphism such that $\Gamma\left(U_i, T_{\mathcal{F}_i} \right)$ is generated by $T_i^1,..., T_i^p\in \Gamma\left(U_i, T_X \right)\otimes A$ which induces $T_{0i}^1,..., T_{0i}^p$. Let $\theta_{ij}=\theta_j^{-1}\theta_i$  correspond to a foliated isomorphism $f_{ij}:\Gamma\left( \Gamma\left( U_{ij}, \mathcal{O}_X \right)\otimes A, T_{\mathcal{F}_j} \right)\to  \left( \Gamma\left(U_{ij}, \mathcal{O}_X \right)\otimes A , T_{\mathcal{F}_i} \right)$. Then $g$ is described by the data $\{g_i\}$, where $g_i:\left(\Gamma\left( U_i, \mathcal{O}_X\right)\otimes A, T_{\mathcal{F}_i} \right)\to \left(  \Gamma\left(U_i, \mathcal{O}_X \right)\otimes A, T_{\mathcal{F}_i} \right)$ which is a foliated isomorphism and $g_i\circ f_{ij}=f_{ij}\circ g_j$. By the induction hypothesis, $g_i$ induces the identity on $\Gamma\left( U_i, \mathcal{O}_X \right)\otimes A/(t)$. Then $g_i$ is of the form $g_i=\textnormal{Id}+ td_i$ such that $d_i\in \Gamma\left(U_i, T_X \right)$ with $d_i=d_j$ and $T_i^\alpha+ t \left[ d_i, T_{0i}^\alpha \right]= \sum_{\beta=1}^p p_{i}^{\alpha\beta} T_i^\beta $ for some $p_i^{\alpha\beta} \in \Gamma\left(U_i, \mathcal{O}_X \right) \otimes A$. On the other hand, we note that by tensoring $0\to (t)\to A\to A/(t)\to 0$ with $T_{\mathcal{F}}$, we have $0\to t T_{\mathcal{F}_0}\to T_{\mathcal{F}} \to T_{\mathcal{F}}\otimes A/(t) \to 0$ ($T_\mathcal{F}$ is locally free and $\mathcal{X}$ is flat over $A$). Then since $T_i^\alpha - \sum_{\beta=1}^p p_i^{\alpha\beta} T_i^\beta$ induces $0$ in $\Gamma\left(U_i, T_{\mathcal{F}} \right)\otimes A/(t)$, we can write $T_i^\alpha - \sum_{\beta=1}^p p_i^{\alpha\beta} T_i^\beta= t W_{i}^\alpha$ for some $W_i^\alpha\in \Gamma\left(U_i, T_{\mathcal{F}_0} \right)$. This implies that $\left[ d_i, T_{0i}^\alpha \right]\in \Gamma\left(U_i, T_{\mathcal{F}_0} \right)$. Hence $d:=\left\{d_i \right\} \in \mathbb{H}^0\left(X, T_{\mathcal{F}_0}^\bullet \right)$. Since $\mathbb{H}^0\left(X, T_{\mathcal{F}_0}^\bullet \right)=0$, we have $d=0$, so that $g=\textnormal{Id}$.
\end{proof}

\begin{proposition}\label{e4}
Let $\left(X, T_{\mathcal{F}_0} \right)$ be a foliated nonsingular variety with $T_{\mathcal{F}_0}$ locally free and $e:0\to (t) \to \tilde{A} \xrightarrow{\mu} A\to 0$ be a small extension in $\bold{Art}$. Let $p:=\textnormal{FDef}_{(X, T_{\mathcal{F}_0})}(\mu):\textnormal{FDef}_{(X, T_{\mathcal{F}_0})}(\tilde{A})\to  \textnormal{FDef}_{(X, T_{\mathcal{F}_0})}(A) $. Then given $\xi= (\mathcal{X}, T_{\mathcal{F}})\in \textnormal{FDef}_{(X, T_{\mathcal{F}_0})}(A)$, there is a transitive action of $\mathbb{H}^1\left( X, T_{\mathcal{F}_0}^\bullet \right)$ on $p^{-1}(\xi)$. Moreover, if $\mathbb{H}^0\left( X, T_{\mathcal{F}_0}^\bullet \right) = 0$, the action is free.
\end{proposition}

\begin{proof}
We will define a map $G:\mathbb{H}^1\left(X, T_{\mathcal{F}_0}^\bullet \right)\times p^{-1}(\xi)\to p^{-1}(\xi)$ which is a group action of $\mathbb{H}^1\left(X, T_{\mathcal{F}_0}^\bullet \right)$ acting on $p^{-1}(\xi)$. Let $\tilde{\xi}=\left(\tilde{\mathcal{X}}, T_{\tilde{\mathcal{F}}} \right)\in p^{-1}\left(\xi \right)$ over $\tilde{A}$. Let $\mathcal{U}=\{U_i\}$ be an affine open cover of $X$ such that $\Gamma\left(U_i, T_{\mathcal{F}_0} \right)$ is generated by $T_{0i}^1,..., T_{0i}^p\in \Gamma\left(U_i, T_X \right)$ with $T_{0i}^\alpha= \sum_{\beta=1}^p r_{0ij}^{\alpha\beta} T_{0j}^\beta$ and $\left[ T_{0i}^\alpha, T_{0i}^\beta\right]= \sum_{\gamma=1}^p g_{0i\alpha\beta}^\gamma T_{0i}^\gamma$, and $\tilde{\theta}_i:\left( U_i\times \textnormal{Spec}(A),  T_{\tilde{\mathcal{F}}_i}\right)\cong \left(\tilde{\mathcal{X}}|_{U_i}, T_{\tilde{\mathcal{F}}}|_{U_i} \right)$ is a foliated isomorphism such that $\Gamma\left( U_i, T_{\tilde{\mathcal{F}}_i}\right)$ is generated by $\tilde{T}_i^1,..., \tilde{T}_i^p\in \Gamma\left( U_i, T_X\right)\otimes \tilde{A}$ which induces $T_{0i}^1,..., T_{0i}^p$, and $\tilde{\theta}_{ji}^* \tilde{T}_i^\alpha= \sum_{\beta=1}^p \tilde{r}_{ij}^{\alpha\beta} \tilde{T}_j^\beta$ where $\tilde{\theta}_{ij}=\tilde{\theta}_j^{-1}\tilde{\theta}_i$ and $\tilde{r}_{ij}^{\alpha\beta}$ induces $r_{0ij}^{\alpha\beta}$, and $\left[ \tilde{T}_i^\alpha, \tilde{T}_i^\beta \right]= \sum_{\gamma=1}^p \tilde{g}_{i\alpha\beta}^\gamma \tilde{T}_i^\gamma$ where $\tilde{g}_{i\alpha\beta}^\gamma$ induces $g_{0i\alpha\beta}^\gamma$.  Let $\tilde{\theta}_{ij}$  correspond to a foliated isomorphism $\tilde{f}_{ij}:\Gamma\left( \Gamma\left( U_i, \mathcal{O}_X \right)\otimes \tilde{A}, T_{\tilde{\mathcal{F}}_j} \right)\to  \left( \Gamma\left(U_{ij}, \mathcal{O}_X \right)\otimes \tilde{A} , T_{\tilde{\mathcal{F}}_i} \right)$. Let $\left( \left\{ \overline{ T_i' } \right\}, \left\{ p_{ij}\right\} \right)\in \mathbb{H}^1\left( X, T_{\mathcal{F}_0}^\bullet \right)$ corresponding to a first-order foliated deformation as in the proof of Proposition \ref{at25}. We keep the notation in the proof of Proposition \ref{at25}. Then $\left\{\tilde{f}_{ij}+ t p_{ij}:\Gamma\left(U_{ij}, \mathcal{O}_X\right)\otimes \tilde{A} \to \Gamma\left(U_{ij}, \mathcal{O}_X\right)\otimes \tilde{A} \right\} $ which corresponds to isomorphisms $\tilde{\theta}_{ij}': U_{ij}\times \textnormal{Spec}\left(\tilde{A} \right)\to U_{ij}\times \textnormal{Spec}\left(\tilde{A}\right)$ satisfies the gluing condition, and since $\tilde{\theta}_{ji}'^*\left( \tilde{T}_i^\alpha + t T_i'^\alpha \right)= \tilde{\theta}_{ji}^*\left( \tilde{T}_i^\alpha \right) +t T_i'^\alpha - t\left[ p_{ij}, T_{0i}^\alpha \right]=\sum_{\beta=1}^p \tilde{r}_{ij}^{\alpha\beta} \tilde{T}_j^\beta +t\sum_{\beta=1}^p r_{0ij}^{\alpha\beta} T_j'^\beta + t \sum_{\beta=1}^p b_{ij}^{\alpha\beta} T_{0j}^\beta = \sum_{\beta=1}^p \left( \tilde{r}_{ij}^{\alpha\beta} + t b_{ij}^{\alpha\beta} \right)\left( \tilde{T}_j^\beta+ t T_j'^\beta \right)$, $\tilde{\theta}_{ij}'$ is a foliated isomorphism.  On the other hand, $\left[\tilde{T}_i^\alpha+ tT_i'^\alpha, \tilde{T}_i^\beta + t T_i'^\beta \right]=\sum_{\gamma=1}^p \tilde{g}_{i\alpha\beta}^\gamma \tilde{T}_i^\gamma + t\sum_{\gamma=1}^p g_{0i\alpha\beta}^\gamma T_i'^\gamma + t\sum_{\gamma=1}^p c_{i\alpha\beta}^\gamma T_{0i}^\gamma= \sum_{\gamma=1}^p \left(\tilde{g}_{i\alpha\beta}^\gamma + t c_{i\alpha\beta}^\gamma \right) \left( \tilde{T}_i^\gamma + t T_i'^\gamma \right)$. Hence this defines an infinitesimal foliated deformation $\tilde{\xi}'=\left( \tilde{\mathcal{X}}', T_{\tilde{\mathcal{F}}'}\right)$ of $\left( X, T_{\mathcal{F}_0} \right)$ over $\tilde{A}$ which induces $\xi=\left(\mathcal{X}, T_{\mathcal{F}} \right)$, so that $\left( \tilde{\mathcal{X}}' , T_{\tilde{\mathcal{F}}'}\right)$ lies in $p^{-1}(\xi)$. Then we define a map $G$ by 
\begin{align*}
G:\mathbb{H}^1\left(X, T_{\mathcal{F}_0}^\bullet \right)\times p^{-1}(\xi) &\to p^{-1}(\xi) \\
                       \left( \left\{ \overline{T_i'}\right\}, \left\{ p_{ij}\right\} \right)\times \tilde{\xi} &\mapsto \tilde{\xi}'
\end{align*}
$G$ is well-defined. In fact, let $\left( \left\{  \overline{T_i''} \right\}, \left\{ p_{ij}' \right\}  \right)$ represent the same cohomology class with $\left( \left\{ \overline{T_i'} \right\}, \left\{ p_{ij} \right\} \right)$. Then we have to show that $\tilde{\xi}'$ defined by $\left( \left\{ \overline{T_i'} \right\} , \left\{ p_{ij} \right\} \right)$ is equivalent to $\tilde{\xi}''$ defined by $\left( \left\{ \overline{T_i''} \right\} , \left\{ p_{ij}' \right\} \right)$. Then there exists $\left\{ a_i \right\}\in C^0\left( \mathcal{U}, T_X \right)$ such that $a_i- a_j= p_{ij}- p_{ij}'$ and $\left[ a_i, T_{0i}^\alpha \right]  = T_i'^\alpha - T_i''^\alpha - \sum_{\beta=1}^p b_i^{\alpha\beta} T_{0i}^\beta$ for $b_i^{\alpha\beta} \in \Gamma\left( U_i, \mathcal{O}_X\right)$. We claim that $\textnormal{Id}- t a_i: \left( \Gamma\left( U_i, \mathcal{O}_X \right)\otimes \tilde{A},  T_{\tilde{\mathcal{F}}_i'} \right)\to \left( \Gamma\left( U_i, \mathcal{O}_X\right)\otimes \tilde{A} , T_{\tilde{\mathcal{F}}_i''}\right)$ is a foliated isomorphism. In fact, $\tilde{T}_i^\alpha  + t T_i'^\alpha - t \left[ a_i, T_{0i}^\alpha \right] = \tilde{T}_i^\alpha + t T_i''^\alpha + t\sum_{\beta=1}^p b_i^{\alpha\beta} T_{0i}^\beta= \sum_{\beta=1}^p \left( \delta_\alpha^\beta + t b_i^{\alpha\beta} \right)\left( \tilde{T}_i^\beta + t T_i''^\beta \right)$. On the other hand, $\left( \tilde{f}_{ij} + tp_{ij}' \right)\circ \left( \textnormal{Id}- t a_j \right) = \left(\textnormal{Id}- t a_i \right)\circ \left( \tilde{f}_{ij} + t p_{ij}\right)$. Hence $\left\{\textnormal{Id} - t a_i \right\}$ glues together to give a foliated isomorphism between $\tilde{\xi}'$ and $\tilde{\xi}''$.

Next we show that the group action $\mathbb{H}^1\left( X, T_{\mathcal{F}_0}^\bullet \right)\times p^{-1}(\xi)\to p^{-1}(\xi)$ is transitive. In fact, let $\tilde{\xi}\in p^{-1}\left( \xi\right)$ be a lifting of $\xi$ as above. Choose arbitrary lifting $\eta\in p^{-1}\left( \xi \right)$ of $\xi$. We have to find $\left( \left\{\overline{T_i'}\right\}, \left\{ p_{ij}\right\} \right)\in \mathbb{H}^1\left(X, T_{\mathcal{F}_0}^\bullet  \right)$ such that $\left( \left( \left\{ \overline{T_i'}\right\}, \left\{ p_{ij} \right\} \right), \tilde{\xi} \right)$ is mapped to $\eta$ under the action $G$. Let $\eta$ consist of the data: $f_{ij}':\left( \Gamma\left( U_{ij}, \mathcal{O}_X \right)\otimes \tilde{A}, T_{\tilde{\mathcal{F}}_j'} \right)\to \left(  \Gamma\left( U_{ij}, \mathcal{O}_X \right)\otimes \tilde{A}, T_{\tilde{\mathcal{F}}_i'} \right)$ such that $\tilde{f}_{ij}'$ and $\tilde{f}_{ij}$ induces the same homomorphism in $\xi$,  and $\Gamma\left( U_i, T_{\tilde{\mathcal{F}}_i'} \right)$ is generated by $\tilde{T}_i'^\alpha\in \Gamma\left( U_i, T_X\right)\otimes \tilde{A},\alpha=1,...,p$ such that $\tilde{T}_i'^\alpha$ and $\tilde{T}_i^\alpha$ induce the same element in $\xi$. Since $\tilde{\xi}$ and $\eta$ induce $\xi$, we have $\tilde{f}_{ij}'= f_{ij}+t p_{ij}$ and $\tilde{T}_i'^\alpha= \tilde{T}_i^\alpha + t V_i^\alpha$. We define $\overline{V_i}:\Gamma\left( U_i, T_{\mathcal{F}_0} \right)\to \Gamma\left( U_i, \frac{T_X}{T_{\mathcal{F}_0}} \right)$ by $T_{0i}^\alpha \mapsto \overline{V_i^\alpha}$. Then as in the same argument in the proof of Proposition \ref{a7}, $\left( \left\{ \overline{V_i} \right\} , \left\{ p_{ij} \right\} \right)$ defines a $1$-cocycle of the \v Cech resolution of $T_{\mathcal{F}_0}^\bullet$.

Now assume that $\mathbb{H}^0\left( X , T_{\mathcal{F}_0}^\bullet \right)=0$. We claim that the group action $G$ is free. Assume that for given $v=\left( \left\{ \overline{T_i'} \right\}, \left\{ p_{ij} \right\} \right) \in \mathbb{H}^1\left( X, T_{\mathcal{F}_0}^\bullet \right)$ and $\tilde{\xi}\in p^{-1}(\xi)$, we have $G\left( v, \tilde{\xi}\right)=\tilde{\xi}$. We have to show that $v=0\in \mathbb{H}^1\left( X, T_{\mathcal{F}_0}^\bullet \right)$. Let $\tilde{\xi}$ and $\tilde{\xi}'=G\left( v, \tilde{\xi}\right)$ as above. Then $\tilde{\xi}'=G\left( v, \tilde{\xi} \right)=\tilde{\xi}$ implies that we have foliated isomorphisms $g_i: \left(\Gamma\left( U_i, \mathcal{O}_X \right)\otimes \tilde{A} , T_{\tilde{\mathcal{F}}_i} \right) \to \left(\Gamma\left( U_i, \mathcal{O}_X \right)\otimes \tilde{A}, T_{\tilde{\mathcal{F}}_i'}  \right)$ such that $g_i \circ f_{ij}= \left( f_{ij} + t p_{ij}\right)\circ g_j$. Since $\left\{ g_i\right\}$ induces a foliated automorphism on $\xi=\left( \mathcal{X}, T_\mathcal{F} \right)$ by modulo $(t)$, which is necessarily the identity by Lemma \ref{a9}, $g_i$ is of the form $\textnormal{Id}+ t q_i$ for some $q_i\in \Gamma\left( U_i, T_X \right)$. Then $\tilde{T}_i^\alpha + t \left[ q_i, T_{0i}^\alpha \right] = \sum_{\beta=1}^p \left(\delta_\alpha^\beta + t b_i^{\alpha\beta} \right) \left( \tilde{T}_i^\beta +t T_i'^\beta \right)$ for some $b_i^{\alpha\beta}\in \Gamma\left(U_i, \mathcal{O}_X \right)$, so that $\left[ q_i, T_{0i}^\alpha \right]= T_i'^\alpha + \sum_{\beta=1}^p b_i^{\alpha\beta} T_{0i}^\beta$. On the other hand, we have $\left(\textnormal{Id}+ t q_i\right)\circ f_{ij}= \left( f_{ij}+ t p_{ij} \right)\circ \left( \textnormal{Id} + t q_j \right)$, so that $q_i- q_j= p_{ij}$. This implies $v=0 \in \mathbb{H}^1\left( X , T_{\mathcal{F}_0}^\bullet \right)$. This completes the proof of Proposition \ref{e4}.

\end{proof}

\begin{proposition}\label{e3}
Let $\left( X, T_{\mathcal{F}_0} \right)$ be a foliated nonsingular projective variety with $T_{\mathcal{F}_0}$ locally free. Then $\textnormal{FDef}_{\left( X , T_{\mathcal{F}_0} \right)}$ has a semiuniversal formal element in the sense of \cite{Ser06} or a miniversal family in the sense of \cite{Har10}.
\end{proposition}

\begin{proof}
We check Schlessinger's criterion for $\textnormal{FDef}_{\left( X, T_{\mathcal{F}_0} \right)}$ (see \cite{Har10} Theorem 16.2 p.113). $\textnormal{FDef}_{\left(X, T_{\mathcal{F}_0} \right)}(k)=\left( X, T_{\mathcal{F}_0} \right)$, so $\left(\textnormal{H}_0 \right)$ is satisfied. Let $\mathcal{U}=\left\{ U_i \right\}$ be an affine open covering such that $\Gamma\left( U_i, T_{\mathcal{F}_0} \right)$ is generated by $T_{0i}^1,..., T_{0i}^p \in \Gamma\left(U_i, T_X \right)$ such that $T_{0i}^\alpha = \sum_{\beta=1}^p r_{0ij}^{\alpha\beta} T_{0j}^\beta$ for some $h_{0ij}^{\alpha\beta} \in \Gamma\left( U_{ij}, \mathcal{O}_X \right)$ and $\left[ T_{0i}^\alpha, T_{0i}^\beta \right]= \sum_{\gamma=1}^p g_{0i\alpha\beta}^\gamma T_{0i}^\gamma$ for some $g_{0i\alpha\beta}^\gamma\in \Gamma\left(U_i, \mathcal{O}_X \right),\alpha,\beta=1,...,p$. Let $A'\to A$ in $\bold{Art}$ and $A''\to A$ be a small extension in $\bold{Art}$. Let $\left(\mathcal{X}, T_\mathcal{F} \right)$ be a foliated deformation of $\left( X, T_{\mathcal{F}_0} \right)$ over $A$, and $\left( \mathcal{X}', T_{\mathcal{F}'} \right)$ be a foliated deformation of $\left( X, T_{\mathcal{F}_0} \right)$ over $A'$ which induces $\left(\mathcal{X}, T_\mathcal{F} \right)$, and $\left(\mathcal{X}'', T_{\mathcal{F}''} \right)$ be a foliated deformation of $\left( X, T_{\mathcal{F}_0}\right)$ over $A''$ which induces $\left( \mathcal{X}, T_\mathcal{F}\right)$. We may assume that $\left( \mathcal{X}, T_\mathcal{F} \right)$ is given by the following data: we have the foliated isomorphism $\theta_i:\left( U_i\times \textnormal{Spec}\left( A \right), T_{\mathcal{F}_i} \right)\to \left( \mathcal{X}|_{U_i}, T_{\mathcal{F}}|_{U_i} \right)$ inducing the identity on $\left(X, T_{\mathcal{F}_0} \right)$ such that $\Gamma\left(U_i, T_{\mathcal{F}_i} \right)$ is generated by $T_i^\alpha  \in \Gamma\left( U_i, T_X\right)\otimes A$ which induces $T_{0i}^\alpha,\alpha=1,...,p$, and  by letting $f_{ij}:\left(\Gamma\left(U_{ij}, \mathcal{O}_X \right)\otimes A, T_{\mathcal{F}_j}\right) \to \left(\Gamma\left( U_{ij}, \mathcal{O}_X \right)\otimes A, T_{\mathcal{F}_i}\right) $ be the ring homomorphism corresponding to $\theta_{ij}:=\theta_j^{-1}\theta_i$, we have $f_{ik}= f_{ij} f_{jk}$ and $\theta_{ji}^*\left(T_i^\alpha \right)= \sum_{\beta=1}^p r_{ij}^{\alpha\beta} T_j^\beta$ with $r_{ij}^{\alpha\beta}\in \Gamma\left(U_{ij}, \mathcal{O}_X \right)\otimes A$ which induces $r_{0ij}^{\alpha\beta}$, and $\left[ T_i^\alpha, T_i^\beta \right]=\sum_{\beta=1}^p g_{i\alpha\beta}^\gamma T_i^\gamma$ for some $g_{i\alpha\beta}^\gamma\in \Gamma\left( U_i, \mathcal{O}_X \right)\otimes A$ which induces $g_{0i\alpha\beta}^\gamma,\alpha,\beta=1,...p$. On the other hand, we may assume that $\left( \mathcal{X}', T_{\mathcal{F}'} \right)$ is given by the following data: we have the foliated isomorphism $\theta_i':\left( U_i\times \textnormal{Spec}\left( A' \right) , T_{\mathcal{F}_i'}\right)\to \left( \mathcal{X}'|_{U_i}, T_{\mathcal{F}'}|_{U_i} \right)$ which induces $\theta_i$ such that $\Gamma\left(U_i, T_{\mathcal{F}_i'} \right)$ is generated by $T_i'^\alpha\in \Gamma\left( U_i, T_X\right)\otimes A',\alpha=1,...,p$ which induces $T_i^\alpha$, and by letting $f_{ij}:\left(\Gamma\left(U_{ij}, \mathcal{O}_X \right)\otimes A', T_{\mathcal{F}_j'}\right) \to \left(\Gamma\left( U_{ij}, \mathcal{O}_X \right)\otimes A, T_{\mathcal{F}_i'}\right) $ be the ring homomorphism corresponding to $\theta_{ij}':=\theta_j'^{-1}\theta_i'$, we have $f_{ik}'= f_{ij}' f_{jk}'$ and $\theta_{ji}'^*\left(T_i'^\alpha \right)= \sum_{\beta=1}^p r_{ij}'^{\alpha\beta} T_j'^\beta$ with $r_{ij}'^{\alpha\beta}\in \Gamma\left(U_{ij}, \mathcal{O}_X \right)\otimes A'$ which induces $r_{ij}^{\alpha\beta}$, and $\left[ T_i'^\alpha, T_i'^\beta \right]=\sum_{\beta=1}^p g_{i\alpha\beta}'^\gamma T_i'^\gamma$ for some $g_{i\alpha\beta}'^\gamma\in \Gamma\left( U_i, \mathcal{O}_X \right)\otimes A'$ which induces $g_{i\alpha\beta}^\gamma,\alpha,\beta=1,...p$, and we may assume similarly for $\left( \mathcal{X}'', T_{\mathcal{F}''} \right)$ (i.e. $\theta_i''$ induces $\theta_i$, and $T_i''^\alpha$ induces $T_i^\alpha$, and $r_{ij}''^{\alpha\beta}$ induces $r_{ij}^{\alpha\beta}$ and $g_{i\alpha\beta}''^\gamma$ induces $g_{i \alpha\beta}^\gamma$). 

Let $A^*=A'\times_A A''$. We define $\mathcal{X}^*$ to be the fibered sum of sheaves of rings $\mathcal{O}_{\mathcal{X}^*}= \mathcal{O}_{\mathcal{X}'}\times_{\mathcal{O}_\mathcal{X}} \mathcal{O}_{\mathcal{X}''}$ which can be described in terms of the covering $\left\{ U_i \right\}$ as follows: $\Gamma\left( U_i, \mathcal{O}_{\mathcal{X}^*} \right)=\Gamma\left(U_i, \mathcal{O}_{\mathcal{X}'} \right)\times_{\Gamma\left( U_i, \mathcal{O}_\mathcal{X} \right)} \Gamma\left( U_i, \mathcal{O}_{\mathcal{X}''} \right)\cong \left( \Gamma\left( U_i, \mathcal{O}_X \right)\otimes A'\right)\times_{\Gamma\left( U_i, \mathcal{O}_X \right)\otimes A} \left(\Gamma\left( U_i, \mathcal{O}_X \right)\otimes A'' \right)=\Gamma\left( U_i , \mathcal{O}_X \right)\otimes A^*$ and $f_{ij}^*=\left( f_{ij}', f_{ij}'' \right):\left(  \Gamma\left( U_{ij}, \mathcal{O}_X \right)\otimes A'  \right) \times_{\Gamma\left( U_{ij} , \mathcal{O}_X \right)\otimes A} \left( \Gamma\left( U_{ij}, \mathcal{O}_X \right) \otimes A'' \right) \to \left(  \Gamma\left( U_{ij}, \mathcal{O}_X \right)\otimes A'  \right) \times_{\Gamma\left( U_{ij} , \mathcal{O}_X \right)\otimes A} \left( \Gamma\left( U_{ij}, \mathcal{O}_X \right) \otimes A'' \right)$ induces $\theta_{ij}^*:U_{ij}\times \textnormal{Spec}\left( A^* \right) \to U_{ij}\times \textnormal{Spec}\left(A^*\right)$ with $f_{ik}^*=f_{ij}^* f_{jk}^*$. Then $\mathcal{X}^*$ is flat over $A^*$. On the other hand, the fibered sum $T_{\mathcal{F}_i^*} := T_{\mathcal{F}_i'} \times_{T_{\mathcal{F}_i}} T_{\mathcal{F}_i''}$ is generated by $T_i^{* \alpha}:=\left( T_i'^\alpha, T_i''^\alpha \right)\in \left(\Gamma\left( U_i, T_X \right)\otimes A' \right) \times_{\Gamma\left(U_i, T_X \right)\otimes A}  \left( \Gamma\left( U_i, T_X \right)\otimes A'' \right)=\Gamma\left( U_i,T_X\right)\otimes A^*$ and by setting $r_{ij}^{*\alpha\beta}:=\left( r_{ij}'^{\alpha\beta}, r_{ij}''^{\alpha\beta} \right)$, we have
\begin{align*}
\left( \theta_{ji}^*\right)^*\left( T_i^{*\alpha} \right)= \left(\theta_{ji}'^*\left( T_i'^\alpha \right), \theta_{ji}''^* \left( T_i''^\alpha \right) \right) = \left( \sum_{\beta=1}^p  r_{ij}'^{\alpha\beta} T_j'^\beta, \sum_{\beta=1}^p r_{ij}''^{\alpha\beta} T_j''^\beta \right) = \sum_{\beta=1}^p \left(r_{ij}'^{\alpha\beta}, r_{ij}''^{\alpha\beta} \right)\left( T_j'^\beta, T_j''^\beta \right)= \sum_{\beta=1}^p r_{ij}^{*\alpha\beta} T_j^{*\beta}
\end{align*}
and the integrability condition holds: by setting $g_{i\alpha\beta}^{*\gamma}:= \left( g_{i\alpha\beta}'^\gamma, g_{i\alpha\beta}''^\gamma \right)\in \Gamma\left( U_i, \mathcal{O}_X \right)\otimes A^*$,
\begin{align*}
\left[ T_i^{*\alpha}, T_i^{*\beta} \right] & = \left[ \left( T_i'^\alpha, T_i''^\alpha \right), \left( T_i'^\beta, T_i''^\beta \right)  \right]=\left( \left[ T_i'^\alpha, T_i'^\beta\right], \left[ T_i''^\alpha, T_i''^\beta   \right] \right) = \left( \sum_{\gamma=1}^p g_{i\alpha\beta}'^\gamma T_i'^\gamma, \sum_{\gamma=1}^p g_{i\alpha\beta}''^\gamma T_i''^\gamma \right) \\
&= \sum_{\gamma=1}^p \left( g_{i\alpha\beta}'^\gamma , g_{i \alpha\beta}''^\gamma \right) \left( T_i'^\gamma, T_i''^\gamma \right)= \sum_{\gamma=1}^p g_{i\alpha\beta}^{*\gamma} T_i^{*\gamma}
\end{align*}
Then $\left\{ T_{\mathcal{F}_i^*} \right\}$ defines a locally free subsheaf $T_{\mathcal{F}^*}$ of $T_{\mathcal{X}^*/A^*}$. Since $T_{\mathcal{F}^*}\hookrightarrow T_{\mathcal{X}^*/A^*}$ induces an injection $T_{\mathcal{F}_0} \hookrightarrow T_X $, we see that $\frac{T_{\mathcal{X}^*/A^*}}{T_{\mathcal{F}^*}}$ is flat over $A^*$. This implies that $\left(\textnormal{H}_1 \right)$ is satisfied. $\left( \textnormal{H}_3 \right)$ follows from Proposition \ref{at25}. It remains to show $\left( \textnormal{H}_2 \right)$. Let $A=k$ and $A''=k[\epsilon]$ and $\left( \mathcal{X}^*, T_{\mathcal{F}^*} \right)$ as above. Let $\left(\mathcal{Y}, T_\mathcal{G} \right)$ be another foliated deformation of $\left(X, T_{\mathcal{F}_0} \right)$ over $A^*$ which induces $\left( \mathcal{X}, T_{\mathcal{F}'}\right)$ and $\left( \mathcal{X}, T_{\mathcal{F}''} \right)$. We may assume that the underlying scheme $\mathcal{Y}=\mathcal{X}^*$ (see \cite{Ser06} p.66-67). We note that we have an isomorphism $\theta_i^*: U_i \times \textnormal{Spec}\left( A^*\right) \to  \mathcal{X}^*|_{U_i}$ such that $\theta_j^{* -1}\theta_i^*=\theta_{ij}^*$. Let $T_{\mathcal{G}_i}^*=\theta_i^*\left( T_\mathcal{G}|_{U_i} \right)$ be the induced foliation from $T_\mathcal{G}$. Then $\Gamma\left(U_i, T_{\mathcal{G}_i} \right)$ induces $\Gamma\left( U_i, T_{\mathcal{F}_i'} \right)$ and $\Gamma\left(U_i, T_{\mathcal{F}_i''} \right)$, i.e. there is a commutative diagram. 
\begin{center}
$\begin{CD}
\Gamma\left( U_i , T_{\mathcal{G}_i} \right)\subset \Gamma\left( U_i, T_X \right)\otimes A^* @>>> \Gamma\left( U_i, T_{\mathcal{F}_i''} \right)\subset \Gamma\left( U_i, T_X \right)\otimes k[\epsilon]\\
@VVV @VVV\\
\Gamma\left( U_i, T_{\mathcal{F}_i'} \right)\otimes A \subset \Gamma\left( U_i, T_X \right)\otimes A' @>>> \Gamma\left( U_i,  T_{\mathcal{F}_0} \right)
\end{CD}$
\end{center}
This induces $\Gamma\left( U_i, T_{\mathcal{G}_i} \right)\to \Gamma\left( U_i, T_{\mathcal{F}_i^*} \right)=\Gamma\left( U_i, T_{\mathcal{F}_i'} \right)\times_{\Gamma\left( U_i, T_{\mathcal{F}_0}\right)}\left( U_i, T_{\mathcal{F}_i''} \right)$ which is an isomorphism since $T_{\mathcal{F}_i^*}$ is free, and $T_{\mathcal{G}_i}$ is flat over $A^*$ and induces the identity $\Gamma\left( U_i, T_{\mathcal{F}_0} \right)$. This implies that $\Gamma\left( U_i, T_{\mathcal{G}_i}\right)= \Gamma\left( U_i, T_{\mathcal{F}_i^*} \right)$ in $\Gamma\left( U_i, T_X  \right)\otimes A^*$. Hence $T_\mathcal{G}= T_{\mathcal{F}^*}$ in $T_{\mathcal{X}^* /A^*}$. This completes the proof of Proposition \ref{e3}.

\end{proof}

\begin{theorem}\label{f2}
Let $\left(X, T_{\mathcal{F}_0}\right)$ be a foliated nonsingular variety with $T_{\mathcal{F}_0}$ locally free, and $\mathbb{H}^0\left(X, T_{\mathcal{F}_0}^\bullet \right)=0$. Then the functor $\textnormal{FDef}_{(X, T_{\mathcal{F}_0})}$ is pro-representable.
\end{theorem}

\begin{proof}
We check Schlessinger's criterion for $\textnormal{FDef}_{\left( X, T_{\mathcal{F}_0} \right)}$ (see \cite{Har10} Theorem 16.2 p.113). $(\textnormal{H}_0),(\textnormal{H}_1),(\textnormal{H}_2)$ and $(\textnormal{H}_3)$ follow from Proposition \ref{e3} and $(\textnormal{H}_4)$ follows from Proposition \ref{e4}.
\end{proof}

\subsection{Stability of a foliated nonsingular variety $\left( X , T_{\mathcal{F}_0} \right)$ with $T_{\mathcal{F}_0}$ locally free}\

\begin{definition}
A foliated nonsingular variety $\left(X, T_{\mathcal{F}_0} \right)$ with $T_{\mathcal{F}_0}$ is called stable if the morphism of functors from $(\ref{e1})$ by forgetting foliated structures
\begin{align*}
\Phi: \textnormal{FDef}_{\left( X, T_{\mathcal{F}_0} \right)}\to \textnormal{Def}_X
\end{align*}
is smooth.
\end{definition}

We note that $\textnormal{FDef}_{\left(X, T_{\mathcal{F}_0} \right)}$ has a semiuniversal formal element from Proposition \ref{e3} and obstruction space $\mathbb{H}^2\left(X, T_{\mathcal{F}_0}^\bullet \right)$ from Proposition \ref{at25}. On the other hand, $\textnormal{Def}_{X}$ has a semiuniversal formal element and obstruction space $\mathbb{H}^2\left( X, T_X\right)$ (see \cite{Ser06}). Assume $\mathbb{H}^1\left(X, \mathscr{H}om_{\mathcal{O}_X}\left( T_{\mathcal{F}_0}, \frac{T_X}{T_{\mathcal{F}_0}} \right)^\bullet \right)=0$. Then from the long exact sequence of $(\ref{e2})$, we have the surjection $d\Phi: \mathbb{H}^1\left(X, T_{\mathcal{F}_0}^\bullet \right)\to \mathbb{H}^1\left( X, T_X \right)$ and an injective obstruction map $\mathbb{H}^2\left( X, T_{\mathcal{F}_0}^\bullet \right)\to \mathbb{H}^2\left( X, T_X \right)$, in other words, $\textnormal{FDef}_{\left( X, T_{\mathcal{F}_0} \right)}$ is less obstructed than $\textnormal{Def}_X$ (see \cite{Ser06} Definition 2.3.5 p.58). This implies that $\Phi$ is smooth. (see \cite{Ser06} Proposition 2.3.6 p.59). However we will directly show the smoothness of $\Phi$.

\begin{proposition}
Let $\left( X, T_{\mathcal{F}_0} \right)$ be a foliated nonsingular variety with $T_{\mathcal{F}_0}$ locally free. Assume that $\mathbb{H}^1\left(X, \mathscr{H}om_{\mathcal{O}_X}\left( T_{\mathcal{F}_0}, \frac{T_X}{T_{\mathcal{F}_0}} \right)^\bullet \right)=0$. Then $\Phi$ is smooth. In other words, for a flat deformation $\mathcal{X}$ of $X$ over $A$, there exists a foliated deformation $\left(\mathcal{X}, T_\mathcal{F} \right)$ over $A$ which induces $\mathcal{X}$.
\end{proposition}

\begin{proof}
Let $0\to (t) \to \tilde{A} \to A \to 0$ be a small extension. We will show that the natural map
\begin{align*}
\textnormal{FDef}_{\left( X, T_{\mathcal{F}_0} \right)} \left( \tilde{A} \right) \to \textnormal{FDef}_{\left( X, T_{\mathcal{F}_0}\right)}\left( A \right)\times_{\textnormal{Def}_X \left(A \right)} \textnormal{Def}_{X} \left( \tilde{A} \right)
\end{align*}
is surjective. Let $\tilde{\mathcal{X}}$ be a flat deformation of $X$ over $\tilde{A}$, and $\left( \mathcal{X}, T_\mathcal{F}\right)$ be a foliated deformation of $\left(X, T_{\mathcal{F}_0}\right)$ over $A$ such that $\tilde{\mathcal{X}}$ induces $\mathcal{X}$. Let $\mathcal{U}=\left\{ U_i \right\}$ be an affine open covering of $X$ such that we have an isomorphism $\tilde{\theta}_i: U_i \times \textnormal{Spec}\left( \tilde{A} \right)\to \tilde{\mathcal{X}}|_{U_i}$ and $\tilde{\theta}_{ij}:=\tilde{\theta}_j^{-1}\tilde{\theta}_i: U_{ij}\times \textnormal{Spec}\left(\tilde{A} \right) \to  U_{ij}\times \textnormal{Spec}\left( \tilde{A} \right)$ corresponds to $\tilde{f}_{ij}: \Gamma\left(U_{ij}, \mathcal{O}_X \right)\otimes \tilde{A} \to \Gamma\left( U_{ij}, \mathcal{O}_X \right)\otimes \tilde{A}$, and we have a foliated isomorphism $\theta_i: \left( U_i\times \textnormal{Spec}\left( A\right), T_{\mathcal{F}_i}\right) \to \left( \mathcal{X}|_{U_i} , T_\mathcal{F}|_{U_i} \right) $ where $\theta_i$ is induced from $\tilde{\theta}_i$ and $T_{\mathcal{F}_i}\subset T_X|_{U_i}\otimes A$ is induced from $T_\mathcal{F}$ via $\theta_i$, and $\theta_{ij}=\theta_j^{-1}\theta_i$ induced from $\tilde{\theta}_{ij}$ corresponds to $f_{ij}:\Gamma\left( U_{ij}, \mathcal{O}_X \right) \otimes A \to \Gamma\left( U_{ij}, \mathcal{O}_X \right) \otimes A $ induced from $\tilde{f}_{ij}$, and $\Gamma\left( U_i, T_{\mathcal{F}_i  }\right)$ is generated by $T_i^\alpha \in \Gamma\left( U_i, T_X \right)\otimes A$ which induces $T_{0i}^\alpha$ modulo $\mathfrak{m}_A$ for $\alpha=1,...,p$, and $\theta_{ji}^* T_i^\alpha= \sum_{\beta=1}^p r_{ij}^{\alpha\beta} T_j^\beta$ for some $r_{ij}^{\alpha\beta}\in \Gamma\left( U_{ij}, \mathcal{O}_X \right)$ which induces $r_{0ij}^{\alpha\beta}\in \Gamma\left( U_{ij}, \mathcal{O}_X \right)$ by modulo $\mathfrak{m}_A$ and $\left[ T_i^\alpha, T_i^\beta \right]= \sum_{\gamma=1}^p g_{i\alpha\beta}^\gamma T_i^\gamma$ for some $g_{i\alpha\beta}^\gamma \in \Gamma\left( U_i, \mathcal{O}_X \right)\otimes A$ which induces $g_{0i\alpha\beta}^\gamma \in \Gamma\left( U_i, \mathcal{O}_X \right)$ by modulo $\mathfrak{m}_A$. Let $\tilde{T}_i^\alpha\in \Gamma\left( U_i, T_X \right)\otimes \tilde{A}$ be an arbitrary lifting of $T_i^\alpha$ and $\tilde{r}_{ij}^{\alpha\beta} \in \Gamma\left( U_{ij}, \mathcal{O}_X \right)\otimes \tilde{A}$ be an arbitrary lifting of $r_{ij}^{\alpha\beta}$ and $\tilde{g}_{i\alpha\beta}^\gamma\in \Gamma\left( U_i, \mathcal{O}_X \right)\otimes \tilde{A}$ be an arbitrary lifting of $g_{i\alpha\beta}^\gamma$. Then as in the proof of Proposition \ref{at25} (we keep the notation of the proof of Proposition \ref{at25}), the obstruction $\left(\left\{P_i \right\}, \left\{- W_{ij} \right\}, \left\{ d_{ijk}=0 \right\} \right)\in C^0\left(\mathcal{U}, \mathscr{H}_{\mathcal{O}_X}\left( \bigwedge^2 T_{\mathcal{F}_0} , \frac{T_X}{T_{\mathcal{F}_0}} \right) \right)\bigoplus C^1\left(\mathcal{U}, \mathscr{H}om_{\mathcal{O}_X}\left( T_{\mathcal{F}_0}, \frac{T_X}{T_{\mathcal{F}_0}} \right) \right)$ lies in $\mathbb{H}^1\left( X, \mathscr{H}om_{\mathcal{O}_X}\left( T_{\mathcal{F}_0}, \frac{T_X}{T_{\mathcal{F}_0} }\right)^\bullet \right)=0$, so that the lifting $\left( \tilde{\mathcal{X}}, T_{\tilde{\mathcal{F}}} \right)$ over $\tilde{A}$ of $\left( X, T_\mathcal{F} \right)$ exits.
\end{proof}

\subsection{Formal foliated deformations of $\left(X, T_{\mathcal{F}_0} \right)$ with $T_{\mathcal{F}_0}$ locally free}\

We only touch upon formal foliated deformations of a foliated nonsingular variety $\left( X, T_{\mathcal{F}_0} \right)$ with $T_{\mathcal{F}_0}$ locally free in a relation to Part I.  The following Theorem \ref{f1} can be considered as the scheme-theoretic version of Theorem of existence for deformations of foliated complex analytic structures in terms of tangent sheaves in Part I.

\begin{theorem}[compare \cite{Ser06} Corollary 2.6.4 p.92] \label{f1}
Let $\left(X , T_{\mathcal{F}_0} \right)$ is a foliated nonsingular projective variety with $T_{\mathcal{F}_0}$ locally free. If $\mathbb{H}^0\left( X, T_{\mathcal{F}_0}^\bullet \right)=\mathbb{H}^2\left( X, T_{\mathcal{F}_0}^\bullet \right)=0$, then $\textnormal{FDef}_{\left( X, T_{\mathcal{F}_0}\right)}$ is pro-represented by a formal power series ring.
\end{theorem}

\begin{proof}
By Theorem \ref{f2}, $\textnormal{FDef}_{\left(X, T_{\mathcal{F}_0} \right)}$ is pro-representable. Let $\textnormal{FDef}_{\left( X, T_{\mathcal{F}_0} \right)}$ be pro-represented by a universal formal couple $\left(R, \hat{u} \right)$ (see \cite{Ser06} p.47).  By Proposition \ref{at25}, $\textnormal{FDef}_{\left(X, T_{\mathcal{F}_0} \right)}$ has an obstruction space $\mathbb{H}^2\left( X, T_{\mathcal{F}_0}^\bullet \right)$ (see \cite{Ser06} Definition 2.2.9 p.53). By Proposition \ref{at25}, we have $\dim_k \mathbb{H}^1\left( X,T_{\mathcal{F}_0}^\bullet \right)\geq \mu(R) \geq  \dim_k \mathbb{H}^1\left( X, T_{\mathcal{F}_0}^\bullet \right)- \dim_k \mathbb{H}^2\left( X , T_{\mathcal{F}_0}^\bullet \right)$, where $\mu(R)$ is the Krull dimension of $R$ (see \cite{Ser06} Corollary 2.2.11 p.54). By the hypothesis $\mathbb{H}^2\left(X, T_{\mathcal{F}_0}^\bullet \right)=0$, we see that $\textnormal{FDef}_{\left(X, T_{\mathcal{F}_0} \right)}$ is smooth, and we have $\mu(R)=\dim_k \mathbb{H}^1\left( X, T_{\mathcal{F}_0}^\bullet \right)=m$. Hence we have $R\cong k\left[ \left[ X_1,..., X_m \right] \right]$.
\end{proof}

We can consider a concept of formal foliated deformations of a foliated nonsingular variety $\left(X, T_{\mathcal{F}_0 } \right)$ with $T_{\mathcal{F}_0}$ locally free as an analogue of formal deformations of a nonsingular variety $X$ as in \cite{Ser06} p.77. We note that we can extend the functor $\textnormal{FDef}_{\left(X, T_{\mathcal{F}_0} \right)}$ to the functor
\begin{align*}
\widehat{\textnormal{FDef}}_{\left( X, T_{\mathcal{F}_0} \right)} : \widehat{\bold{Art}}&\to \left(\textnormal{Sets} \right)\\
     \bar{A} &\to \varprojlim \textnormal{FDef}_{\left(X, T_{\mathcal{F}_0} \right)}\left(\bar{A}/\mathfrak{m}_{\bar{A}}^{n+1} \right)
\end{align*}
(see \cite{Ser06} p.46). An element $\hat{\eta}\in \widehat{\textnormal{FDef}}_{\left( X, T_{\mathcal{F}_0} \right)}\left(\bar{A} \right)$ is called a formal element of $\textnormal{FDef}_{\left( X, T_{\mathcal{F}_0} \right)}$. By definition $\hat{\eta}$ can be represented as a sysemt of elements $\left\{ \eta_n \in \textnormal{FDef}_{\left( X, T_{\mathcal{F}_0} \right)}\left(\bar{A}/\mathfrak{m}_{\bar{A}}^{n+1} \right) \right\}_{n\geq 0}$  such that for every $n\geq 1$ the map $\textnormal{FDef}_{\left( X, T_{\mathcal{F}_0}\right)}\left( \bar{A}/\mathfrak{m}_{\bar{A}}^{n+1}\right)\to \textnormal{FDef}_{\left( X, T_{\mathcal{F}_0} \right)}\left( \bar{A}/\mathfrak{m}_{\bar{A}}^n \right)$ induced by the projection $\bar{A}/\mathfrak{m}_{\bar{A}}^{n+1}\to \bar{A}/{\mathfrak{m}}_{\bar{A}}^n$ sends $\eta_n\mapsto \eta_{n-1}$. A formal foliated deformation of $\left( X, T_{\mathcal{F}_0}\right)$ over $\bar{A}$ is a formal element $\hat{\eta}=\left\{\eta_n \right\} \in \widehat{\textnormal{FDef}}_{\left( X, T_{\mathcal{F}_0} \right)}\left( \bar{A}\right)$ for $\textnormal{FDef}_{\left(X, T_{\mathcal{F}_0} \right)}$, where $\eta_n$ is an infinitesimal foliated deformation of $\left( X, T_{\mathcal{F}_0}\right)$ over $\bar{A}_n$
\begin{center}
$\eta_n:\begin{CD}
\left(X, T_{\mathcal{F}_0} \right) @>>> \left( \mathcal{X}_n, T_{\mathcal{F}_n} \right) \\
@VVV @VV \pi_n V\\
\textnormal{Spec}(k) @>>> \textnormal{Spec}\left( \bar{A}_n \right)
\end{CD}$
\end{center}
where $\bar{A}_n=\bar{A}/\mathfrak{m}_{\bar{A}}^{n+1}$, such that for all $n\geq 1$, $\eta_n$ induces $\eta_{n-1}$ by pullback under the natural inclusion $\textnormal{Spec}\left(\bar{A}_{n-1} \right)\to \textnormal{Spec}\left( \bar{A}_n \right)$, i.e. we have
\begin{center}
$\eta_{n-1}:\begin{CD}
\left(X, T_{\mathcal{F}_0} \right) @>>> \left( \mathcal{X}_{n-1}=\mathcal{X}_n \otimes_{\bar{A}_n} \bar{A}_{n-1}, T_{\mathcal{F}_n}\otimes_{\bar{A}_n} \bar{A}_{n-1}  \right)\\
@VVV @VV\pi_{n-1}V\\
\textnormal{Spec}(k) @>>> \textnormal{Spec}\left(\bar{A}_{n-1} \right)
\end{CD}$
\end{center}
We will denote such a formal foliated deformation by $\left( \bar{A}, \left\{ \eta_n \right\} \right)$.

Let $\left(X, T_{\mathcal{F}_0}\right)$ be a foliated nonsingular projective variety with $T_{\mathcal{F}_0}$ locally free and consider a flat family of foliated deformations of $\left( X, T_{\mathcal{F}_0} \right)$ parametrized by an affine scheme $S=\textnormal{Spec}\left( B \right)$ where $B$ is a $k$-algebras of finite type:
\begin{center}
$\eta:\begin{CD}
\left( X, T_{\mathcal{F}_0} \right) @>f>> \left( \mathcal{X}, T_\mathcal{F} \right)\\
@VVV @VV\pi V\\
\textnormal{Spec}(k) @>s >> S
\end{CD}$
\end{center}
where $\pi$ is projective, and $\mathcal{X}$ is smooth over $S$, $s$ is a $k$-rational point and $T_\mathcal{F}$ is a locally free subsheaf of $T_{\mathcal{X}/S}$ such that $\frac{T_{\mathcal{X}/S}}{T_\mathcal{F}}$ is flat over $S$ and $\left[ T_{\mathcal{F}}, T_{\mathcal{F}}\right] \subset T_{\mathcal{F}}$ under the Lie bracket. We will denote the foliated deformation $\eta$ by $\left( S, s , \eta \right)$. Given such a foliated deformation $\left(  S, s, \eta \right)$ let $\eta_n$ be the infinitesimal foliated deformation induced pulling back of $\eta$ under the closed embedding $\textnormal{Spec}\left( \mathcal{O}_{S,s}/\mathfrak{m}^{n+1} \right) \to S$. Then $\left( \widehat{\mathcal{O}}_{S,s}, \left\{ \eta_n \right\} \right)$ is a formal foliated deformation of $\left(X, T_{\mathcal{F}_0} \right)$ which we call the formal foliated deformation defined by $\eta$.

\begin{definition}
A foliated deformation $\left(S, s, \eta \right)$ of a foliated nonsingular variety $\left( X, T_{\mathcal{F}_0} \right)$ with $T_{\mathcal{F}_0}$ locally free is called formally universal $($formally semiuniversal, formally versal$)$ if the formal foliated deformation $\left( \widehat{\mathcal{O}}_{S,s}, \left\{ \eta_n\right\}_{n \geq 0} \right)$ defined by $\eta$ is universal $($resp. semiuniversal, versal$)$ \textnormal{(see \cite{Ser06} Definition 2.2.3 p.47 and Definition 2.2.6 p.49)}. A flat foliated family $\pi:\left(\mathcal{X}, T_\mathcal{F} \right)\to S$, where $\mathcal{X}$ is smooth over $S$, and $T_{\mathcal{F}}$ is a locally free subsehaf of $T_{\mathcal{X}/S}$, is called formally universal $($resp. formally semiuniversal, formally versal$)$ at a $k$-rational point $s\in S$ if
\begin{center}
$\eta:\begin{CD}
\left(\mathcal{X}(s), T_{\mathcal{F}(s)}\right) @>\hookrightarrow>> \left( \mathcal{X}, T_\mathcal{F} \right)\\
@VVV @VV\pi V \\
\textnormal{Spec}(k) @>s>> S
\end{CD}$
\end{center}
is a formally universal $($resp. formally semiuniversal, formally versal$)$ foliated deformation of $\left( \mathcal{X}(s), T_{\mathcal{F}(s)}\right)$.
\end{definition}

\begin{definition}[foliated Kodaira-Spencer map in an algebraic foliated family]
Let $\left( X, T_{\mathcal{F}_0} \right)$ be a foliated nonsingular variety with $T_{\mathcal{F}_0}$ locally free. Let us consider a foliated deformation $(S, s, \eta)$
\begin{center}
$\eta: \begin{CD}
\left( X, T_{\mathcal{F}_0} \right) @>>> \left( \mathcal{X}, T_\mathcal{F} \right)\\
@VVV @VV\pi V\\
\textnormal{Spec}(k) @>s>> S
\end{CD}$
\end{center}
By pulling back $(S, s, \eta)$ by morphisms $\textnormal{Spec}\left( k[\epsilon] \right)\to S$ with image $s$ and applying \textnormal{Proposition \ref{at25}} we define a linear map \textnormal{(see \cite{MO15} p.200)}
\begin{align*}
\kappa_{\pi, s}: T_{S,s}\to  \textnormal{FDef}_{\left( X, T_{\mathcal{F}_0} \right)} \cong \mathbb{H}^1\left( X, T_{\mathcal{F}_0}^\bullet \right)
\end{align*}
which is called the foliated Kodaira-Spencer map of $(S,s, \eta)$ in terms of tangent sheaves.
\end{definition}

\begin{proposition}[compare \cite{Ser06} Proposition 2.5.8 p.79]\label{f3}
Let $\left(S , s , \eta \right)$ be a foliated deformation of a foliated nonsingular variety $\left( X, T_{\mathcal{F}_0} \right)$ with $T_{\mathcal{F}_0}$ locally free.
\begin{enumerate}
\item If $\eta$ is formally versal $($resp. formally semiuniversal or formally universal$)$, then the foliated Kodaira-Spencer map $\kappa_{\pi, s}: T_s S\to \textnormal{FDef}_{\left(X, T_{\mathcal{F}_0} \right)}\left( k[\epsilon]\right)\cong \mathbb{H}^1 \left(X, T_{\mathcal{F}_0}^\bullet \right)$ is surjective $($resp. an isomorphism$)$.
\item If $S$ is nonsingular at $s$ and the foliated Kodaira-Spencer map $\kappa_{\pi,s}$ is surjective $($resp. an isomorphism$)$, then $\eta$ is formally versal $($resp. formally semiuniversal$)$ and $\left(X, T_{\mathcal{F}_0} \right)$ is unobstructed, i.e. $\textnormal{FDef}_{\left( X, T_{\mathcal{F}_0} \right)}$ is smooth.
\end{enumerate}
\end{proposition}
\begin{proof}
Let $h_{\widehat{\mathcal{O}}_{S,s}}: \widehat{\bold{Art}}\to \left( \textnormal{Sets}\right)$ be a functor defined by $\bar{A} \to \textnormal{Hom}_{\widehat{\bold{Art}}}\left(\bar{A},  \widehat{\mathcal{O}}_{S,s} \right) $. Then there is a $1-1$ correspondence between $\widehat{\textnormal{FDef}}_{\left(X, T_{\mathcal{F}_0} \right)}\left(\widehat{\mathcal{O}}_{S, s} \right)$ and the set of morphisms of functors $h_{\widehat{\mathcal{O}}_{S,s}}\to \textnormal{FDef}_{\left(X, T_{\mathcal{F}_0} \right)}$ (see \cite{Ser06} Lemma 2.2.2 p.46). So $\left(\widehat{\mathcal{O}}_{S,s}\, \left\{\eta_n \right\} \right)$ defines a morphism $h_{\widehat{\mathcal{O}}_{S,s}}\to \textnormal{FDef}_{\left(X, T_{\mathcal{F}_0} \right)}$. Then $(1)$ follows from the definitions of versality and semiuniversality of the formal couple $\left(\widehat{\mathcal{O}}_{S,s}, \left\{\eta_n \right\} \right)$ applied to $h_{\widehat{\mathcal{O}}_{S,s}}\to \textnormal{FDef}_{\left(X, T_{\mathcal{F}_0} \right)}$. Let us prove $(2)$. Since $h_{\widehat{\mathcal{O}}_{S,s}}$ is smooth, i.e unobstructed, and $\textnormal{FDef}_{\left(X, T_{\mathcal{F}_0} \right)}$ has a semiuniversal formal element by Proposition \ref{e3} and $\textnormal{FDef}_{\left( X, T_{\mathcal{F}_0} \right)}$ has an obstruction space $\mathbb{H}^2\left( X, T_{\mathcal{F}_0}^\bullet \right)$ by Proposition \ref{at25}, and the foliated Kodaira-Spencer map $\kappa_{S,s}$ is surjective by the assumption, the morphism $h_{\widehat{\mathcal{O}}_{S,s}}\to \textnormal{FDef}_{\left( X, T_{\mathcal{F}_0} \right)}$ is smooth and $\textnormal{FDef}_{\left( X , T_{\mathcal{F}_0} \right)}$ is unobstructed. (see \cite{Ser06} Corollary 2.3.7 p.59).
\end{proof}

Proposition \ref{f3} is the algebraic version of Theorem of completeness for deformations of foliated complex analytic structures in terms of tangent sheaves in Part I. We finish this section by mentioning the concept of algebraizability in foliated deformations.

\begin{definition}[compare \cite{Ser06} p.80]
A formal foliated deformation $\left( \bar{A}, \left\{ \eta_n \right\} \right)$ of a foliated nonsingular variety $\left( X, T_{\mathcal{F}_0}\right)$ with $T_{\mathcal{F}_0}$ locally free is called algebraizable if there exists an algebraic foliated deformation $\left( S, s, \xi \right)$ of $\left( X, T_{\mathcal{F}_0} \right)$ and an isomorphism $\widehat{\mathcal{O}}_{S,s}\cong \bar{A}$ sending $\eta_n\to \xi_n$ for all $n$, i.e. $\left( \bar{A}, \left\{ \eta_n \right\} \right)$ is isomorphic to the formal foliated deformation defined by $\xi$. The foliated deformation $\left( S, s, \xi \right)$ is called an algebraization of $\left( \bar{A} , \left\{ \eta_n \right\}  \right)$.
\end{definition}

The existence of algebraizations can be considered as the counterpart of the convergence step in the construction of local families of foliated deformations of compact foliated complex manifolds in terms of tangent sheaves in Part I.

\section{Generalities on foliated deformations in terms of cotangent sheaves}

Let $X$ be a nonsingular variety over $k=\mathbb{C}$ and $\mathcal{N}_{\mathcal{F}_0}^*$ be a singular foliation on $X$. In other words,  we have an exact sequence
\begin{align*}
0\to \mathcal{N}_{\mathcal{F}_0}^* \to \Omega_X^1 \to \Omega_X^1/ \mathcal{N}_{\mathcal{F}_0}^* \to 0
\end{align*}
where $\mathcal{N}_{\mathcal{F}_0}^*$ is a coherent subsheaf of the cotangent sheaf $\Omega_X^1$ such that $d\mathcal{N}_{\mathcal{F}_0,x}^* \subset \mathcal{N}_{\mathcal{F}_0, x}^*\bigwedge \Omega_{X,x}^1$ for $x\in X- \textnormal{Sing}\left(\mathcal{F}_0 \right)$ where $\textnormal{Sing}\left( \mathcal{F}_0\right)$ is the singular locus of $\mathcal{F}_0$, and $\Omega_X^1/ \mathcal{N}_{\mathcal{F}_0}^*$ is a torsion-free $\mathcal{O}_X$-module, so that $\mathcal{N}_{\mathcal{F}_0}^*$ is reflexive.  A cartesian diagram of morphism of algebraic schemes
\begin{center}
$\eta:\begin{CD}
\left(X, \mathcal{N}_{\mathcal{F}_0}^* \right)@>>> \left(\mathcal{X}, \mathcal{N}_\mathcal{F}^* \right)\\
@VVV @VV\pi V\\
\textnormal{Spec}(k)@>s>> S
\end{CD}$
\end{center}
where $\pi$ is flat, smooth and surjective, and $S$ is a connected algebraic scheme, is called a family of (singularly) foliated deformations of $\left(X, \mathcal{N}_{\mathcal{F}_0}^* \right)$ in terms of cotangent sheaves parametrized by $S$ if there exists an exact sequence
\begin{align}\label{b1}
0\to  \mathcal{N}_\mathcal{F}^* \to \Omega_{\mathcal{X}/S}^1 \to \Omega_{\mathcal{X}/S}^1 /\mathcal{N}_\mathcal{F}^* \to 0
\end{align}
where  $d_{\mathcal{X}/S} \left(\mathcal{N}_{\mathcal{F},x}^*\right) \subset \mathcal{N}_{\mathcal{F},x}^* \bigwedge \Omega_{\mathcal{X}/S,x}^1 $ for $x\in \mathcal{X} - \textnormal{Sing}\left( \mathcal{F}\right)$ where $\textnormal{Sing}\left(\mathcal{F}\right)$ is the locus of points of $\mathcal{X}$ at which $\frac{ \Omega_{\mathcal{X}/S}^1 }{\mathcal{N}_\mathcal{F}^*}$ is not locally free, and $\Omega_{\mathcal{X}/S}^1/\mathcal{N}_\mathcal{F}^*$ is flat over $S$, and for each $k$-rational point $t\in S$ the exact sequence $(\ref{a1})$ induces an exact sequence on the fiber $X_t$
\begin{align*}
0\to  \mathcal{N}_{\mathcal{F}_t}^*\to \Omega_{X_t/k}^1 \to \Omega_{X_t/k}^1/ \mathcal{N}_{\mathcal{F}_t}^* \to 0
\end{align*}
such that $\Omega_{X_t}^1/ \mathcal{N}_{\mathcal{F}_t}^*$ is a torsion-free $\mathcal{O}_{X_t}$-module. When $S=\textnormal{Spec}(A)$ with $A\in \textnormal{ob}\left(\bold{Art}^*\right)$ and $s\in S$ is the closed point we have a local family of foliated deformations of $\left(X, \mathcal{N}_{\mathcal{F}}^* \right)$ in terms of cotangent sheaves over $A$. The foliated deformation $\eta$ in terms of cotangent sheaves will be denoted by $\left( S, \eta\right)$ or $\left( A, \eta\right)$ when $S=\textnormal{Spec}(A)$. The local foliated deformation $\left( A, \eta\right)$ is infinitesimal (resp. first-order) if $A\in \bold{Art}$ (resp. $A=k[\epsilon]$). Given another foliated deformation of $\left( X, \mathcal{N}_{\mathcal{F}_0}^*\right)$ over $S$ in terms of cotangent sheaves
\begin{center}
$\xi: \begin{CD}
\left( X, \mathcal{N}_{\mathcal{F}_0}^* \right) @>>>  \left( \mathcal{Y}, \mathcal{N}_\mathcal{G}^*\right)\\
@VVV @VVV\\
\textnormal{Spec}(k) @>>> S
\end{CD}$
\end{center}
An isomorphism between foliated deformations $\eta$ and $\xi$ is an foliated $S$-isomorphism $\phi:\left(\mathcal{X}, \mathcal{N}_{\mathcal{F}}^* \right)\to \left(  \mathcal{Y}, \mathcal{N}_{\mathcal{G}}^* \right)$ with $\mathcal{N}_\mathcal{F}^*=\phi^* \mathcal{N}_{\mathcal{G}}^*$ as a subsheaf of $\Omega_{\mathcal{X}/S}^1$ inducing the identity on $\left( X, \mathcal{N}_{\mathcal{F}_0}^*\right)$, or equivalently the commutative diagram
\begin{center}
$\begin{CD}
0 @>>> \mathcal{N}_\mathcal{F}^* @>>> \Omega_{\mathcal{X}/S}^1 \\
@. @V\cong VV @| \\
0@>>> \phi^* \mathcal{N}_\mathcal{G} @>>> \phi^* \Omega_{\mathcal{Y}/S}^1=\Omega_{\mathcal{X}/S}^1
\end{CD}$
\end{center}

\begin{definition}[Trivial foliated deformations in terms of cotangent sheaves]
Let $\left( X, \mathcal{N}_{\mathcal{F}_0}^* \right)$ be a foliated nonsingular variety over $k$ and $(S,s)$ be a $k$-pointed scheme. We define the trivial foliated family induced by $\left( X, \mathcal{N}_{\mathcal{F}_0}^* \right)$ and $(S,s)$ to be the following foliated deformation of $\left( X, \mathcal{N}_{\mathcal{F}_0}^* \right)$ over $S$,
\begin{center}
$\begin{CD}
\left( X, \mathcal{N}_{\mathcal{F}_0}^* \right) @>>> \left(  X\times_k S    , \pi_1^* \mathcal{N}_{\mathcal{F}_0}^* \right)\\
@VVV @VV\pi_2 V\\
\textnormal{Spec}(k) @>s>> S
\end{CD}$
\end{center}
where $\pi_1:X\times S\to X$ and $\pi_2:X\times S \to S$ are the projections. A flat foliated deformation of $\left( X, \mathcal{N}_{\mathcal{F}_0}^* \right)$ over $S$ is called trivial if it is isomorphic to the trivial foliated family over $S$.
\end{definition}

\begin{definition}[Rigid foliated deformations in terms of cotangent sheaves]
Let $\left( X, \mathcal{N}_{\mathcal{F}_0}^* \right)$ be a foliated nonsingular variety. $\left(X, \mathcal{N}_{\mathcal{F}_0}^*\right)$ is called rigid in foliated deformations in terms of cotangent sheaves if every infinitesimal foliated deformation of $\left( X, \mathcal{N}_{\mathcal{F}_0}^*\right)$ over $A$ is trivial for every $A \in \bold{Art}$.

\end{definition}

\section{Foliated deformations of $\left(X, \mathcal{N}_{\mathcal{F}_0}^* \right)$ with $\mathcal{N}_{\mathcal{F}_0}^*$ locally free such that the underlying deformations of $\left(X, \mathcal{N}_{\mathcal{F}_0}^* \right)$ are trivial}

Let $\left(X, \mathcal{N}_{\mathcal{F}_0}^*\right)$ be a foliated nonsingular variety such that $\mathcal{N}_{\mathcal{F}_0}^*$ is a locally free subsheaf of $\Omega_X^1$, so that we have an exact sequence $0\to \mathcal{N}_{\mathcal{F}_0}^*\to \Omega_X^1 \to \Omega_X^1/  \mathcal{N}_{\mathcal{F}_0}^* \to 0$ where $\Omega_X^1 / \mathcal{N}_{\mathcal{F}_0}^* $ is a torsion-free $\mathcal{O}_X$-module. We start from the most simple foliated deformations in terms of cotangent sheaves. In this section, we study flat foliated deformations of $\left(X, \mathcal{N}_{\mathcal{F}_0}^* \right)$ such that the underlying deformations of a pair of a nonsingular variety $X$ and a locally free sheaf $\mathcal{N}_{\mathcal{F}_0}^*$ are trivial.

\begin{definition}
Let $(X, \mathcal{N}_{\mathcal{F}_0}^*)$ be a foliated nonsingular variety as above. An infinitesimal foliated deformations $\left( \mathcal{X}= X\times \textnormal{Spec}(A), \mathcal{N}_\mathcal{F}^* \right)$ over $A\in \bold{Art}$ of $\left(X, \mathcal{N}_{\mathcal{F}_0}^*\right)$  with $\mathcal{N}_{\mathcal{F}_0}^*$ and $X$ fixed is a locally free subsheaf $\mathcal{N}_{\mathcal{F}}^*$ of $\Omega_{\frac{X\times \textnormal{Spec}(A)}{\textnormal{Spec}(A)}}^1=\Omega_X^1\otimes A$ such that $\mathcal{N}_{\mathcal{F}}^* \cong \mathcal{N}_{\mathcal{F}_0}^*\otimes A$ and $d_{\mathcal{X}/A} \left(\mathcal{N}_{\mathcal{F},x}^* \right) \subset \mathcal{N}_{\mathcal{F},x}^* \bigwedge \Omega_{\mathcal{X}/A,x}^1$ for $x \in X-\textnormal{Sing}\left(\mathcal{F}_0 \right)$ and $\frac{\Omega_{X\times \textnormal{Spec}(A)/\textnormal{Spec}(A)}}{\mathcal{N}_\mathcal{F}^*}$ is flat over $A$. Equivalently an infinitesimal foliated deformations over $A$ of $\left(X, \mathcal{N}_{\mathcal{F}_0}^*\right)$  with $\mathcal{N}_{\mathcal{F}_0}^*$ and $X$ fixed is an exact sequence on $\mathcal{X}= X\times_k \textnormal{Spec}(A)$
\begin{align}
0\to \mathcal{N}_{\mathcal{F}_0}^* \otimes A \xrightarrow{l} \Omega_{\mathcal{X}/\textnormal{Spec}(A)}^1=\Omega_X^1\otimes A \to \frac{\Omega_X^1\otimes A}{ l\left(\mathcal{N}_{\mathcal{F}_0}^* \otimes A\right)}\to 0
\end{align}
which induces $0\to \mathcal{N}_{\mathcal{F}_0}^* \to \Omega_X^1 \to \frac{\Omega_X^1}{\mathcal{N}_{\mathcal{F}_0}^*}\to 0$ over $\textnormal{Spec}(k)$ such that $d_{\mathcal{X}/A}\left(l\left(\mathcal{N}_{\mathcal{F}_0}^*\otimes A \right)_x \right)\subset l\left( \mathcal{N}_{\mathcal{F}_0}^*\otimes A\right)_x \bigwedge \Omega_{\mathcal{X}/A,x}^1$ for $x\in X-\textnormal{Sing}\left( \mathcal{F}_0 \right)$, and $\Omega_{\mathcal{X}/A}^1/l\left(\mathcal{N}_{\mathcal{F}_0}^* \otimes A \right)$ is flat over $\textnormal{Spec}(A)$. Then we have a functor of Artin rings
\begin{align*}
\textnormal{FDef}_{\left(X, \mathcal{N}_{\mathcal{F}_0}^*\right)_{\textnormal{fixed}}}:\bold{Art}&\to \textnormal{Sets}\\
                                                                              A &\mapsto        \{\textnormal{infinitesimal foliated deformations of $\left(X, \mathcal{N}_{\mathcal{F}_0}^*\right)$ over $A$  with $\mathcal{N}_{\mathcal{F}_0}^*$ and $X$ fixed} \}
\end{align*}
\end{definition}

Then we will show that 
{\small{\begin{align*}
\mathscr{H}om_{\mathcal{O}_X}\left(\mathcal{N}_{\mathcal{F}_0}^*, \Omega_X^1 \right)^\bullet: \mathscr{H}om_{\mathcal{O}_X}\left(\mathcal{N}_{\mathcal{F}_0}^*, \Omega_X^1 \right) \xrightarrow{E_1'}\mathscr{H}om_{\mathcal{O}_X}\left( \mathcal{N}_{\mathcal{F}_0}^*, \tilde{\mathcal{S}}^2 \right) \xrightarrow{E_2} \mathscr{H}om_{\mathcal{O}_X}\left( \mathcal{N}_{\mathcal{F}_0}^*, \tilde{\mathcal{S}}^3 \right) \xrightarrow{E_3} \cdots
\end{align*}}}
controls infinitesimal foliated deformations of  $\left( X, \mathcal{N}_{\mathcal{F}_0}^* \right)$ with $X$ and $\mathcal{N}_{\mathcal{F}_0}^*$ fixed. Here $E_1'$ is defined by the composition of the natural map $\mathscr{H}om_{\mathcal{O}_X}\left( \mathcal{N}_{\mathcal{F}_0}^*, \Omega_X^1  \right)\to \mathscr{H}om_{\mathcal{O}_M}\left( \mathcal{N}_{\mathcal{F}_0}^* , \frac{\Omega_X^1}{\mathcal{N}_{\mathcal{F}_0}^*}\right)$ induced by the quotient $\Omega_X^1\to \frac{\Omega_X^1}{\mathcal{N}_{\mathcal{F}_0}^*}$ with $E_1$ (see Part I). We will denote the $i$-th cohomology group of $\mathscr{H}om_{\mathcal{O}_X}\left( \mathcal{N}_{\mathcal{F}_0}^* , \Omega_X^1 \right)^\bullet$ by $\mathbb{H}^i\left( X, \mathscr{H}om_{\mathcal{O}_X}\left( \mathcal{N}_{\mathcal{F}_0}^* , \Omega_X^1 \right)^\bullet\right)$.

\begin{proposition}\label{b2}
Let $\left( X, \mathcal{N}_{\mathcal{F}_0}^* \right)$ be a foliated nonsingular variety over $k$ with $\mathcal{N}_{\mathcal{F}_0}^*$ locally free. Then
\begin{enumerate}
\item there is a $1-1$ correspondence 
\begin{align*}
\textnormal{FDef}_{\left( X, \mathcal{N}_{\mathcal{F}_0}^*\right)_{\textnormal{fixed}}} \left(k[\epsilon] \right) \cong \mathbb{H}^0\left( X, \mathscr{H}om_{\mathcal{O}_X}\left(\mathcal{N}_{\mathcal{F}_0}^*, \Omega_X^1 \right)^\bullet \right)
\end{align*}
\item Let $A\in \bold{Art}$ and $\eta=\left( \mathcal{X}, \mathcal{N}_\mathcal{F}^* \right)$ be an infinitesimal foliated deformations over $A$ of $\left( X, \mathcal{N}_{\mathcal{F}_0}^* \right)$ with $X$ and $\mathcal{N}_{\mathcal{F}_0}^*$ fixed. To every small extension $0\to (t)\to \tilde{A} \to A \to 0$, we can associate an element $o_\eta(e)\in \mathbb{H}^1\left( X, \mathscr{H}om_{\mathcal{O}_X}\left( \mathcal{N}_{\mathcal{F}_0}^*, \Omega_X^1 \right)^\bullet \right)$ called the obstruction lifting of $\eta$ to $\tilde{A}$, which is $0$ if and only if a lifting of $\eta$ to $\tilde{A}$ exists.
\end{enumerate}
\end{proposition}

\begin{proof}
Given a first-order deformation $\left(X\times \textnormal{Spec}\left(k[\epsilon ]\right) , \mathcal{N}_{\mathcal{F}}^*\right)$ such that the underlying deformation of $\left(X, \mathcal{N}_{\mathcal{F}_0}^* \right)$ is trivial, we choose an affine open covering $\mathcal{U}=\{U_i\}$ of $X$ such that $\Gamma\left( U_i, \mathcal{N}_{\mathcal{F}_0}^* \right)$ is generated by $w_{0i}^1,..., w_{0i}^q\in \Gamma\left( U_i, \Omega_X^1 \right)$ with $w_{0i}^\alpha= \sum_{\beta=1}^q h_{0ij}^{\alpha\beta} w_{0j}^\beta$ for some $h_{0ij}^{\alpha\beta}\in \Gamma\left( U_{ij}, \mathcal{O}_X \right),\alpha=1,...,q$ and $w_{0i}^1 \wedge \cdots \wedge w_{0i}^q \wedge d w_{0i}^\alpha=0,\alpha=1,...,q$, and $\Gamma\left(U_i, \mathcal{N}_\mathcal{F}^* \right)$ is generated by $w_{0i}^\alpha + \epsilon w_i'^\alpha$ for some $w_i'^\alpha \in \Gamma\left( U_i, \Omega_X^1 \right)$ with $w_{0i}^\alpha + \epsilon w_i'^\alpha = \sum_{\beta=1}^q h_{0ij}^{\alpha\beta}\left( w_{0j}^\beta + \epsilon w_j'^\beta \right)$ and $\left( w_{0i}^1+ \epsilon w_i'^1\right) \wedge \cdots \wedge \left( w_{0i}^q + \epsilon w_i'^q \right) \wedge d \left(w_{0i}^\alpha + \epsilon w_i'^\alpha \right)=0,\alpha=1,...,q$. Then we have
\begin{align}
&w_i'^\alpha = \sum_{\beta=1}^q h_{0ij}^{\alpha\beta} w_j'^\beta \label{b3}\\
\sum_{\beta=1}^q w_{0i}^1 \wedge \cdots \wedge \overbrace{w_{i}'^\beta}^{\beta-\textnormal{th}}& \wedge \cdots  \wedge  w_{0i}^q \wedge dw_{0i}^\alpha + w_{0i}^1 \wedge \cdots \wedge w_{0i}^q \wedge d w_i'^\alpha \label{b4}
\end{align}
We define $w_i'\in \Gamma\left(U_i, \mathscr{H}om_{\mathcal{O}_X}\left( \mathcal{N}_{\mathcal{F}_0}^*, \Omega_X^1 \right) \right)$ in the following way:
\begin{align*}
w_i':\Gamma\left( U_i, \mathcal{N}_{\mathcal{F}_0}^*  \right) &\to \Gamma\left( U_i, \Omega_X^1 \right) \\ 
w_{0i}^\alpha &\mapsto w_i'^\alpha
\end{align*}
and linearly extends to $\Gamma\left(U_i, \mathcal{N}_{\mathcal{F}_0}^* \right)$. Then $(\ref{b3})$ implies that $w':=\left\{ w_i'\right\}\in H^0\left(X, \mathscr{H}om_{\mathcal{O}_X}\left( \mathcal{N}_{\mathcal{F}_0}^*, \Omega_X^1 \right) \right)$ and $(\ref{b4})$ implies that $E_1'\left(w_i' \right)=0$, so that $w'\in \mathbb{H}^0\left( X, \mathscr{H}om_{\mathcal{O}_X}\left( \mathcal{N}_{\mathcal{F}_0}^*, \Omega_X^1 \right)^\bullet \right)$.

We define an inverse map. Given a $w'\in \mathbb{H}^0\left( X, \mathscr{H}om_{\mathcal{O}_X}\left( \mathcal{N}_{\mathcal{F}_0}^*, \Omega_X^1 \right)^\bullet \right)$, we have $w_i':\Gamma\left( U_i, \mathcal{N}_{\mathcal{F}_0}^* \right) \to \Gamma\left(U_i, \Omega_X^1 \right), w_{0i}^\alpha \mapsto w_i'^\alpha$ and $\sum_{\beta=1}^q w_{0i}^1 \wedge \cdots \wedge w_i'^\beta \wedge \cdots \wedge w_i'^q \wedge dw_{0i}^\alpha + w_{0i}^1 \wedge \cdots \wedge w_{0i}^q \wedge dw_i'^\alpha=0,\alpha=1,...,q$. By reversing the above process, we obtain a first-oder foliated deformation of $\left( X, \mathcal{N}_{\mathcal{F}_0}^* \right)$ with $X$ and $\mathcal{N}_{\mathcal{F}_0}^*$ fixed.

Next we identify obstructions. Given a small extension $e:0\to (t)\to \tilde{A} \to A \to 0$ in $\bold{Art}$ and an infinitesimal foliated deformation $\eta=\left( \mathcal{X}, \mathcal{N}_\mathcal{F}^* \right)$ over $A$ of $\left( X, \mathcal{N}_{\mathcal{F}_0}^*\right)$ with $X$ and $\mathcal{N}_{\mathcal{F}_0}^*$ fixed, we choose an affine open covering $\mathcal{U}=\{U_i\}$ of $X$ such that $\Gamma\left( U_i, \mathcal{N}_\mathcal{F}^* \right)$ is generated by $w_i^1 ,..., w_i^q\in \Gamma\left( U_i, \Omega_X^1 \right) \otimes A$ which induces $w_{0i}^\alpha$ by modulo $\mathfrak{m}_A$ and $w_i^\alpha = \sum_{\beta=1}^q h_{0ij}^{\alpha\beta} w_j^\beta$ and $w_i^1\wedge \cdots \wedge w_i^q \wedge dw_i^\alpha =0,\alpha=1,...,q$. In order to see if there exists a lifting of $\left( \mathcal{X}, \mathcal{N}_\mathcal{F}^* \right)$ to $\tilde{A}$, we choose an arbitrary lifting $\tilde{w}_i^\alpha \in \Gamma\left(U_i, \Omega_X \right)\otimes \tilde{A}$. Then we have
\begin{align*}
\tilde{w}_i^\alpha - \sum_{\beta=1}^q h_{0ij}^{\alpha\beta} \tilde{w}_j^\beta = t R_{ij}^\alpha\,\,\,\,\,\,\textnormal{for some}\,\,\,\,\, R_{ij}^\alpha \in \Gamma\left( U_{ij}, \Omega_X^1 \right)\\
\tilde{w}_i^\alpha \wedge \cdots \wedge \tilde{w}_i^q \wedge  d \tilde{w}_i^\alpha= t B_i^\alpha  \,\,\,\,\,\,\,\textnormal{for some} \,\,\,\,\, B_i^\alpha\in \Gamma\left( U_i, \Omega_X^1 \right)
\end{align*}
We define $R_{ij}\in \Gamma\left(U_{ij}, \mathscr{H}om_{\mathcal{O}_X}\left( \mathcal{N}_{\mathcal{F}_0}^*, \Omega_X^1 \right) \right)$ by
\begin{align*}
R_{ij}: \Gamma\left( U_{ij}, \mathcal{N}_{\mathcal{F}_0}^* \right) &\to \Gamma\left( U_{ij}, \Omega_X^1  \right)\\
 w_{0i}^\alpha &\mapsto R_{ij}^\alpha
\end{align*}
and linearly extends to $\Gamma\left( U_{ij}, \mathcal{N}_{\mathcal{F}_0}^* \right)$. On the other hand, we show that $B_i^\alpha \in \Gamma\left( U_i, \tilde{\mathcal{S}}^2 \right),\alpha=1,...,q$. For $x\in U_i-S$ with $S=\textnormal{Sing}\left(\mathcal{F}_0 \right)$, we choose a neighborhood $U_x\subset U_i-S$ of $x$ such that $dw_{0i}^\alpha- \sum_{\beta=1}^q a_{i_x}^{\alpha\beta} \wedge w_{0i}^\beta$ on $U_x$ for some $a_{i_x}^{\alpha\beta}\in \Gamma\left(U_i, \Omega_X^1 \right)$. Then 
\begin{align*}
dw_i^\alpha= \sum_{\beta=1}^q b_{i_x}^{\alpha\beta} \wedge w_i^\beta\,\,\,\,\,\,\,\textnormal{for some}\,\,\,b_{i_x}^{\alpha\beta}\in \Gamma\left( U_x, \Omega_X^1 \right)\otimes A
\end{align*}
such that $b_{i_x}^{\alpha\beta}$ induces $a_{i_x}^{\alpha\beta}$ by modulo $\mathfrak{m}_A$. Choose a lifting $\tilde{b}_{i_x}^{\alpha\beta}\in \Gamma\left( U_x, \Omega_X^1 \right)\otimes \tilde{A} $ of $b_{i_x}^{\alpha\beta}$. Then we have
\begin{align}\label{b5}
d\tilde{w}_i^\alpha-\sum_{\beta=1}^q \tilde{b}_{i_x}^{\alpha\beta}\wedge \tilde{w}_i^\beta = t C_{i_x}^\alpha \,\,\,\,\,\,\,\,\,\,\textnormal{for some}\,\,\, C_{i_x}^{\alpha}\in \Gamma\left(U_{i_x}, \bigwedge^2 \Omega_X^1\right)
\end{align}
Then we have
\begin{align}\label{si42}
tB_i^\alpha =\tilde{w}_i^1\wedge \cdots \wedge \tilde{w}_i^q\wedge \left(tC_{i_x}^\alpha+ \sum_{\beta=1}^q \tilde{b}_{i_x}^{\alpha\beta}\wedge \tilde{w}_i^\beta \right)= tw_{0i}^1\wedge \cdots \wedge w_{0i}^q\wedge C_{i_x}^\alpha
\end{align}
so that we have
\begin{align*}
B_i^\alpha=w_i^1\wedge \cdots \wedge w_i^q\wedge C_{i_x}^\alpha
\end{align*}
This implies that $B_i^\alpha\in \Gamma \left(U_i, \tilde{\mathcal{S}}^2\right)$. We define $B_i\in \Gamma\left( U_i, \mathscr{H}om_{\mathcal{O}_X}\left( \mathcal{N}_{\mathcal{F}_0}^*, \tilde{\mathcal{S}}^2 \right)\right)$ by
\begin{align}\label{nb7}
B_i:\Gamma\left(U_i, \mathcal{N}_{\mathcal{F}_0}^* \right) &\to  \Gamma\left( U_i, \tilde{\mathcal{S}}^2 \right)\\
 w_{0i}^\alpha &\mapsto B_i^\alpha \notag
\end{align}
Then we claim that $\left( \{B_i \}, \{ R_{ij}\} \right) \in C^0\left( \mathcal{U}, \mathscr{H}om_{\mathcal{O}_X}\left( \mathcal{N}_{\mathcal{F}_0}^*, \tilde{\mathcal{S}}^2 \right) \right) \bigoplus C^1\left( \mathscr{H}om_{\mathcal{O}_X}\left( \mathcal{N}_{\mathcal{F}_0}^* , \Omega_X^1 \right) \right)$ define a $1$-cocycle in the following \v Cech resolution of $\mathscr{H}om_{\mathcal{O}_X}\left( \mathcal{N}_{\mathcal{F}_0}, \Omega_X^1 \right)^\bullet$:
\begin{center}
$\begin{CD}
\cdots \\
@AE_3AA \\
C^0\left( \mathcal{U}, \mathscr{H}om_{\mathcal{O}_X}\left( \mathcal{N}_{\mathcal{F}_0}^*,   \tilde{\mathcal{S}}^3 \right) \right) @>>> \cdots\\
@AE_2AA @AE_2AA\\
C^0\left( \mathcal{U}, \mathscr{H}om_{\mathcal{O}_X}\left( \mathcal{N}_{\mathcal{F}_0}^*, \tilde{\mathcal{S}}^2 \right) \right) @>-\delta >> C^1\left( \mathcal{U}, \mathscr{H}om_{\mathcal{O}_X}\left( \mathcal{N}_{\mathcal{F}_0}^*, \tilde{\mathcal{S}}^2 \right) \right) @>\delta>> \cdots \\
@AE_1'AA @AE_1'AA @AE_1'AA \\
C^0\left( \mathcal{U}, \mathscr{H}om_{\mathcal{O}_X}\left( \mathcal{N}_{\mathcal{F}_0}^*,  \Omega_X^1 \right) \right) @>\delta>> C^1\left(\mathcal{U}, \mathscr{H}om_{\mathcal{O}_X}\left( \mathcal{N}_{\mathcal{F}_0}^*, \Omega_X^1 \right) \right) @>-\delta>> C^2\left( \mathcal{U}, \mathscr{H}om_{\mathcal{O}_X}\left( \mathcal{N}_{\mathcal{F}_0}^*, \Omega_X^1 \right) \right) @>>> \cdots
\end{CD}$
\end{center}

We note that by taking $d$ on $(\ref{b5})$, we have
\begin{align}\label{b6}
-\sum_{\beta=1}^q d\tilde{b}_{i_x}^{\alpha\beta}\wedge \tilde{w}_i^\beta+\sum_{\beta=1}^q \tilde{b}_{i_x}^{\alpha\beta} \wedge d \tilde{w}_i^\beta = t dC_{i_x}^\alpha
\end{align}
Then by multiplying $\tilde{w}_i^1 \wedge \cdots \wedge \tilde{w}_i^q$ on $(\ref{b6})$, we have
\begin{align*}
\sum_{\beta=1}^q \tilde{w}_i^1\wedge \cdots \wedge \tilde{w}_i^q\wedge \tilde{b}_{i_x}^{\alpha\beta} \wedge \left( t C_{i_x}^\beta + \sum_{\gamma=1}^q \tilde{b}_{i_x}^{\beta \gamma} \wedge \tilde{w}_i^\gamma \right) = t w_{0i}^1\wedge \cdots \wedge w_{0i}^q\wedge d C_{i_x}^\alpha
\end{align*}
so that we have
\begin{align*}
w_{0i}^1\wedge \cdots w_{0i}^q\wedge \left(dC_{i_x}^\alpha-\sum_{\beta=1}^q a_{i_x}^{\alpha\beta}\wedge C_{i_x}^\beta\right)=0
\end{align*}
This implies that $E_2\left( \left\{ B_i \right\} \right)=0$. We note that
{\small{\begin{align*}
tB_i^\alpha&=\tilde{w}_i^1\wedge \cdots \wedge \tilde{w}_i^q\wedge d\tilde{w}_i^\alpha=\left(tR_{ij}^1+ \sum_{\beta=1}^q h_{0ij}^{1\beta} \tilde{w}_j^\beta\right) \wedge \cdots \wedge \left( t R_{ij}^q+\sum_{\beta=1}^q h_{0ij}^{q\beta} \tilde{w}_j^\beta \right) \wedge d\left(tR_{ij}^\alpha+\sum_{\beta=1}^q h_{0ij}^{\alpha\beta} \tilde{w}_j^\beta \right)\\
                 &= t\sum_{\beta=1}^q w_{0i}^1\wedge \cdots \wedge R_{ij}^\beta \wedge \cdots \wedge w_{0i}^q \wedge dw_{0i}^\alpha+t w_{0i}^1\wedge \cdots \wedge w_{0i}^q\wedge dR_{ij}^\alpha +\det\left( h_{0ij}\right) \tilde{w}_j^1\wedge \cdots \wedge \tilde{w}_j^q \wedge \left(\sum_{\beta=1}^q dh_{0ij}^{\alpha\beta}\wedge \tilde{w}_j^\beta+\sum_{\beta=1}^q h_{0ij}^{\alpha\beta} d\tilde{w}_j^\beta\right)\\
                                  &= t\sum_{\beta=1}^q w_{0i}^1\wedge \cdots \wedge R_{ij}^\beta \wedge \cdots \wedge w_{0i}^q \wedge dw_{0i}^\alpha+t w_{0i}^1\wedge \cdots \wedge w_{0i}^q\wedge d R_{ij}^\alpha +\sum_{\beta=1}^q h_{0ij}^{\alpha\beta}\det\left( h_{0ij}\right) \tilde{w}_j^1\wedge \cdots \wedge \tilde{w}_j^q \wedge  d\tilde{w}_j^\beta
\end{align*}}}
so that we have
\begin{align*}
B_i^\alpha- \sum_{\beta=1}^q h_{0ij}^{\alpha\beta} \det \left(h_{0ij}\right) B_j^\beta =\sum_{\beta=1}^q w_{0i}^1\wedge \cdots \wedge R_{ij}^\beta\wedge \cdots \wedge w_{0i}^q\wedge dw_{0i}^\alpha+ w_{0i}^1\wedge \cdots \wedge w_{0i}^q\wedge dR_{ij}^\alpha
\end{align*}
This implies that $-\delta\left( \left\{ B_i \right\} \right) + E_1'\left( \left\{ -R_{ij} \right\} \right)=0$. We note that
\begin{align*}
tR_{ij}^\alpha-tR_{ik}^\alpha+t\sum_{\beta=1}^q h_{0ij}^{\alpha\beta} R_{jk}^\beta&=\tilde{w}_i^\alpha-\sum_{\beta=1}^q h_{0ij}^{\alpha\beta} \tilde{w}_j^\beta -\tilde{w}_i^\alpha+\sum_{\beta=1}^q h_{0ik}^{\alpha\beta} \tilde{w}_k^\beta+\sum_{\beta=1}^q h_{0ij}^{\alpha\beta} \tilde{w}_j^\beta-\sum_{\beta,\gamma=1}^q h_{0ij}^{\alpha\beta} h_{0jk}^{\beta\gamma} \tilde{w}_k^\gamma=0
\end{align*}
This implies that $-\delta\left( \left\{ - R_{ij} \right\} \right)=0$. This proves the claim. We show that for another lifting $\tilde{w}_i'^\alpha$ of $w_i^\alpha,\alpha=1,...,q$, the associated $1$-cocycle $\beta:=\left( \left\{ B_i'\right\},\left\{ - R_{ij}' \right\} \right)$ is cohomologous to $\alpha:=\left( \left\{ B_i\right\}, \left\{ - R_{ij} \right\} \right)$ as above. Then $\tilde{w}_i'^\alpha = \tilde{w}_i^\alpha + t V_i^\alpha$ for some $V_i^\alpha \in \Gamma\left( U_i, \Omega_X^1 \right)$. We define $V_i\in \Gamma\left( U_i, \mathscr{H}om_{\mathcal{O}_X} \left(  \mathcal{N}_{\mathcal{F}_0}^*, \Omega_X^1   \right) \right)$ by
\begin{align*}
V_i:\Gamma\left( U_i, \mathcal{N}_{\mathcal{F}_0}^* \right) &\to \Gamma\left( U_i, \Omega_X^1 \right) \\
 w_{0i}^\alpha &\mapsto V_i^\alpha
\end{align*}
Then we see that
\begin{align}\label{b10}
t\left(B_i'^\alpha - B_i^\alpha   \right)&= \left( \tilde{w}_i^1 + t V_i^\alpha \right) \wedge \cdots \wedge \left( \tilde{w}_i^q + t V_i^q \right) \wedge d \left( \tilde{w}_i^\alpha + t V_i^\alpha \right) - \tilde{w}_i^1 \wedge \cdots \wedge w_i^q \wedge d \tilde{w}_i^\alpha \\
 &= t\sum_{\beta=1}^q w_{0i}^1 \wedge \cdots \wedge V_i^\beta \wedge \cdots \wedge w_{0i}^q \wedge d w_{0i}^\alpha + t w_{0i}^1 \wedge \cdots \wedge w_{0i}^q \wedge d V_i^\alpha \notag
\end{align}
This implies that $E_1'\left( \left\{ V_i \right\} \right)=\left\{ B_i'- B_i\right\}$. On the other hand, we see that
\begin{align*}
t\left(- R_{ij}' + R_{ij} \right) = -\tilde{w}_i^\alpha - t V_i^\alpha + \sum_{\beta=1}^q h_{0ij}^{\alpha\beta} \left(\tilde{w}_j^\beta + t V_j^\beta \right) + \tilde{w}_i^\alpha - \sum_{\beta=1}^q h_{0ij}^{\alpha\beta} \tilde{w}_j^\beta =  t\left( \sum_{\beta=1}^q h_{0ij}^{\alpha\beta} V_j^\beta- V_i^\alpha \right)
\end{align*}
This implies that $\delta\left(\left\{ V_i\right\} \right)=\left\{ - R_{ij}'-\left( - R_{ij}\right) \right\}$. This proves that $\alpha$ is cohomologous to $\beta$. So given a small extension $e:0\to (t)\to \tilde{A} \to A \to 0$ and an infinitesimal foliated deformation $\eta$ of $\left( X, \mathcal{N}_{\mathcal{F}_0}^* \right)$ with $X$ and $\mathcal{N}_{\mathcal{F}_0}^*$ fixed, we can associate an element $o_\eta(e):=$ the cohomology class of $\alpha \in \mathbb{H}^1\left( X, \mathscr{H}om_{\mathcal{O}_X}\left( \mathcal{N}_{\mathcal{F}_0}^* , \Omega_X^1 \right)^\bullet \right)$. We note that $o_\eta(e)=0$ if and only if there exists $\left\{ N_i \right\} \in C^0\left( \mathcal{U}, \mathscr{H}om_{\mathcal{O}_X}\left(  \mathcal{N}_{\mathcal{F}_0}^*, \Omega_X^1 \right) \right)$ where $N_i:\Gamma\left( U_i, \mathcal{N}_{\mathcal{F}_0}^* \right) \to \Gamma\left( U_i, \Omega_X^1 \right), w_{0i}^\alpha \mapsto N_i^\alpha$ such that $N_j- N_i=- R_{ij}$ and $E_1'\left(N_i \right)= B_i$, so that
{\small{\begin{align}
&t\left(\sum_{\beta=1}^q h_{0ij}^{\alpha\beta} N_j^\beta - N_i^\alpha\right)= - \tilde{w}_i^\alpha + \sum_{\beta=1}^q h_{0ij}^{\alpha\beta} \tilde{w}_j^\beta \Longrightarrow \left(\tilde{w}_i^\alpha - tN_i^\alpha\right) = \sum_{\beta=1}^q h_{0ij}^{\alpha\beta}\left(\tilde{w}_j^\beta - t N_j^\beta \right)\notag\\
& t\left( \sum_{\beta=1}^q w_{0i}^1\wedge \cdots \wedge N_i^\beta \wedge \cdots \wedge w_{0i}^q \wedge dw_{0i}^\alpha + w_{0i}^1 \wedge \cdots \wedge w_{0i}^q \wedge dN_i^\alpha  \right)= \tilde{w}_i^1 \wedge \cdots \wedge \tilde{w}_i^q \wedge d \tilde{w}_i^\alpha \notag\\
&\Longrightarrow \left(\tilde{w}_i^1 - t N_i^1 \right) \wedge \cdots \left( \tilde{w}_i^q - t N_i^q \right) \wedge d \left( \tilde{w}_i^\alpha - t N_i^\alpha \right)= 0 \label{b12}
\end{align}}}
Hence a lifting $\tilde{\eta}$ of $\eta$ exists. This completes the proof of Proposition \ref{b2}.
\end{proof}

\section{Foliated deformations of $\left(X, \mathcal{N}_{\mathcal{F}_0}^*\right)$ with $\mathcal{N}_{\mathcal{F}_0}^*$ locally free under the trivial deformation of $X$}

Let $\left(X, \mathcal{N}_{\mathcal{F}_0}^* \right)$ be a foliated nonsingular variety over $k$. In the previous section we study flat foliated deformations of $\left(X, \mathcal{N}_{\mathcal{F}_0}^* \right)$ such that the underling flat deformations of $\left(X,\mathcal{N}_{\mathcal{F}_0}^* \right)$ is trivial. In this section, in addition, we vary $\mathcal{N}_{\mathcal{F}_0}^*$ with $X$ fixed.

\begin{definition}
Let $\left( X, \mathcal{N}_{\mathcal{F}_0}^* \right)$ be a foliated nonsingular variety over $k$ with $\mathcal{N}_{\mathcal{F}_0}^*$ locally free, so that we have an exact sequence $0\to \mathcal{N}_{\mathcal{F}_0}^* \to \Omega_X^1 \to \frac{\Omega_X^1}{\mathcal{N}_{\mathcal{F}_0}^*}\to 0$.  An infinitesimal foliated deformations of $\left(X, \mathcal{N}_{\mathcal{F}_0}^* \right)$ under
 the trivial deformation of $X$ over $A\in \bold{Art}$ is a  locally free subsheaf  $\mathcal{N}_{\mathcal{F}}^*$ of $\Omega_{\mathcal{X}/A}^1=\Omega_{X\times \textnormal{Spec}(A)/\textnormal{Spec}(A)}^1$ such that $ d_{\mathcal{X}/A}\left( \mathcal{N}_{\mathcal{F},x}^* \right) \subset \mathcal{N}_{\mathcal{F},x}^*\bigwedge \Omega_{\mathcal{X}/A,x}^1 $ for $x\in X- \textnormal{Sing}\left( \mathcal{F}_0 \right)$ and $\Omega_{\mathcal{X}/A}^1/ \mathcal{N}_{\mathcal{F}}^*$ is flat over $\textnormal{Spec}(A)$ and the exact sequence on $\mathcal{X}:=X\times_{\textnormal{Spec}(k)} \textnormal{Spec}(A)$ 
\begin{align}\label{}
0\to \mathcal{N}_{\mathcal{F}}^* \to \Omega_{\mathcal{X}/A}^1=\Omega_X^1\otimes A \to \frac{\Omega_X^1\otimes A}{\mathcal{N}_{\mathcal{F}}^*}\to 0
\end{align}
induces $0\to \mathcal{N}_{\mathcal{F}_0}^* \to \Omega_X^1 \to \frac{\Omega_X^1}{\mathcal{N}_{\mathcal{F}_0}^*}\to 0$ over $\textnormal{Spec}(k)$.
Then we have a functor of Artin rings
\begin{align*}
\textnormal{FDef}_{\mathcal{N}_{\mathcal{F}_0}}: \bold{Art} &\to \textnormal{(Sets)} \\
                             A &\to \left\{ \textnormal{Infinitesimal foliated deformation of $\left( X, \mathcal{N}_{\mathcal{F}_0}^*\right)$ with $X$ fixed}   \right\}
\end{align*}
\end{definition}

We will show that the following complex of sheaves truncated from the dual leaf complex (see Part I)
{\small{\begin{align*}
\mathscr{H}om_{\mathcal{O}_X}\left( \mathcal{N}_{\mathcal{F}_0}, \Omega_X^1 /\mathcal{N}_{\mathcal{F}_0}^* \right)^\bullet :\mathscr{H}om_{\mathcal{O}_X}\left(\mathcal{N}_{\mathcal{F}_0}^*, \Omega_X^1/ \mathcal{N}_{\mathcal{F}_0}^* \right) \xrightarrow{E_1} \mathscr{H}om_{\mathcal{O}_X}\left(\mathcal{N}_{\mathcal{F}_0}^* , \tilde{\mathcal{S}}^2 \right) \xrightarrow{E_2}  \mathscr{H}om_{\mathcal{O}_X}\left( \mathcal{N}_{\mathcal{F}_0}^* , \tilde{\mathcal{S}}^3 \right) \xrightarrow{E_3} \cdots
\end{align*}}}
controls foliated deformations of $\left(X, \mathcal{N}_{\mathcal{F}_0}^* \right)$ under the trivial deformation of $X$ in the language of functors of Artin rings. We will denote the $i$-th cohomology group of $\mathscr{H}om_{\mathcal{O}_X}\left(  \mathcal{N}_{\mathcal{F}_0}^* , \frac{\Omega_X^1}{\mathcal{N}_{\mathcal{F}_0}^*} \right)^\bullet$ by $\mathbb{H}^i\left( X, \mathscr{H}om_{\mathcal{O}_X}\left(  \mathcal{N}_{\mathcal{F}_0}^* , \frac{\Omega_X^1}{\mathcal{N}_{\mathcal{F}_0}^*} \right)^\bullet\right)$.

\begin{proposition}\label{b11}
Let $\left(X, \mathcal{N}_{\mathcal{F}_0}^* \right)$ be a foliated nonsingular variety with $\mathcal{N}_{\mathcal{F}_0}^*$ locally free.
\begin{enumerate}
\item there is a $1-1$ correspondence
\begin{align*}
\textnormal{FDef}_{\mathcal{N}_{\mathcal{F}_0}^*}\left(k[\epsilon] \right)\cong \mathbb{H}^0\left(X, \mathscr{H}om_{\mathcal{O}_X}\left(  \mathcal{N}_{\mathcal{F}_0}^*, \frac{\Omega_X^1}{ \mathcal{N}_{\mathcal{F}_0}^*} \right)^\bullet \right)
\end{align*}
\item Let $A\in \bold{Art}$ and $\eta=\left( \mathcal{X}, \mathcal{N}_\mathcal{F}^* \right)$ be an infinitesimal foliated deformations over $A$ of $\left( X, \mathcal{N}_{\mathcal{F}_0}^* \right)$ with $X$ fixed. To every small extension $0\to (t)\to \tilde{A} \to A \to 0$, we can associate an element $o_\eta(e)\in \mathbb{H}^1\left( X, \mathscr{H}om_{\mathcal{O}_X}\left( \mathcal{N}_{\mathcal{F}_0}^*, \frac{\Omega_X^1}{\mathcal{N}_{\mathcal{F}_0}^*} \right)^\bullet \right)$ called the obstruction lifting of $\eta$ to $\tilde{A}$, which is $0$ if and only if a lifting of $\eta$ to $\tilde{A}$ exists.
\end{enumerate}
\end{proposition}

\begin{proof}
Given a first-order foliated deformation $\left(X\times \textnormal{Spec}\left( k[\epsilon]\right), \mathcal{N}_\mathcal{F}^*\right)$ of $\left( X, \mathcal{N}_{\mathcal{F}_0}^* \right)$ under the trivial deformation of $X$, we choose an affine open covering $\mathcal{U}=\{U_i\}$ of $X$ such that $\Gamma\left( U_i, \mathcal{N}_{\mathcal{F}_0}^* \right)$ is generated by $w_{0i}^1,..., w_{0i}^q\in \Gamma\left( U_i, \Omega_X^1 \right)$ with $w_{0i}^\alpha= \sum_{\beta=1}^q h_{0ij}^{\alpha\beta} w_{0j}^\beta$ for some $h_{0ij}^{\alpha\beta}\in \Gamma\left( U_{ij}, \mathcal{O}_X \right)$ and $w_{0i}^1 \wedge \cdots \wedge w_{0i}^q \wedge d w_{0i}^\alpha =0, \alpha=1,...,q$, and $\Gamma\left( U_i, \mathcal{N}_\mathcal{F}^* \right)$ is generated by $w_{0i}^\alpha + \epsilon w_i'^\alpha$ for some $w_i'^\alpha \in \Gamma\left( U_i, \Omega_X^1 \right),\alpha=1,...,q$ with $ w_{0i}^\alpha + \epsilon w_i'^\alpha = \sum_{\beta=1}^q \left(h_{0ij}^{\alpha \beta} + \epsilon g_{ij}^{\alpha\beta} \right)\left( w_{0j}^\alpha + \epsilon w_j'^\beta\right)$ for some $g_{ij}^{\alpha\beta}\in \Gamma\left( U_{ij}, \mathcal{O}_X \right)$ and $\left(w_{0i}^1+ \epsilon w_i'^1 \right)\wedge \cdots \left( w_{0i}^q + \epsilon w_i'^q \right)\wedge d \left(w_{0i}^\alpha + \epsilon w_i'^\alpha \right)=0,\alpha=1,...q$. Then we have
\begin{align}
&w_i'^\alpha= \sum_{\beta=1}^q h_{0ij}^{\alpha\beta} w_j'^\beta + \sum_{\beta=1}^q g_{ij}^{\alpha\beta} w_{0j}^\beta \label{b7}\\
\sum_{\beta=1}^q w_{0i}^1 \wedge \cdots \wedge w_i'^\beta &\wedge \cdots \wedge w_{0i}^q \wedge dw_{0i}^\alpha+ w_{0i}^1 \wedge \cdots \wedge w_{0i}^q \wedge dw_i'^\alpha =0 \label{b8}
\end{align}
We define $\overline{w_i'} \in \Gamma\left(U_i, \mathscr{H}om_{\mathcal{O}_X}\left( \mathcal{N}_{\mathcal{F}_0}^*, \frac{\Omega_X^1}{\mathcal{N}_{\mathcal{F}_0}^*} \right) \right)$ in the following way:
\begin{align}\label{nb4}
\overline{w_i'}:\Gamma\left( U_i, \mathcal{N}_{\mathcal{F}_0}^* \right) &\to \Gamma\left(  U_i, \frac{\Omega_X^1}{\mathcal{N}_{\mathcal{F}_0}^*}\right)\\
 w_{0i}^\alpha &\mapsto \overline{w_i'^\alpha} \notag
\end{align}
and linearly extends to $\Gamma\left( U_i , \mathcal{N}_{\mathcal{F}_0}^* \right)$. Here $\overline{w_i'^\alpha}$ is the image of $w_i'^\alpha$ in $\Gamma\left( U_i, \frac{\Omega_X^1}{\mathcal{N}_{\mathcal{F}_0}^*} \right)$. Then $(\ref{b7})$ and $(\ref{b8})$ implies that $w' :=\left\{ \overline{w_i' }\right\} \in \mathbb{H}^0\left( X, \mathscr{H}om_{\mathcal{O}_X}\left( \mathcal{N}_{\mathcal{F}_0}^*, \frac{\Omega_X^1}{\mathcal{N}_{\mathcal{F}_0}^*} \right)^\bullet \right)$.

We define an inverse map. Given a $w'\in \mathbb{H}^0\left(X, \mathscr{H}om_{\mathcal{O}_X}\left(\mathcal{N}_{\mathcal{F}_0}^*, \frac{\Omega_X^1}{\mathcal{N}_{\mathcal{F}_0}^*} \right)^\bullet \right)$, we have $w_i':\Gamma\left(U_i, \mathcal{N}_{\mathcal{F}_0}^*\right)\to \Gamma\left( U_i, \frac{\Omega_X^1}{\mathcal{N}_{\mathcal{F}_0}^*}\right), w_{0i}^\alpha \mapsto \overline{w_i'^\alpha}$ for some $w_i'^\alpha\in \Gamma\left(U_i, \Omega_X^1 \right)$ and $\sum_{\beta=1}^q w_{0i}^1 \wedge \cdots \wedge w_i'^\beta \wedge \cdots \wedge w_{0i}^q \wedge dw_{0i}^\alpha + w_{0i}^1 \wedge \cdots \wedge w_{0i}^q \wedge d w_i'^\alpha=0 , \alpha=1,...,q$. Then by reversing the above process, we obtain a first-order foliated deformation $\left( \mathcal{X}, \mathcal{N}_{\mathcal{F}}^*\right)$ with $X$ fixed.

Next we identify obstructions. Given a small extension $0\to (t) \to \tilde{A} \to A \to 0$ in $\bold{Art}$ and a foliated deformation $\left(\mathcal{X}, \mathcal{N}_\mathcal{F}^* \right)$ over $A$ of $\left( X, \mathcal{N}_{\mathcal{F}_0}^* \right)$ with $X$ fixed, we choose an affine open covering $\mathcal{U}=\{U_i\}$ of $X$  such that $\Gamma\left( U_i, \mathcal{N}_{\mathcal{F}}^* \right)$ is generated by $w_i^1,..., w_i^q\in \Gamma\left(U_i, \Omega_X^1 \right)\otimes A$ which induces $w_{0i}^\alpha$ by modulo $\mathfrak{m}_A$ and $w_i^\alpha= \sum_{\beta=1}^q h_{ij}^{\alpha\beta} w_j^\beta$ for some $h_{ij}^{\alpha\beta}\in \Gamma\left( U_{ij}, \mathcal{O}_X \right) \otimes A$ which induces $h_{0ij}^{\alpha\beta}$ by modulo $\mathfrak{m}_A$ and $w_i^1\wedge \cdots \wedge w_i^q \wedge dw_i^\alpha=0,\alpha=1,...,q$. In order to see if there exists a lifting of $\left( \mathcal{X},  \mathcal{N}_\mathcal{F}^* \right)$ to $\tilde{A}$, we choose an arbitrary lifting $\tilde{w}_i^\alpha\in \Gamma\left( U_i, \Omega_X^1 \right)\otimes \tilde{A}$ of $w_i^\alpha$ and arbitrary lifting of $\tilde{h}_{ij}^{\alpha\beta}\in \Gamma\left( U_{ij}, \mathcal{O}_X \right)\otimes \tilde{A}$ of $h_{ij}^{\alpha\beta}$. Then we have
\begin{align*}
\tilde{w}_i^\alpha - \sum_{\beta=1}^q \tilde{h}_{ij} \tilde{w}_j^\beta= tR_{ij}^{\alpha}\,\,\,\,\,&\textnormal{for some}\,\,\,R_{ij}^\alpha \in \Gamma\left( U_{ij}, \Omega_X^1 \right)\\
\tilde{w}_i^1\wedge \cdots \wedge \tilde{w}_i^q \wedge d\tilde{w}_i^\alpha = t B_i^\alpha\,\,\,\,\,&\textnormal{for some}\,\,\, B_i^\alpha \in \Gamma\left( U_i, \tilde{\mathcal{S}}^2 \right)\\
\tilde{h}_{ik}^{\alpha \gamma}- \sum_{\beta=1}^q \tilde{h}_{ij}^{\alpha\beta} \tilde{h}_{jk}^{\beta \gamma}= t G_{ijk}^{\alpha \gamma}\,\,\,\,\,\,&\textnormal{for some}\,\,\, G_{ijk}^{\alpha \gamma} \in \Gamma\left( U_{ijk}, \mathcal{O}_X \right)
\end{align*}
Here we note that we can show $B_i^\alpha \in \Gamma\left(U_i, \tilde{\mathcal{S}}^2 \right)$ as in the proof of Proposition \ref{b2}. Then we define $\overline{R_{ij}}\in \Gamma\left( U_{ij}, \mathscr{H}om_{\mathcal{O}_X}\left( \mathcal{N}_{\mathcal{F}_0}^*, \frac{\Omega_X^1}{ \mathcal{N}_{\mathcal{F}_0}^*} \right) \right)$ in the following way:
\begin{align}\label{nb6}
\overline{R_{ij}}:\Gamma\left( U_{ij}, \mathcal{N}_{\mathcal{F}_0}^* \right) &\to \Gamma\left( U_{ij}, \frac{\Omega_X^1}{\mathcal{N}_{\mathcal{F}_0}^*} \right)\\
w_{0i}^\alpha &\mapsto \overline{ R_{ij}^\alpha } \notag
\end{align}
and linearly extends to $\Gamma\left(U_{ij}, \mathcal{N}_{\mathcal{F}_0}^* \right)$. Here $\overline{R_{ij}^\alpha }$ is the image of $R_{ij}^\alpha$ in $\Gamma\left( U_{ij}, \frac{\Omega_X^1}{\mathcal{N}_{\mathcal{F}_0}^*} \right)$, and we define $B_i\in \Gamma\left( U_i,  \mathscr{H}om_{\mathcal{O}_X}\left( \mathcal{N}_{\mathcal{F}_0}^*,  \tilde{\mathcal{S}}^2 \right) \right)$ in the following way:
\begin{align*}
B_i:\Gamma\left(  U_i, \mathcal{N}_{\mathcal{F}_0}^* \right) &\to  \Gamma\left( U_i, \tilde{\mathcal{S}}^2 \right) \\
 w_{0i}^\alpha &\mapsto B_i^\alpha
\end{align*}
and linearly extends to $\Gamma\left( U_i, \mathcal{N}_{\mathcal{F}_0}^* \right)$. Then we claim that $\left( \left\{ B_i \right\}, \left\{ - \overline{R_{ij}} \right\} \right)\in C^0 \left( \mathcal{U}, \mathscr{H}om_{\mathcal{O}_X}\left( \mathcal{N}_{\mathcal{F}_0}^*, \tilde{\mathcal{S}}^2 \right)   \right)$ $  \bigoplus C^1\left( \mathcal{U}, \mathscr{H}om_{\mathcal{O}_X}\left( \mathcal{N}_{\mathcal{F}_0}^*, \frac{\Omega_X^1}{\mathcal{N}_{\mathcal{F}_0}^*} \right) \right)$ defines a $1$-cocycle in the following \v Cech resolution of $\mathscr{H}om_{\mathcal{O}_X}\left( \mathcal{N}_{\mathcal{F}_0}^*, \frac{\Omega_X^1}{\mathcal{N}_{\mathcal{F}_0}^*} \right)^\bullet$:
\begin{center}
$\begin{CD}
\cdots \\
@AD_2AA\\
C^0\left( \mathcal{U}, \mathscr{H}om_{\mathcal{O}_X}\left( \mathcal{N}_{\mathcal{F}_0}^* , \tilde{\mathcal{S}}^2 \right) \right) @>-\delta>> \cdots \\
@AD_1'AA @AD_1' AA \\
C^0\left(\mathcal{U}, \mathscr{H}om_{\mathcal{O}_X}\left( \mathcal{N}_{\mathcal{F}_0}^* , \frac{\Omega_X^1}{ \mathcal{N}_{\mathcal{F}_0}^*} \right) \right) @>\delta>> C^1\left( \mathcal{U}, \mathscr{H}om_{\mathcal{O}_X}\left(  \mathcal{N}_{\mathcal{F}_0}^* , \frac{\Omega_X^1}{\mathcal{N}_{\mathcal{F}_0}^*} \right) \right) @>-\delta>> 
\end{CD}$
\end{center}
We note that
\begin{align}\label{b13}
tR_{ij}^\alpha-tR_{ik}^\alpha+t\sum_{\beta=1}^q h_{0ij}^{\alpha\beta} R_{jk}^\beta&=\tilde{w}_i^\alpha-\sum_{\beta=1}^q \tilde{h}_{ij}^{\alpha\beta} \tilde{w}_j^\beta -\tilde{w}_i^\alpha+\sum_{\beta=1}^q \tilde{h}_{ik}^{\alpha\beta} \tilde{w}_k^\beta+\sum_{\beta=1}^q \tilde{h}_{ij}^{\alpha\beta} \tilde{w}_j^\beta-\sum_{\beta,\gamma=1}^q \tilde{h}_{ij}^{\alpha\beta} \tilde{h}_{jk}^{\beta\gamma} \tilde{w}_k^\gamma\\
&=\sum_{\gamma=1}^q \left(\tilde{h}_{ik}^{\alpha\gamma}-\sum_{\beta=1}^q \tilde{h}_{ij}^{\alpha\beta}\tilde{h}_{jk}^{\beta\gamma} \right) \tilde{w}_k^\gamma=t\sum_{\gamma=1}^q t G_{ijk}^{\alpha\gamma} w_{0k}^\gamma \notag
\end{align}
This implies that $-\delta\left( \left\{ - R_{ij} \right\} \right)=0$. As in the proof of Proposition \ref{b2}, $-\delta\left( \left\{ B_i \right\} \right) + E_1'\left( \left\{ - R_{ij} \right\} \right)=0$ and $E_2\left( \left\{ B_i \right\} \right)=0$. We show that for another lifting $\tilde{w}_i'^\alpha$ of $w_i^\alpha$, and another lifting $\tilde{h}_{ij}'^{\alpha\beta}$ of $h_{ij}^{\alpha\beta},\alpha,\beta=1,...,q$, the associated $1$-cocycle $\beta:=\left( \left\{ B_i'\right\} , \left\{ - R_{ij}'\right\} \right)$ is cohomologous to $\alpha:=\left(  \left\{ B_i \right\}, \left\{ - R_{ij} \right\} \right)$. Then $\tilde{w}_i'^\alpha= \tilde{w}_i^\alpha + t V_i^\alpha$ for some $V_i^\alpha\in \Gamma\left(U_i, \Omega_X^1 \right)$, and $\tilde{h}_{ij}'^{\alpha\beta}= \tilde{h}_{ij}^{\alpha\beta} + t y_{ij}^{\alpha\beta}$ for some $y_{ij}^{\alpha\beta}\in \Gamma\left( U_{ij}, \mathcal{O}_X \right)$. We define $V_i\in \Gamma\left( U_i, \mathscr{H}om_{\mathcal{O}_X}\left(\mathcal{N}_{\mathcal{F}_0}^*, \frac{\Omega_X^1}{\mathcal{N}_{\mathcal{F}_0}^*} \right) \right)$ by $V_i:\Gamma\left( U_i, \mathcal{N}_{\mathcal{F}_0}^* \right)\to \Gamma\left( U_i, \frac{\Omega_X^1}{\mathcal{N}_{\mathcal{F}_0}^*} \right), w_{0i}^\alpha \mapsto \overline{V_i^\alpha}$. Then as in $(\ref{b10})$, we have $E_1\left(\left\{V_i \right\} \right)= \left\{ B_i'- B_i\right\}$. On the other hand, we have
\begin{align}\label{b14}
t\left(- R_{ij}'^\alpha + R_{ij}^\alpha \right) &= -\tilde{w}_i^\alpha - t V_i^\alpha + \sum_{\beta=1}^q \left(\tilde{h}_{ij}^{\alpha\beta}+ t y_{ij}^{\alpha\beta} \right) \left(\tilde{w}_j^\beta + t V_j^\beta \right) + \tilde{w}_i^\alpha - \sum_{\beta=1}^q \tilde{h}_{ij}^{\alpha\beta} \tilde{w}_j^\beta  \\
&=  t\left( \sum_{\beta=1}^q h_{0ij}^{\alpha\beta} V_j^\beta- V_i^\alpha \right) + t\sum_{\beta=1}^q  y_{ij}^{\alpha\beta} w_{0j}^\beta \notag
\end{align}
This implies that $\delta\left(\left\{ V_i\right\} \right)=\left\{ -R_{ij}'-\left(- R_{ij}\right) \right\}$. This proves that $\alpha$ is cohomologous to $\beta$. So given a small extension $e:0\to (t)\to \tilde{A} \to A \to 0$ and an infinitesimal foliated deformation $\eta$ of $\left( X, \mathcal{N}_{\mathcal{F}_0}^*\right)$ with $X$ fixed, we can associate an element $o_\eta(e):=$ the cohomology class of $\alpha \in \mathbb{H}^1\left( X, \mathscr{H}om_{\mathcal{O}_X}\left( \mathcal{N}_{\mathcal{F}_0}^* , \frac{\Omega_X^1}{\mathcal{N}_{\mathcal{F}_0}^*} \right)^\bullet \right)$. We note that $o_\eta(e)=0$ if and only if there exists $\left\{ N_i \right\} \in C^0\left( \mathcal{U}, \mathscr{H}om_{\mathcal{O}_X}\left( \mathcal{N}_{\mathcal{F}_0}^* , \frac{\Omega_X}{\mathcal{N}_{\mathcal{F}_0}^*} \right)\right)$ where $N_i:\Gamma\left(U_i, \mathcal{N}_{\mathcal{F}_0}^* \right)\to \Gamma\left( U_i, \Omega_X^1 \right), w_{0i}^\alpha \mapsto \overline{N_i^\alpha}$ with $N_i^\alpha \in \Gamma\left( U_i, \Omega_X^1 \right)$ such that $N_j- N_i= -R_{ij}$ and $E_1\left(N_i \right)=B_i$, so that
{\small{\begin{align}
&t\left(\sum_{\beta=1}^q h_{0ij}^{\alpha\beta} N_j^\beta - N_i^\alpha + t \sum_{\beta=1}^q u_{ij}^{\alpha\beta} w_{0j}^\beta\right)= - \tilde{w}_i^\alpha + \sum_{\beta=1}^q \tilde{h}_{ij}^{\alpha\beta} \tilde{w}_j^\beta \notag \\
& \Longrightarrow \tilde{w}_i^\alpha - tN_i^\alpha  = \sum_{\beta=1}^q \left(\tilde{h}_{0ij}^{\alpha\beta}- tu_{ij}^{\alpha\beta} \right)\left(\tilde{w}_j^\beta - t N_j^\beta \right)\,\,\,\textnormal{for some} \,\,\, u_{ij}^{\alpha\beta}\in \Gamma\left( U_{ij}, \mathcal{O}_X \right) \label{b16}
\end{align}}}
and as in $(\ref{b12})$, we have  $\left(\tilde{w}_i^1 - t N_i^1 \right) \wedge \cdots \left( \tilde{w}_i^q - t N_i^q \right) \wedge d \left( \tilde{w}_i^\alpha - t N_i^\alpha \right)= 0 , \alpha=1,...,q $. Hence a lifting $\tilde{\eta}$ of $\eta$ exits. This completes the proof of Proposition \ref{b11}.

\end{proof}

\subsection{Foliated deformations of $\mathcal{N}_{\mathcal{F}_0}^*$ with $X$ fixed in terms of vector bundles} \

Let $\left( X, \mathcal{N}_{\mathcal{F}_0}^* \right)$ be a foliated nonsingular variety with $\mathcal{N}_{\mathcal{F}_0}^*$ locally free. The underlying deformations of locally free sheaf $\mathcal{N}_{\mathcal{F}_0}^*$ is implicit in the complex $\mathscr{H}om_{\mathcal{O}_X}\left(\mathcal{N}_{\mathcal{F}_0}^*, \frac{\Omega_X^1}{\mathcal{N}_{\mathcal{F}_0}^* } \right)^\bullet$. In order to make the underlying deformations of $\mathcal{N}_{\mathcal{F}_0}^*$ explicit, we define another complex of sheaves also controlling flat foliated deformations of $\left(X, \mathcal{N}_{\mathcal{F}_0}^* \right)$ under the trivial deformation of $X$ and present foliated deformations of $\left(X, \mathcal{N}_{\mathcal{F}_0}^* \right)$ with $X$ fixed by using the complex. As is known, deformations of a locally free sheaf $\mathcal{N}_{\mathcal{F}_0}^*$ on $X$ is controlled by $\mathscr{H}om_{\mathcal{O}_X}\left( \mathcal{N}_{\mathcal{F}_0}^* , \mathcal{N}_{\mathcal{F}_0}^* \right)$ (see \cite{Har10}). We note that we have an injection $j:\mathscr{H}om_{\mathcal{O}_X}\left( \mathcal{N}_{\mathcal{F}_0}^* , \mathcal{N}_{\mathcal{F}_0}^* \right) \to \mathscr{H}om_{\mathcal{O}_X}\left( \mathcal{N}_{\mathcal{F}_0}^* , \Omega_X^1 \right)$ induced by $\mathcal{N}_{\mathcal{F}_0}^* \hookrightarrow \Omega_X^1$. Then let us consider an exact sequence of complexes of sheaves
{\small{\begin{center}
$\begin{CD}
@. \cdots @. \cdots @. \cdots \\
@. @AAA @AE_3AA @AE_3AA \\
0@>>>0 @>>> \mathscr{H}om_{\mathcal{O}_X}\left( \mathcal{N}_{\mathcal{F}_0}^* , \tilde{\mathcal{S}}^3 \right) @>>> \mathscr{H}om_{\mathcal{O}_X}\left( \mathcal{N}_{\mathcal{F}_0}^* ,  \tilde{\mathcal{S}}^2 \right)@>>> 0 \\
@.@AAA @AE_2AA @AE_2AA \\
0@>>>0 @>>> \mathscr{H}om_{\mathcal{O}_X}\left( \mathcal{N}_{\mathcal{F}_0}^* , \tilde{\mathcal{S}}^2 \right) @>>> \mathscr{H}om_{\mathcal{O}_X}\left(  \mathcal{N}_{\mathcal{F}_0}^* ,  \tilde{\mathcal{S}}^2 \right) @>>> 0\\
@. @AAA @AE_1'AA @AE_1AA \\
0@>>> \mathscr{H}om_{\mathcal{O}_X}\left( \mathcal{N}_{\mathcal{F}_0}^*, \mathcal{N}_{\mathcal{F}_0}^* \right) @>>> \mathscr{H}om_{\mathcal{O}_X}\left( \mathcal{N}_{\mathcal{F}_0}^* , \Omega_X^1 \right) @>>> \mathscr{H}om_{\mathcal{O}_X}\left( \mathcal{N}_{\mathcal{F}_0}^* ,\frac{\Omega_X^1}{\mathcal{N}_{\mathcal{F}_0}^*} \right) @>>> 0 \\
@. @A\textnormal{Id}AA @AjAA @AAA \\
0@>>> \mathscr{H}om_{\mathcal{O}_X}\left(\mathcal{N}_{\mathcal{F}_0}^*, \mathcal{N}_{\mathcal{F}_0}^* \right) @>>> \mathscr{H}om_{\mathcal{O}_X}\left( \mathcal{N}_{\mathcal{F}_0}^*, \mathcal{N}_{\mathcal{F}_0}^* \right) @>>> 0@>>> 0
\end{CD}$
\end{center}}}
Let us denote the middle complex of sheaves by $\mathscr{H}om_{\mathcal{O}_X}\left(\mathcal{N}_{\mathcal{F}_0}^*, \mathcal{N}_{\mathcal{F}_0}^* \right)^\bullet$:
{\Small{\begin{align}\label{b9}
\mathscr{H}om_{\mathcal{O}_X}\left( \mathcal{N}_{\mathcal{F}_0}^* , \mathcal{N}_{\mathcal{F}_0}^* \right)^\bullet: \mathscr{H}om_{\mathcal{O}_X}\left(\mathcal{N}_{\mathcal{F}_0}^*, \mathcal{N}_{\mathcal{F}_0}^* \right) \to  \mathscr{H}om_{\mathcal{O}_X}\left(\mathcal{N}_{\mathcal{F}_0}^* , \Omega_X^1 \right) \to \mathscr{H}om_{\mathcal{O}_X}\left( \mathcal{N}_{\mathcal{F}_0}^* , \tilde{\mathcal{S}}^2 \right)\to \mathscr{H}om_{\mathcal{O}_X}\left( \mathcal{N}_{\mathcal{F}_0}^* , \tilde{\mathcal{S}}^3 \right)\to \cdots
\end{align}}}
We will denote the $i$-th cohomology group by $\mathbb{H}^i\left( X, \mathscr{H}om_{\mathcal{O}_X}\left( \mathcal{N}_{\mathcal{F}_0}^* , \mathcal{N}_{\mathcal{F}_0}^* \right)^\bullet \right)$. From the exact sequence we see that
\begin{align*}
\mathbb{H}^i\left(X, \mathscr{H}om_{\mathcal{O}_X}\left( \mathcal{N}_{\mathcal{F}_0}^* , \mathcal{N}_{\mathcal{F}_0}^* \right)^\bullet \right)\cong \mathbb{H}^{i-1}\left( X, \mathscr{H}om_{\mathcal{O}_X}\left( \mathcal{N}_{\mathcal{F}_0}^*, \Omega_X^1/ \mathcal{N}_{\mathcal{F}_0}^* \right)^\bullet\right),\,\,\,\,\,i\geq 0 . 
\end{align*}
This shows that $\mathscr{H}om_{\mathcal{O}_X}\left( \mathcal{N}_{\mathcal{F}_0}^* , \mathcal{N}_{\mathcal{F}_0}^* \right)^\bullet $ also controls flat foliated deformations of $\left( X, \mathcal{N}_{\mathcal{F}_0}^* \right)$ with $X$ fixed. Moreover the complex $(\ref{tf1})$ explicitly contains $\mathscr{H}om_{\mathcal{O}_X}\left( \mathcal{N}_{\mathcal{F}_0}^* ,  \mathcal{N}_{\mathcal{F}_0}^* \right)$ as its component. We will reinterpret flat foliated deformations of $\left(X, \mathcal{N}_{\mathcal{F}_0}^* \right)$ with $X$ fixed by using the complex $(\ref{b9})$.

\begin{remark}
To emphasize deformations of the underlying locally free sheaf $\mathcal{N}_{\mathcal{F}_0}^*$, we can define infinitesimal deformations of $\left(X, \mathcal{N}_{\mathcal{F}_0}^* \right)$ with $X$ fixed equivalently in the following way: an infinitesimal foliated deformation $\left(\mathcal{X}=X\times \textnormal{Spec}(A), \mathcal{N}_{\mathcal{F}}^* \right)$ of $\left(X, \mathcal{N}_{\mathcal{F}_0}^* \right)$ under the trivial deformation of $X$ over $A\in \bold{Art}$ is a locally free sheaf $\mathcal{N}_{\mathcal{F}}^*$ and an exact sequence on $\mathcal{X}=X\times \textnormal{Spec}(A)$
\begin{align*}
0 \to \mathcal{N}_{\mathcal{F}}^* \xrightarrow{l} \Omega_{\mathcal{X}/A}=\Omega_X^1 \otimes A \to \frac{\Omega_X^1\otimes A}{l\left(\mathcal{N}_{\mathcal{F}}^* \right)}\to 0
\end{align*}
such that $d_{\mathcal{X}/A}\left( l\left(\mathcal{N}_{\mathcal{F},x}^* \right) \right) \subset l\left( \mathcal{N}_{\mathcal{F}}^* \right)_x \wedge \Omega_{\mathcal{X}/A,x}^1$ for $x\in X-\textnormal{Sing}\left( \mathcal{F}_0 \right)$ and $\frac{\Omega_{\mathcal{X}/A}^1}{l\left(\mathcal{N}_{\mathcal{F}}^* \right)}$ is flat over $A$, and induces $0\to \mathcal{N}_{\mathcal{F}_0}^* \to \Omega_X^1 \to \frac{\Omega_X^1}{\mathcal{N}_{\mathcal{F}_0}^*}\to 0$ over $\textnormal{Spec}(k)$. Two infinitesimal deformations $\left(\mathcal{X}, \mathcal{N}_{\mathcal{F}}^*\right)$ and $\left(\mathcal{X}, \mathcal{N}_{\mathcal{F}'}^*\right)$ of $\left(X, \mathcal{N}_{\mathcal{F}_0}^* \right)$ with $X$ fixed are equivalent  if there is an isomorphism $\alpha: \mathcal{N}_{\mathcal{F}}^* \to \mathcal{N}_{\mathcal{F}'}^* $ over $\mathcal{X}$ such that $l= l'\circ \alpha$ where $0 \to \mathcal{N}_{\mathcal{F}'}^* \xrightarrow{l'} \Omega_{\mathcal{X}/A}^1 \to \frac{\Omega_{\mathcal{X}/A}^1}{\mathcal{N}_{\mathcal{F}'}^*}\to 0$. Then we have a functor of Artin rings
\begin{align*}
\textnormal{FDef}_{\mathcal{N}_{\mathcal{F}_0}}: \bold{Art} &\to \textnormal{(Sets)} \\
                             A &\to \left\{ \textnormal{Infinitesimal foliated deformation of $\left( X, \mathcal{N}_{\mathcal{F}_0}^*\right)$ with $X$ fixed}   \right\}/\textnormal{isomorphisms}
\end{align*}
\end{remark}

\begin{proposition}\label{b15}
Let $\left(X, \mathcal{N}_{\mathcal{F}_0}^* \right)$ be a foliated nonsingular variety with $\mathcal{N}_{\mathcal{F}_0}^*$ locally free.
\begin{enumerate}
\item there is a $1-1$ correspondence
\begin{align*}
\textnormal{FDef}_{\mathcal{N}_{\mathcal{F}_0}^*}\left(k[\epsilon] \right)\cong \mathbb{H}^1\left(X, \mathscr{H}om_{\mathcal{O}_X}\left(  \mathcal{N}_{\mathcal{F}_0}^* , \mathcal{N}_{\mathcal{F}_0}^* \right)^\bullet \right)
\end{align*}
\item Let $A\in \bold{Art}$ and $\eta=\left( \mathcal{X}, \mathcal{N}_\mathcal{F}^* \right)$ be an infinitesimal foliated deformations over $A$ of $\left( X, \mathcal{N}_{\mathcal{F}_0}^* \right)$ with $X$ fixed. To every small extension $0\to (t)\to \tilde{A} \to A \to 0$, we can associate an element $o_\eta(e)\in \mathbb{H}^2\left( X, \mathscr{H}om_{\mathcal{O}_X}\left( \mathcal{N}_{\mathcal{F}_0}^* , \mathcal{N}_{\mathcal{F}_0}^* \right)^\bullet \right)$ called the obstruction lifting of $\eta$ to $\tilde{A}$, which is $0$ if and only if a lifting of $\eta$ to $\tilde{A}$ exists.
\end{enumerate}
\end{proposition}

\begin{proof}
We keep the notations in the proof of Proposition \ref{b11}. Given a first-order foliated deformation $\left(\mathcal{X}, \mathcal{N}_{\mathcal{F}}^* \right)$ of $\left( X, \mathcal{N}_{\mathcal{F}_0}^* \right)$ with $X$ fixed, let us define $w_i' \in \Gamma\left( U_i, \mathscr{H}om_{\mathcal{O}_X}\left( \mathcal{N}_{\mathcal{F}_0}^* , \Omega_X^1 \right)\right)$ by $w_i':\Gamma\left(U_i, \mathcal{N}_{\mathcal{F}_0}^* \right)\to \Gamma\left(U_i, \Omega_X^1 \right), w_{0i}^\alpha \mapsto w_i'^\alpha$ and define $G_{ij}\in \Gamma\left(U_{ij}, \mathscr{H}om_{\mathcal{O}_X}\left( \mathcal{N}_{\mathcal{F}_0}^* , \mathcal{N}_{\mathcal{F}_0}^* \right) \right)$ by
\begin{align*}
G_{ij}: \Gamma\left( U_{ij}, \mathcal{N}_{\mathcal{F}_0}^*\right) &\to \Gamma\left( U_{ij},  \mathcal{N}_{\mathcal{F}_0}^* \right) \\
 w_{0i}^\alpha &\mapsto \sum_{\beta=1}^q g_{ij}^{\alpha\beta} w_{0j}^\beta
\end{align*}
Then we claim that $\left( \left\{ w_i' \right\} , \left\{ G_{ij}  \right\} \right) \in C^0\left( \mathcal{U}, \mathscr{H}om_{\mathcal{O}_X}\left( \mathcal{N}_{\mathcal{F}_0}^* ,  \Omega_X^1  \right) \right) \bigoplus C^1\left(\mathcal{U}, \mathscr{H}om_{\mathcal{O}_X}\left(  \mathcal{N}_{\mathcal{F}_0}^* , \mathcal{N}_{\mathcal{F}_0}^* \right) \right)$ is a $1$-cocycle in the following \v Cech resolution of $\mathscr{H}om_{\mathcal{O}_X}\left( \mathcal{N}_{\mathcal{F}_0}^* , \mathcal{N}_{\mathcal{F}_0}^* \right)^\bullet$:{\small{\begin{center}
$\begin{CD}
\cdots \\
@AE_2AA \\
C^0\left(\mathcal{U}, \mathscr{H}om_{\mathcal{O}_X}\left(  \mathcal{N}_{\mathcal{F}_0}^* ,  \tilde{\mathcal{S}}^2  \right) \right) @>-\delta>> \cdots \\
@AE_1'AA @AE_1'AA \\
C^0\left( \mathcal{U}, \mathscr{H}om_{\mathcal{O}_X}\left( \mathcal{N}_{\mathcal{F}_0}^* , \Omega_X^1 \right) \right) @>\delta>> C^1\left( \mathcal{U}, \mathscr{H}om_{\mathcal{O}_X}\left( \mathcal{N}_{\mathcal{F}_0}^* , \Omega_X^1 \right) \right) @>-\delta>> \cdots \\
@AjAA @AjAA @AjAA \\
C^0\left( \mathcal{U}, \mathscr{H}om_{\mathcal{O}_X}\left( \mathcal{N}_{\mathcal{F}_0}^* , \mathcal{N}_{\mathcal{F}_0}^* \right) \right) @>-\delta>> C^1\left( \mathcal{U}, \mathscr{H}om_{\mathcal{O}_X}\left( \mathcal{N}_{\mathcal{F}_0}^* , \mathcal{N}_{\mathcal{F}_0}^* \right) \right) @>\delta>> C^2\left( \mathcal{U}, \mathscr{H}om_{\mathcal{O}_X}\left( \mathcal{N}_{\mathcal{F}_0}^* , \mathcal{N}_{\mathcal{F}_0}^* \right) \right) @>-\delta>> \cdots
\end{CD}$
\end{center}}}
As in the proof of Proposition  \ref{b11}, we have $E_1'\left(\left\{  w_i' \right\} \right)=0$ and $(\ref{b7})$ implies that $\delta\left( \left\{ w_i' \right\} \right)+ j\left(\left\{ G_{ij} \right\} \right)=0$ and we note that
\begin{align*}
&h_{0ik}^{\alpha \gamma} + \epsilon g_{ik}^{\alpha \gamma} = \sum_{\beta=1}^q \left( h_{0ij}^{\alpha \beta} + \epsilon g_{ij}^{\alpha\beta} \right) \left( h_{0jk}^{\beta\gamma} + \epsilon g_{jk}^{\beta \gamma} \right) \Longrightarrow g_{ik}^{\alpha \gamma}- \sum_{\beta=1}^q h_{0ij}^{\alpha \beta} g_{jk}^{\beta \gamma} - \sum_{\beta=1}^q g_{ij}^{\alpha\beta} h_{0jk}^{\beta \gamma} =0 \\
&\Longrightarrow \left(G_{ij} - G_{ik} + G_{jk} \right)\left( w_{0i}^\alpha \right) = \sum_{\beta=1}^q g_{ij}^{\alpha\beta} w_{0j}^\beta - \sum_{\beta=1}^q g_{ik}^{\alpha\beta} w_{0k}^\beta + \sum_{\beta,\gamma=1}^q h_{0ij}^{\alpha\beta} g_{jk}^{\beta \gamma} w_{0k}^\gamma = 0
\end{align*}
This implies that $\delta\left(\left\{ G_{ij} \right\} \right)=0$. This proves the claim. For another equivalent foliated deformations $\beta:=\left( \left\{ w_i''\right\}, \left\{ G_{ij}'\right\} \right)$ of $\left(X, \mathcal{N}_{\mathcal{F}_0}^* \right)$ with $X$ fixed, we show that $\beta$ is cohomologous to $\alpha:=\left( \left\{ w_i'\right\}, \left\{ G_{ij}\right\} \right)$. We note that we have $w_{0i}^\alpha+ \epsilon w_i''^\alpha= \sum_{\beta=1}^q \left( \delta_\alpha^\beta + \epsilon b_{i}^{\alpha\beta} \right)\left(w_{0i}^\beta + \epsilon w_i'^\beta \right)$ for some $b_i^{\alpha\beta} \in \Gamma\left( U_i, \mathcal{O}_X \right)$, so that we have $w_i''^\alpha= \sum_{\beta=1}^q b_i^{\alpha\beta} w_{0i}^\beta + w_i'^\alpha$.  Then we define $B_i\in \Gamma\left( U_i, \mathscr{H}om_{\mathcal{O}_X}\left( \mathcal{N}_{\mathcal{F}_0}^* , \mathcal{N}_{\mathcal{F}_0}^* \right) \right)$ by 
\begin{align*}
B_i:\Gamma\left( U_i, T_{\mathcal{F}_0}\right) \to \Gamma\left( U_i, T_{\mathcal{F}_0}\right)\\
 w_{0i}^\alpha \mapsto \sum_{\beta=1}^q b_i^{\alpha\beta} w_{0i}^\beta
\end{align*}
Then $j(B_i)(w_{0i}^\alpha)= w_i''^\alpha- w_i'^\alpha$ and since we have from $(\ref{b7})$
{\small{\begin{align*}
\left(G_{ij}'- G_{ij}\right)\left( w_{0i}^\alpha \right)= \sum_{\beta=1}^q g_{ij}'^{\alpha\beta} w_{0j}^\beta- \sum_{\beta=1}^q g_{ij}^{\alpha\beta}w_{0j}^\beta = w_i''^\alpha - w_i'^\alpha - \sum_{\beta=1}^q h_{0ij}^{\alpha\beta}\left( w_j''^\beta - w_j'^\beta\right)= \sum_{\beta=1}^q b_i^{\alpha\beta} w_{0i}^\beta - \sum_{\beta,\gamma=1}^q h_{0ij}^{\alpha\beta} b_j^{\beta \gamma} w_{0j}^\gamma
\end{align*}}}
we see that $-\delta\left(\left\{ B_i \right\} \right)= G_{ij}'-G_{ij}$. Hence $\alpha$ is cohomologus to $\beta$. As in Proposition $\ref{b11}$ we define an inverse map.

Next we identify obstructions. We define $R_{ij}\in \Gamma\left( U_{ij}, \mathscr{H}om_{\mathcal{O}_X}\left( \mathcal{N}_{\mathcal{F}_0}^* , \Omega_X^1 \right) \right)$ by $R_{ij}:\Gamma\left(U_{ij}, \mathcal{N}_{\mathcal{F}_0}^* \right) \to \Gamma\left( U_{ij}, \Omega_X^1 \right), w_{0i}^\alpha \mapsto R_{ij}^\alpha$ and $B_i\in \Gamma\left( U_i, \mathscr{H}om_{\mathcal{O}_X}\left(  \mathcal{N}_{\mathcal{F}_0}^* , \tilde{\mathcal{S}}^2 \right) \right)$ as in the proof of Proposition \ref{b11}. We define $G_{ijk}\in \Gamma\left( U_{ijk}, \mathscr{H}om_{\mathcal{O}_X}\left( \mathcal{N}_{\mathcal{F}_0}^* ,  \mathcal{N}_{\mathcal{F}_0}^* \right) \right)$ by
\begin{align*}
G_{ijk}:\Gamma\left(U_{ijk}, \mathcal{N}_{\mathcal{F}_0}^* \right)&\to \Gamma\left( U_{ijk}, \mathcal{N}_{\mathcal{F}_0}^* \right) \\
 w_{0i}^\alpha &\mapsto \sum_{\gamma=1}^q G_{ijk}^{\alpha \gamma} w_{0k}^\gamma
\end{align*}
Then $\left\{ G_{ijk}   \right\}\in C^2\left(\mathcal{U}, \mathscr{H}om_{\mathcal{O}_X}\left( \mathcal{N}_{\mathcal{F}_0}^*, \mathcal{N}_{\mathcal{F}_0}^* \right) \right)$. Then in a similar way to $(\ref{ot46})$, we have $-\delta\left(\left\{ -G_{ijk} \right\} \right)=0$. From $(\ref{b13})$, we have $-\delta\left( \left\{ -R_{ij} \right\} \right) + j\left(  \left\{ - G_{ijk} \right\} \right)=0$. As in the proof of Proposition $(\ref{b11})$, $-\delta\left(\left\{B_i\right\} \right)+ E_1'\left( \left\{ - R_{ij}  \right\} \right)=0$ and $E_2\left( \left\{ B_i \right\} \right)=0$. Hence we see that
{\small{\begin{align*}
\left( \left\{ B_i \right\}, \left\{ - R_{ij} \right\}, \left\{ - G_{ijk} \right\} \right) \in C^0\left( \mathcal{U}, \mathscr{H}om_{\mathcal{O}_X}\left(  \mathcal{N}_{\mathcal{F}_0}^* , \tilde{\mathcal{S}}^2 \right) \right) \bigoplus C^1\left( \mathcal{U}, \mathscr{H}om_{\mathcal{O}_X}\left( \mathcal{N}_{\mathcal{F}_0}^* , \Omega_X^1 \right) \right) \bigoplus C^2\left( \mathcal{U}, \mathscr{H}om_{\mathcal{O}_X}\left( \mathcal{N}_{\mathcal{F}_0}^* ,  \mathcal{N}_{\mathcal{F}_0}^*  \right) \right)
\end{align*}}}
defines a $2$-cocycle in the above \v Cech resolution of $\mathscr{H}om_{\mathcal{O}_X}\left( \mathcal{N}_{\mathcal{F}_0}^*, \mathcal{N}_{\mathcal{F}_0}^* \right)^\bullet$. We show that for another lifting $\tilde{w}_i'^\alpha$ of $w_i^\alpha$, and another lifting $\tilde{h}_{ij}'^{\alpha\beta}$ of $h_{ij}^{\alpha\beta},\alpha,\beta=1,...,p$, the associated $2$-cocycle $\beta:=\left(\left\{ B_i' \right\}, \left\{ - R_{ij}' \right\} , \left\{ - G_{ijk}' \right\} \right)$ is cohomologous to $\alpha:=\left( \left\{ B_i \right\}, \left\{ -R_{ij} \right\} , \left\{ -G_{ijk} \right\}\right)$. Then $\tilde{w}_i'^\alpha= \tilde{w}_i^\alpha + t V_i^\alpha$ for some $Y_i^\alpha$ and $\tilde{h}_{ij}'^{\alpha\beta}= \tilde{h}_{ij}^{\alpha\beta} + t y_{ij}^{\alpha\beta}$ for some $y_{ij}^{\alpha\beta}\in \Gamma\left( U_{ij}, \mathcal{O}_X \right)$. We define $Y_i \in \Gamma\left( U_i, \mathscr{H}om_{\mathcal{O}_X}\left( T_{\mathcal{F}_0}, T_X \right) \right)$ by $V_i:\Gamma\left( U_i, \mathcal{N}_{\mathcal{F}_0}^* \right)\to \Gamma\left( U_i, \Omega_X^1 \right), w_{0i}^\alpha \mapsto V_i^\alpha$ and define $y_{ij}\in \Gamma\left(U_{ij}, \mathscr{H}om_{\mathcal{O}_X}\left(\mathcal{N}_{\mathcal{F}_0}^*, \mathcal{N}_{\mathcal{F}_0}^* \right)  \right)$ by
\begin{align*}
y_{ij}:\Gamma\left( U_{ij}, \mathcal{N}_{\mathcal{F}_0}^* \right) &\to  \Gamma\left( U_{ij} , \mathcal{N}_{\mathcal{F}_0}^* \right) \\
  w_{0i}^\alpha &\mapsto \sum_{\beta=1}^q y_{ij}^{\alpha\beta} w_{0j}^\beta
\end{align*}
Then as in $(\ref{b10})$, we have $E_1'\left( \left\{ V_i \right\} \right)= B_i'- B_i$ and $(\ref{b14})$ implies that $\delta\left(\left\{ V_i \right\} \right) + j\left(\left\{ y_{ij} \right\} \right)= - R_{ij}'+ R_{ij}$. We note that in a similar way to $(\ref{ot47})$, we can show $\delta\left( \left\{ y_{ij} \right\} \right) = \left\{ - G_{ijk}' + G_{ijk} \right\}$. This proves that $\alpha$ is cohomologous to $\beta$. So given a small extension $e:0\to (t)\to \tilde{A} \to A \to 0$ and an infinitesimal foliated deformation $\eta$ of $\left( X, \mathcal{N}_{\mathcal{F}_0}^* \right)$ with $X$ fixed, we can associate an element $o_\eta(e):=$ the cohomology class of $\alpha \in \mathbb{H}^2\left( X, \mathscr{H}om_{\mathcal{O}_X}\left( \mathcal{N}_{\mathcal{F}_0}^*,  \mathcal{N}_{\mathcal{F}_0}^* \right)^\bullet \right)$. We note that $o_\eta(e)=0$ if and only if there exists $\left\{ N_i \right\} \in C^0\left( \mathcal{U}, \mathscr{H}om_{\mathcal{O}_X}\left( \mathcal{N}_{\mathcal{F}_0}^* , \Omega_X^1 \right)\right)$ where $N_i:\Gamma\left(U_i, \mathcal{N}_{\mathcal{F}_0}^* \right)\to \Gamma\left( U_i, \Omega_X^1 \right), w_{0i}^\alpha \mapsto N_i^\alpha$ and $u_{ij} \in \Gamma\left( U_{ij}, \mathscr{H}om_{\mathcal{O}_X}\left( \mathcal{N}_{\mathcal{F}_0}^* , \mathcal{N}_{\mathcal{F}_0}^* \right) \right)$ where $u_{ij}:\Gamma\left(U_{ij}, \mathcal{N}_{\mathcal{F}_0}^* \right) \to \Gamma\left( U_{ij}, \mathcal{N}_{\mathcal{F}_0}^* \right), w_{0i}^\alpha \mapsto \sum_{\beta=1}^q u_{ij}^{\alpha\beta} w_{0j}^\beta$ such that $N_j- N_i + j\left(u_{ij}\right)= - R_{ij}$ and $E_1'\left( N_i\right)=B_i$ and $u_{ij}- u_{ik} + u_{jk}= - G_{ijk}$. Then as in $(\ref{b16})$ we have $\tilde{w}_i^\alpha - tN_i^\alpha  = \sum_{\beta=1}^q \left(\tilde{h}_{ij}^{\alpha\beta}- tu_{ij}^{\alpha\beta} \right)\left(\tilde{w}_j^\beta - t N_j^\beta \right)$ and as in $(\ref{b12})$, we have  $\left(\tilde{w}_i^1 - t N_i^1 \right) \wedge \cdots \left( \tilde{w}_i^q - t N_i^q \right) \wedge d \left( \tilde{w}_i^\alpha - t N_i^\alpha \right)= 0, \alpha=1,...,q$, and in a similar way to $(\ref{ot32})$, we have $\tilde{h}_{ik}^{\alpha \gamma}- tu_{ik}^{\alpha \gamma} = \sum_{\beta=1}^q \left(\tilde{h}_{ij}^{\alpha \beta} - t u_{ij}^{\alpha \beta} \right)\left( \tilde{h}_{jk}^{\beta\gamma} - t u_{jk}^{\beta \gamma} \right)$. Hence a lifting $\tilde{\eta}$ of $\eta$ exits. This completes the proof of Proposition \ref{b15}.

\end{proof}

\subsection{Relationship between foliated deformations of $\left(X, \mathcal{N}_{\mathcal{F}_0}^* \right)$ with $X$ fixed and the underlying deformations of $\mathcal{N}_{\mathcal{F}_0}^*$} \

In the previous section we have shown that $\mathscr{H}om_{\mathcal{O}_X}\left( \mathcal{N}_{\mathcal{F}_0}^*,  \mathcal{N}_{\mathcal{F}_0}^* \right)^\bullet$ controls flat foliated deformations of $\left(X, \mathcal{N}_{\mathcal{F}_0}^* \right)$ with $X$ fixed. As is known, the underlying flat deformations of the locally free sheaf $\mathcal{N}_{\mathcal{F}_0}^*$ with $X$ fixed is controlled by $\mathscr{H}om_{\mathcal{O}_X}\left(\mathcal{N}_{\mathcal{F}_0}^*, \mathcal{N}_{\mathcal{F}_0}^* \right)$. Let us denote by $\textnormal{Def}_{\mathcal{N}_{\mathcal{F}_0}^*}$ the functor of Artin rings defined by
\begin{align*}
\textnormal{Def}_{\mathcal{N}_{\mathcal{F}_0}^*}:\bold{Art} &\to \left(\textnormal{Sets} \right) \\
 A &\mapsto \left\{ \textnormal{locally free sheaves on $X\times \textnormal{Spec}(A)$ inducing  $\mathcal{N}_{\mathcal{F}_0}^*$}   \right\}/\textnormal{isomorphisms}
\end{align*}
Then we have a morphism of deformation functors
\begin{align*}
\textnormal{FDef}_{\mathcal{N}_{\mathcal{F}_0}^*} \to \textnormal{Def}_{\mathcal{N}_{\mathcal{F}_0}^*}
\end{align*}
by forgetting foliated structures. We also have a canonical morphism of deformation functors 
\begin{align*}
\textnormal{FDef}_{\left( X, \mathcal{N}_{\mathcal{F}_0}^* \right)_{\textnormal{fixed}}}\to \textnormal{FDef}_{\mathcal{N}_{\mathcal{F}_0}^*}
\end{align*}
Then the relationship of those three functors $\textnormal{FDef}_{\left(X, \mathcal{N}_{\mathcal{F}_0}^* \right)}, \textnormal{FDef}_{ \mathcal{N}_{\mathcal{F}_0}^* }$ and $\textnormal{Def}_{\mathcal{N}_{\mathcal{F}_0}^*}$ is expressed by the following exact sequence of complexes of sheaves
{\small{\begin{center}
$\begin{CD}
@. \cdots @. \cdots @. \cdots \\
@. @AE_3AA @AE_3AA @AAA \\
0 @>>> \mathscr{H}om_{\mathcal{O}_X}\left( \mathcal{N}_{\mathcal{F}_0}^* , \tilde{\mathcal{S}}^3 \right) @>>> \mathscr{H}om_{\mathcal{O}_X}\left(  \mathcal{N}_{\mathcal{F}_0}^*, \tilde{\mathcal{S}}^3 \right) @>>> 0 @>>> 0\\
@. @AE_2AA @AE_2AA @AAA \\
0 @>>> \mathscr{H}om_{\mathcal{O}_X}\left( \mathcal{N}_{\mathcal{F}_0}^* , \tilde{\mathcal{S}}^2 \right) @>>> \mathscr{H}om_{\mathcal{O}_X}\left( \mathcal{N}_{\mathcal{F}_0}^* , \tilde{\mathcal{S}}^2 \right) @>>> 0 @>>> 0\\
@. @AE_1'AA @AE_1'AA @AAA \\
0 @>>> \mathscr{H}om_{\mathcal{O}_X}\left(\mathcal{N}_{\mathcal{F}_0}^*, \Omega_X^1  \right) @>>>  \mathscr{H}om_{\mathcal{O}_X}\left( \mathcal{N}_{\mathcal{F}_0}^* , \Omega_X^1 \right) @>>> 0 @>>> 0\\
@. @AAA @A jAA @AAA \\
0 @>>> 0 @>>> \mathscr{H}om_{\mathcal{O}_X}\left( \mathcal{N}_{\mathcal{F}_0}^*, \mathcal{N}_{\mathcal{F}_0}^* \right) @>>> \mathscr{H}om_{\mathcal{O}_X}\left( \mathcal{N}_{\mathcal{F}_0}^*,  \mathcal{N}_{\mathcal{F}_0}^* \right) @>>> 0
\end{CD}$
\end{center}}}
Then from the long exact sequence of cohomology groups, we have
\begin{align*}
 \mathbb{H}^0\left(X, \mathscr{H}om_{\mathcal{O}_X}\left(\mathcal{N}_{\mathcal{F}_0}^*, \Omega_X^1 \right)^\bullet \right) \to \mathbb{H}^1\left( X, \mathscr{H}om_{\mathcal{O}_X}\left( \mathcal{N}_{\mathcal{F}_0}^* , \mathcal{N}_{\mathcal{F}_0}^* \right)^\bullet \right) \to H^1\left(X, \mathscr{H}om_{\mathcal{O}_X}\left(\mathcal{N}_{\mathcal{F}_0}^*, \mathcal{N}_{\mathcal{F}_0}^* \right) \right)
\end{align*}
which gives a morphism of first-order deformations 
\begin{align*}
\textnormal{FDef}_{\left( X, \mathcal{N}_{\mathcal{F}_0}^* \right)_{\textnormal{fixed}}}\left(k[\epsilon] \right)\to \textnormal{FDef}_{ \mathcal{N}_{\mathcal{F}_0}^*} \left( k[\epsilon] \right) \to \textnormal{Def}_{\mathcal{N}_{\mathcal{F}_0}^*} \left( k[\epsilon ] \right)
\end{align*}

and the long exact sequence also gives an exact sequence 
\begin{align*}
\mathbb{H}^1\left(X, \mathscr{H}om_{\mathcal{O}_X}\left( \mathcal{N}_{\mathcal{F}_0}^*, \Omega_X^1 \right)^\bullet \right) \to \mathbb{H}^2\left( X, \mathscr{H}om_{\mathcal{O}_X}\left( \mathcal{N}_{\mathcal{F}_0}^* , \mathcal{N}_{\mathcal{F}_0}^* \right)^\bullet \right) \to H^2\left(X, \mathscr{H}om_{\mathcal{O}_X}\left( \mathcal{N}_{\mathcal{F}_0}^* , \mathcal{N}_{\mathcal{F}_0}^* \right) \right)
\end{align*}
which gives obstruction maps of morphisms of deformation functors $\textnormal{FDef}_{\left( X, \mathcal{N}_{\mathcal{F}_0}^* \right)_{\textnormal{fixed}}} \to \textnormal{FDef}_{\mathcal{N}_{\mathcal{F}_0}^*} \to \textnormal{Def}_{\mathcal{N}_{\mathcal{F}_0}^*}$. 

\subsection{Relationship between foliated deformations of $\left(X, \mathcal{N}_{\mathcal{F}_0}^* \right)$ with $X$ fixed and the underlying flat deformations of $\frac{\Omega_X^1}{ \mathcal{N}_{\mathcal{F}_0}^* }$}\

By definition of foliated deformations of $\left(X, \mathcal{N}_{\mathcal{F}_0}^* \right)$ with $\mathcal{N}_{\mathcal{F}_0}^*$ locally free under the trivial deformation of $X$, it also involves flat deformations of the underlying coherent sheaf $\frac{\Omega_X^1}{\mathcal{N}_{\mathcal{F}_0}^*}$. We denote by $\textnormal{Def}_{\frac{\Omega_X^1}{\mathcal{N}_{\mathcal{F}_0}^*}}$ the functor of Artin rings of flat deformations of $\Omega_X^1/ \mathcal{N}_{\mathcal{F}_0}^*$ on the fixed underlying nonsingular variety $X$. Then we have a natural morphism of deformation functor 
\begin{align*}
\textnormal{FDef}_{\mathcal{N}_{\mathcal{F}_0}^*} \to \textnormal{Def}_{\frac{\Omega_X^1}{ \mathcal{N}_{\mathcal{F}_0}^*}}
\end{align*}
We will identify the relationships between first-order deformations and obstructions of two deformation functors. This subsection is a direct consequence of \cite{Har10} Theorem 7.2 and Theorem 7.3 p.54-p.56.

We note that we have an exact sequence of complex of sheaves
{\small{\begin{center}
$\begin{CD}
@. \cdots @. \cdots @. \cdots \\
@. @AAA @AAA @AAA \\
0@>>> \mathscr{H}om_{\mathcal{O}_X}\left( \mathcal{N}_{\mathcal{F}_0}^*, \tilde{\mathcal{S}}^3 \right) @>>> \mathscr{H}om_{\mathcal{O}_X}\left(  \mathcal{N}_{\mathcal{F}_0}^* ,  \tilde{\mathcal{S}}^3 \right) @>>> 0 @>>> 0 \\
@. @AE_2AA @AE_2AA @AAA \\
0@>>> \mathscr{H}om_{\mathcal{O}_X}\left( \mathcal{N}_{\mathcal{F}_0}^* ,  \tilde{\mathcal{S}}^2 \right) @>>> \mathscr{H}om_{\mathcal{O}_X}\left(  \mathcal{N}_{\mathcal{F}_0}^* ,  \tilde{\mathcal{S}}^2 \right) @>>> 0 @>>> 0 \\
@. @AAA @AE_1AA @AAA \\
0@>>> 0 @>>> \mathscr{H}om_{\mathcal{O}_X}\left(\mathcal{N}_{\mathcal{F}_0}^*, \frac{\Omega_X^1}{ \mathcal{N}_{\mathcal{F}_0}^*} \right) @>>> \mathscr{H}om_{\mathcal{O}_X}\left( \mathcal{N}_{\mathcal{F}_0}^*,  \frac{\Omega_X^1}{\mathcal{N}_{\mathcal{F}_0}^*}\right) @>>> 0
\end{CD}$
\end{center}}}  

Then we have the natural indued morphisms
\begin{align}
\mathbb{H}^0\left(X, \mathscr{H}om_{\mathcal{O}_X}\left(  \mathcal{N}_{\mathcal{F}_0}^*, \Omega_X^1 /  \mathcal{N}_{\mathcal{F}_0}^* \right)^\bullet \right) \to H^0\left(X, \mathscr{H}om_{\mathcal{O}_X}\left( \mathcal{N}_{\mathcal{F}_0}^*, \Omega_X^1 / \mathcal{N}_{\mathcal{F}_0}^* \right) \right)\cong \textnormal{Hom}\left( \mathcal{N}_{\mathcal{F}_0}^* , \Omega_X^1 / \mathcal{N}_{\mathcal{F}_0}^* \right) \label{b20} \\
\mathbb{H}^1\left(X, \mathscr{H}om_{\mathcal{O}_X}\left( \mathcal{N}_{\mathcal{F}_0}^* ,  \Omega_X^1 / \mathcal{N}_{\mathcal{F}_0}^* \right)^\bullet \right) \to H^1\left(X, \mathscr{H}om_{\mathcal{O}_X}\left( \mathcal{N}_{\mathcal{F}_0}^*, \Omega_X^1 / \mathcal{N}_{\mathcal{F}_0}^* \right) \right)\cong \textnormal{Ext}^1\left( \mathcal{N}_{\mathcal{F}_0}^* , \Omega_X^1 / \mathcal{N}_{\mathcal{F}_0}^* \right) \notag
\end{align}
Here we note that since $\mathcal{N}_{\mathcal{F}_0}^*$ is locally free, we have
\begin{align*}
\textnormal{Ext}^i\left( \mathcal{N}_{\mathcal{F}_0}^* , \Omega_X^1 / \mathcal{N}_{\mathcal{F}_0}^* \right)\cong \textnormal{Ext}^i\left(\mathcal{O}_X, \left( \mathcal{N}_{\mathcal{F}_0}^*  \right)^*\otimes \Omega_X^1 / \mathcal{N}_{\mathcal{F}_0}^* \right) \cong H^i\left(X, \mathscr{H}om_{\mathcal{O}_X}\left( \mathcal{N}_{\mathcal{F}_0}^*, \Omega_X^1/  \mathcal{N}_{\mathcal{F}_0}^* \right) \right)
\end{align*}

On the other hand, from the exact sequence $0\to \mathcal{N}_{\mathcal{F}_0}^* \to  \Omega_X^1 \to  \Omega_X^1 / \mathcal{N}_{\mathcal{F}_0}^* \to 0$, we have the exact sequence
\begin{align}
0&\to \textnormal{Hom}\left( \Omega_X^1  /  \mathcal{N}_{\mathcal{F}_0}^* , \Omega_X^1 /  \mathcal{N}_{\mathcal{F}_0}^*  \right)\to  \textnormal{Hom}\left( \Omega_X^1 ,  \Omega_X^1 /  \mathcal{N}_{\mathcal{F}_0}^*  \right)\to  \textnormal{Hom}\left( \mathcal{N}_{\mathcal{F}_0}^* , \Omega_X^1 / \mathcal{N}_{\mathcal{F}_0}^* \right)\notag\\
 &\to \textnormal{Ext}^1\left( \Omega_X^1 / \mathcal{N}_{\mathcal{F}_0}^* , \Omega_X^1 / \mathcal{N}_{\mathcal{F}_0}^* \right) \to \textnormal{Ext}^1\left( \Omega_X^1, \Omega_X^1 / \mathcal{N}_{\mathcal{F}_0}^* \right) \to \textnormal{Ext}^1\left( \mathcal{N}_{\mathcal{F}_0}^*  , \Omega_X^1 / \mathcal{N}_{\mathcal{F}_0}^* \right) \label{b21}\\
 &\to \textnormal{Ext}^2\left( \Omega_X^1 / \mathcal{N}_{\mathcal{F}_0}^* , \Omega_X^1 / \mathcal{N}_{\mathcal{F}_0}^* \right) \to \textnormal{Ext}^2\left( \Omega_X^1 , \Omega_X^1 / \mathcal{N}_{\mathcal{F}_0}^* \right) \to \textnormal{Ext}^2\left( \mathcal{N}_{\mathcal{F}_0}^* , \Omega_X^1 / \mathcal{N}_{\mathcal{F}_0}^* \right)\to \cdots \notag
\end{align}
Now we assume that $X$ is projective. Then we know that $\textnormal{Def}_{\Omega_X^1 / \mathcal{N}_{\mathcal{F}_0}^*}\left(k[\epsilon]\right)\cong \textnormal{Ext}^1\left( \Omega_X^1 / \mathcal{N}_{\mathcal{F}_0}^* , \Omega_X^1 / \mathcal{N}_{\mathcal{F}_0}^* \right)$ and the obstruction space for $\textnormal{Def}_{\Omega_X^1 / \mathcal{N}_{\mathcal{F}_0}^* }$ is  $\textnormal{Ext}^2\left( \Omega_X^1 /  \mathcal{N}_{\mathcal{F}_0}^* , \Omega_X^1 / \mathcal{N}_{\mathcal{F}_0}^* \right)$. By combining $(\ref{b20})$ and $(\ref{b21})$, we have
\begin{align*}
\mathbb{H}^0\left(X, \mathscr{H}om_{\mathcal{O}_X}\left( \mathcal{N}_{\mathcal{F}_0}^*, \Omega_X^1 / \mathcal{N}_{\mathcal{F}_0}^* \right)^\bullet \right) \to \textnormal{Hom}\left( \mathcal{N}_{\mathcal{F}_0}^* , \Omega_X^1  / \mathcal{N}_{\mathcal{F}_0}^* \right) \to \textnormal{Ext}^1\left( \Omega_X^1 / \mathcal{N}_{\mathcal{F}_0}^* , \Omega_X^1 / \mathcal{N}_{\mathcal{F}_0}^* \right)
\end{align*}
which gives a morphism of first order deformations between two deformation functors
\begin{align*}
\textnormal{FDef}_{\mathcal{N}_{\mathcal{F}_0}^*}\left(k[\epsilon]\right)\to \textnormal{Def}_{\Omega_X^1 / \mathcal{N}_{\mathcal{F}_0}^* }\left(k[\epsilon] \right)
\end{align*}

On the other hand, we have
\begin{align*}
\mathbb{H}^1\left(X, \mathscr{H}om_{\mathcal{O}_X}\left( \mathcal{N}_{\mathcal{F}_0}^* , \Omega_X^1 / \mathcal{N}_{\mathcal{F}_0}^* \right)^\bullet \right)\to \textnormal{Ext}^1\left( \mathcal{N}_{\mathcal{F}_0}^* , \Omega_X^1 / \mathcal{N}_{\mathcal{F}_0}^*  \right) \to \textnormal{Ext}^2\left( \Omega_X^1 / \mathcal{N}_{\mathcal{F}_0}^* , \Omega_X^1  / \mathcal{N}_{\mathcal{F}_0}^* \right)
\end{align*}
which gives the obstruction map between two deformation functors.

\section{Flat foliated deformations of $\left(X, \mathcal{N}_{\mathcal{F}_0}^* \right)$ with $\mathcal{N}_{\mathcal{F}_0}^*$ locally free}

Let $\left( X, \mathcal{N}_{\mathcal{F}_0}^*\right)$ be a foliated nonsingular variety over $k$ with $\mathcal{N}_{\mathcal{F}_0}^*$ locally free. In this section we now vary the underlying nonsingular variety $X$. We recall that

\begin{definition}[Flat foliated deformations in terms of cotangent sheaves]
Let $\left( X, \mathcal{N}_{\mathcal{F}_0}^* \right)$ be a foliated nonsingular variety with $\mathcal{N}_{\mathcal{F}_0}^*$ locally free. An infinitesimal foliated deformation $\eta$ of $\left( X, \mathcal{N}_{\mathcal{F}_0}^* \right)$ over $A\in \bold{Art}$ is a cartesian diagram of schemes
\begin{center}
$\eta:\begin{CD}
\left( X, \mathcal{N}_{\mathcal{F}_0}^* \right) @>>>\left( \mathcal{X}, \mathcal{N}_\mathcal{F}^* \right)\\
@VVV @VV \pi V\\
\textnormal{Spec}(k)@>>> \textnormal{Sec}(A)
\end{CD}$
\end{center}
where $\pi$ is flat, $\mathcal{N}_\mathcal{F}^*$ is a locally free subsheaf of $\Omega_{\mathcal{X}/\textnormal{Spec}(A)}^1$ such that $d\mathcal{N}_{\mathcal{F}}^* \subset \mathcal{N}_{\mathcal{F}}^*\bigwedge \Omega_{\mathcal{X}/A}^1$ at $x\in X-\textnormal{Sing}\left(\mathcal{F}_0 \right)$, and $\left( \mathcal{X}, \mathcal{N}_\mathcal{F}^* \right)$ induces $\left( X, \mathcal{N}_{\mathcal{F}_0}^* \right)$. Two infinitesimal foliated deformations $\left( \mathcal{X}, \mathcal{N}_\mathcal{F}^* \right)$ and $\left( \mathcal{X}', \mathcal{N}_{\mathcal{F}'}^* \right)$ of $\left( X, \mathcal{N}_{\mathcal{F}_0}^* \right)$ over $A$ are isomorphic if there is an isomorphism $\phi:\mathcal{X}\to \mathcal{X}'$ over $A$ inducing the identity on $\left( X, \mathcal{N}_{\mathcal{F}_0}^*\right)$ such that $\phi^* \mathcal{N}_{\mathcal{F}'}^*= \mathcal{N}_\mathcal{F}^*$. Then we have a functor of Artin rings
\begin{align*}
\textnormal{FDef}_{\left( X, \mathcal{N}_{\mathcal{F}_0}^* \right)}: \bold{Art} &\to \left(\textnormal{Sets} \right) \\
 A &\mapsto \left\{\textnormal{Infinitesimal foliated deformations of $\left( X , \mathcal{N}_{\mathcal{F}_0}^* \right)$ } \right\}/\textnormal{isomorphisms}
\end{align*}
\end{definition}

We will show that the dual leaf complex (see Part I)
\begin{align}\label{b22}
\mathcal{N}_{\mathcal{F}_0}^{*\bullet}: T_X \xrightarrow{E_0} \mathscr{H}om_{\mathcal{O}_X}\left( \mathcal{N}_{\mathcal{F}_0}^*, \frac{\Omega_X^1}{\mathcal{N}_{\mathcal{F}_0}^* } \right) \xrightarrow{E_1} \mathscr{H}om_{\mathcal{O}_X}\left( \mathcal{N}_{\mathcal{F}_0}^*, \tilde{\mathcal{S}}^2 \right) \xrightarrow{E_2}  \mathscr{H}om_{\mathcal{O}_X}\left( \mathcal{N}_{\mathcal{F}_0}^*, \tilde{\mathcal{S}}^3 \right)  \xrightarrow{E_3} \cdots
\end{align}
controls foliated deformations of $\left( X, \mathcal{N}_{\mathcal{F}_0}^* \right)$ in the language of functors of Artin rings.

\begin{lemma}[compare \cite{Ser06} Lemma 1.2.6 p.26]\label{nb2}
Let $B_0$ be a $k$-algebra with an element $w_0\in \Omega_{B_0/k}^1$, and $e: 0\to (t) \to \tilde{A} \to A  \to 0$ be a small extension in $\bold{Art}$. Let $w\in \Omega_{B_0/k}^1\otimes A$ which induces $w_0$, and $w_1,w_2 \in \Omega_{B_0/k}^1\otimes \tilde{A}$ which induces $w$, so that there exists $w' \in \Omega_{B_0/k}^1$ such that $w_1-w_2= tw'$. Then there is one to one correspondence
\begin{align*}
\left\{\textnormal{isomorphisms on $\left( B_0\otimes_k \tilde{A}, w_1 \right)$ and $\left( B_0\otimes_k \tilde{A}, w_2 \right)$ inducing the identity on $\left(B_0\otimes_k A, w \right)$} \right\} \\
\to \left\{ \textnormal{ $P\in \textnormal{Der}_k\left( B_0, B_0 \right)=\textnormal{Hom}_{B_0}\left( \Omega_{B_0/k}^1, B_0 \right) | w' +\mathcal{L}_P\left(w_0 \right) = 0 $ } \right\}
\end{align*}
In particular, when $w_1=w_2$, there is a canonical isomorphism of groups
\begin{align*}
\left\{ \textnormal{automorphisms on $\left( B_0\otimes_k \tilde{A}, w_1 \right)$ inducing the identity on $\left(B_0\otimes_k A, w \right)$} \right\} \to \left\{ P\in \textnormal{Der}_k\left( B_0 , B_0 \right) | \mathcal{L}_P \left( w_0 \right)=0  \right\}
\end{align*}
\end{lemma}

\begin{proof}
Let $\theta:\left( B_0\otimes_k \tilde{A} , w_1 \right)\to \left( B_0 \otimes_k \tilde{A} , w_2 \right)$ be an isomorphism sending $w_1$ to $w_2$ inducing the identity on $\left( B_0\otimes_k A, T \right)$. Then $\theta(x)=x + t P(dx)$ for some $P\in \textnormal{Der}_k\left( B_0, B_0 \right)= \textnormal{Hom}_{B_0}\left( \Omega_{B_0/k}^1, B_0 \right)$. We define the correspondence by $\theta\to P$. On the other hand, let $w_1= \sum_{i=1}^r \tilde{g}_i d \tilde{f}_i$ for some $g_i, f_i \in B_0\otimes_k \tilde{A}$. Then $w_0= \sum_{i=1}^r g_i df_i$ where $g_i,f_i\in B_0$ which are the restriction of $\tilde{g}_i, \tilde{f}_i$ on $B_0$. Since $w_1$ corresponds to $w_2$ via $\theta=\textnormal{Id}+ t P$, we have $w_2=\sum_{i=1}^r  \left( \tilde{g}_i+ t P\left(\tilde{g}_i \right) \right) d \left( \tilde{f}_i + t P\left( \tilde{f}_i \right) \right)= \sum_{i=1}^r\left( \tilde{g}_i d \tilde{f}_i + t P(g_i)df_i + t g_i P(f_i) \right)$, so that $w_2= w_1 + t \mathcal{L}_P\left( \sum_{i=1}^r g_idf_i \right)$. Hence $w'+ \mathcal{L}_P\left( w_0\right)=0$. Since $\theta$ is determined by $P$, the correspondence is one to one. Now assume that $w_1=w_2$. Then $\theta$ corresponds to $P$ with $\mathcal{L}_P (w_0)=0$ and the correspondence $\theta\to P$ is a group isomorphism.

\end{proof}

\begin{proposition}\label{nb1}
Let $\left( X, \mathcal{N}_{\mathcal{F}_0}^* \right)$ be a foliated nonsingular variety with $\mathcal{N}_{\mathcal{F}_0}^*$ locally free.
\begin{enumerate}
\item there is a $1-1$ correspondence
\begin{align*}
\textnormal{FDef}_{\left(X, \mathcal{N}_{\mathcal{F}_0}^*\right)}\left(k[\epsilon] \right)\cong \mathbb{H}^1\left(X, \mathcal{N}_{\mathcal{F}_0}^{* \bullet } \right)
\end{align*}
\item Let $A\in \bold{Art}$ and $\eta=\left( \mathcal{X}, \mathcal{N}_\mathcal{F}^* \right)$ be an infinitesimal foliated deformations of $\left( X, \mathcal{N}_{\mathcal{F}_0}^* \right)$ over $A$. To every small extension $0\to (t)\to \tilde{A} \to A \to 0$, we can associate an element $o_\eta(e)\in \mathbb{H}^2\left( X,  \mathcal{N}_{\mathcal{F}_0}^{* \bullet } \right)$ called the obstruction lifting of $\eta$ to $\tilde{A}$, which is $0$ if and only if a lifting of $\eta$ to $\tilde{A}$ exists.
\end{enumerate}
\end{proposition}

\begin{proof}
Given a first-order foliated deformation $\left( \mathcal{X}, \mathcal{N}_\mathcal{F}^* \right)$ of $\left( X, \mathcal{N}_{\mathcal{F}_0}^* \right)$ over $k[\epsilon]$, we choose an affine open covering $\mathcal{U}=\left\{ U_i \right\}$ of $X$ such that $\Gamma\left( U_i, \mathcal{N}_{\mathcal{F}_0}^* \right)$ is generated by $w_{0i}^1,..., w_{0i}^q\in \Gamma\left( U_i , \Omega_X^1 \right)$ with $w_{0i}^\alpha = \sum_{\beta=1}^q h_{0ij}^{\alpha\beta} w_{0j}^\beta$ for some $h_{0ij}^{\alpha\beta}\in \Gamma\left( U_{ij}, \mathcal{O}_X \right)$ and $w_{0i}^1\wedge \cdots \wedge w_{0i}^q\wedge dw_{0i}^\alpha=0,\alpha=1,...,q$, and we have $\mathcal{X}|_{U_i}\cong U_i\times \textnormal{Spec}\left( k[\epsilon] \right)$ is trivial for all $i$ with the induced foliated structure $\mathcal{N}_{\mathcal{F}_i}^*$ on $U_i\times \textnormal{Spec}\left( k[\epsilon ] \right)$, and $\Gamma\left( U_i, \mathcal{N}_{\mathcal{F}_i}^* \right)$ is generated by $w_{0i}^\alpha + \epsilon w_i'^\alpha,\alpha=1,...,q$, where $w_i'^\alpha\in \Gamma\left( U_i, \Omega_X^1 \right)$, and $\left( w_{0i}^1 + \epsilon w_i'^1 \right)\wedge \cdots \wedge \left( w_{0i}^q + \epsilon w_i'^q \right)\wedge d\left( w_{0i}^\alpha + \epsilon w_i'^\alpha \right)=0$. For all $i$, we have foliated isomorphisms
\begin{align*}
\theta_i : \left( U_i\times \textnormal{Spec}\left(k[\epsilon] \right) , \mathcal{N}_{\mathcal{F}_i}^* \right)\to \left( \mathcal{X}|_{U_i} ,  \mathcal{N}_{\mathcal{F}}^*|_{U_i} \right)
\end{align*}
Then for each $i,j, \theta_{ij}:=\theta_j^{-1}\theta_i:\left(U_{ij}\times \textnormal{Spec}\left(k[\epsilon] \right) , \mathcal{N}_{\mathcal{F}_i}^* \right)\to \left(U_{ij}\times \textnormal{Spec}\left( k[\epsilon ] \right) , \mathcal{N}_{\mathcal{F}_j}^* \right)$ is a foliated isomorphism inducing the identity on $\left( U_{ij} , \mathcal{N}_{\mathcal{F}_0}^* \right)$ by modulo $\epsilon$. Then by Lemma \ref{nb2}, $\theta_{ij}$ corresponds to $\textnormal{Id}+ \epsilon p_{ij}: \left( \Gamma\left( U_{ij}, \mathcal{O}_X \right), \mathcal{N}_{\mathcal{F}_j}^* \right)\to \left( \Gamma\left( U_{ij}, \mathcal{O}_X\right) , \mathcal{N}_{\mathcal{F}_i}^* \right)$ for some $p_{ij}\in \Gamma\left( U_{ij}, T_X\right)$ such that
\begin{align}\label{nb5}
\theta_{ji}^*\left( w_{0i}^\alpha + \epsilon w_i'^\alpha \right)= w_{0i}^\alpha + \epsilon w_i'^\alpha - \epsilon \mathcal{L}_{p_{ij}}\left( w_{0i}^\alpha \right) = \sum_{\beta=1}^q \left( h_{0ij}^{\alpha\beta} +\epsilon g_{ij}^{\alpha\beta} \right)\left( w_{0j}^\beta + \epsilon w_j'^\beta \right)
\end{align}
for some $g_{ij}^{\alpha\beta}\in \Gamma\left( U_{ij} , \mathcal{O}_X \right)$. Then $(\ref{nb5})$ implies that
\begin{align}\label{nb3}
w_i'^\alpha -\mathcal{L}_{p_{ij}}\left( w_{0i}^\alpha \right) = \sum_{\beta=1}^q h_{0ij}^{\alpha\beta} w_j'^\beta + \sum_{\beta=1}^q g_{ij}^{\alpha\beta} w_{0j}^\beta
\end{align}
We define $\overline{w_i'}\in \Gamma\left( U_i, \mathscr{H}om_{\mathcal{O}_X}\left( \mathcal{N}_{\mathcal{F}_0}^*, \frac{\Omega_M^1}{\mathcal{N}_{\mathcal{F}_0}^*} \right) \right)$ as in $(\ref{nb4})$. Then $(\ref{nb3})$ and $(\ref{b8})$ implies that $\left( \left\{ \overline{w_i'} \right\}, \left\{ p_{ij} \right\} \right)\in C^0\left( \mathcal{U}, \mathscr{H}om_{\mathcal{O}_X}\left( \mathcal{N}_{\mathcal{F}_0}^*, \frac{\Omega_X^1}{\mathcal{N}_{\mathcal{F}_0}^*} \right) \right) \bigoplus C^1\left( \mathcal{U}, T_X \right)$ defines a $1$-cocycle in the following \v Cech resolution of $\mathcal{N}_{\mathcal{F}_0}^{*\bullet}$:
{\Small{\begin{center}
$\begin{CD}
\cdots\\
@AE_3AA \\
C^0\left(\mathcal{U}, \mathscr{H}om_{\mathcal{O}_X}\left(\mathcal{N}_{\mathcal{F}_0}^*, \tilde{\mathcal{S}}^3 \right) \right)@>\delta>> \cdots\\
@AE_2AA @AE_2AA \\
C^0\left(\mathcal{U}, \mathscr{H}om_{\mathcal{O}_X}\left(\mathcal{N}_{\mathcal{F}_0}^*, \tilde{\mathcal{S}}^2 \right) \right) @>-\delta>> C^1\left(\mathcal{U}, \mathscr{H}om_{\mathcal{O}_X}\left( \mathcal{N}_{\mathcal{F}_0}^*, \tilde{\mathcal{S}}^2\right) \right) @>\delta>> \cdots\\ 
@AE_1AA @AE_1AA @AE_1AA \\
C^0\left(\mathcal{U}, \mathscr{H}om_{\mathcal{O}_X}\left( \mathcal{N}_{\mathcal{F}_0}^*, \frac{\Omega_X^1}{\mathcal{N}_{\mathcal{F}_0}^*}\right) \right)@>\delta>> C^1\left(\mathcal{U}, \mathscr{H}om_{\mathcal{O}_X}\left( \mathcal{N}_{\mathcal{F}_0}^*, \frac{\Omega_X^1}{\mathcal{N}_{\mathcal{F}_0}^*}\right) \right) @>-\delta>> C^2\left(\mathcal{U}, \mathscr{H}om_{\mathcal{O}_X}\left( \mathcal{N}_{\mathcal{F}_0}^*, \frac{\Omega_X^1}{\mathcal{N}_{\mathcal{F}_0}^*}\right) \right) \\
@AE_0AA @AE_0AA @AE_0AA \\
C^0\left(\mathcal{U},T_X\right)@>-\delta>> C^1\left(\mathcal{U},T_X \right) @>\delta>> C^2 \left(\mathcal{U},T_X\right) @>-\delta>> \cdots
\end{CD}$
\end{center}}}
We show that for two equivalent first-order foliated deformations of $\left( X, \mathcal{N}_{\mathcal{F}_0}^* \right)$, the associated $1$-cocycles are equivalent. Let $\left( \mathcal{X}', \mathcal{N}_{\mathcal{F}'}^* \right)$ be another foliated deformation of $\left( X, \mathcal{N}_{\mathcal{F}_0}^* \right)$ over $k[\epsilon]$ which induces $\left(\left\{\overline{w_i''}\right\}, \left\{ p_{ij}' \right\} \right)$ and $\Phi:\left( \mathcal{X}, \mathcal{N}_\mathcal{F}^* \right)\to \left( \mathcal{X}', \mathcal{N}_{\mathcal{F}'}^* \right)$ is a foliated isomorphism inducing the identity on $\left( X, \mathcal{N}_{\mathcal{F}_0}^* \right)$. Then there is a foliated isomorphism $\alpha_i$ on $U_i\times \textnormal{Spec}\left(k[\epsilon ] \right)$ given by 
\begin{align*}
\alpha_i:\left(U_i\times \textnormal{Spec}\left(k[\epsilon]\right) , \mathcal{N}_{\mathcal{F}_i}^* \right) \xrightarrow{\theta_i} \left( \mathcal{X}|_{U_i}, \mathcal{N}_{\mathcal{F}}^*|_{U_i} \right) \xrightarrow{\Phi |_{U_i}} \left( \mathcal{X}'|_{U_i}, \mathcal{N}_{\mathcal{F}'}^*|_{U_i} \right) \xrightarrow{\theta_i'^{-1}} \left( U_i \times \textnormal{Spec}\left(k[\epsilon] \right) , \mathcal{N}_{\mathcal{F}_i'}^* \right)
\end{align*}
which corresponds to $1+\epsilon a_i:\left(\Gamma\left( U_i, \mathcal{O}_X  \right)\otimes k[\epsilon], \mathcal{N}_{\mathcal{F}_i'}^* \right)\to \left( \Gamma\left(U_i, \mathcal{O}_X \right) \otimes k[\epsilon] , \mathcal{N}_{\mathcal{F}_i}^* \right)$ for some $a_i\in \Gamma\left( U_i, T_X \right)$ such that $w_{0i}^\alpha + \epsilon w_i'^\alpha - \epsilon \mathcal{L}_{a_i}\left( w_{0i}^\alpha \right) = \sum_{\beta=1}^q \left(\delta_\alpha^\beta + \epsilon g_i^{\alpha\beta} \right) \left( w_{0i}^\beta + \epsilon w_i''^\beta \right)$ for some $g_i^{\alpha\beta} \in \Gamma\left( U_i, \mathcal{O}_X \right)$ by Lemma \ref{nb1}. Then we have $\mathcal{L}_{ a_i,}\left( w_{0i}^\alpha\right)= w_i'^\alpha- w_i''^\alpha - \sum_{\beta=1}^q  g_i^{\alpha\beta} w_{0i}^\beta $. On the other hand, since $\theta_i'\alpha_i= \Phi |_{U_i} \theta_i$, we have $\alpha_j^{-1} \theta_{ij}' \alpha_i=\left(\theta_j' \alpha_j \right)^{-1}\left( \theta_i' \alpha_i \right) = \theta_j^{-1} \Phi |_{U_{ij}}^{-1} \Phi|_{U_{ij}} \theta_i= \theta_{ij}$ which corresponds to $\left(\textnormal{Id}+\epsilon a_i \right)\left(\textnormal{Id}+ \epsilon p_{ij}' \right)\left( \textnormal{Id} - \epsilon a_j \right)= \textnormal{Id}+t p_{ij}$, so that $a_i- a_j=  p_{ij}- p_{ij}'$. This implies that $\left( \left\{ \overline{w_i'} \right\}, \left\{ p_{ij} \right\} \right)$ and $\left( \left\{ \overline{w_i''} \right\}, \left\{  p_{ij}' \right\} \right)$ are cohomologous.

We define an inverse map. Given an element in $\mathbb{H}^1\left(X, \mathcal{N}_{\mathcal{F}_0}^{ * \bullet } \right)$, we represent it by a \v Cech 1-cocycle $\left(\left\{ \overline{w_i'} \right\}, \left\{ p_{ij} \right\} \right)$ for an affine open covering $\mathcal{U}=\{U_i\}$ of $X$ satisfying $(\ref{nb3})$ and $(\ref{b8})$ and $p_{ij}- p_{ik}+ p_{jk}=0$. By reversing the above process, we obtain a first-order foliated deformation of $\left( X, \mathcal{N}_{\mathcal{F}_0}^* \right)$.

Next we identify obstructions. Let $e:0\to (t)\to \tilde{A} \to A\to 0$ be a small extension in $\bold{Art}$ and $\eta:=\left( \mathcal{X}, \mathcal{N}_\mathcal{F}^* \right)$ be a foliated deformations of $\left( X, \mathcal{N}_{\mathcal{F}_0}^* \right)$ over $A$. Let $\mathcal{U}= \left\{U_i \right\}$ be an affine open covering of  $X$ such that we have foliated isomorphisms $\theta_i:\left( U_i \times\textnormal{Spec}(A), \mathcal{N}_{\mathcal{F}_i}^* \right)\to \left( \mathcal{X}|_{U_i}, \mathcal{N}_\mathcal{F}^*|_{U_i} \right)$, where $\mathcal{N}_{\mathcal{F}_i}^* \subset \Omega_X^1|_{U_i}\otimes A$ is the induced foliation on $U_i\times \textnormal{Spec}(A)$ from $\mathcal{N}_{\mathcal{F}}^*$ via $\theta_i$. Let $\Gamma\left( U_i, \mathcal{N}_{\mathcal{F}_i}^*\right)\subset \Gamma\left(U_i, \Omega_X^1 \right)\otimes A$ be generated by $w_i^\alpha,\alpha=1,..., q$ which induces $w_{0i}^\alpha$ modulo $\mathfrak{m}_A$ with $w_i^1\wedge \cdots \wedge w_i^q \wedge dw_i^\alpha=0,\alpha=1,...,q$, and let $\theta_{ij}:=\theta_j^{-1} \theta_i:\left( U_{ij}\times \textnormal{Spec}(A), \mathcal{N}_{\mathcal{F}_j}^* \right)\to \left( U_{ij}\times \textnormal{Spec}(A), \mathcal{N}_{\mathcal{F}_i}^* \right)$ is a foliated isomorphism with $\theta_{jk}\theta_{ij}=\theta_{ik}$ on $U_{ijk}\times \textnormal{Spec}(A)$ and $\theta_{ij}^*\left( \mathcal{N}_{\mathcal{F}_j}^* \right)=\mathcal{N}_{\mathcal{F}_i}^*$, and $\theta_{ji}^*\left( w_i^\alpha \right)= \sum_{\beta=1}^q h_{ij}^{\alpha\beta} w_j^\beta$ for some $h_{ij}^{\alpha\beta}\in \Gamma\left( U_{ij}, \mathcal{O}_X \right)\otimes A$ which induces $w_{0i}^\alpha= \sum_{\beta=1}^q h_{0ij}^{\alpha\beta} w_{0j}^\beta$ modulo $\mathfrak{m}_A$. We let $f_{ij}:\left(\Gamma\left( U_{ij}, \mathcal{O}_X \right)\otimes A, \mathcal{N}_{\mathcal{F}_j}^* \right) \to \left( \Gamma\left( U_{ij}, \mathcal{O}_X \right)\otimes A,  \mathcal{N}_{\mathcal{F}_i}^* \right)$ be the ring homomorphism corresponding to $\theta_{ij}$. Then we have $f_{ij}f_{jk}=f_{ik}$. To give a lifting $\tilde{\eta}$ of $\eta$ to $\tilde{A}$ is equivalent to give a collection $\left\{\mathcal{N}_{\tilde{\mathcal{F}}_i}^* \right\}$, where $\Gamma\left( U_i, \mathcal{N}_{\tilde{\mathcal{F}}_i}^* \right)\subset \Gamma\left( U_i, \Omega_X^1 \right)\otimes \tilde{A}$ is generated by $\tilde{w}_i^\alpha\in \Gamma\left( U_i , \Omega_X^1 \right)\otimes \tilde{A},\alpha=1,...,q$ which induces ${w}_i^\alpha$ such that $\tilde{w}_i^1\wedge \cdots \wedge \tilde{w}_i^q \wedge d\tilde{w}_i^\alpha=0,\alpha=1,...,q$, and a collection of foliated isomorphisms $\left\{ \tilde{\theta}_{ij}\right\}$ which induces $\theta_{ij}$, where $\tilde{\theta}_{ij}:\left( U_{ij}\times \textnormal{Spec}\left( \tilde{A} \right), \mathcal{N}_{\tilde{\mathcal{F}}_i}^* \right)\to \left( U_{ij}\times \textnormal{Spec}\left( \tilde{A} \right), \mathcal{N}_{\tilde{\mathcal{F}}_j}^* \right)$ such that $\tilde{\theta}_{ji}^*\left( \tilde{w}_i^\alpha \right) = \sum_{\beta=1}^q \tilde{h}_{ij}^{\alpha\beta} \tilde{w}_j^\beta$ for some $\tilde{h}_{ij}^{\alpha\beta}\in \Gamma\left( U_{ij}, \mathcal{O}_X \right)\otimes \tilde{A}$ which induces $h_{ij}^{\alpha\beta}$, and $\tilde{\theta}_{jk}\tilde{\theta}_{ij}=\tilde{\theta}_{ik}$ (equivalently, $\tilde{f}_{ij}\tilde{f}_{jk}= \tilde{f}_{ik}$). From such data, we can glue together $\left(U_i\times \textnormal{Spec}(\tilde{A}), \mathcal{N}_{\tilde{\mathcal{F}}_i}^*\right)$ to make a foliated deformation $\left(\tilde{\mathcal{X}}, \mathcal{N}_{\tilde{\mathcal{F}}}^* \right)$ inducing $\left(\mathcal{X}, \mathcal{N}_\mathcal{F}^*\right)$. Now given a foliated deformation $\eta=\left(\mathcal{X}, \mathcal{N}_\mathcal{F}^*\right)$ over $A$ as above, we associate an element $o_\xi(e)\in \mathbb{H}^2 \left(X, \mathcal{N}_{\mathcal{F}_0}^{* \bullet } \right)$. In order to see if there exists a lifting of $\left(\mathcal{X}, \mathcal{N}_\mathcal{F}^* \right)$ to $\tilde{A}$, we choose arbitrary automorphism $\left\{\tilde{\theta}_{ij}\right\}$ which restricts to $\theta_{ij}$ on $U_{ij}\times \textnormal{Spec}\left(A\right)$ and arbitrary lifting $\tilde{w}_i^\alpha\in \Gamma\left(U_i, \Omega_X^1 \right)\otimes \tilde{A}$ which restricts to $w_i^\alpha$ and arbitrary lifting $\tilde{h}_{ij}^{\alpha\beta} \in \Gamma\left( U_{ij}, \mathcal{O}_X \right)\otimes \tilde{A}$ of $h_{ij}^{\alpha\beta},\alpha,\beta=1,...,q$. Since $\tilde{\theta}_{ijk}=\tilde{\theta}_{ki}\tilde{\theta}_{jk}\tilde{\theta}_{ij}$ (equivalently $\tilde{f}_{ijk}=\tilde{f}_{ij}\tilde{f}_{jk}\tilde{f}_{ki}$) is an automorphism on $U_{ijk}\times \textnormal{Spec}\left(\tilde{A}\right)$ inducing the identity on $U_{ijk}\times \textnormal{Spec}(A)$, $\tilde{\theta}_{ijk}$ corresponds to $\tilde{d}_{ijk}\in \Gamma\left(U_{ijk}, T_X\right)$ with $\tilde{d}_{jkl}-\tilde{d}_{ikl}+\tilde{d}_{ijl}-\tilde{d}_{jkl}=0$, so that $\delta(\{\tilde{d}_{ijk}\})=0$. On the other hand, we have
\begin{align}
\tilde{\theta}_{ji}^*\left( \tilde{w}_i^\alpha \right) - \sum_{\beta=1}^q \tilde{h}_{ij}^{\alpha\beta} \tilde{w}_j^\beta & = tR_{ij}^\alpha\,\,\,\,\,\,\,\textnormal{for some}\,\,\, R_{ij}^\alpha\in \Gamma\left( U_{ij}, \Omega_X^1 \right)\label{si34}\\
\tilde{w}_i^1\wedge \cdots \wedge \tilde{w}_i^q\wedge d\tilde{w}_i^\alpha&=t B_i^\alpha \,\,\,\,\,\,\,\textnormal{for some}\,\,\, B_i^\alpha\in \Gamma\left(U_i, \tilde{\mathcal{S}}^2 \right) \notag
\end{align}
We define $\overline{R_{ij}}\in \Gamma\left(U_{ij}, \mathscr{H}om_{\mathcal{O}_X}\left(\mathcal{N}_{\mathcal{F}_0}^*, \frac{\Omega_X^1}{\mathcal{N}_{\mathcal{F}_0}^*} \right) \right)$ as in a similar way to $(\ref{nb6})$ and $B_i \in \Gamma\left(U_i, \mathscr{H}om_{\mathcal{O}_X}\left( \mathcal{N}_{\mathcal{F}_0}^*, \tilde{\mathcal{S}}^2 \right) \right)$ in a similar way to $(\ref{nb7})$. We note that by applying $\theta_{ij}^*$ to $\theta_{ji}^*\left( w_i^\alpha \right)= \sum_{\beta=1}^q h_{ij}^{\alpha\beta} w_j^\beta$, we have $w_i^\alpha= \sum_{\beta=1}^q f_{ij}\left(h_{ij}^{\alpha\beta} \right)\left(\sum_{\gamma=1}^q h_{ji}^{\beta \gamma} w_i^\gamma \right)$, so that $\sum_{\beta=1}^q f_{ij}\left( h_{ij}^{\alpha\beta} \right)h_{ji}^{\beta \gamma}= \delta_\alpha^\gamma$. Then we have $\sum_{\beta=1}^q \tilde{f}_{ij}\left( \tilde{h}_{ij}^{\alpha\beta} \right) \tilde{h}_{ji}^{\beta \gamma}= \delta_\alpha^\gamma + t c_{ij\alpha}^\gamma$ for some $c_{ij\alpha}^\gamma \in \Gamma\left( U_{ij}, \mathcal{O}_X \right)$. Then we can show $\left\{ \overline{R_{ij}} \right\} \in C^1\left( \mathcal{U}, \mathscr{H}om_{\mathcal{O}_X}\left( \mathcal{N}_{\mathcal{F}_0}^* , \frac{\Omega_X^1}{\mathcal{N}_{\mathcal{F}_0}^*} \right) \right)$ in a similar way to $(\ref{atb1})$. We claim that $\left( \left\{  B_i \right\}, \left\{ - \overline{R_{ij}} \right\} , \left\{ d_{ijk} \right\}\right)\in C^0\left( \mathcal{U}, \mathscr{H}om_{\mathcal{O}_X}\left( \mathcal{N}_{\mathcal{F}_0}^* , \tilde{\mathcal{S}}^2 \right) \right) $ $ \bigoplus C^1\left( \mathcal{U}, \mathscr{H}om_{\mathcal{O}_X}\left( \mathcal{N}_{\mathcal{F}_0}^* , \frac{\Omega_X^1}{\mathcal{N}_{\mathcal{F}_0}^*} \right) \right) $ $\bigoplus C^2\left(\mathcal{U}, T_X \right)$ defines a $2$-cocycle in the above \v Cech resolution of $\mathcal{N}_{\mathcal{F}_0}^{* \bullet }$. We have $-\delta\left( \left\{ d_{ijk} \right\} \right)=0$ as above and we have $E_2\left( B_i \right)=0$ as in the proof of Proposition \ref{b2}. We show that $-\delta\left( \left\{B_i \right\} \right) + D_1\left( \left\{ - \overline{R_{ij}} \right\} \right)=0$. In fact, we note that
{\small{\begin{align*}
&tB_i^\alpha= \tilde{\theta}_{ji}^*\left(tB_i^\alpha\right)= \tilde{\theta}_{ji}^*\left(\tilde{w}_i^1 \right) \wedge  \cdots \wedge \tilde{\theta}_{ji}^*\left(\tilde{w}_i^q \right) \wedge d \tilde{\theta}_{ji}^*\left( \tilde{w}_i^\alpha \right) \\
&= \left( tR_{ij}^1+\sum_{\beta_1=1}^q \tilde{h}_{ij}^{1\beta_1} \tilde{w}_j^{\beta_1} \right)\wedge \cdots \wedge \left( tR_{ij}^q+\sum_{\beta_q=1}^q \tilde{h}_{ij}^{q\beta_q} \tilde{w}_j^{\beta_q} \right) \wedge d\left( tR_{ij}^1+\sum_{\beta=1}^q \tilde{h}_{ij}^{\alpha\beta} \tilde{w}_j^{\beta}   \right)\\
&=t \sum_{\beta=1}^q w_{0i}^1 \wedge \cdots \wedge R_{ij}^\beta \wedge \cdots \wedge w_{0i}^q \wedge dw_{0i}^\alpha + t w_{0i}^1\wedge \cdots \wedge w_{0i}^q \wedge d R_{ij}^\alpha+ \det\left(\tilde{h}_{ij} \right)\tilde{w}_j^1\wedge \cdots \wedge \tilde{w}_j^q\wedge \left(\sum_{\beta=1}^q d \tilde{h}_{ij}^{\alpha\beta} \wedge \tilde{w}_j^\beta + \sum_{\beta=1}^q \tilde{h}_{ij}^{\alpha\beta} d\tilde{w}_j^\beta \right)\\
&=t \left( \sum_{\beta=1}^q w_{0i}^1 \wedge \cdots \wedge R_{ij}^\beta \wedge \cdots \wedge w_{0i}^q \wedge dw_{0i}^\alpha + w_{0i}^1\wedge \cdots \wedge w_{0i}^q \wedge d R_{ij}^\alpha + \sum_{\beta=1}^q h_{0ij}^{\alpha\beta}\det\left( h_{0ij} \right) B_j^\beta \right)
\end{align*}}}
This implies $-\delta\left( \left\{B_i \right\} \right) + D_1\left( \left\{ - \overline{R_{ij}} \right\} \right)=0$. On the other hand, we show that $-\delta\left(\left\{ - \overline{R_{ij}} \right\} \right) + E_0\left( \left\{ d_{ijk}\right\}\right)=0$. In fact, first we note that
\begin{align}
w_i^\alpha =\theta_{ij}^* \theta_{jk}^* &\theta_{ki}^*\left( w_i^\alpha \right) = \theta_{ij}^*\theta_{jk}^*\left(\sum_{\beta=1}^q h_{ik}^{\alpha\beta} w_k^\beta \right)= \theta_{ij}^*\left( \sum_{\beta, \gamma=1}^q f_{jk}\left(h_{ik}^{\alpha\beta} \right)  h_{kj}^{\beta \gamma} w_j^\gamma \right) = \sum_{\beta, \gamma, \eta=1}^q \left(f_{ij} f_{jk} \right)\left(h_{ik}^{\alpha\beta} \right)f_{ij}\left( h_{kj}^{\beta \gamma} \right) h_{ji}^{\gamma \eta} w_i^\eta \notag \\
&\Longrightarrow \sum_{\beta, \gamma=1}^q  \left(\tilde{f}_{ij} \tilde{f}_{jk} \right)\left(\tilde{h}_{ik}^{\alpha\beta} \right) \tilde{f}_{ij}\left( \tilde{h}_{kj}^{\beta \gamma} \right) \tilde{h}_{ji}^{\gamma \eta} = \delta_\alpha^\eta + t u_{ijk}^{\alpha \eta}\,\,\,\,\,\,\,\textnormal{for some}\,\,\,  u_{ijk}^{\alpha \eta} \in \Gamma\left( U_{ijk}, \mathcal{O}_X \right) \label{tn6}
\end{align}
Then we note that we have an isomorphism $\textnormal{Id}+ td_{ijk}:\Gamma\left( U_{ijk}, \mathcal{O}_X \right)\otimes \tilde{A}\to \Gamma\left(U_{ijk}, \mathcal{O}_X \right)\otimes\tilde{A}$ corresponding to $\tilde{\theta}_{ki}\tilde{\theta}_{jk} \tilde{\theta}_{ij}$ and so $\tilde{f}_{ij}\tilde{f}_{jk}\tilde{f}_{ki}=\textnormal{Id}+t d_{ijk}$. Then by Lemma \ref{nb2}, 
\begin{align*}
\tilde{w}_{i}^\alpha+t \mathcal{L}_{ d_{ijk}}\left( w_{0i}^\alpha \right)&=\tilde{\theta}_{ij}^*\tilde{\theta}_{jk}^*\tilde{\theta}_{ki}^* \left(\tilde{w}_i^\alpha\right)=\tilde{\theta}_{ij}^*\tilde{\theta}_{jk}^*\left(tR_{ik}^\alpha +\sum_{\beta=1}^q \tilde{h}_{ik}^{\alpha\beta} \tilde{w}_k^\beta \right)=\tilde{\theta}_{ij}^*\left(tR_{ik}^\alpha+\sum_{\beta=1}^q\tilde{h}_{jk}\left(\tilde{h}_{ik}^{\alpha\beta} \right)\left(tR_{kj}^\beta+\sum_{\gamma=1}^q \tilde{h}_{kj}^{\beta\gamma} \tilde{w}_j^\gamma \right) \right)\\
&=tR_{ik}^\alpha+t\sum_{\beta=1}^q h_{0ik}^{\alpha\beta} R_{kj}^\beta +\sum_{\beta,\gamma=1}^q \left(\tilde{f}_{ij}\tilde{f}_{jk}\right)\left( \tilde{h}_{ik}^{\alpha\beta}\right) \tilde{f}_{ij}\left( \tilde{h}_{kj}^{\beta\gamma}\right)\left(tR_{ji}^\gamma +\sum_{\eta=1}^q \tilde{h}_{ji}^{\gamma\eta} \tilde{w}_i^\eta \right)\\
&=t\left(R_{ik}^\alpha+ \sum_{\beta=1}^q h_{0ik}^{\alpha\beta} R_{kj}^\beta+\sum_{\gamma=1}^p r_{0ij}^{\alpha\gamma} R_{ji}^\gamma \right)+ \tilde{w}_i^\alpha+ t\sum_{\gamma=1}^q u_{ijk}^{\alpha\eta} w_{0i}^\eta
\end{align*}
This proves our claim. We show that for another lifting $\tilde{\theta}_{ij}'$ of $\theta_{ij}$, and another lifting $\tilde{w}_i'^\alpha$ of $w_i^\alpha$, and another lifting of $\tilde{h}_{ij}'^{\alpha\beta}$ of $h_{ij}^{\alpha\beta}$, the associated $2$-cocycle $\beta:=\left(\left\{ B_i' \right\}, \left\{ - \overline{R_{ij}'} \right\}, \left\{d_{ijk}' \right\} \right)$ is cohomologus to $\alpha:=\left( \left\{ B_i \right\}, \left\{ - \overline{ R_{ij}} \right\}, \left\{ d_{ijk} \right\} \right)$. Then $\tilde{f}_{ij}'=\tilde{f}_{ij} + t p_{ij}$ for some $p_{ij}\in \Gamma\left( U_{ij}, T_X \right)$ and $\tilde{w}_i'^\alpha= \tilde{w}_i^\alpha + t V_i^\alpha$ for some $V_i^\alpha \in \Gamma\left( U_i, \Omega_X^1 \right)$, and $\tilde{h}_{ij}'^{\alpha\beta}= \tilde{h}_{ij}^{\alpha\beta} + t y_{ij}^{\alpha\beta}$ for some $y_{ij}^{\alpha\beta}\in \Gamma\left( U_{ij}, \mathcal{O}_X \right)$. We define $\overline{V_i}\in \Gamma\left( U_i,  \mathscr{H}om_{\mathcal{O}_X}\left( \mathcal{N}_{\mathcal{F}_0}^*, \frac{\Omega_X^1}{\mathcal{N}_{\mathcal{F}_0}^*} \right) \right), w_{0i}^\alpha \mapsto \overline{V_i^\alpha}$. Then as in $(\ref{b10})$, we have $E_1\left(\overline{V_i} \right) =  B_i'-B_i$. We note that
\begin{align}\label{si38}
&t\left(- R_{ij}'^\alpha + R_{ij}^\alpha \right) =  - \tilde{\theta}_{ji}'^*\left( \tilde{w}_i'^\alpha \right) + \sum_{\beta=1}^q \tilde{h}_{ij}'^{\alpha\beta} \tilde{w}_j'^\beta + \tilde{\theta}_{ji}^*\left( \tilde{w}_i^\alpha \right) - \sum_{\beta=1}^q \tilde{h}_{ij}^{\alpha\beta} \tilde{w}_j^\beta  \\
&=\left(- \tilde{\theta}_{ji}'^* +\tilde{\theta}_{ji}^*\right) \left( \tilde{w}_i^\alpha \right)- t V_{i}^\alpha + \sum_{\beta=1}^q \left( \tilde{h}_{ij}^{\alpha\beta} +t y_{ij}^{\alpha\beta} \right)\left( \tilde{w}_j^\beta + t V_j^\beta \right) - \sum_{\beta=1}^q \tilde{h}_{ij}^{\alpha\beta} \tilde{w}_j^\beta \notag \\
&= - t \mathcal{L}_{p_{ji}}\left( w_{0i}^\alpha \right)+ t  \left(  \sum_{\beta=1}^q h_{0ij}^{\alpha\beta} V_j^\beta - V_i^\alpha \right) + t \sum_{\beta=1}^q y_{ij}^{\alpha\beta} w_{0j}^\beta = t\left( \mathcal{L}_{p_{ij}}\left( w_{0i}^\alpha \right) + \sum_{\beta=1}^q h_{0ij}^{\alpha\beta} V_j^\beta - V_i^\alpha + \sum_{\beta=1}^q y_{ij}^{\alpha\beta} w_{0j}^\beta \right) \notag
\end{align}
This implies that $\delta\left(\left\{ \overline{V_i}\right\} \right) + D_0\left( \left\{ p_{ij} \right\}\right)= \left\{ -R_{ij}' + R_{ij} \right\}$ and $\delta\left(\left\{ p_{ij} \right\} \right)= \left\{ d_{ijk}' - d_{ijk}\right\}$ (see \cite{Ser06}). This proves that $\alpha$ is cohomologous to $\beta$. So given a small extension $e:0\to (t)\to \tilde{A} \to A \to 0$ and an infinitesimal foliated deformation $\eta$ of $\left( X, \mathcal{N}_{\mathcal{F}_0}^* \right)$ over $A$, we can associate an element $o_\eta(e):=$ the cohomology class of $\alpha \in \mathbb{H}^2\left( X, \mathcal{N}_{\mathcal{F}_0}^{* \bullet } \right)$. We note that $o_\eta(e)=0$ if and only if there exists $\left(\left\{ \overline{ N_i} \right\}, \left\{ d_{ij} \right\} \right) \in C^0\left( \mathcal{U}, \mathscr{H}om_{\mathcal{O}_X}\left( \mathcal{N}_{\mathcal{F}_0}^* , \frac{ \Omega_X^1 }{ \mathcal{N}_{\mathcal{F}_0}^*} \right) \right) \bigoplus C^1\left( \mathcal{U}, T_X \right)$ where $\overline{N_i}:\Gamma\left( U_i, \mathcal{N}_{\mathcal{F}_0}^* \right) \to \Gamma\left( U_i, \frac{\Omega_X^1}{ \mathcal{N}_{\mathcal{F}_0}^*}  \right), w_{0i}^\alpha \mapsto \overline{N_i^\alpha}$ such that $E_1\left(N_i\right)= B_i$ and $\overline{N_j}- \overline{N_i} + E_0\left( d_{ij} \right)= - \overline{ R_{ij}}$ and $d_{ij}- d_{ik} + d_{jk}= d_{ijk}$, so that $\left(\tilde{f}_{ij}- t d_{ij} \right)\left(\tilde{f}_{jk}- t d_{jk} \right)= \tilde{f}_{ik} - td_{ik}$ as in $(\ref{atb2})$. Let $\Phi_{ij}:U_{ij}\times \textnormal{Spec}\left(\tilde{A}\right) \to  U_{ij} \times \textnormal{Spec}\left( \tilde{A}\right)$ be an isomorphism corresponding to $\tilde{f}_{ij} - t d_{ij}$. Then $\Phi_{jk}\Phi_{ij} = \Phi_{ik}$. We also note that
{\small{\begin{align}
&t\left( \mathcal{L}_{d_{ij}}\left(w_{0i}^\alpha \right) +\sum_{\beta=1}^q h_{0ij}^{\alpha\beta} N_j^\beta - N_i^\alpha +  \sum_{\beta=1}^q u_{ij}^{\alpha\beta} w_{0j}^\beta\right)= - \tilde{\theta}_{ji}^* \left(\tilde{w}_i^\alpha \right) + \sum_{\beta=1}^q \tilde{h}_{ij}^{\alpha\beta} \tilde{w}_j^\beta \notag \\
& \Longrightarrow \tilde{\theta}_{ji}^*\left( \tilde{w}_i^\alpha \right) - tN_i^\alpha + t \mathcal{L}_{d_{ij}}\left( w_{0i}^\alpha \right)  = \sum_{\beta=1}^q \left(\tilde{h}_{ij}^{\alpha\beta}- tu_{ij}^{\alpha\beta} \right)\left(\tilde{w}_j^\beta - t N_j^\beta \right)\,\,\,\textnormal{for some} \,\,\, u_{ij}^{\alpha\beta}\in \Gamma\left( U_{ij}, \mathcal{O}_X \right) \notag\\
&\Longrightarrow \Phi_{ji}^*\left( \tilde{w}_i^\alpha- tN_i^\alpha \right) = \sum_{\beta=1}^q \left(\tilde{h}_{ij}^{\alpha\beta}- tu_{ij}^{\alpha\beta} \right)\left(\tilde{w}_j^\beta - t N_j^\beta \right)
\end{align}}}
and as in $(\ref{b12})$, we have  $\left(\tilde{w}_i^1 - t N_i^1 \right) \wedge \cdots \left( \tilde{w}_i^q - t N_i^q \right) \wedge d \left( \tilde{w}_i^\alpha - t N_i^\alpha \right)= 0 , \alpha=1,...,q $. Hence a lifting $\tilde{\eta}$ of $\eta$ exits. This completes the proof of Proposition \ref{nb1}.
\end{proof}

\subsection{Relationship between foliated deformations of $\left(X, \mathcal{N}_{\mathcal{F}_0}^* \right)$ with $X$ fixed and flat foliated deformation of $\left(X, \mathcal{N}_{\mathcal{F}_0}^* \right)$ and underlying flat deformations of $X$}\

Let $\left( X, \mathcal{N}_{\mathcal{F}_0}^* \right)$ be a foliated nonsingular variety with $\mathcal{N}_{\mathcal{F}_0}^*$ locally free. Let $\textnormal{Def}_X$ be the functor of Artin rings of flat deformations of the underlying nonsingular variety $X$ (see \cite{Ser06}). As is known, flat deformations of the underlying nonsingular variety $X$ is controlled by the tangent sheaf $T_X$. We have a canonical morphism of deformation functors
\begin{align*}
\textnormal{FDef}_{\mathcal{N}_{\mathcal{F}_0}^*}\to \textnormal{FDef}_{\left(X, \mathcal{N}_{\mathcal{F}_0}^*\right)}
\end{align*}
and we have a morphism of deformation functors
\begin{align}\label{a15}
\textnormal{FDef}_{\left(X, \mathcal{N}_{\mathcal{F}_0}^*\right)}\to \textnormal{Def}_X
\end{align}
by forgetting foliated structures. Then the relationship of three deformation functors $\textnormal{FDef}_{\mathcal{N}_{\mathcal{F}_0}^*}, \textnormal{FDef}_{\left( X, \mathcal{N}_{\mathcal{F}_0}^* \right)}$ and $\textnormal{Def}_X$ is expressed by the following exact sequence of complexes of sheaves
{\Small{\begin{equation} \label{a18}
\begin{CD}
@. \cdots @. \cdots @. \cdots \\
@. @AE_3AA @AE_3AA @AAA \\
0@>>> \mathscr{H}om_{\mathcal{O}_X}\left( \mathcal{N}_{\mathcal{F}_0}^*,  \tilde{\mathcal{S}}^3 \right) @>>> \mathscr{H}om_{\mathcal{O}_X}\left( \mathcal{N}_{\mathcal{F}_0}^*, \tilde{\mathcal{S}}^3 \right) @>>> 0@>>> 0 \\
@. @AE_2AA @AE_2AA @AAA \\
0@>>> \mathscr{H}om_{\mathcal{O}_X}\left( \mathcal{N}_{\mathcal{F}_0}^* , \tilde{\mathcal{S}}^2 \right) @>>> \mathscr{H}om_{\mathcal{O}_X}\left( \mathcal{N}_{\mathcal{F}_0}^* ,  \tilde{\mathcal{S}}^2 \right) @>>> 0@>>> 0 \\
@. @AE_1AA @AE_1AA @AAA \\
0@>>> \mathscr{H}om_{\mathcal{O}_X}\left( \mathcal{N}_{\mathcal{F}_0}^*, \frac{\Omega_X^1}{\mathcal{N}_{\mathcal{F}_0}^*}\right) @>>> \mathscr{H}om_{\mathcal{O}_X}\left( \mathcal{N}_{\mathcal{F}_0}^* , \frac{\Omega_X^1}{\mathcal{N}_{\mathcal{F}_0}^*} \right) @>>> 0@>>> 0 \\
@. @AAA @AE_0AA @AAA \\
0 @>>> 0 @>>> T_X @>>> T_X @>>> 0
\end{CD}
\end{equation}}}

Then the long exact sequence induces maps of first-order deformations
{\small{\begin{align*}
\textnormal{FDef}_{\mathcal{N}_{\mathcal{F}_0}^*}\left(k[\epsilon]\right)\cong \mathbb{H}^0\left(X,\mathscr{H}om_{\mathcal{O}_X}\left( \mathcal{N}_{\mathcal{F}_0}^*, \frac{\Omega_X^1}{\mathcal{N}_{\mathcal{F}_0}^*} \right)^\bullet  \right) \to \textnormal{FDef}_{\left(X, \mathcal{N}_{\mathcal{F}_0}^*\right)}\left(k[\epsilon]\right)\cong  \mathbb{H}^1\left(X, \mathcal{N}_{\mathcal{F}_0}^{* \bullet } \right)\to \textnormal{Def}_X\left(k[\epsilon]\right)\cong H^1\left(X,T_X\right)
\end{align*}}}
and obstruction maps
\begin{align*}
\mathbb{H}^1\left( X, \mathscr{H}om_{\mathcal{O}_X}\left(\mathcal{N}_{\mathcal{F}_0}^*, \Omega_X^1/ \mathcal{N}_{\mathcal{F}_0}^* \right)^\bullet \right) \to \mathbb{H}^2\left(X, \mathcal{N}_{\mathcal{F}_0}^{*\bullet}\right) \to H^2\left(X, T_X\right)
\end{align*}
of morphisms of deformation functors $\textnormal{FDef}_{\mathcal{N}_{\mathcal{F}_0}^*}\to \textnormal{FDef}_{\left(X, \mathcal{N}_{\mathcal{F}_0}^* \right)} \to \textnormal{Def}_X$.

\subsection{Flat foliated deformations of $\left(X, \mathcal{N}_{\mathcal{F}_0}^* \right)$ with $\mathcal{N}_{\mathcal{F}_0}^*$ locally free in terms of vector bundles}\

The underlying deformation of locally free sheaf $\mathcal{N}_{\mathcal{F}_0}^*$ is implicit in the dual leaf complex $(\ref{b22})$. In order to make the underlying deformations of $\mathcal{N}_{\mathcal{F}_0}^*$ explicit, we define a complex of sheaves also controlling flat foliated deformations of $\left(X, \mathcal{N}_{\mathcal{F}_0}^* \right)$ and present foliated deformations of $\left(X, \mathcal{N}_{\mathcal{F}_0}^* \right)$ by using the complex. As is known, deformations of a pair of a nonsingular variety $X$ and a locally free sheaf $\mathcal{N}_{\mathcal{F}_0}^*$ is controlled by the Atiyah extension (see \cite{Hua95}) which is the exact sequence 
\begin{align}\label{si2}
0\to \mathscr{H}om_{\mathcal{O}_X}\left(\mathcal{N}_{\mathcal{F}_0}^*, \mathcal{N}_{\mathcal{F}_0}^* \right) \to \mathcal{E}_{\mathcal{N}_{\mathcal{F}_0}^*}\to T_X \to 0
\end{align}
where $\mathcal{E}_{\mathcal{N}_{\mathcal{F}_0}^*}$ is defined in the following way: let $\{U_i\}$ be an affine open covering of $X$ such that $\Gamma\left(U_i,\mathcal{N}_{\mathcal{F}_0}^* \right)\cong \bigoplus^q \Gamma \left(U_i, \mathcal{O}_X\right)$ is trivialized in terms of the basis $\left\{w_{0i}^\alpha \right\}_{\alpha=1,...,q}$ with the relation $w_{0i}^\alpha= \sum_{\beta=1}^q h_{0ij}^{\alpha \beta} w_{0j}^\beta,\alpha=1,...,q$. Then we define $\mathcal{E}_{\mathcal{N}_{\mathcal{F}_0}^*}$ locally on $U_i$
\begin{align*}
\mathcal{E}_{\mathcal{N}_{\mathcal{F}_0}^*}|_{U_i}\cong \mathscr{H}om_{\mathcal{O}_X}\left( \mathcal{N}_{\mathcal{F}_0}^*, \mathcal{N}_{\mathcal{F}_0}^* \right)|_{U_i} \bigoplus T_X |_{U_i}
\end{align*}
such that $\left(\phi_i, X_i\right)\in \mathscr{H}om_{\mathcal{O}_X}\left( \mathcal{N}_{\mathcal{F}_0}^*, \mathcal{N}_{\mathcal{F}_0}^* \right) |_{U_i} \bigoplus T_X |_{U_i}$ and $\left(\phi_j, X_j\right)\in \mathscr{H}om_{\mathcal{O}_X}\left(\mathcal{N}_{\mathcal{F}_0}^*, \mathcal{N}_{\mathcal{F}_0}^* \right) |_{U_j} \bigoplus T_X |_{U_j} $ are identified on $U_{ij}$ if $X_i=X_j=X$ and
\begin{align*}
\phi_i \left(w_{0i}^\alpha\right)=\phi_j \left(w_{0i}^\alpha\right)+\sum_{\beta=1}^q \left[X, h_{0ij}^{\alpha\beta}\right]w_{0j}^\beta
\end{align*}
Then $\mathcal{E}_{\mathcal{N}_{\mathcal{F}_0}^*}$ is well-defined. In fact, since $\phi_j \left(w_{0j}^\beta\right)=\phi_k \left(w_{0j}^\beta\right)+\sum_{\gamma=1}^q \left[X, h_{0jk}^{\beta \gamma} \right]w_{0k}^\gamma$, we have
\begin{align}\label{si4}
\phi_i \left(w_{0i}^\alpha \right)&=\phi_j \left(w_{0i}^\alpha \right)+\sum_{\beta=1}^q \left[ X,h_{0ij}^{\alpha\beta} \right] w_{0j}^\beta =\sum_{\beta=1}^q h_{0ij}^{\alpha\beta} \phi_j\left(w_{0j}^\beta\right) + \sum_{\beta,\gamma=1}^q \left[X, h_{0ij}^{\alpha\beta} \right]h_{0jk}^{\beta \gamma}w_{0k}^\gamma \\
&=\sum_{\beta=1}^q h_{0ij}^{\alpha\beta}\left( \phi_k(w_{0j}^\beta)+\sum_{\gamma=1}^q \left[X, h_{0jk}^{\beta\gamma} \right] w_{0k}^\gamma \right) + \sum_{\beta,\gamma=1}^q \left[X, h_{0ij}^{\alpha\beta} \right]h_{0jk}^{\beta \gamma}w_{0k}^\gamma =\phi_k \left(w_{0i}^\alpha\right)+ \sum_{\gamma=1}^q \left[X, h_{0ik}^{\alpha \gamma}\right]w_{0k}^\gamma \notag
\end{align}
Then we shall construct a complex of sheaves
\begin{align}
\mathcal{E}_{\mathcal{N}_{\mathcal{F}_0}^*}^\bullet:\mathcal{E}_{\mathcal{N}_{\mathcal{F}_0}^*}\xrightarrow{E_0'} \mathscr{H}om_{\mathcal{O}_X}\left( \mathcal{N}_{\mathcal{F}_0}^*, \Omega_X^1 \right) \xrightarrow{E_1'} \mathscr{H}om_{\mathcal{O}_X}\left(\mathcal{N}_{\mathcal{F}_0}^*, \tilde{\mathcal{S}}^2 \right) \xrightarrow{E_2} \mathscr{H}om_{\mathcal{O}_X}\left(\mathcal{N}_{\mathcal{F}_0}^*, \tilde{\mathcal{S}}^3 \right) \xrightarrow{E_3} \cdots \label{cv1}
\end{align}
which also controls foliated deformations of $\left( X, \mathcal{N}_{\mathcal{F}_0}^* \right)$. We will denote the $i$-th cohomoloy group of $\mathcal{E}_{\mathcal{N}_{\mathcal{F}_0}^*}^\bullet$ by $\mathbb{H}^i\left(X, \mathcal{E}_{\mathcal{N}_{\mathcal{F}_0}}^{\bullet}\right)$. We define $E_0': \mathcal{E}_{\mathcal{N}_{\mathcal{F}_0}^*}\to \mathscr{H}om_{\mathcal{O}_X}\left( \mathcal{N}_{\mathcal{F}_0}^*, \Omega_X^1 \right) $ in $(\ref{cv1})$ locally  in the following way:
\begin{align*}
E_0':\mathcal{E}_{\mathcal{N}_{\mathcal{F}_0}^*} |_{U_i} \cong \mathscr{H}om_{\mathcal{O}_X}\left( \mathcal{N}_{\mathcal{F}_0}^*, \mathcal{N}_{\mathcal{F}_0}^* \right) |_{U_i}\bigoplus T_X |_{U_i} \to \mathscr{H}om_{\mathcal{O}_X}\left( \mathcal{N}_{\mathcal{F}_0}^*, \Omega_X^1 \right) |_{U_i}
\end{align*}
On $U_i$, we define $E_0' \left( \left(\phi_i,X \right) \right)$ in terms of the basis $\left\{w_{0i}^\alpha\right\}$
\begin{align}\label{si8}
E_0'\left(\left(\phi_i,X \right)\right)\left(w_{0i}^\alpha \right)= \left(\mathcal{L}_X-\phi_i \right)\left(w_{0i}^\alpha\right)
\end{align}
Then we linearly  extend to $\mathcal{N}_{\mathcal{F}_0}^*|_{U_i}$. We claim that $E_0'$ is well-defined. In fact, if $(\phi_i,X)=(\phi_j,X)$ on $U_{ij}$, then we have
\begin{align}\label{si10}
&\left( \mathcal{L}_X-\phi_i \right)\left(w_{0i}^\alpha\right)-\sum_{\beta=1}^q h_{0ij}^{\alpha\beta}\left(\mathcal{L}_X-\phi_j\right)\left(w_{0j}^\beta\right)= \mathcal{L}_X\left(w_{0i}^\alpha \right)- \phi_i\left(w_{0i}^\alpha\right)-\sum_{\beta=1}^q h_{0ij}^{\alpha\beta}\mathcal{L}_X\left(w_{0j}^\beta\right)+ \sum_{\beta=1}^q h_{0ij}^{\alpha\beta}\phi_j\left(w_{0j}^\beta\right)\\
&=\mathcal{L}_X\left(\sum_{\beta=1}^q h_{0ij}^{\alpha\beta} w_{0j}^\beta \right) -\phi_i\left(w_{0i}^\alpha\right) -\sum_{\beta=1}^q h_{0ij}^{\alpha\beta}\mathcal{L}_X\left(w_{0j}^\beta\right)+\phi_j\left(w_{0i}^\alpha\right)=\sum_{\beta=1}^q \left[X, h_{0ij}^{\alpha\beta} \right]w_{0j}^\beta -\phi_i\left(w_{0i}^\alpha\right)+\phi_j\left(w_{0i}^\alpha\right)=0 \notag
\end{align} 
Next we define $E_1':\mathscr{H}om_{\mathcal{O}_X}\left(\mathcal{N}_{\mathcal{F}_0}^*, \Omega_X^1 \right) \to \mathscr{H}om_{\mathcal{O}_X}\left( \mathcal{N}_{\mathcal{F}_0}^*, \tilde{\mathcal{S}}^2 \right)$ in $(\ref{cv1})$ to be the composition of the natural quotient map $\Omega_X^1\to \Omega_X^1/\mathcal{N}_{\mathcal{F}_0}$ and $E_1$.  Then we have the exact sequence of complex of sheaves
{\Small{\begin{center}
$\begin{CD}
@. \cdots @. \cdots @. \cdots \\
@. @AAA @AE_3AA @AE_3 AA \\
0@>>> 0 @>>> \mathscr{H}om_{\mathcal{O}_X}\left(\mathcal{N}_{\mathcal{F}_0}^*, \tilde{\mathcal{S}}^3 \right) @>>> \mathscr{H}om_{\mathcal{O}_X}\left( \mathcal{N}_{\mathcal{F}_0}^*, \tilde{\mathcal{S}}^3 \right) @>>> 0\\
@. @AAA @AE_2AA @AE_2AA \\
0 @>>> 0@>>> \mathscr{H}om_{\mathcal{O}_X}\left(\mathcal{N}_{\mathcal{F}_0}^*, \tilde{\mathcal{S}}^2 \right) @>>>  \mathscr{H}om_{\mathcal{O}_X}\left( \mathcal{N}_{\mathcal{F}_0}^*, \tilde{\mathcal{S}}^2\right) @>>>0 \\
@. @AAA @AE_1' AA @AE_1AA \\
0@>>> \mathscr{H}om_{\mathcal{O}_X}\left(\mathcal{N}_{\mathcal{F}_0}^*, \mathcal{N}_{\mathcal{F}_0}^* \right) @>>> \mathscr{H}om_{\mathcal{O}_X}\left(\mathcal{N}_{\mathcal{F}_0}^*, \Omega_X^1 \right) @>>> \mathscr{H}om_{\mathcal{O}_X}\left(\mathcal{N}_{\mathcal{F}_0}^*, \frac{\Omega_X^1}{\mathcal{N}_{\mathcal{F}_0}^* }\right) @>>> 0\\
@. @A-\textnormal{Id}AA @AE_0' AA @AE_0AA \\
0@>>> \mathscr{H}om_{\mathcal{O}_X}\left(\mathcal{N}_{\mathcal{F}_0}^*,  \mathcal{N}_{\mathcal{F}_0}^* \right) @>>> \mathcal{E}_{\mathcal{N}_{\mathcal{F}_0}^*}@>>> T_X @>>> 0
\end{CD}$
\end{center}}}
The above exact sequence implies that
\begin{align*}
\mathbb{H}^i\left(X, \mathcal{E}_{\mathcal{N}_{\mathcal{F}_0}^*}^\bullet \right) \cong \mathbb{H}^i\left( X, \mathcal{N}_{\mathcal{F}_0}^{*\bullet}\right)
\end{align*}

We reinterpret Proposition $\ref{nb1}$ by using $\mathcal{E}_{\mathcal{N}_{\mathcal{F}_0}^*}^\bullet$ instead of $\mathcal{N}_{\mathcal{F}_0}^{* \bullet}$.

\begin{proposition}\label{si20}
Let $\left( X, \mathcal{N}_{\mathcal{F}_0}^* \right)$ be a foliated nonsingular variety with $\mathcal{N}_{\mathcal{F}_0}^*$ locally free.
\begin{enumerate}
\item there is a $1-1$ correspondence
\begin{align*}
\textnormal{FDef}_{\left(X, \mathcal{N}_{\mathcal{F}_0}^* \right)}\left(k[\epsilon] \right)\cong \mathbb{H}^1\left(X, \mathcal{E}_{\mathcal{N}_{\mathcal{F}_0}^*}^\bullet \right)
\end{align*}
\item Let $A\in \bold{Art}$ and $\eta=\left( \mathcal{X}, \mathcal{N}_\mathcal{F}^* \right)$ be an infinitesimal foliated deformations of $\left( X, \mathcal{N}_{\mathcal{F}_0}^* \right)$ over $A$. To every small extension $0\to (t)\to \tilde{A} \to A \to 0$, we can associate an element $o_\eta(e)\in \mathbb{H}^2\left( X,  \mathcal{E}_{ \mathcal{N}_{\mathcal{F}_0}^*}^\bullet \right)$ called the obstruction lifting of $\eta$ to $\tilde{A}$, which is $0$ if and only if a lifting of $\eta$ to $\tilde{A}$ exists.
\end{enumerate}
\end{proposition}

\begin{proof}
Let $\left( \mathcal{X}, \mathcal{N}_\mathcal{F}^* \right)$ be a first-order foliated deformation of $\left( X, \mathcal{N}_{\mathcal{F}_0}^* \right)$. We keep the notations in the proof of Proposition \ref{nb1}. We define an element $w_i'\in \Gamma\left( U_i, \mathscr{H}om_{\mathcal{O}_X}\left( \mathcal{N}_{\mathcal{F}_0}^*, \Omega_X^1 \right) \right)$ by $w_i':\Gamma\left( U_i, \mathcal{N}_{\mathcal{F}_0}^* \right) \to \Gamma\left( U_i, \Omega_X^1 \right), w_{0i}^\alpha \mapsto w_i'^\alpha$. We define an element $g_{ij}\in \Gamma\left( U_{ij}, \mathscr{H}om_{\mathcal{O}_X}\left( \mathcal{N}_{\mathcal{F}_0}^*, \mathcal{N}_{\mathcal{F}_0}^* \right) \right)$ by
\begin{align*}
g_{ij}: \Gamma\left( U_{ij}, \mathcal{N}_{\mathcal{F}_0}^* \right) &\to \Gamma\left( U_{ij}, \mathcal{N}_{\mathcal{F}_0}^* \right) \\
 w_{0i}^\alpha & \mapsto \sum_{\beta=1}^q g_{0ij}^{\alpha\beta}  w_{0j}^\beta
\end{align*}
By considering $\left( p_{ij}, - g _{ij} \right)\in T_X|_{U_{ij}} \bigoplus \mathscr{H}om_{\mathcal{O}_X}\left( \mathcal{N}_{\mathcal{F}_0}^*, \mathcal{N}_{\mathcal{F}_0}^* \right)|_{U_{ij}} \subset T_X|_{U_i} \bigoplus \mathscr{H}om_{\mathcal{O}_X}\left( \mathcal{N}_{\mathcal{F}_0}^* , \mathcal{N}_{\mathcal{F}_0}^* \right)|_{U_i} $, we can show $\left\{  \left( p_{ij}, - g_{ij} \right) \right\} \in C^1\left(\mathcal{U}, \mathcal{E}_{\mathcal{N}_{\mathcal{F}_0}^*} \right)$ in a similar way to $(\ref{tn1})$. Then we claim that $\left(\left\{ w_i' \right\}, \left\{ \left( p_{ij}, -g_{ij} \right) \right\}\right)\in C^0\left( \mathcal{U}, \mathscr{H}om_{\mathcal{O}_X}\left( \mathcal{N}_{\mathcal{F}_0}^*, \Omega_X^1 \right) \right) \bigoplus C^1\left( \mathcal{U}, \mathcal{E}_{\mathcal{N}_{\mathcal{F}_0}^*}^\bullet \right)$ defines a $1$-cocycle in the following \v Cech resolution of $\mathcal{E}_{\mathcal{N}_{\mathcal{F}_0}^*}^\bullet$:
{\Small{\begin{center}
$\begin{CD}
\cdots\\
@AE_3AA \\
C^0\left(\mathcal{U}, \mathscr{H}om_{\mathcal{O}_X}\left( \mathcal{N}_{\mathcal{F}_0}^* ,  \tilde{\mathcal{S}}^3 \right) \right)@>\delta>> \cdots\\
@AE_2AA @AE_2AA \\
C^0\left(\mathcal{U}, \mathscr{H}om_{\mathcal{O}_X}\left(\mathcal{N}_{\mathcal{F}_0}^* , \tilde{\mathcal{S}}^2 \right) \right) @>-\delta>> C^1\left(\mathcal{U},  \mathscr{H}om_{\mathcal{O}_X}\left( \mathcal{N}_{\mathcal{F}_0}^* , \tilde{\mathcal{S}}^2  \right) \right)@>\delta>> \cdots\\ 
@AE_1' AA @AE_1' AA @AE_1'AA \\
C^0\left(\mathcal{U}, \mathscr{H}om_{\mathcal{O}_X}\left( \mathcal{N}_{\mathcal{F}_0}^* , \Omega_X^1 \right) \right)@>\delta>> C^1\left(\mathcal{U}, \mathscr{H}om_{\mathcal{O}_X}\left( \mathcal{N}_{\mathcal{F}_0}^*, \Omega_X^1 \right) \right) @>-\delta>> C^2\left(\mathcal{U}, \mathscr{H}om_{\mathcal{O}_X}\left( \mathcal{N}_{\mathcal{F}_0}^* , \Omega_X^1 \right) \right) \\
@AE_0' AA @AE_0' AA @AE_0' AA \\
C^0\left(\mathcal{U}, \mathcal{E}_{\mathcal{N}_{\mathcal{F}_0}^*}\right)@>-\delta>> C^1\left(\mathcal{U}, \mathcal{E}_{\mathcal{N}_{\mathcal{F}_0}^*}\right) @>\delta>> C^2 \left(\mathcal{U}, \mathcal{E}_{\mathcal{N}_{\mathcal{F}_0}^*}\right) @>-\delta>> \cdots
\end{CD}$
\end{center}}}
It is sufficient to show that $\delta\left( \left\{ \left( p_{ij}, -g_{ij}\right)  \right\} \right)=0$ which can be shown in a similar way to $(\ref{tn2})$. We show that for two equivalent first-order foliated deformations of $\left( X, \mathcal{N}_{\mathcal{F}_0}^* \right)$, the associated $1$-cocycles are equivalent. Let $\left( \mathcal{X}', \mathcal{N}_{\mathcal{F}'}^* \right)$ be another foliated deformations of $\left( X, \mathcal{N}_{\mathcal{F}_0}^* \right)$ over $k[\epsilon]$ which induces $\left( \left\{ w_i'' \right\}, \left\{ \left( p_{ij}' , - g_{ij}'\right)\right\}  \right)$ and $\Phi:\left( \mathcal{X}, \mathcal{N}_{\mathcal{F}}^* \right) \to \left( \mathcal{X}', \mathcal{N}_{\mathcal{F}'}^* \right)$ is a foliated isomorphism inducing the identity on $\left(X, \mathcal{N}_{\mathcal{F}_0}^* \right)$. Then we define $g_i\in \Gamma\left( U_i, \mathscr{H}om_{\mathcal{O}_X}\left( \mathcal{N}_{\mathcal{F}_0}^* , \mathcal{N}_{\mathcal{F}_0}^* \right) \right)$ by $g_i:\Gamma\left( U_i, \mathcal{N}_{\mathcal{F}_0}^*\right) \to \Gamma\left(U_i, \mathcal{N}_{\mathcal{F}_0}^* \right), w_{0i}^\alpha \mapsto \sum_{\beta=1}^p  g_i^{\alpha\beta} w_{0i}^\beta$. We consider $\left\{\left\{ \left( a_i, -g_i \right)\right\} \right\}\in C^0\left( \mathcal{U}, \mathcal{E}_{\mathcal{N}_{\mathcal{F}_0}^*}\right)$. It is sufficient to show that $- g_i+ g_j^i= - g_{ij} + g_{ij}'$, i.e. $-g_i\left( w_{0i}^\alpha \right) + g_j\left( w_{0i}^\alpha \right) - \sum_{\beta=1}^p\left[ a_j, h_{0ij}^{\alpha\beta} \right] w_{0j}^\beta = - g_{ij}\left( w_{0i}^\alpha\right) + g_{ij}'\left( w_{0i}^\alpha \right)$. In fact, since $ \theta_{ji} \alpha_j^{-1}= \alpha_i^{-1} \theta_{ji}'$ which can be shown in a similar way to $(\ref{tn3})$ and $(\ref{tn4})$.

Next we identify obstructions. Let $e: 0 \to (t) \to \tilde{A} \to A \to 0$ be a small extension in $\bold{Art}$ and $\eta:=\left( \mathcal{X}, \mathcal{N}_\mathcal{F}^* \right)$ be a foliated deformation of $\left( X, \mathcal{N}_{\mathcal{F}_0}^* \right)$ over $A$.  We keep the notations in the proof of $(2)$ of Proposition \ref{nb1}. We may assume that $\sum_{\beta=1}^q \tilde{f}_{ij}\left( \tilde{h}_{ij}^{\alpha\beta} \right)\tilde{h}_{ji}^{\beta \gamma}= \delta_\alpha^\gamma,\alpha=1,...,q$ in a similar way to $(\ref{ts8})$. We define an element $R_{ij}\in \Gamma\left( U_{ij}, \mathscr{H}om_{\mathcal{O}_X}\left( \mathcal{N}_{\mathcal{F}_0}^*, \Omega_X^1 \right)  \right)$ by $R_{ij}:\Gamma\left( U_{ij}, \mathcal{N}_{\mathcal{F}_0}^* \right)\to \Gamma\left( U_{ij},  \Omega_X^1 \right), w_{0i}^\alpha \mapsto R_{ij}^\alpha$. Then we can show that $\left\{ R_{ij}\right\}\in C^1\left( \mathcal{U}, \mathscr{H}om_{\mathcal{O}_X}\left( \mathcal{N}_{\mathcal{F}_0}^*, \Omega_X^1 \right) \right)$ in a similar way to $(\ref{tn5})$. On the other hand, from $(\ref{tn6})$ we define an element $u_{ijk}\in \Gamma\left(U_{ijk}, \mathscr{H}om_{\mathcal{O}_X}\left(  \mathcal{N}_{\mathcal{F}_0}^*, \mathcal{N}_{\mathcal{F}_0}^* \right) \right)$ by 
\begin{align*}
u_{ijk}:\Gamma\left( U_{ijk}, \mathcal{N}_{\mathcal{F}_0}^* \right) &\to \Gamma\left( U_{ijk}, \mathcal{N}_{\mathcal{F}_0}^* \right) \\
 w_{0i}^\alpha &\mapsto \sum_{\eta=1}^q u_{ijk}^{\alpha \eta} w_{0i}^\eta
\end{align*}
and linearly extends to $\Gamma\left( U_{ijk}, \mathcal{N}_{\mathcal{F}_0}^* \right)$. We consider $\left( \left\{ d_{ijk} \right\}, \left\{ u_{ijk} \right\} \right)\in T_X|_{U_{ijk}} \bigoplus \mathscr{H}om_{\mathcal{O}_X}\left( \mathcal{N}_{\mathcal{F}_0}^*, \mathcal{N}_{\mathcal{F}_0}^* \right)|_{U_{ijk}} \subset T_X|_{U_i}\bigoplus \mathscr{H}om_{\mathcal{O}_X}\left( \mathcal{N}_{\mathcal{F}_0}^*, \mathcal{N}_{\mathcal{F}_0}^* \right)|_{U_i}$ as an element of $\Gamma\left( U_{ijk}, \mathcal{E}_{\mathcal{N}_{\mathcal{F}_0}^*} \right)$. Then we claim that
\begin{align*}
\left( \left\{ B_i\right\}, \left\{- R_{ij} \right\}   ,\left\{ \left( d_{ijk}, u_{ijk}\right) \right\}\right) \in C^0\left( \mathcal{U}, \mathscr{H}om_{\mathcal{O}_X}\left( \mathcal{N}_{\mathcal{F}_0}^*, \tilde{\mathcal{S}}^2 \right) \right) \bigoplus C^1\left(\mathcal{U}, \mathscr{H}om_{\mathcal{O}_X}\left( \mathcal{N}_{\mathcal{F}_0}^*, \Omega_X^1 \right) \right) \bigoplus C^2\left( \mathcal{U}, \mathcal{E}_{\mathcal{N}_{\mathcal{F}_0}^*} \right)
\end{align*}
defines a $2$-cocycle in the above \v Cech resolution of $\mathcal{E}_{\mathcal{N}_{\mathcal{F}_0}^*}^\bullet$. It is sufficient to show that $-\delta\left(\left\{ \left( d_{ijk}, u_{ijk} \right) \right\} \right)=0$ which follows in a similar way to show $-\delta\left( \left\{ \left( d_{ijk}, g_{ijk}\right) \right\} \right)$ in the proof of Proposition \ref{at27}. We show that for another lifting $\tilde{\theta}_{ij}'$ of $\theta_{ij}$, and another lifting $\tilde{w}_i'^\alpha$ of $w_i^\alpha$, and another lifting of $\tilde{h}_{ij}'^{\alpha\beta}$ of $h_{ij}^{\alpha\beta}$, the associated $2$-cocycle $\beta:=\left( \left\{ B_i' \right\}, \left\{ - R_{ij} \right\}, \left\{\left(d_{ijk}', u_{ijk}' \right) \right\} \right)$ is cohomologous to $\alpha:=\left( \left\{ B_i\right\}, \left\{- R_{ij} \right\}, \left\{ \left(d_{ijk}, u_{ijk} \right) \right\} \right)$. We define $y_{ij}\in \Gamma\left( U_{ij}, \mathscr{H}om_{\mathcal{O}_X}\left( \mathcal{N}_{\mathcal{F}_0}^*, \mathcal{N}_{\mathcal{F}_0}^* \right) \right), w_{0i}^\alpha \mapsto \sum_{\beta=1}^q y_{ij}^{\alpha\beta} w_{0j}^\beta$. Then we consider $\left( p_{ij}, - y_{ij}\right)\in T_X|_{U_{ij}}$ $\bigoplus \mathscr{H}om_{\mathcal{O}_X}\left( \mathcal{N}_{\mathcal{F}_0}^*, \mathcal{N}_{\mathcal{F}_0}^* \right)|_{U_{ij}}\subset T_X|_{U_i}\bigoplus \mathscr{H}om_{\mathcal{O}_X}\left( \mathcal{N}_{\mathcal{F}_0}^*, \mathcal{N}_{\mathcal{F}_0}^*\right)|_{U_i}$ as an element of $\Gamma\left( U_{ij}, \mathcal{E}_{\mathcal{N}_{\mathcal{F}_0}^*} \right)$. Then we have $E_1'\left( V_i \right)=B_i'- B_i$ and $\delta\left( \left\{ V_i \right\} \right)+ E_0'\left(\left\{\left( p_{ij}, - y_{ij}\right) \right\} \right)=\left\{\left( - R_{ij}'+ R_{ij}  \right) \right\}$ and, $\delta\left(\left\{\left(  p_{ij}, - y_{ij} \right) \right\} \right)= \left\{ \left(d_{ijk}'- d_{ijk}, u_{ijk}' - u_{ijk} \right) \right\}$ which follows in a similar way to show $\delta\left(\left\{ \left( p_{ij}, - s_{ij} \right) \right\} \right)=\left\{\left(d_{ijk}' - d_{ijk}, g_{ijk}'- g_{ijk} \right) \right\}$ in the proof of Proposition \ref{at27}. This proves that $\alpha$ is cohomologous to $\beta$. So given a small extension $e: 0 \to (t) \to \tilde{A} \to A \to 0$ and an infinitesimal foliated deformation $\eta$ of $\left(X, \mathcal{N}_{\mathcal{F}_0}^* \right)$ over $A$, we can associate an element $o_\eta(e):=$ the cohomology class of $\alpha\in \mathbb{H}^2\left( X , \mathcal{E}_{\mathcal{N}_{\mathcal{F}_0}^*}^\bullet \right)$. We note that $o_\eta(e)=0$ if and only if there exists $\left(\left\{ N_i \right\}, \left\{\left(d_{ij}, - u_{ij} \right) \right\} \right)\in C^0\left( \mathcal{U}, \mathscr{H}om_{\mathcal{O}_X}\left( \mathcal{N}_{\mathcal{F}_0}^*, \Omega_X^1 \right) \right) \bigoplus C^1\left( \mathcal{U},  \mathcal{E}_{\mathcal{N}_{\mathcal{F}_0}^*} \right)$ where $N_i:\Gamma\left( U_i, \mathcal{N}_{\mathcal{F}_0}^* \right) \to \Gamma\left( U_i , \Omega_X^1 \right), w_{0i}^\alpha \mapsto N_i^\alpha$ and $u_{ij}:\Gamma\left( U_{ij}, \mathcal{N}_{\mathcal{F}_0}^* \right) \to \Gamma\left( U_{ij}, \mathcal{N}_{\mathcal{F}_0}^* \right), w_{0i}^\alpha \mapsto \sum_{\beta=1}^q u_{ij}^{\alpha\beta} w_{0j}^\beta$ such that $E_1'\left( N_i \right)=B_i$ and $N_j- N_i+ E_0'\left(\left( d_{ij}, - u_{ij} \right) \right)= - R_{ij}$ and $d_{ij}- d_{ik} + d_{jk}= d_{ijk}$ and $-u_{ij}+ u_{ik}- u_{jk}^i= u_{ijk}$. This implies that a lifting $\tilde{\eta}$ of $\eta$ exists. This completes the proof of Proposition \ref{si20}.

\end{proof}

\subsection{Relationships between foliated deformations of $\left(X, \mathcal{N}_{\mathcal{F}_0}^* \right)$ and the underlying flat deformations of a pair $\left(X, \mathcal{N}_{\mathcal{F}_0} \right)$ of a nonsingular variety $X$ and a locally free shea $\mathcal{N}_{\mathcal{F}_0}^*$} \

Let $\left( X, \mathcal{N}_{\mathcal{F}_0}^* \right)$ be a foliated nonsingular variety with $\mathcal{N}_{\mathcal{F}_0}^*$ locally free. Let $\textnormal{Def}_{\left(X, \mathcal{N}_{\mathcal{F}_0}^* \right)}$ be the functor of Artin rings of flat deformations of the pair of the underlying nonsingular variety $X$ and the underlying locally free sheaf $\mathcal{N}_{\mathcal{F}_0}^*$. As mentioned above, $\textnormal{Def}_{\left( X, T_{\mathcal{F}_0} \right)}$ is controlled by the Atiyah extension $\mathcal{E}_{\mathcal{N}_{\mathcal{F}_0}^*}$. We have a morphism of deformation functors
\begin{align*}
\textnormal{FDef}_{\left( X, \mathcal{N}_{\mathcal{F}_0}^* \right)}\to \textnormal{Def}_{\left( X, \mathcal{N}_{\mathcal{F}_0}^* \right)}
\end{align*}
by forgetting foliated structures and we have a canonical morphism of deformation functors
\begin{align*}
\textnormal{FDef}_{\left( X, \mathcal{N}_{\mathcal{F}_0}^* \right)_{\textnormal{fixed}}} \to \textnormal{FDef}_{\left( X, \mathcal{N}_{\mathcal{F}_0}^* \right)}
\end{align*}
Then the relationship of three deformation functors $\textnormal{FDef}_{\left( X, \mathcal{N}_{\mathcal{F}_0}^* \right)_{\textnormal{fixed}}}, \textnormal{FDef}_{\left( X, \mathcal{N}_{\mathcal{F}_0}^* \right)}$ and $\textnormal{Def}_{\left(X, \mathcal{N}_{\mathcal{F}_0}^* \right)}$ is expressed by the following exact sequence of complexes of sheaves
{\Small{\begin{center}
$\begin{CD}
@. \cdots @. \cdots @. \cdots \\
@. @AE_3AA @AE_3AA @AAA \\
0@>>> \mathscr{H}om_{\mathcal{O}_X}\left(\mathcal{N}_{\mathcal{F}_0}^*,  \tilde{\mathcal{S}}^3 \right) @>>> \mathscr{H}om_{\mathcal{O}_X}\left( \mathcal{N}_{\mathcal{F}_0}^* , \tilde{\mathcal{S}}^3 \right) @>>> 0 @>>> 0 \\
@. @AE_2AA @AE_2AA @AAA \\
0@>>> \mathscr{H}om_{\mathcal{O}_X}\left( \mathcal{N}_{\mathcal{F}_0}^*, \tilde{\mathcal{S}}^2 \right) @>>> \mathscr{H}om_{\mathcal{O}_X}\left( \mathcal{N}_{\mathcal{F}_0}^* ,  \tilde{\mathcal{S}}^2 \right) @>>> 0 @>>> 0 \\
@. @AE_1'AA @AE_1AA @AAA \\
0@>>> \mathscr{H}om_{\mathcal{O}_X}\left( \mathcal{N}_{\mathcal{F}_0}^* , \Omega_X^1 \right) @>>> \mathscr{H}om_{\mathcal{O}_X}\left( \mathcal{N}_{\mathcal{F}_0}^* , \Omega_X^1 \right) @>>> 0 @>>> 0  \\
@. @AAA @AE_0'AA @AAA \\
0 @>>> 0 @>>> \mathcal{E}_{\mathcal{N}_{\mathcal{F}_0}^*} @>>> \mathcal{E}_{\mathcal{N}_{\mathcal{F}_0}^*} @>>> 0
\end{CD}$
\end{center}}}
Then the long exact sequence induces maps of first-order deformations
{\small{\begin{align*}
\textnormal{FDef}_{\left( X, \mathcal{N}_{\mathcal{F}_0}^* \right)_{\textnormal{fixed}}}\left(k[\epsilon] \right)\cong \mathbb{H}^0\left( X, \mathscr{H}om_{\mathcal{O}_X}\left( \mathcal{N}_{\mathcal{F}_0}^* , \Omega_X^1 \right)^\bullet \right) \to \textnormal{FDef}_{\left( X, \mathcal{N}_{\mathcal{F}_0}^*\right)} \left( k[\epsilon] \right) \cong \mathbb{H}^1\left( X, \mathcal{E}_{\mathcal{N}_{\mathcal{F}_0}^*}^\bullet  \right) \to \textnormal{Def}_{\left( X, \mathcal{N}_{\mathcal{F}_0}^* \right)}\left(k[\epsilon] \right)\cong H^1\left( X, \mathcal{E}_{ \mathcal{N}_{\mathcal{F}_0}^*} \right)
\end{align*}}}
and obstruction maps
\begin{align*}
\mathbb{H}^1\left( X, \mathscr{H}om_{\mathcal{O}_X}\left( \mathcal{N}_{\mathcal{F}_0}^* , \Omega_X^1 \right)^\bullet \right) \to \mathbb{H}^2\left( X, \mathcal{E}_{ \mathcal{N}_{\mathcal{F}_0}^* }^\bullet \right) \to H^2\left( X, \mathcal{E}_{ \mathcal{N}_{\mathcal{F}_0}^* } \right)
\end{align*}
of morphisms of deformation functors $\textnormal{FDef}_{\left( X, \mathcal{N}_{\mathcal{F}_0}^* \right)_{\textnormal{fixed}}} \to \textnormal{FDef}_{\left( X, \mathcal{N}_{\mathcal{F}_0}^* \right)} \to \textnormal{Def}_{\left( X , \mathcal{N}_{\mathcal{F}_0}^* \right) }$.

\section{Flat foliated deformations of $\left(X, \mathcal{N}_{\mathcal{F}_0}^* \right)$ with $\mathcal{N}_{\mathcal{F}_0}^*$ locally free in terms of deformations of sections}\label{nt5}

Let $\left( X, \mathcal{N}_{\mathcal{F}_0}^* \right)$ be a foliated nonsingular variety with $\mathcal{N}_{\mathcal{F}_0}^*$ locally free of rank $q$. The exact sequence $0\to \mathcal{N}_{\mathcal{F}_0}^* \to \Omega_X^1 \to \frac{\Omega_X^1}{\mathcal{N}_{\mathcal{F}_0}^*} \to 0$ induces a section $w_0\in H^0\left( X, \bigwedge^q \Omega_X^1 \otimes \mathcal{L}_0 \right)$ where $\mathcal{L}_0:=\left( \bigwedge^q \mathcal{N}_{\mathcal{F}_0}^* \right)^*$. An infinitesimal flat foliated deformation $\left( \mathcal{X}, \mathcal{N}_\mathcal{F}^* \right)$ of $\left( X, \mathcal{N}_{\mathcal{F}_0}^* \right)$ over $A\in \bold{Art}$ defined by  the exact sequence $0 \to \mathcal{N}_\mathcal{F}^* \to \Omega_{\mathcal{X}/A}^1 \to \frac{\Omega_{\mathcal{X}/A}^1}{\mathcal{N}_\mathcal{F}^*} \to 0$ induces a section $w\in H^0\left( \mathcal{X}, \bigwedge^q \Omega_{\mathcal{X}/A}^1 \otimes \mathcal{L} \right)$ where $\mathcal{L}:=\left( \bigwedge^q \mathcal{N}_\mathcal{F}^* \right)^*$ which is an infinitesimal deformation of the section $w_0$ and a pair of the nonsingular variety $X$ and a line bundle $\mathcal{L}_0$. We will describe infinitesimal foliated deformation of $\left( X, \mathcal{N}_{\mathcal{F}_0}^* \right)$ in terms of deformations of the section $w_0$ and the pair of $\left( X, \mathcal{L}_0 \right)$. Let $\mathcal{U}= \left\{U_i \right\}$ be an affine open covering of $X$ such that $\Gamma\left(U_i, \mathcal{N}_{\mathcal{F}_0}^* \right)$ is generated by $w_{0i}^1,..., w_{0i}^q \in \Gamma\left(U_i, \Omega_X^1 \right)$ with $w_{0i}^\alpha = \sum_{\beta=1}^q h_{0ij}^{\alpha\beta} w_{0j}^\beta$ for some $h_{0ij}^{\alpha\beta}\in \Gamma\left( U_{ij}, \mathcal{O}_X \right)$. Let $\det\left( h_{0ij}^{\alpha\beta} \right)= h_{0ij}$. Then we have
\begin{align*}
w_{0i}^1 \wedge \cdots \wedge w_{0i}^q = h_{0ij} w_{0j}^1 \wedge \cdots \wedge w_{0j}^q
\end{align*}
Then $\mathcal{L}_0= \left( \bigwedge^q \mathcal{N}_{\mathcal{F}_0}^* \right)^*$ is defined by the transition functions $\left\{ h_{0ij}\right\}$ with $h_{0ij}\in \Gamma\left( U_{ij} , \mathcal{O}_X^* \right)$ and $w_0$ is defined locally on $U_i$ by $w_0:=\left\{ w_{0i}^1 \wedge \cdots \wedge w_{0i}^q \right\} \in H^0\left( X, \bigwedge^q \Omega_X^1 \otimes \mathcal{L}_0 \right)$. As is known, deformations of a pair of a nonsingular variety $X$ and a line bundle $\mathcal{L}_0$ is controlled by the Atiyah extension (see \cite{Ser06}) which is the exact sequence
\begin{align}\label{si22}
0\to \mathcal{O}_X \to \mathcal{E}_{\mathcal{L}_0} \to T_X \to 0
\end{align}
where $\mathcal{E}_{\mathcal{L}_0}$ is defined locally on $U_i$
\begin{align*}
\mathcal{E}_{\mathcal{L}_0}|_{U_i} := \mathcal{O}_X|_{U_i} \bigoplus T_X|_{U_i}
\end{align*}
and $\left( a_i, b_i \right)\in \Gamma\left( U_i, \mathcal{O}_X \right) \bigoplus \Gamma\left( U_i, T_X \right)$ is identified with $\left(a_j, b_j \right) \in \Gamma\left( U_j , \mathcal{O}_X \right) \bigoplus \Gamma\left( U_j, T_X \right)$ if and only if $b_i=b_j$ and $a_j = a_i + \left[ b_i, \log h_{0ij} \right]$. In the following we will construct a complex of sheaves also controlling foliated deformations of $\left( X, \mathcal{N}_{\mathcal{F}_0}^* \right)$ which explicitly shows the underlying deformations of the section $w_0$ and the line bundle $\mathcal{L}_0$:
\begin{align}\label{si21}
\mathcal{E}_{\mathcal{L}_0}^\bullet : \mathcal{L}_0 \xrightarrow{\tilde{E}_0} \mathcal{S}^1 \xrightarrow{\tilde{E}_1} \mathscr{H}om_{\mathcal{O}_X}\left( \mathcal{N}_{\mathcal{F}_0}^*,  \tilde{\mathcal{S}}^2 \right) \xrightarrow{E_2} \mathscr{H}om_{\mathcal{O}_X}\left( \mathcal{N}_{\mathcal{F}_0}^*, \tilde{\mathcal{S}}^3 \right) \xrightarrow{E_3} \cdots
\end{align}
We will denote the $i$-th cohomology group of $\mathcal{E}_{\mathcal{L}_0}^\bullet$ by $\mathbb{H}^i\left( X, \mathcal{E}_{\mathcal{L}_0}^\bullet \right)$. We define a subsheaf $\mathcal{S}^1\subset \bigwedge^q \Omega_X^1 \otimes \mathcal{L}_0$ which is defined in the following way: $\Gamma\left( U_i, \mathcal{S}^1 \right)\subset \Gamma\left( U_i, \bigwedge^q \Omega_X^1\otimes \mathcal{L}_0 \right)$ is generated by
\begin{align*}
w_{0i}^1\wedge \cdots \wedge  \widehat{w_{0i}^\alpha}^{\alpha-\textnormal{th}}  \wedge \cdots \wedge w_{0i}^q,\,\,\,\,\,\,\,\,\,\,\alpha=1,..., q
\end{align*}
This implies that any element in $\Gamma\left(U_i, \mathcal{S}^1\right)$ is of the form
\begin{align*}
a_i =\sum_{\alpha=1}^q w_{0i}^1\wedge \cdots \wedge \overbrace{a_i^\alpha}^{\alpha-\textnormal{th}} \wedge \cdots \wedge w_{0i}^q,\,\,\,\,\,\,\,\,\, a_i^\alpha  \in \Gamma\left(U_i, \Omega_X^1\right)
\end{align*}
and $a_i$ is transformed into $\Gamma\left(U_j, \mathcal{S}^1 \right)$ in the following way:
\begin{align*}
 a_j:=\frac{1}{h_{0ij}} a_i= \frac{1}{h_{0ij}}\sum_{\alpha=1}^q w_{0i}^1\wedge \cdots \wedge \overbrace{a_i^\alpha}^{\alpha-\textnormal{th}} \wedge \cdots \wedge w_{0i}^q = \sum_{\beta=1}^q w_{0j}^1\wedge \cdots \wedge \overbrace{a_j^\beta}^{\beta-\textnormal{th}} \wedge \cdots \wedge w_{0j}^q,\,\,\,\,\,\,\,\,\, \textnormal{for some}\,\,\, a_j^\beta  \in \Gamma\left(U_i\cap U_j , \Omega_X^1\right),
\end{align*}
so that $\mathcal{S}^1$ is well-defined.

We define a map
\begin{align*}
\alpha: \mathcal{S}^1 \to \mathscr{H}om_{\mathcal{O}_X}\left( \mathcal{N}_{\mathcal{F}_0}^*, \frac{\Omega_X^1}{ \mathcal{N}_{\mathcal{F}_0}^*} \right)
\end{align*}
which is defined locally in the following way: for $a_i\in \Gamma\left(U_i, \mathcal{S}^1 \right)\subset \Gamma\left(U_i, \bigwedge^q \Omega_X^1\right)$ as above,  we define
\begin{align*}
\alpha: \Gamma\left(U_i, \mathcal{S}^1\right)&\to \Gamma\left( U_i, \mathscr{H}om_{\mathcal{O}_X}\left( \mathcal{N}_{\mathcal{F}_0}^*, \Omega_X^1/\mathcal{N}_{\mathcal{F}_0}^* \right)  \right)\\
a_i&\mapsto \left( w_{0i}^\alpha \mapsto \overline{a_i^\alpha}     \right)
\end{align*}
It is well-defined: if we have 
\begin{align*}
a_i=\sum_{\alpha=1}^q w_{0i}^1\wedge \cdots \wedge \overbrace{a_i^\alpha}^{\alpha-\textnormal{th}} \wedge \cdots \wedge w_{0i}^q = \sum_{\alpha=1}^q w_{0i}^1 \wedge \cdots \wedge \overbrace{a_i'^\alpha }^{\alpha-\textnormal{th}} \wedge \cdots \wedge w_{0i}^q
\end{align*}
Then we have $\left(a_i^\alpha-a_i'^\alpha\right)\wedge w_{0i}^1\wedge \cdots \wedge w_{0i}^q=0$. Then for $x\in U_{i_x}\subset U_i-S$, we can write $a_i^\alpha- a_i'^\alpha=\sum_{\alpha=1}^q A_{i_x}^\alpha w_{0i}^\alpha\in \Gamma\left(U_{i_x}, \mathcal{N}_{\mathcal{F}_0}^*\right)$, so that $a_i^\alpha - a_i'^\alpha \in \Gamma\left(U_i-S, \mathcal{N}_{\mathcal{F}_0}^*\right)$. Since $\mathcal{N}_{\mathcal{F}_0}^*$ is reflexive, we have $a_i^\alpha - a_i'^\alpha \in \Gamma\left(U_i, \mathcal{N}_{\mathcal{F}_0}^*\right)$. On the other hand, $a_j$ on $U_j$ is transformed into $a_i=h_{0ij} a_j$ on $U_i$. Then
\begin{align*}
&-\frac{1}{h_{0ij}} w_{0i}^1\wedge \cdots \wedge w_{0i}^q \wedge a_i^\alpha= \frac{1}{h_{0ij}} a_i\wedge w_{0i}^\alpha= a_j \wedge w_{0i}^\alpha\\
&= \sum_{\beta=1}^q w_{0j}^1 \wedge \cdots \wedge \overbrace{a_j^\beta}^{\beta-\textnormal{th}} \wedge \cdots w_{0j}^q \wedge \left( \sum_{\gamma=1}^q h_{0ij}^{\alpha \gamma} w_{0j}^\gamma  \right)  = \sum_{\beta=1}^q w_{0j}^1 \wedge \cdots \wedge \overbrace{a_j^\beta}^{\beta-\textnormal{th}} \wedge \cdots \wedge w_{0j}^q \wedge \left( h_{0ij}^{\alpha\beta} w_{0j}^\beta \right)\\
&= - w_{0j}^1\wedge \cdots \wedge w_{0j}^q \wedge \left(\sum_{\beta=1}^q h_{0ij}^{\alpha\beta} a_j^\beta \right) = -\frac{1}{h_{0ij}} w_{0i}^1\wedge \cdots \wedge w_{0i}^q \wedge \left( \sum_{\beta=1}^q h_{0ij}^{\alpha\beta} a_j^\beta   \right)
\end{align*}
This implies that $a_i^\alpha- \sum_{\beta=1}^q h_{0ij}^{\alpha\beta} a_j^\beta \in \Gamma\left( U_i \cap U_j, \mathcal{N}_{\mathcal{F}_0}^* \right)$. Hence $\mathcal{S}^1 \xrightarrow{\alpha}  \mathscr{H}om_{\mathcal{O}_X}\left( \mathcal{N}_{\mathcal{F}_0}^*, \frac{\Omega_X^1}{ \mathcal{N}_{\mathcal{F}_0}^* } \right)$ is well-defined. Then we define $\tilde{E}_1$ in $(\ref{si21})$ by the composition of $\alpha$ and $E_1:\mathscr{H}om_{\mathcal{O}_X}\left( \mathcal{N}_{\mathcal{F}_0}^*, \frac{\Omega_X^1}{\mathcal{N}_{\mathcal{F}_0}^*} \right)\to \mathscr{H}om_{\mathcal{O}_X}\left(  \mathcal{N}_{\mathcal{F}_0}^*, \tilde{\mathcal{S}}^2 \right)$. 

We have an exact sequence
\begin{align*}
0\to \mathcal{O}_X &\to \mathcal{S}^1 \xrightarrow{\alpha}  \mathscr{H}om_{\mathcal{O}_X}\left( \mathcal{N}_{\mathcal{F}_0}^*, \frac{\Omega_X^1}{\mathcal{N}_{\mathcal{F}_0}^*} \right) \to 0 \\
a &\mapsto aw_0
\end{align*}
where $w_0=\left\{ w_{0i}^1\wedge \cdots \wedge w_{0i}^q \right\}\in H^0\left( X,  \bigwedge^q \Omega_X^1 \otimes \mathcal{L}_0 \right)$. In fact, it is clear that $\alpha$ is surjective. Assume that $\alpha(a_i)=\alpha\left( \sum_{\alpha=1}^q w_{0i}^1 \wedge \cdots \wedge \overbrace{a_i^\alpha}^{\alpha-\textnormal{th}} \wedge \cdots \wedge w_{0i}^q \right)=0$. This implies that $a_i^\alpha \in \Gamma\left(U_i, \mathcal{N}_{\mathcal{F}_0}^*\right)$. Then $a_i= A_i w_{0i}^1\wedge \cdots \wedge w_{0i}^q, A_i\in \Gamma(U_i, \mathcal{O}_X)$.

Next let us define $\tilde{E}_0$ in $(\ref{si21})$. We define $\tilde{E}_0:\mathcal{E}_{\mathcal{L}_0} \to \mathcal{S}^1$ locally on $U_i$ by
\begin{align*}
\mathcal{O}_X|_{U_i} \bigoplus T_X|_{U_i} \to &\mathcal{S}^1|_{U_i}\\
 (a_i, b_i) \mapsto &\mathcal{L}_{b_i}\left(w_{0i}^1\wedge \cdots \wedge w_{0i}^q \right) + a_i w_{0i}^1\wedge \cdots \wedge w_{0i}^q \\
                            &= \sum_{\alpha=1}^q w_{0i}^1\wedge \cdots \wedge \overbrace{ \mathcal{L}_{b_i}\left( w_{0i}^\alpha \right) }^{\alpha-\textnormal{th}} \wedge \cdots \wedge w_{0i}^q + a_i w_{0i}^1\wedge \cdots \wedge w_{0i}^q
\end{align*}
It is clear that the image is in $\Gamma\left(U_i, \mathcal{S}^1\right)$. We show that it is well-defined. In fact,  $\left(a_i, b_i \right)$ on $U_i$ is translated into $\left(a_i+ \left[b_i, \log h_{0ij}\right], b_i \right)$ on $U_j$, and we have
\begin{align*}
&\frac{1}{h_{0ij}}\left(\mathcal{L}_{b_i}\left( w_{0i}^1\wedge \cdots \wedge w_{0i}^q \right) + a_i w_{0i}^1\wedge \cdots \wedge w_{0i}^q \right) = \frac{1}{h_{0ij}}\left( \mathcal{L}_{b_i}\left( h_{0ij} w_{0j}^1\wedge \cdots \wedge w_{0j}^q  \right) \right) + a_i w_{0j}^1\wedge \cdots \wedge w_{0j}^q\\
&=\frac{1}{h_{0ij}}[b_i, h_{0ij}] \left( w_{0j}^1\wedge \cdots \wedge w_{0j}^q \right) + \mathcal{L}_{b_i}\left(w_{0j}^1\wedge \cdots \wedge w_{0j}^q \right) + a_i w_{0j}^1\wedge \cdots \wedge w_{0j}^q\\
& = \mathcal{L}_{b_i}\left( w_{0j}^1\wedge \cdots \wedge w_{0j}^q \right) +\left(a_i+\left[ b_i, \log h_{0ij} \right] \right) w_{0j}^1\wedge \cdots \wedge w_{0j}^q
\end{align*}
Then we have the commutative diagram
\begin{center}
$\begin{CD}
\mathcal{S}^1 @>\alpha>> \mathscr{H}om_{\mathcal{O}_X}\left( \mathcal{N}_{\mathcal{F}_0}, \frac{\Omega_X^1}{\mathcal{N}_{\mathcal{F}_0}^* } \right)\\
@A\tilde{E}_0 AA @AE_0AA \\
 \mathcal{E}_{\mathcal{L}_0} @>>> T_X
\end{CD}$
\end{center}
where $\mathcal{E}_{\mathcal{L}_0}\to T_X$ is the map in the Atiyah extension $(\ref{si22})$. In fact, for $\left(a_i,b_i\right)\in \Gamma\left(U_i, \mathcal{E}_{\mathcal{L}_0}\right)$, we have
\begin{align*}
\alpha\left(\tilde{E}_0\left(a_i,b_i\right) \right) \left(w_{0i}^\alpha\right)=  \alpha\left(\sum_{\alpha=1}^q w_{0i}^1\wedge \cdots \wedge \overbrace{\mathcal{L}_{b_i}(w_{0i}^\alpha)}^{\alpha-\textnormal{th}} \wedge \cdots \wedge w_{0i}^q + a_i w_{0i}^1\wedge \cdots \wedge w_{0i}^q\right) (w_{0i}^\alpha) =\overline{\mathcal{L}_{b_i}(w_{0i}^\alpha)}=E_0(b_i)(w_{0i}^\alpha)
\end{align*}
This shows that $\mathcal{E}_{\mathcal{L}_0}^\bullet$ in $(\ref{si21})$ is a complex. Moreover we have an exact sequence of complexes of sheaves
\begin{center}
$\begin{CD}
@. \cdots @. \cdots @. \cdots \\
@. @AAA @AE_3AA @AE_3 AA \\
0 @>>> 0 @>>> \mathscr{H}om_{\mathcal{O}_X}\left( \mathcal{N}_{\mathcal{F}_0}^*, \tilde{\mathcal{S}}^3 \right) @>>> \mathscr{H}om_{\mathcal{O}_X}\left(\mathcal{N}_{\mathcal{F}_0}^*, \tilde{\mathcal{S}}^3 \right) @>>> 0 \\
@. @AAA @AE_2AA @AE_2AA \\
0@>>> 0 @>>> \mathscr{H}om_{\mathcal{O}_X}\left(\mathcal{N}_{\mathcal{F}_0}^*, \tilde{\mathcal{S}}^2 \right) @>>> \mathscr{H}om_{\mathcal{O}_X}\left( \mathcal{N}_{\mathcal{F}_0}^*, \tilde{\mathcal{S}}^2 \right) @>>> 0\\
@. @AAA @A\tilde{E}_1AA @AE_1AA \\
0@>>> \mathcal{O}_X @>>> \mathcal{S}^1 @>\alpha >> \mathscr{H}om_{\mathcal{O}_X}\left( \mathcal{N}_{\mathcal{F}_0}^*, \frac{\Omega_X^1}{ \mathcal{N}_{\mathcal{F}_0}^* } \right) @>>> 0 \\
@. @A\textnormal{Id} AA @A\tilde{E}_0 AA @AE_0AA \\
0 @>>> \mathcal{O}_X @>>> \mathcal{E}_{\mathcal{L}_0}  @>>> T_X @>>> 0 
\end{CD}$
\end{center}
This implies that we have
\begin{align*}
\mathbb{H}^i\left(X, \mathcal{E}_{\mathcal{L}_0}^\bullet \right) \cong \mathbb{H}^1\left( X, \mathcal{N}_{\mathcal{F}_0}^{*\bullet} \right)
\end{align*}

We reinterpret  Proposition \ref{nb1} in terms of $\mathcal{E}_{\mathcal{L}_0}^\bullet$ instead of $\mathcal{N}_{\mathcal{F}_0}^{*\bullet}$.

\begin{proposition}\label{si39}
Let $\left( X, \mathcal{N}_{\mathcal{F}_0}^* \right)$ be a foliated nonsingular variety with $\mathcal{N}_{\mathcal{F}_0}^*$ locally free.
\begin{enumerate}
\item there is a $1-1$ correspondence
\begin{align*}
\textnormal{FDef}_{\left( X, \mathcal{N}_{\mathcal{F}_0}^* \right)} \left( k[\epsilon] \right) \cong \mathbb{H}^1\left( X, \mathcal{E}_{\mathcal{L}_0}^\bullet \right)
\end{align*}
\item Let $A\in \bold{Art}$ and $\eta=\left( \mathcal{X}, \mathcal{N}_\mathcal{F}^* \right)$ be an infinitesimal foliated deformation of $\left( X, \mathcal{N}_{\mathcal{F}_0}^* \right)$ over $A$. To every small extension $0\to (t) \to \tilde{A} \to A \to 0$, we can associate an element $o_\eta(e)\in \mathbb{H}^2\left( X, \mathcal{E}_{\mathcal{L}_0}^\bullet \right)$ called the obstruction lifting of $\eta$ to $\tilde{A}$, which is $0$ if and only if a lifting of $\eta$ to $\tilde{A}$ exists.
\end{enumerate}

\end{proposition}

\begin{proof}
Let $\left( \mathcal{X}, \mathcal{N}_\mathcal{F}^* \right)$ be a first-order foliated deformation of $\left( X, \mathcal{N}_{\mathcal{F}_0}^* \right)$ over $k[ \epsilon]$. We keep the notations in the proof of $(1)$ of Proposition \ref{nb1}. We define an element $w_i'^*:= \sum_{\alpha=1}^q w_{0i}^1 \wedge \cdots \wedge w_i'^\alpha \wedge \cdots \wedge w_{0i}^q\in \Gamma\left( U_i, \mathcal{S}^1 \right)$. By taking the wedge product of $(\ref{nb5})$, $\alpha=1,...,q$, we have
\begin{align*}
&\theta_{ji}^*\left( \left(w_{0i}^1+ \epsilon w_i'^1 \right) \wedge \cdots \wedge \left( w_{0i}^q + \epsilon w_i'^q \right) \right) \\
&= \left( w_{0i}^1+ \epsilon w_i'^1 -\epsilon \mathcal{L}_{p_{ij}}\left( w_{0i}^1 \right) \right) \wedge \cdots \wedge \left( w_{0i}^q + \epsilon w_i'^q - \epsilon \mathcal{L}_{p_{ij}} \left( w_{0i}^q \right) \right)= \det\left( h_{0ij}^{\alpha\beta} + \epsilon g_{ij}^{\alpha\beta} \right)\left( w_{0j}^1 + \epsilon w_j'^1 \right) \wedge \cdots \wedge \left( w_{0j}^q + \epsilon w_j'^q \right) \\
&\Longrightarrow \sum_{\alpha=1}^q w_{0i}^1 \wedge \cdots \wedge \left( w_i'^\alpha - \mathcal{L}_{p_{ij}}\left( w_{0i}^\alpha \right) \right) \wedge \cdots \wedge w_{0i}^q = h_{0ij}\left( \sum_{\alpha=1}^q w_{0j}^1 \wedge \cdots \wedge w_j'^\alpha \wedge \cdots \wedge w_{0j}^q \right) + G_{ij} w_{0j}^1 \wedge \cdots \wedge w_{0j}^q \\
&\Longrightarrow w_i'^* - \mathcal{L}_{p_{ij}}\left( w_{0i}^1 \wedge \cdots \wedge w_{0i}^q \right) = h_{0ij} w_j'^* + \frac{G_{ij}}{h_{0ij}} w_{0i}^1 \wedge \cdots \wedge w_{0i}^q
\end{align*}
Here we write $\det\left( h_{0ij}^{\alpha\beta} + \epsilon g_{ij}^{\alpha\beta} \right) : = h_{0ij} + \epsilon G_{ij}$ for some $G_{ij}\in \Gamma\left( U_{ij}, \mathcal{O}_X \right)$. Then in a similar way to $(\ref{si31})$, we can show that $\left\{ \left( p_{ij} , \frac{G_{ij}}{h_{0ij}} \right)  \right\}\in C^1\left(\mathcal{U}, \mathcal{E}_{\mathcal{L}_0} \right)$. On the other hand, in a similar way to $(\ref{si32})$, we can show that $\frac{G_{ij}}{h_{0ij}}- \frac{G_{ik}}{ h_{0ik}} + \frac{G_{jk}}{h_{0jk}} -\left[ p_{ij}, \log h_{0ij} \right]=$, so that $\delta\left( \left\{\left( p_{ij}, \frac{G_{ij}}{h_{0ij}} \right) \right\} \right)=0$. This implies that $\left( \left\{ w_i'^* \right\}, \left\{ \left( p_{ij}, \frac{G_{ij}}{h_{0ij}} \right)  \right\} \right)\in C^0\left( \mathcal{U}, \mathcal{S}^1 \right) \bigoplus C^1\left( \mathcal{U}, \mathcal{E}_{\mathcal{L}_0} \right)$ defines a $1$-cocycle in the following \v Cech resolution of $\mathcal{E}_{\mathcal{L}_0}^\bullet$:
{\Small{\begin{center}
$\begin{CD}
\cdots\\
@AE_3AA \\
C^0\left(\mathcal{U}, \mathscr{H}om_{\mathcal{O}_X}\left( \mathcal{N}_{\mathcal{F}_0}^*,  \tilde{\mathcal{S}}^3 \right) \right)@>\delta>> \cdots\\
@AE_2AA @AE_2AA \\
C^0\left(\mathcal{U}, \mathscr{H}om_{\mathcal{O}_X}\left( \mathcal{N}_{\mathcal{F}_0}^* ,  \tilde{\mathcal{S}}^2 \right) \right) @>-\delta>> C^1\left(\mathcal{U}, \mathscr{H}om_{\mathcal{O}_X}\left(  \mathcal{N}_{\mathcal{F}_0}^* , \tilde{\mathcal{S}}^2 \right) \right) @>\delta>> \cdots\\ 
@A\tilde{E}_1AA @A\tilde{E}_1AA @A\tilde{E}_1AA \\
C^0\left(\mathcal{U}, \mathcal{S}^1 \right)@>\delta>> C^1\left(\mathcal{U}, \mathcal{S}^1 \right) @>-\delta>> C^2\left(\mathcal{U},  \mathcal{S}^1 \right) \\
@A \tilde{E}_0AA @A\tilde{E}_0AA @A\tilde{E}_0AA \\
C^0\left(\mathcal{U}, \mathcal{E}_{\mathcal{L}_0}\right)@>-\delta>> C^1\left(\mathcal{U}, \mathcal{E}_{\mathcal{L}_0} \right) @>\delta>> C^2 \left(\mathcal{U}, \mathcal{E}_{\mathcal{L}_0} \right) @>-\delta>> \cdots
\end{CD}$
\end{center}}}

Moreover we can show that for any two equivalent first-order foliated deformations of $\left( X , \mathcal{N}_{\mathcal{F}_0}^* \right)$, the associated $1$-cocycles are equivalent as in the proof of Proposition \ref{ts11}.

Next we identify obstructions. Let $e: 0\to (t) \to \tilde{A} \to A \to 0$ be a small extension in $\bold{Art}$ and $\eta:=\left( \right)$ be a foliated deformation of $\left( X , \mathcal{N}_{\mathcal{F}_0}^* \right)$ over $A$. We keep the notations in the proof of $(2)$ in Proposition \ref{nb1}. We set $h_{ij}:=\det\left( h_{ij}^{\alpha\beta} \right)$ and $\tilde{h}_{ij}:=\det\left( \tilde{h}_{ij}^{\alpha\beta} \right)$. Since $\theta_{ji}^*\left( w_i^1 \wedge \cdots \wedge w_i^q \right)= h_{ij} w_j^1 \wedge \cdots \wedge w_j^q$, we can show that $\left( \tilde{f}_{ij}\tilde{f}_{jk} \right)\left( \tilde{h}_{ik} \right)\cdot \tilde{f}_{ij}\left( \tilde{h}_{kj} \right)\cdot \tilde{h}_{ji}= 1 + t h_{ijk}^*$ for some $h_{ijk}^*\in \Gamma\left( U_{ijk} , \mathcal{O}_X \right)$ in a similar way to $(\ref{si33})$. In a similar way to show  $(\ref{ts8})$, we may assume that $\tilde{f}_{ij}\left( \tilde{h}_{ji} \right) \tilde{h}_{ji}=1$. From $(\ref{si34})$, we have
\begin{align*}
\tilde{\theta}_{ji}^*\left( \tilde{w}_i^1 \wedge \cdots \wedge \tilde{w}_i^q \right) & = \left( \sum_{\beta=1}^q \tilde{h}_{ij}^{1\beta} \tilde{w}_j^\beta + t R_{ij}^1 \right) \wedge \cdots \wedge \left( \sum_{\beta=1}^q \tilde{h}_{ij}^{q\beta} \tilde{w}_j^\beta + t R_{ij}^q \right) \\
&= \tilde{h}_{ij} \tilde{w}_j^1 \wedge \cdots \wedge \tilde{w}_j^q + t \sum_{\alpha=1}^q w_{0i}^1 \wedge \cdots \wedge R_{ij}^\alpha \wedge \cdots \wedge w_{0i}^q
\end{align*}

We define $R_{ij}^*\in \Gamma\left( U_{ij}, \mathcal{S}^1 \right)$ by
\begin{align*}
R_{ij}^* : = \sum_{\alpha=1}^q w_{0i}^1 \wedge \cdots R_{ij}^\alpha \wedge \cdots \wedge w_{0i}^q \in \Gamma\left( U_{ij}, \mathcal{S}^1 \right) \subset \Gamma\left( U_{ij}, \bigwedge^q \Omega_X^1 \otimes \mathcal{L}_0 \right)
\end{align*}
In a similar way to $(\ref{si35})$, we can show that $\left\{ R_{ij}^* \right\}\in C^1\left( \mathcal{U}, \mathcal{S}^1 \right)$. Then we can show that $\left( \left\{ B_i \right\}, \left\{ - R_{ij}^* \right\}, \left\{\left( d_{ijk}, - h_{ijk}^* \right) \right\} \right)$ $\in C^0\left( \mathcal{U}, \mathscr{H}om_{\mathcal{O}_X}\left( \mathcal{N}_{\mathcal{F}_0}^*, \tilde{\mathcal{S}}^2 \right) \right) \bigoplus C^1\left( \mathcal{U}, \mathcal{S}^1 \right) \bigoplus C^2 \left( \mathcal{U}, \mathcal{E}_{\mathcal{L}_0} \right)$ in the above \v Cech resolution of $\mathcal{E}_{\mathcal{L}_0}^\bullet$ in a similar way to show $\left(\left\{ P_i \right\}, \left\{  - W_{ij}^* \right\}, \left\{\left( d_{ijk}, - g_{ijk}^* \right) \right\} \right)$ is a $2$-cocycle as in the proof of Proposition \ref{ts11}. Here we consider $\left( d_{ijk}, - h_{ijk}^* \right)\in T_X|_{U_{ijk}}\bigoplus \mathcal{O}_X|_{U_{ijk}} \subset T_X|_{U_i}\bigoplus \mathcal{O}|_{U_i}$. For another lifting $\tilde{\theta}_{ij}'$ of $\theta_{ij}$, and another lifting $\tilde{w}_i'^\alpha$ of $w_i^\alpha$, and another lifting $\tilde{h}_{ij}'^{\alpha\beta}$ of $h_{ij}^{\alpha\beta}$, the associated $2$-cocycle $\beta:=\left(\left\{ B_i'\right\}, \left\{ - R_{ij}'^*\right\}, \left\{ \left( d_{ijk}' , - h_{ijk}'^* \right) \right\} \right)$ is cohomologous to $\alpha:=\left( \left\{ B_i \right\}, \left\{ - R_{ij}^* \right\}, \left\{ \left( d_{ijk}, - h_{ijk}^* \right) \right\} \right)$. In fact, we define $V_i^*=\sum_{\alpha=1}^q w_{0i}^1 \wedge \cdots \wedge V_i^\alpha \wedge \cdots \wedge w_{0i}^q\in \Gamma\left( U_i, \mathcal{S}^1 \right)$. We can show that $\tilde{h}_{ij}'= \tilde{h}_{ij}+ t h_{0ij} \sum_{\alpha, \beta=1}^q y_{ij}^{\alpha\beta} h_{0ji}^{\beta \alpha}:= \tilde{h}_{ij}+ t h_{0ij} B_{ij}$, where $B_{ij}:=\sum_{\alpha, \beta=1}^q y_{ij}^{\alpha\beta} h_{0ji}^{\beta \alpha}$ as in a similar way to $(\ref{ts9})$ and $\left\{ \left( p_{ij},  B_{ij} \right) \right\} \in C^1\left(\mathcal{U}, \mathcal{E}_{\mathcal{L}_0} \right)$ as in a similar way to $(\ref{si36})$ and $\delta\left( \left\{\left( p_{ij}, B_{ij} \right) \right\} \right)=$ as in a similar way to $(\ref{si37})$. On the other hand, from $(\ref{si38})$, we have
\begin{align*}
- R_{ij}'^* + R_{ij}^* &= \mathcal{L}_{p_{ij}}\left( w_{0i}^1 \wedge \cdots \wedge w_{0i}^q \right) + h_{0ij} V_j^*- V_i^* + \sum_{\alpha=1}^q w_{0i}^1 \wedge \cdots \wedge \left( \sum_{\beta, \gamma=1}^q y_{ij}^{\alpha \beta} h_{0ji}^{\beta \gamma} w_{0i}^\gamma \right) \wedge \cdots \wedge w_{0i}^q \\
& = \mathcal{L}_{p_{ij}}\left( w_{0i}^1 \wedge \cdots \wedge w_{0i}^q \right) + h_{0ij} V_j^* - V_i^* + B_{ij} w_{0i}^1 \wedge \cdots \wedge w_{0i}^q
\end{align*}
This implies that $\delta\left( V_i^* \right) + \tilde{E}_0\left( \left\{\left( p_{ij}, B_{ij} \right) \right\}\right)= \left\{ - R_{ij}'^* + R_{ij}^* \right\}$. This proves that $\alpha$ is cohomologous to $\beta$. So given a small extension $0\to (t) \to \tilde{A} \to A \to 0$ and an infinitesimal foliated deformation $\eta$ of $\left( X, \mathcal{N}_{\mathcal{F}_0}^* \right)$ over $A$, we can associate an element $o_\eta(e):=$ the cohomology class of $\alpha \in \mathbb{H}^2\left( X, \mathcal{E}_{\mathcal{L}_0}^\bullet \right)$.  We note that $o_\eta(e)=0$ if and only if there exists $\left( \left\{ N_i^* \right\}, \left\{ \left(  d_{ij}, U_{ij} \right) \right\} \right)\in  C^0\left(\mathcal{U}, \mathcal{S}^1 \right) \bigoplus C^1\left(  \mathcal{U}, \mathcal{E}_{\mathcal{L}_0} \right)$ where $N_i^*=\sum_{\alpha=1}^q w_{0i}^1 \wedge \cdots \wedge N_i^\alpha \wedge \cdots \wedge w_{0i}^q$ such that $\tilde{E}_1\left( N_i^* \right)= B_i$ and $N_j^*- N_i^*+ \tilde{E}_0\left( \left\{\left( d_{ij}, U_{ij} \right) \right\} \right)= - R_{ij}^*$ and $d_{ij}- d_{ik} + d_{jk}= d_{ijk}$ and $U_{ij}- U_{ik} + U_{jk}^i= - h_{ijk}^*$. Then we have
{\small{\begin{align*}
&t\sum_{\alpha=1}^q w_{0i}^1 \wedge \cdots \wedge \left( \mathcal{L}_{d_{ij}}\left( w_{0i}^\alpha \right) + \sum_{\beta=1}^q h_{0ij}^{\alpha\beta} N_j^\beta - N_i^\alpha \right) \wedge \cdots \wedge w_{0i}^q + t U_{ij} w_{0i}^1 \wedge \cdots \wedge w_{0i}^q= - \tilde{\theta}_{ji}^* \left( \tilde{w}_i^1 \wedge \cdots \wedge \tilde{w}_i^q \right)+ \tilde{h}_{ij} \tilde{w}_j^1 \wedge \cdots \wedge \tilde{w}_j^q \\
&\Longrightarrow \Phi_{ji}^*\left( \left( \tilde{w}_i^1 -t N_i^1 \right)\wedge \cdots \wedge \left( \tilde{w}_i^q - t N_i^q \right) \right)= \left( \tilde{h}_{ij}- t h_{0ij} U_{ij} \right)\left( \tilde{w}_j^1 - t N_j^1\right) \wedge \cdots \left( \tilde{w}_j^q - t N_j^q\right)
\end{align*}}}
and we can show that $\left( \tilde{f}_{ij} - t d_{ij}\right)\left( \tilde{f}_{jk} - t d_{jk}\right)\left(\tilde{h}_{ik}- t h_{0ik} U_{ik} \right) \cdot \left( \tilde{f}_{ij} - t d_{ij} \right)\left(\tilde{h}_{kj} -t h_{0kj} U_{kj} \right)\cdot \left(\tilde{h}_{ji} - t h_{0ji} U_{ji} \right) = 1$ in a similar way to $(\ref{si40})$. Hence a lifting $\tilde{\eta}$ of $\eta$ exists. This completes the proof of Proposition \ref{si39}.
\end{proof}

\subsection{Relationship between foliated deformations of $\left( X, \mathcal{N}_{\mathcal{F}_0}^* \right)$ in terms of vector bundles and in terms of deformations of sections (in terms of line bundles)}\

Let $\left( X, \mathcal{N}_{\mathcal{F}_0}^* \right)$ be a foliated nonsingular variety with $\mathcal{N}_{\mathcal{F}_0}^*$ locally free. We will show that the following diagram commutes and is exact.
\begin{center}
$\begin{CD}
@. \cdots @. \cdots @. \cdots \\
@. @AAA @AE_3AA @AE_3AA \\
0@>>> 0 @>>> \mathscr{H}om_{\mathcal{O}_X}\left( \mathcal{N}_{\mathcal{F}_0}^*, \tilde{\mathcal{S}}^3 \right) @>>> \mathscr{H}om_{\mathcal{O}_X}\left( \mathcal{N}_{\mathcal{F}_0}^*, \tilde{\mathcal{S}}^3 \right) @>>> 0 \\
@. @AAA @AE_2AA @AE_2AA \\
0@>>> 0 @>>> \mathscr{H}om_{\mathcal{O}_X}\left( \mathcal{N}_{\mathcal{F}_0}^*, \tilde{\mathcal{S}}^2\right) @>>> \mathscr{H}om_{\mathcal{O}_X}\left( \mathcal{N}_{\mathcal{F}_0}^*, \tilde{\mathcal{S}}^2 \right) @>>> 0 \\
@. @AAA @AE_1'AA @A\tilde{E}_1AA \\
0@>>> \mathscr{H}om_{\mathcal{O}_X}\left(\mathcal{N}_{\mathcal{F}_0}^*, \mathcal{N}_{\mathcal{F}_0}^* \right)|_{\textnormal{Tr}=0} @>>> \mathscr{H}om_{\mathcal{O}_X}\left(\mathcal{N}_{\mathcal{F}_0}^*, \Omega_X^1 \right) @>\alpha_1>> \mathcal{S}^1\subset \bigwedge^q \Omega_X^1\otimes \mathcal{L}_0 @>>> 0 \\
@. @A-\textnormal{Id}AA @AE_0'AA @A\tilde{E}_0AA \\
0 @>>> \mathscr{H}om_{\mathcal{O}_X}\left(\mathcal{N}_{\mathcal{F}_0}^*, \mathcal{N}_{\mathcal{F}_0}^*\right)|_{\textnormal{Tr}=0} @>>> \mathcal{E}_{\mathcal{N}_{\mathcal{F}_0}^*} @>\alpha_0 >> \mathcal{E}_{\mathcal{L}_0} @>>> 0
\end{CD}$
\end{center}
We define tha map $\alpha_1$ in the following way:
\begin{align*}
\alpha_1:\mathscr{H}om_{\mathcal{O}_X}\left(\mathcal{N}_{\mathcal{F}_0}^*, \Omega_X^1 \right) &\to \mathcal{S}^1 \subset \bigwedge^q \Omega_X^1\otimes \mathcal{L}_0^* \\
                       \phi& \mapsto \phi_i:=\sum_{\alpha=1}^q w_{0i}^1\wedge \cdots \wedge \phi(w_{0i}^\alpha)\wedge \cdots \wedge w_{0i}^q\,\,\,\,\,\,\,\,\,\textnormal{on}\,\,\, U_i
\end{align*}
This is well-defined since
\begin{align*}
\phi_i=\sum_{\alpha=1}^q w_{0i}^1\wedge \cdots \wedge \phi(w_{0i}^\alpha) \wedge \cdots \wedge w_{0i}^q = h_{0ij} \sum_{\alpha=1}^q w_{0j}^1\wedge \cdots \wedge \phi(w_{0j}^\alpha)\wedge \cdots \wedge w_{0j}^q = h_{0ij}\phi_j
\end{align*}
This is clearly surjective. Let us find the kernel. Assume that $\alpha_1(\phi)=0$. Then 
\begin{align*}
\sum_{\alpha=1}^q w_{0i}^1\wedge \cdots \wedge \phi(w_{0i}^\alpha) \wedge \cdots \wedge w_{0i}^q = 0
\end{align*}
By multiplying $w_{0i}^\alpha$, we have
\begin{align*}
w_{0i}^1\wedge \cdots \wedge w_{0i}^q \wedge \phi(w_{0i}^\alpha)=0
\end{align*}
Then $\phi(w_{0i}^\alpha)\in \Gamma\left(U_i-S, \mathcal{N}_{\mathcal{F}_0}^*\right)$, so that by the reflexiveness of $\mathcal{N}_{\mathcal{F}_0}^*$, we have $\phi(w_{0i}^\alpha)\in \Gamma(U_i, \mathcal{N}_{\mathcal{F}_0}^*)$ and so $\phi(w_{0i}^\alpha)= \sum_{\beta=1}^q f_i^{\alpha\beta} w_{0i}^\beta$ for some $f_i^{\alpha\beta}\in \Gamma(U_i, \mathcal{O}_X)$. Then
\begin{align*}
0=\sum_{\alpha=1}^p w_{0i}^1\wedge \cdots \wedge \phi(w_{0i}^\alpha) \wedge \cdots \wedge w_{0i}^q=\sum_{\alpha=1}^p w_{0i}^1\wedge \cdots \wedge f_i^{\alpha\alpha} w_{0i}^\alpha \wedge \cdots \wedge  w_{0i}^q= \left( f_i^{11} + \cdots + f_i^{qq} \right) w_{0i}^1\wedge \cdots \wedge w_{0i}^q
\end{align*}
so that $\sum_{\alpha=1}^p f_i^{\alpha\alpha}=0$. This implies that $\phi\in \mathscr{H}om_{\mathcal{O}_X}\left( \mathcal{N}_{\mathcal{F}_0}^*, \mathcal{N}_{\mathcal{F}_0}^* \right)|_{\textnormal{Tr}=0}$.

We define $\alpha_0 :\mathcal{E}_{\mathcal{N}_{\mathcal{F}_0}^*}\to \mathcal{E}_{\mathcal{L}_0^*}$ locally in the following way:
\begin{align*}
\alpha_0: \mathcal{E}_{\mathcal{N}_{\mathcal{F}_0}^*}|_{U_i} \cong  \mathscr{H}om_{\mathcal{O}_X}\left( \mathcal{N}_{\mathcal{F}_0}^*, \mathcal{N}_{\mathcal{F}_0}^* \right) |_{U_i} \bigoplus T_X|_{U_i}  &\to \mathcal{E}_{\mathcal{L}_0}|_{U_i} \cong  \mathcal{O}_X|_{U_i} \bigoplus T_X|_{U_i} \\
(\phi_i, X) &\mapsto ( -\textnormal{Tr}(\phi_i), X)
\end{align*}
Then the kernel of $\alpha_0$ on $U_i$ is $\mathscr{H}om_{\mathcal{O}_X}\left(\mathcal{N}_{\mathcal{F}_0}^* , \mathcal{N}_{\mathcal{F}_0}^* \right)|_{\textnormal{Tr}=0}|_{U_i}$. Next we show that  $\alpha_0$ is well-defined, i.e. the following diagram commutes
\begin{center}
$\begin{CD}
(\phi_j, X) \in \mathscr{H}om_{\mathcal{O}_X}\left( \mathcal{N}_{\mathcal{F}_0}^*, \mathcal{N}_{\mathcal{F}_0}^*  \right) |_{U_i} \bigoplus T_X |_{U_i} @>>> (- \textnormal{Tr}(\phi_j), X)\in \mathcal{O}_X|_{U_i} \bigoplus T_X |_{U_i} \\
@AAA @AAA \\
(\phi_i, X)\in \mathscr{H}om_{\mathcal{O}_X}\left(\mathcal{N}_{\mathcal{F}_0}^*, \mathcal{N}_{\mathcal{F}_0}^*  \right) |_{U_i} \bigoplus T_X |_{U_i} @>>>  (- \textnormal{Tr}(\phi_i), X) \in  \mathcal{O}_X |_{U_i} \bigoplus  T_X|_{U_i}
\end{CD}$
\end{center}
where $\phi_i(w_{0i}^\alpha)= \phi_j( w_{0i}^\alpha) + \sum_{\beta=1}^q \left[ X, h_{0ij}^{\alpha\beta} \right] w_{0j}^\beta=\phi_j(w_{0i}^\alpha)+ \mathcal{L}_X\left(w_{0i}^\alpha\right) - \sum_{\beta=1}^q h_{0ij}^{\alpha\beta} \mathcal{L}_X\left(w_{0j}^\beta\right)$. Then
{\small{\begin{align*}
&\textnormal{Tr}(\phi_i) w_{0i}^1\wedge \cdots \wedge w_{0i}^q =\sum_{\alpha=1}^q w_{0i}^1\wedge \cdots \wedge \phi_i(w_{0i}^\alpha) \wedge \cdots \wedge w_{0i}^q =\sum_{\alpha=1}^q w_{0i}^1\wedge \cdots \wedge \left( \phi_j(w_{0i}^\alpha) + \mathcal{L}_X\left(w_{0i}^\alpha \right) - \sum_{\beta=1}^q h_{0ij}^{\alpha\beta} \mathcal{L}_X\left( w_{0j}^\beta \right)   \right) \wedge \cdots \wedge w_{0i}^q\\
&= h_{0ij} \textnormal{Tr}(\phi_j)w_{0j}^1\wedge \cdots \wedge w_{0j}^q +\mathcal{L}_X\left( w_{0i}^1\wedge \cdots \wedge w_{0i}^q \right)  - h_{0ij} \mathcal{L}_X\left( w_{0j}^1 \wedge \cdots \wedge w_{0j}^q  \right) =  \textnormal{Tr}(\phi_j) w_{0i}^1\wedge \cdots \wedge w_{0i}^q +\left[ X, h_{0ij}\right] w_{0j}^1\wedge \cdots \wedge  w_{0j}^q
\end{align*}}}
Then we have $\textnormal{Tr}(\phi_i)= \textnormal{Tr}(\phi_j)+  \left[ X,  \log h_{0ij} \right]$. Hence the diagram commutes.

Next we show that the following diagram commutes
\begin{center}
$\begin{CD}
\mathscr{H}om_{\mathcal{O}_X}\left(\mathcal{N}_{\mathcal{F}_0}^*, \Omega_X^1 \right) @>\alpha_1>> \mathcal{S}^1 \\
@AE_0'AA @AA\tilde{E}_0A \\
\mathcal{E}_{\mathcal{N}_{\mathcal{F}_0}^*} @>\alpha_0 >> \mathcal{E}_{\mathcal{L}_0}
\end{CD}$ 
\end{center}
In fact, on $U_i$, $E_0'(\phi_i, X)(w_{0i}^\alpha)= \left(\mathcal{L}_X-\phi_i\right)(w_{0i}^\alpha)=\mathcal{L}_X(w_{0i}^\alpha) -\phi_i(w_{0i}^\alpha)$ for $(\phi_i, X)\in \Gamma(U_i, \mathcal{E}_{\mathcal{N}_{\mathcal{F}_0}^*})$. Then we have
\begin{align*}
\alpha_1E_0'(\phi_i, X)&= \sum_{\alpha=1}^q w_{0i}^1\wedge \cdots \wedge \left( \mathcal{L}_X(w_{0i}^\alpha)- \phi_i(w_{0i}^\alpha)    \right) \wedge \cdots \wedge w_{0i}^q = \mathcal{L}_X \left( w_{0i}^1\wedge \cdots \wedge w_{0i}^q  \right) - \textnormal{Tr}(\phi_i) w_{0i}^1\wedge \cdots \wedge w_{0i}^q \\
 &=\tilde{E}_0(-\textnormal{Tr}(\phi_i), X)=\tilde{E}_0\alpha_0(\phi_i, X)
\end{align*}
and by definition $E_1'=\tilde{E}_1\circ \alpha_1$, and for $\phi\in \mathscr{H}om_{\mathcal{O}_X}\left(\mathcal{N}_{\mathcal{F}_0}^*, \mathcal{N}_{\mathcal{F}_0}^* \right)|_{\textnormal{Tr}=0}$, we have $E_0'(\phi, 0)=-\phi$. Hence the above diagram is exact. This implies that we have
\begin{align*}
\mathbb{H}^i\left( X, \mathcal{E}_{\mathcal{N}_{\mathcal{F}_0}^*}^\bullet   \right) \cong \mathbb{H}^i\left( X , \mathcal{E}_{\mathcal{L}_0}^\bullet \right)
\end{align*}
The map $\mathcal{E}_{\mathcal{N}_{\mathcal{F}_0}^*}^\bullet \to \mathcal{E}_{\mathcal{L}_0}^\bullet$ above gives an isomorphism of first-order deformations and obstructions of foliated deformations of $\left(X, \mathcal{N}_{\mathcal{F}_0}^* \right)$ in terms of cotangent sheaves between interpretations in terms of the vector bundle and in terms of line bundle. In fact, we keep the notations in the proof of Proposition \ref{si20} and Proposition \ref{si39}. Then we have to show that
\begin{align*}
\mathbb{H}^1\left( X, \mathcal{E}_{\mathcal{N}_{\mathcal{F}_0}^*}^\bullet \right)&\xrightarrow{\alpha} \mathbb{H}^1\left( X, \mathcal{E}_{\mathcal{L}_0}^\bullet  \right) \\
\left(\left\{ w_i' \right\}, \left\{ p_{ij}, - g_{ij} \right\} \right) &\mapsto \left( \left\{ w_i'^*\right\}, \left\{\left( p_{ij}, \frac{G_{ij}}{h_{0ij}} \right) \right\} \right)
\end{align*}
In fact, it is clear that $w_i'\mapsto w_i'^*$ under $\alpha_1$, and we have $-\textnormal{Tr}\left( - g_{ij} \right) = \frac{G_{ij}}{h_{0ij}}$ as in a similar way with $(\ref{si43})$. 

Next we show that $\alpha:\mathcal{E}_{\mathcal{N}_{\mathcal{F}_0}^*}^\bullet \to \mathcal{E}_{\mathcal{L}_0}^\bullet$ gives the obstruction map for a small extension $e:0\to (t) \to \tilde{A} \to A \to 0$. We have to show that
\begin{align*}
\mathbb{H}^2\left(X, \mathcal{E}_{\mathcal{N}_{\mathcal{F}_0}^*}^\bullet \right) &\xrightarrow{\alpha} \mathbb{H}^2\left( X, \mathcal{E}_{\mathcal{L}_0}^\bullet \right) \\
\left(  \left\{ B_i \right\}, \left\{- R_{ij} \right\} , \left\{\left(d_{ijk}, - h_{ijk} \right) \right\} \right) &\mapsto \left(  \left\{ B_i \right\}, \left\{ - R_{ij}^* \right\}, \left\{ \left( d_{ijk}, - h_{ijk}^* \right) \right\}\right)
\end{align*} 
In fact, it is clear that $- W_{ij}\mapsto - W_{ij}^*$ under $\alpha_1$ and $- \textnormal{Tr}(h_{ijk})= - h_{ijk}^*$ in a similar way with $(\ref{si45})$.

We also remark  the case of foliated deformations of $\left(X, \mathcal{N}_{\mathcal{F}_0}^* \right)$ with $X$ fixed (i.e. without deformations of complex structures of $X$). We note that we have shown that $\textnormal{FDef}_{\mathcal{N}_{\mathcal{F}_0}^*}$ is controlled by $\mathscr{H}om_{\mathcal{O}_X}\left( \mathcal{N}_{\mathcal{F}_0}^*, \mathcal{N}_{\mathcal{F}_0}^* \right)^\bullet$ from $(\ref{b9})$. We consider foliated deformations of $T_{\mathcal{F}_0}$ with $X$ fixed in terms of the line bundle $\mathcal{L}_0$. We have the following exact sequence of complexes of sheaves

{\small{\begin{center}
$\begin{CD}
@. \cdots @. \cdots @. \cdots \\
@. @AAA @AAA @AAA \\
0 @>>> 0 @>>> \mathscr{H}om_{\mathcal{O}_X } \left( \mathcal{N}_{\mathcal{F}_0}^*, \tilde{\mathcal{S}}^3 \right) @>>> \mathscr{H}om_{\mathcal{O}_X}\left( \mathcal{N}_{\mathcal{F}_0}^*, \tilde{\mathcal{S}}^3 \right) @>>> 0 \\
@. @AAA @AAA @AAA \\
0 @>>> 0 @>>> \mathscr{H}om_{\mathcal{O}_X}\left( \mathcal{N}_{\mathcal{F}_0}^*, \tilde{\mathcal{S}}^2  \right) @>>>  \mathscr{H}om_{\mathcal{O}_X}\left( \mathcal{N}_{\mathcal{F}_0}^*, \tilde{\mathcal{S}}^2 \right) @>>> 0 \\
@. @AAA @AAA @AAA \\
0@>>> \mathscr{H}om_{\mathcal{O}_X}\left( \mathcal{N}_{\mathcal{F}_0}^*, \mathcal{N}_{\mathcal{F}_0}^* \right)|_{\textnormal{Tr}=0} @>>> \mathscr{H}om_{\mathcal{O}_X}\left( \mathcal{N}_{\mathcal{F}_0}^*, \Omega_X^1   \right) @>>> \mathcal{S}^1 @>>> 0 \\
@. @A-\textnormal{Id}AA @AAA @AAA \\
0 @>>> \mathscr{H}om_{\mathcal{O}_X}\left( \mathcal{N}_{\mathcal{F}_0}^*, \mathcal{N}_{\mathcal{F}_0}^* \right)|_{\textnormal{Tr}=0} @>>> \mathscr{H}om_{\mathcal{O}_X}\left( \mathcal{N}_{\mathcal{F}_0}^*, \mathcal{N}_{\mathcal{F}_0}^* \right) @>-\textnormal{Tr}>> \mathcal{O}_X @>>> 0
\end{CD}$
\end{center}}}
We denote the complex of sheave $\mathcal{O}_X \to \mathcal{S}^1\to \mathscr{H}om_{\mathcal{O}_X}\left(\mathcal{N}_{\mathcal{F}_0}^*, \tilde{\mathcal{S}}^2 \right) \to \mathscr{H}om_{\mathcal{O}_X}\left( \mathcal{N}_{\mathcal{F}_0}^*, \tilde{\mathcal{S}}^3 \right) \to \cdots $ by $\mathcal{L}_0^{\bullet}$. Then we have
\begin{align*}
\mathbb{H}^i\left( X, \mathscr{H}om_{\mathcal{O}_X}\left( \mathcal{N}_{\mathcal{F}_0}^*, \mathcal{N}_{\mathcal{F}_0}^* \right)^\bullet \right) \cong \mathbb{H}^i\left( X, \mathcal{L}_0^{\bullet} \right)
\end{align*}

\section{Foliated deformations of $\left( X, \mathcal{N}_{\mathcal{F}_0}^* \right)$ with $\mathcal{N}_{\mathcal{F}_0}^*$ locally free by using the complex $\mathcal{N}_{\mathcal{F}_0}'^{*\bullet}$ }\label{nt9}

We note that we have an isomorphism of complexes of sheaves (see Part I):
{\Tiny{\begin{equation}\label{si41}
\begin{CD}
\mathcal{N}_{\mathcal{F}_0}^{*\bullet}:@. T_X @>E_0>> \mathscr{H}om_{\mathcal{O}_X}\left( \mathcal{N}_{\mathcal{F}_0}^*, \frac{\Omega_X^1}{\mathcal{N}_{\mathcal{F}_0}^*} \right) @>E_1>> \mathscr{H}om_{\mathcal{O}_X}\left( \mathcal{N}_{\mathcal{F}_0}^*, \tilde{\mathcal{S}}^2 \right)@>E_2>> \mathscr{H}om_{\mathcal{O}_X}\left( \mathcal{N}_{\mathcal{F}_0}^*, \tilde{\mathcal{S}}^3 \right) @>E_3>> \cdots \\
@. @| @| @V\alpha_2 VV @V \alpha_3 VV \\
\mathcal{N}_{\mathcal{F}_0}'^{*\bullet} :@. T_X @>E_0>> \mathscr{H}om_{\mathcal{O}_X}\left( \mathcal{N}_{\mathcal{F}_0}^*, \frac{\Omega_X^1}{\mathcal{N}_{\mathcal{F}_0}^*} \right) @>E_1'>> \mathscr{H}om_{\mathcal{O}_X}\left( \mathcal{N}_{\mathcal{F}_0}^*, \mathscr{H}om_{\mathcal{O}_X}\left(\bigwedge^2 T_{\mathcal{F}_0} , \mathcal{O}_X \right) \right)@>E_2'>>  \mathscr{H}om_{\mathcal{O}_X}\left( \mathcal{N}_{\mathcal{F}_0}^*, \mathscr{H}om_{\mathcal{O}_X}\left(\bigwedge^3 T_{\mathcal{F}_0} , \mathcal{O}_X \right) \right) @>E_3'>> \cdots
\end{CD}
\end{equation}}}
Hence we have
\begin{align*}
\mathbb{H}^i\left(X, \mathcal{N}_{\mathcal{F}_0}^{*\bullet} \right) \cong \mathbb{H}^i\left( X, \mathcal{N}_{\mathcal{F}_0}'^{*\bullet}\right)
\end{align*}

We present foliated deformations of $\left(X, \mathcal{N}_{\mathcal{F}_0}^* \right)$ by using $\mathcal{N}_{\mathcal{F}_0}'^{*\bullet}$, and we will describe the correspondences of first-oder deformations and obstructions from isomorphisms $\mathbb{H}^1\left( X, \mathcal{N}_{\mathcal{F}_0}^{*\bullet} \right) \to \mathbb{H}^1\left(X, \mathcal{N}_{\mathcal{F}_0}'^{*\bullet} \right)$ and $\mathbb{H}^2\left( X, \mathcal{N}_{\mathcal{F}_0}^{*\bullet} \right) \to \mathbb{H}^2\left(X, \mathcal{N}_{\mathcal{F}_0}'^{*\bullet} \right)$.

\begin{definition}
To emphasize the integrability condition for using the complex of sheaves $\mathcal{N}_{\mathcal{F}_0}'^{*\bullet}$, we define foliated deformations of $\left(X, \mathcal{N}_{\mathcal{F}_0}^* \right)$ with locally free equivalently in the following way$:$ let $\left(X, \mathcal{N}_{\mathcal{F}_0}^* \right)$ be a foliated nonsingular variety with $\mathcal{N}_{\mathcal{F}_0}^*$ locally free. An infinitesimal foliated deformation $\eta$ of $\left( X, \mathcal{N}_{\mathcal{F}_0}^* \right)$ over $A \in \bold{Art}$ is a cartesian diagram of schemes
\begin{center}
$\eta:\begin{CD}
\left( X, \mathcal{N}_{\mathcal{F}_0}^* \right) @>>> \left( \mathcal{X}, \mathcal{N}_\mathcal{F}^* \right)\\
@VVV @VV \pi V\\
\textnormal{Spec}(k) @>>> \textnormal{Spec}(A)
\end{CD}$
\end{center}
where $\pi$ is flat, $\mathcal{N}_{\mathcal{F}}^*$ is a locally free subsheaf of $\Omega_{\mathcal{X}/A}^1$ and $\frac{\Omega_{\mathcal{X}/A}^1}{\mathcal{N}_\mathcal{F}^*}$ is flat over $A$ such that $d_{\mathcal{X}/A}\left( \mathcal{N}_\mathcal{F}^* \right)$ vanishes on $\bigwedge^2 T_{\mathcal{F}}$ where $T_{\mathcal{F}}$ is defined by the dual of $\frac{\Omega_{\mathcal{X}/A}^1}{\mathcal{N}_\mathcal{F}^*}$ and $d_{\mathcal{X}/A}$ is the relative differential over $A$, and $\left( \mathcal{X}, \mathcal{N}_\mathcal{F}^* \right)$ induces $\left(X, \mathcal{N}_{\mathcal{F}_0}^* \right)$. Two infinitesimal foliated deformations $\left( \mathcal{X}, \mathcal{N}_\mathcal{F}^* \right)$ and $\left(\mathcal{X}', \mathcal{N}_{\mathcal{F}'}^* \right)$ of $\left( X, \mathcal{N}_{\mathcal{F}_0}^* \right)$ over $A$ are isomorphic if
\end{definition}

\begin{proposition}\label{si30}
Let $\left( X, \mathcal{N}_{\mathcal{F}_0}^* \right)$ be a foliated nonsingular variety with $\mathcal{N}_{\mathcal{F}_0}^*$ locally free.
\begin{enumerate}
\item there is a $1-1$ correspondence
\begin{align*}
\textnormal{FDef}_{\left(X, \mathcal{N}_{\mathcal{F}_0}^* \right)}\left(k[\epsilon] \right)\cong \mathbb{H}^1\left(X, \mathcal{N}_{\mathcal{F}_0}'^{*\bullet} \right)
\end{align*}
\item Let $A\in \bold{Art}$ and $\eta=\left( \mathcal{X}, \mathcal{N}_\mathcal{F}^* \right)$ be an infinitesimal foliated deformations of $\left( X, \mathcal{N}_{\mathcal{F}_0}^* \right)$ over $A$. To every small extension $0\to (t)\to \tilde{A} \to A \to 0$, we can associate an element $o_\eta(e)\in \mathbb{H}^2\left( X,  \mathcal{N}_{\mathcal{F}_0}^{ * \bullet} \right)$ called the obstruction lifting of $\eta$ to $\tilde{A}$, which is $0$ if and only if a lifting of $\eta$ to $\tilde{A}$ exists.
\end{enumerate}
\end{proposition}

\begin{proof}
We note that descriptions of first-oder deformations by using $\mathcal{N}_{\mathcal{F}_0}'^{*\bullet}$ is exactly same to the descriptions by using $\mathcal{N}_{\mathcal{F}_0}^{*\bullet}$ in Proposition \ref{nb1} since the first two terms of $\mathcal{N}_{\mathcal{F}_0}^{*\bullet}$ and $\mathcal{N}_{\mathcal{F}_0}'^{*\bullet}$ are identical.

We identify obstructions. Let $0\to (t)\to \tilde{A} \to A \to 0$ be a small extension in $\bold{Art}$. Let $\left( \mathcal{X}, \mathcal{N}_{\mathcal{F}}^* \right)$ be an infinitesimal foliated deformation of $\left(X, \mathcal{N}_{\mathcal{F}_0}^* \right)$ over $\textnormal{Spec}(A)$, so that we have the exact sequence $0 \to \mathcal{N}_{\mathcal{F}}^*\to \Omega_{\mathcal{X}/\textnormal{Spec}(A)}^1 \to \frac{\Omega_{\mathcal{X}/\textnormal{Spec}(A)}^1}{\mathcal{N}_{\mathcal{F}}^*}\to 0$. We set $T_{\mathcal{F}}:=\mathscr{H}om_{\mathcal{O}_{\mathcal{X}}}\left( \frac{\Omega_{\mathcal{X}/\textnormal{Spec}(A)}^1}{\mathcal{N}_{\mathcal{F}}^*},  \mathcal{O}_{\mathcal{X}} \right) \subset T_{\mathcal{X}/\textnormal{Spec}(A)}$. Let $\mathcal{U}= \left\{U_i \right\}$ be an affine open covering of $X$ such that we have the foliated isomorphisms $\theta_i:\left( U_i\times \textnormal{Spec}(A), \mathcal{N}_{\mathcal{F}_i}^* \right)\to \left(\mathcal{X}|_{U_i}, \mathcal{N}_{\mathcal{F}}^*|_{U_i} \right)$, where $\mathcal{N}_{\mathcal{F}_i}^*\subset \Omega_X^1|_{U_i}\otimes A$ is the induced foliation on $U_i\times \textnormal{Spec}(A)$ from $\mathcal{N}_{\mathcal{F}}^*$ via $\theta_i$. Let $\Gamma\left( U_i, \mathcal{N}_{\mathcal{F}_i}^* \right) \subset \Gamma\left( U_i, \Omega_X^1 \right)\otimes A$ be generated by $w_i^\alpha,\alpha=1,...,q$ which induces $w_{0i}^\alpha$ modulo $\mathfrak{m}_A$ as in the proof of Proposition \ref{nb1}. Then the integrability condition can be written as $dw_i^\alpha\left(T_1, T_2     \right)=i_{T_1\wedge T_2}\left(dw_i^\alpha \right)=0 $ for $T_1,T_2\in \Gamma\left( U_i, T_{\mathcal{F}_i}\right)$ and for $\alpha=1,...,q$, where $T_{\mathcal{F}_i}:=\mathscr{H}om_{\mathcal{O}_X|_{U_i}\otimes A}\left( \frac{\Omega_X^1|_{U_i}\otimes A}{\mathcal{N}_{\mathcal{F}_i}^*}, \mathcal{O}_X|_{U_i}\otimes A \right)$. Let $\theta_{ij}:= \theta_j^{-1}\theta_i:\left(  U_{ij} \times \textnormal{Spec}(A), \mathcal{N}_{\mathcal{F}_i}^* \right) \to \left( U_{ij}\times \textnormal{Spec}(A) , \mathcal{N}_{\mathcal{F}_j}^* \right)$ is a foliated isomorphism with $\theta_{jk}\theta_{ij}= \theta_{ik}$ on $U_{ijk}\times \textnormal{Spec}(A)$ and $\theta_{ij}^*\left( \mathcal{N}_{\mathcal{F}_j}^* \right)= \mathcal{N}_{\mathcal{F}_i}^*$, and $\theta_{ji}^*\left( w_i^\alpha \right)= \sum_{\beta=1}^q h_{ij}^{\alpha\beta} w_j^\beta$ for some $h_{ij}^{\alpha\beta}\in \Gamma\left( U_i, \mathcal{O}_X \right)\otimes A$ which induces $w_{0i}^\alpha= \sum_{\beta=1}^q h_{0ij}^{\alpha\beta} w_{0j}^\beta$ modulo $\mathfrak{m}_A$. We let $f_{ij}:\left( \Gamma\left( U_{ij}, \mathcal{O}_X \right)\otimes A , \mathcal{N}_{\mathcal{F}_j}^* \right) \to \left( \Gamma\left( U_{ij}, \mathcal{O}_X \right), \mathcal{N}_{\mathcal{F}_i}^* \right)$ be the ring homomorphism corresponding to $\theta_{ij}$.

In oder to see if there exists a lifting of $\left( \mathcal{X}, \mathcal{N}_\mathcal{F}^* \right)$ to $\tilde{A}$, we choose arbitrary automorphisms $\left\{ \tilde{\theta}_{ij} \right\}$ which restricts to $\theta_{ij}$ on $U_{ij}\times \textnormal{Spec}(A)$ and arbitrary lifting $\tilde{w}_i^\alpha \in \Gamma\left( U_i, \Omega_X^1 \right) \otimes \tilde{A}$ of $w_i^\alpha$, and arbitrary lifting $\tilde{h}_{ij}^{\alpha\beta}$ of $h_{ij}^{\alpha\beta}$. Then  $\tilde{\theta}_{ijk}:=\tilde{\theta}_{ki}\tilde{\theta}_{jk}\tilde{\theta}_{ij}$ on $U_{ijk}$ corresponds to $d_{ijk}\in \Gamma\left(U_{ijk}, T_X \right)$, and we have $\tilde{\theta}_{ji}^*\left( \tilde{w}_i^\alpha \right)- \sum_{\beta=1}^q \tilde{h}_{ij}^{\alpha\beta} \tilde{w}_j^\beta= t R_{ij}^\alpha$ for some $R_{ij}^\alpha \in \Gamma\left( U_{ij}, \Omega_X^1 \right),\alpha=1,...,q$. We define $\overline{R_{ij}}:\Gamma\left( U_{ij}, \mathcal{N}_{\mathcal{F}_0}^* \right) \to \Gamma\left( U_{ij}, \frac{\Omega_X^1}{\mathcal{N}_{\mathcal{F}_0}^*} \right), w_{0i}^\alpha \mapsto \overline{R_{ij}^\alpha}$. Let $\mathcal{N}_{\tilde{\mathcal{F}}_i}^*$ be a subsheaf of $\Omega_X^1|_{U_i}\otimes \tilde{A}$ generated by $\tilde{w}_i^1,..., \tilde{w}_i^q$. We set $T_{\tilde{\mathcal{F}}_i} : =\mathscr{H}om_{\mathcal{O}_X|_{U_i}\otimes \tilde{A}}\left( \frac{\Omega_X^1|_{U_i}\otimes \tilde{A}}{\mathcal{N}_{\tilde{\mathcal{F}}_i}^*} , \mathcal{O}_X|_{U_i}\otimes  \tilde{A} \right)\subset T_X|_{U_i}\otimes \tilde{A}$. We note that since $\frac{\Omega_{\mathcal{X}/\textnormal{Spec}(A)}^1}{\mathcal{N}_{\mathcal{F}}^*}$ is flat over $\textnormal{Spec}(A)$, $\frac{\Omega_{\mathcal{X}/\textnormal{Spec}(A)}^1}{\mathcal{N}_{\mathcal{F}}^*}$ is locally free on $X-S$. This implies $\frac{\Omega_X^1|_{U_i}\otimes A}{\mathcal{N}_{\mathcal{F}_i}^*}$ is locally  free on $U_i-S$, so that  $T_{\mathcal{F}_i}|_{U_i- S}$ is locally free. Then for $x\in U_i-S$, we can choose an affine neighborhood $U_x \subset U_i-S$ such that $T_{\mathcal{F}_i}$ on $U_x$ is generated by $T_{x}^1,..., T_{x}^p\in \Gamma\left(U_x, T_{\mathcal{F}_i}\right)$ and the induced $T_{0x}^1,..., T_{0x}^p$ by modulo $\mathfrak{m}_A$ generates $\Gamma(U_x, T_{\mathcal{F}_0})$, and for another $y\in U_i-S$ such that $U_x\cap U_y\ne \emptyset$, we have $T_{x}^\alpha= \sum_{\beta=1}^p r_{xy}^{\alpha\beta} T_y^\beta$ for some $r_{xy}^{\alpha\beta}\in \Gamma\left(U_x\cap U_y, \mathcal{O}_{\mathcal{X}}\right)=\Gamma\left( U_x\cap U_y, \mathcal{O}_X\right)\otimes A$. We choose arbitrary lifting $\tilde{T}_x^\alpha \in \Gamma\left( U_x, T_{\tilde{\mathcal{F}}_i }\right)$ of $T_x^\alpha$ \footnote{We note that $\mathcal{N}_{\tilde{\mathcal{F}}_i}^*$ is free. In fact, if $\sum_{\alpha=1}^q a_i^\alpha \tilde{w}_i^\alpha=0$ for some $a_i^\alpha \in \Gamma\left( U_i, \mathcal{O}_X\right) \otimes \tilde{A}$, we see that $\sum_{\alpha=1}^q \overline{a_i^\alpha} w_i^\alpha=0$ by modulo $(t)$. Then $\overline{a_i^\alpha}=0$, so that $a_i^\alpha = t a_{0i}^\alpha$ for some $a_{0i}^\alpha \in \Gamma\left(U_i, \mathcal{O}_X \right)$. This implies that $a_i^\alpha=0$. Then since the exact sequence $0\to \mathcal{N}_{\tilde{\mathcal{F}}_i,x}^*  \to \Omega_{X,x}^1|_{U_i}\otimes \tilde{A}\to \frac{\Omega_{X,x}^1|_{U_i}\otimes \tilde{A}}{\mathcal{N}_{\tilde{\mathcal{F}}_i,x}^*} \to 0  $ induces the exact sequence $0\to \mathcal{N}_{\mathcal{F}_0,x}^* \to \Omega_{X,x}^1 \to \frac{\Omega_{X,x}^1}{\mathcal{N}_{\mathcal{F}_0,x}^*}\to 0$, we see that $\frac{\Omega_{X}^1|_{U_i}\otimes \tilde{A}}{\mathcal{N}_{\tilde{\mathcal{F}}_i}^*}$ is flat over $\tilde{A}$. In particular $\frac{\Omega_X^1|_{U_i}\otimes \tilde{A}}{\mathcal{N}_{\tilde{\mathcal{F}}_i^*} }$ is locally free on $U_x$.  Then we have $T_{\tilde{\mathcal{F}}_i}\otimes_{\tilde{A}} A = \mathscr{H}om_{\mathcal{O}_X|_{U_i}\otimes \tilde{A}}\left( \frac{\Omega_X^1|_{U_i}\otimes \tilde{A}}{\mathcal{N}_{\tilde{\mathcal{F}}_i}^*} , \mathcal{O}_X|_{U_i}\otimes  \tilde{A} \right)\otimes_{\tilde{A}} A = T_{\mathcal{F}_i}$ on $U_x$. Hence a lifting $\tilde{T}_x^\alpha$ of $T_x^\alpha$ exists.}, $\alpha=1,...,p$. Then $i_{\tilde{T}_x^\alpha}\left( \tilde{w}_i^\beta \right) =0$ by our construction for $\beta=1,...,q$. We note that since $i_{T_x^\alpha \wedge T_x^\beta}\left( dw_i^\gamma\right)=0$, we have $i_{\tilde{T}_x^\alpha \wedge \tilde{T}_x^\beta}\left( d\tilde{w}_i^\gamma \right)= t f_{\alpha\beta x}^\gamma$ for some $f_{\alpha\beta x}^\gamma\in \Gamma(U_x, \mathcal{O}_X)$. Then we define
\begin{align*}
\phi_{ix}^\gamma : \Gamma\left(U_x, \bigwedge^2 T_{\mathcal{F}_0}\right) &\to \Gamma\left(U_x, \mathcal{O}_X \right) \\
                        T_{0x}^\alpha \wedge T_{0x}^\beta &\mapsto  t f_{\alpha\beta x}^\gamma
\end{align*}
and then linearly extends to $\Gamma\left(U_x, \bigwedge^2 T_{\mathcal{F}_0}\right)$. This is independent of choice of liftings $\tilde{T}_x^\alpha$ of $T_x^\alpha$. In fact, for another choice $\tilde{T}_x'^\alpha \in \Gamma\left( U_x, T_{\tilde{\mathcal{F}}_i} \right)$, we have $\tilde{T}_x^\alpha - \tilde{T}_x'^\alpha = t V_x^\alpha $ for some $V_x^\alpha \in \Gamma\left(U_x, T_X\right)$. Since $i_{\tilde{T}_x^\alpha}\left(\tilde{w}_i^\beta\right)= i_{ \tilde{T}_x'^\alpha}\left(\tilde{w}_i^\beta\right)=0$, we see that $i_{tV_x^\alpha}\left(w_{0i}^\beta\right)=0,\beta=1,...,q$, so that $V_x^\alpha \in \Gamma(U_x, T_{\mathcal{F}_0})$. Then
\begin{align*}
tf_{\alpha\beta x}^\gamma&= i_{\tilde{T}_x^\alpha \wedge \tilde{T}_x^\beta}\left(d\tilde{w}_i^\gamma\right)= i_{\left(\tilde{T}_x'^\alpha + t V_x^\alpha\right) \wedge \left( \tilde{T}_x'^\beta + t V_x^\beta \right) } \left( d \tilde{w}_i^\gamma \right) = i_{\tilde{T}_x'^\alpha \wedge \tilde{T}_x'^\beta}\left( d\tilde{w}_i^\gamma \right) + t i_{V_x^\alpha \wedge T_{0x}^\beta}\left( dw_{0i}^\gamma \right) + i_{T_{0x}^\alpha \wedge V_x^\beta}\left( d w_{0i}^\gamma \right) \\
&=  i_{\tilde{T}_x'^\alpha \wedge \tilde{T}_x'^\beta}\left( d\tilde{w}_i^\gamma \right)= t f_{\alpha\beta x}'^\gamma \Longrightarrow  f_{\alpha\beta x}^\gamma= f_{\alpha\beta x}'^\gamma
\end{align*}
On the other hand, let $\tilde{r}_{xy}^{\alpha\beta}\in \Gamma\left( U_x \cap U_y, \mathcal{O}_X \right)\otimes \tilde{A}$ be an arbitrary lifting of $r_{xy}^{\alpha\beta}$. Then $\tilde{T}_x^\alpha- \sum_{\beta=1}^p \tilde{r}_{xy}^{\alpha\beta} \tilde{T}_y^\beta = t B_{xy\alpha}$ for some $B_{xy\alpha}\in \Gamma\left( U_x \cap U_y , \mathcal{O}_X \right)$. Then $ti_{B_{xy\alpha}}\left(w_{0i}^\gamma \right)= i_{\tilde{T}_x^\alpha- \sum_{\beta=1}^p \tilde{r}_{xy}^{\alpha\beta} \tilde{T}_y^\beta }\left(\tilde{w}_i^\gamma \right)=0,\gamma=1,...,q$. This implies that $B_{xy\alpha}\in \Gamma(U_x\cap U_y, T_{\mathcal{F}_0})$. Then we have
\begin{align*}
t f_{\alpha\beta x}^\gamma &=i_{\tilde{T}_x^\alpha \wedge \tilde{T}_y^\beta}\left(d \tilde{w}_i^\gamma \right) = i_{\left( t B_{xy\alpha} + \sum_{\delta=1}^p \tilde{r}_{xy}^{\alpha\delta} \tilde{T}_y^\delta \right) \wedge \left( t B_{xy\beta}  + \sum_{\eta=1}^p \tilde{r}_{xy}^{\beta \eta} \tilde{T}_y^\eta \right)}\left( d\tilde{w}_i^\gamma \right) \\
&=t i_{B_{xy\alpha} \wedge T_{0x}^\beta }(dw_{0i}^\gamma) + t i_{ T_{0x}^\alpha \wedge B_{xy\beta}}(dw_{0i}^\gamma) + t\sum_{\delta, \eta=1}^p r_{0xy}^{\alpha \delta} r_{0xy}^{\beta \eta} f_{\delta \eta y}^\gamma\\
\Longrightarrow & f_{\alpha\beta x}^\gamma = \sum_{\delta, \eta=1}^p r_{0xy}^{\alpha \delta} r_{0xy}^{\beta \eta} f_{\delta \eta y}^\gamma
\end{align*}

Then we can glue together $\{\phi_{ix}^\gamma \}$ to define a $\phi_i^\gamma \in \Gamma\left( U_i-S, \mathscr{H}om_{\mathcal{O}_X}\left( \bigwedge^2 T_{\mathcal{F}_0}, \mathcal{O}_X \right) \right)$, and hence $\phi_i^\gamma \in \Gamma\left( U_i, \mathscr{H}om_{\mathcal{O}_X} \left( \bigwedge^2 T_{\mathcal{F}_0}, \mathcal{O}_X   \right) \right)$. Then we define $\phi_i\in \Gamma\left( U_i, \mathscr{H}om_{\mathcal{O}_X}\left( \mathcal{N}_{\mathcal{F}_0}^*, \left(\bigwedge^2 T_{\mathcal{F}_0} \right)^* \right) \right)$ by
\begin{align*}
\phi_i :\Gamma\left( U_i, \mathcal{N}_{\mathcal{F}_0}^* \right) &\to \Gamma\left( U_i, \mathscr{H}om_{\mathcal{O}_X}\left( \bigwedge^2 T_{\mathcal{F}_0}, \mathcal{O}_X \right) \right) \\
 w_{0i}^\gamma &\mapsto \phi_i^\gamma
\end{align*}
and then linearly extends to $\Gamma\left( U_i, \mathcal{N}_{\mathcal{F}_0}^* \right)$. Then we claim that
\begin{align*}
\left( \left\{  \phi_i \right\},  \left\{  - \overline{R_{ij}}\right\} , \left\{ d_{ijk} \right\} \right) \in C^0\left( \mathcal{U}, \mathscr{H}om_{\mathcal{O}_X}\left(\mathcal{N}_{\mathcal{F}_0}^*, \left(\bigwedge^2 T_{\mathcal{F}_0} \right)^* \right) \right) \bigoplus C^1\left(\mathcal{U}, \mathscr{H}om_{\mathcal{O}_X}\left( \mathcal{N}_{\mathcal{F}_0}^*, \frac{ \Omega_X^1 }{\mathcal{N}_{\mathcal{F}_0}^*} \right) \right) \bigoplus C^2\left( \mathcal{U}, T_X \right)
\end{align*}
defines a $2$-cocycle in the following \v Cech resolution of $\mathcal{N}_{\mathcal{F}_0}'^{*\bullet}$:
\begin{center}
$\begin{CD}
\cdots \\
@AAA \\
C^0\left( \mathcal{U}, \mathscr{H}om_{\mathcal{O}_X}\left( \mathcal{N}_{\mathcal{F}_0}^*, \left( \bigwedge^2 T_{\mathcal{F}_0} \right)^* \right) \right) @>>> \cdots \\
@AAA @AAA \\
C^0\left( \mathcal{U}, \mathscr{H}om_{\mathcal{O}_X}\left( \mathcal{N}_{\mathcal{F}_0}^*, \frac{\Omega_X^1}{\mathcal{N}_{\mathcal{F}_0}^*} \right) \right) @ >>> C^0\left( \mathcal{U},  \mathscr{H}om_{\mathcal{O}_X}\left( \mathcal{N}_{\mathcal{F}_0}^*, \frac{\Omega_X^1}{\mathcal{N}_{\mathcal{F}_0}^*}  \right) \right) @>>> \cdots \\
@AAA @AAA @AAA \\
C^0\left( \mathcal{U}, T_X \right) @>-\delta >> C^0\left( \mathcal{U}, T_X \right) @>\delta >> C^0\left( \mathcal{U}, T_X \right) @>-\delta >> \cdots
\end{CD}$
\end{center}
As in the proof of Proposition \ref{nb1}, $-\delta\left(\left\{ d_{ijk} \right\} \right)=0$ and $-\delta\left( \left\{ - R_{ij}\right\} \right) + E_0\left( \left\{ d_{ijk} \right\} \right)=0$. We show that $E_2'\left(\left\{ \phi_i\right\} \right)=0$. In fact, for $x\in U_i-S$, choose an affine neighborhood of $U_x\subset U_i-S$ as above. Then we check locally on $U_x\subset U_i-S$. We note that $dw_{i}^\alpha= \sum_{\beta=1}^q a_{x}^{\alpha\beta} \wedge w_{i}^\beta$ for some $a_x^{\alpha\beta}\in \Gamma\left(U_x, \Omega_X^1 \right)\otimes A$. Then $i_{T_{x}^\eta}\left( dw_{i}^\alpha\right)= \sum_{\beta=1}^q i_{T_{x}^\eta}\left( a_{x}^{\alpha\beta} \right) w_{i}^\beta $. We let $\tilde{a}_x^{\alpha\beta}\in \Gamma\left(U_x, \Omega_X^1 \right)\otimes \tilde{A}$ be a lifting of $a_{x}^{\alpha\beta}$. We also note that $\left[ T_{x}^\alpha, T_{x}^\beta \right]= \sum_{\gamma=1}^p g_{\alpha\beta x}^\gamma T_{x}^\gamma$ for some $g_{\alpha\beta x}^\gamma\in \Gamma\left(U_x, \mathcal{O}_X \right)\otimes A$. We let $\tilde{g}_{\alpha\beta x}^\gamma\in \Gamma\left( U_x , \mathcal{O}_X \right)\otimes \tilde{A}$ be a lifting of $g_{\alpha\beta x}^\gamma$. Then we have 
\begin{align*}
d\tilde{w}_i^\alpha - \sum_{\beta=1}^q \tilde{a}_x^{\alpha\beta} \wedge \tilde{w}_i^\beta= t b_{x}^\alpha \Longrightarrow i_{\tilde{T}_x^\gamma}\left(d\tilde{w}_i^\alpha \right)- \sum_{\beta=1}^q i_{\tilde{T}_x^\gamma}\left( \tilde{a}_x^{\alpha\beta} \right) \tilde{w}_i^\beta = t i_{T_{0x}^\gamma}(b_x^\alpha) \\
\Longrightarrow \mathcal{L}_{\tilde{T}_x^\gamma} \left(\tilde{w}_i^\alpha \right)- \sum_{\beta=1}^q i_{\tilde{T}_x^\gamma}\left( \tilde{a}_x^{\alpha\beta} \right) \tilde{w}_i^\beta = t i_{T_{0x}^\gamma}\left(b_x^\alpha \right)
\end{align*}
for some $b_x^\alpha\in \Gamma\left( U_x, \bigwedge^2 \Omega_X^1 \right)$ and we have $i_{\tilde{T}_x^\alpha \wedge \tilde{T}_x^\beta}\left(d\tilde{w}_i^\gamma  \right)= i_{\tilde{T}_x^\beta} i_{\tilde{T}_x^\alpha}\left( d\tilde{w}_i^\gamma \right) =i_{\tilde{T}_x^\beta} \mathcal{L}_{\tilde{T}_x^\alpha}\left(\tilde{w}_i^\gamma \right)  = - i_{\left[ \tilde{T}_x^\alpha, \tilde{T}_x^\beta \right]}\left(\tilde{w}_i^\gamma\right)$. Then we have
{\small{\begin{align*}
&tE_2'(\phi_i)\left(T_{0x}^\alpha \wedge T_{0x}^\beta \wedge T_{0x}^\gamma \right)\left(w_{0i}^\delta \right) \\
&= t\left[ T_{0x}^\alpha, \phi_i\left(T_{0x}^\beta \wedge T_{0x}^\gamma\right)\left(w_{0i}^\delta\right)  \right] - t\phi\left(T_{0x}^\beta \wedge T_{0x}^\gamma \right)\left(\mathcal{L}_{T_{0x}^\alpha}\left(w_{0i}^\delta \right) \right) - t\left[ T_{0x}^\beta, \phi_i\left( T_{0x}^\alpha \wedge T_{0x}^\gamma \right)\left(w_{0i}^\delta \right)   \right] + t\phi\left(T_{0x}^\alpha \wedge T_{0x}^\gamma \right)\left(\mathcal{L}_{T_{0x}^\beta}\left(w_{0i}^\delta \right)\right)\\
&+t\left[ T_{0x}^\gamma, \phi_i\left( T_{0x}^\alpha \wedge T_{0x}^\beta \right)\left(w_{0i}^\delta \right)   \right] - t\phi_i\left(T_{0x}^\alpha \wedge T_{0x}^\beta \right)\left(\mathcal{L}_{T_{0x}^\gamma}\left(w_{0i}^\delta \right)\right) - t \phi_i\left(\left[ T_{0x}^\alpha, T_{0x}^\beta\right], T_{0x}^\gamma    \right)\left(w_{0i}^\delta\right) + t\phi\left(\left[ T_{0x}^\alpha, T_{0x}^\gamma \right], T_{0x}^\beta \right) \\
&- t \phi_i \left(\left[ T_{0x}^\beta, T_{0x}^\gamma\right] , T_{0x}^\alpha \right)(w_{0i}^\delta)\\
&= \left[ \tilde{T}_x^\alpha, i_{\tilde{T}_x^\beta \wedge \tilde{T}_x^\gamma}\left(d \tilde{w}_i^\delta \right) \right] -  \sum_{\eta=1}^q i_{\tilde{T}_{x}^\alpha}\left( \tilde{a}_x^{\delta\eta} \right) i_{\tilde{T}_x^\beta \wedge \tilde{T}_x^\gamma}\left(d\tilde{w}_i^\eta \right) -  \left[ \tilde{T}_x^\beta, i_{\tilde{T}_x^\alpha \wedge \tilde{T}_x^\gamma}\left(d \tilde{w}_i^\delta \right) \right] + \sum_{\eta=1}^q i_{\tilde{T}_{x}^\beta}\left( \tilde{a}_x^{\delta\eta} \right) i_{\tilde{T}_x^\alpha \wedge \tilde{T}_x^\gamma}\left(d\tilde{w}_i^\eta \right) \\
&+  \left[ \tilde{T}_x^\gamma, i_{\tilde{T}_x^\alpha \wedge \tilde{T}_x^\beta}\left(d \tilde{w}_i^\delta \right) \right] -  \sum_{\eta=1}^q i_{\tilde{T}_{x}^\gamma}\left( \tilde{a}_x^{\delta\eta} \right) i_{\tilde{T}_x^\alpha \wedge \tilde{T}_x^\beta}\left(d\tilde{w}_i^\eta \right) - \sum_{\eta=1}^p \tilde{g}_{\alpha \beta x}^\eta i_{\tilde{T}_{x}^\eta \wedge \tilde{T}_{x}^\gamma}\left(d\tilde{w}_{i}^\delta \right) + \sum_{\eta=1}^p \tilde{g}_{\alpha \gamma x}^\eta i_{\tilde{T}_x^\eta \wedge \tilde{T}_x^\beta }\left( d \tilde{w}_i^\delta\right) -  \sum_{\eta=1}^p \tilde{g}_{\beta \gamma x}^\eta i_{\tilde{T}_x^\eta \wedge \tilde{T}_x^\alpha }\left( d \tilde{w}_i^\delta\right) \\
&=- \mathcal{L}_{\tilde{T}_x^\alpha} i_{\left[ \tilde{T}_x^\beta, \tilde{T}_x^\gamma  \right]}\left(\tilde{w}_i^\delta \right) + \sum_{\eta=1}^q i_{\tilde{T}_x^\alpha}\left( \tilde{a}_x^{\delta \eta} \right) i_{\left[ \tilde{T}_x^\beta , \tilde{T}_x^\gamma \right]}\left( \tilde{w}_i^\eta \right) + \mathcal{L}_{\tilde{T}_x^\beta} i_{\left[  \tilde{T}_x^\alpha   , \tilde{T}_x^\gamma \right]} \left(\tilde{w}_i^\delta \right) - \sum_{\eta=1}^q i_{\tilde{T}_x^\beta}\left(\tilde{a}_x^{\delta \eta} \right) i_{\left[ \tilde{T}_x^\alpha,  \tilde{T}_x^\gamma \right]}\left( \tilde{w}_i^\eta \right)\\
&- \mathcal{L}_{\tilde{T}_x^\gamma} i_{\left[  \tilde{T}_x^\alpha   , \tilde{T}_x^\beta \right]} \left(\tilde{w}_i^\delta\right) + \sum_{\eta=1}^q i_{\tilde{T}_x^\gamma}\left(\tilde{a}_x^{\delta \eta} \right) i_{\left[ \tilde{T}_x^\alpha,  \tilde{T}_x^\beta \right]}\left( \tilde{w}_i^\eta \right) - \sum_{\eta=1}^p \tilde{g}_{\alpha \beta x}^\eta i_{\tilde{T}_{x}^\eta \wedge \tilde{T}_{x}^\gamma}\left(tb_x^\delta  \right) + \sum_{\eta=1}^p \tilde{g}_{\alpha \gamma x}^\eta i_{\tilde{T}_x^\eta \wedge \tilde{T}_x^\beta }\left(t b_x^\delta \right) -  \sum_{\eta=1}^p \tilde{g}_{\beta \gamma x}^\eta i_{\tilde{T}_x^\eta \wedge \tilde{T}_x^\alpha }\left( t b_x^\delta \right)\\
&=- \mathcal{L}_{\tilde{T}_x^\alpha} i_{\left[ \tilde{T}_x^\beta, \tilde{T}_x^\gamma  \right]}\left(\tilde{w}_i^\delta \right) +  i_{\left[ \tilde{T}_x^\beta , \tilde{T}_x^\gamma \right]} \mathcal{L}_{\tilde{T}_x^\alpha}\left( \tilde{w}_i^\delta \right) - t i_{\left[ T_{0x}^\beta, T_{0x}^\gamma \right]} i_{ T_{0x}^\alpha }\left( b_i^\delta \right) + \mathcal{L}_{\tilde{T}_x^\beta} i_{\left[  \tilde{T}_x^\alpha   , \tilde{T}_x^\gamma \right]}\left(\tilde{w}_i^\delta\right) - i_{\left[ \tilde{T}_x^\alpha , \tilde{T}_x^\gamma \right]} \mathcal{L}_{\tilde{T}_x^\beta}\left( \tilde{w}_i^\delta \right) + t i_{\left[ T_{0x}^\alpha, T_{0x}^\gamma \right]} i_{T_{0x}^\beta} \left( b_i^\delta \right) \\
& - \mathcal{L}_{\tilde{T}_x^\gamma} i_{\left[  \tilde{T}_x^\alpha   , \tilde{T}_x^\beta \right]}(\tilde{w}_i^\delta) + i_{\left[ \tilde{T}_x^\alpha , \tilde{T}_x^\beta \right]} \mathcal{L}_{\tilde{T}_x^\gamma}\left( \tilde{w}_i^\delta \right) - t i_{\left[ T_{0x}^\alpha, T_{0x}^\beta \right]} i_{T_{0x}^\gamma} \left( b_i^\delta \right) - ti_{\left[ T_{0x}^\alpha, T_{0x}^\beta  \right]\wedge T_{0x}^\gamma}\left( b_x^\delta \right) + ti_{\left[ T_{0x}^\alpha, T_{0x}^\gamma  \right]\wedge T_{0x}^\beta}\left( b_x^\delta \right) - ti_{\left[ T_{0x}^\beta, T_{0x}^\gamma  \right]\wedge T_{0x}^\alpha }\left( b_x^\delta \right)\\
&= - \mathcal{L}_{\left(\left[ \tilde{T}_x^\alpha, \left[ \tilde{T}_x^\beta, \tilde{T}_x^\gamma \right] \right]  + \left[ \tilde{T}_x^\beta, \left[ \tilde{T}_x^\gamma, \tilde{T}_x^\alpha \right] \right] + \left[ \tilde{T}_x^\gamma, \left[ \tilde{T}_x^\alpha, \tilde{T}_x^\beta \right] \right] \right) }\left(\tilde{w}_i^\delta \right) =0.
\end{align*}}}
Hence $E_1'\left(\phi_i \right)=0$ on $U_i-S$, and so on $U_i$. Next we claim that $\phi_i-\phi_j+ E_1'\left( \overline{ - R_{ij}}\right)=0$ on $U_i\cap U_j$. It is sufficient to check locally on $U_i\cap U_j - S$. Let $x\in U_x\subset U_i-S$ and $y\in U_y\subset U_j-S$ with $U_x\cap U_y\ne \emptyset$. Since $\tilde{\theta}_{ji}^* \left(\tilde{w}_i^\alpha \right) - \sum_{\beta=1}^q \tilde{h}_{ij}^{\alpha\beta} \tilde{w}_j^\beta= t R_{ij}^\alpha$ for some $R_{ij}^\alpha \in \Gamma\left( U_{ij}, \Omega_X^1 \right)$, we have $\tilde{w}_i^\alpha - \sum_{\beta=1}^q \tilde{f}_{ij}\left( \tilde{h}_{ij}^{\alpha\beta} \right) \theta_{ij}^*\left(\tilde{w}_j^\beta \right)= t R_{ij}^\alpha$.  On the other hand,  since $\theta_{ji}^*\left( T_{i_x}^\alpha \right)- \sum_{\beta=1}^p r_{i_x j_y}^{\alpha\beta} T_{j_y}^\beta=0$ for some $r_{i_x j_y}^{\alpha\beta}\in \Gamma\left( U_x\cap U_y, \mathcal{O}_X \right)\otimes A$, by choosing arbitrary lifting of $\tilde{r}_{i_x j_y}^{\alpha\beta}\in \Gamma\left( U_x \cap U_y, \mathcal{O}_X \right)\otimes \tilde{A}$ of $r_{i_x j_y}^{\alpha\beta}$ and liftings $\tilde{T}_{i_x}^\alpha, \tilde{T}_{j_y}^\beta$ of $T_{i_x}^\alpha, T_{j_y}$ as above. Then we can write $\tilde{\theta}_{ji}^*\left(\tilde{T}_{i_x}^\alpha\right) - \sum_{\beta=1}^p \tilde{r}_{i_xj_y}^{\alpha\beta} \tilde{T}_{j_y}^\beta = t V_{i_xi_y}^\alpha$ for some $V_{i_x j_y}^\alpha\in \Gamma\left( U_x\cap U_y, T_X\right)$, equivalently $\tilde{T}_{i_x}^\alpha - \sum_{\beta=1}^p \tilde{f}_{ij}\left(\tilde{r}_{i_xi_y}^{\alpha\beta} \right)\tilde{\theta}_{ij}^*\left( \tilde{T}_{j_y}^\beta \right)= tV_{i_xi_y}^\alpha$. Then we have
{\small{\begin{align*}
t i_{V_{i_xj_y}^\alpha}\left( w_{0i}^\gamma \right) & = i_{ \tilde{T}_{i_x}^\alpha - \sum_{\beta=1}^p \tilde{f}_{ij}\left(\tilde{r}_{i_xi_y}^{\alpha\beta} \right)\tilde{\theta}_{ij}^*\left( \tilde{T}_{j_y}^\beta \right) }\left(\tilde{w}_i^\gamma \right) = - i_{\sum_{\beta=1}^p \tilde{f}_{ij}\left(\tilde{r}_{i_xi_y}^{\alpha\beta} \right)\tilde{\theta}_{ij}^*\left( \tilde{T}_{j_y}^\beta \right) }\left( t R_{ij}^\gamma + \sum_{\delta=1}^q \tilde{f}_{ij}\left(\tilde{h}_{ij}^{\gamma\delta} \right)  \tilde{\theta}_{ij}^*\left( \tilde{w}_j^\delta \right) \right)= - ti_{T_{0i_x}^\alpha}\left( R_{ij}^\gamma \right),
\end{align*}}}
so that $i_{V_{i_xj_y}^\alpha}\left( w_{0i}^\gamma \right)= -  i_{T_{0i_x}^\alpha}\left( R_{ij}^\gamma \right)$. We note that
\begin{align*}
&E_1'\left( R_{ij} \right)\left(T_{0i_x}^\alpha \wedge T_{0i_x}^\beta\right)\left(w_{0i}^\gamma\right) \\
&= \left[ T_{0i_x}^\alpha, i_{T_{0i_x}^\beta} \left( R_{ij}^\gamma \right) \right] - \sum_{\eta=1}^p i_{T_{0i_x}^\alpha} \left( a_{0i_x}^{\gamma \eta} \right) i_{T_{0i_x}^\beta}\left( R_{ij}^\eta \right) - \left[ T_{0i_x}^\beta, i_{T_{0i_x}^\alpha}\left( R_{ij}^\gamma \right) \right]  + \sum_{\eta=1}^p i_{T_{0i_x}^\beta} \left( a_{0i_x}^{\gamma \eta} \right) i_{T_{0i_x}^\alpha}\left( R_{ij}^\eta \right) - i_{\left[ T_{0i_x}^\alpha, T_{0i_x}^\beta \right]}\left( R_{ij}^\gamma \right)
\end{align*}
Then we have
{\small{\begin{align*}
&t\phi_i\left(T_{0i_x}^\alpha \wedge T_{0i_x}^\beta\right)(w_{0i}^\gamma)= i_{\tilde{T}_{i_x}^\alpha \wedge \tilde{T}_{i_x}^\beta}\left(d \tilde{w}_i^\gamma \right) = -i_{\left[\tilde{T}_{i_x}^\alpha , \tilde{T}_{i_x}^\beta \right]}\left(\tilde{w}_i^\gamma \right)\\
&= - i_{\left[ t V_{i_xi_y}^\alpha + \sum_{\delta=1}^p \tilde{f}_{ij}\left( \tilde{r}_{i_xi_y}^{\alpha \delta} \right) \tilde{\theta}_{ij}^*\left( \tilde{T}_{j_x}^\delta \right),  t V_{i_xi_y}^\beta + \sum_{\sigma=1}^p \tilde{f}_{ij}\left( \tilde{r}_{i_xi_y}^{\alpha \sigma} \right) \tilde{\theta}_{ij}^*\left( \tilde{T}_{j_x}^\sigma \right)  \right]} \left( t R_{ij}^\gamma + \sum_{\eta=1}^q \tilde{f}_{ij}\left(\tilde{h}_{ij}^{\gamma \eta} \right) \tilde{\theta}_{ij}^*\left( \tilde{w}_j^\eta  \right) \right)\\
&= -ti_{\left[ V_{i_xi_y}^\alpha, T_{0i_x}^\beta \right]} \left( w_{0i}^\gamma \right) - t i_{\left[ T_{0i_x}^\alpha,  V_{i_xi_y}^\beta \right]}\left(w_{0i}^\gamma \right) -t i_{\left[ T_{0i_x}^\alpha, T_{0i_x}^\beta \right]}\left( R_{ij}^\gamma \right) \\
&+ i_{\left( \sum_{\delta=1}^p \tilde{f}_{ij}\left( \tilde{r}_{i_xi_y}^{\alpha \delta} \right) \tilde{\theta}_{ij}^*\left( \tilde{T}_{j_x}^\delta \right) \right) \wedge \left(\sum_{\sigma=1}^q \tilde{f}_{ij}\left(\tilde{r}_{i_xj_y}^{\beta \sigma} \right) \tilde{\theta}_{ij}^*\left( \tilde{T}_{j_y}^\sigma  \right)  \right)}d\left(  \sum_{\eta=1}^q \tilde{f}_{ij}\left(\tilde{h}_{ij}^{\gamma \eta} \right) \tilde{\theta}_{ij}^*\left( \tilde{w}_j^\eta  \right) \right) \\
&= t \mathcal{L}_{T_{0_{i_x}}^\beta} i_{V_{i_xi_y}^\alpha}\left(w_{0i}^\gamma \right) - t i_{V_{i_x i_y}^\alpha}\mathcal{L}_{T_{0i_x}^\beta}\left( w_{0i}^\gamma \right)-  t \mathcal{L}_{T_{0_{i_x}}^\alpha} i_{V_{i_xi_y}^\beta}\left(w_{0i}^\gamma \right) + i_{V_{i_x i_y}^\beta}\mathcal{L}_{T_{0i_x}^\alpha}\left( w_{0i}^\gamma \right) - t i_{\left[ T_{0i_x}^\alpha, T_{0i_x}^\beta \right]} \left( R_{ij}^\gamma \right) \\
& + \sum_{\eta=1}^q \tilde{f}_{ij}\left(\tilde{h}_{ij}^{\gamma \eta} \right)\theta_{ij}^*\left( i_{\left( \sum_{\delta=1}^p \tilde{r}_{i_x j_y}^{\alpha \delta} \tilde{T}_{j_y}^\delta \right)\wedge \left( \sum_{\sigma=1}^p \tilde{r}_{i_x j_y}^{\beta \sigma} \tilde{T}_{j_y}^\sigma \right)}\left( d  \tilde{w}_j^\eta \right)      \right) \\
&= - t\left[ T_{0i_x}^\beta, i_{T_{0i_x}^\alpha} \left( R_{ij}^\gamma \right) \right] + t\sum_{\eta=1}^q i_{T_{0i_x}^\beta}\left( a_{0i_x}^{\gamma \eta} \right) i_{T_{0i_x}^\alpha}\left( R_{ij}^\eta \right)  +  t\left[ T_{0i_x}^\alpha, i_{T_{0i_x}^\beta} \left( R_{ij}^\gamma \right) \right] - t\sum_{\eta=1}^q i_{T_{0i_x}^\alpha}\left( a_{0i_x}^{\gamma \eta} \right) i_{T_{0i_x}^\beta}\left( R_{ij}^\eta \right) \\
& - t i_{\left[ T_{0i_x}^\alpha, T_{0i_x}^\beta \right]}\left( R_{ij}^\gamma \right) + t \sum_{\eta=1}^q h_{0ij}^{\gamma \eta}\phi_j\left( T_{0i_x}^\alpha \wedge T_{0i_x}^\beta \right)\left(w_{0j}^\eta \right) = tE_1'\left( R_{ij} \right)\left( T_{0i_x}^\alpha \wedge T_{0i_x}^\beta \right) \left(w_{0i}^\gamma \right) +t \phi_j\left( T_{0i_x}^\alpha \wedge T_{0i_x}^\beta \right)  \left(w_{0i}^\gamma \right)
\end{align*}}}
This implies that $\phi_i - \phi_j - E_1'\left( \overline{R_{ij}} \right)=0$ on $U_i \cap U_j - S$, and so $U_i \cap U_j$. This proves the claim.

We show that for another lifting $\tilde{\theta}_{ij}'$ of $\theta_{ij}$, and another lifting $\tilde{w}_i'^\alpha$ of $w_i^\alpha$, and another lifting of $\tilde{h}_{ij}'^{\alpha\beta}$ of $h_{ij}^{\alpha\beta}$, the associated $2$-cocycle $\beta:=\left( \left\{ \phi_i '\right\} , \left\{ \overline{- R_{ij}' }\right\} , \left\{ d_{ijk}'\right\} \right)$ is cohomologous to $\alpha:=\left( \left\{ \phi_i \right\}, \left\{ - R_{ij} \right\}, \left\{ d_{ijk} \right\} \right)$. Then $\tilde{f}_{ij}'= \tilde{f}_{ij} + t p_{ij}$ for some $p_{ij} \in \Gamma\left(U_{ij}, T_X \right)$ and $\tilde{w}_i'^\alpha = \tilde{w}_i^\alpha + t V_i^\alpha$ for some $V_i^\alpha \in \Gamma\left( U_i, \Omega_X^1 \right)$, and $\tilde{h}_{ij}'^{\alpha\beta}= \tilde{h}_{ij}^{\alpha\beta} + t y_{ij}^{\alpha\beta}$ for some $y_{ij}^{\alpha\beta} \in \Gamma\left( U_{ij}, \mathcal{O}_X \right)$. We define $\overline{V_i}\in  \Gamma\left( U_i, \mathscr{H}om_{\mathcal{O}_X} \left( \mathcal{N}_{\mathcal{F}_0}^*, \frac{\Omega_X^1}{\mathcal{N}_{\mathcal{F}_0}^* } \right) \right), w_{0i}^\alpha \mapsto \overline{V_i^\alpha}$. Then $\delta\left( \left\{ \overline{V_i} \right\}\right)= \left\{- R_{ij}'+ R_{ij} \right\}$ and $\delta\left( \left\{ p_{ij} \right\} \right)=\left\{ d_{ijk}'  - d_{ijk} \right\}$ as in the proof of Proposition \ref{nb1}. On the other hand, we claim that $E_1'\left(\overline{V_i} \right)=\left\{ \phi_i'- \phi_i \right\}$. In fact, first we note that we have $T_x'^\alpha= T_x^\alpha + t W_x^\alpha$ for some $W_i^\alpha \in \Gamma\left( U_x, T_X\right)$. Then
\begin{align*}
0= i_{\tilde{T}_x'^\alpha}\left( \tilde{w}_i'^\delta \right) = i_{\tilde{T}_x^\alpha + t W_x^\alpha}\left( \tilde{w}_i^\delta + t V_i^\delta \right) = t i_{T_{0x}^\alpha}\left( V_i^\delta \right) + t i_{W_x^\alpha}\left( w_{0i}^\delta \right) \Longrightarrow  i_{W_x^\alpha}\left( w_{0i}^\delta \right) = - i_{T_{0x}^\alpha}\left(  V_i^\delta \right)
\end{align*}
Then we have
\begin{align*}
t f_{\alpha\beta x}'^\gamma =i_{\tilde{T}_x'^\alpha\wedge \tilde{T}_x'^\beta}\left( d \tilde{w}_i'^\gamma \right) &= i_{\left( \tilde{T}_x^\alpha + t W_x^\alpha \right) \wedge \left( \tilde{T}_x^\beta + t W_x^\beta \right) } \left( d\left( \tilde{w}_i^\gamma + t V_i^\gamma \right) \right)\\
& =t f_{\alpha\beta x}^\gamma + t i_{T_{0x}^\alpha \wedge W_x^\beta}\left( dw_{0i}^\gamma \right) + ti_{W_x^\alpha \wedge T_{0x}^\beta}\left( d w_{0i}^\gamma \right) + t i_{T_{0x}^\alpha \wedge T_{0x}^\beta}\left( dV_i^\gamma \right) \\
&= tf_{\alpha\beta x}^\gamma+ \sum_{\delta=1}^q a_{0\alpha x}^{\gamma \delta} i_{W_x^\beta}\left( w_{0i}^\delta \right) - \sum_{\delta=1}^q a_{0\beta x}^{\gamma \delta} i_{W_x^\alpha}\left(  w_{0i}^\delta \right) + i_{T_{0x}^\alpha \wedge T_{0x}^\beta}\left( dV_i^\gamma \right) \\
&= t f_{\alpha\beta x}^\gamma - \sum_{\delta=1}^q a_{0\alpha x}^{\gamma \delta} i_{T_{0x}^\beta}\left( V_i^\delta \right) + \sum_{\delta=1}^q a_{0\beta x}^{\gamma \delta} i_{T_{0x}^\alpha}\left( V_i^\delta  \right) + i_{ T_{0x}^\alpha \wedge T_{0x}^\beta} \left(  d V_i^\gamma \right)
\end{align*} 
Then we have
\begin{align*}
f_{\alpha\beta x}'^\gamma- f_{\alpha \beta x}^\gamma =  - \sum_{\delta=1}^q a_{0\alpha x}^{\gamma \delta} i_{T_{0x}^\beta}\left( V_i^\delta \right) + \sum_{\delta=1}^q a_{0\beta x}^{\gamma \delta} i_{T_{0x}^\alpha}\left( V_i^\delta  \right) + i_{ T_{0x}^\alpha \wedge T_{0x}^\beta} \left(  d V_i^\gamma \right)
\end{align*}
This implies that $E_1'\left( \left\{ \overline{V_i} \right\} \right)= \left\{ \phi_i' - \phi_i \right\}$. This proves that $\alpha$ is cohomologous to $\beta$. So given a small extension $0 \to (t)\to \tilde{A} \to A \to 0$ and an infinitesimal foliated deformation $\eta$ of $\left( X , \mathcal{N}_{\mathcal{F}_0}^* \right)$ over $A$, we can associate an element $o_\eta(e):=$ the cohomology class of $\alpha\in \mathbb{H}^2\left( X, \mathcal{N}_{\mathcal{F}_0}'^{*\bullet} \right)$. We note that $o_\eta(e)=0$ if and only if there exists $\left(\left\{ \overline{N_i}\right\}, \left\{ d_{ij} \right\} \right)\in C^0\left( \mathcal{U}, \mathscr{H}om_{\mathcal{O}_X}\left( \mathcal{N}_{\mathcal{F}_0}^*, \frac{\Omega_X^1}{\mathcal{N}_{\mathcal{F}_0}^*} \right) \right) \bigoplus C^1\left( \mathcal{U}, T_X \right)$ where $\overline{N_i}:\Gamma\left( U_i, \mathcal{N}_{\mathcal{F}_0}^* \right) \to \Gamma\left( U_i, \frac{\Omega_X^1}{\mathcal{N}_{\mathcal{F}_0}^*} \right), w_{0i}^\alpha \mapsto \overline{N_i^\alpha}$ such that $E_1'\left(N_i \right)=\phi_i$ and $\overline{N_j}- \overline{N_i}+ E_0\left(d_{ij} \right)= - \overline{R_{ij}}$ and $d_{ij}- d_{ik} + d_{jk}= d_{ijk}$. Then $\Phi_{jk}\Phi_{ij}= \Phi_{ik}$ where $\Phi_{ij}:U_{ij}\times \textnormal{Spec}\left(\tilde{A}\right) \to U_{ij}\times \textnormal{Spec}\left( \tilde{A} \right)$ be an isomorphism corresponding $\tilde{f}_{ij}- t d_{ij}$ and $\Phi_{ji}^*\left( \tilde{w}_i^\alpha - t N_i^\alpha \right)= \sum_{\beta=1}^q \left( \tilde{h}_{ij}^{\alpha\beta} - t u_{ij}^{\alpha\beta} \right)\left(  \tilde{w}_j^\beta -t N_j^\beta\right)$ as in the proof of Proposition \ref{nb1}. Then $\left\{\Phi_{ij} \right\}$ glues together to define $\tilde{\mathcal{X}}$ over $\textnormal{Spec}\left(\tilde{A} \right)$ inducing $\mathcal{X}$ and $\mathcal{N}_{\tilde{\mathcal{F}}_i}^*$ generated by $\left\{ \tilde{w}_i^\alpha - t N_i^\alpha \right\}_{\alpha=1,...,q}$ glues together to define a locally free sheaf $\mathcal{N}_{\tilde{\mathcal{F}}}^*$ on $\tilde{\mathcal{X}}$ which induces $\mathcal{N}_\mathcal{F}^*$ on $\mathcal{X}$. Let $T_{\tilde{\mathcal{F}}}:=\mathscr{H}om_{\mathcal{O}_{\tilde{\mathcal{X}}} }\left( \frac{\Omega_{\tilde{\mathcal{X}}/\tilde{A}}^1}{\mathcal{N}_{\tilde{\mathcal{F}}}^1}, \mathcal{O}_{\tilde{\mathcal{X}}} \right)$. Then $T_{\tilde{\mathcal{F}}}$ is locally the dual of $\frac{\Omega_X^1 |_{U_i} \otimes \tilde{A}}{\mathcal{N}_{\tilde{\mathcal{F}}_i}^*}$ on $U_i$ which we denote by $T_{\tilde{\mathcal{F}}_i}$. Let $A_1,A_2\in \Gamma\left( U_i, T_{\tilde{\mathcal{F}}_i} \right)$. For $x\in U_i-S$, choose an affine neighborhood $U_x\subset U_i-S$. Then the image $\overline{A_1},\overline{A_2}$ of $A_1,A_2$ in $\Gamma\left( U_x, T_{\mathcal{F}_i} \right)$ by modulo $(t)$ can be written as $\overline{A_\eta}=\sum_{\alpha=1}^q b_\eta^\alpha T_x^\alpha$ for some $b_\eta^\alpha \in \Gamma\left( U_x, \mathcal{O}_X\right)\otimes A,\eta=1,2$. Let $\tilde{b}_\eta^\alpha \in \Gamma\left(U_x, \mathcal{O}_X \right)\otimes \tilde{A}$ be liftings of $b_\eta^\alpha$.  Then $A_\eta- \sum_{\alpha=1}^p \tilde{b}_\eta^\alpha \tilde{T}_x^\alpha= t P_\eta$ for some $P_\eta \in \Gamma\left(U_x,  T_X \right), \eta=1,2$. Then for $\eta=1,2$ we have on $U_x$
\begin{align*}
0= i_{A_\eta}\left( \tilde{w}_i^\gamma - t N_i^\gamma\right)&= i_{\sum_{\alpha=1}^p \tilde{b}_\eta^\alpha \tilde{T}_x^\alpha + tP_\eta }\left(\tilde{w}_i^\gamma - t N_i^\gamma \right) = -t\sum_{\alpha=1}^p b_{0\eta}^\alpha i_{T_{0x}^\alpha}\left(N_i^\gamma \right)+ t i_{P_\eta}\left( w_{0i}^\gamma \right) \\
&\Longrightarrow \sum_{\alpha=1}^p b_{0\eta}^\alpha i_{T_{0x}^\alpha}\left(N_i^\gamma \right) = i_{P_\eta}\left(w_{0i}^\gamma \right)
\end{align*}
Then since $E_1'(N_i)=\phi_i$, we have on $U_x$
\begin{align*}
&i_{A_1\wedge A_2}\left( d\left(\tilde{w}_i^\gamma - t N_i^\gamma \right) \right)= i_{\left(\sum_{\alpha=1}^p \tilde{b}_1^\alpha \tilde{T}_x^\alpha + t P_1\right) \wedge \left( \sum_{\beta=1}^p \tilde{b}_2^\beta \tilde{T}_x^\beta + t P_2 \right)}\left( d\tilde{w}_i^\gamma - t d N_i^\gamma \right)\\
&= t\sum_{\alpha,\beta=1}^p {b}_{01}^\alpha {b}_{02}^\beta f_{\alpha\beta x}^\gamma +t \sum_{\alpha=1}^p b_{01}^\alpha i_{T_{0x}^\alpha\wedge P_2}\left( dw_{0i}^\gamma \right) +t \sum_{\beta=1}^p b_{02}^\beta i_{P_1 \wedge T_{0x}^\beta}\left( d w_{0i}^\gamma \right) - t \sum_{\alpha,\beta=1}^p b_{01}^\alpha b_{02}^\beta i_{T_{0x}^\alpha \wedge T_{0x}^\beta}\left( d N_i^\gamma \right) \\
&= t\sum_{\alpha,\beta=1}^p b_{01}^\alpha b_{02}^\beta  f_{\alpha\beta x}^\gamma + t \sum_{\alpha=1}^p \sum_{\delta=1}^q b_{01}^\alpha a_{0x\alpha}^\delta i_{P_2}\left( w_{0i}^\delta \right) - t \sum_{\beta=1}^p \sum_{\delta=1}^q b_{02}^\beta a_{0x\beta}^{\delta} i_{P_1}\left(w_{0i}^\delta \right) - t \sum_{\alpha, \beta=1}^p b_{01}^\alpha b_{02}^\beta i_{T_{0x}^\alpha \wedge T_{0x}^\beta}\left( d N_i^\gamma \right)\\
&= t \sum_{\alpha, \beta=1}^p b_{01}^\alpha b_{02}^\beta \left( f_{\alpha \beta x}^\gamma + \sum_{\delta=1}^q a_{0x\alpha}^\delta i_{T_{0x}^\beta}\left( N_i^\delta \right) - \sum_{\delta=1}^q a_{0x\beta}^\delta i_{T_{0x}^\alpha}\left( w_{0i}^\delta \right) - i_{T_{0x}^\alpha \wedge T_{0x}^\beta}\left( d N_i^\gamma \right)   \right) = 0
\end{align*}
This implies that $i_{A_1\wedge A_2}\left(d\left(\tilde{w}_i^\gamma - t N_i^\gamma \right) \right)=0$ on $U_i-S$, and so on $U_i$, so that the integrability condition is satisfied, i.e. $d_{\tilde{\mathcal{X}}/\tilde{A}}\left( \mathcal{N}_{\tilde{\mathcal{F}}}^* \right)$ vanishes on $\bigwedge^2 T_{\tilde{\mathcal{F}}}$. This completes the proof of Proposition \ref{si30}.

\end{proof}

The morphism of complexes of sheaves $(\ref{si41})$ induces the morphism of first-order foliated deformations and obstructions.  In fact, we describe the correspondence of first-order deformations between $\mathbb{H}^1\left( X, \mathcal{N}_{\mathcal{F}_0}^{*\bullet} \right)$ and $\mathbb{H}^1\left( X, \mathcal{N}_{\mathcal{F}_0}'^{*\bullet}\right)$. Since the first-two terms of $\mathcal{N}_{\mathcal{F}_0}^{*\bullet}$ and $\mathcal{N}_{\mathcal{F}_0}'^{*\bullet}$ are identical, we have the identity map (see the proof of $(1)$ in Proposition \ref{nb1})
\begin{align*}
\mathbb{H}^1\left( X, \mathcal{N}_{\mathcal{F}_0}^{*\bullet} \right) &\to \mathbb{H}^1\left( X, \mathcal{N}_{\mathcal{F}_0}'^{*\bullet} \right)\\
\left( \left\{ \overline{w_i'}\right\}, \left\{ p_{ij} \right\} \right) &\mapsto \left( \left\{\overline{w_i'} \right\},  \left\{ p_{ij} \right\} \right)
\end{align*}
We describe the obstruction map. Given a foliated deformation of $\left( X, \mathcal{N}_{\mathcal{F}_0}^* \right)$ over $A\in \bold{Art}$ and a small extension $0\to (t)\to \tilde{A} \to A \to 0$, the obstruction lifting to $\tilde{A}$ in terms of $\mathbb{H}^2\left(X, \mathcal{N}_{\mathcal{F}_0}^{*\bullet} \right)$ is given by $\left( \left\{ B_i \right\}, \left\{ \overline{- R_{ij} } \right\} , \left\{ d_{ijk} \right\} \right)$ as in the proof of Proposition \ref{nb1}. On the other hand, the obstruction lifting to $\tilde{A}$ in terms of $\mathbb{H}^2\left( X, \mathcal{N}_{\mathcal{F}_0}'^{*\bullet} \right)$ is given by $\left( \left\{ \phi_i \right\}, \left\{ \overline{- R_{ij}}\right\}, \left\{ d_{ijk}\right\} \right)$ as above. Hence it is sufficient show that $\alpha_2\left( B_i \right)=\phi_i$ in $(\ref{si41})$. Indeed, we note that for $x\in U_x \subset U_i-S$, $B_i$ is of the form $B_i\left( w_{0i}^\gamma  \right)= w_{0i}^1 \wedge \cdots \wedge w_{0i}^q \wedge C_{i_x}^\gamma$ from $(\ref{si42})$, where $t C_{i_x}^\gamma= d\tilde{w}_i^\gamma - \sum_{\beta=1}^q \tilde{a}_{i_x}^{\gamma\beta} \wedge \tilde{w}_i^\beta$. Then we have on $U_x$
\begin{align*}
t i_{T_{0x}^\alpha \wedge T_{0x}^\beta}\left( C_{i_x}^\gamma \right)= i_{\tilde{T}_x^\alpha \wedge \tilde{T}_x^\beta}\left(t C_{i_x}^\gamma \right)= i_{\tilde{T}_x^\alpha \wedge \tilde{T}_x^\beta}\left( d\tilde{w}_i^\gamma \right)=t f_{\alpha\beta x}^\gamma \Longrightarrow i_{T_{0x}^\alpha \wedge T_{0x}^\beta}\left( C_{i_x}^\gamma \right) = f_{\alpha\beta x}^\gamma =\phi_{i}^\gamma\left(T_{0x}^\alpha \wedge T_{0x}^\beta \right)\left( w_{0i}^\gamma \right)
\end{align*}
Hence we have on $U_x$
\begin{align*}
i_{T_{0x}^\alpha \wedge T_{0x}^\beta}\left( w_{0i}^1 \wedge \cdots \wedge w_{0i}^q \wedge C_{i_x}^\alpha \right) = i_{T_{0x}^\alpha\wedge T_{0x}^\beta}\left( C_{i_x}^\alpha \right) w_{0i}^1 \wedge \cdots \wedge w_{0i}^q= f_{\alpha\beta x}^\gamma w_{0i}^1 \wedge \cdots \wedge w_{0i}^q,
\end{align*}
so that $\alpha_2\left( B_i\right)\left( T_{0x}^\alpha \wedge T_{0x}^\beta \right)\left(w_{0i}^\gamma \right)= f_{\alpha\beta x}^\gamma= \phi_i\left(T_{0x}^\alpha \wedge T_{0x}^\beta \right)\left( w_{0i}^\gamma \right)$. Hence $\alpha_2(B_i)= \phi_i$. This proves the claim.

\section{Rigid foliated deformations, Schlessinger's criterion for $\textnormal{FDef}_{\left( X , \mathcal{N}_{\mathcal{F}_0}^* \right)}$ and stability }

\subsection{Rigid foliated deformations}\
 
We recall that a foliated nonsingular variety $\left( X, T_{\mathcal{F}_0} \right)$ with $\mathcal{N}_{\mathcal{F}_0}^*$ locally free is called rigid in foliated deformations if every infinitesimal foliated deformation $\left(\mathcal{X}, \mathcal{N}_\mathcal{F}^* \right)$ of $\left(X, \mathcal{N}_{\mathcal{F}_0}^* \right)$ over $A$ is trivial for every $A$ in $\bold{Art}$, in other words, $\left(\mathcal{X}, \mathcal{N}_\mathcal{F}^* \right) \cong \left( X\times \textnormal{Spec}(A) , \mathcal{N}_{\mathcal{F}_0}^* \right)$.

\begin{proposition}\label{a11}
A foliated nonsingular variety $\left( X, \mathcal{N}_{\mathcal{F}_0}^* \right)$ with $\mathcal{N}_{\mathcal{F}_0}^*$ locally free is rigid in foliated deformations if and only if $\mathbb{H}^1\left(X, \mathcal{N}_{\mathcal{F}_0}^{*\bullet} \right)=0$.
\end{proposition}

\begin{proof}
Assume that $\left(X, \mathcal{N}_{\mathcal{F}_0}^* \right)$ is rigid in foliated deformations. Since any infinitesimal foliated deformation (in particular, any first-order foliated deformation) is trivial, we have $\mathbb{H}^1\left( X, \mathcal{N}_{\mathcal{F}_0}^{*\bullet} \right)=0$ by Proposition \ref{nb1}. Assume that $\mathbb{H}^1\left( X, \mathcal{N}_{\mathcal{F}_0}^{* \bullet} \right)=0$. First we claim that given an infinitesimal foliated deformation $\eta$ of $\left( X, \mathcal{N}_{\mathcal{F}_0}^* \right)$ over $A\in \bold{Art}$ and a small extension $0\to (t) \to \tilde{A} \to A \to 0$, any two liftings $\xi, \xi'$ of $\eta$ to $\tilde{A}$ are equivalent. Let $\mathcal{U}=\left\{ U_i \right\}$ be an affine open covering of $X$ such that $\xi= \left(\mathcal{X}, \mathcal{N}_\mathcal{F}^*  \right)$ is represented by the following data: we have foliated isomorphisms $\theta_i:\left( U_i \times \textnormal{Spec}\left( \tilde{A} \right), \mathcal{N}_{\mathcal{F}_i}^*\right)\to \left( \mathcal{X}|_{U_i}, \mathcal{N}_\mathcal{F}^*|_{U_i} \right)$ where $\mathcal{N}_{\mathcal{F}_i}^*$ is the foliated structure on $U_i\times \textnormal{Spec}\left( \tilde{A} \right)$ induced from $\mathcal{N}_\mathcal{F}^*$ via $\theta_i$ and $\Gamma\left(U_i , \mathcal{N}_{\mathcal{F}_i}^* \right)$ is generated by $w_i^1,..., w_i^q\in \Gamma\left( U_i, \Omega_X^1 \right)\otimes \tilde{A}$ with $w_i^1\wedge \cdots \wedge w_i^q \wedge dw_i^\alpha=0$, and letting $\theta_{ij}= \theta_j^{-1}\theta_i$, we have foliated isomorphisms $\theta_{ij}: \left( U_{ij}\times \textnormal{Spec}\left( \tilde{A} \right) , \mathcal{N}_{\mathcal{F}_i}^* \right)\to \left( U_{ij}\times \textnormal{Spec}\left( \tilde{A} \right), \mathcal{N}_{\mathcal{F}_j}^* \right)$ such that $\theta_{ji}^* w_i^\alpha= \sum_{\beta=1}^q h_{ij}^{\alpha\beta} w_j^\beta$ for $h_{ij}^{\alpha\beta}\in \Gamma\left(U_{ij}, \mathcal{O}_X \right)\otimes \tilde{A}$, and on the other hand, $\xi'=\left( \mathcal{X}', \mathcal{N}_{\mathcal{F}'}^* \right)$ is represented by the following data: similarly $\theta_i', \theta_{ij}'=\theta_j'^{-1}\theta_i'$ and $\left\{ \mathcal{N}_{\mathcal{F}_i'}^* \right\}, w_i'^\alpha, h_{ij}'^{\alpha\beta}$. We may assume that $w_i^\alpha$ and $w_i'^\alpha$ induce the same element in $\Gamma\left( U_i, \Omega_X^1 \right)\otimes A$. Let $f_{ij}:\left( \Gamma\left( U_{ij} , \mathcal{O}_X \right)\otimes \tilde{A}, \mathcal{N}_{\mathcal{F}_j}^* \right) \to \left( \Gamma\left( U_{ij}, \mathcal{O}_X \right)\otimes \tilde{A}, \mathcal{N}_{\mathcal{F}_i}^* \right)$ be the ring homomorphism corresponding to $\theta_{ij}$ and similarly $f_{ij}'$ corresponding to $\theta_{ij}'$. Since $\xi,\xi'$ induce $\eta$ over $A$, we have $f_{ij}'=f_{ij} +  t p_{ij}$ for some $p_{ij}\in \Gamma\left( U_{ij}, T_X \right)$ and $w_i'^\alpha= w_i^\alpha + t B_i^\alpha$ for some $B_i^\alpha \in \Gamma\left( U_i, \Omega_X^1 \right),\alpha=1,...,q$, and $h_{ij}'^{\alpha\beta}= h_{ij}^{\alpha\beta} + t g_{ij}^{\alpha\beta}$ for some $g_{ij}^{\alpha\beta}\in \Gamma\left( U_{ij}, \mathcal{O}_X \right)$. We define $\overline{B_i}\in \Gamma\left( U_i , \mathscr{H}om_{\mathcal{O}_X} \left( \mathcal{N}_{\mathcal{F}_0}^*, \frac{\Omega_X^1}{\mathcal{N}_{\mathcal{F}_0}^*} \right) \right)$ by $\overline{B_i}:\Gamma\left( U_i, \mathcal{N}_{\mathcal{F}_0}^*  \right) \to \Gamma\left(  U_i, \frac{\Omega_X^1}{\mathcal{N}_{\mathcal{F}_0}^*}\right) ,  w_{0i}^\alpha \mapsto \overline{B_i^\alpha}$. Then we claim that $\left( \left\{ \overline{B_i}\right\}, \left\{ p_{ij} \right\} \right)\in C^0\left( \mathcal{U}, \mathscr{H}om_{\mathcal{O}_X}\left( \mathcal{N}_{\mathcal{F}_0}^*, \frac{\Omega_X^1}{\mathcal{N}_{\mathcal{F}_0}^*} \right) \right) \bigoplus C^1\left( \mathcal{U}, T_X \right)$ defines a $1$-cocycle of the \v Cech resolution of $\mathcal{N}_{\mathcal{F}_0}^{*\bullet}$. In fact, we have $p_{ij}- p_{ik} + p_{jk}=0$, and
\begin{align*}
0&=w_i'^1 \wedge \cdots \wedge w_i'^q \wedge d w_i'^\alpha= \left( w_i^1+ t B_i^1 \right) \wedge \cdots \wedge \left( w_i^q + t B_i^q \right) \wedge d \left( w_i^\alpha + t B_i^\alpha \right) \\
&= w_i^1 \wedge \cdots \wedge w_i^q \wedge d w_i^\alpha +t \sum_{\beta=1}^q w_{0i}^1 \wedge \cdots \wedge B_i^\beta \wedge \cdots \wedge w_{0i}^q \wedge dw_{0i}^\alpha + t w_{0i}^q \wedge \cdots \wedge w_{0i}^q \wedge d B_i^\alpha\\
&\Longrightarrow   \sum_{\beta=1}^q w_{0i}^1 \wedge \cdots \wedge B_i^\beta \wedge \cdots \wedge w_{0i}^q \wedge dw_{0i}^\alpha +  w_{0i}^q \wedge \cdots \wedge w_{0i}^q \wedge d B_i^\alpha = 0
\end{align*}
This implies that $E_1\left( \overline{B}_i \right)=0$. On the other hand, we note that
\begin{align*}
\theta_{ji}'^*\left( w_i'^\alpha \right)&= \theta_{ji}'^*\left( w_i^\alpha + t B_i^\alpha \right) = \theta_{ji}^* w_i^\alpha + t B_i^\alpha  - \mathcal{L}_{p_{ij}}\left( w_{0i}^\alpha \right)\\
\theta_{ji}'^*\left( w_i'^\alpha \right)&= \sum_{\beta=1}^q h_{ij}'^{\alpha\beta} w_j'^\beta= \sum_{\beta=1}^q \left( h_{ij}^{\alpha\beta} + t g_{ij}^{\alpha\beta} \right) \left( w_j^\beta + t B_j^\beta \right)= \sum_{\beta=1}^q h_{ij}^{\alpha\beta} w_j^\beta + t\sum_{\beta=1}^q h_{0ij}^{\alpha\beta} B_j^\beta + t \sum_{\beta=1}^q g_{ij}^{\alpha\beta} w_{0j}^\beta
\end{align*}
This implies that $\delta\left( \left\{ \overline{B_i}\right\} \right) + E_0\left( \left\{ p_{ij} \right\} \right)=0$. This proves the claim. Since $\mathbb{H}^1\left( X, \mathcal{N}_{\mathcal{F}_0}^{*\bullet} \right)=0$, there exists $\left\{a_i \right\}\in C^0\left( \mathcal{U}, T_X  \right)$ such that $a_i- a_j= p_{ij}$ and $\mathcal{L}_{a_i}\left( w_{0i}^\alpha \right) - B_i^\alpha= \sum_{\beta=1}^q e_i^{\alpha\beta} w_{0i}^\beta$ for some $e_i^{\alpha\beta}\in \Gamma\left(  U_i, \mathcal{O}_X \right)$. We explicitly construct a foliated isomorphism $\left( \mathcal{X}, \mathcal{N}_\mathcal{F}^* \right) \cong \left( \mathcal{X}' , \mathcal{N}_{\mathcal{F}'}^* \right)$. We claim that $\textnormal{Id}- t a_i:\left( \Gamma\left( U_i, \mathcal{O}_X \right)\otimes \tilde{A}, \mathcal{N}_{\mathcal{F}_i'}^* \right)\to \left( \Gamma\left( U_i, \mathcal{O}_X \right)\otimes  \tilde{A} , \mathcal{N}_{\mathcal{F}_i}^* \right)$ is a foliated isomorphism. In fact, $w_i'^\alpha- t\mathcal{L}_{a_i}\left( w_{0i}^\alpha \right)= w_i^\alpha +t B_i^\alpha - t \mathcal{L}_{a_i}\left( w_{0i}^\alpha \right)= w_i^\alpha - \sum_{\beta=1}^q e_i^{\alpha\beta} w_{0i}^\beta= \sum_{\beta=1}^q \left( \delta_\alpha^\beta - t e_i^{\alpha\beta} \right) w_i^\beta $. On the other hand, since $p_{ij}=a_i - a_j$, we have $\left( \textnormal{Id}- t a_i \right)\circ f_{ij}'=f_{ij} \circ \left( \textnormal{Id} - t a_j \right)$. Hence $\left\{ \textnormal{Id}- t a_i \right\}$ glues together to give a foliated isomorphism $\left( \mathcal{X}, \mathcal{N}_\mathcal{F}^* \right)\cong \left( \mathcal{X}' , \mathcal{N}_{\mathcal{F}'}^* \right)$.

Now we prove if $\mathbb{H}^1\left( X, \mathcal{N}_{\mathcal{F}_0}^{* \bullet} \right)=0$, then $\left( X, \mathcal{N}_{\mathcal{F}_0}^* \right)$ is rigid in foliated deformations. We prove by induction on the dimension of $A\in \bold{Art}$. Assume that any infinitesimal foliated deformation of $\left( X, \mathcal{N}_{\mathcal{F}_0}^* \right)$ over $A$ with $\dim_k A \leq n-1$ is trivial. Let $\xi$ be an infinitesimal foliated deformation of $\left(X, T_{\mathcal{F}_0} \right)$ over $A$ with $\dim_k A =n$ such that $\mathfrak{m}^{p-1}\ne 0$ and $\mathfrak{m}^p=0$. Choose an element $t\ne 0 \in \mathfrak{m}^{p-1}$. Then $0\to (t)\to A \to A/(t)\to 0$ is a small extension and $\dim_k A/(t)\leq n-1$. Hence the induced foliated deformation $\overline{\xi}$ over $A/(t)$ from $\xi$ is trivial by induction hypothesis. Since $\xi$ is a lifting of $\overline{\xi}$ and the trivial foliated deformation over $A$ is also a lifting of $\overline{\xi}$, $\xi$ is trivial since we have at most one lifting of $\overline{\xi}$. This completes the proof of Proposition \ref{a11}.
\end{proof}

\subsection{Schlessinger's criterion for $\textnormal{FDef}_{\left( X, \mathcal{N}_{\mathcal{F}_0}^* \right)}$}\

\begin{proposition}\label{a12}
Let $\left( X, \mathcal{N}_{\mathcal{F}_0}^* \right)$ be a foliated nonsingular projective variety with $\mathcal{N}_{\mathcal{F}_0}^*$ locally free. Assume that $\mathbb{H}^0\left( X, \mathcal{N}_{\mathcal{F}_0}^{*\bullet}\right)=0$. Then for any infinitesimal foliated deformation $\left( \mathcal{X}, \mathcal{N}_\mathcal{F}^* \right)$ of $\left( X, \mathcal{N}_{\mathcal{F}_0}^*\right)$ over $A\in \bold{Art}$, 
\begin{align*}
\textnormal{Aut}\left( \frac{\left( \mathcal{X}, \mathcal{N}_{\mathcal{F}}^* \right)}{\left( X, \mathcal{N}_{\mathcal{F}_0}^* \right)} \right)= \textnormal{Id}
\end{align*}
where $\textnormal{Aut}\left( \frac{\left( \mathcal{X}, \mathcal{N}_\mathcal{F}^*\right)}{\left( X, \mathcal{N}_{\mathcal{F}_0}^*\right)} \right)$ is  the set of foliated automorphisms of $\left( \mathcal{X}, \mathcal{N}_\mathcal{F}^* \right)$ restricting to the identity on $\left( X, \mathcal{N}_{\mathcal{F}_0}^* \right)$. 
\end{proposition}

\begin{proof}
We prove by the induction on the dimension of $A$. Let's assume that Proposition \ref{a12} holds for $A$ with $\dim_k A\leq n-1$. Let $\dim_k A= n$ and $\left( \mathcal{X}, \mathcal{N}_\mathcal{F}^* \right)$ be an infinitesimal foliated deformation of $\left( X, \mathcal{N}_{\mathcal{F}_0}^* \right)$ over $A$. Assume that the maximal ideal $\mathfrak{m}$ of $A$ satisfies $\mathfrak{m}^{p-1}\ne 0 $ and $\mathfrak{m}^p=0$. Choose $t\ne 0\in \mathfrak{m}^{p-1}$. Then $A/(t)\in \bold{Art}$ with $\dim_k A/(t)\leq n-1$ and $0\to (t) \to A \to A/(t)\to 0$ is a small extension. Now let $g:\left( \mathcal{X}, \mathcal{N}_\mathcal{F}^* \right)\to \left( \mathcal{X}, \mathcal{N}_\mathcal{F}^* \right)$ over $A$ be a foliated automorphism with $g^* \mathcal{N}_\mathcal{F}^*= \mathcal{N}_\mathcal{F}^*$ restricting to the identity on $\left(X, \mathcal{N}_{\mathcal{F}_0}^* \right)$. Let $\mathcal{U}=\left\{ U_i \right\}$ be an affine open cover of $X$ such that $\Gamma\left( U_i , \mathcal{N}_{\mathcal{F}_0}^* \right)$ is generated by $w_{0i}^1,..., w_{0i}^q \in \Gamma\left( U_i , \Omega_X^1 \right)$ with $w_{0i}^\alpha= \sum_{\beta=1}^q h_{0ij}^{\alpha\beta} w_{0j}^\beta$  and $w_{0i}^1\wedge \cdots \wedge w_{0i}^q \wedge dw_{0i}^\alpha= 0,\alpha = 1,...,q$, and $\theta_i:\left( U_i \times \textnormal{Spec}(A), \mathcal{N}_{\mathcal{F}_i}^*  \right) \to \left( \mathcal{X}|_{U_i}, \mathcal{N}_\mathcal{F}^*|_{U_i} \right)$ is a foliated isomorphism such that $\Gamma\left( U_i, \mathcal{N}_\mathcal{F}^* \right)$ is generated by $w_i^1,..., w_i^q\in \Gamma\left( U_i,\Omega_X^1\right)\otimes A$ which induces $w_{0i}^1,..., w_{0i}^q$. Let $\theta_{ij}= \theta_j^{-1}\theta_i$ correspond to a foliated isomorphism $f_{ij}:\left( \Gamma\left( U_{ij}, \mathcal{O}_X \right)\otimes A, \mathcal{N}_{\mathcal{F}_j}^* \right) \to \left( \Gamma\left( U_{ij}, \mathcal{O}_X \right) \otimes A , \mathcal{N}_{\mathcal{F}_i}^*     \right)$. Then $g$ is described by the data $\left\{ g_i \right\}$, where $g_i:\left( \Gamma\left( U_i, \mathcal{O}_X \right)\otimes A, \mathcal{N}_{\mathcal{F}_i}^* \right)\to \left( \Gamma\left( U_i, \mathcal{O}_X \right)\otimes A, \mathcal{N}_{\mathcal{F}_i}^* \right)$ which is a foliated isomorphism and $g_i\circ f_{ij}= f_{ij}\circ g_i$. By the induction hypothesis, $g_i$ induces the identity on $\Gamma\left( U_i, \mathcal{O}_X \right)\otimes A/(t)$. Then $g_i$ is of the form $g_i = \textnormal{Id} + t d_i$ such that $d_i\in \Gamma\left( U_i, T_X \right)$ with $d_i= d_j$ and $w_i^\alpha + t \mathcal{L}_{d_i}\left( w_{0i}^\alpha \right)= \sum_{\beta=1}^q s_i^{\alpha\beta} w_i^\beta$ for some $s_i^{\alpha\beta}\in \Gamma\left( U_i, \mathcal{O}_X \right)\otimes A$. Then $\mathcal{L}_{d_i}\left( w_{0i}^\alpha \right) \in \Gamma\left( U_i, \mathcal{N}_{\mathcal{F}_0}^* \right)$. Hence $d:=\left\{ d_i \right\} \in \mathbb{H}^0\left( X, \mathcal{N}_{\mathcal{F}_0}^{*\bullet} \right)$. Since $\mathbb{H}^0\left( X, \mathcal{N}_{\mathcal{F}_0}^{*\bullet} \right)=0$, we have $d=0$, so that $g= \textnormal{Id}$.
\end{proof}

\begin{proposition} \label{e7}
Let $\left(  X, \mathcal{N}_{\mathcal{F}_0}^* \right)$ be a foliated nonsingular variety with $\mathcal{N}_{\mathcal{F}_0}^*$ locally free, and $e: 0 \to (t) \to \tilde{A} \xrightarrow{\mu} A \to 0$ be a small extension in $\bold{Art}$. Let $p:=\textnormal{FDef}_{\left( X, \mathcal{N}_{\mathcal{F}_0}^* \right)}(\mu): \textnormal{FDef}_{\left( X, \mathcal{N}_{\mathcal{F}_0}^* \right)}\left(\tilde{A} \right)\to \textnormal{FDef}_{\left( X, \mathcal{N}_{\mathcal{F}_0}^* \right)}\left( A \right)$. Then given $\xi=\left( \mathcal{X}, \mathcal{N}_\mathcal{F}^*\right)$, there is a transitive action of $\mathbb{H}^1\left(X, \mathcal{N}_{\mathcal{F}_0}^{*\bullet} \right)$ on $p^{-1}\left( \xi\right)$. Moreover, if $\mathbb{H}^0\left( X, \mathcal{N}_{\mathcal{F}_0}^{*\bullet} \right)=0$, then the action is free.
\end{proposition}

\begin{proof}
We will define a map $G:\mathbb{H}^1\left( X, \mathcal{N}_{\mathcal{F}_0}^{*\bullet} \right)\times p^{-1}(\xi)\to p^{-1}(\xi)$ which is a group action of $\mathbb{H}^1\left( X, \mathcal{N}_{\mathcal{F}_0}^{*\bullet} \right)$ acting on $p^{-1}(\xi)$. Let $\tilde{\xi}=\left( \tilde{\mathcal{X}}, \mathcal{N}_{\tilde{\mathcal{F}}}^*  \right) \in p^{-1}(\xi)$ over $\tilde{A}$. Let $\mathcal{U}=\left\{ U_i \right\}$ be an affine open cover of $X$ such that $\Gamma\left( U_i, \mathcal{N}_{\mathcal{F}_0}^* \right)$ is generated by $w_{0i}^1,..., w_{0i}^q\in \Gamma\left( U_i, \Omega_X^1 \right)$ with $w_{0i}^\alpha= \sum_{\beta=1}^q h_{0ij}^{\alpha\beta} w_{0j}^\beta$ and $w_{0i}^1\wedge \cdots \wedge w_{0i}^q \wedge dw_{0i}^\alpha=0$, and $\tilde{\theta}_i: \left( U_i\times \textnormal{Spec}\left( \tilde{A} \right), \mathcal{N}_{\tilde{\mathcal{F}}_i}^* \right) \to \Gamma\left( \mathcal{X}|_{U_i}, \mathcal{N}_\mathcal{F}^*|_{U_i} \right)$ is a foliated isomorphism such that $\Gamma\left( U_i, \mathcal{N}_{\tilde{\mathcal{F}}_i}^* \right)$ is generated by $\tilde{w}_i^1,..., \tilde{w}_i^q\in \Gamma\left( U_i , \Omega_X^1 \right)\otimes \tilde{A}$ which induces $w_{0i}^1,..., w_{0i}^q$, and $\tilde{\theta}_{ji}^* \tilde{w}_i^\alpha = \sum_{\beta=1}^q \tilde{h}_{ij}^{\alpha\beta} \tilde{w}_j^\beta$ where $\tilde{\theta}_{ij}= \tilde{\theta}_j^{-1}\tilde{\theta}_i$ and $\tilde{h}_{ij}^{\alpha\beta}$ induces $h_{0ij}^{\alpha\beta}$, and $\tilde{w}_i^1 \wedge \cdots \wedge \tilde{w}_i^q \wedge d\tilde{w}_i^\alpha=0$. Let $\tilde{\theta}_{ij}$ correspond to a foliated isomorphism $\tilde{f}_{ij}: \left(  \Gamma\left( U_{ij}, \mathcal{O}_X \right) \otimes \tilde{A}, \mathcal{N}_{\tilde{\mathcal{F}}_j}^* \right)\to \left( \Gamma\left( U_{ij}, \mathcal{O}_X \right)\otimes \tilde{A}, \mathcal{N}_{\tilde{\mathcal{F}}_i}^* \right)$. Let $\left( \left\{ \overline{w_i'}\right\}, \left\{ p_{ij} \right\} \right)\in \mathbb{H}^1\left( X, \mathcal{N}_{\mathcal{F}_0}^{*\bullet} \right)$ corresponding to a first-order foliated deformation as in the proof of Proposition \ref{nb1}. We keep the notations in the proof of Proposition \ref{nb1}. Then $\left\{ \tilde{f}_{ij} +t p_{ij}: \Gamma\left(U_{ij}, \mathcal{O}_X \right)\otimes \tilde{A} \to \Gamma\left( U_{ij}, \mathcal{O}_X \right) \otimes \tilde{A} \right\}$ which corresponds to isomorphisms $\tilde{\theta}_{ij}':U_{ij}\times \textnormal{Spec}\left( \tilde{A} \right)\to U_{ij}\times \textnormal{Spec}\left(  \tilde{A} \right)$ satisfies the gluing condition, and  since $\tilde{\theta}_{ji}'^*\left( \tilde{w}_i^\alpha + t w_i'^\alpha \right)= \tilde{\theta}_{ji}^*\left( \tilde{w}_i^\alpha \right) + t w_i'^\alpha - t\mathcal{L}_{p_{ij}}\left( w_{0i}^\alpha \right)= \sum_{\beta=1}^q \tilde{h}_{ij}^{\alpha\beta} \tilde{w}_j^\beta + t \sum_{\beta=1}^q h_{0ij}^{\alpha\beta} w_j'^\beta + t \sum_{\beta=1}^q g_{ij}^{\alpha\beta} w_{0j}^\beta = \sum_{\beta=1}^q \left( \tilde{h}_{ij}^{\alpha\beta} + t g_{ij}^{\alpha\beta} \right)\left(  \tilde{w}_j^\beta + t w_j'^\beta \right)$, $\tilde{\theta}_{ij}'$ is a foliated isomorphism. On the other hand, $\left( \tilde{w}_i^1+ t w_i'^\alpha \right) \wedge \cdots \wedge \left( \tilde{w}_i^q + t w_i'^q \right) \wedge d \left( \tilde{w}_i^\alpha + t w_i'^\alpha \right)= \tilde{w}_i^1 \wedge \cdots \wedge \tilde{w}_i^q \wedge d \tilde{w}_i^\alpha + t \sum_{\beta=1}^q w_{0i}^1 \wedge \cdots {w}_i'^\beta \wedge \cdots \wedge w_{0i}^q \wedge dw_{0i}^\alpha + t w_{0i}^1 \wedge \cdots \wedge w_{0i}^q \wedge d w_i'^\alpha =0$. This defines an infinitesimal foliated deformation $\tilde{\xi}'=\left( \tilde{\mathcal{X}}', \mathcal{N}_{\tilde{\mathcal{F}}'}^* \right)$ of $\left( X, \mathcal{N}_{\mathcal{F}_0}^* \right)$ over $\tilde{A}$ which induces $\xi=\left( \mathcal{X}, \mathcal{N}_\mathcal{F}^* \right)$, so that $\left( \tilde{\mathcal{X}}', \mathcal{N}_{\tilde{\mathcal{F}}'}^* \right)$ lies in $p^{-1}(\xi)$, so that we define a map $G$ by
\begin{align*}
G: \mathbb{H}^1\left( X, \mathcal{N}_{\mathcal{F}_0}^{*\bullet} \right) \times p^{-1}(\xi)&\to p^{-1}(\xi) \\
\left(\left\{ \overline{w_i'}\right\}, \left\{ p_{ij} \right\} \right) \times \tilde{\xi} &\mapsto \tilde{\xi}'
\end{align*}
$G$ is well-defined. In fact, let $\left( \left\{ \overline{w_i''} \right\}, \left\{ p_{ij}' \right\} \right)$ represent the same cohomology class with $\left( \left\{ \overline{w_i'}\right\}, \left\{ p_{ij} \right\} \right)$. Then we have to show that $\tilde{\xi}'=G\left(  \left( \left\{ \overline{w_i'}\right\}, \left\{ p_{ij} \right\} \right), \tilde{\xi} \right)$ is equivalent to $\tilde{\xi}''= G\left(  \left( \left\{ \overline{w_i''}\right\}, \left\{ p_{ij}'\right\} \right)  , \tilde{\xi} \right)$. Then there exists $\left\{ a_i \right\} \in C^0\left( \mathcal{U}, T_X \right)$ such that $a_i- a_j= p_{ij}-  p_{ij}'$ and $\mathcal{L}_{a_i}\left( w_{0i}^\alpha \right)= w_i'^\alpha - w_i''^\alpha - \sum_{\beta=1}^q c_i^{\alpha\beta} w_{0i}^\beta$ for $c_i^{\alpha\beta} \in \Gamma\left( U_i, \mathcal{O}_X \right)$. We claim that $\textnormal{Id} - t a_i : \left(\Gamma\left( U_i, \mathcal{O}_X \right)\otimes \tilde{A}, \mathcal{N}_{\tilde{\mathcal{F}}_i'}^* \right)\to \left( \Gamma\left( U_i, \mathcal{O}_X \right)\otimes \tilde{A}, \mathcal{N}_{\tilde{\mathcal{F}}_i''}^* \right)$ is a foliated isomorphism. In fact, $\tilde{w}_i^\alpha+ t w_i'^\alpha- t \mathcal{L}_{a_i}\left( w_{0i}^\alpha \right)= \tilde{w}_i^\alpha +t w_i''^\alpha + t\sum_{\beta=1}^q c_i^{\alpha\beta} w_{0i}^\beta = \sum_{\beta=1}^q \left( \delta_\alpha^\beta + t c_i^{\alpha\beta} \right) \left( \tilde{w}_i^\beta + t w_i''^\beta \right)$. On the other hand, $\left(\tilde{f}_{ij} + t p_{ij}' \right) \circ \left( \textnormal{Id}- t a_j\right)=\left( \textnormal{Id} -t a_i \right)\circ \left(  \tilde{f}_{ij} + t p_{ij} \right)$. Hence $\left\{ \textnormal{Id}- t a_i \right\}$ glues together to give a foliated isomorphism between $\tilde{\xi}'$ and $\tilde{\xi}''$.  

Next we show that the group action $\mathbb{H}^1\left( X, \mathcal{N}_{\mathcal{F}_0}^{*\bullet} \right)\times p^{-1}(\xi)\to p^{-1}(\xi)$ is transitive. In fact, let $\tilde{\xi}\in p^{-1}(\xi)$ be a lifting of $\xi$ as above. Choose arbitrary lifting $\eta\in p^{-1}(\xi)$ of $\xi$. We have to find $\left( \left\{\overline{w_i'}\right\}, \left\{ p_{ij} \right\} \right)\in \mathbb{H}^1\left( X, \mathcal{N}_{\mathcal{F}_0}^{*\bullet} \right)$ such that $\left(\left( \left\{ \overline{w_i'}\right\}, \left\{ p_{ij} \right\} \right), \tilde{\xi} \right)$ is mapped to $\eta$ under the action $G$. Let $\eta$ consist of the data: $\tilde{f}_{ij}':\left( \Gamma\left( U_{ij}, \mathcal{O}_X \right)\otimes \tilde{A}, \mathcal{N}_{\tilde{\mathcal{F}}_j'}^*  \right)\to \left( \Gamma\left( U_{ij}, \mathcal{O}_X \right)\otimes \tilde{A}, \mathcal{N}_{\tilde{\mathcal{F}}_i'}^* \right)$ such that $\tilde{f}_{ij}'$ and $\tilde{f}_{ij}$ induce the same homomorphism in $\xi$, and $\Gamma\left( U_i, \mathcal{N}_{\tilde{\mathcal{F}}_i'}^* \right)$ is generated by $\tilde{w}_i'^\alpha\in \Gamma\left( U_i, \Omega_X^1 \right)\otimes \tilde{A}, \alpha=1,...,q$ such that $\tilde{w}_i'^\alpha$ and $\tilde{w}_i^\alpha$ induce the same element in $\Gamma\left( U_i, \Omega_X^1\right) \otimes A$. Since $\tilde{\xi}$ and $\eta$ induce $\xi$, we have $\tilde{f}_{ij}'= \tilde{f}_{ij}+ t p_{ij}$ and $\tilde{w}_i'^\alpha = \tilde{w}_i^\alpha + t B_i^\alpha$. We define $\overline{B_i}:\Gamma\left( U_i , \mathcal{N}_{\mathcal{F}_0}^* \right)\to \Gamma\left( U_i, \frac{\Omega_X^1}{\mathcal{N}_{\mathcal{F}_0}^*} \right)$ by $w_{0i}^\alpha \mapsto \overline{B_i^\alpha}$. Then as in the same argument in the proof of Proposition \ref{a11}, $\left( \left\{ \overline{B_i} \right\} , \left\{ p_{ij} \right\} \right)$ defines a $1$-cocycle of the \v Cech resolution of $\mathcal{N}_{\mathcal{F}_0}^{*\bullet}$.

Now assume that $\mathbb{H}^0\left( X, \mathcal{N}_{\mathcal{F}_0}^{*\bullet} \right)=0$. We claim that the group action $G$ is free. Assume that for given $v=\left( \left\{ \overline{w_i'} \right\}, \left\{ p_{ij} \right\} \right)\in \mathbb{H}^1\left( X, \mathcal{N}_{\mathcal{F}_0}^{*\bullet} \right)$ and $\tilde{\xi}\in p^{-1}(\xi)$, we have $G\left( v, \tilde{\xi}\right)=\tilde{\xi}$. We have to show that $v=0\in \mathbb{H}^1\left( X,  \mathcal{N}_{\mathcal{F}_0}^{* \bullet} \right)$. Let $\tilde{\xi}$ and $\tilde{\xi}'= G\left(v, \tilde{\xi} \right)$ as above. Then $\tilde{\xi}'=G\left( v , \tilde{\xi} \right)=\tilde{\xi}$ implies that we have foliated isomorphisms $g_i:\left( \Gamma\left( U_i, \mathcal{O}_X \right)\otimes \tilde{A}, \mathcal{N}_{\tilde{\mathcal{F}}_i}^* \right)\to \left( \Gamma\left( U_i, \mathcal{O}_X \right)\otimes \tilde{A}, \mathcal{N}_{\tilde{\mathcal{F}}_i'}^* \right)$ such that $g_i\circ f_{ij}= \left( f_{ij} + t p_{ij} \right) \circ g_j$. Since $\left\{ g_i \right\}$ induces a foliated automorphism on $\xi=\left( \mathcal{X}, \mathcal{N}_\mathcal{F}^* \right)$ by modulo $(t)$, which is necessarily the identity by Lemma \ref{a12}, $g_i$ is of the from $\textnormal{Id} + tq_i$ for some $q_i\in \Gamma\left( U_i, T_X \right)$. Then $\tilde{w}_i^\alpha + t \mathcal{L}_{q_i}\left( w_{0i}^\alpha \right)= \sum_{\beta=1}^q \tilde{c}_i^{\alpha\beta}\left( \tilde{w}_i^\beta + t w_i'^\beta \right)$ for some $\tilde{c}_i^{\alpha\beta}\in \Gamma\left( U_i, \mathcal{O}_X \right)\otimes \tilde{A}$ which induces $c_i^{\alpha\beta}\in \Gamma\left(  U_i, \mathcal{O}_X\right)\otimes A$. By restricting to $A$, we have $w_i^\alpha=\sum_{\beta=1}^q c_i^{\alpha\beta} w_i^\beta$, so that $c_i^{\alpha\beta}= \delta_\alpha^\beta$. Then $\tilde{c}_i^{\alpha\beta}:= \delta_\alpha^\beta + t e_i^{\alpha\beta}$ for some $e_i^{\alpha\beta}\in \Gamma\left( U_i, \mathcal{O}_X \right)$. Then we have $\tilde{w}_i^\alpha + t \mathcal{L}_{q_i}\left( w_{0i}^\alpha \right)=\sum_{\beta=1}^q\left( \delta_\alpha^\beta+ t e_i^{\alpha\beta}\right)\left( \tilde{w}_i^\beta + t w_i'^\beta \right)$, so that $\mathcal{L}_{q_i}\left( w_{0i}^\alpha \right)= w_i'^\alpha + \sum_{\beta=1}^q e_i^{\alpha\beta} w_{0i}^\beta$. On the other hand, we have $\left( \textnormal{Id}+ t q_i \right) \circ f_{ij} = \left( f_{ij}+ t p_{ij} \right) \circ \left( \textnormal{Id} + t q_j \right)$, so that $q_i- q_j = p_{ij} $. This implies that $v=0\in \mathbb{H}^1\left( X , \mathcal{N}_{\mathcal{F}_0}^{*\bullet} \right)$. This completes the proof of Proposition \ref{e7}.
\end{proof}

\begin{proposition}\label{e5}
Let $\left( X, \mathcal{N}_{\mathcal{F}_0}^* \right)$ be a foliated nonsingular projective variety with $\mathcal{N}_{\mathcal{F}_0}^*$ locally free. Then $\textnormal{FDef}_{\left( X, \mathcal{N}_{\mathcal{F}_0}^* \right)}$ has a semiuniversal formal element in the sense of \cite{Ser06} or a miniversal family in the sense of \cite{Har10}.
\end{proposition}

\begin{proof}
We check Schlessinger's criterion for $\textnormal{FDef}_{\left( X, \mathcal{N}_{\mathcal{F}_0}^* \right)}$ (see \cite{Har10} Theorem 16.2 p.113). $\textnormal{FDef}_{\left(X, \mathcal{N}_{\mathcal{F}_0}^* \right)}(k)=\left( X, \mathcal{N}_{\mathcal{F}_0}^* \right)$, so $\left(\textnormal{H}_0 \right)$ is satisfied. Let $\mathcal{U}=\left\{ U_i \right\}$ be an affine open covering such that $\Gamma\left(U_i, \mathcal{N}_{\mathcal{F}_0}^* \right)$ is generated by $w_{0i}^1,..., w_{0i}^q\in \Gamma\left( U_i, \Omega_X^1\right)$ such that $w_{0i}^\alpha= \sum_{\beta=1}^q h_{0ij}^{\alpha \beta} w_{0j}^\beta$ for some $h_{0ij}^{\alpha\beta}\in \Gamma\left( U_{ij}, \mathcal{O}_X \right)$ and $w_{0i}^1\wedge \cdots \wedge w_{0i}^q \wedge dw_{0i}^\alpha=0,\alpha=1,...,q$. Let $A'\to A$ in $\bold{Art}$ and $A''\to A$ be a small extension in $\bold{Art}$. Let $\left(\mathcal{X}, \mathcal{N}_\mathcal{F}^* \right)$ be a foliated deformation of $\left(X, \mathcal{N}_{\mathcal{F}_0}^* \right)$ over $A$, and $\left(\mathcal{X}', \mathcal{N}_{\mathcal{F}'}^*  \right)$ be a foliated deformation of $\left(X, \mathcal{N}_{\mathcal{F}_0}^* \right)$ over $A'$ which induces $\left(\mathcal{X}, \mathcal{N}_\mathcal{F}^* \right)$, and $\left( \mathcal{X}'', \mathcal{N}_{\mathcal{F}''}^* \right)$ be a foliated deformation of $\left( X, \mathcal{N}_{\mathcal{F}_0}^*\right)$ over $A''$ which induces $\left( \mathcal{X}, \mathcal{N}_\mathcal{F}^* \right)$. We may assume that $\left(\mathcal{X}, \mathcal{N}_\mathcal{F}^* \right)$ is given by the following data: we have the foliated isomorphism $\theta_i:\left( U_i \times \textnormal{Spec}(A), \mathcal{N}_{\mathcal{F}_i}^* \right)\to \left( \mathcal{X}|_{U_i}, \mathcal{N}_{\mathcal{F}}^*|_{U_i} \right)$  inducing the identity on $\left( X, \mathcal{N}_{\mathcal{F}_0}^*\right)$ such that $\Gamma\left( U_i, \mathcal{N}_{\mathcal{F}_i}^*\right)$ is generated by $w_i^\alpha\in \Gamma\left( U_i, \Omega_X^1\right)\otimes A$ which induces $w_{0i}^\alpha,\alpha=1,...,q$, and by letting $f_{ij}:\left(\Gamma\left(U_{ij}, \mathcal{O}_X \right)\otimes A, \mathcal{N}_{\mathcal{F}_j}^* \right) \to \left(\Gamma\left( U_{ij}, \mathcal{O}_X \right)\otimes A, \mathcal{N}_{\mathcal{F}_i}^* \right) $ be the ring homomorphism corresponding to $\theta_{ij}:=\theta_j^{-1}\theta_i$, we have $f_{ik}= f_{ij} f_{jk}$ and $\theta_{ji}^*\left(w_i^\alpha \right)= \sum_{\beta=1}^q h_{ij}^{\alpha\beta} w_j^\beta$ with $h_{ij}^{\alpha\beta}\in \Gamma\left(U_{ij}, \mathcal{O}_X \right)\otimes A$ which induces $h_{0ij}^{\alpha\beta}$, and $w_i^1\wedge \cdots \wedge w_i^q \wedge d w_i^\alpha =0 \in \Gamma\left( U_i, \bigwedge^{q+2} \Omega_X^1\right)\otimes A ,\alpha=1,...,q$. On the other hand, we may assume that $\left( \mathcal{X}', \mathcal{N}_{\mathcal{F}'}^*\right)$ is given by the following data: we have the foliated isomorphism $\theta_i': \left( U_i\times \textnormal{Spec}\left( A' \right), \mathcal{N}_{\mathcal{F}_i'}^*\right)\to \left( \mathcal{X}'|_{U_i}, \mathcal{N}_{\mathcal{F}'}^*|_{U_i} \right)$ which induces $\theta_i$ such that $\Gamma\left( U_i, \mathcal{N}_{\mathcal{F}_i'}^* \right)$ is generated by $w_i'^\alpha\in \Gamma\left( U_i, \Omega_X^1\right) \otimes A' , \alpha=1,...,q$ which induces $w_i^\alpha$ and by letting $f_{ij}':\left(\Gamma\left(U_{ij}, \mathcal{O}_X \right)\otimes A, \mathcal{N}_{\mathcal{F}_j'}^* \right) \to \left(\Gamma\left( U_{ij}, \mathcal{O}_X \right)\otimes A, \mathcal{N}_{\mathcal{F}_i'}^* \right) $ be the ring homomorphism corresponding to $\theta_{ij}':=\theta_j'^{-1}\theta_i'$, which induces $f_{ij}$, we have $f_{ik}'= f_{ij}' f_{jk}'$ and $\theta_{ji}'^*\left( w_i'^\alpha \right)= \sum_{\beta=1}^q h_{ij}'^{\alpha\beta} w_j'^\beta$ for some $h_{ij}'^{\alpha\beta}\in \Gamma\left( U_{ij}, \mathcal{O}_X \right) \otimes A'$ which induces $h_{ij}^{\alpha\beta}$, and $w_i'^1\wedge \cdots \wedge w_i'^q \wedge dw_i'^\alpha=0\in \Gamma\left( U_i, \bigwedge^{q+2} \Omega_X^1\right)\otimes A',\alpha=1,...,q$, and we may assume similarly for $\left(\mathcal{X}'', \mathcal{N}_{\mathcal{F}''}^* \right)$ (i.e. $\theta_i''$ induces $\theta_i$, and $w_i''^\alpha$ induces $w_i^\alpha$ and $h_{ij}''^{\alpha\beta}$ induces $h_{ij}^{\alpha\beta}$). Let $A^*=A'\times_A A''$. We define $\mathcal{X}^*$ to be the fibered sum of sheaves of rings $\mathcal{O}_{\mathcal{X}^*}= \mathcal{O}_{\mathcal{X}'}\times_{\mathcal{O}_\mathcal{X}} \mathcal{O}_{\mathcal{X}''}$ which can be described in terms of the covering $\left\{ U_i \right\}$ as follows: $\Gamma\left( U_i, \mathcal{O}_{\mathcal{X}^*} \right)=\Gamma\left(U_i, \mathcal{O}_{\mathcal{X}'} \right)\times_{\Gamma\left( U_i, \mathcal{O}_\mathcal{X} \right)} \Gamma\left( U_i, \mathcal{O}_{\mathcal{X}''} \right)\cong \left( \Gamma\left( U_i, \mathcal{O}_X \right)\otimes A'\right)\times_{\Gamma\left( U_i, \mathcal{O}_X \right)\otimes A} \left(\Gamma\left( U_i, \mathcal{O}_X \right)\otimes A'' \right)=\Gamma\left( U_i , \mathcal{O}_X \right)\otimes A^*$ and $f_{ij}^*=\left( f_{ij}', f_{ij}'' \right):\left(  \Gamma\left( U_{ij}, \mathcal{O}_X \right)\otimes A'  \right) \times_{\Gamma\left( U_{ij} , \mathcal{O}_X \right)\otimes A} \left( \Gamma\left( U_{ij}, \mathcal{O}_X \right) \otimes A'' \right) \to \left(  \Gamma\left( U_{ij}, \mathcal{O}_X \right)\otimes A'  \right) \times_{\Gamma\left( U_{ij} , \mathcal{O}_X \right)\otimes A} \left( \Gamma\left( U_{ij}, \mathcal{O}_X \right) \otimes A'' \right)$ induces $\theta_{ij}^*:U_{ij}\times \textnormal{Spec}\left( A^* \right) \to U_{ij}\times \textnormal{Spec}\left(A^*\right)$ with $f_{ik}^*=f_{ij}^* f_{jk}^*$. Then $\mathcal{X}$ is flat over $A^*$. On the other hand, the fibered sum $\mathcal{N}_{\mathcal{F}_i^*}^* = \mathcal{N}_{\mathcal{F}_i'}^* \times_{\mathcal{N}_{\mathcal{F}_i}^*} \mathcal{N}_{\mathcal{F}_i''}^*$ is generated by $w_i^{* \alpha}:=\left( w_i'^\alpha, w_i''^\alpha \right)\in \left(\Gamma\left( U_i, \Omega_X^1 \right)\otimes A' \right) \times_{\Gamma\left(U_i, \Omega_X^1 \right)\otimes A}  \left( \Gamma\left( U_i, \Omega_X^1 \right)\otimes A'' \right)=\Gamma\left( U_i, \Omega_X^1 \right)\otimes A^*$ and by setting $h_{ij}^{*\alpha\beta}=\left( h_{ij}'^{\alpha\beta}, h_{ij}''^{\alpha\beta} \right)$, we have
\begin{align*}
\left( \theta_{ji}^*\right)^*\left( w_i^{*\alpha} \right)= \left(\theta_{ji}'^*\left( w_i'^\alpha \right), \theta_{ji}''^* \left( w_i''^\alpha \right) \right) = \left( \sum_{\beta=1}^q  h_{ij}'^{\alpha\beta} w_j'^\beta, \sum_{\beta=1}^q h_{ij}''^{\alpha\beta} w_j''^\beta \right) = \sum_{\beta=1}^q \left(h_{ij}'^{\alpha\beta}, h_{ij}''^{\alpha\beta} \right)\left( w_j'^\beta, w_j''^\beta \right)= \sum_{\beta=1}^q h_{ij}^{*\alpha\beta} w_j^{*\beta}
\end{align*}
and the integrability condition holds:
\begin{align*}
&\left(w_i'^1,w_i''^1 \right) \wedge \cdots \wedge \left( w_i'^q , w_i''^q \right) \wedge d \left( w_i'^\alpha, w_i''^\alpha \right)= \left(w_i'^1 \wedge \cdots \wedge w_i'^q \wedge dw_i'^\alpha, w_i''^1 \wedge \cdots \wedge w_i''^q \wedge dw_i''^\alpha    \right)= 0 \\
& \in \left( \Gamma\left( U_i, \bigwedge^{q+2}\Omega_X^1\right)\otimes A' \right)\times_{ \Gamma\left( U_i, \bigwedge^{q+2}\Omega_X^1\right)\otimes A} \left(  \Gamma\left( U_i, \bigwedge^{q+2}\Omega_X^1\right)\otimes A'' \right) = \Gamma\left( U_i, \bigwedge^{q+2} \Omega_X^1 \right)\otimes A^*
\end{align*}
Then $\left\{ \mathcal{N}_{\mathcal{F}_i^*}^* \right\}$ defines a locally free subsheaf $\mathcal{N}_{\mathcal{F}^*}^*$ of $\Omega_{\mathcal{X}^*/A^*}^1$. Since $\mathcal{N}_{\mathcal{F}^*}^*\hookrightarrow \Omega_{\mathcal{X}^*/A^*}^1$ induces an injection $\mathcal{N}_{\mathcal{F}_0}^* \hookrightarrow \Omega_{X}^1$, we see that $\frac{\Omega_{\mathcal{X}^*/A^*}}{\mathcal{N}_{\mathcal{F}^*}^*}$ is flat over $A^*$. This implies that $\left(\textnormal{H}_1 \right)$ is satisfied. $\left( \textnormal{H}_3 \right)$ follows from Proposition \ref{nb1}. It remains to show $\left( \textnormal{H}_2 \right)$. Let $A=k$ and $A''=k[\epsilon]$ and $\left( \mathcal{X}^*, \mathcal{N}_{\mathcal{F}^*}^* \right)$ as above. Let $\left(\mathcal{Y}, \mathcal{N}_\mathcal{G}^* \right)$ be another foliated deformation of $\left(X, \mathcal{N}_{\mathcal{F}_0}^* \right)$ over $A^*$ which induces $\left( \mathcal{X}, \mathcal{N}_{\mathcal{F}'}^*\right)$ and $\left( \mathcal{X}, \mathcal{N}_{\mathcal{F}''}^* \right)$. We may assume that the underlying scheme $\mathcal{Y}=\mathcal{X}^*$ (see \cite{Ser06} p.66-67). We note that we have an isomorphism $\theta_i^*: U_i \times \textnormal{Spec}\left( A^*\right) \to  \mathcal{X}^*|_{U_i}$ such that $\theta_j^{* -1}\theta_i^*=\theta_{ij}^*$. Let $\mathcal{N}_{\mathcal{G}_i}^*=\theta_i^*\left( \mathcal{N}_\mathcal{G}^*|_{U_i} \right)$ be the induced foliation from $\mathcal{N}_\mathcal{G}^*$. Then $\Gamma\left(U_i, \mathcal{N}_{\mathcal{G}_i}^*\right)$ induces $\Gamma\left( U_i, \mathcal{N}_{\mathcal{F}_i'}^* \right)$ and $\Gamma\left(U_i, \mathcal{N}_{\mathcal{F}_i''}^* \right)$, i.e. there is a commutative diagram. 
\begin{center}
$\begin{CD}
\Gamma\left( U_i , \mathcal{N}_{\mathcal{G}_i}^* \right)\subset \Gamma\left( U_i, \Omega_X^1 \right)\otimes A^* @>>> \Gamma\left( U_i, \mathcal{N}_{\mathcal{F}_i''}^* \right)\subset \Gamma\left( U_i, \Omega_X^1 \right)\otimes k[\epsilon]\\
@VVV @VVV\\
\Gamma\left( U_i, \mathcal{N}_{\mathcal{F}_i'}^*\right)\otimes A \subset \Gamma\left( U_i, \Omega_X^1 \right)\otimes A' @>>> \Gamma\left( U_i, \mathcal{N}_{\mathcal{F}_0}^* \right)
\end{CD}$
\end{center}
This induces $\Gamma\left( U_i, \mathcal{N}_{\mathcal{G}_i}^* \right)\to \Gamma\left( U_i, \mathcal{N}_{\mathcal{F}_i^*}^* \right)=\Gamma\left( U_i, \mathcal{N}_{\mathcal{F}_i'}^* \right)\times_{\Gamma\left( U_i, \mathcal{N}_{\mathcal{F}_0}^*\right)}\left( U_i, \mathcal{N}_{\mathcal{F}_i''}^* \right)$ which is an isomorphism since $\mathcal{N}_{\mathcal{F}_i^*}^*$ is free, and $\mathcal{N}_{\mathcal{G}_i}^*$ is flat over $A^*$ and induces the identity $\Gamma\left( U_i, \mathcal{N}_{\mathcal{F}_0}^* \right)$. This implies that $\Gamma\left( U_i, \mathcal{N}_{\mathcal{G}_i}^*\right)= \Gamma\left( U_i, \mathcal{N}_{\mathcal{F}_i^*}^* \right)$ in $\Gamma\left( U_i, \Omega_X^1 \right)\otimes A^*$. Hence $\mathcal{N}_\mathcal{G}^*= \mathcal{N}_{\mathcal{F}^*}^*$ in $\Omega_{\mathcal{X}^* /A^*}^1$. This completes the proof of Proposition \ref{e5}.
\end{proof}

\begin{theorem}
Let $\left( X ,  \mathcal{N}_{\mathcal{F}_0}^* \right)$ be a foliated nonsingular projective variety with $\mathcal{N}_{\mathcal{F}_0}^*$ locally free, and $\mathbb{H}^0\left( X, \mathcal{N}_{\mathcal{F}_0}^{*\bullet} \right)=0$. Then the functor $\textnormal{FDef}_{\left( X, \mathcal{N}_{\mathcal{F}_0}^* \right)}$ is pro-representable.
\end{theorem}

\begin{proof}
We check Schlessinger's criterion for $\textnormal{FDef}_{\left( X, \mathcal{N}_{\mathcal{F}_0}^* \right)}$ (see \cite{Har10} Theorem 16.2 p.113). $(\textnormal{H}_0)-(\textnormal{H}_3)$ follow from Proposition \ref{e5}, and $(\textnormal{H}_4)$ follows from Proposition \ref{e7}.
\end{proof}

\subsection{Stability of a foliated nonsingular variety $\left( X, \mathcal{N}_{\mathcal{F}_0}^* \right)$ with $\mathcal{N}_{\mathcal{F}_0}^*$ locally free }\

\begin{definition}
A foliated nonsingular variety $\left(X, \mathcal{N}_{\mathcal{F}_0}^* \right)$ with $\mathcal{N}_{\mathcal{F}_0}^*$ is called stable if the morphism of functors from $(\ref{a15})$ by forgetting foliated structures
\begin{align*}
\Phi: \textnormal{FDef}_{\left( X, \mathcal{N}_{\mathcal{F}_0}^* \right)}\to \textnormal{Def}_X
\end{align*}
is smooth.
\end{definition}

We note that $\textnormal{FDef}_{\left(X, \mathcal{N}_{\mathcal{F}_0}^* \right)}$ has a semiuniversal formal element from Proposition \ref{e5} and obstruction space $\mathbb{H}^2\left(X, \mathcal{N}_{\mathcal{F}_0}^{*\bullet} \right)$ from Proposition \ref{nb1}. On the other hand, $\textnormal{Def}_{X}$ has a semiuniversal formal element and obstruction space $H^2\left( X, T_X\right)$ (see \cite{Ser06}). Assume $\mathbb{H}^1\left(X, \mathscr{H}om_{\mathcal{O}_X}\left( \mathcal{N}_{\mathcal{F}_0}^*, \frac{\Omega_X^1}{\mathcal{N}_{\mathcal{F}_0}^*} \right)^\bullet \right)=0$. Then from the long exact sequence of $(\ref{a18})$, we have the surjection $d\Phi: \mathbb{H}^1\left(X, \mathcal{N}_{\mathcal{F}_0}^{* \bullet }\right)\to H^1\left( X, T_X \right)$ and an injective obstruction map $\mathbb{H}^2\left( X, \mathcal{N}_{\mathcal{F}_0}^{ * \bullet } \right)\to H^2\left( X, T_X \right)$, in other words, $\textnormal{FDef}_{\left( X, \mathcal{N}_{\mathcal{F}_0}^* \right)}$ is less obstructed than $\textnormal{Def}_X$ (see \cite{Ser06} Definition 2.3.5 p.58). This implies that $\Phi$ is smooth. (see \cite{Ser06} Proposition 2.3.6 p.59). We summarize as follows.

\begin{proposition}
Let $\left( X, \mathcal{N}_{\mathcal{F}_0}^* \right)$ be a foliated nonsingular variety with $\mathcal{N}_{\mathcal{F}_0}^*$ locally free. Assume that $\mathbb{H}^1\left(X, \mathscr{H}om_{\mathcal{O}_X}\left( \mathcal{N}_{\mathcal{F}_0}^*, \frac{\Omega_X^1}{\mathcal{N}_{\mathcal{F}_0}^*} \right)^\bullet \right)=0$. Then $\Phi$ is smooth. In other words, for a flat deformation $\mathcal{X}$ of $X$ over $A$, there exists a foliated deformation $\left(\mathcal{X}, \mathcal{N}_\mathcal{F}^* \right)$ over $A$ which induces $\mathcal{X}$.
\end{proposition}

\subsection{Formal foliated deformations of $\left( X, \mathcal{N}_{\mathcal{F}_0}^* \right)$ with $\mathcal{N}_{\mathcal{F}_0}^*$ locally free} \

We only touch upon formal foliated deformations of a foliated nonsingular variety $\left( X, \mathcal{N}_{\mathcal{F}_0}^* \right)$ with $\mathcal{N}_{\mathcal{F}_0}^*$ locally free in a relation to Part I. The following Theorem \ref{a16} can be considered as the scheme-theoretic version of Theorem of existence for deformations of foliated complex analytic structures in terms of cotangent sheaves in Part I.

\begin{theorem}[compare \cite{Ser06} Corollary 2.6.4 p.92]\label{a16}
Let $\left( X, \mathcal{N}_{\mathcal{F}_0}^* \right)$ be a foliated nonsingular projective variety with $\mathcal{N}_{\mathcal{F}_0}^*$ locally free. If $\mathbb{H}^0\left( X , \mathcal{N}_{\mathcal{F}_0}^{*\bullet} \right)=\mathbb{H}^2\left( X, \mathcal{N}_{\mathcal{F}_0}^{*\bullet} \right)=0$, then $\textnormal{FDef}_{\left( X, \mathcal{N}_{\mathcal{F}_0}^* \right)}$ is pro-represented by a formal power series ring.
\end{theorem}

\begin{proof}
By Theorem \ref{e5}, $\textnormal{FDef}_{\left(X, \mathcal{N}_{\mathcal{F}_0}^* \right)}$ is pro-representable. Let $\textnormal{FDef}_{\left( X, \mathcal{N}_{\mathcal{F}_0}^* \right)}$ be pro-represented by a universal formal couple (see \cite{Ser06} p.47). By Proposition \ref{nb1}, $\textnormal{FDef}_{\left(X, \mathcal{N}_{\mathcal{F}_0}^* \right)}$ has an obstruction space $\mathbb{H}^2\left( X , \mathcal{N}_{\mathcal{F}_0}^{*\bullet} \right)$ (see \cite{Ser06} Definition 2.2.9 p.53).  By Proposition \ref{nb1}, we have $\dim_k \mathbb{H}^1\left( X, \mathcal{N}_{\mathcal{F}_0}^{*\bullet} \right) \geq \mu(R) \geq \dim_k \mathbb{H}^1\left( X, \mathcal{N}_{\mathcal{F}_0}^{*\bullet} \right)- \dim_k \mathbb{H}^2\left( X, \mathcal{N}_{\mathcal{F}_0}^{*\bullet} \right)$, where $\mu(R)$ is the Krull dimension of $R$ (see \cite{Ser06} Corollarly 2.2.11 p.54). By the hypothesis $\mathbb{H}^2\left( X, \mathcal{N}_{\mathcal{F}_0}^{*\bullet} \right)=0$, we see that $\textnormal{FDef}_{\left( X , \mathcal{N}_{\mathcal{F}_0}^* \right)}$ is smooth, and we have $\mu(R)=\dim_k \mathbb{H}^1\left(X, \mathcal{N}_{\mathcal{F}_0}^{*\bullet} \right)$. Hence we have $R\cong k[[X_1,..., X_m]]$.
\end{proof}

We can consider a concept of formal foliated deformations of a foliated nonsingular variety $\left( X, \mathcal{N}_{\mathcal{F}_0}^* \right)$ with $\mathcal{N}_{\mathcal{F}_0}^*$ locally free as an analogue of formal deformations of a nonsingular variety $X$ as in \cite{Ser06} p.77. We note that we can extend the functor $\textnormal{FDef}_{\left(X, \mathcal{N}_{\mathcal{F}_0}^* \right)}$ to the functor
\begin{align*}
\widehat{\textnormal{FDef}}_{\left( X, \mathcal{N}_{\mathcal{F}_0}^* \right)} : \widehat{\bold{Art}}&\to \left(\textnormal{Sets} \right)\\
     \bar{A} &\to \varprojlim \textnormal{FDef}_{\left(X, \mathcal{N}_{\mathcal{F}_0}^* \right)}\left(\bar{A}/\mathfrak{m}_{\bar{A}}^{n+1} \right)
\end{align*}
(see \cite{Ser06} p.46). An element $\hat{\eta}\in \widehat{\textnormal{FDef}}_{\left( X, \mathcal{N}_{\mathcal{F}_0}^* \right)}\left(\bar{A} \right)$ is called a formal element of $\textnormal{FDef}_{\left( X, \mathcal{N}_{\mathcal{F}_0}^* \right)}$. By definition $\hat{\eta}$ can be represented as a sysemt of elements $\left\{ \eta_n \in \textnormal{FDef}_{\left( X, \mathcal{N}_{\mathcal{F}_0}^* \right)}\left(\bar{A}/\mathfrak{m}_{\bar{A}}^{n+1} \right) \right\}_{n\geq 0}$  such that for every $n\geq 1$ the map $\textnormal{FDef}_{\left( X, \mathcal{N}_{\mathcal{F}_0}^* \right)}\left( \bar{A}/\mathfrak{m}_{\bar{A}}^{n+1}\right)\to \textnormal{FDef}_{\left( X, \mathcal{N}_{\mathcal{F}_0}^* \right)}\left( \bar{A}/\mathfrak{m}_{\bar{A}}^n \right)$ induced by the projection $\bar{A}/\mathfrak{m}_{\bar{A}}^{n+1}\to \bar{A}/{\mathfrak{m}}_{\bar{A}}^n$ sends $\eta_n\mapsto \eta_{n-1}$. A formal foliated deformation of $\left( X, \mathcal{N}_{\mathcal{F}_0}^* \right)$ over $\bar{A}$ is a formal element $\hat{\eta}=\left\{\eta_n \right\} \in \widehat{\textnormal{FDef}}_{\left( X, \mathcal{N}_{\mathcal{F}_0}^* \right)}\left( \bar{A}\right)$ for $\textnormal{FDef}_{\left(X, \mathcal{N}_{\mathcal{F}_0}^* \right)}$, where $\eta_n$ is an infinitesimal foliated deformation of $\left( X, \mathcal{N}_{\mathcal{F}_0}^* \right)$ over $\bar{A}_n$
\begin{center}
$\eta_n:\begin{CD}
\left(X, \mathcal{N}_{\mathcal{F}_0}^* \right) @>>> \left( \mathcal{X}_n, \mathcal{N}_{\mathcal{F}_n}^* \right) \\
@VVV @VV \pi_n V\\
\textnormal{Spec}(k) @>>> \textnormal{Spec}\left( \bar{A}_n \right)
\end{CD}$
\end{center}
where $\bar{A}_n=\bar{A}/\mathfrak{m}_{\bar{A}}^{n+1}$, such that for all $n\geq 1$, $\eta_n$ induces $\eta_{n-1}$ by pullback under the natural inclusion $\textnormal{Spec}\left(\bar{A}_{n-1} \right)\to \textnormal{Spec}\left( \bar{A}_n \right)$, i.e. we have
\begin{center}
$\eta_{n-1}:\begin{CD}
\left(X, \mathcal{N}_{\mathcal{F}_0}^* \right) @>>> \left( \mathcal{X}_{n-1}=\mathcal{X}_n \otimes_{\bar{A}_n} \bar{A}_{n-1}, \mathcal{N}_{\mathcal{F}_n}^*\otimes_{\bar{A}_n} \bar{A}_{n-1}  \right)\\
@VVV @VV\pi_{n-1}V\\
\textnormal{Spec}(k) @>>> \textnormal{Spec}\left(\bar{A}_{n-1} \right)
\end{CD}$
\end{center}
We will denote such a formal foliated deformation by $\left( \bar{A}, \left\{ \eta_n \right\} \right)$.

Let $\left(X, \mathcal{N}_{\mathcal{F}_0}^*\right)$ be a foliated nonsingular projective variety with $\mathcal{N}_{\mathcal{F}_0}^*$ locally free and consider a flat family of foliated deformations of $\left( X, \mathcal{N}_{\mathcal{F}_0}^* \right)$ parametrized by an affine scheme $S=\textnormal{Spec}\left( B \right)$ where $B$ is a $k$-algebras of finite type:
\begin{center}
$\eta:\begin{CD}
\left( X, \mathcal{N}_{\mathcal{F}_0}^* \right) @>f>> \left( \mathcal{X}, \mathcal{N}_\mathcal{F}^* \right)\\
@VVV @VV\pi V\\
\textnormal{Spec}(k) @>s >> S
\end{CD}$
\end{center}
where $\pi$ is projective, and $\mathcal{X}$ is smooth over $S$, $s$ is a $k$-rational point and $\mathcal{N}_\mathcal{F}^*$ is a locally free subsheaf of $\Omega_{\mathcal{X}/S}^1$ such that $\frac{\Omega_{\mathcal{X}/S}^1}{\mathcal{N}_\mathcal{F}^*}$ is flat over $S$ and $d\left( \mathcal{N}_{\mathcal{F}, x} ^*  \right) \subset \mathcal{N}_{\mathcal{F} , x} ^* \bigwedge \Omega_{\mathcal{X}/S , x }^1$ for $x\in \mathcal{X}-\textnormal{Sing}\left( \mathcal{F}\right)$, where $\textnormal{Sing}\left( \mathcal{F}\right)$ is the locus of points of $\mathcal{X}$ at which $\frac{\Omega_{\mathcal{X}/S}^1}{\mathcal{N}_\mathcal{F}^*}$ is not locally free. We will denote the foliated deformation $\eta$ by $\left( S, s , \eta \right)$. Given such a foliated deformation $\left(  S, s, \eta \right)$ let $\eta_n$ be the infinitesimal foliated deformation induced pulling back of $\eta$ under the closed embedding $\textnormal{Spec}\left( \mathcal{O}_{S,s}/\mathfrak{m}^{n+1} \right) \to S$. Then $\left( \widehat{\mathcal{O}}_{S,s}, \left\{ \eta_n \right\} \right)$ is a formal foliated deformation of $\left(X, \mathcal{N}_{\mathcal{F}_0}^* \right)$ which we call the formal foliated deformation defined by $\eta$.

\begin{definition}
A foliated deformation $\left(S, s, \eta \right)$ of a foliated nonsingular variety $\left( X, \mathcal{N}_{\mathcal{F}_0}^* \right)$ with $\mathcal{N}_{\mathcal{F}_0}^*$ locally free is called formally universal $($formally semiuniversal, formally versal$)$ if the formal foliated deformation $\left( \widehat{\mathcal{O}}_{S,s}, \left\{ \eta_n\right\}_{n \geq 0} \right)$ defined by $\eta$ is universal $($resp. semiuniversal, versal$)$ \textnormal{(see \cite{Ser06} Definition 2.2.3 p.47 and Definition 2.2.6 p.49)}. A flat foliated family $\pi:\left(\mathcal{X}, \mathcal{N}_\mathcal{F}^* \right)\to S$, where $\mathcal{X}$ is smooth over $S$, and $\mathcal{N}_{\mathcal{F}}^*$ is a locally free subsehaf of $\Omega_{\mathcal{X}/S}^1$, is called formally universal $($resp. formally semiuniversal, formally versal$)$ at a $k$-rational point $s\in S$ if
\begin{center}
$\eta:\begin{CD}
\left(\mathcal{X}(s), \mathcal{N}_{\mathcal{F}(s)}^*\right) @>\hookrightarrow>> \left( \mathcal{X}, \mathcal{N}_\mathcal{F}^* \right)\\
@VVV @VV\pi V \\
\textnormal{Spec}(k) @>s>> S
\end{CD}$
\end{center}
is a formally universal $($resp. formally semiuniversal, formally versal$)$ foliated deformation of $\left( \mathcal{X}(s), \mathcal{N}_{\mathcal{F}(s)}^* \right)$.
\end{definition}

\begin{definition}[foliated Kodaira-Spencer map in an algebraic foliated family]
Let $\left( X, \mathcal{N}_{\mathcal{F}_0}^* \right)$ be a foliated nonsingular variety with $\mathcal{N}_{\mathcal{F}_0}^*$ locally free. Let us consider a foliated deformation $(S, s, \eta)$
\begin{center}
$\eta: \begin{CD}
\left( X, \mathcal{N}_{\mathcal{F}_0}^* \right) @>>> \left( \mathcal{X}, \mathcal{N}_\mathcal{F}^* \right)\\
@VVV @VV\pi V\\
\textnormal{Spec}(k) @>s>> S
\end{CD}$
\end{center}
By pulling back $(S, s, \eta)$ by morphisms $\textnormal{Spec}\left( k[\epsilon] \right)\to S$ with image $s$ and applying \textnormal{Proposition \ref{at25}} we define a linear map \textnormal{(see \cite{MO15} p.200)}
\begin{align*}
\kappa_{\pi, s}: T_{S,s}\to  \textnormal{FDef}_{\left( X, \mathcal{N}_{\mathcal{F}_0}^* \right)} \cong \mathbb{H}^1\left( X, \mathcal{N}_{\mathcal{F}_0}^{* \bullet} \right)
\end{align*}
which is called the foliated Kodaira-Spencer map of $(S,s, \eta)$ in terms of cotangent sheaves.
\end{definition}

\begin{proposition}[compare \cite{Ser06} Proposition 2.5.8 p.79]\label{a17}
Let $\left(S , s , \eta \right)$ be a foliated deformation of a foliated nonsingular variety $\left( X, \mathcal{N}_{\mathcal{F}_0}^* \right)$ with $\mathcal{N}_{\mathcal{F}_0}^*$ locally free.
\begin{enumerate}
\item If $\eta$ is formally versal $($resp. formally semiuniversal or formally universal$)$, then the foliated Kodaira-Spencer map $\kappa_{\pi, s}: T_s S\to \textnormal{FDef}_{\left(X, \mathcal{N}_{\mathcal{F}_0}^* \right)}\left( k[\epsilon]\right)\cong \mathbb{H}^1 \left(X, \mathcal{N}_{\mathcal{F}_0}^{* \bullet} \right)$ is surjective $($resp. an isomorphism$)$.
\item If $S$ is nonsingular at $s$ and the foliated Kodaira-Spencer map $\kappa_{\pi,s}$ is surjective $($resp. an isomorphism$)$, then $\eta$ is formally versal $($resp. formally semiuniversal$)$ and $\left(X, \mathcal{N}_{\mathcal{F}_0}^* \right)$ is unobstructed, i.e. $\textnormal{FDef}_{\left( X, \mathcal{N}_{\mathcal{F}_0}^* \right)}$ is smooth.
\end{enumerate}
\end{proposition}
\begin{proof}
Let $h_{\widehat{\mathcal{O}}_{S,s}}: \widehat{\bold{Art}}\to \left( \textnormal{Sets}\right)$ be a functor defined by $\bar{A} \to \textnormal{Hom}_{\widehat{\bold{Art}}}\left(\bar{A},  \widehat{\mathcal{O}}_{S,s} \right) $. Then there is a $1-1$ correspondence between $\widehat{\textnormal{FDef}}_{\left(X, \mathcal{N}_{\mathcal{F}_0}^* \right)}\left(\widehat{\mathcal{O}}_{S, s} \right)$ and the set of morphisms of functors $h_{\widehat{\mathcal{O}}_{S,s}}\to \textnormal{FDef}_{\left(X, \mathcal{N}_{\mathcal{F}_0}^* \right)}$ (see \cite{Ser06} Lemma 2.2.2 p.46). So $\left(\widehat{\mathcal{O}}_{S,s}\, \left\{\eta_n \right\} \right)$ defines a morphism $h_{\widehat{\mathcal{O}}_{S,s}}\to \textnormal{FDef}_{\left(X, T_{\mathcal{F}_0} \right)}$. Then $(1)$ follows from the definitions of versality and semiuniversality of the formal couple $\left(\widehat{\mathcal{O}}_{S,s}, \left\{\eta_n \right\} \right)$ applied to $h_{\widehat{\mathcal{O}}_{S,s}}\to \textnormal{FDef}_{\left(X, T_{\mathcal{F}_0} \right)}$. Let us prove $(2)$. Since $h_{\widehat{\mathcal{O}}_{S,s}}$ is smooth, i.e unobstructed, and $\textnormal{FDef}_{\left(X, \mathcal{N}_{\mathcal{F}_0}^* \right)}$ has a semiuniversal formal element by Proposition \ref{e5} and $\textnormal{FDef}_{\left( X, \mathcal{N}_{\mathcal{F}_0}^* \right)}$ has an obstruction space $\mathbb{H}^2\left( X, \mathcal{N}_{\mathcal{F}_0}^{* \bullet} \right)$ by Proposition \ref{nb1}, and the foliated Kodaira-Spencer map $\kappa_{S,s}$ is surjective by the assumption, the morphism $h_{\widehat{\mathcal{O}}_{S,s}}\to \textnormal{FDef}_{\left( X, \mathcal{N}_{\mathcal{F}_0}^* \right)}$ is smooth and $\textnormal{FDef}_{\left( X , \mathcal{N}_{\mathcal{F}_0}^* \right)}$ is unobstructed. (see \cite{Ser06} Corollary 2.3.7 p.59).
\end{proof}

Proposition \ref{a17} is the algebraic version of Theorem of completeness for deformations of foliated complex analytic structures in terms of cotangent sheaves in Part I. We finish this section by mentioning the concept of algebraizability in foliated deformations.

\begin{definition}[compare \cite{Ser06} p.80]
A formal foliated deformation $\left( \bar{A}, \left\{ \eta_n \right\} \right)$ of a foliated nonsingular variety $\left( X, \mathcal{N}_{\mathcal{F}_0}^* \right)$ with $\mathcal{N}_{\mathcal{F}_0}^*$ locally free is called algebraizable if there exists an algebraic foliated deformation $\left( S, s, \xi \right)$ of $\left( X, \mathcal{N}_{\mathcal{F}_0}^* \right)$ and an isomorphism $\widehat{\mathcal{O}}_{S,s}\cong \bar{A}$ sending $\eta_n\to \xi_n$ for all $n$, i.e. $\left( \bar{A}, \left\{ \eta_n \right\} \right)$ is isomorphic to the formal foliated deformation defined by $\xi$. The foliated deformation $\left( S, s, \xi \right)$ is called an algebraization of $\left( \bar{A} , \left\{ \eta_n \right\}  \right)$.
\end{definition}

The existence of algebraizations can be considered as the counterpart of the convergence step in the construction of local families of foliated deformations of compact foliated complex manifolds in terms of cotangent sheaves in Part I.

\section{Generalities on foliated deformations in terms of both tangent and cotangent sheaves}

Let $\left(X, T_{\mathcal{F}_0}, \mathcal{N}_{\mathcal{F}_0}^* \right)$ be a foliated nonsingular variety. In other words, we have exact sequences
\begin{align*}
0\to T_{\mathcal{F}_0} &\to T_X \to T_X/T_{\mathcal{F}_0} \to 0\\
0\to \mathcal{N}_{\mathcal{F}_0}^* &\to \Omega_X^1 \to \Omega_X^1/\mathcal{N}_{\mathcal{F}_0}^* \to 0
\end{align*}
where $T_{\mathcal{F}_0}$ is a coherent subsheaf of $T_X$ and $\mathcal{N}_{\mathcal{F}_0}^*$ is a coherent subsheaf of $\Omega_X^1$ such that $\left( \frac{T_X}{T_{\mathcal{F}_0}} \right)^*= \mathcal{N}_{\mathcal{F}_0}^*$ and $\left(\frac{\Omega_X^1}{\mathcal{N}_{\mathcal{F}_0}^*} \right)^*= T_{\mathcal{F}_0}$ and $\left[ T_{\mathcal{F}_0}, T_{\mathcal{F}_0}\right] \subset T_{\mathcal{F}_0}$ and $d\left(\mathcal{N}_{\mathcal{F}_0,x}^* \right) \subset \mathcal{N}_{\mathcal{F}_0,x}^* \bigwedge \Omega_X^1 $ for $x\in X- \textnormal{Sing}\left( \mathcal{F}_0\right)$, and $T_X/T_{\mathcal{F}_0}$ is a torsion-free $\mathcal{O}_X$-module and $\Omega_X^1/\mathcal{N}_{\mathcal{F}_0}^*$ is a torsion-free $\mathcal{O}_X$-module. A cartesian diagram of morphism of algebraic schemes
\begin{center}
$\begin{CD}
\left( X, T_{\mathcal{F}_0}, \mathcal{N}_{\mathcal{F}_0}^*\right) @>>> \left( \mathcal{X}, T_\mathcal{F}, \mathcal{N}_\mathcal{F}^* \right)\\
@VVV @VVV\\
\textnormal{Spec}(k) @>>> S
\end{CD}$
\end{center}
where $\pi$ is flat, smooth and surjective, and $S$ is a connected algebraic scheme, is called a family of foliated deformations of $\left( X, T_{\mathcal{F}_0}, \mathcal{N}_{\mathcal{F}_0}^* \right)$ in terms of both tangent and cotangent sheaves parameterized by $S$ if there exist exact sequences
\begin{align}
0\to T_{\mathcal{F}} &\to T_{\mathcal{X}/S} \to T_{\mathcal{X}/S}/T_{\mathcal{F}}\to 0\label{c1}\\
0 \to \mathcal{N}_{\mathcal{F}}^* & \to \Omega_{\mathcal{X}/S}^1 \to \Omega_{\mathcal{X}/S}^1/\mathcal{N}_{\mathcal{F}}^* \to 0 \label{c2}
\end{align}
such that $\left( \frac{T_{\mathcal{X}/S}}{T_\mathcal{F}} \right)^*=\mathcal{N}_\mathcal{F}^*$ and $\left( \frac{\Omega_{\mathcal{X}/S}^1}{ \mathcal{N}_\mathcal{F}^* } \right)^*= T_\mathcal{F}^*$, and  $\left[ T_\mathcal{F}, T_\mathcal{F} \right] \subset T_{\mathcal{F}}$ and $T_{\mathcal{X}/S}/T_{\mathcal{F}}$ is flat over $S$, and $d_{\mathcal{X}/S}\left( \mathcal{N}_{\mathcal{F},x}^*\right) \subset \mathcal{N}_{\mathcal{F},x}^* \bigwedge \Omega_{\mathcal{X}/S,x}^1$ for $x\in \mathcal{X} - \textnormal{Sing}\left( \mathcal{F} \right)$ and $\Omega_{\mathcal{X}/S}^1/\mathcal{N}_\mathcal{F}^*$ is flat over $S$, and for each $k$-rational point $t\in S$ the exact sequence $(\ref{c1})$ and $(\ref{c2})$ induce exact sequences on the fiber $X_t$
\begin{align*}
0\to T_{\mathcal{F}_t} &\to T_{X_t} \to T_{X_t}/T_{\mathcal{F}_t} \to 0\\
0 \to \mathcal{N}_{\mathcal{F}_t}^* & \to \Omega_{X_t}^1 \to \Omega_{X_t}^1/\mathcal{N}_{\mathcal{F}_t}^* \to 0
\end{align*}
such that $\left( \frac{ T_{X_t} }{T_{\mathcal{F}_t}} \right)^* = \mathcal{N}_{\mathcal{F}_t}^*$ and $\left( \frac{\Omega_{X_t}^1}{\mathcal{N}_{\mathcal{F}_t}^*} \right)^*= T_{\mathcal{F}_t}$, and $T_{X_t}/T_{\mathcal{F}_t}$ is a torsion-free $\mathcal{O}_{X_t}$-module and $\Omega_{X_t}^1/ \mathcal{N}_{\mathcal{F}_t}^*$ is a torsion-free $\mathcal{O}_{X_t}$-module. When $S=\textnormal{Spec}(A)$ with $A\in \textnormal{ob}\left(\bold{Art}^*\right)$ and $s\in S$ is the closed point we have a local family of foliated deformations of $\left(X, T_{\mathcal{F}}\right)$ in terms of tangent sheaves over $A$. The foliated deformation $\eta$ in terms of both tangent and cotangent sheaves will be denoted by $\left( S, \eta\right)$ or $\left( A, \eta\right)$ when $S=\textnormal{Spec}(A)$. The local foliated deformation $\left( A, \eta\right)$ is infinitesimal (resp. first-order) if $A\in \bold{Art}$ (resp. $A=k[\epsilon]$). Given another foliated deformation of $\left( X, T_{\mathcal{F}_0}, \mathcal{N}_{\mathcal{F}_0}^*\right)$ over $S$ in terms of both tangent and cotangent sheaves
\begin{center}
$\xi: \begin{CD}
\left( X, T_{\mathcal{F}_0}, \mathcal{N}_{\mathcal{F}_0}^*\right) @>>>  \left( \mathcal{Y}, T_\mathcal{G} , \mathcal{N}_{\mathcal{G}}^* \right)\\
@VVV @VVV\\
\textnormal{Spec}(k) @>>> S
\end{CD}$
\end{center}
An isomorphism between foliated deformations $\eta$ and $\xi$ is an foliated $S$-isomorphism $\phi:\left(\mathcal{X}, T_{\mathcal{F}} \right)\to \left(  \mathcal{Y}, T_{\mathcal{G}} \right)$ with $T_\mathcal{F}=\phi^* T_{\mathcal{G}}$ as a subsheaf of $T_{\mathcal{X}/S}$ inducing the identity on $\left( X, T_{\mathcal{F}_0}\right)$ and $\mathcal{N}_\mathcal{F}^*=\phi^* \mathcal{N}_\mathcal{G}^*$ as a subsheaf of $\Omega_{\mathcal{X}/S}^1$, or equivalently the commutative diagrams
\begin{center}
$\begin{CD}
0 @>>> T_\mathcal{F} @>>> T_{\mathcal{X}/S} \\
@. @V\cong VV @| \\
0@>>> \phi^* T_\mathcal{G} @>>> \phi^*T_{\mathcal{Y}/S}=T_{\mathcal{X}/S}
\end{CD}$\,\,\,\,\,\,\,\,\,\,
$\begin{CD}
0 @>>> \mathcal{N}_\mathcal{F}^* @>>>  \Omega_{\mathcal{X}/S}^1 \\
@. @V\cong VV @| \\
0@>>> \phi^* \mathcal{N}_\mathcal{G}^* @>>> \phi^* \Omega_{\mathcal{Y}/S}^1= \Omega_{\mathcal{X}/S}^1
\end{CD}$
\end{center}

\begin{definition}[Trivial foliated deformations in terms of tangent and cotangent sheaves]
Let $\left( X, T_{\mathcal{F}_0}, \mathcal{N}_{\mathcal{F}_0}^* \right)$ be a foliated nonsingular variety over $k$, and $(S,s)$ be a $k$-pointed scheme. We define the trivial foliated family induced by $\left( X, T_{\mathcal{F}_0},\mathcal{N}_{\mathcal{F}_0}^* \right)$ and $(S,s)$ to be the following foliated deformation of $\left( X, T_{\mathcal{F}_0}, \mathcal{N}_{\mathcal{F}_0}^* \right)$ over $S$,
\begin{center}
$\begin{CD}
\left( X, T_{\mathcal{F}_0} , \mathcal{N}_{\mathcal{F}_0}^* \right) @>>> \left(  X\times_k S    , \pi_1^* T_{\mathcal{F}_0} , \pi^* \mathcal{N}_{\mathcal{F}_0}^* \right)\\
@VVV @VV\pi_2 V\\
\textnormal{Spec}(k) @>s>> S
\end{CD}$
\end{center}
where $\pi_1:X\times S\to X$ and $\pi_2:X\times S \to S$ are the projections. A flat foliated deformation of $\left( X, T_{\mathcal{F}_0} , \mathcal{N}_{\mathcal{F}_0}^* \right)$ over $S$ is called trivial if it is isomorphic to the trivial foliated family over $S$.
\end{definition}

\begin{definition}[Rigid foliated deformations in terms of both tangent and cotangent sheaves] \label{a10}
Let $\left( X, T_{\mathcal{F}_0}\right)$ be a foliated nonsingular variety. $\left(X, T_{\mathcal{F}_0}, \mathcal{N}_{\mathcal{F}_0}^* \right)$ is called rigid in foliated deformations in terms of both tangent and cotangent sheaves if every infinitesimal foliated deformation of $\left( X, T_{\mathcal{F}_0} , \mathcal{N}_{\mathcal{F}_0}^* \right)$ over $A$ is trivial for every $A \in \bold{Art}$.

\end{definition}

\section{Foliated deformations of $\left(X, T_{\mathcal{F}_0}, \mathcal{N}_{\mathcal{F}_0}^* \right)$ with both $T_{\mathcal{F}_0}$ and $\mathcal{N}_{\mathcal{F}_0}^*$ locally free}

Let $\left( X, T_{\mathcal{F}_0},  \mathcal{N}_{\mathcal{F}_0}^* \right)$ be a foliated nonsingular variety with $T_{\mathcal{F}_0}$ and $\mathcal{N}_{\mathcal{F}_0}^*$ locally free. In this section we study infinitesimal simultaneous deformations of $\left( X, T_{\mathcal{F}_0}, \mathcal{N}_{\mathcal{F}_0}^* \right)$ in the language of functors of Artin rings.

\begin{definition}[flat foliated deformations in terms of both tangent and cotangent sheaves]\label{ac1}
Let $\left( X, T_{\mathcal{F}_0}, \mathcal{N}_{\mathcal{F}_0}^* \right)$ be a foliated nonsingular variety with $T_{\mathcal{F}_0}$ and $\mathcal{N}_{\mathcal{F}_0}^*$ locally free. An infinitesimal simultaneous foliated deformation $\eta$ of $\left(X, T_{\mathcal{F}_0}, \mathcal{N}_{\mathcal{F}_0}^* \right)$ over $A\in \bold{Art}$
\begin{center}
$\eta:\begin{CD}
\left(  X, T_{\mathcal{F}_0}, \mathcal{N}_{\mathcal{F}_0}^* \right) @>>> \left( \mathcal{X}, T_\mathcal{F}, \mathcal{N}_\mathcal{F}^* \right) \\
@VVV @VV\pi V\\
\textnormal{Spec}(k) @>>> \textnormal{Spec}(A)
\end{CD}$
\end{center}
where $\pi$ is flat, $T_{\mathcal{F}}$ is a locally free subsheaf of $T_{\mathcal{X}/A}$ and $\mathcal{N}_{\mathcal{F}}^*$ is a locally free subsheaf of $\Omega_{\mathcal{X}/A}^1$ such that $T_{\mathcal{F}}=\left( \frac{\Omega_{\mathcal{X}/A}^1}{\mathcal{N}_\mathcal{F}^*} \right)^*$ and $\mathcal{N}_\mathcal{F}^*=\left( \frac{T_{\mathcal{X}/A}}{T_{\mathcal{F}}}\right)^*$, and $\left[ T_{\mathcal{F}}, T_{\mathcal{F}}\right] \subset T_{\mathcal{F}}$ and $d \left( \mathcal{N}_{\mathcal{F},x}^* \right)\subset  \mathcal{N}_{\mathcal{F},x}^* \bigwedge \Omega_{\mathcal{X}/A, x}^1$ for $x\in X-\textnormal{Sing}\left( \mathcal{F}_0 \right)$, and $\frac{T_{\mathcal{X}/A}}{T_\mathcal{F}}$ is flat over $A$ and $\frac{\Omega_{\mathcal{X}/A}^1}{\mathcal{N}_\mathcal{F}^*}$ is flat over $A$. Two infinitesimal simultaneous deformations $\left( \mathcal{X}, T_\mathcal{F}, \mathcal{N}_\mathcal{F}^* \right)$ and $\left(\mathcal{X}', T_{\mathcal{F}'}, \mathcal{N}_{\mathcal{F}'}^* \right)$ of $\left( X,T_{\mathcal{F}_0}, \mathcal{N}_{\mathcal{F}_0}^* \right)$ over $A$ are isomorphic if there is an isomorphism $\phi:\mathcal{X}\to \mathcal{X}'$ over $A$ inducing the identity on $\left( X, T_{\mathcal{F}_0} , \mathcal{N}_{\mathcal{F}_0}^* \right)$ such that $\phi^* T_{\mathcal{F}'}= T_\mathcal{F}$ and $\phi^* \mathcal{N}_{\mathcal{F}'}^* = \mathcal{N}_\mathcal{F}^*$. Then we have a functor of Artin rings
\begin{align*}
\textnormal{FDef}_{\left( X, T_{\mathcal{F}_0}, \mathcal{N}_{\mathcal{F}_0}^* \right)}:\bold{Art} &\to \left( \textnormal{Sets} \right) \\
 A &\mapsto \left\{\textnormal{Infinitesimal simultaneous foliated deformations of $\left( X, T_{\mathcal{F}_0}, \mathcal{N}_{\mathcal{F}_0}^* \right)$} \right\}/\textnormal{isomorphisms}
\end{align*}
\end{definition}

We will show that the following complex of sheaves $\mathcal{F}_0^\bullet$ (see Part I)
{\small{\begin{center}
$\begin{CD}
\cdots \\
@AF_3AA \\
 \mathscr{H}om_{\mathcal{O}_X}\left( \bigwedge^2 T_{\mathcal{F}_0}, \frac{T_X}{T_{\mathcal{F}_0}} \right) \bigoplus \mathscr{H}om_{\mathcal{O}_X}\left(\mathcal{N}_{\mathcal{F}_0}^*, \tilde{\mathcal{S}}^2 \right) \bigoplus \mathscr{H}om_{\mathcal{O}_X}\left(\mathcal{N}_{\mathcal{F}_0}^*, T_{\mathcal{F}_0}^* \right) \\
@AF_2AA \\
 \mathscr{H}om_{\mathcal{O}_X}\left( \bigwedge^2 T_{\mathcal{F}_0}, \frac{T_X}{T_{\mathcal{F}_0}} \right) \bigoplus \mathscr{H}om_{\mathcal{O}_X}\left(\mathcal{N}_{\mathcal{F}_0}^*, \tilde{\mathcal{S}}^2 \right) \bigoplus \mathscr{H}om_{\mathcal{O}_X}\left(\mathcal{N}_{\mathcal{F}_0}^*, T_{\mathcal{F}_0}^* \right) \\
@AF_1AA \\
\mathscr{H}om_{\mathcal{O}_X}\left( T_{\mathcal{F}_0}, \frac{T_X}{T_{\mathcal{F}_0}} \right) \bigoplus \mathscr{H}om_{\mathcal{O}_X}\left( \mathcal{N}_{\mathcal{F}_0}^*, \frac{\Omega_X^1}{\mathcal{N}_{\mathcal{F}_0}^*} \right)\\
@AAA\\
T_X
\end{CD}$
\end{center}}}

\begin{proposition} \label{si11}
Let $\left( X, T_{\mathcal{F}_0} , \mathcal{N}_{\mathcal{F}_0}^* \right)$ be a foliated nonsingular variety with $T_{\mathcal{F}_0}$ and $\mathcal{N}_{\mathcal{F}_0}^*$ locally free.
\begin{enumerate}
\item there is a $1-1$ correspondence
\begin{align*}
\textnormal{FDef}_{\left( X, T_{\mathcal{F}_0} , \mathcal{N}_{\mathcal{F}_0}^* \right)}\left( k[\epsilon] \right) \cong \mathbb{H}^1\left( X, \mathcal{F}_0^\bullet \right)
\end{align*}
\item Let $A\in \bold{Art}$ and $\eta=\left( \mathcal{X}, T_\mathcal{F}, \mathcal{N}_\mathcal{F}^* \right)$ be an infinitesimal foliated deformation of $\left( X, T_{\mathcal{F}_0}, \mathcal{N}_{\mathcal{F}_0}^* \right)$ over $A$. To every small extension $0\to (t) \to \tilde{A} \to A \to 0$, we can associate an element $o_\eta(e)\in \mathbb{H}^2\left( X, \mathcal{F}_0^\bullet \right)$ called the obstruction lifting of $\eta$ to $\tilde{A}$, which is $0$ if and  only if  a lifting of $\eta$ to $\tilde{A}$ exists.
\end{enumerate}
\end{proposition}

\begin{proof}
Let $\left( \mathcal{X}, T_\mathcal{F}, \mathcal{N}_\mathcal{F}^* \right)$ a first-order simultaneous foliated deformation of $\left( X , T_{\mathcal{F}_0} , \mathcal{N}_{\mathcal{F}_0}^* \right)$ over $k[\epsilon]$. We keep the notations of the proof of (1) in Proposition \ref{at25} and Proposition \ref{nb1}. Then we need to consider the additional condition
\begin{align}\label{sd1}
i_{T_{0i}^\alpha+ \epsilon T_i'^\alpha}\left( w_{0i}^\beta + \epsilon w_i'^\beta \right) \Longrightarrow i_{T_{0i}^\alpha}\left(  w_i'^\beta \right) + i_{T_i'^\alpha}\left( w_{0i}^\beta \right)=0
\end{align}

Then $(\ref{sd1})$ implies that
\begin{align*}
\left( \left\{\overline{T_i'}\right\}, \left\{ \overline{w_i' }\right\}, \left\{ p_{ij}\right\}  \right)\in C^0\left( \mathcal{U}, \mathscr{H}om_{\mathcal{O}_X} \left( T_{\mathcal{F}_0}  , \frac{T_X}{T_{\mathcal{F}_0}} \right) \right) \bigoplus C^0\left( \mathcal{U}, \mathscr{H}om_{\mathcal{O}_X}\left( \mathcal{N}_{\mathcal{F}_0}^*, \frac{\Omega_X^1}{\mathcal{N}_{\mathcal{F}_0}^*} \right) \right) \bigoplus C^1\left( \mathcal{U}, T_X \right)
\end{align*}
defines a $1$-cocycle in the following \v Cech resolution of $\mathcal{F}_0^\bullet$:
{\tiny{\begin{center}
$\begin{CD}
\cdots \\
@AF_2AA \\
C^0\left( \mathcal{U}, \bigwedge^2 T_{\mathcal{F}_0}\otimes \frac{T_X}{T_{\mathcal{F}_0}}\right) \bigoplus C^0\left(\mathcal{U}, \left( \mathcal{N}_{\mathcal{F}_0}^* \right)^*\otimes \tilde{\mathcal{S}}^2 \right) \bigoplus C^0\left(\mathcal{U}, \left( \mathcal{N}_{\mathcal{F}_0}^* \right)\otimes T_{\mathcal{F}_0}^* \right) @>-\delta>> \cdots\\
@AF_1AA @AF_1AA \\
C^0\left(\mathcal{U}, T_{\mathcal{F}_0}^*\otimes \frac{T_X}{T_{\mathcal{F}_0}} \right) \bigoplus C^0\left(\mathcal{U}, \mathcal{N}_{\mathcal{F}_0}^*\otimes \frac{\Omega_X^1}{\mathcal{N}_{\mathcal{F}_0}^*} \right)@ >\delta >> C^1\left(\mathcal{U}, T_{\mathcal{F}_0}^*\otimes \frac{T_X}{T_{\mathcal{F}_0}} \right) \bigoplus C^1\left(\mathcal{U}, \mathcal{N}_{\mathcal{F}_0}^*\otimes \frac{\Omega_X^1}{\mathcal{N}_{\mathcal{F}_0}^*} \right) @>-\delta>> \cdots \\
@AF_0AA @AF_0AA @AF_0AA \\
C^0\left( \mathcal{U}, T_X \right) @>-\delta>> C^1\left( \mathcal{U}, T_X \right) @>\delta>> C^2\left( \mathcal{U}, T_X \right) @>-\delta>> \cdots
\end{CD}$
\end{center}}}
For two equivalent first-order deformations of $\left( X, T_{\mathcal{F}_0}, \mathcal{N}_{\mathcal{F}_0}^*\right)$, the associated $1$-cocycles are equivalent.

Next we identify obstructions. Let $e: 0 \to (t) \to \tilde{A} \to A \to 0$ be a small extension in $\bold{Art}$ and $\eta:=\left( \mathcal{X}, T_{\mathcal{F}}, \mathcal{N}_{\mathcal{F}}^*  \right)$ be an infinitesimal simultaneous foliated deformation of $\left( X, T_{\mathcal{F}_0}, \mathcal{N}_{\mathcal{F}_0}^* \right)$ over $A$. We keep the notations of the proof of (2) in Proposition \ref{at25} and Proposition \ref{nb1}. We need to consider the additional condition: since $i_{T_i^\alpha}\left( w_i^\beta \right)=0$ for $\alpha=1,...,p,\beta=1,...,q$, we have
\begin{align*}
i_{\tilde{T}_i^\alpha}\left( \tilde{w}_i^\beta \right) = t G_i^{\alpha \beta}, \,\,\,\,\,\,\,\,\,\textnormal{for some}\,\,\,\, G_i^{\alpha \beta}\in \Gamma\left( U_i, \mathcal{O}_X \right)
\end{align*}
We define $G_i\in \Gamma\left( U_i, \mathscr{H}om_{\mathcal{O}_X}\left( \mathcal{N}_{\mathcal{F}_0}^*, T_{\mathcal{F}_0}^* \right)\right)$ by
\begin{align*}
G_i:\Gamma\left( U_i, \mathcal{N}_{\mathcal{F}_0}^* \right)&\to \Gamma\left(U_i, T_{\mathcal{F}_0}^* \right)\\
 w_{0i}^\beta &\mapsto \left( T_{0i}^\alpha \mapsto G_i^{\alpha\beta} \right)
\end{align*}
Then we claim that 
{\Small{\begin{align*}
&\left( \left\{ \overline{P_i} \right\}, \left\{ B_i \right\}, \left\{ G_i \right\},  \left\{ - \overline{ W_{ij} } \right\}, \left\{ - \overline{R_{ij}} \right\}, \left\{ d_{ijk} \right\} \right)\\
 &\in C^0\left(\mathcal{U}, \mathscr{H}om_{\mathcal{O}_X}\left(\bigwedge^2 T_{\mathcal{F}_0}, \frac{T_X}{T_{\mathcal{F}_0}} \right) \right) \bigoplus C^0\left(\mathcal{U}, \mathscr{H}om_{\mathcal{O}_X}\left( \mathcal{N}_{\mathcal{F}_0}^*, \tilde{\mathcal{S}}^2 \right) \right) \bigoplus C^0\left( \mathcal{U}, \mathscr{H}om_{\mathcal{O}_X}\left( \mathcal{N}_{\mathcal{F}_0}^*, T_{\mathcal{F}_0}^* \right) \right) \\
& \bigoplus C^1\left( \mathcal{U}, \mathscr{H}om_{\mathcal{O}_X}\left( T_{\mathcal{F}_0}, \frac{T_X}{T_{\mathcal{F}_0}} \right) \right) \bigoplus C^1\left( \mathcal{U}, \mathscr{H}om_{\mathcal{O}_X}\left( \mathcal{N}_{\mathcal{F}_0}^*, \frac{\Omega_X^1}{\mathcal{N}_{\mathcal{F}_0}^*} \right) \right) \bigoplus C^2\left( \mathcal{U}, T_X \right)
\end{align*}}}
defines a $2$-cocycle in the above \v Cech resolution of $\mathcal{F}_0^\bullet$. It remains to show that $\beta_2\left(P_i \right)+ \alpha_2\left(B_i \right)- E_1'\left( G_i \right)=0$ and $G_i- G_j+ \beta_1\left(- \overline{W}_{ij} \right) + \gamma_1\left(- \overline{R}_{ij} \right)=0$. First we show that $\beta_2\left(P_i \right)+ \alpha_2\left(B_i \right)- E_1'\left( G_i \right)=0$. In fact, for $x\in U_i-S$, we check locally on $U_x\subset U_i-S$. We note that
\begin{align}\label{ad10}
E_1'\left( G_i\right) \left( T_{0i}^\alpha \wedge T_{0i}^\beta \right)\left( w_{0i}^\gamma \right) = \left[ T_{0i}^\alpha, G_i^{\beta \gamma} \right] - \sum_{\delta=1}^q a_{0i\alpha}^{\gamma \delta} G_i^{\beta \delta}- \left[ T_{0i}^\beta, G_i^{\alpha \gamma}\right] + \sum_{\delta=1}^q a_{0i\beta}^{\gamma \delta} G_i^{\alpha \delta} - \sum_{\eta=1}^p g_{0i\alpha\beta}^\eta G_i^{\eta \gamma}
\end{align}
Then from $(\ref{ad10})$ we have
\begin{align*}
&ti_{P_i\left( T_{0i}^\alpha \wedge T_{0i}^\beta \right)}\left( w_{0i}^\gamma \right)= i_{\left[\tilde{T}_i^\alpha, \tilde{T}_i^\beta \right]- \sum_{\eta=1}^p \tilde{g}_{i\alpha\beta}^\eta \tilde{T}_i^\eta }\left(\tilde{w}_i^\eta \right) =\mathcal{L}_{\tilde{T}_i^\alpha}i_{\tilde{T}_i^\beta}\left( \tilde{w}_i^\gamma \right) - i_{\tilde{T}_i^\beta}\mathcal{L}_{\tilde{T}_i^\alpha}\left( \tilde{w}_i^\gamma \right) - \sum_{\eta=1}^p \tilde{g}_{i\alpha\beta}^\eta i_{\tilde{T}_i^\eta}\left( \tilde{w}_i^\gamma \right)\\
&= t \mathcal{L}_{T_{0i}^\alpha}\left( G_i^{\beta \gamma} \right)- i_{\tilde{T}_i^\beta} i_{\tilde{T}_i^\alpha} d \tilde{w}_i^\gamma - i_{\tilde{T}_i^\beta} d i_{\tilde{T}_i^\alpha } \left( \tilde{w}_i^\gamma \right) - t \sum_{\eta=1}^p g_{0i\alpha\beta}^\eta G_i^{\eta \gamma} \\
&= t \left[ T_{0i}^\alpha, G_i^{\beta \gamma}\right] - i_{\tilde{T}_i^\beta} i_{\tilde{T}_i^\alpha}\left( t C_{i_x}^{\gamma} + \sum_{\delta=1}^q \tilde{a}_{i_x}^{\gamma \delta} \wedge \tilde{w}_i^\delta    \right) - t i_{T_{0i}^\beta} d G_i^{\alpha \gamma} - t \sum_{\eta=1}^p  g_{0i\alpha\beta}^\eta G_i^{\eta \gamma}\\
&= t \left[ T_{0i}^\alpha, G_i^{\beta \gamma}\right] -  ti_{T_{0i}^\alpha \wedge T_{0i}^\beta}\left( C_{i_x}^\gamma \right) - t\sum_{\delta=1}^q i_{\tilde{T}_{0i}^\alpha}\left( a_{0i_x}^{\gamma \delta} \right)G_i^{\beta\delta} + t \sum_{\delta=1}^q i_{T_{0i}^\beta} \left( a_{0i_x}^{\gamma \delta}\right) G_i^{\alpha \delta} - t \left[ T_{0i}^\beta, G_i^{\alpha \gamma}\right]  - t \sum_{\eta=1}^p g_{0i\alpha\beta}^\eta G_i^{\eta \gamma}\\
&=- i_{T_{0i}^\alpha \wedge T_{0i}^\beta}\left( C_{i_x}^\gamma \right) + t E_i'\left( G_i\right)\left( T_{0i}^\alpha \wedge T_{0i}^\beta \right)\left( w_{0i}^\gamma \right)
\end{align*}
This implies that $\beta_2(P_i)+ \alpha_2(B_i)- E_i'(G) = 0$ on $U_i-S$, and so $U_i$. Next we show that $G_i- G_j+ \beta_1\left(- \overline{W}_{ij} \right) + \gamma_1\left(- \overline{R}_{ij} \right)=0$. In fact, we note that
\begin{align} \label{ad5}
t G_i^{\alpha \beta}= i_{ \tilde{\theta}_{ji}^*\left(\tilde{T}_i^\alpha \right)}\left( \tilde{\theta}_{ji}^* \left( \tilde{w}_i^\beta \right)\right)= i_{tW_{ij}^\alpha+\sum_{\eta=1}^p \tilde{r}_{ij}^{\alpha\eta} \tilde{T}_j^\eta} \left( t R_{ij}^\beta + \sum_{\delta=1}^q \tilde{h}_{ij}^{\beta \delta} \tilde{w}_j^\delta \right) = t i_{W_{ij}^\alpha}\left( w_{0i}^\beta \right) + t i_{T_{0i}^\alpha}\left( R_{ij}^\beta \right) +  t  \sum_{\eta=1}^p \sum_{\delta=1}^q r_{0ij}^{\alpha \eta} h_{0ij}^{\beta \delta} G_j^{\eta \delta}
\end{align}
This implies that $-\delta\left(\left\{ G_i \right\} \right) + \beta_1\left( \left\{- \overline{W_{ij}} \right\} \right) + \gamma_1\left( \left\{- \overline{R_{ij}} \right\} \right)=0$.

For another liftings $\tilde{\theta}_{ij}', \tilde{T}_i^\alpha, \tilde{r}_{ij}^{\alpha\beta}, \tilde{g}_{i\alpha\beta}, \tilde{w}_i^\alpha, \tilde{h}_{ij}^{\alpha\beta}$ of $\theta_{ij}, r_{ij}^{\alpha\beta}, g_{i\alpha\beta}, w_i^\alpha, h_{ij}^{\alpha\beta}$ respectively, we will show that the associated $2$-cocycle $\beta:=\left( \left\{\overline{P_i'} \right\}, \left\{ B_i' \right\}, \left\{ G_i'\right\}, \left\{ - \overline{W}'_{ij}  \right\} , \left\{ - \overline{R}'_{ij} \right\}, \left\{ d_{ijk}' \right\}\right)$ is cohomologous to $\alpha:=\left(\left\{\overline{P_i} \right\}, \left\{ B_i \right\}, \left\{ G_i\right\}, \left\{ - \overline{W}_{ij}  \right\} , \left\{ - \overline{R}_{ij} \right\}, \left\{ d_{ijk} \right\} \right)$. In fact, it remains to show that $i_{T_{0i}^\alpha }\left( V_i^\beta \right) + i_{Y_i^\alpha}\left( w_{0i}^\beta \right)= G_i'^{\alpha\beta} - G_i^{\alpha\beta}$. In fact,
\begin{align}\label{ad7}
tG_i'^{\alpha\beta}= i_{\tilde{T}_i'^\alpha}\left( \tilde{w}_i'^\beta \right)= i_{\tilde{T}_i^\alpha + t Y_i^\alpha}\left( \tilde{w}_i^\beta + t V_i^\beta \right) = t G_i^{\alpha\beta} + t i_{T_{0i}^\alpha}\left( V_i^\beta \right) + t i_{Y_i^\alpha} \left( w_{0i}^\beta \right).
\end{align}
So given a small extension $0 \to (t) \to \tilde{A} \to A \to 0$ and an infinitesimal simultaneous deformation $\eta$ of $\left( X, T_{\mathcal{F}_0} , \mathcal{N}_{\mathcal{F}_0}^* \right)$ over $A$, we can associate an element $o_\eta(e):=$ the cohomology class of $\alpha\in \mathbb{H}^2\left( X, \mathcal{F}_0^\bullet \right)$. We note that $o_\eta(e)=0$ if and only if there exists 
\begin{align*}
\left( \left\{ X_i \right\}, \left\{ N_i \right\}, \left\{ d_{ij}\right\} \right)\in C^0\left( \mathscr{H}om_{\mathcal{O}_X}\left( T_{\mathcal{F}_0}, \frac{T_X}{T_{\mathcal{F}_0}} \right)\right) \bigoplus C^0\left( \mathcal{U}, \mathscr{H}om_{\mathcal{O}_X}\left( \mathcal{N}_{\mathcal{F}_0}^*, \frac{\Omega_X^1}{\mathcal{N}_{\mathcal{F}_0}^*} \right) \right) \bigoplus C^1\left( \mathcal{U}, T_X \right)
\end{align*}
such that in particular, $G_i^{\alpha\beta}= i_{X_i^\alpha}\left( w_{0i}^\beta  \right)+ i_{T_{0i}^\alpha}\left( N_i^\beta \right)$. In order to show a lifting of $\eta$ exists, it remains to show that $i_{\tilde{T}_i^\alpha- t X_i^\alpha}\left( \tilde{w}_i^\beta - t N_i^\beta \right)= 0$ which is equivalent to $G_i^{\alpha\beta}= i_{X_i^\alpha}\left( w_{0i}^\beta  \right)+ i_{T_{0i}^\alpha}\left( N_i^\beta \right)$, and we will show that the lifting $\tilde{\eta}=\left( \tilde{\mathcal{X}}, T_{\tilde{\mathcal{F}}}, \mathcal{N}_{\tilde{\mathcal{F}}}^* \right)$ satisfies $T_{\tilde{\mathcal{F}}} = \left( \frac{\Omega_{\tilde{\mathcal{X}}/\textnormal{Spec}\left(\tilde{A}\right) }^1}{ \mathcal{N}_{\tilde{\mathcal{F}}}^* }\right)^*$ and $\mathcal{N}_{\tilde{\mathcal{F}}}^* = \left( \frac{T_{\tilde{\mathcal{X}}/\textnormal{Spec}(A)}}{T_{\tilde{\mathcal{F}}} }  \right)^*$. In fact, first we note that $T_{\tilde{\mathcal{F}}} \subset \left( \frac{\Omega_{\tilde{\mathcal{X}}/\textnormal{Spec}\left( \tilde{A} \right)}^1 }{ \mathcal{N}_{\tilde{\mathcal{F}}}^*} \right)^*$ by our construction. Conversely let $W \in \Gamma\left(U_i, \left( \frac{\Omega_{\tilde{\mathcal{X}}/\textnormal{Spec}\left(\tilde{A}\right) }^1}{ \mathcal{N}_{\tilde{\mathcal{F}}}^* }\right)^* \right)\subset \Gamma\left( U_i, T_{\tilde{\mathcal{X}}/\textnormal{Spec}\left( \tilde{A} \right)} \right)$. Then $i_W\left( \tilde{w}_i'^\beta \right)$ for all $\beta=1,...,q$, where $\tilde{w}_i'^\beta= \tilde{w}_i^\beta - t N_i^\beta$. By modulo $(t)$,  we have $i_{\overline{W}} \left( w_i^\beta \right)=0$, so that $\overline{W}\in \Gamma\left( U_i, T_{\mathcal{F}} \right)$ since $T_{\mathcal{F}}=\left( \frac{\Omega_{\mathcal{X}/\textnormal{Spec}\left(A \right) }^1}{\mathcal{N}_\mathcal{F}^* } \right)^*$. Let $W'\in \Gamma\left( U_i, T_{\tilde{\mathcal{F}}} \right)$ be a lifting of $\overline{W}$. Then $W-W'= t W_0$ for some $W_0\in \Gamma\left(U_i, T_X \right)$. Since $t i_{W_0}\left( w_{0i}^\beta \right)= i_{W-W'}\left( \tilde{w}_i'^\beta  \right)=0$. Then $W_0\in \Gamma\left(U_i, T_{\mathcal{F}_0} \right)$. This implies that $W= W'+ t W_0\in \Gamma\left( U_i, T_{\tilde{\mathcal{F}}} \right)$. Similarly we can show $\mathcal{N}_{\tilde{\mathcal{F}}}^* = \left( \frac{T_{\tilde{\mathcal{X}}/\textnormal{Spec}(A)}}{T_{\tilde{\mathcal{F}}} }  \right)^*$. Hence a lifting of $\eta$ to $\tilde{A}$ exists.

\end{proof}

\subsection{Relationship between $\textnormal{FDef}_{\left(T_{\mathcal{F}_0}, \mathcal{N}_{\mathcal{F}_0}^* \right)}$ and $\textnormal{FDef}_{\left( X, T_{\mathcal{F}_0}, \mathcal{N}_{\mathcal{F}_0}^* \right)}$ and $\textnormal{Def}_X$} \

Let $\left(X, T_{\mathcal{F}_0}, \mathcal{N}_{\mathcal{F}_0}^* \right)$ be a foliated nonsingular variety with $T_{\mathcal{F}_0}$ and $\mathcal{N}_{\mathcal{F}_0}^*$ locally free. Let us consider simultaneous foliated deformations of $\left( X, T_{\mathcal{F}_0}, \mathcal{N}_{\mathcal{F}_0}^* \right)$ with $X$ fixed.

\begin{definition}
Let $\left( X, T_{\mathcal{F}_0}, \mathcal{N}_{\mathcal{F}_0}^* \right)$ be a foliated nonsingular variety with $T_{\mathcal{F}_0}$ and $\mathcal{N}_{\mathcal{F}_0}^*$ locally free. An infinitesimal simultaneous foliated deformation $\eta$ under the trivial deformation of $X$ over $A\in \bold{Art}$ is a pair of a locally free subsheaf $T_\mathcal{F}$ of $T_{\mathcal{X}/A}=T_{X\times \textnormal{Spec}(A)/A}$ and a locally free subsheaf $\mathcal{N}_{\mathcal{F}}^*$ of $\Omega_{\mathcal{X}/A}^1=\Omega_{X\times \textnormal{Spec}(A)/A}^1$ such that $\left[ T_{\mathcal{F}}, T_{\mathcal{F}}\right] \subset T_{\mathcal{F}}$ and $\frac{T_{\mathcal{X}A}}{T_\mathcal{F}}$ is flat over $A$, and $d_{\mathcal{X}/A}\left( \mathcal{N}_{\mathcal{F}, x}^* \right)\subset \mathcal{N}_{\mathcal{F},x}^* \bigwedge \Omega_{\mathcal{X}/A,x}^1$ for $x\in X-\textnormal{Sing}\left( \mathcal{F}_0 \right)$ and $\frac{\Omega_{\mathcal{X}/A}^1}{\mathcal{N}_\mathcal{F}^*}$ is flat over $A$ and $T_\mathcal{F}=\left( \frac{\Omega_{\mathcal{X}/A}^1}{\mathcal{N}_\mathcal{F}^*} \right)^*$ and $\mathcal{N}_\mathcal{F}^* =\left( \frac{T_{\mathcal{X}/A}}{T_\mathcal{F}} \right)^*$, and $T_{\mathcal{F}}$ induces $T_{\mathcal{F}_0}$ and $\mathcal{N}_\mathcal{F}^*$ induces $\mathcal{N}_{\mathcal{F}_0}^*$ on $X$. Then we have a functor of Artin rings
\begin{align*}
\textnormal{FDef}_{\left(T_{\mathcal{F}_0}, \mathcal{N}_\mathcal{F}^* \right)}:\bold{Art} &\to \left(\textnormal{Sets} \right) \\
 A &\mapsto \left\{\textnormal{Infinitesimal simultaneous foliated deformation of $\left(X, T_{\mathcal{F}_0}, \mathcal{N}_{\mathcal{F}_0}^* \right)$ with $X$ fixed} \right\}
\end{align*}
\end{definition}

Then the following truncated complex of sheaves $\mathcal{F}_0^{\bullet \geq 1}$ from $\mathcal{F}_0^\bullet$
{\small{\begin{center}
$\begin{CD}
\cdots \\
@AF_3AA \\
 \mathscr{H}om_{\mathcal{O}_X}\left( \bigwedge^2 T_{\mathcal{F}_0}, \frac{T_X}{T_{\mathcal{F}_0}} \right) \bigoplus \mathscr{H}om_{\mathcal{O}_X}\left(\mathcal{N}_{\mathcal{F}_0}^*, \tilde{\mathcal{S}}^2 \right) \bigoplus \mathscr{H}om_{\mathcal{O}_X}\left(\mathcal{N}_{\mathcal{F}_0}^*, T_{\mathcal{F}_0}^* \right) \\
@AF_2AA \\
 \mathscr{H}om_{\mathcal{O}_X}\left( \bigwedge^2 T_{\mathcal{F}_0}, \frac{T_X}{T_{\mathcal{F}_0}} \right) \bigoplus \mathscr{H}om_{\mathcal{O}_X}\left(\mathcal{N}_{\mathcal{F}_0}^*, \tilde{\mathcal{S}}^2 \right) \bigoplus \mathscr{H}om_{\mathcal{O}_X}\left(\mathcal{N}_{\mathcal{F}_0}^*, T_{\mathcal{F}_0}^* \right) \\
@AF_1AA \\
\mathscr{H}om_{\mathcal{O}_X}\left( T_{\mathcal{F}_0}, \frac{T_X}{T_{\mathcal{F}_0}} \right) \bigoplus \mathscr{H}om_{\mathcal{O}_X}\left( \mathcal{N}_{\mathcal{F}_0}^*, \frac{\Omega_X^1}{\mathcal{N}_{\mathcal{F}_0}^*} \right)
\end{CD}$
\end{center}}}
controls simultaneous foliated deformations of $\left(X, T_{\mathcal{F}_0},  \mathcal{N}_{\mathcal{F}_0}^* \right)$ with $X$ fixed in the language of functors of Artin rings. We will denote the $i$-th cohomology group of $\mathcal{F}_0^{\bullet \geq 1}$ by $\mathbb{H}^i\left( X, \mathcal{F}_0^{\bullet \geq 1} \right)$. Now we note that we have the natural morphisms of functors
\begin{align}\label{a24}
\textnormal{FDef}_{\left( T_{\mathcal{F}_0}, \mathcal{N}_{\mathcal{F}_0}^* \right)}\to \textnormal{FDef}_{\left( X, T_{\mathcal{F}_0}, \mathcal{N}_{\mathcal{F}_0}^* \right)} \to \textnormal{Def}_X
\end{align}
The relationships between functors are described by the following exact sequence of complexes of sheaves
{\small{\begin{equation}\label{a23}
\begin{CD}
\cdots  @. \cdots @. \cdots \\
@AF_3 AA @AF_3AA @AAA \\
\bigwedge^3 T_{\mathcal{F}_0} \otimes \frac{T_X}{T_{\mathcal{F}_0}} \bigoplus \left( \mathcal{N}_{\mathcal{F}_0}^* \right)^*\otimes \tilde{\mathcal{S}}^3   \bigoplus \left( \mathcal{N}_{\mathcal{F}_0}^*  \right)^* \otimes \bigwedge^2 T_{\mathcal{F}_0}^* @>>> \bigwedge^3 T_{\mathcal{F}_0}^* \otimes \frac{T_X}{T_{\mathcal{F}_0}}  \bigoplus \left(\mathcal{N}_{\mathcal{F}_0}^*\right)^* \otimes \tilde{\mathcal{S}}^3  \bigoplus \left(\mathcal{N}_{\mathcal{F}_0}^*\right)^*\otimes  \bigwedge^2 T_{\mathcal{F}_0}^*  @>>> 0 \\
@AF_2 AA @AF_2AA @AAA \\
\bigwedge^2 T_{\mathcal{F}_0} \otimes \frac{T_X}{T_{\mathcal{F}_0}}  \bigoplus \left( \mathcal{N}_{\mathcal{F}_0}^* \right)^*\otimes \tilde{\mathcal{S}}^2 \bigoplus \left( \mathcal{N}_{\mathcal{F}_0}^*\right)^* \otimes T_{\mathcal{F}_0}^*  @>>> \bigwedge^2 T_{\mathcal{F}_0}^*\otimes \frac{T_X}{T_{\mathcal{F}_0}} \bigoplus \left( \mathcal{N}_{\mathcal{F}_0}^* \right)^* \otimes \tilde{\mathcal{S}}^2  \bigoplus \left(\mathcal{N}_{\mathcal{F}_0}^* \right)^* \otimes T_{\mathcal{F}_0}^* @>>>  0  \\
@AF_1 AA @AF_1AA @AAA \\
T_{\mathcal{F}_0}^* \otimes \frac{T_X}{ T_{\mathcal{F}_0}} \bigoplus \left(\mathcal{N}_{\mathcal{F}_0}^*\right)^* \otimes \frac{\Omega_X^1}{\mathcal{N}_{\mathcal{F}_0}^*}  @>>> T_{\mathcal{F}_0}^* \otimes \frac{T_X}{T_{\mathcal{F}_0}} \bigoplus \left(\mathcal{N}_{\mathcal{F}_0}^* \right)^* \otimes \frac{\Omega_X^1}{\mathcal{N}_{\mathcal{F}_0}^*}@>>> 0 \\
 @AAA @AF_0 AA @A AA\\
 0 @>>> T_X @>>> T_X \\
\end{CD}
\end{equation}}}

Then we have the long exact sequence
\begin{align*}
0\to \mathbb{H}^0\left(X, \mathcal{F}_0^\bullet \right)\to  H^0\left( X, T_X \right) \to &\mathbb{H}^0\left( X, \mathcal{F}_0^{\bullet \geq 1}\right) \to \mathbb{H}^1\left(X, \mathcal{F}_0^\bullet \right) \to H^1 \left( X,  T_X \right) \to \mathbb{H}^1\left( X, \mathcal{F}_0^{\bullet \geq 1} \right)  \\
&\to \mathbb{H}^2\left( X, \mathcal{F}_0^\bullet \right) \to \mathbb{H}^2\left( X, T_{\mathcal{F}_0}^\bullet \right) \to \mathbb{H}^2\left( X, \mathcal{F}_0^{\bullet \geq 1} \right) \to \cdots
\end{align*}

This induces morphisms of first-oder deformations
\begin{align*}
\textnormal{FDef}_{\left( T_{\mathcal{F}_0}, \mathcal{N}_{\mathcal{F}_0}^*\right)}\left( k[\epsilon]\right) \cong \mathbb{H}^0\left(X, \mathcal{F}_0^{\bullet \geq 1} \right) \to \textnormal{FDef}_{\left( X, T_{\mathcal{F}_0}, \mathcal{N}_{\mathcal{F}_0}^*\right) } \left( k[\epsilon] \right) \cong \mathbb{H}^1\left( X, \mathcal{F}_0^\bullet \right) &  \to \textnormal{Def}_X\left( k[\epsilon]\right)\cong H^1\left(X, T_X \right)\\
\left(\left\{ \overline{T_i'} \right\},\left\{ \overline{w_i'} \right\}, \left\{ p_{ij} \right\}  \right) &\mapsto \left\{ p_{ij}\right\} 
\end{align*}
and induces the obstruction maps
\begin{align*}
\mathbb{H}^1\left( X, \mathcal{F}_0^{\bullet \geq 1} \right)\to \mathbb{H}^2\left( X, \mathcal{F}_0^\bullet \right) &\to \mathbb{H}^2\left( X, T_{\mathcal{F}_0}^\bullet \right)\\
\left( \left\{ \overline{P_i}\right\}, \left\{ B_i \right\}, \left\{ G_i \right\}, \left\{ - \overline{W_{ij}} \right\}, \left\{ - \overline{R_{ij}} \right\}, \left\{ d_{ijk}\right\} \right) &\mapsto \left( \left\{ \overline{P_i}  \right\} ,\left\{ - \overline{W_{ij}} \right\} , \left\{ d_{ijk} \right\} \right)
\end{align*}

\subsection{Relationship between $\textnormal{FDef}_{\left( X, T_{\mathcal{F}_0}, \mathcal{N}_{\mathcal{F}_0}^* \right)}$ and $\textnormal{FDef}_{\left( X, T_{\mathcal{F}_0} \right)}$}\

Let $\left( X, T_{\mathcal{F}_0}, \mathcal{N}_{\mathcal{F}_0}^* \right)$ be a foliated nonsingular variety with $T_{\mathcal{F}_0}$ and $\mathcal{N}_{\mathcal{F}_0}^*$ locally free. Then we have the natural morphism of functors by forgetting $\mathcal{N}_{\mathcal{F}_0}^*$
\begin{align}\label{a21}
\textnormal{FDef}_{\left(X, T_{\mathcal{F}_0}, \mathcal{N}_{\mathcal{F}_0}^* \right)} \to \textnormal{FDef}_{\left( X, T_{\mathcal{F}_0} \right)}
\end{align}

The relationship between two functors is described by the following exact sequence of complexes of sheaves by the canonical projection $\mathcal{F}_0^\bullet \to T_{\mathcal{F}_0}^\bullet$:
{\small{\begin{equation}\label{a19}
\begin{CD}
\cdots  @. \cdots @. \cdots \\
@AF_3^\sharp AA @AF_3AA @AD_3AA \\
\left( \mathcal{N}_{\mathcal{F}_0}^* \right)^*\otimes \tilde{\mathcal{S}}^3   \bigoplus \left( \mathcal{N}_{\mathcal{F}_0}^*  \right)^* \otimes \bigwedge^3 T_{\mathcal{F}_0}^* @>>> \bigwedge^3 T_{\mathcal{F}_0}^* \otimes \frac{T_X}{T_{\mathcal{F}_0}}  \bigoplus \left(\mathcal{N}_{\mathcal{F}_0}^*\right)^* \otimes \tilde{\mathcal{S}}^3  \bigoplus \left(\mathcal{N}_{\mathcal{F}_0}^*\right)^*\otimes  \bigwedge^2 T_{\mathcal{F}_0}^*  @>>> \bigwedge^3 T_{\mathcal{F}_0} \otimes \frac{T_X}{T_{\mathcal{F}_0}} \\
@AF_2^\sharp AA @AF_2AA @AD_2AA \\
\left( \mathcal{N}_{\mathcal{F}_0}^* \right)^*\otimes \tilde{\mathcal{S}}^2 \bigoplus \left( \mathcal{N}_{\mathcal{F}_0}^*\right)^* \otimes T_{\mathcal{F}_0}^*  @>>> \bigwedge^2 T_{\mathcal{F}_0}^*\otimes \frac{T_X}{T_{\mathcal{F}_0}} \bigoplus \left( \mathcal{N}_{\mathcal{F}_0}^* \right)^* \otimes \tilde{\mathcal{S}}^2  \bigoplus \left(\mathcal{N}_{\mathcal{F}_0}^* \right)^* \otimes T_{\mathcal{F}_0}^* @>>>  \bigwedge^2 T_{\mathcal{F}_0} \otimes \frac{T_X}{T_{\mathcal{F}_0}}   \\
@AF_1^\sharp AA @AF_1AA @AD_1AA \\
\left(\mathcal{N}_{\mathcal{F}_0}^*\right)^* \otimes \frac{\Omega_X^1}{\mathcal{N}_{\mathcal{F}_0}^*}  @>>> T_{\mathcal{F}_0}^* \otimes \frac{T_X}{T_{\mathcal{F}_0}} \bigoplus \left(\mathcal{N}_{\mathcal{F}_0}^* \right)^* \otimes \frac{\Omega_X^1}{\mathcal{N}_{\mathcal{F}_0}^*}@>>> T_{\mathcal{F}_0}^* \otimes \frac{T_X}{ T_{\mathcal{F}_0}} \\
 @AAA @AF_0 AA @AD_0 AA\\
 0 @>>> T_X @>>> T_X \\
\end{CD}
\end{equation}}}
Let us denote the right complex of sheaves by $A^\bullet$. Then we have the long exact sequence
\begin{align*}
0\to \mathbb{H}^0\left(X, \mathcal{F}_0^\bullet \right)\to  \mathbb{H}^0\left( X, T_{\mathcal{F}_0}^\bullet \right) \to &\mathbb{H}^0\left( X, A^\bullet \right)=0 \to \mathbb{H}^1\left(X, \mathcal{F}_0^\bullet \right) \to \mathbb{H}^1\left( X,  T_{\mathcal{F}_0}^\bullet \right) \to \mathbb{H}^1\left( X, A^\bullet \right)  \\
&\to \mathbb{H}^2\left( X, \mathcal{F}_0^\bullet \right) \to \mathbb{H}^2\left( X, T_{\mathcal{F}_0}^\bullet \right) \to \mathbb{H}^2\left( X, A^\bullet \right) \to \cdots
\end{align*}
Then we have the isomorphism between the spaces of infinitesimal foliated automorphisms of $\left( X, T_{\mathcal{F}_0}\right)$ and $\left( X, T_{\mathcal{F}_0}, \mathcal{N}_{\mathcal{F}_0}^* \right)$:
\begin{align*}
\mathbb{H}^0\left( X, T_{\mathcal{F}_0}^\bullet \right) \cong \mathbb{H}^0\left( X , \mathcal{F}_0^\bullet \right)
\end{align*}
and we have the injective morphism of first-oder deformations
\begin{align*}
\textnormal{FDef}_{\left( X, T_{\mathcal{F}_0}, \mathcal{N}_{\mathcal{F}_0}^*\right) } \left( k[\epsilon] \right) \cong \mathbb{H}^1\left( X, \mathcal{F}_0^\bullet \right) & \hookrightarrow \textnormal{FDef}_{\left( X, T_{\mathcal{F}_0} \right)}  \left( k[\epsilon] \right) \cong  \mathbb{H}^1\left( X, T_{\mathcal{F}_0}^\bullet \right)\\
\left(\left\{ \overline{T_i'} \right\},\left\{ \overline{w_i'} \right\}, \left\{ p_{ij} \right\}  \right) &\mapsto \left( \left\{ \overline{T_i'} \right\}, \left\{ p_{ij}\right\} \right)
\end{align*}
and the above long exact sequence induces the obstruction map
\begin{align*}
\mathbb{H}^2\left( X, \mathcal{F}_0^\bullet \right) &\to \mathbb{H}^2\left( X, T_{\mathcal{F}_0}^\bullet \right)\\
\left( \left\{ \overline{P_i}\right\}, \left\{ B_i \right\}, \left\{ G_i \right\}, \left\{ - \overline{W_{ij}} \right\}, \left\{ - \overline{R_{ij}} \right\}, \left\{ d_{ijk}\right\} \right) &\mapsto \left( \left\{ \overline{P_i}  \right\} ,\left\{ - \overline{W_{ij}} \right\} , \left\{ d_{ijk} \right\} \right)
\end{align*}

\subsection{Relationship between $\textnormal{FDef}_{\left(X, T_{\mathcal{F}_0}, \mathcal{N}_{\mathcal{F}_0}^* \right)}$ and $\textnormal{FDef}_{\left( X, \mathcal{N}_{\mathcal{F}_0}^*\right)}$}\

Let $\left( X, T_{\mathcal{F}_0}, \mathcal{N}_{\mathcal{F}_0}^* \right)$ be a foliated nonsingular variety with $T_{\mathcal{F}_0}$ and $\mathcal{N}_{\mathcal{F}_0}^*$ locally free. Then we have the natural morphism of functors by forgetting $T_{\mathcal{F}_0}$
\begin{align}\label{a22}
\textnormal{FDef}_{\left( X , T_{\mathcal{F}_0}, \mathcal{N}_{\mathcal{F}_0}^* \right)} \to \textnormal{FDef}_{\left(X, \mathcal{N}_{\mathcal{F}_0}^* \right)}
\end{align}

Then the relationship between two functors is described by the following exact sequences of complexes of sheaves by the canonical projection $\mathcal{F}_0^\bullet \to \mathcal{N}_{\mathcal{F}_0}^{*\bullet}$:

{\small{\begin{equation}\label{a20}
\begin{CD}
\cdots  @. \cdots @. \cdots \\
@AF_3^\flat AA @AF_3AA @AE_3AA \\
\bigwedge^3 T_{\mathcal{F}_0}^*\otimes \frac{T_X}{T_{\mathcal{F}_0}}   \bigoplus \left( \mathcal{N}_{\mathcal{F}_0}^*  \right)^* \otimes \bigwedge^2 T_{\mathcal{F}_0}^* @>>> \bigwedge^2 T_{\mathcal{F}_0}^* \otimes \frac{T_X}{T_{\mathcal{F}_0}}  \bigoplus \left(\mathcal{N}_{\mathcal{F}_0}^*\right)^* \otimes \tilde{\mathcal{S}}^2  \bigoplus \left(\mathcal{N}_{\mathcal{F}_0}^*\right)^*\otimes \bigwedge^2 T_{\mathcal{F}_0}^*  @>>> \left( \mathcal{N}_{\mathcal{F}_0}^*\right)^* \otimes \tilde{\mathcal{S}}^3 \\
@AF_2^\flat AA @AF_2AA @AE_2AA \\
\bigwedge^2 T_{\mathcal{F}_0}^* \otimes \frac{T_X}{T_{\mathcal{F}_0}} \bigoplus \left( \mathcal{N}_{\mathcal{F}_0}^*\right)^* \otimes T_{\mathcal{F}_0}^*  @>>> \bigwedge^2 T_{\mathcal{F}_0}^*\otimes \frac{T_X}{T_{\mathcal{F}_0}} \bigoplus \left( \mathcal{N}_{\mathcal{F}_0}^* \right)^* \otimes \tilde{\mathcal{S}}^2  \bigoplus \left(\mathcal{N}_{\mathcal{F}_0}^* \right)^* \otimes T_{\mathcal{F}_0}^* @>>>  \left( \mathcal{N}_{\mathcal{F}_0}^* \right)^* \otimes \tilde{\mathcal{S}}^2  \\
@AF_1^\flat AA @AF_1AA @AE_1AA \\
T_{\mathcal{F}_0}^* \otimes \frac{T_X}{T_{\mathcal{F}_0}}  @>>> T_{\mathcal{F}_0}^* \otimes \frac{T_X}{T_{\mathcal{F}_0}} \bigoplus \left(\mathcal{N}_{\mathcal{F}_0}^* \right)^* \otimes \frac{\Omega_X^1}{\mathcal{N}_{\mathcal{F}_0}^*}@>>> \left(\mathcal{N}_{\mathcal{F}_0}^*\right)^* \otimes \frac{\Omega_X^1}{\mathcal{N}_{\mathcal{F}_0}^*} \\
 @AAA @AF_0 AA @AE_0 AA\\
 0 @>>> T_X @>>> T_X \\
\end{CD}
\end{equation}}}
Let us denote the right complex of sheaves by $B^\bullet$. Then we have the long exact sequence
\begin{align*}
0\to \mathbb{H}^0\left(X, \mathcal{F}_0^\bullet \right)\to  \mathbb{H}^0\left( X, \mathcal{N}_{\mathcal{F}_0}^{*\bullet} \right) \to &\mathbb{H}^0\left( X, A^\bullet \right)=0 \to \mathbb{H}^1\left(X, \mathcal{F}_0^\bullet \right) \to \mathbb{H}^1\left( X,  \mathcal{N}_{\mathcal{F}_0}^{* \bullet } \right) \to \mathbb{H}^1\left( X, B^\bullet \right)  \\
&\to \mathbb{H}^2\left( X, \mathcal{F}_0^\bullet \right) \to \mathbb{H}^2\left( X, \mathcal{N}_{\mathcal{F}_0}^{* \bullet} \right) \to \mathbb{H}^2\left( X, B^\bullet \right) \to \cdots
\end{align*}
Then we have the isomorphism between the spaces of infinitesimal foliated automorphisms of $\left( X, \mathcal{N}_{\mathcal{F}_0}^*\right)$ and $\left( X, T_{\mathcal{F}_0}, \mathcal{N}_{\mathcal{F}_0}^* \right)$:
\begin{align*}
\mathbb{H}^0\left( X, \mathcal{N}_{\mathcal{F}_0}^{*\bullet} \right) \cong \mathbb{H}^0\left( X , \mathcal{F}_0^\bullet \right)
\end{align*}
and we have the injective morphism of first-oder deformations
\begin{align*}
\textnormal{FDef}_{\left( X, T_{\mathcal{F}_0}, \mathcal{N}_{\mathcal{F}_0}^*\right) } \left( k[\epsilon] \right) \cong \mathbb{H}^1\left( X, \mathcal{F}_0^\bullet \right) & \hookrightarrow \textnormal{FDef}_{\left( X, \mathcal{N}_{\mathcal{F}_0}^* \right)}  \left( k[\epsilon] \right) \cong  \mathbb{H}^1\left( X, \mathcal{N}_{\mathcal{F}_0}^{*\bullet} \right)\\
\left(\left\{ \overline{T_i'} \right\},\left\{ \overline{w_i'} \right\}, \left\{ p_{ij} \right\}  \right) &\mapsto \left( \left\{ \overline{w_i'} \right\}, \left\{ p_{ij}\right\} \right)
\end{align*}
and the above long exact sequence induces the obstruction map
\begin{align*}
\mathbb{H}^2\left( X, \mathcal{F}_0^\bullet \right) &\to \mathbb{H}^2\left( X, \mathcal{N}_{\mathcal{F}_0}^{* \bullet } \right)\\
\left( \left\{ \overline{P_i}\right\}, \left\{ B_i \right\}, \left\{ G_i \right\}, \left\{ - \overline{W_{ij}} \right\}, \left\{ - \overline{R_{ij}} \right\}, \left\{ d_{ijk}\right\} \right) &\mapsto \left( \left\{ \overline{B_i}  \right\} ,\left\{ - \overline{R_{ij}} \right\} , \left\{ d_{ijk} \right\} \right)
\end{align*}

\subsection{Foliated deformations of $\left( X, T_{\mathcal{F}_0 } , \mathcal{N}_{\mathcal{F}_0}^* \right)$ in terms of vector bundles}\

The underlying deformations of locally free sheaves $T_{\mathcal{F}_0}$ and $\mathcal{N}_{\mathcal{F}_0}^*$  are implicit in the complex of sheaves $\mathcal{F}_0^\bullet$. In order to make the underlying deformations of $T_{\mathcal{F}_0}$ and $\mathcal{N}_{\mathcal{F}_0}^*$ explicit,  we  define a complex of sheaves also controlling simultaneous flat foliated deformations of $\left(X, T_{\mathcal{F}_0}, \mathcal{N}_{\mathcal{F}_0}^* \right)$ and present foliated deformations of $\left( X, T_{\mathcal{F}_0}, \mathcal{N}_{\mathcal{F}_0}^* \right)$ by using the complex. By combining the Atiyah extension $(\ref{si1})$ associated with $T_{\mathcal{F}_0}$ and the Atiyah extension $(\ref{si2})$ associated with $\mathcal{N}_{\mathcal{F}_0}^*$, we define the generalized Atiyah extension which is the exact sequence
\begin{align}
0\to  \mathscr{H}om_{\mathcal{O}_X}\left( T_{\mathcal{F}_0}, T_{\mathcal{F}_0} \right) \bigoplus \mathscr{H}om_{\mathcal{O}_X}\left( \mathcal{N}_{\mathcal{F}_0}^* , \mathcal{N}_{\mathcal{F}_0}^* \right) \to \mathcal{E}_{\mathcal{F}_0} \to T_X \to 0
\end{align}
where $\mathcal{E}_{\mathcal{F}_0}$ is defined in the following way: let $\left\{ U_i \right\}$ be an affine open covering of $X$ such that $\Gamma\left(U_i, T_{\mathcal{F}_0}\right)\cong \bigoplus^p \Gamma\left( U_i, \mathcal{O}_X \right)$ is trivialized in terms of the basis $\left\{ T_{0i}^\alpha \right\}_{\alpha=1,...,p}$ with the relation $T_{0i}^\alpha= \sum_{\beta=1}^q r_{0ij}^{\alpha\beta} T_{0j}^\beta,\alpha=1,...p$ and $\Gamma\left(U_i, \mathcal{N}_{\mathcal{F}_0}^* \right)\cong \bigoplus^q \Gamma\left( U_i, \mathcal{O}_X \right)$ is trivialized in terms of the basis $\left\{ w_{0i}^\alpha \right\}_{\alpha=1,...,q}$ with the relation $w_{0i}^\alpha= \sum_{\beta=1}^q h_{0ij}^{\alpha\beta} w_{0j}^\beta , \alpha=1,...,q$. Then we define $\mathcal{E}_{\mathcal{F}_0}$ locally on $U_i$
\begin{align*}
\mathcal{E}_{\mathcal{F}_0}|_{U_i} \cong \mathscr{H}om_{\mathcal{O}_X}\left( T_{\mathcal{F}_0}, T_{\mathcal{F}_0} \right)|_{U_i} \bigoplus \mathscr{H}om_{\mathcal{O}_X}\left( \mathcal{N}_{\mathcal{F}_0}^*, \mathcal{N}_{\mathcal{F}_0}^* \right)|_{U_i} \bigoplus T_X |_{U_i}
\end{align*}
such that $\left(\phi_i, \psi_i, X_i \right)\in \mathscr{H}om_{\mathcal{O}_X}\left( T_{\mathcal{F}_0}, T_{\mathcal{F}_0} \right) \bigoplus \mathscr{H}om_{\mathcal{O}_X}\left( \mathcal{N}_{\mathcal{F}_0}^*, \mathcal{N}_{\mathcal{F}_0}^* \right)$ and $\left( \phi_j, \psi_j, X_i \right)\in \mathscr{H}om_{\mathcal{O}_X}\left( T_{\mathcal{F}_0}, T_{\mathcal{F}_0} \right)|_{U_j}$ $\bigoplus \mathscr{H}om_{\mathcal{O}_X}\left( \mathcal{N}_{\mathcal{F}_0}^*, \mathcal{N}_{\mathcal{F}_0}^* \right) |_{U_j} \bigoplus T_X|_{U_j}$ are identified on $U_{ij}$ if and only if $X_i=X_j= X$ and
\begin{align*}
\phi_i\left( T_{0i}^\alpha \right) = \phi_j\left( T_{0i}^\alpha \right) + \sum_{\beta=1}^p \left[ X, r_{0ij}^{\alpha\beta} \right] T_{0j}^\beta,\,\,\,\,\,\,\,\,\,\,\,\,\,\,\,\,\,\psi_i\left( w_{0i}^\alpha \right)  = \psi_j\left(w_{0i}^\alpha \right) + \sum_{\beta=1}^q \left[ X, h_{0ij}^{\alpha\beta} \right] w_{0j}^\beta
\end{align*}
Then $\mathcal{E}_{\mathcal{F}_0}$ is well-defined by $(\ref{si3})$ and $(\ref{si4})$. Then we shall construct a complex of sheaves $\mathcal{E}_{\mathcal{F}_0}^\bullet$:
{\small{\begin{equation} \label{si5}
\begin{CD}
\cdots \\
@AF_3AA \\
 \mathscr{H}om_{\mathcal{O}_X}\left( \bigwedge^2 T_{\mathcal{F}_0}, \frac{T_X}{T_{\mathcal{F}_0}} \right) \bigoplus \mathscr{H}om_{\mathcal{O}_X}\left(\mathcal{N}_{\mathcal{F}_0}^*, \tilde{\mathcal{S}}^2 \right) \bigoplus \mathscr{H}om_{\mathcal{O}_X}\left(\mathcal{N}_{\mathcal{F}_0}^*, T_{\mathcal{F}_0}^* \right) \\
@AF_2AA \\
 \mathscr{H}om_{\mathcal{O}_X}\left( \bigwedge^2 T_{\mathcal{F}_0}, \frac{T_X}{T_{\mathcal{F}_0}} \right) \bigoplus \mathscr{H}om_{\mathcal{O}_X}\left(\mathcal{N}_{\mathcal{F}_0}^*, \tilde{\mathcal{S}}^2 \right) \bigoplus \mathscr{H}om_{\mathcal{O}_X}\left(\mathcal{N}_{\mathcal{F}_0}^*, T_{\mathcal{F}_0}^* \right) \\
@AF_1''AA \\
\mathscr{H}om_{\mathcal{O}_X}\left( T_{\mathcal{F}_0}, T_X \right) \bigoplus \mathscr{H}om_{\mathcal{O}_X}\left( \mathcal{N}_{\mathcal{F}_0}^*,  \Omega_X^1 \right)\\
@AF_0'AA \\
\mathcal{E}_{\mathcal{F}_0}
\end{CD}
\end{equation}}}
which also controls foliated deformations of $\left( X, T_{\mathcal{F}_0}, \mathcal{N}_{\mathcal{F}_0}^* \right)$. We will denote the $i$-th cohomology group of $\mathcal{E}_{\mathcal{F}_0}^\bullet$ by $\mathbb{H}^i\left( X , \mathcal{E}_{\mathcal{F}_0}^\bullet \right)$. We define $F_0': \mathcal{E}_{\mathcal{F}_0} \to \mathscr{H}om_{\mathcal{O}_X}\left( T_{\mathcal{F}_0}, T_X \right) \bigoplus \mathscr{H}om_{\mathcal{O}_X}\left( \mathcal{N}_{\mathcal{F}_0}^*, \Omega_X^1 \right)$ in $(\ref{si5})$ locally in the following way:
{\small{\begin{align*}
F_0':\mathcal{E}_{\mathcal{F}_0}|_{U_i}\cong \mathscr{H}om_{\mathcal{O}_X} \left( T_{\mathcal{F}_0}, T_{\mathcal{F}_0} \right)|_{U_i} \bigoplus \mathscr{H}om_{\mathcal{O}_X}\left( \mathcal{N}_{\mathcal{F}_0}^* , \mathcal{N}_{\mathcal{F}_0}^* \right)|_{U_i} \bigoplus T_X|_{U_i} \to \mathscr{H}om_{\mathcal{O}_X}\left( T_{\mathcal{F}_0}, T_X \right)|_{U_i} \bigoplus \mathscr{H}om_{\mathcal{O}_X}\left( \mathcal{N}_{\mathcal{F}_0}^*, \Omega_X^1 \right)|_{U_i}
\end{align*}}}
On $U_i$, we define $F_0'\left(\left( \phi_i, \psi_i, X \right) \right)$ by $(\ref{si7})$ and $(\ref{si8})$
\begin{align*}
F_0'\left(\left( \phi_i, \psi_i, X \right) \right)= \left( D_0'\left(\phi_i, X \right) , E_0'\left(\psi_i, X \right)  \right)
\end{align*}
Then $F_0'$ is well-defined by $(\ref{si9})$ and $(\ref{si10})$. Next we define $F_1''$ by the composition of the natural quotient map $\mathscr{H}om_{\mathcal{O}_X}\left( T_{\mathcal{F}_0}, T_X \right) \bigoplus \mathscr{H}om_{\mathcal{O}_X}\left( \mathcal{N}_{\mathcal{F}_0}^*, \Omega_X^1 \right)\to \mathscr{H}om_{\mathcal{O}_X}\left(  T_{\mathcal{F}_0}, \frac{T_X}{T_{\mathcal{F}_0}} \right) \bigoplus \mathscr{H}om_{\mathcal{O}_X}\left( \mathcal{N}_{\mathcal{F}_0}^*, \frac{\Omega_X^1}{\mathcal{N}_{\mathcal{F}_0}^*} \right)$ with $F_1$.

Then we have the exact sequence of complexes of sheaves
{\tiny{\begin{equation}\label{ad8}
\begin{CD}
\cdots @. \cdots @. \cdots \\
@AAA @AF_3AA @AF_3AA \\
0 @>>> \bigwedge^2 T_{\mathcal{F}_0}^* \otimes \frac{T_X}{T_{\mathcal{F}_0}} \bigoplus \left( \mathcal{N}_{\mathcal{F}_0}^* \right)^* \otimes \tilde{\mathcal{S}}^2 \bigoplus \left(\mathcal{N}_{\mathcal{F}_0}^* \right)^*\otimes \bigwedge^2 T_{\mathcal{F}_0}^* @>>> \bigwedge^2 T_{\mathcal{F}_0}^* \otimes \frac{T_X}{T_{\mathcal{F}_0}} \bigoplus \left( \mathcal{N}_{\mathcal{F}_0}^* \right)^* \otimes \tilde{\mathcal{S}}^2 \bigoplus \left(\mathcal{N}_{\mathcal{F}_0}^* \right)^*\otimes  \bigwedge^2 T_{\mathcal{F}_0}^* \\
@AAA @AF_2AA @AF_2AA \\
 0 @>>> \bigwedge^3 T_{\mathcal{F}_0}^* \otimes \frac{T_X}{T_{\mathcal{F}_0}} \bigoplus \left( \mathcal{N}_{\mathcal{F}_0}^* \right)^* \otimes \tilde{\mathcal{S}}^3 \bigoplus \left(\mathcal{N}_{\mathcal{F}_0}^* \right)^*\otimes T_{\mathcal{F}_0}^* @>>> \bigwedge^2 T_{\mathcal{F}_0}^* \otimes \frac{T_X}{T_{\mathcal{F}_0}} \bigoplus \left( \mathcal{N}_{\mathcal{F}_0}^* \right)^* \otimes \tilde{\mathcal{S}}^2 \bigoplus \left(\mathcal{N}_{\mathcal{F}_0}^* \right)^*\otimes T_{\mathcal{F}_0}^* \\
@AAA @AF_1''AA @AF_1AA \\
T_{\mathcal{F}_0}^* \otimes T_{\mathcal{F}_0} \bigoplus \left( \mathcal{N}_{\mathcal{F}_0}^* \right)^*\otimes \mathcal{N}_{\mathcal{F}_0}^* @>>> T_{\mathcal{F}_0}^*\otimes T_X \bigoplus \left( \mathcal{N}_{\mathcal{F}_0}^* \right)^* \otimes \Omega_X^1 @>>>  T_{\mathcal{F}_0}^*\otimes \frac{T_X}{T_{\mathcal{F}_0}}\bigoplus \left( \mathcal{N}_{\mathcal{F}_0}^* \right)^* \otimes \frac{\Omega_X^1}{\mathcal{N}_{\mathcal{F}_0}^*} \\
@A -\textnormal{Id} AA @AF_0 AA @AF_0'AA \\
T_{\mathcal{F}_0}^* \otimes T_{\mathcal{F}_0} \bigoplus \left( \mathcal{N}_{\mathcal{F}_0}^* \right)^*\otimes \mathcal{N}_{\mathcal{F}_0}^* @>>> \mathcal{E}_{\mathcal{F}_0} @>>> T_X \end{CD}
\end{equation}}}
The above exact sequence implies that
\begin{align*}
\mathbb{H}^i\left( X, \mathcal{E}_{\mathcal{F}_0}^\bullet \right) \cong \mathbb{H}^i\left( X, \mathcal{F}_0^\bullet \right)
\end{align*}
By combining Proposition \ref{at27} and Proposition \ref{si20} and Proposition \ref{si11}, we obtain

\begin{proposition}\label{ad9}
Let $\left(X, T_{\mathcal{F}_0}, \mathcal{N}_{\mathcal{F}_0}^* \right)$ be a foliated nonsingular variety with $T_{\mathcal{F}_0}$ and $\mathcal{N}_{\mathcal{F}_0}^*$ locally free.
\begin{enumerate}
\item there is a $1-1$ correspondence
\begin{align*}
\textnormal{FDef}_{\left(X, T_{\mathcal{F}_0}, \mathcal{N}_{\mathcal{F}_0}^* \right)}\left( k[\epsilon] \right) \cong \mathbb{H}^1\left( X, \mathcal{E}_{\mathcal{F}_0}^\bullet \right)
\end{align*}
\item Let $A\in \bold{Art}$ and $\eta=\left( \mathcal{X}, T_{\mathcal{F}}, \mathcal{N}_{\mathcal{F}}^* \right)$ be an infinitesimal simultaneous foliated deformation of $\left( X, T_{\mathcal{F}_0}, \mathcal{N}_{\mathcal{F}_0}^* \right)$ over $A$. To every small extension $0 \to (t) \to \tilde{A} \to A \to 0$, we can associate an element $o_\eta(e)\in \mathbb{H}^2\left( X, \mathcal{E}_{\mathcal{F}_0}^\bullet \right)$ called the obstruction lifting of $\eta$ to $\tilde{A}$, which is $0$ if and only if a lifting of $\eta$ to $\tilde{A}$ exists.
\end{enumerate}
\end{proposition}

\subsection{Relationship between $\textnormal{FDef}_{\left( X, T_{\mathcal{F}_0}, \mathcal{N}_{\mathcal{F}_0}^* \right)}$ and $\textnormal{FDef}_{\left( X, T_{\mathcal{F}_0} \right)}$ in terms of vector bundles,  and relationship between $\textnormal{FDef}_{\left(X, T_{\mathcal{F}_0}, \mathcal{N}_{\mathcal{F}_0}^* \right)}$ and $\textnormal{FDef}_{\left( X, \mathcal{N}_{\mathcal{F}_0}^*\right)}$ in terms of vector bundles}\

Let $\left( X, T_{\mathcal{F}_0}, \mathcal{N}_{\mathcal{F}_0}^* \right)$ be a foliated nonsingular variety with $T_{\mathcal{F}_0}$ and $\mathcal{N}_{\mathcal{F}_0}^*$ locally free. Then we have the exact sequence of complexes of sheaves by the canonical projection $\mathcal{F}_0^\bullet \to T_{\mathcal{F}_0}^\bullet$:
{\small{\begin{center}
$\begin{CD}
\cdots  @. \cdots @. \cdots \\
@AF_3^\sharp AA @AF_3AA @AD_3AA \\
\left( \mathcal{N}_{\mathcal{F}_0}^* \right)^*\otimes \tilde{\mathcal{S}}^3   \bigoplus \left( \mathcal{N}_{\mathcal{F}_0}^*  \right)^* \otimes \bigwedge^2 T_{\mathcal{F}_0}^* @>>> \bigwedge^2 T_{\mathcal{F}_0}^* \otimes \frac{T_X}{T_{\mathcal{F}_0}}  \bigoplus \left(\mathcal{N}_{\mathcal{F}_0}^*\right)^* \otimes \tilde{\mathcal{S}}^2  \bigoplus \left(\mathcal{N}_{\mathcal{F}_0}^*\right)^*\otimes T_{\mathcal{F}_0}^*  @>>> \bigwedge^3 T_{\mathcal{F}_0} \otimes \frac{T_X}{T_{\mathcal{F}_0}} \\
@AF_2^\sharp AA @AF_2AA @AD_2AA \\
\left( \mathcal{N}_{\mathcal{F}_0}^* \right)^*\otimes \tilde{\mathcal{S}}^2 \bigoplus \left( \mathcal{N}_{\mathcal{F}_0}^*\right)^* \otimes T_{\mathcal{F}_0}^*  @>>> \bigwedge^2 T_{\mathcal{F}_0}^*\otimes \frac{T_X}{T_{\mathcal{F}_0}} \bigoplus \left( \mathcal{N}_{\mathcal{F}_0}^* \right)^* \otimes \tilde{\mathcal{S}}^2  \bigoplus \left(\mathcal{N}_{\mathcal{F}_0}^* \right)^* \otimes T_{\mathcal{F}_0}^* @>>>  \bigwedge^2 T_{\mathcal{F}_0} \otimes \frac{T_X}{T_{\mathcal{F}_0}}   \\
@AF_1'^\sharp AA @AF_1''AA @AD_1'AA \\
\mathscr{H}om_{\mathcal{O}_X}\left(\mathcal{N}_{\mathcal{F}_0}^*, \Omega_X^1 \right)  @>>> T_{\mathcal{F}_0}^* \otimes T_X \bigoplus \left(\mathcal{N}_{\mathcal{F}_0}^* \right)^* \otimes \Omega_X^1 @>>> T_{\mathcal{F}_0}^* \otimes T_X \\
 @AF_0^\sharp AA @AF_0 AA @AD_0' AA\\
 \mathscr{H}om_{\mathcal{O}_X}\left( \mathcal{N}_{\mathcal{F}_0}^* , \mathcal{N}_{\mathcal{F}_0}^* \right) @>>> \mathcal{E}_{\mathcal{F}_0} @>>> \mathcal{E}_{T_{\mathcal{F}_0}} \\
\end{CD}$
\end{center}}}
Let us denote the right complex of sheaves by $B^\bullet$. Then we have the long exact sequence
\begin{align*}
\mathbb{H}^0 \left( X, B^\bullet \right) =0 \to \mathbb{H}^0\left(X, \mathcal{E}_{\mathcal{F}_0}^\bullet \right)\to  \mathbb{H}^0\left( X, \mathcal{E}_{T_{\mathcal{F}_0}}^\bullet \right) \to &\mathbb{H}^1\left( X, B^\bullet \right)=0 \to \mathbb{H}^1\left(X, \mathcal{E}_{\mathcal{F}_0}^\bullet \right) \to \mathbb{H}^1\left( X,  \mathcal{E}_{T_{\mathcal{F}_0}}^\bullet \right) \to \mathbb{H}^2\left( X, B^\bullet \right)  \\
\to &\mathbb{H}^2\left( X, \mathcal{E}_{\mathcal{F}_0}^\bullet \right) \to \mathbb{H}^2\left( X, \mathcal{E}_{T_{\mathcal{F}_0}}^\bullet \right) \to \cdots
\end{align*}
Then we have the isomorphism between the spaces of infinitesimal foliated automorphisms of $\left( X, T_{\mathcal{F}_0}\right)$ and $\left( X, T_{\mathcal{F}_0}, \mathcal{N}_{\mathcal{F}_0}^* \right)$:
\begin{align*}
\mathbb{H}^0\left( X, \mathcal{E}_{T_{\mathcal{F}_0}}^\bullet \right) \cong \mathbb{H}^0\left( X , \mathcal{E}_{\mathcal{F}_0}^\bullet \right)
\end{align*}
and we have the injective morphism of first-oder deformations
\begin{align*}
\textnormal{FDef}_{\left( X, T_{\mathcal{F}_0}, \mathcal{N}_{\mathcal{F}_0}^*\right) } \left( k[\epsilon] \right) \cong \mathbb{H}^1\left( X, \mathcal{E}_{\mathcal{F}_0}^\bullet \right) & \hookrightarrow \textnormal{FDef}_{\left( X, T_{\mathcal{F}_0} \right)}  \left( k[\epsilon] \right) \cong  \mathbb{H}^1\left( X, \mathcal{E}_{T_{\mathcal{F}_0}}^\bullet \right)\\
\left(\left\{ T_i'  \right\},\left\{ w_i'  \right\}, \left\{ \left( p_{ij} , - b_{ij} ,- g_{ij}  \right) \right\}  \right) &\mapsto \left( \left\{ T_i'  \right\}, \left\{ \left( p_{ij} , - b_{ij} \right)\right\} \right)
\end{align*}
and the above long exact sequence induces the obstruction map
\begin{align*}
\mathbb{H}^2\left( X, \mathcal{E}_{\mathcal{F}_0}^\bullet \right) &\to \mathbb{H}^2\left( X, \mathcal{E}_{T_{\mathcal{F}_0}}^\bullet \right)\\
\left( \left\{ \overline{P_i}\right\}, \left\{ B_i \right\}, \left\{ G_i \right\}, \left\{ - W_{ij} \right\}, \left\{ - R_{ij} \right\}, \left\{ \left( d_{ijk} , g_{ijk} , u_{ijk} \right) \right\} \right) &\mapsto \left( \left\{ \overline{P_i}  \right\} ,\left\{ -  W_{ij} \right\} , \left\{ \left( d_{ijk} , g_{ijk} \right) \right\} \right)
\end{align*}
On the other hand, we have a similar relationship between $\textnormal{FDef}_{\left(X, T_{\mathcal{F}_0}, \mathcal{N}_{\mathcal{F}_0}^* \right)}$ and $\textnormal{FDef}_{\left( X, \mathcal{N}_{\mathcal{F}_0}^*\right)}$ in terms of vector bundles.

\section{Foliated deformations of $\left( X, T_{\mathcal{F}_0}, \mathcal{N}_{\mathcal{F}_0}^*  \right)$ II}

We have another complex of sheaves $\mathcal{F}_0'^\bullet$ which also controls simultaneous foliated deformations of $\left( X, T_{\mathcal{F}_0}, \mathcal{N}_{\mathcal{F}_0}^* \right)$ (see Part I):
{\small{\begin{center}
$\begin{CD}
\cdots \\
@AF_3AA \\
 \mathscr{H}om_{\mathcal{O}_X}\left( \bigwedge^3 T_{\mathcal{F}_0}, \frac{T_X}{T_{\mathcal{F}_0}} \right) \bigoplus \mathscr{H}om_{\mathcal{O}_X}\left(\mathcal{N}_{\mathcal{F}_0}^*, \bigwedge^3 T_{\mathcal{F}_0}^* \right) \bigoplus \mathscr{H}om_{\mathcal{O}_X}\left(\mathcal{N}_{\mathcal{F}_0}^*, \bigwedge^2 T_{\mathcal{F}_0}^* \right) \\
@AF_2AA \\
 \mathscr{H}om_{\mathcal{O}_X}\left( \bigwedge^2 T_{\mathcal{F}_0}, \frac{T_X}{T_{\mathcal{F}_0}} \right) \bigoplus \mathscr{H}om_{\mathcal{O}_X}\left(\mathcal{N}_{\mathcal{F}_0}^*, \bigwedge^2 T_{\mathcal{F}_0}^* \right) \bigoplus \mathscr{H}om_{\mathcal{O}_X}\left(\mathcal{N}_{\mathcal{F}_0}^*, T_{\mathcal{F}_0}^* \right) \\
@AF_1AA \\
\mathscr{H}om_{\mathcal{O}_X}\left( T_{\mathcal{F}_0}, \frac{T_X}{T_{\mathcal{F}_0}} \right) \bigoplus \mathscr{H}om_{\mathcal{O}_X}\left( \mathcal{N}_{\mathcal{F}_0}^*, \frac{\Omega_X^1}{\mathcal{N}_{\mathcal{F}_0}^*} \right)\\
@AF_0 AA\\
T_X
\end{CD}$
\end{center}}}
We will denote the $i$-th cohomology group of $\mathcal{F}_0'^\bullet$ by $\mathbb{H}^i\left( X, \mathcal{F}_0'^\bullet \right)$. Then we have $\mathbb{H}^i\left(X, \mathcal{F}_0'^\bullet  \right)\cong \mathbb{H}^i\left( X, \mathcal{F}_0^\bullet \right)$. We will reinterpret Proposition \ref{si11} by using the complex $\mathcal{F}_0'^\bullet$. For this, we modify Definition \ref{ac1} to make the following equivalent definition of simultaneous foliated deformations by changing the form of the integrability condition in terms of cotangent sheaves.

\begin{definition}[Flat foliated deformation in terms of both tangent and cotangent sheaves]
Let $\left(X, T_{\mathcal{F}_0}, \mathcal{N}_{\mathcal{F}_0}^* \right)$ be a foliated nonsingular variety with $T_{\mathcal{F}_0}$ and $\mathcal{N}_{\mathcal{F}_0}^*$ locally free. An infinitesimal simultaneous foliated deformation $\eta$ of $\left( X, T_{\mathcal{F}_0}, \mathcal{N}_{\mathcal{F}_0}^* \right)$ over $A\in \bold{Art}$ is a cartesian diagram of schemes
\begin{center}
$\eta:\begin{CD}
\left(X, T_{\mathcal{F}_0}, \mathcal{N}_{\mathcal{F}_0}^*  \right) @>>> \left( \mathcal{X}, T_\mathcal{F}, \mathcal{N}_\mathcal{F}^* \right)\\
@VVV @VV\pi V \\
\textnormal{Spec}(k) @>>> \textnormal{Spec}(A)
\end{CD}$
\end{center}
where $\pi$ is flat, $T_{\mathcal{F}}$ is a locally free subsheaf of $T_{\mathcal{X}/A}$ and $\mathcal{N}_{\mathcal{F}}^*$ is a locally free subsheaf of $\Omega_{\mathcal{X}/A}^1$ such that $T_{\mathcal{F}}\cong \left(\frac{\Omega_{\mathcal{X}}^1}{\mathcal{N}_{\mathcal{F}}^*} \right)^*$ and $\mathcal{N}_\mathcal{F}^*\cong \left( \frac{T_{\mathcal{X}/A}}{T_{\mathcal{F}}} \right)^*$ and $\left[ T_{\mathcal{F}}, T_{\mathcal{F}}\right] \subset T_{\mathcal{F}}$ and $d_{\mathcal{X}/A} \left( \mathcal{N}_\mathcal{F}^* \right) $ vanishes on $\bigwedge^2 T_{\mathcal{F}}$ where $d_{\mathcal{X}/A}$ is the relative differential over $A$, and $\left( \mathcal{X}, T_\mathcal{F}, \mathcal{N}_\mathcal{F}^* \right)$ induces $\left( X, T_{\mathcal{F}_0}, \mathcal{N}_{\mathcal{F}_0}^* \right)$ over $\textnormal{Spec}(k)$. Two infinitesimal simultaneous foliated deformations $\left( \mathcal{X}, T_{\mathcal{F}}, \mathcal{N}_{\mathcal{F}}^* \right)$ and $\left( \mathcal{X}' , T_{\mathcal{F}'}, \mathcal{N}_{\mathcal{F}'}^* \right)$ of $\left( X, T_{\mathcal{F}_0}, \mathcal{N}_{\mathcal{F}_0}^* \right)$ over $A$ are isomorphic if there is an isomorphism $\phi:\mathcal{X}\to \mathcal{X}'$ over $A$ inducing the identity on $\left( X, T_{\mathcal{F}_0}\right)$ such that $\phi^* T_{\mathcal{F}'}= T_\mathcal{F}$ and $\phi^* \mathcal{N}_{\mathcal{F}'}^*= \mathcal{N}_\mathcal{F}^*$.
\end{definition}

\begin{proposition}\label{ad6}
Let $\left( X, T_{\mathcal{F}_0}, \mathcal{N}_{\mathcal{F}_0}^* \right)$ be a foliated nonsingular variety with $T_{\mathcal{F}_0}$ and $\mathcal{N}_{\mathcal{F}_0}^*$ locally free.
\begin{enumerate}
\item there is a $1-1$ correspondence
\begin{align*}
 \textnormal{FDef}_{\left(X, T_{\mathcal{F}_0}, \mathcal{N}_{\mathcal{F}_0}^* \right)} \cong \mathbb{H}^1\left( X, \mathcal{F}_0'^\bullet \right) 
\end{align*}
\item Let $A\in \bold{Art}$ and $\eta=\left( \mathcal{X}, T_{\mathcal{F}}, \mathcal{N}_{\mathcal{F}}^* \right)$ be an infinitesimal simultaneous foliated deformation of $\left( X, T_{\mathcal{F}_0}, \mathcal{N}_{\mathcal{F}_0}^* \right)$ over $A$. To every small extension $0 \to (t) \to \tilde{A} \to A \to 0$, we can associate an element $o_\eta(e)\in \mathbb{H}^2\left( X, \mathcal{F}_0'^\bullet \right)$ called the obstruction lifting of $\eta$ to $\tilde{A}$, which is $0$ if and only if a lifting of $\eta$ to $\tilde{A}$ exists.
\end{enumerate}
\end{proposition}

\begin{proof}
Let $\left(\mathcal{X}, T_{\mathcal{F}}, \mathcal{N}_\mathcal{F}^* \right)$ be a first-order infinitesimal simultaneous foliated deformation of $\left(X, T_{\mathcal{F}_0}, \mathcal{N}_{\mathcal{F}_0}^*  \right)$ over $k[\epsilon]$. We keep the notations in the proof of (1) in Proposition \ref{at25} and Proposition \ref{nb1} and Proposition \ref{si11}. Then we need to consider the additional conditions:
\begin{align}
&i_{T_i^\alpha \wedge T_i^\beta}\left( d w_i^\gamma \right)=0 \iff i_{\left(T_{0i}^\alpha + \epsilon T_i'^\alpha\right) \wedge \left(  T_{0i}^\beta + \epsilon T_i'^\beta \right)}\left(  d\left( w_{0i}^\gamma + \epsilon w_i'^\gamma \right) \right) = 0 \notag \\
&\Longrightarrow i_{T_{0i}^\alpha \wedge T_i'^\beta}\left( dw_{0i}^\gamma\right) + i_{T_i'^\alpha \wedge T_{0i}^\beta}\left( dw_{0i}^\gamma \right) + i_{T_{0i}^\alpha \wedge T_{0i}^\beta}\left( dw_i'^\gamma \right) = 0 \notag \\
&\Longrightarrow  \sum_{\delta=1}^q a_{0i\alpha}^{\gamma\delta} i_{T_i'^\beta}\left( w_{0i}^\delta \right) - \sum_{\delta=1}^q a_{0i\beta}^{\gamma \delta} i_{T_i'^\alpha}\left( w_{0i}^\delta \right) + i_{T_{0i}^\alpha \wedge T_{0i}^\beta}\left( dw_i'^\gamma \right) = 0 \label{sd2}
\end{align}
Since from $(\ref{sd1})$ we have $i_{T_{0i}^\alpha}\left( w_i'^\beta \right)+ i_{T_{i}'^\alpha}\left( w_{0i}^\beta \right)=0$, $(\ref{sd2})$ implies that
\begin{align}\label{sd3}
-\sum_{\delta=1}^q a_{0i\alpha}^{\gamma\delta} i_{T_{0i}^\beta}\left( w_{i}'^\delta \right) + \sum_{\delta=1}^q a_{0i\beta}^{\gamma \delta} i_{T_i'^\alpha}\left( w_{0i}^\delta \right)+ i_{T_{0i}^\alpha \wedge T_{0i}^\beta}\left( dw_i'^\gamma \right) = 0
\end{align}
Then $(\ref{sd3})$ implies that 
\begin{align*}
\left( \left\{\overline{T_i'}\right\}, \left\{ \overline{w_i' }\right\}, \left\{ p_{ij}\right\}  \right)\in C^0\left( \mathcal{U}, \mathscr{H}om_{\mathcal{O}_X} \left( T_{\mathcal{F}_0}  , \frac{T_X}{T_{\mathcal{F}_0}} \right) \right) \bigoplus C^0\left( \mathcal{U}, \mathscr{H}om_{\mathcal{O}_X}\left( \mathcal{N}_{\mathcal{F}_0}^*, \frac{\Omega_X^1}{\mathcal{N}_{\mathcal{F}_0}^*} \right) \right) \bigoplus C^1\left( \mathcal{U}, T_X \right)
\end{align*}
 defines a $1$-cocycle in the following \v Cech resolution of $\mathcal{F}_0'^\bullet$:
{\tiny{\begin{center}
$\begin{CD}
\cdots \\
@AF_2AA \\
C^0\left( \mathcal{U}, \bigwedge^2 T_{\mathcal{F}_0}\otimes \frac{T_X}{T_{\mathcal{F}_0}}\right) \bigoplus C^0\left(\mathcal{U}, \left( \mathcal{N}_{\mathcal{F}_0}^* \right)^*\otimes \bigwedge^2 T_{\mathcal{F}_0}^* \right) \bigoplus C^0\left(\mathcal{U}, \left( \mathcal{N}_{\mathcal{F}_0}^* \right)\otimes T_{\mathcal{F}_0}^* \right) @>-\delta>> \cdots\\
@AF_1AA @AF_1AA \\
C^0\left(\mathcal{U}, T_{\mathcal{F}_0}^*\otimes \frac{T_X}{T_{\mathcal{F}_0}} \right) \bigoplus C^0\left(\mathcal{U}, \mathcal{N}_{\mathcal{F}_0}^*\otimes \frac{\Omega_X^1}{\mathcal{N}_{\mathcal{F}_0}^*} \right)@ >\delta >> C^1\left(\mathcal{U}, T_{\mathcal{F}_0}^*\otimes \frac{T_X}{T_{\mathcal{F}_0}} \right) \bigoplus C^1\left(\mathcal{U}, \mathcal{N}_{\mathcal{F}_0}^*\otimes \frac{\Omega_X^1}{\mathcal{N}_{\mathcal{F}_0}^*} \right) @>-\delta>> \cdots \\
@AF_0AA @AF_0AA @AF_0AA \\
C^0\left( \mathcal{U}, T_X \right) @>-\delta>> C^1\left( \mathcal{U}, T_X \right) @>\delta>> C^2\left( \mathcal{U}, T_X \right) 
\end{CD}$
\end{center}}}
For two equivalent first-order deformations of $\left( X, T_{\mathcal{F}_0}, \mathcal{N}_{\mathcal{F}_0}^*\right)$, the associated $1$-cocycles are equivalent.

Next we identify obstructions. Let $e: 0\to (t) \to \tilde{A} \to A \to 0$ be a small extension in $\bold{Art}$ and $\eta=\left( \right)$ be an infinitesimal simultaneous foliated deformation of $\left( \mathcal{X}, T_{\mathcal{F}}, \mathcal{N}_\mathcal{F}^* \right)$ over $A$. We keep the notations in the proof of (2) in Proposition \ref{at25} and Proposition \ref{nb1} and Proposition \ref{si11}. We need to consider the additional condition: since $i_{T_i^\alpha \wedge T_i^\beta }\left( dw_i^\gamma \right)=0$ for $\alpha,\beta=1,...p$ and $\gamma=1,...,q$, we have
\begin{align*}
i_{\tilde{T}_i^\alpha \wedge \tilde{T}_i^\beta}\left( d \tilde{w}_i^\gamma \right) = t f_{i\alpha\beta}^\gamma,\,\,\,\,\,\,\,\,\,\,\textnormal{for some}\,\,\, f_{i\alpha\beta}^\gamma \in \Gamma\left( U_i, \mathcal{O}_X \right)
\end{align*}
We define an element $\phi_i\in \Gamma\left( U_i, \mathscr{H}om_{\mathcal{O}_M}\left( \mathcal{N}_{\mathcal{F}_0}^*, \bigwedge^2 T_{\mathcal{F}_0}^* \right) \right)$ by
\begin{align*}
\phi_i : \Gamma\left( U_i, \mathcal{N}_{\mathcal{F}_0}^* \right) &\to \Gamma\left( U_i, \bigwedge^2 T_{\mathcal{F}_0}^* \right) \\
 w_{0i}^\gamma  & \mapsto \left( T_{0i}^\alpha \wedge T_{0i}^\beta \mapsto f_{i\alpha\beta}^\gamma - \sum_{\delta=1}^q a_{0i\alpha}^{\gamma \delta} G_i^{\beta \delta} + \sum_{\delta=1}^q a_{0i\beta}^{\gamma \delta} G_i^{\alpha \delta} \right)
\end{align*}
Then we claim that 
{\Small{\begin{align*}
&\left( \left\{ \overline{P_i} \right\}, \left\{ \phi_i \right\}, \left\{ G_i \right\},  \left\{ - \overline{ W_{ij} } \right\}, \left\{ - \overline{R_{ij}} \right\}, \left\{ d_{ijk} \right\} \right)\\
 &\in C^0\left(\mathcal{U}, \mathscr{H}om_{\mathcal{O}_X}\left(\bigwedge^2 T_{\mathcal{F}_0}, \frac{T_X}{T_{\mathcal{F}_0}} \right) \right) \bigoplus C^0\left(\mathcal{U}, \mathscr{H}om_{\mathcal{O}_X}\left( \mathcal{N}_{\mathcal{F}_0}^*, \tilde{\mathcal{S}}^2 \right) \right) \bigoplus C^0\left( \mathcal{U}, \mathscr{H}om_{\mathcal{O}_X}\left( \mathcal{N}_{\mathcal{F}_0}^*, T_{\mathcal{F}_0}^* \right) \right) \\
& \bigoplus C^1\left( \mathcal{U}, \mathscr{H}om_{\mathcal{O}_X}\left( T_{\mathcal{F}_0}, \frac{T_X}{T_{\mathcal{F}_0}} \right) \right) \bigoplus C^1\left( \mathcal{U}, \mathscr{H}om_{\mathcal{O}_X}\left( \mathcal{N}_{\mathcal{F}_0}^*, \frac{\Omega_X^1}{\mathcal{N}_{\mathcal{F}_0}^*} \right) \right) \bigoplus C^2\left( \mathcal{U}, T_X \right)
\end{align*}}}
defines a $2$-cocycle in the above \v Cech resolution of $\mathcal{F}_0'^\bullet$. It remains to show that $E_2'\left( \phi_i \right)=0$, and $\beta_2\left(P_i  \right)+ \phi_i- E_1'\left(G_i \right)=0$ and $\phi_i-\phi_j + E_1'' \left( -\overline{ R_{ij}}\right)=0$. First we show that $\alpha_2\left(P_i  \right)+ \phi_i- E_1'\left(G_i \right)=0$. In fact,

\begin{align*}
&ti_{P_i\left( T_{0i}^\alpha \wedge T_{0i}^\beta \right)}\left( w_{0i}^\gamma \right)= i_{\left[\tilde{T}_i^\alpha, \tilde{T}_i^\beta \right] - \sum_{\eta=1}^p \tilde{g}_{i\alpha\beta}^\eta \tilde{T}_i^\eta}\left( \tilde{w}_i^\gamma \right) = \mathcal{L}_{\tilde{T}_i^\alpha} i_{\tilde{T}_i^\beta}\left( \tilde{w}_i^\gamma \right) - i_{\tilde{T}_i^\beta}i_{\tilde{T}_i^\alpha}\left( d\tilde{w}_i^\gamma \right) - i_{\tilde{T}_i^\beta }d i_{\tilde{T}_i^\alpha}\left( \tilde{w}_i^\gamma \right)- \sum_{\eta=1}^p \tilde{g}_{i\alpha\beta}^\eta i_{\tilde{T}_i^\eta}\left( \tilde{w}_i^\gamma \right) \\
&= t\left[ T_{0i}^\alpha, G_i^{\beta \gamma} \right]- t f_{i\alpha\beta}^\gamma - ti_{T_{0i}^\beta} d G_i^{\alpha \gamma} - t \sum_{\eta=1}^p g_{0i\alpha\beta}^\eta G_i^{\eta \gamma}=  t\left[ T_{0i}^\alpha, G_i^{\beta \gamma} \right]- t f_{i\alpha\beta}^\gamma - t \left[ T_{0i}^\beta,   G_i^{\alpha \gamma} \right] - t \sum_{\eta=1}^p g_{0i\alpha\beta}^\eta G_i^{\eta \gamma}
\end{align*}
Then we have
\begin{align*}
i_{P_i^{\alpha\beta}}\left( w_{0i}^\gamma \right) + f_{i\alpha\beta}^\gamma - \sum_{\delta=1}^q a_{0i\alpha}^{\gamma \delta} G_i^{\beta \delta} + \sum_{\delta=1}^q a_{0i\beta}^{\gamma \delta} G_i^{\alpha \delta} = \left[ T_{0i}^\alpha, G_i^{\beta \gamma} \right] - \sum_{\delta=1}^q a_{0i\alpha}^{\gamma \delta} G_i^{\beta \delta}  - \left[ T_{0i}^\beta, G_i^{\alpha \gamma} \right] + \sum_{\delta=1}^q a_{0i\beta}^{\gamma \delta} G_i^{\alpha \delta} - \sum_{\eta=1}^p  g_{0i\alpha\beta}^\eta G_i^{\eta \gamma}
\end{align*}
This implies that $\beta_2(P_i)+ \phi_i- E_1'(G_i)=0$, which also implies that $E_2' \left(\phi_i \right)= E_2' E_1'(G_i)- E_2'\circ \alpha_2(P_i)= - \alpha_3 D_2\left( P_i \right)=0$. We show that $\phi_i - \phi_j + E_1'' \left(- R_{ij} \right)=0$. In fact, first we note that from Part I, we have $\sum_{\eta=1}^q a_{0i\alpha}^{\gamma \eta} h_{0ij}^{\eta \xi}= i_{T_{0i}^\alpha}\left( dh_{0ij}^{\gamma \xi} \right) + \sum_{\delta=1}^q \sum_{\beta=1}^p h_{0ij}^{\gamma \delta} r_{0ij}^{\alpha\beta} a_{0j\beta}^{\delta \xi}$, and recall from $(\ref{ad5})$ that $G_i^{\alpha\beta}= i_{W_{ij}^\alpha} \left( w_{0i}^\beta \right) + i_{T_{0i}^\alpha} \left( R_{ij}^\beta \right) + \sum_{\eta=1}^p \sum_{\delta=1}^q r_{0ij}^{\alpha \eta} h_{0ij}^{\beta \delta} G_j^{\eta \delta}$. Then we have
{\small{\begin{align*}
&tf_{i\alpha\beta}^\gamma=t i_{T_{0i}^\alpha \wedge T_{0ii}^\beta}\left( dw_{0i}^\gamma \right) = ti_{ \tilde{\theta}_{ji}^*\left(\tilde{T}_i^\alpha \right) \wedge \tilde{\theta}_{ji}^*\left( \tilde{T}_i^\beta \right)}  \left(d \tilde{\theta}_{ji}^*\left( \tilde{w}_i^\gamma \right) \right) =  i_{\left( tW_{ij}^\alpha + \sum_{\eta=1}^p \tilde{r}_{ij}^{\alpha \eta} \tilde{T}_j^\eta  \right)\wedge \left( tW_{ij}^\beta + \sum_{\sigma=1}^p \tilde{r}_{ij}^{\beta \sigma} \tilde{T}_j^\sigma \right)}\left( d \left( t R_{ij}^\gamma + \sum_{\delta=1}^q \tilde{h}_{ij}^{\gamma \delta} \tilde{w}_j^\delta  \right) \right)\\
&= \sum_{\eta, \sigma=1}^p \tilde{r}_{ij}^{\alpha\eta} \tilde{r}_{ij}^{\beta \sigma} i_{\tilde{T}_j^\eta \wedge \tilde{T}_j^\sigma}\left( \sum_{\delta=1}^q \tilde{h}_{ij}^{\gamma \delta} d \tilde{w}_j^\delta +  \sum_{\delta=1}^q d \tilde{h}_{ij}^{\gamma \delta} \wedge \tilde{w}_j^\delta \right) + t i_{W_{ij}^\alpha \wedge T_{0i}^\beta}\left( dw_{0i}^\gamma \right) + t i_{T_{0i}^\alpha \wedge W_{ij}^\beta}\left( dw_{0i}^\gamma \right) + t i_{T_{0i}^\alpha \wedge T_{0i}^\beta}\left( d R_{ij}^\gamma \right)\\
&= t \sum_{\eta, \sigma=1}^p \sum_{\delta=1}^q r_{0ij}^{\alpha \eta} r_{0ij}^{\beta \sigma} h_{0ij}^{\gamma \delta} f_{j\eta \sigma}^\delta +t \sum_{\sigma=1}^p \sum_{\delta=1}^q r_{0ij}^{\beta \sigma} i_{T_{0i}^\alpha} \left( d h_{0ij}^{\gamma \delta} \right) G_j^{\sigma \delta} - t \sum_{\eta=1}^p \sum_{\delta=1}^q r_{0ij}^{\alpha \eta} i_{T_{0i}^\beta}\left( d h_{0ij}^{\gamma \delta} \right) G_j^{\eta \delta} \\
&-t \sum_{\delta=1}^q a_{0i\beta}^{\gamma \delta} i_{W_{ij}^\alpha}\left( w_{0i}^\delta \right) + t \sum_{\delta=1}^q a_{0i\alpha}^{\gamma \delta} i_{W_{ij}^\beta}\left( w_{0i}^\delta \right) + t i_{T_{0i}^\alpha \wedge T_{0i}^\beta} \left( d R_{ij}^\gamma \right) \\
&=  t \sum_{\eta, \sigma=1}^p \sum_{\delta=1}^q r_{0ij}^{\alpha \eta} r_{0ij}^{\beta \sigma} h_{0ij}^{\gamma \delta} f_{j\eta \sigma}^\delta + t \sum_{\sigma=1}^p \sum_{\delta=1}^q r_{0ij}^{\beta \sigma}\left( \sum_{\xi=1}^q a_{0i\alpha}^{\gamma \xi} h_{0ij}^{\xi \delta} - \sum_{\xi=1}^q \sum_{\eta=1}^p  h_{0ij}^{\gamma \xi} r_{0ij}^{\alpha \eta} a_{0j\eta}^{\xi \delta}  \right) G_j^{\sigma \delta} \\
& -  t \sum_{\eta=1}^p \sum_{\delta=1}^q r_{0ij}^{\alpha \eta}\left( \sum_{\xi=1}^q a_{0i\beta}^{\gamma \xi} h_{0ij}^{\xi \delta} - \sum_{\xi=1}^q \sum_{\sigma=1}^p  h_{0ij}^{\gamma \xi} r_{0ij}^{\beta \sigma} a_{0j\sigma}^{\xi \delta}  \right) G_j^{\eta \delta} -t \sum_{\delta=1}^q a_{0i\beta}^{\gamma \delta} \left( G_i^{\alpha \delta} - i_{T_{0i}^\alpha}\left( R_{ij}^\delta \right)  - \sum_{\eta=1}^p\sum_{\xi=1}^q r_{0ij}^{\alpha \eta} h_{0ij}^{\delta \xi} G_j^{\eta \xi} \right) \\
& + t \sum_{\delta=1}^q a_{0i\alpha}^{\gamma \delta} \left( G_i^{\beta \delta} - i_{T_{0i}^\beta}\left( R_{ij}^\delta \right)  - \sum_{\eta=1}^p\sum_{\xi=1}^q r_{0ij}^{\beta \eta} h_{0ij}^{\delta \xi} G_j^{\eta \xi} \right) + t i_{T_{0i}^\alpha \wedge T_{0i}^\beta} \left( d R_{ij}^\gamma \right)
\end{align*}}}
Then we have
\begin{align*}
&f_{i\alpha\beta}^\gamma - \sum_{\delta=1}^q a_{0i\alpha}^{\gamma \delta} G_i^{\beta \delta} + \sum_{\delta=1}^q a_{0i\beta}^{\gamma \delta} G_i^{\alpha \delta} \\
&= \sum_{\eta, \sigma=1}^p \sum_{\delta=1}^q r_{0ij}^{\alpha \eta} r_{0ij}^{\beta \sigma} h_{0ij}^{\gamma \delta}\left( f_{j\eta\sigma}^\delta  - \sum_{\xi=1}^q a_{0j\eta}^{\delta \xi} G_j^{\sigma \xi} + \sum_{\xi=1}^q a_{0j\sigma}^{\delta \xi} G_j^{\eta \xi}    \right) - \sum_{\delta=1}^q a_{0i\alpha}^{\gamma \delta} i_{T_{0i}^\beta}\left( R_{ij}^\delta \right) + \sum_{\delta=1}^q a_{0i\beta}^{\gamma \delta} i_{T_{0i}^\alpha}\left( R_{ij}^\delta \right) + i_{T_{0i}^\alpha \wedge T_{0i}^\beta} \left( d R_{ij}^\gamma \right)
\end{align*}
This implies that $\phi_i- \phi_j +E_1''\left( - R_{ij} \right)=0$. For another lifting $\tilde{\theta}_{ij}', \tilde{T}_i^\alpha, \tilde{r}_{ij}^{\alpha\beta}, \tilde{g}_{i\alpha\beta}^\gamma, \tilde{w}_i^\alpha, \tilde{h}_{ij}^{\alpha\beta}$ of $\theta_{ij}, T_i^\alpha, r_{ij}^{\alpha\beta}, g_{i\alpha\beta}^\gamma, w_i^\alpha, h_{ij}^{\alpha\beta}$ respectively, we show that the associated $2$-cocycle $\beta:=\left( \left\{ \overline{P_i'} \right\} , \left\{ \phi_i'\right\}, \left\{ G_i'\right\}, \left\{ - \overline{W_{ij}'}  \right\}, \left\{ - \overline{R_{ij}'}\right\}, \left\{ d_{ijk}'\right\}\right)$ is cohomoglous to $\alpha:=\left( \left\{ \overline{P_i} \right\} , \left\{ \phi_i\right\}, \left\{ G_i\right\}, \left\{ - \overline{W_{ij}}  \right\}, \left\{ - \overline{R_{ij}}\right\}, \left\{ d_{ijk}\right\} \right)$. It remains to show that 
\begin{align*}
f_{i\alpha\beta}'^\gamma - f_{i\alpha\beta}^\gamma - \sum_{\delta=1}^q a_{0i\alpha}^{\gamma \delta} i_{T_{0i}^\beta}\left( G_i'^{\alpha \delta} - G_i^{\alpha \delta}  \right)+ \sum_{\delta=1}^q a_{0i\beta}^{\gamma \delta}\left( G_i'^{\beta \delta} - G_i^{\beta \delta} \right)= - \sum_{\delta=1}^q a_{0i\alpha}^{\gamma \delta} i_{T_{0i}^\beta}\left( V_i^\delta \right) + \sum_{\delta=1}^q a_{0i\beta}^{\gamma \delta} i_{T_{0i}^\alpha}\left( V_i^\delta \right) + i_{T_{0i}^\alpha \wedge T_{0i}^\beta}\left( d V_i^\gamma \right)
\end{align*}
 In fact, since we have  $i_{T_{0i}^\alpha}\left( V_i^\beta \right)+i_{Y_i^\alpha}\left( w_{0i}^\beta  \right)= G_i'^{\alpha \beta}- G_i^{\alpha\beta}$ from $(\ref{ad7})$,
\begin{align*}
&t f_{i\alpha\beta}'^\gamma= i_{T_{0i}'^\alpha \wedge T_{0i}'^\beta}\left( d \tilde{w}_i'^\gamma \right) = i_{\left( T_i^\alpha + t Y_i^\alpha \right) \wedge \left( T_i^\beta + t Y_i^\beta \right) } \left( d\left( w_i^\gamma + t V_i^\gamma \right) \right) \\
&= t f_{i\alpha\beta}^\gamma + t i_{T_{0i}^\alpha \wedge Y_i^\beta}\left( dw_{0i}^\gamma \right) + t i_{Y_i^\alpha \wedge T_{0i}^\beta}\left( d w_{0i}^\gamma \right) + t i_{ T_{0i}^\alpha \wedge T_{0i}^\beta}\left( dV_i^\gamma \right)\\
&= t f_{i\alpha\beta}^\gamma + t \sum_{\delta=1}^q a_{0i\alpha}^{\gamma \delta} \left(G_i'^{\beta \delta} - G_i^{\beta \delta} - i_{T_{0i}^\beta}\left( V_i^\delta \right) \right) - t \sum_{\delta=1}^q a_{0i\beta}^{\gamma \delta} \left( G_i'^{\alpha \delta} - G_i^{\alpha \delta} - i_{T_{0i}^\alpha}\left( V_i^\delta \right)  \right) + t i_{T_{0i}^\alpha \wedge T_{0i}^\beta}\left( dV_i^\gamma \right) 
\end{align*}
So given a small extension $0\to (t)\to \tilde{A} \to A \to 0$ and an infinitesimal simultaneous foliated deformation $\eta$ of $\left( X, T_{\mathcal{F}_0}, \mathcal{N}_{\mathcal{F}_0}^* \right)$ over $A$, we can associate an element $o_\eta(e):=$ the cohomology class of $\alpha \in \mathbb{H}^2\left( X, \mathcal{F}_0'^\bullet \right)$. We note that $o_\eta(e)=0$ if and only if there exists 
\begin{align*}
\left( \left\{ X_i \right\} , \left\{ N_i \right\}, \left\{ d_{ij} \right\} \right)\in C^0\left(\mathcal{U}, \mathscr{H}om_{\mathcal{O}_X}\left( T_{\mathcal{F}_0}, \frac{T_X}{T_{\mathcal{F}_0}} \right) \right) \bigoplus C^0\left(\mathcal{U}, \mathscr{H}om_{\mathcal{O}_X}\left( \mathcal{N}_{\mathcal{F}_0}^*, \frac{\Omega_X^1}{\mathcal{N}_{\mathcal{F}_0}^*} \right) \right) \bigoplus C^1\left(\mathcal{U} , T_X \right)
\end{align*}
such that in particular, $E_1''(N_i)= \phi_i$. In order to show a lifting of $\eta$ exists, it remains to show that $i_{\left( \tilde{T}_i^\alpha -t X_i^\alpha \right)\wedge \left( \tilde{T}_i^\beta - t X_i^\beta \right)} \left( d\left(  \tilde{w}_i^\gamma -t N_i^\gamma \right) \right)=0$, In fact, since $G_i^{\alpha\beta}= i_{X_i^\alpha}\left( w_{0i}^\beta\right) + i_{T_{0i}^\alpha}\left( N_i^\beta \right)$ from $(\ref{ad7})$, we have
\begin{align*}
&i_{\left( \tilde{T}_i^\alpha -t X_i^\alpha \right)\wedge \left( \tilde{T}_i^\beta - t X_i^\beta \right)} \left( d\left(  \tilde{w}_i^\gamma -t N_i^\gamma \right) \right)=0 \\ 
&\iff t f_{i\alpha\beta}^\gamma - t i_{T_{0i}^\alpha \wedge X_i^\beta}\left( dw_{0i}^\gamma \right) - t i_{X_i^\alpha \wedge T_{0i}^\beta}\left( dw_{0i}^\gamma \right) - t i_{T_{0i}^\alpha \wedge T_{0i}^\beta }\left( d N_i^\gamma \right)= 0 \\
&\iff  f_{i\alpha\beta}^\gamma - \sum_{\delta=1}^q a_{0i\alpha}^{\gamma \delta} i_{X_i^\beta}\left( w_{0i}^\gamma \right) + \sum_{\delta=1}^q a_{0i\beta}^{\gamma \delta} i_{X_i^\alpha}\left( w_{0i}^\delta \right) - i_{T_{0i}^\alpha \wedge T_{0i}^\beta}\left( d N_i^\gamma \right) = 0 \\
& \iff  f_{i\alpha\beta}^\gamma - \sum_{\delta=1}^q a_{0i\alpha}^{\gamma \delta}\left( G_i^{\beta \delta} - i_{T_{0i}^\beta}\left(N_i^\delta \right) \right) + \sum_{\delta=1}^q a_{0i\beta}^{\gamma \delta} \left( G_i^{\alpha\delta} - i_{T_{0i}^\alpha}\left(N_i^\delta \right)\right) - i_{T_{0i}^\alpha \wedge T_{0i}^\beta}\left( d N_i^\gamma \right) = 0 \\
&\iff \phi_i\left( T_{0i}^\alpha \wedge T_{0i}^\beta \right)\left( w_{0i}^\gamma \right) = E_1''\left( N_i \right)\left( T_{0i}^\alpha \wedge T_{0i}^\beta \right) \left( w_{0i}^\gamma \right)
\end{align*}
This proves the claim. Hence a lifting of $\eta$ to $\tilde{A}$ exists.

\end{proof}

We have the natural morphism of complexes of sheaves $\mathcal{F}^\bullet \to \mathcal{F}'^\bullet$. We will show the the morphism of complexes of sheaves induces isomorphisms of first-order deformations and obstructions. In fact, we describe the correspondence of first-order deformations between $\mathbb{H}^1\left( X, \mathcal{F}_0^\bullet \right)$ and $\mathbb{H}^1\left( X, \mathcal{F}_0'^\bullet \right)$. Since the first-two terms of $\mathcal{F}_0^\bullet$ and $\mathcal{F}_0'^\bullet$ are identical,  we have the identity map (see the proof of (1) in Proposition \ref{si11})
\begin{align*}
\mathbb{H}^1\left( X, \mathcal{F}_0^\bullet \right) &\to \mathbb{H}^1\left( X, \mathcal{F}_0'^\bullet \right) \\
\left( \left\{ \overline{T_i'} \right\} , \left\{ \overline{w_i'} \right\} , \left\{ p_{ij} \right\} \right) &\mapsto \left(  \left\{ \overline{T_i' } \right\} , \left\{ \overline{w_i'}  \right\} , \left\{ p_{ij} \right\} \right)
\end{align*}
We describe the obstruction map. Given a simultaneous foliated deformation of $\left( X, T_{\mathcal{F}_0}, \mathcal{N}_{\mathcal{F}_0}^* \right)$ over $A\in \bold{Art}$ and a small extension $e:0\to (t) \to \tilde{A} \to A \to 0$, the obstruction lifting to $\tilde{A}$ in terms of $\mathbb{H}^2\left( X, \mathcal{F}_0^\bullet \right)$ is given by $\left(\left\{ \overline{P_i} \right\}, \left\{  B_i \right\} , \left\{ G_i \right\}, \left\{ - \overline{W_{ij} }\right\}, \left\{ - \overline{R_{ij}}\right\}, \left\{ d_{ijk} \right\} \right)$ as in the proof of Proposition \ref{si11}. On the other hands, the obstruction lifting to $\tilde{A}$ in terms of $\mathbb{H}^2\left( X, \mathcal{F}_0'^\bullet \right)$ is given by $\left(\left\{ \overline{P_i} \right\} , \left\{ \phi_i \right\} ,  \left\{ G_i \right\}, \left\{ - \overline{W_{ij} }\right\}, \left\{ - \overline{R_{ij}}\right\}, \left\{ d_{ijk} \right\} \right)$ as in the proof of Proposition \ref{ad6}. Hence it is sufficient to show that $\alpha_2(B_i)= \phi_i$. Indeed,  we note that for $x\in U_x \subset U_i-S$, $B_i\left( w_{0i}^\gamma \right)= w_{0i}^1 \wedge \cdots \wedge w_{0i}^q \wedge C_{i_x}^\gamma$ from $(\ref{si42})$, where $t C_{i_x}^\gamma= d\tilde{w}_i^\gamma - \sum_{\delta=1}^q \tilde{a}_{i_x}^{\gamma\delta} \wedge \tilde{w}_i^\delta$. Then we have on $U_x$
{\small{\begin{align*} 
ti_{T_{0i}^\alpha \wedge T_{0i}^\beta}\left( C_{i_x}^\gamma \right) &= i_{\tilde{T}_i^\alpha \wedge \tilde{T}_i^\beta}\left( t C_{i_x}^\gamma \right)= i_{\tilde{T}_i^\alpha \wedge \tilde{T}_i^\beta}\left( d \tilde{w}_i^\gamma - \sum_{\delta=1}^q \tilde{a}_{i_x}^{\gamma\delta} \wedge \tilde{w}_i^\delta \right)= t f_{i\alpha\beta}^\gamma - t \sum_{\delta=1}^q i_{T_{0i}^\alpha}\left( a_{0i_x}^{\gamma \delta} \right) G_i^{\beta \delta} + t \sum_{\delta=1}^q i_{T_{0i}^\beta}\left( a_{0i_x}^{\gamma \delta} \right) G_i^{\alpha \delta} \\
&= t f_{i\alpha\beta}^\gamma - t \sum_{\delta=1}^q a_{0i\alpha}^{\gamma \delta} G_i^{\beta \delta} + t \sum_{\delta=1}^q a_{0i\beta}^{\gamma \delta} G_i^{\alpha\delta} = t\phi_i\left( T_{0i}^\alpha \wedge T_{0i}^\beta \right)\left( w_{0i}^\gamma \right)
\end{align*}}}
This implies that $\alpha_2 (B_i)= \phi_i$. This proves the claim.

\subsection{Simultaneous foliated deformations of $\left(X, T_{\mathcal{F}_0}, \mathcal{N}_{\mathcal{F}_0}^* \right)$ in terms of vector bundles II}\

Let $\left( X, T_{\mathcal{F}_0}, \mathcal{N}_{\mathcal{F}_0}^* \right)$ be a foliated nonsingular variety. In order to make the underlying deformations of $T_{\mathcal{F}_0}$ and $\mathcal{N}_{\mathcal{F}_0}^*$ explicit, we define the following complex of sheaves $\mathcal{F}_0'^\bullet$:
{\small{\begin{center}
$\begin{CD}
\cdots \\
@AF_3'AA \\
 \mathscr{H}om_{\mathcal{O}_X}\left( \bigwedge^3 T_{\mathcal{F}_0}, \frac{T_X}{T_{\mathcal{F}_0}} \right) \bigoplus \mathscr{H}om_{\mathcal{O}_X}\left(\mathcal{N}_{\mathcal{F}_0}^*, \bigwedge^3 T_{\mathcal{F}_0}^* \right) \bigoplus \mathscr{H}om_{\mathcal{O}_X}\left(\mathcal{N}_{\mathcal{F}_0}^*, \bigwedge^2 T_{\mathcal{F}_0}^* \right) \\
@AF_2'AA \\
 \mathscr{H}om_{\mathcal{O}_X}\left( \bigwedge^2 T_{\mathcal{F}_0}, \frac{T_X}{T_{\mathcal{F}_0}} \right) \bigoplus \mathscr{H}om_{\mathcal{O}_X}\left(\mathcal{N}_{\mathcal{F}_0}^*, \bigwedge^2 T_{\mathcal{F}_0}^* \right) \bigoplus \mathscr{H}om_{\mathcal{O}_X}\left(\mathcal{N}_{\mathcal{F}_0}^*, T_{\mathcal{F}_0}^* \right) \\
@AF_1'''AA \\
\mathscr{H}om_{\mathcal{O}_X}\left( T_{\mathcal{F}_0}, T_X \right) \bigoplus \mathscr{H}om_{\mathcal{O}_X}\left( \mathcal{N}_{\mathcal{F}_0}^*, \Omega_X^1 \right)\\
@AF_0'AA\\
\mathcal{E}_{\mathcal{F}_0}
\end{CD}$
\end{center}}}
It is sufficient to define $F_1''$. We define $F_1'''$ by the natural quotient map $\mathscr{H}om_{\mathcal{O}_X}\left( T_{\mathcal{F}_0}, T_X \right) \bigoplus \mathscr{H}om_{\mathcal{O}_X}\left( \mathcal{N}_{\mathcal{F}_0}^*, \Omega_X^1 \right) \to \mathscr{H}om_{\mathcal{O}_X}\left( T_{\mathcal{F}_0}, \frac{T_X}{T_{\mathcal{F}_0}} \right) \bigoplus \mathscr{H}om_{\mathcal{O}_X}\left( \mathcal{N}_{\mathcal{F}_0}^*, \frac{\Omega_X^1}{\mathcal{N}_{\mathcal{F}_0}^*} \right) $ with $F_1'$. We will denote $i$-th cohomology group of $\mathcal{E}_{\mathcal{F}_0}'^\bullet$ by $\mathbb{H}^i\left( X, \mathcal{E}_{\mathcal{F}_0}'^\bullet  \right)$.  As in a similar exact sequence with $(\ref{ad8})$, we have $\mathbb{H}^i\left( X, \mathcal{F}_0'^\bullet \right) \cong \mathbb{H}^i\left( X, \mathcal{E}_{\mathcal{F}_0}'^\bullet \right)$. By combining Proposition \ref{ad9} and Proposition \ref{ad6}, we obtain
\begin{proposition}
Let $\left(X, T_{\mathcal{F}_0}, \mathcal{N}_{\mathcal{F}_0}^* \right)$ be a foliated nonsingular variety with $T_{\mathcal{F}_0}$ and $\mathcal{N}_{\mathcal{F}_0}^*$ locally free.
\begin{enumerate}
\item there is a $1-1$ correspondence
\begin{align*}
\textnormal{FDef}_{\left(X, T_{\mathcal{F}_0}, \mathcal{N}_{\mathcal{F}_0}^* \right)}\left( k[\epsilon] \right) \cong \mathbb{H}^1\left( X, \mathcal{E}_{\mathcal{F}_0}'^\bullet \right)
\end{align*}
\item Let $A\in \bold{Art}$ and $\eta=\left( \mathcal{X}, T_{\mathcal{F}}, \mathcal{N}_{\mathcal{F}}^* \right)$ be an infinitesimal simultaneous foliated deformation of $\left( X, T_{\mathcal{F}_0}, \mathcal{N}_{\mathcal{F}_0}^* \right)$ over $A$. To every small extension $0 \to (t) \to \tilde{A} \to A \to 0$, we can associate an element $o_\eta(e)\in \mathbb{H}^2\left( X, \mathcal{E}_{\mathcal{F}_0}'^\bullet \right)$ called the obstruction lifting of $\eta$ to $\tilde{A}$, which is $0$ if and only if a lifting of $\eta$ to $\tilde{A}$ exists.
\end{enumerate}
\end{proposition}

\section{Rigid foliated deformations, Schlessinger's criterion for $\textnormal{FDef}_{\left(X, T_{\mathcal{F}_0} , \mathcal{N}_{\mathcal{F}_0}^* \right) }$ and stability }\

\subsection{Rigid simultaneous foliated deformations}\

We recall Definition \ref{a10} that a foliated nonsingular variety $\left( X, T_{\mathcal{F}_0}, \mathcal{N}_{\mathcal{F}_0}^* \right)$ is called rigid in simultaneous foliated deformations if every infinitesimal simultaneous foliated deformation $\left(\mathcal{X}, T_\mathcal{F}, \mathcal{N}_\mathcal{F}^* \right)$ of $\left( X, T_{\mathcal{F}_0}, \mathcal{N}_{\mathcal{F}_0}^* \right)$ over $A$ is trivial for every $A$ in $\bold{Art}$, in other words, $\left( \mathcal{X}, T_\mathcal{F}, \mathcal{N}_\mathcal{F}^* \right)\cong \left(X\times \textnormal{Spec}\left( A \right), T_{\mathcal{F}_0}, \mathcal{N}_{\mathcal{F}_0}^* \right)$.

\begin{proposition}
A foliated nonsingular variety $\left( X, T_{\mathcal{F}_0}, \mathcal{N}_{\mathcal{F}_0}^* \right)$ with $T_{\mathcal{F}_0}$ and $\mathcal{N}_{\mathcal{F}_0}^*$ locally free is rigid in simultaneous foliated deformations if and only if $\mathbb{H}^1\left( X, \mathcal{F}_0^\bullet \right)=0$.
\end{proposition}

\begin{proof}
Assume that $\left( X, T_{\mathcal{F}_0}, \mathcal{N}_{\mathcal{F}_0}^* \right)$ is rigid in simultaneous foliated deformations. Since any infinitesimal simultaneous foliated deformation (in particular, any first-order deformation) is trivial, we have $\mathbb{H}^1\left( X, \mathcal{F}_0^\bullet \right)=0$ by Proposition \ref{si11}. Assume that $\mathbb{H}^1\left( X, \mathcal{F}_0^\bullet \right)$. We claim that given an infinitesimal simultaneous foliated deformation $\eta$ of $\left( X , T_{\mathcal{F}_0} , \mathcal{N}_{\mathcal{F}_0}^* \right)$ over $A\in \bold{Art}$ and a small extension $0\to (t) \to \tilde{A} \to A \to 0$, any two liftings $\xi,\xi'$ of $\eta$ to $\tilde{A}$ are equivalent. We keep the notations in the proof of Proposition \ref{a7} and Proposition \ref{a11}. Then we note that
\begin{align*}
0=i_{T_i'^\alpha}\left( w_i'^\beta \right)= i_{T_i^\alpha +t V_i^\alpha } \left( w_i^\beta + t B_i^\beta \right) \Longrightarrow i_{T_{0i}^\alpha }\left( B_i^\beta \right)+ i_{V_i^\alpha}\left( w_{0i}^\alpha \right)=0
\end{align*}
This implies that $\left( \left\{ \overline{V_i} \right\} , \left\{ \overline{B_i} \right\}, \left\{ p_{ij}\right\}  \right)\in C^0\left( \mathcal{U}, \mathscr{H}om_{\mathcal{O}_X}\left( T_{\mathcal{F}_0}, \frac{T_X}{T_{\mathcal{F}_0}} \right) \right) \bigoplus C^0\left( \mathcal{U}, \mathscr{H}om_{\mathcal{O}_X}\left( \mathcal{N}_{\mathcal{F}_0}^*, \frac{\Omega_X^1}{\mathcal{N}_{\mathcal{F}_0}^*}  \right) \right) \bigoplus C^1\left( \mathcal{U}, T_X \right)$ defines a $1$-cocycle of the \v Cech resolution of $\mathcal{F}_0^\bullet$. Since $\mathbb{H}^1\left(X, \mathcal{F}_0^\bullet \right)=0$, we see that $\xi$ and $\xi'$ are equivalent, which implies that $\left(X, T_{\mathcal{F}_0}, \mathcal{N}_{\mathcal{F}_0}^* \right)$ is rigid in simultaneous foliated deformation.
\end{proof}

We note that if $\mathbb{H}^1\left( X, T_{\mathcal{F}_0}^\bullet \right)=0$ or $\mathbb{H}^1\left( X, \mathcal{N}_{\mathcal{F}_0}^{*\bullet} \right)=0$, then we have $\mathbb{H}^1\left( X, \mathcal{F}^\bullet \right)=0$. From Proposition \ref{a7} and Proposition \ref{a11}, if $\left( X, T_{\mathcal{F}_0}\right)$ is rigid in foliated deformation in terms of tangent sheaves or $\left(X, \mathcal{N}_{\mathcal{F}_0}^* \right)$ is rigid in foliated deformations in cotangent sheaves, then $\left(X, \mathcal{F}_0 \right)$ is rigid in simultaneous foliated deformations.

\subsection{Schlessinger's criterion for $\textnormal{FDef}_{\left( X, T_{\mathcal{F}_0}, \mathcal{N}_{\mathcal{F}_0}^* \right)}$}\

Let $\left( X, T_{\mathcal{F}_0}, \mathcal{N}_{\mathcal{F}_0}^* \right)$ be a foliated nonsingular variety with $T_{\mathcal{F}_0}$ and $\mathcal{N}_{\mathcal{F}_0}^*$ locally free. Since $\mathbb{H}^0\left( X, T_{\mathcal{F}_0}^\bullet \right)=\mathbb{H}^0\left( X, \mathcal{N}_{\mathcal{F}_0}^{*\bullet} \right)=\mathbb{H}^0\left( X, \mathcal{F}_0^\bullet \right)$, we see that from Proposition \ref{a9} and Proposition \ref{a12} we have
\begin{proposition}
Let $\left( X, T_{\mathcal{F}_0}, \mathcal{N}_{\mathcal{F}_0}^* \right)$ be a foliated nonsingular variety with $T_{\mathcal{F}_0}$ and $\mathcal{N}_{\mathcal{F}_0}^*$ locally free. Assume that $\mathbb{H}^0\left( X, \mathcal{F}_0^\bullet \right)=0$. Then for any infinitesimal simultaneous foliated deformation $\left(\mathcal{X}, T_\mathcal{F}, \mathcal{N}_\mathcal{F}^* \right)$ of $\left( X, T_{\mathcal{F}_0}, \mathcal{N}_{\mathcal{F}_0}^* \right)$ over $A\in \bold{Art}$,
\begin{align*}
\textnormal{Aut}\left( \frac{ \left( \mathcal{X}, T_{\mathcal{F}}, \mathcal{N}_\mathcal{F}^*\right)}{\left( X, T_{\mathcal{F}_0}, \mathcal{N}_{\mathcal{F}_0}^* \right) } \right) = \textnormal{Id}
\end{align*}
where $\textnormal{Aut}\left( \frac{ \left( \mathcal{X}, T_{\mathcal{F}}, \mathcal{N}_\mathcal{F}^*\right)}{\left( X, T_{\mathcal{F}_0}, \mathcal{N}_{\mathcal{F}_0}^* \right) }  \right)$ is the set of foliated automorphisms of $\left( \mathcal{X}, T_\mathcal{F}, \mathcal{N}_\mathcal{F}^* \right)$ restricting the the identity on $\left( X, T_{\mathcal{F}_0}, \mathcal{N}_{\mathcal{F}_0}^* \right)$.
\end{proposition}

By combining Proposition \ref{e4} and Proposition \ref{e7}, we can show

\begin{proposition}
Let $\left( X, T_{\mathcal{F}_0}, \mathcal{N}_{\mathcal{F}_0}^* \right)$ be a foliated nonsingular variety with $T_{\mathcal{F}_0}$ and $\mathcal{N}_{\mathcal{F}_0}^*$ locally free, and $e:0\to (t) \to \tilde{A}\xrightarrow{\mu} A \to 0$ be a small extension in $\bold{Art}$. Let $p:=\textnormal{FDef}_{\left(X, T_{\mathcal{F}_0}, \mathcal{N}_{\mathcal{F}_0}^* \right)}(\mu): \textnormal{FDef}_{\left( X, T_{\mathcal{F}_0}, \mathcal{N}_{\mathcal{F}_0}^* \right)} \left( \tilde{A} \right) \to \textnormal{FDef}_{\left( X, T_{\mathcal{F}_0}, \mathcal{N}_{\mathcal{F}_0}^* \right)}\left( A \right)$. Then given $\xi=\left(\mathcal{X}, T_\mathcal{F}, \mathcal{N}_\mathcal{F}^* \right)$, there is a transitive action of $\mathbb{H}^1\left(X, \mathcal{F}_0^\bullet \right)$ on $p^{-1}(\xi)$. Moreover, if $\mathbb{H}^0\left(X, \mathcal{F}_0^\bullet \right)=0$, then the action is free.
\end{proposition}

By combining Proposition \ref{e3} and Proposition \ref{e5} and Proposition \ref{si11} and considering the additional conditions $\left(\frac{T_{\mathcal{X}/A}}{T_{\mathcal{F}}} \right)^*=\mathcal{N}_\mathcal{F}^*$ and $\left (\frac{\Omega_{\mathcal{X}/A}^1}{\mathcal{N}_\mathcal{F}^*} \right)^* = T_{\mathcal{F}}$ for a foliated deformation $\left( \mathcal{X}, T_\mathcal{F}, \mathcal{N}_\mathcal{F}^* \right)$ over $A \in \bold{Art}$, we obtain

\begin{proposition}\label{e6}
Let $\left(X, T_{\mathcal{F}_0}, \mathcal{N}_{\mathcal{F}_0}^* \right)$ be a foliated nonsingular projective variety with $T_{\mathcal{F}_0}$ and $\mathcal{N}_{\mathcal{F}_0}^*$ locally free. Then $\textnormal{FDef}_{\left( X, T_{\mathcal{F}_0}, \mathcal{N}_{\mathcal{F}_0}^* \right)}$ has a semiuniversal formal element in the sense of \cite{Ser06} or a miniversal family in the sense of \cite{Har10}.
\end{proposition}

\begin{theorem}
Let $\left( X, T_{\mathcal{F}_0}, \mathcal{N}_{\mathcal{F}_0}^* \right)$ be a foliated nonsingular projective variety with $T_{\mathcal{F}_0}$ and $\mathcal{N}_{\mathcal{F}_0}^*$ locally free, and $\mathbb{H}^0\left( X, \mathcal{F}_0^\bullet \right)=0$. Then the functor $\textnormal{FDef}_{\left( X , T_{\mathcal{F}_0}  , \mathcal{N}_{\mathcal{F}_0}^* \right)}$ is pro-representable.
\end{theorem}

\subsection{Stability of a foliated nonsingular variety $\left( X, T_{\mathcal{F}_0}, \mathcal{N}_{\mathcal{F}_0}^* \right)$ in $X$ }\

\begin{definition}
A foliated nonsingular variety $\left(X, T_{\mathcal{F}_0}, \mathcal{N}_{\mathcal{F}_0}^* \right)$ with $T_{\mathcal{F}_0}$ and $\mathcal{N}_{\mathcal{F}_0}^*$ locally free is called stable in $X$ if the morphism of functors from $(\ref{a24})$ by forgetting foliated structures
\begin{align*}
\Phi: \textnormal{FDef}_{\left( X, T_{\mathcal{F}_0}, \mathcal{N}_{\mathcal{F}_0}^* \right)}\to \textnormal{Def}_X
\end{align*}
is smooth.
\end{definition}

We note that $\textnormal{FDef}_{\left(X, T_{\mathcal{F}_0}, \mathcal{N}_{\mathcal{F}_0}^* \right)}$ has a semiuniversal formal element from Proposition \ref{e6} and obstruction space $\mathbb{H}^2\left(X,\mathcal{F}_0^\bullet \right)$ from Proposition \ref{si11}. On the other hand, $\textnormal{Def}_{X}$ has a semiuniversal formal element and obstruction space $\mathbb{H}^2\left( X, T_X\right)$ (see \cite{Ser06}). Assume $\mathbb{H}^1\left(X, \mathcal{F}_0^{\bullet \geq 1} \right)=0$. Then from the long exact sequence of $(\ref{a23})$, we have the surjection $d\Phi: \mathbb{H}^1\left(X, \mathcal{N}_{\mathcal{F}_0}^{* \bullet }\right)\to \mathbb{H}^1\left( X, T_X \right)$ and an injective obstruction map $\mathbb{H}^2\left( X, \mathcal{N}_{\mathcal{F}_0}^{ * \bullet } \right)\to \mathbb{H}^2\left( X, T_X \right)$, in other words, $\textnormal{FDef}_{\left( X, \mathcal{N}_{\mathcal{F}_0}^* \right)}$ is less obstructed than $\textnormal{Def}_X$ (see \cite{Ser06} Definition 2.3.5 p.58). This implies that $\Phi$ is smooth. (see \cite{Ser06} Proposition 2.3.6 p.59). We summarize as follows.

\begin{proposition}
Let $\left( X, T_{\mathcal{F}_0},  \mathcal{N}_{\mathcal{F}_0}^* \right)$ be a foliated nonsingular variety with $T_{\mathcal{F}_0}$ and $\mathcal{N}_{\mathcal{F}_0}^*$ locally free. Assume that $\mathbb{H}^1\left(X, \mathcal{F}_0^{\bullet \geq 1} \right)=0$. Then $\Phi$ is smooth. In other words, for a flat deformation $\mathcal{X}$ of $X$ over $A$, there exists a simultaneous foliated deformation $\left(\mathcal{X}, T_{\mathcal{F}}, \mathcal{N}_\mathcal{F}^* \right)$ over $A$ which induces $\mathcal{X}$.
\end{proposition}

\subsection{Stability of a foliated nonsingular variety $\left( X, T_{\mathcal{F}_0}, \mathcal{N}_{\mathcal{F}_0}^* \right)$ in $\left(X, T_{\mathcal{F}_0}\right)$ }\

\begin{definition}
A foliated nonsingular variety $\left(X, T_{\mathcal{F}_0}, \mathcal{N}_{\mathcal{F}_0}^* \right)$ with $T_{\mathcal{F}_0}$ and $\mathcal{N}_{\mathcal{F}_0}^*$ locally free is called stable in $\left(X, T_{\mathcal{F}_0}\right)$ if the morphism of functors from $(\ref{a21})$ by forgetting $\mathcal{N}_{\mathcal{F}_0}^*$
\begin{align*}
\Phi: \textnormal{FDef}_{\left( X, T_{\mathcal{F}_0}, \mathcal{N}_{\mathcal{F}_0}^* \right)}\to \textnormal{FDef}_{\left(X, T_{\mathcal{F}_0}\right)}
\end{align*}
is smooth.
\end{definition}

We note that $\textnormal{FDef}_{\left(X, T_{\mathcal{F}_0}, \mathcal{N}_{\mathcal{F}_0}^* \right)}$ has a semiuniversal formal element from Proposition \ref{e6} and obstruction space $\mathbb{H}^2\left(X,\mathcal{F}_0^\bullet \right)$ from Proposition \ref{si11}. On the other hand, $\textnormal{FDef}_{\left(X, T_{\mathcal{F}_0} \right)}$ has a semiuniversal formal element and obstruction space $\mathbb{H}^2\left( X, T_{\mathcal{F}_0}^\bullet \right)$ from Proposition \ref{at25}. Assume $\mathbb{H}^1\left(X, A^\bullet \right)=0$. Then from the long exact sequence of $(\ref{a19})$, we have the surjection $d\Phi: \mathbb{H}^1\left(X,  \mathcal{F}_0^\bullet \right)\to \mathbb{H}^1\left( X, T_{\mathcal{F}_0}^\bullet  \right)$ and an injective obstruction map $\mathbb{H}^2\left( X,  \mathcal{F}_0^\bullet \right)\to \mathbb{H}^2\left( X, T_{\mathcal{F}_0}^\bullet \right)$, in other words, $\textnormal{FDef}_{\left( X, T_{\mathcal{F}_0} ,  \mathcal{N}_{\mathcal{F}_0}^* \right)}$ is less obstructed than $\textnormal{FDef}_{\left(X, T_{\mathcal{F}_0}\right)}$ (see \cite{Ser06} Definition 2.3.5 p.58). This implies that $\Phi$ is smooth. (see \cite{Ser06} Proposition 2.3.6 p.59). We summarize as follows.

\begin{proposition}
Let $\left( X, T_{\mathcal{F}_0},  \mathcal{N}_{\mathcal{F}_0}^* \right)$ be a foliated nonsingular variety with $T_{\mathcal{F}_0}$ and $\mathcal{N}_{\mathcal{F}_0}^*$ locally free. Assume that $\mathbb{H}^1\left(X, A^\bullet \right)=0$. Then $\Phi$ is smooth. In other words, for a flat foliated deformation $\left(\mathcal{X}, T_{\mathcal{F}}\right)$ of $\left(X, T_{\mathcal{F}_0}\right)$ over $A$, there exists a simultaneous foliated deformation $\left(\mathcal{X}, T_{\mathcal{F}}, \mathcal{N}_\mathcal{F}^* \right)$ over $A$ which induces $\left(\mathcal{X}, T_\mathcal{F} \right) $.
\end{proposition}

\subsection{Stability of a foliated nonsingular variety $\left( X, T_{\mathcal{F}_0}, \mathcal{N}_{\mathcal{F}_0}^* \right)$ in $\left(X, \mathcal{N}_{\mathcal{F}_0}^* \right)$ }\

\begin{definition}
A foliated nonsingular variety $\left(X, T_{\mathcal{F}_0}, \mathcal{N}_{\mathcal{F}_0}^* \right)$ with $T_{\mathcal{F}_0}$ and $\mathcal{N}_{\mathcal{F}_0}^*$ locally free is called stable in $\left(X, \mathcal{N}_{\mathcal{F}_0}^*\right)$ if the morphism of functors from $(\ref{a22})$ by forgetting $T_{\mathcal{F}_0}$
\begin{align*}
\Phi: \textnormal{FDef}_{\left( X, T_{\mathcal{F}_0}, \mathcal{N}_{\mathcal{F}_0}^* \right)}\to \textnormal{FDef}_{\left(X, \mathcal{N}_{\mathcal{F}_0}^*\right)}
\end{align*}
is smooth.
\end{definition}

We note that $\textnormal{FDef}_{\left(X, T_{\mathcal{F}_0}, \mathcal{N}_{\mathcal{F}_0}^* \right)}$ has a semiuniversal formal element from Proposition \ref{e6} and obstruction space $\mathbb{H}^2\left(X,\mathcal{F}_0^\bullet \right)$ from Proposition \ref{si11}. On the other hand, $\textnormal{FDef}_{\left(X, \mathcal{N}_{\mathcal{F}_0}^* \right)}$ has a semiuniversal formal element and obstruction space $\mathbb{H}^2\left( X, \mathcal{N}_{\mathcal{F}_0}^{ * \bullet} \right)$ from Proposition \ref{nb1}. Assume $\mathbb{H}^1\left(X, B^\bullet \right)=0$. Then from the long exact sequence of $(\ref{a20})$, we have the surjection $d\Phi: \mathbb{H}^1\left(X,  \mathcal{F}_0^\bullet \right)\to \mathbb{H}^1\left( X, \mathcal{N}_{\mathcal{F}_0}^{*\bullet}  \right)$ and an injective obstruction map $\mathbb{H}^2\left( X,  \mathcal{F}_0^\bullet \right)\to \mathbb{H}^2\left( X, \mathcal{N}_{\mathcal{F}_0}^{*\bullet} \right)$, in other words, $\textnormal{FDef}_{\left( X, T_{\mathcal{F}_0} ,  \mathcal{N}_{\mathcal{F}_0}^* \right)}$ is less obstructed than $\textnormal{FDef}_{\left(X, \mathcal{N}_{\mathcal{F}_0}^*\right)}$ (see \cite{Ser06} Definition 2.3.5 p.58). This implies that $\Phi$ is smooth. (see \cite{Ser06} Proposition 2.3.6 p.59). We summarize as follows.

\begin{proposition}
Let $\left( X, T_{\mathcal{F}_0},  \mathcal{N}_{\mathcal{F}_0}^* \right)$ be a foliated nonsingular variety with $T_{\mathcal{F}_0}$ and $\mathcal{N}_{\mathcal{F}_0}^*$ locally free. Assume that $\mathbb{H}^1\left(X, B^\bullet \right)=0$. Then $\Phi$ is smooth. In other words, for a flat foliated deformation $\left(\mathcal{X}, \mathcal{N}_{\mathcal{F}}^*\right)$ of $\left(X, \mathcal{N}_{\mathcal{F}_0}^*\right)$ over $A$, there exists a simultaneous foliated deformation $\left(\mathcal{X}, T_{\mathcal{F}}, \mathcal{N}_\mathcal{F}^* \right)$ over $A$ which induces $\left(\mathcal{X}, \mathcal{N}_\mathcal{F}^* \right) $.
\end{proposition}

\subsection{Formal foliated deformations of $\left(X, T_{\mathcal{F}_0} , \mathcal{N}_{\mathcal{F}_0}^* \right)$ with $T_{\mathcal{F}_0}$ and $\mathcal{N}_{\mathcal{F}_0}^*$ locally free}\

We only touch upon formal foliated deformations of a foliated nonsingular variety $\left( X, T_{\mathcal{F}_0}, \mathcal{N}_{\mathcal{F}_0}^* \right)$ with $T_{\mathcal{F}_0}$ and $\mathcal{N}_{\mathcal{F}_0}^*$ locally free in a relation to Part I.  The following Theorem \ref{a25} can be considered as the scheme-theoretic version of Theorem of existence for deformations of foliated complex analytic structures in terms of tangent sheaves in Part I. 

\begin{theorem}[compare \cite{Ser06} Corollary 2.6.4 p.92] \label{a25}
Let $\left(X , T_{\mathcal{F}_0}, \mathcal{N}_{\mathcal{F}_0}^* \right)$ is a foliated nonsingular projective variety with $T_{\mathcal{F}_0}$ and $\mathcal{N}_{\mathcal{F}_0}^*$ locally free. If $\mathbb{H}^0\left( X, \mathcal{F}_0^\bullet \right)=\mathbb{H}^2\left( X, \mathcal{F}_0^\bullet \right)=0$, then $\textnormal{FDef}_{\left( X, T_{\mathcal{F}_0} , \mathcal{N}_{\mathcal{F}_0}^* \right)}$ is pro-represented by a formal power series ring.
\end{theorem}

We can consider a concept of formal foliated deformations of a foliated nonsingular variety $\left(X, T_{\mathcal{F}_0 } , \mathcal{N}_{\mathcal{F}_0}^* \right)$ with $T_{\mathcal{F}_0}$ and $\mathcal{N}_{\mathcal{F}_0}^*$ locally free as an analogue of formal deformations of a nonsingular variety $X$ as in \cite{Ser06} p.77. We note that we can extend the functor $\textnormal{FDef}_{\left(X, T_{\mathcal{F}_0} , \mathcal{N}_{\mathcal{F}_0}^* \right)}$ to the functor
\begin{align*}
\widehat{\textnormal{FDef}}_{\left( X, T_{\mathcal{F}_0} , \mathcal{N}_{\mathcal{F}_0}^* \right)} : \widehat{\bold{Art}}&\to \left(\textnormal{Sets} \right)\\
     \bar{A} &\to \varprojlim \textnormal{FDef}_{\left(X, T_{\mathcal{F}_0} , \mathcal{N}_{\mathcal{F}_0}^* \right)}\left(\bar{A}/\mathfrak{m}_{\bar{A}}^{n+1} \right)
\end{align*}
(see \cite{Ser06} p.46). An element $\hat{\eta}\in \widehat{\textnormal{FDef}}_{\left( X, T_{\mathcal{F}_0} , \mathcal{N}_{\mathcal{F}_0}^* \right)}\left(\bar{A} \right)$ is called a formal element of $\textnormal{FDef}_{\left( X, T_{\mathcal{F}_0} , \mathcal{N}_{\mathcal{F}_0}^* \right)}$. By definition $\hat{\eta}$ can be represented as a system of elements $\left\{ \eta_n \in \textnormal{FDef}_{\left( X, T_{\mathcal{F}_0} , \mathcal{N}_{\mathcal{F}_0}^* \right)}\left(\bar{A}/\mathfrak{m}_{\bar{A}}^{n+1} \right) \right\}_{n\geq 0}$  such that for every $n\geq 1$ the map $\textnormal{FDef}_{\left( X, T_{\mathcal{F}_0} , \mathcal{N}_{\mathcal{F}_0}^* \right)}\left( \bar{A}/\mathfrak{m}_{\bar{A}}^{n+1}\right)\to \textnormal{FDef}_{\left( X, T_{\mathcal{F}_0}  , \mathcal{N}_{\mathcal{F}_0}^* \right)}\left( \bar{A}/\mathfrak{m}_{\bar{A}}^n \right)$ induced by the projection $\bar{A}/\mathfrak{m}_{\bar{A}}^{n+1}\to \bar{A}/{\mathfrak{m}}_{\bar{A}}^n$ sends $\eta_n\mapsto \eta_{n-1}$. A formal foliated deformation of $\left( X, T_{\mathcal{F}_0} , \mathcal{N}_{\mathcal{F}_0}^* \right)$ over $\bar{A}$ is a formal element $\hat{\eta}=\left\{\eta_n \right\} \in \widehat{\textnormal{FDef}}_{\left( X, T_{\mathcal{F}_0} , \mathcal{N}_{\mathcal{F}_0}^* \right)}\left( \bar{A}\right)$ for $\textnormal{FDef}_{\left(X, T_{\mathcal{F}_0} , \mathcal{N}_{\mathcal{F}_0}^* \right)}$, where $\eta_n$ is an infinitesimal foliated deformation of $\left( X, T_{\mathcal{F}_0} , \mathcal{N}_{\mathcal{F}_0}^* \right)$ over $\bar{A}_n$
\begin{center}
$\eta_n:\begin{CD}
\left(X, T_{\mathcal{F}_0} , \mathcal{N}_{\mathcal{F}_0}^* \right) @>>> \left( \mathcal{X}_n, T_{\mathcal{F}_n} , \mathcal{N}_{\mathcal{F}_n}^* \right) \\
@VVV @VV \pi_n V\\
\textnormal{Spec}(k) @>>> \textnormal{Spec}\left( \bar{A}_n \right)
\end{CD}$
\end{center}
where $\bar{A}_n=\bar{A}/\mathfrak{m}_{\bar{A}}^{n+1}$, such that for all $n\geq 1$, $\eta_n$ induces $\eta_{n-1}$ by pullback under the natural inclusion $\textnormal{Spec}\left(\bar{A}_{n-1} \right)\to \textnormal{Spec}\left( \bar{A}_n \right)$, i.e. we have
\begin{center}
$\eta_{n-1}:\begin{CD}
\left(X, T_{\mathcal{F}_0} , \mathcal{N}_{\mathcal{F}_0}^* \right) @>>> \left( \mathcal{X}_{n-1}=\mathcal{X}_n \otimes_{\bar{A}_n} \bar{A}_{n-1}, T_{\mathcal{F}_n}\otimes_{\bar{A}_n} \bar{A}_{n-1} , \mathcal{N}_{\mathcal{F}_n}^*\otimes_{\bar{A}_n} \bar{A}_{n-1} \right)\\
@VVV @VV\pi_{n-1}V\\
\textnormal{Spec}(k) @>>> \textnormal{Spec}\left(\bar{A}_{n-1} \right)
\end{CD}$
\end{center}
We will denote such a formal foliated deformation by $\left( \bar{A}, \left\{ \eta_n \right\} \right)$.

Let $\left(X, T_{\mathcal{F}_0}, \mathcal{N}_{\mathcal{F}_0}^* \right)$ be a foliated nonsingular projective variety with $T_{\mathcal{F}_0}$ and $\mathcal{N}_{\mathcal{F}_0}^*$ locally free and consider a flat family of foliated deformations of $\left( X, T_{\mathcal{F}_0} , \mathcal{N}_{\mathcal{F}_0}^* \right)$ parametrized by an affine scheme $S=\textnormal{Spec}\left( B \right)$ where $B$ is a $k$-algebras of finite type:
\begin{center}
$\eta:\begin{CD}
\left( X, T_{\mathcal{F}_0} , \mathcal{N}_{\mathcal{F}_0}^* \right) @>f>> \left( \mathcal{X}, T_\mathcal{F} , \mathcal{N}_\mathcal{F}^* \right)\\
@VVV @VV\pi V\\
\textnormal{Spec}(k) @>s >> S
\end{CD}$
\end{center}
where $\pi$ is projective, and $\mathcal{X}$ is smooth over $S$, $s$ is a $k$-rational point and $T_\mathcal{F}$ is a locally free subsheaf of $T_{\mathcal{X}/S}$ such that $\frac{T_{\mathcal{X}/S}}{T_\mathcal{F}}$ is flat over $S$ and $\left[ T_{\mathcal{F}}, T_{\mathcal{F}}\right] \subset T_{\mathcal{F}}$ under the Lie bracket, and $\mathcal{N}_\mathcal{F}^*$ is a locally free subsheaf of $\Omega_{\mathcal{X}/S}^1$ such that $\frac{\Omega_{\mathcal{X}/S}^1}{\mathcal{N}_\mathcal{F}^*}$ is flat over $S$ and $d\left( \mathcal{N}_{\mathcal{F}, x}^*  \right) \subset \mathcal{N}_{\mathcal{F} , x}^* \bigwedge \Omega_{\mathcal{X}/S, x}^1$ for  $x\in \mathcal{X}- \textnormal{Sing}\left( \mathcal{F}\right)$ where $\textnormal{Sing}\left( \mathcal{F}\right)$ is the locus of points of $\mathcal{X}$ at which $\frac{\Omega_{\mathcal{X}/S}^1}{\mathcal{N}_\mathcal{F}^*}$ is not locally free, and $\left( \frac{\Omega_{\mathcal{X}/S}^1}{\mathcal{N}_\mathcal{F}^*} \right)^*=T_{\mathcal{F}}$ and $\left( \frac{T_{\mathcal{X}/S}}{T_\mathcal{F}} \right)^*=\mathcal{N}_\mathcal{F}^*$. We will denote the foliated deformation $\eta$ by $\left( S, s , \eta \right)$. Given such a foliated deformation $\left(  S, s, \eta \right)$ let $\eta_n$ be the infinitesimal foliated deformation induced pulling back of $\eta$ under the closed embedding $\textnormal{Spec}\left( \mathcal{O}_{S,s}/\mathfrak{m}^{n+1} \right) \to S$. Then $\left( \widehat{\mathcal{O}}_{S,s}, \left\{ \eta_n \right\} \right)$ is a formal foliated deformation of $\left(X, T_{\mathcal{F}_0} , \mathcal{N}_{\mathcal{F}_0}^* \right)$ which we call the formal foliated deformation defined by $\eta$.

\begin{definition}
A foliated deformation $\left(S, s, \eta \right)$ of a foliated nonsingular variety $\left( X, T_{\mathcal{F}_0} , \mathcal{N}_{\mathcal{F}_0}^* \right)$ with $T_{\mathcal{F}_0}$ and $\mathcal{N}_{\mathcal{F}_0}^*$ locally free is called formally universal $($formally semiuniversal, formally versal$)$ if the formal foliated deformation $\left( \widehat{\mathcal{O}}_{S,s}, \left\{ \eta_n\right\}_{n \geq 0} \right)$ defined by $\eta$ is universal $($resp. semiuniversal, versal$)$ \textnormal{(see \cite{Ser06} Definition 2.2.3 p.47 and Definition 2.2.6 p.49)}. A flat foliated family $\pi:\left(\mathcal{X}, T_\mathcal{F} , \mathcal{N}_\mathcal{F}^* \right)\to S$, where $\mathcal{X}$ is smooth over $S$, and $T_{\mathcal{F}}$ is a locally free subsehaf of $T_{\mathcal{X}/S}$, and $\mathcal{N}_\mathcal{F}^*$ is a locally free subsheaf of $\Omega_{\mathcal{X}/S}^1$, is called formally universal $($resp. formally semiuniversal, formally versal$)$ at a $k$-rational point $s\in S$ if
\begin{center}
$\eta:\begin{CD}
\left(\mathcal{X}(s), T_{\mathcal{F}(s)} , \mathcal{N}_{\mathcal{F}(s)}^* \right) @>\hookrightarrow>> \left( \mathcal{X}, T_\mathcal{F} , \mathcal{N}_\mathcal{F}^*\right)\\
@VVV @VV\pi V \\
\textnormal{Spec}(k) @>s>> S
\end{CD}$
\end{center}
is a formally universal $($resp. formally semiuniversal, formally versal$)$ foliated deformation of $\left( \mathcal{X}(s), T_{\mathcal{F}(s)}\right)$.
\end{definition}

\begin{definition}[foliated Kodaira-Spencer map in an algebraic foliated family]
Let $\left( X, T_{\mathcal{F}_0} \right)$ be a foliated nonsingular variety with $T_{\mathcal{F}_0}$ locally free. Let us consider a foliated deformation $(S, s, \eta)$
\begin{center}
$\eta: \begin{CD}
\left( X, T_{\mathcal{F}_0},  \mathcal{N}_{\mathcal{F}_0}^* \right) @>>> \left( \mathcal{X}, T_\mathcal{F} , \mathcal{N}_\mathcal{F}^* \right)\\
@VVV @VV\pi V\\
\textnormal{Spec}(k) @>s>> S
\end{CD}$
\end{center}
By pulling back $(S, s, \eta)$ by morphisms $\textnormal{Spec}\left( k[\epsilon] \right)\to S$ with image $s$ and applying \textnormal{Proposition \ref{si11}} we define a linear map
\begin{align*}
\kappa_{\pi, s}: T_{S,s}\to  \textnormal{FDef}_{\left( X, T_{\mathcal{F}_0}  , \mathcal{N}_{\mathcal{F}_0}^* \right)} \cong \mathbb{H}^1\left( X, \mathcal{F}_0^\bullet \right)
\end{align*}
which is called the foliated Kodaira-Spencer map of $(S,s, \eta)$ in terms of tangent sheaves.
\end{definition}

\begin{proposition}[compare \cite{Ser06} Proposition 2.5.8 p.79]\label{a26}
Let $\left(S , s , \eta \right)$ be a foliated deformation of a foliated nonsingular variety $\left( X, T_{\mathcal{F}_0}, \mathcal{N}_{\mathcal{F}_0}^* \right)$ with $T_{\mathcal{F}_0}$ and $\mathcal{N}_{\mathcal{F}_0}^*$ locally free.
\begin{enumerate}
\item If $\eta$ is formally versal $($resp. formally semiuniversal or formally universal$)$, then the foliated Kodaira-Spencer map $\kappa_{\pi, s}: T_s S\to \textnormal{FDef}_{\left(X, T_{\mathcal{F}_0} , \mathcal{N}_{\mathcal{F}_0}^* \right)}\left( k[\epsilon]\right)\cong \mathbb{H}^1 \left(X, \mathcal{F}_0^\bullet \right)$ is surjective $($resp. an isomorphism$)$.
\item If $S$ is nonsingular at $s$ and the foliated Kodaira-Spencer map $\kappa_{\pi,s}$ is surjective $($resp. an isomorphism$)$, then $\eta$ is formally versal $($resp. formally semiuniversal$)$ and $\left(X, T_{\mathcal{F}_0} \right)$ is unobstructed, i.e. $\textnormal{FDef}_{\left( X, T_{\mathcal{F}_0} , \mathcal{N}_{\mathcal{F}_0}^* \right)}$ is smooth.
\end{enumerate}
\end{proposition}

Proposition \ref{a26} is the algebraic version of Theorem of completeness for deformations of foliated complex analytic structures in terms of both tangent and cotangent sheaves in Part I. We finish this section by mentioning the concept of algebraizability in foliated deformations.

\begin{definition}[compare \cite{Ser06} p.80]
A formal foliated deformation $\left( \bar{A}, \left\{ \eta_n \right\} \right)$ of a foliated nonsingular variety $\left( X, T_{\mathcal{F}_0}, \mathcal{N}_{\mathcal{F}_0}^*\right)$ with $T_{\mathcal{F}_0}$ and $\mathcal{N}_{\mathcal{F}_0}^*$ locally free is called algebraizable if there exists an algebraic foliated deformation $\left( S, s, \xi \right)$ of $\left( X, T_{\mathcal{F}_0} , \mathcal{N}_{\mathcal{F}_0}^* \right)$ and an isomorphism $\widehat{\mathcal{O}}_{S,s}\cong \bar{A}$ sending $\eta_n\to \xi_n$ for all $n$, i.e. $\left( \bar{A}, \left\{ \eta_n \right\} \right)$ is isomorphic to the formal foliated deformation defined by $\xi$. The foliated deformation $\left( S, s, \xi \right)$ is called an algebraization of $\left( \bar{A} , \left\{ \eta_n \right\}  \right)$.
\end{definition}

The existence of algebraizations can be considered as the counterpart of the convergence step in the construction of local families of foliated deformations of compact foliated complex manifolds in terms of tangent sheaves in Part I.

\section{Relationship between foliated deformations of $\left(X, T_{\mathcal{F}_0}\right)$ of dimension $n-1$ with $T_{\mathcal{F}_0}$ locally free and foliated deformations of $\left(X, \mathcal{N}_{\mathcal{F}_0}^* \right)$ of codimension $1$  }\

Let $\left( X, T_{\mathcal{F}_0} \right)$ be a foliated nonsingular complex variety with $\dim_\mathbb{C} X= n$ and $T_{\mathcal{F}_0}$ be locally free with rank $n-1$. Let $\mathcal{U}=\{U_i\}$ be an Stein open covering of $X$ with local coordinates $\left(z_i^1,..., z_i^n \right)$ on $U_i$ such that $\Gamma\left( U_i, T_{\mathcal{F}_0} \right)$ is generated by $T_{0i}^1,..., T_{0i}^{n-1}\in \Gamma\left( U_i, T_X \right)$ with $T_{0i}^\alpha= \sum_{\beta=1}^{n-1} r_{0ij}^{\alpha\beta} w_{0j}^\beta$ with $r_{0ij}^{\alpha\beta}\in \Gamma\left( U_{ij}, \mathcal{O}_X \right),\alpha=1,...,{n-1}$ and $\left[ T_{0i}^\alpha, T_{0i}^\beta \right]= \sum_{\gamma=1}^{n-1} g_{0i\alpha\beta}^\gamma T_{0i}^\gamma$ for some $g_{0i\alpha\beta}^\gamma \in \Gamma\left( U_i , \mathcal{O}_X \right),\alpha,\beta=1,..., n-1$. Then $T_{0i}^1\wedge \cdots \wedge T_{0i}^{n-1}= r_{0ij} T_{0j}^1 \wedge \cdots \wedge T_{0j}^{n-1}$, where $r_{0ij}=\det\left(r_{0ij}^{\alpha\beta} \right)$. Then $\mathcal{Q}_0=\left(\bigwedge^p T_{\mathcal{F}_0} \right)^*$ is defined by the transition functions $\{r_{0ij}\}$ with $r_{0ij}\in \Gamma(U_{ij}, \mathcal{O}_X^*)$ and $T_0$ is defined locally on $U_i$ by $T_0:=\left\{ T_{0i}^1\wedge \cdots \wedge T_{0i}^{n-1} \right\}\in H^0\left(X, \bigwedge^{n-1} T_X \otimes \mathcal{Q}_0 \right)$ as in section \ref{nt1}.

We will describe the associated codimension $1$-foliation $\mathcal{N}_{\mathcal{F}_0}^*$. Since $\bigwedge^{n-1} T_X \times T_X \to \bigwedge^n T_X=K_X^*$ is a perfect pairing, we have
\begin{align*}
\bigwedge^{n-1} T_X \cong \mathscr{H}om_{\mathcal{O}_X}\left(T_X, K_X^* \right)\cong \Omega_X^1\otimes K_X^*
\end{align*}
so that we have
\begin{align*}
\bigwedge^{n-1} T_X \otimes \mathcal{Q}_0 &\cong \Omega_X^1\otimes K_X^* \otimes \mathcal{Q}_0 \\
 T_{0i}^1\wedge \cdots \wedge T_{0i}^{n-1} &\mapsto \left(\frac{\partial}{\partial z_i^\alpha}\mapsto  T_{0i}^1\wedge \cdots \wedge T_{0i}^{n-1}\wedge \frac{\partial}{\partial z_i^\alpha}=a_i^\alpha \frac{\partial}{\partial z_i^1}\wedge \cdots  \wedge \frac{\partial}{\partial z_i^n}  \right) \cong \sum_{\alpha=1}^n a_i^\alpha dz_i^\alpha\,\,\,\,\,\,\,\textnormal{for}\,\,\, a_i^\alpha \in \Gamma\left( U_i, \mathcal{O}_X \right)
\end{align*}

We note that
\begin{align*}
i_{T_{0i}^1\wedge \cdots \wedge T_{0i}^{n-1} \wedge \frac{\partial}{\partial z_i^\alpha}} \left(dz_i^1\wedge \cdots \wedge dz_i^n \right)= a_i^\alpha \iff i_{\frac{\partial}{\partial z_i^\alpha}}\left(i_{T_{0i}^1\wedge \cdots \wedge T_{0i}^{n-1}}\left( dz_i^1\wedge \cdots \wedge dz_i^n \right) \right)= a_i^\alpha
\end{align*}

This implies that
\begin{align*}
i_{T_{0i}^1\wedge \cdots \wedge T_{0i}^{n-1}}\left(  dz_i^1\wedge \cdots \wedge dz_i^n \right) = \sum_{\alpha=1}^n a_i^\alpha dz_i^\alpha : = w_{0i}
\end{align*}

Then we see that
\begin{align*}
i_{T_{0i}^\beta}\left( w_{0i} \right) = i_{T_{0i}^\beta}\left(\sum_{\alpha=1}^n a_i^\alpha dz_i^\alpha \right) = i_{T_{0i}^\beta}i_{T_{0i}^1\wedge \cdots \wedge T_{0i}^{n-1} \wedge \frac{\partial}{\partial z_i^\alpha}} \left(dz_i^1\wedge \cdots \wedge dz_i^n \right)=0 
\end{align*}

On the other hand, we note that
\begin{align*}
dz_i^1\wedge \cdots \wedge dz_i^n =K_{0ij} dz_j^1\wedge \cdots \wedge dz_j^n \,\,\,\,\,\,\,\,\textnormal{for some}\,\,\,K_{0ij}\in \Gamma\left( U_{ij}, \mathcal{O}_X^* \right)
\end{align*}
Then $K_X^*$ is defined by the transition functions $\left\{K_{0ij} \right\}\in C^1\left( \mathcal{U},  \mathcal{O}_X^*\right)$. Then
\begin{align*}
w_{0i}= i_{T_{0i}^1\wedge \cdots \wedge T_{0i}^{n-1}}\left( dz_i^1\wedge \cdots \wedge dz_i^n\right)= r_{0ij} K_{0ij} i_{T_{0j}^1\wedge \cdots \wedge T_{0j}^{n-1}}\left(dz_j^1\wedge \cdots \wedge dz_j^n \right)= r_{0ij} K_{0ij} w_{0j}
\end{align*}
This implies that
\begin{align*}
w=\{w_{0i}\}\in H^0\left(X,  \Omega_X^1\otimes K_X^*\otimes \mathcal{Q}_0 \right) = H^0\left( X , \Omega_X^1 \otimes \mathcal{L}_0 \right),\,\,\,\,\,\,\,\textnormal{where}\,\,\,\mathcal{L}_0:= K_X^*\otimes \mathcal{Q}_0 = \left( \mathcal{N}_{\mathcal{F}_0}^* \right)^*
\end{align*}
Moreover, for $\alpha, \beta= 1,..., n-1$,
\begin{align*}
 i_{\left[T_{0i}^\alpha, T_{0i}^\beta\right]}\left(w_{0i} \right) =  i_{\sum_{\gamma=1}^{n-1} g_{0i\alpha\beta}^\gamma T_{0i}^\gamma}i_{T_{0i}^1\wedge \cdots \wedge T_{0i}^{n-1}}(dz_i^1\wedge \cdots \wedge dz_i^n)= 0
\end{align*}
Then we have 
\begin{align} \label{nt2}
0 = i_{\left[T_{0i}^\alpha ,  T_{0i}^\beta \right]}(w_{0i})=\mathcal{L}_{T_{0i}^\alpha} i_{T_{0i}^\beta}(w_{0i}) - i_{T_{0i}^\beta} \mathcal{L}_{T_{0i}^\alpha}(w_{0i}) = -i_{T_{0i}^\beta}\left(i_{T_{0i}^\alpha} d + d i_{T_{0i}^\alpha}    \right)(w_{0i})= - i_{T_{0i}^\alpha\wedge T_{0i}^\beta}(dw_{0i})
\end{align}
For $x\in U_i-S$, there exists a neighborhood $U_x \subset U_i-S$ of $x$ such that $T_{0i}^1,..., T_{0i}^{n-1}, V_{i_x}^n$ is a basis of $\Gamma(U_x, T_X)$ for some $V_{i_x}^n  \in \Gamma\left( U_x, T_X \right)$. Then
\begin{align*}
i_{V_{i_x}^n}\left( w_{0i} \right) = i_{V_{i_x}^n} i_{T_{0i}^1\wedge \cdots \wedge T_{0i}^n}\left( dz_i^1\wedge \cdots \wedge dz_i^n\right)= i_{b_i \frac{\partial}{\partial z_i^1}\wedge \cdots \wedge \frac{\partial}{\partial z_i^n}}\left( dz_i^1 \wedge \cdots \wedge dz_i^n \right) = b_i \,\,\,\,\,\textnormal{for some}\,\,\, b_i\in \Gamma\left( U_i, \mathcal{O}_X^* \right) 
\end{align*}
This implies that $w_{0i}$ does not vanish on $U_x$. Then by considering the dual basis of $T_{0i}^1,..., T_{0i}^{n-1}, V_{i_x}^n$ on $U_x$ which extend $w_{0i}$, we can deduce from $(\ref{nt2})$ that $dw_{0i}= a_{i_x}\wedge w_{0i}$ on $U_x$ for some $a_{i_x}\in \Gamma(U_{i_x}, \Omega_X^1)$. Then $w_{0i}\wedge dw_{0i}=0$ on $U_i-S$, so that $w_{0i}\wedge d w_{0i} = 0$.

\subsection{Foliated deformations of $\left(X, T_{\mathcal{F}_0}\right)$ with $T_{\mathcal{F}_0}$ locally free of rank $n-1$ in terms of tangent sheaves and foliated deformations of $\left( X, \mathcal{N}_{\mathcal{F}_0}^*\right)$ with $\mathcal{N}_{\mathcal{F}_0}^*$ locally free of rank $1$ in terms of cotangent sheaves } \

We keep the notations in the previous subsection. We will show that a foliated deformation $(\mathcal{X}, T_{\mathcal{F}})$ of $(X, T_{\mathcal{F}_0})$ with $T_{\mathcal{F}_0}$ locally free of rank $n-1$ over $A\in \bold{Art}$ in terms of tangent sheaves also defines a foliated deformation $(\mathcal{X}, \mathcal{N}_\mathcal{F}^*)$ of $(X, \mathcal{N}_{\mathcal{F}_0}^*)$ of codimension $1$ in terms of cotangent sheaves, equivalently in terms of deformations of the section $w_0$. Let $\mathcal{U}=\left\{ U_i \right\}$ be a Stein open covering of $X$ with local coordinates $\left( z_i^1,..., z_i^n \right)$ on $U_i$ such that for all $i$, we have an isomorphism
\begin{align*}
\theta_i: (U_i\times \textnormal{Spec}(A), T_{\mathcal{F}_i})\to (\mathcal{X}|_{U_i}, T_\mathcal{F}|_{U_i})
\end{align*}
and $\Gamma(U_i, T_{\mathcal{F}_i})$ is generatesecd by $T_i^\alpha \in \Gamma(U_i, T_X)\otimes A$ which induces $T_{0i}^\alpha, \alpha=1,...,n-1$, and $\theta_{ij}=\theta_j^{-1}\theta_i$ is a foliated isomorphism with $\theta_{jk}\theta_{ij}= \theta_{ik}$ and
\begin{align*}
\theta_{ji}^*(T_i^\alpha)=\sum_{\beta=1}^{n-1} r_{ij}^{\alpha\beta} T_j^\beta,\,\,\,\,\,\,\,\,\,r_{ij}^{\alpha\beta} \in \Gamma(U_{ij}, \mathcal{O}_X)\otimes A
\end{align*}
Let $\mathcal{Q}= \left(\bigwedge^{n-1}T_\mathcal{F}\right)^*$. Then $\left\{T_i^1\wedge \cdots \wedge T_i^{n-1}\right\}\in H^0\left(\mathcal{X}, \bigwedge^{n-1} T_{\mathcal{X}/A}\otimes \mathcal{Q} \right)$. We note that
\begin{align*}
\theta_{ji}^*\left( T_i^1\wedge \cdots \wedge T_i^{n-1} \right) =\theta_{ji}^*(T_i^1)\wedge \cdots \wedge \theta_{ji}^*(T_i^{n-1})= r_{ij} T_j^1\wedge \cdots \wedge T_j^{n-1},\,\,\,\,\,\,\, r_{ij}=  \det\left( r_{ij}^{\alpha\beta} \right)
\end{align*}
We note that since $\bigwedge^{n-1} T_{\mathcal{X}/A}\times T_{\mathcal{X}/A}\to \bigwedge^n T_{\mathcal{X}/A}\cong K_{\mathcal{X}/A}^*$ is a perfect pairing, we have
\begin{align*}
\bigwedge^{n-1} T_{\mathcal{X}/A} \cong \mathscr{H}om_{\mathcal{O}_\mathcal{X}}\left( T_{\mathcal{X}/A}, K_{\mathcal{X}/A}^* \right) \cong \Omega_{\mathcal{X}/A}^1\otimes K_{\mathcal{X}/A}^*
\end{align*}
so that we have
\begin{align*}
\bigwedge^{n-1} T_{\mathcal{X}/A}\otimes \mathcal{Q} \cong \Omega_{\mathcal{X}/A}^1\otimes K_{\mathcal{X}/A}^*\otimes \mathcal{Q} 
\end{align*}

Locally on $U_i$, we have an isomorphism
\begin{align*}
\Gamma\left(U_i, \bigwedge^{n-1} T_X \right)\otimes A &\cong \Gamma(U_i, \Omega_X^1)\otimes A \\
                         T_i^1\wedge \cdots \wedge T_i^{n-1} &\mapsto i_{T_i^1\wedge \cdots \wedge T_i^{n-1}} (dz_i^1\wedge \cdots \wedge dz_i^n)=\sum_{\alpha=1}^n A_i^\alpha dz_i^\alpha = w_i ,\,\,\,\,\,\,\,\,\,A_i^\alpha\in \Gamma(U_i, \mathcal{O}_X)\otimes A
\end{align*}
Then we have $\left\{ w_i = \sum_{\alpha=1}^n A_i^\alpha dz_i^\alpha \right\} \in H^0\left( \mathcal{X}, \Omega_X^1\otimes K_{\mathcal{X}/A}^*\otimes \mathcal{Q} \right)$ which induces $\left\{ w_{0i} = \sum_{\alpha=1}^n a_i^\alpha dz_i^\alpha \right\}$. This is well-defined. In fact let $f_{ji}: \left(\Gamma\left( U_{ii}, \mathcal{O}_X \right), T_{\mathcal{F}_i} \right) \to  \left( \Gamma\left(U_{ij}, \mathcal{O}_X \right), T_{\mathcal{F}_j} \right) $ which corresponds to $\theta_{ji}:\left(U_{ij}, T_{\mathcal{F}_j}\right) \to \left(U_{ij}, T_{\mathcal{F}_i}\right)$. $K_{\mathcal{X}/A}^*$ is described by the transition functions $\left\{ K_{ij} \right\}\in C^1 \left( \mathcal{U}, \mathcal{O}_\mathcal{X}^* \right)$ in terms of $dz_i^1\wedge \cdots \wedge dz_i^n$ as
\begin{align*}
K_{ij} dz_j^1\wedge \cdots \wedge dz_j^n= d\left( f_{ji}(z_i^1) \right) \wedge \cdots \wedge d \left( f_{ji}(z_i^n) \right)
\end{align*}
Then
\begin{align*}
\theta_{ji}^*\left( i_{T_i^1\wedge \cdots \wedge T_i^{n-1}} \left(dz_i^1\wedge \cdots \wedge dz_i^n \right) \right)= i_{\theta_{ji}^* T_i^1\wedge \cdots \wedge \theta_{ji}^* T_i^{n-1}}\left( df_{ji}(z_i^1)\wedge \cdots df_{ij}(z_i^n)   \right) = r_{ij} K_{ij}  i_{T_j^1\wedge \cdots \wedge T_j^{n-1}} \left( dz_j^1\wedge \cdots \wedge d z_j^n \right)
\end{align*}
We also note that for $\alpha,\beta=1,...,n-1$,
\begin{align*}
0=i_{\left[ T_i^\alpha, T_i^\beta \right]}\left( w_i\right) = i_{T_i^\alpha \wedge T_i^\beta}\left( dw_i \right)  \Longrightarrow  w_i \wedge dw_i = 0
\end{align*}

On the other hand, let $\mathcal{N}_\mathcal{F}^*$ be the locally free subsheaf of $\Omega_{\mathcal{X}/A}^1$ on $\mathcal{X}$ defined by $\left\{ w=\sum_{\alpha=1}^n A_i^\alpha dz_i^\alpha \right\}$, i.e. defined by the transition function $\left\{ \frac{1}{r_{ij} K_{ij}} \right\}$. Then we have the exact sequence $0\to \mathcal{N}_\mathcal{F}^* \to \Omega_{\mathcal{X}/A}^1 \to \frac{\Omega_{\mathcal{X}/A}^1}{\mathcal{N}_\mathcal{F}^*}\to 0$. We note that since $\mathcal{N}_{\mathcal{F}}^*$ and $\Omega_{\mathcal{X}/A}^1$ are locally free, they are flat over $A$, and the exact sequence induces an injection $0\to \mathcal{N}_{\mathcal{F}_0}^* \to \Omega_X^1$ since $\mathcal{N}_{\mathcal{F}_0}^*$ is the locally free sheaf defined by $\left\{ w_{0i}=\sum_{\alpha=1}^n a_i^\alpha dz_i^\alpha \right\}$. This implies that $\frac{\Omega_{\mathcal{X}/A}^1}{\mathcal{N}_\mathcal{F}^*}$ is flat over $A$. It remains to show that $\mathcal{N}_\mathcal{F}^* =\left( \frac{T_{\mathcal{X}/A}}{T_{\mathcal{F}}} \right)^*$. First we claim that 
\begin{align}\label{nt4}
T_\mathcal{F}=\left( \frac{\Omega_{\mathcal{X}/A}^1}{\mathcal{N}_\mathcal{F}^*} \right)^*.
\end{align}
 In fact, we have the exact sequence $0\to T_\mathcal{F} \to T_{\mathcal{X}/A}\to \frac{T_{\mathcal{X}/A}}{T_\mathcal{F}}\to 0$, and $\frac{T_{\mathcal{X}/A}}{T_{\mathcal{F}}}$ is flat over $A$, and since $\frac{T_X}{T_{\mathcal{F}_0}}$ is locally free on $X-S$, we see that $\frac{T_{\mathcal{X}/A } }{T_\mathcal{F}}$ is locally free on $X-S$. In particular on $U_i-S$, we have the exact sequence of locally free sheaves
\begin{align}\label{nt3}
0 \to T_\mathcal{F}|_{U_i-S} \to T_{\mathcal{X}/A}|_{U_i-S} \to \frac{T_{\mathcal{X}/A}}{T_\mathcal{F}}|_{U_i-S}\to 0
\end{align}
We see that for $x\in U_i-S$, there is a neighborhood $U_x\subset U_i-S$ such that $\Gamma\left(U_x, \frac{T_{\mathcal{X}/A}}{T_\mathcal{F}} \right)$ is generated by $\bar{V}_x$ such that $V_x \in \Gamma(U_x, T_X)\otimes A$ where $\bar{V_x}$ is the image of $V_x$ under $T_{\mathcal{X}/A}\to \frac{T_{\mathcal{X}/A}}{T_\mathcal{F}}$. Then $T_{\mathcal{X}/A}|_{U_x}$ is generated by $T_i^1,..., T_i^{n-1}, V_x$. Let $P\in \Gamma\left( U_i, T_{\mathcal{X}/A}\right)= \Gamma\left( U_i, T_X \right)\otimes A$ with $i_{P}\left( \sum_{\alpha=1}^n A_i^\alpha dz_i^\alpha \right)=0$, i.e. $P \in \left(\frac{\Omega_{\mathcal{X}/A}^1}{\mathcal{N}_\mathcal{F}^*} \right)^*$. We can write $P=\sum_{\alpha=1}^{n-1} b_x^\alpha T_i^\alpha + a_x V_x$ on $U_x$. Then
\begin{align*}
0 &= i_{P}\left( \sum_{\alpha=1}^n A_i^\alpha dz_i^\alpha \right)= i_{T_i^1\wedge \cdots \wedge T_i^{n-1}\wedge P}\left(dz_i^1\wedge \cdots \wedge dz_i^n \right)=0 \Longrightarrow T_i^1\wedge \cdots \wedge T_i^{n-1} \wedge P=0 \\
&\Longrightarrow a_x=0  \Longrightarrow P\in \Gamma\left(U_x , T_\mathcal{F} \right)
\end{align*}
Then $P\in \Gamma(U_i-S, T_\mathcal{F})$. We note that $\Gamma\left(U_i, T_\mathcal{F} \right)\cong \bigoplus^{n-1} \Gamma\left(U_i, \mathcal{O}_X \right)\otimes A$, so that $\Gamma\left(U_i-S, T_\mathcal{F} \right)\cong \bigoplus^{n-1} \Gamma\left(U_i-S, \mathcal{O}_X \right)\otimes A= \bigoplus^{n-1} \Gamma\left(U_i, \mathcal{O}_X \right)\otimes A \cong \Gamma(U_i, T_\mathcal{F})$. Hence $P\in\Gamma(U_i, T_\mathcal{F})$. This implies that $\left(\frac{\Omega_{\mathcal{X}/A}^1}{\mathcal{N}_\mathcal{F}^*}\right)^* \subset T_\mathcal{F}$. Clearly we have $T_\mathcal{F}\subset \left( \frac{\Omega_{\mathcal{X}/A}^1}{\mathcal{N}_\mathcal{F}^*} \right)^*$ since $i_{T_i^1\wedge \cdots \wedge T_i^{n-1}\wedge V}(dz_i^1\wedge \cdots \wedge dz_i^n)=0$ for $V= \sum_{\alpha=1}^{n-1} a_i^\alpha T_i^\alpha$. Hence $T_\mathcal{F}=\left(\frac{\Omega_{\mathcal{X}/A}^1}{\mathcal{N}_\mathcal{F}^*} \right)^*$. This proves the claim.

By definition of $\mathcal{N}_{\mathcal{F}}^*$, we have $\mathcal{N}_\mathcal{F}^*\subset \left(\frac{T_{\mathcal{X}/A}}{T_\mathcal{F}} \right)^*$. Conversely, assume that $W\in \Gamma\left(U_i,  \Omega_{\mathcal{X}/A}^1 \right)$ such that $i_{T_i^\alpha} (W)=0$ for $\alpha=1,...,p$. In other words, $W\in \left( \frac{T_{\mathcal{X}/A}}{T_\mathcal{F}} \right)^* \subset \Omega_{\mathcal{X}/A}^1$. Then we will show that $W\in \Gamma\left( U_i, \mathcal{N}_\mathcal{F}^* \right)$. In fact, since $\frac{\Omega_{\mathcal{X}/A}^1}{\mathcal{N}_\mathcal{F}^*}$ is flat over $A$ and $\frac{\Omega_X^1}{\mathcal{N}_{\mathcal{F}_0}^*}$ is locally free, we see that $\frac{\Omega_{\mathcal{X}/A}^1}{\mathcal{N}_{\mathcal{F}}^* }$ is locally free on $U_i-S$. Then from $(\ref{nt4})$ and $(\ref{nt3})$, we have the exact sequence
\begin{align*}
0 \to \left( \frac{T_{\mathcal{X}/A}}{ T_\mathcal{F}} \right)^* |_{U_i-S} \to \Omega_{\mathcal{X}/A}^1|_{U_i-S}   \to T_\mathcal{F}^*|_{U_i-S} = \frac{\Omega_{\mathcal{X}/A}^1}{\mathcal{N}_\mathcal{F}^*} |_{U_i-S} \to 0
\end{align*}
which shows that $\mathcal{N}_\mathcal{F}^*|_{U_i-S} = \left(\frac{T_{\mathcal{X}/A}}{T_\mathcal{F}}\right)^* |_{U_i-S}$. Then $W\in \Gamma\left( U_i-S, \mathcal{N}_\mathcal{F}^* \right)$, so that $W\in \Gamma\left( U_i, \mathcal{N}_\mathcal{F}^*  \right)$. This implies that $\left( \frac{T_{\mathcal{X}/A}}{T_\mathcal{F}} \right)^*\subset \mathcal{N}_\mathcal{F}^*$. Hence we have $\mathcal{N}_\mathcal{F}^* =\left( \frac{T_{\mathcal{X}/A}}{T_\mathcal{F}} \right)^*$.

Hence a foliated deformation $\left(\mathcal{X}, T_{\mathcal{F}} \right)$ of $\left(X, T_{\mathcal{F}_0}\right)$ of dimension $n-1$ also defines a foliated deformation $\left(\mathcal{X}, \mathcal{N}_\mathcal{F}^* \right)$ of $\left(X, \mathcal{N}_{\mathcal{F}_0}^*\right)$ of codimension $1$, so that we have a morphism of functors
\begin{align*}
\textnormal{FDef}_{\left(X, T_{\mathcal{F}_0}\right)} \to  \textnormal{FDef}_{\left(X, \mathcal{N}_{\mathcal{F}_0}^* \right)}
\end{align*}

Now we describe the relationship between $\textnormal{FDef}_{\left(X, T_{\mathcal{F}_0}\right)} $ and $\textnormal{FDef}_{\left(X, \mathcal{N}_{\mathcal{F}_0}^* \right)}$. We shall describe a morphism of first-order deformations and obstructions between two functors. We  recall section \ref{nt1}. We keep the notations in section \ref{nt1}.  We note that foliated deformations of $\left( X, T_{\mathcal{F}_0} \right)$ is controlled by the following complex of sheaves in terms of deformations of the section $T_0:=\left\{T_{0i}^1\wedge \cdots \wedge T_{0i}^{n-1} \right\} \in H^0\left( X, \bigwedge^{n-1} T_X \otimes \mathcal{Q}_0\right)$
\begin{align*}
\mathcal{Q}_0^\bullet : \mathcal{E}_{\mathcal{Q}_0} \xrightarrow{\tilde{D}_0 }\mathcal{R} \xrightarrow{\tilde{D}_1} \mathscr{H}om_{\mathcal{O}_X}\left( \bigwedge^2 T_{\mathcal{F}_0}, \frac{T_X}{T_{\mathcal{F}_0}} \right) \xrightarrow{D_2} \mathscr{H}om_{\mathcal{O}_X}\left( \bigwedge^3 T_{\mathcal{F}_0}, \frac{T_X}{T_{\mathcal{F}_0}} \right) \xrightarrow{D_3} \cdots
\end{align*}

On the other hand, we recall section \ref{nt5}. We keep the notations in section \ref{nt5}. We note that foliated deformations of $\left( X, \mathcal{N}_{\mathcal{F}_0}^* \right)$ is controlled by the complex of sheaves in terms of deformations of the section $w_0\in H^0\left( X, \Omega_X^1 \otimes \mathcal{L}_0 \right)$ where $\mathcal{L}_0= K_X^*\otimes \mathcal{Q}_0$
\begin{align}\label{nt10}
\mathcal{E}_{\mathcal{L}_0}^\bullet : \mathcal{E}_{K_X^*\otimes \mathcal{Q}_0}\xrightarrow{\tilde{E}_0} \tilde{\mathcal{S}}^1= \Omega_X^1\otimes \mathcal{Q}_0\otimes K_X^* \xrightarrow{\tilde{E}_1} \mathscr{H}om_{\mathcal{O}_X}\left( \mathcal{N}_{\mathcal{F}_0}^* , \tilde{\mathcal{S}}^2 \right) \xrightarrow{E_2} \mathscr{H}om_{\mathcal{O}_X}\left( \mathcal{N}_{\mathcal{F}_0}^*, \tilde{\mathcal{S}}^3 \right) \xrightarrow{ E_3} \cdots 
\end{align}

We also recall section \ref{nt9}. We note that since $\tilde{\mathcal{S}}^r \cong \bigwedge^r T_{\mathcal{F}_0}^*$, we see that from $( \ref{si41})$, the complex of sheaves $\mathcal{E}_{\mathcal{L}_0}^\bullet$ (\ref{nt10}) is isomorphic to the complex of sheaves which we still denote by $\mathcal{E}_{K_X^*\otimes \mathcal{Q}_0}^\bullet$.
\begin{align} \label{nt11}
\mathcal{E}_{\mathcal{Q}_0\otimes K_X^*}^\bullet: \mathcal{E}_{\mathcal{Q}_0\otimes K_X^*} \xrightarrow{\tilde{E}_0} \tilde{S}^1=\Omega_X^1\otimes \mathcal{Q}_0\otimes K_X^* \xrightarrow{\alpha_2\circ \tilde{E}_1} \mathscr{H}om_{\mathcal{O}_X}\left( \mathcal{N}_{\mathcal{F}_0}^*, \bigwedge^2 T_{\mathcal{F}_0}^* \right) \xrightarrow{ E_2'} \mathscr{H}om_{\mathcal{O}_X}\left( \mathcal{N}_{\mathcal{F}_0}^* , \bigwedge^3 T_{\mathcal{F}_0}^* \right) \xrightarrow{E_3'} \cdots
\end{align}

Then relationship between $\textnormal{FDef}_{\left(X, T_{\mathcal{F}_0} \right)}$ and $\textnormal{FDef}_{\left( X, \mathcal{N}_{\mathcal{F}_0}^*\right)}$ is described by the following exact sequence of complexes of sheaves by defining a morphism $\mathcal{E}_{\mathcal{Q}_0}^\bullet\to \mathcal{E}_{\mathcal{L}_0}^\bullet$.
\begin{center}
$\begin{CD}
@. \cdots @. \cdots @. \cdots \\
@. @AD_3AA @AE_3'AA @AAA \\
0@>>> \mathscr{H}om_{\mathcal{O}_X}\left( \bigwedge^3 T_{\mathcal{F}_0}, \frac{T_X}{T_{\mathcal{F}_0}} \right) @> - \textnormal{dual}>> \mathscr{H}om_{\mathcal{O}_X}\left( \mathcal{N}_{\mathcal{F}_0}^*, \bigwedge^3 T_{\mathcal{F}_0}^* \right) @>>> \mathcal{P}^3 @>>> 0 \\
@. @AD_2AA @AE_2'AA @AAA \\
0 @>>> \mathscr{H}om_{\mathcal{O}_X}\left(\bigwedge^2 T_{\mathcal{F}_0}, \frac{T_X}{T_{\mathcal{F}_0}} \right) @> - \textnormal{dual}>> \mathscr{H}om_{\mathcal{O}_X}\left(\mathcal{N}_{\mathcal{F}_0}^*, \bigwedge^2 T_{\mathcal{F}_0}^* \right) @>>> \mathcal{P}^2 @>>> 0 \\
@. @A\tilde{D}_1AA @A\tilde{E}_1AA @AAA \\
0@>>> \mathcal{R}\subset \bigwedge^{n-1}T_X \otimes \mathcal{Q}_0 @>i_{(-)}(dz_i^1\wedge \cdots \wedge dz_i^n)>> \Omega_X^1 \otimes  \mathcal{Q}_0\otimes K_X^* @>>> \mathcal{P}^1 @>>> 0 \\
@. @A\tilde{D}_0AA @A\tilde{E}_0AA @AAA \\
0 @>>> \mathcal{E}_{\mathcal{Q}_0} @>\gamma>> \mathcal{E}_{\mathcal{Q}_0\otimes K_X^*} @>>> 0 @>>> 0
\end{CD}$
\end{center}
We will show that the above diagram commutes and is exact. First we note that from the injection
\begin{align*}
0\to \frac{T_X}{T_{\mathcal{F}_0}} \to \left( \mathcal{N}_{\mathcal{F}_0}^* \right)^* \to \mathcal{D} \to 0
\end{align*}
we have the exact sequence by tensoring the locally free sheaf $\bigwedge^r T_{\mathcal{F}_0}^* (r\geq 2)$
\begin{align*}
0\to \bigwedge^r T_{\mathcal{F}_0}^* \otimes \frac{T_X}{T_{\mathcal{F}_0}}\cong \mathscr{H}om_{\mathcal{O}_X}\left(\bigwedge^r T_{\mathcal{F}_0}, \frac{T_X}{T_{\mathcal{F}_0}} \right) \to \bigwedge^r T_{\mathcal{F}_0}^*\otimes \left( \mathcal{N}_{\mathcal{F}_0}^* \right)^* \cong \mathscr{H}om_{\mathcal{O}_X}\left( \mathcal{N}_{\mathcal{F}_0}^*, \bigwedge^r T_{\mathcal{F}_0}^* \right) \to \bigwedge^r T_{\mathcal{F}_0}^*\otimes \mathcal{D}=\mathcal{P}^r \to 0
\end{align*}

 We define $\gamma:\mathcal{E}_{\mathcal{Q}_0}\to \mathcal{E}_{\mathcal{Q}_0\otimes K_X^*}$ locally in the following way:
\begin{align*}
\mathcal{O}_X|_{U_i}\bigoplus T_X|_{U_i} &\xrightarrow{\gamma} \mathcal{O}_X |_{U_i} \bigoplus T_X|_{U_i}\\
\left(a_i, b_i=\sum_{\gamma=1}^n b_i^\gamma \frac{\partial}{\partial z_i^\gamma} \right) &\mapsto \left(a_i -\sum_{\alpha=1}^n \frac{\partial b_i^\alpha}{\partial z_i^\alpha}  , b_i  \right)
\end{align*}
We show that $\gamma$ is well-defined. In other words, it is independent of local coordinates, i.e. the following diagram commutes
\begin{center}
$\begin{CD}
(a_i +[b_i, \log r_{0ij}],     b_i)\in \mathcal{O}_X|_{U_j} \bigoplus T_X|_{U_j} @>\gamma >>  (a_i'+ [b_i, \log r_{0ij}] + \left[b_i, \log K_{ij}\right] )\in \mathcal{O}_X|_{U_j} \bigoplus T_X|_{U_j} \\
@AAA @AAA \\
(a_i, b_i)\in\mathcal{O}_X|_{U_i} \bigoplus T_X|_{U_i} @>\gamma>> (a_i', b_i)\in \mathcal{O}_X|_{U_i}\bigoplus T_X|_{U_i}
\end{CD}$
\end{center}
We set that $b_i=\sum_{\gamma=1}^n b_i^\gamma\frac{\partial}{\partial z_i^\gamma}= \sum_{\gamma, \beta=1}^n b_i^\gamma \frac{\partial z_j^\beta}{\partial z_i^\gamma} \frac{\partial}{\partial z_j^\beta}= \sum_{\beta=1}^n b_j^\beta \frac{\partial}{\partial z_j^\beta}$. Then we have to show that 
{\small{\begin{align}\label{nt8}
a_i-\sum_{\alpha=1}^n \frac{\partial b_i^\alpha}{\partial z_i^\alpha} +\left[ b_i , \log r_{0ij} \right] +\left[b_i , \log K_{0ij}\right]  = a_i+\left[ b_i,  \log r_{0ij} \right]  - \sum_{\beta=1}^n \frac{\partial b_j^\beta}{\partial z_j^\beta} \iff -\sum_{\alpha=1}^n \frac{\partial b_i^\alpha}{\partial z_i^\alpha} + \left[ b_i, \log K_{0ij} \right]=-\sum_{\beta=1}^n \frac{\partial b_j^\beta}{\partial z_j^\beta}
\end{align}}}
We note that
{\small{\begin{align*}
d\left(i_{b_i}\left(dz_i^1\wedge \cdots \wedge dz_i^n \right) \right)&=d\left(\sum_{\alpha=1}^n (-1)^{\alpha-1} b_i^\alpha  dz_i^1\wedge \cdots \wedge \widehat{dz_i^\alpha} \wedge \cdots \wedge dz_i^n\right)=\sum_{\alpha=1}^n (-1)^{\alpha-1} \frac{\partial b_i^\alpha}{\partial  z_i^\alpha} dz_i^\alpha \wedge dz_i^1 \wedge \cdots \wedge \widehat{d z_i^\alpha}\wedge \cdots \wedge dz_i^n\\
&= \left(\sum_{\alpha=1}^n  \frac{\partial b_i^\alpha}{\partial z_i^\alpha} \right) dz_i^1\wedge \cdots \wedge dz_i^n
\end{align*}}}
Then we have
\begin{align*}
&\left(\sum_{\alpha=1}^n  \frac{\partial b_i^\alpha}{\partial z_i^\alpha} \right) dz_i^1\wedge \cdots \wedge dz_i^n= d\left(i_{b_i}\left(dz_i^1\wedge \cdots \wedge dz_i^n \right) \right)=d\left(K_{ij} i_{b_i}\left(dz_j^1\wedge \cdots \wedge dz_j^n \right) \right)\\
&= dK_{ij}\wedge i_{b_i}\left(dz_j^1\wedge \cdots \wedge dz_j^n \right) + K_{ij} \left(\sum_{\beta=1}^n \frac{\partial b_j^\beta}{\partial z_j^\beta} \right)dz_j^1\wedge \cdots \wedge dz_j^n\\
&=\left( \sum_{\beta=1}^n \frac{\partial K_{ij}}{\partial z_j^\beta} dz_j^\beta \right) \wedge \left(\sum_{\alpha=1}^n (-1)^{\alpha-1} b_j^\alpha dz_i^1 \wedge \cdots \wedge \widehat{dz_j^\alpha} \wedge \cdots \wedge dz_j^n \right)+ \left( \sum_{\beta=1}^n \frac{\partial b_j^\beta}{\partial z_j^\beta} \right) dz_i^1\wedge \cdots \wedge dz_i^n\\
&= \sum_{\alpha=1}^n b_j^\alpha\frac{\partial K_{ij}}{\partial z_j^\alpha}dz_j^1\wedge \cdots \wedge dz_j^n +\left(\sum_{\beta=1}^n \frac{\partial b_j^\beta}{\partial z_j^\beta} \right)dz_i^1\wedge \cdots \wedge dz_i^n =\left(  \frac{1}{K_{ij}}\left[ b_i, K_{ij}  \right]  + \sum_{\beta=1}^n \frac{\partial b_j^\beta}{\partial z_j^\beta}  \right)dz_i^1\wedge \cdots \wedge dz_i^n
\end{align*}
Hence $\gamma$ is well-defined. Next we show that the following diagram commutes
\begin{center}
$\begin{CD}
\bigwedge^{n-1}T_X\otimes \mathcal{Q}_0 @>i_{(-)}(dz_i^1\wedge \cdots \wedge dz_i^n) >\cong> \Omega_X^1\otimes \mathcal{Q}_0\otimes K_X^* \\
@A\tilde{D}_0AA @AA\tilde{E}_0A \\
\mathcal{E}_{\mathcal{Q}_0} @>\gamma >\cong > \mathcal{E}_{\mathcal{Q}_0 \otimes K_X^*}
\end{CD}$
\end{center}
In fact, 
\begin{align*}
&i_{\tilde{D}_0(a_i,b_i)}\left( dz_i^1\wedge \cdots \wedge dz_i^n \right)= i_{\left[b_i, T_{0i}^1\wedge \cdots \wedge T_{0i}^{n-1} \right] + a_i T_{0i}^1\wedge \cdots \wedge T_{0i}^{n-1} }(dz_i^1\wedge \cdots \wedge dz_i^n)\\
&= i_{\left[ b_i, T_{0i}^1\wedge \cdots \wedge T_{0i}^{n-1}\right]}( dz_i^1\wedge \cdots \wedge dz_i^n) + a_i w_{0i} = \left( \mathcal{L}_{b_i}i_{T_{0i}^1\wedge \cdots \wedge T_{0i}^{n-1}} - i_{T_{0i}^1\wedge \cdots \wedge T_{0i}^{n-1}}\mathcal{L}_{b_i}  \right)(dz_i^1\wedge \cdots \wedge dz_i^n) + a_i w_{0i}\\
&=\mathcal{L}_{b_i}(w_{0i})+ a_i w_{0i}- i_{T_{0i}^1\wedge \cdots \wedge T_{0i}^{n-1}}\left(\sum_{\alpha=1}^n dz_i^1\wedge \cdots \wedge \mathcal{L}_{b_i}(dz_i^\alpha)\wedge \cdots \wedge dz_i^n \right) =\mathcal{L}_{b_i} (w_{0i}) + a_iw_{0i}- \left( \sum_{\alpha=1}^n \frac{\partial b_i^\alpha}{\partial z_i^\alpha} \right)w_{0i}\\
&=\tilde{E}_0(\gamma(a_i,b_i))
\end{align*}

We will show that other diagram also commutes.
\begin{center}
$\begin{CD}
\mathscr{H}om_{\mathcal{O}_X}\left(\bigwedge^2 T_{\mathcal{F}_0}, \frac{T_X}{T_{\mathcal{F}_0}} \right) @> - \textnormal{dual} >> \mathscr{H}om_{\mathcal{O}_X}\left( \mathcal{N}_{\mathcal{F}_0}^*, \bigwedge^2 T_{\mathcal{F}_0}^* \right) \\
@AD_2 AA @AAA \\
\mathscr{H}om_{\mathcal{O}_X}\left(T_{\mathcal{F}_0}, \frac{T_X}{T_{\mathcal{F}_0}} \right) @>  - \textnormal{dual} >> \mathscr{H}om_{\mathcal{O}_X}\left(\mathcal{N}_{\mathcal{F}_0}^*, T_{\mathcal{F}_0}^* \right) \\
@A\tilde{D}_1 AA @AAA \\
\mathcal{R} @. \mathscr{H}om_{\mathcal{O}_X}\left( \mathcal{N}_{\mathcal{F}_0}^*, \frac{\Omega_X^1}{\mathcal{Q}_0\otimes K_X} \right) \\
@A\textnormal{Id} AA @AAA \\
\mathcal{R}\subset \bigwedge^{n-1}T_X\otimes \mathcal{Q}_0 @>>> \Omega_X^1\otimes \mathcal{Q}_0\otimes K_X^*
\end{CD}$
\end{center}
We note that for $a_i=\sum_{\alpha=1}^{n-1} T_{0i}^1\wedge \cdots \wedge a_i^\alpha \wedge \cdots \wedge T_{0i}^\alpha \in \Gamma(U_i , \mathcal{R})$, we have
\begin{align*}
\mathcal{R} &\to \mathscr{H}om_{\mathcal{O}_X}\left(T_{\mathcal{F}_0}, \frac{T_X}{T_{\mathcal{F}_0}} \right) \xrightarrow{  - \textnormal{dual}}  \mathscr{H}om_{\mathcal{O}_X}\left(\mathcal{N}_{\mathcal{F}_0}^*, T_{\mathcal{F}_0}^* \right)\\
a_i &\mapsto \phi_i (T_{0i}^\alpha)= \overline{a_i^\alpha} \mapsto -\phi_i^*(w_{0i})(T_{0i}^\alpha)= - i_{\phi_i(T_{0i}^\alpha)} (w_{0i})=- i_{a_i^\alpha}(w_{0i})
\end{align*}
On the other hand, $a_i$ is mapped to $i_{a_i}(dz_i^1\wedge \cdots \wedge dz_i^n)\in \Gamma(U_i, \Omega_X^1\otimes \mathcal{Q}_0 \otimes K_X^*)$, and
{\begin{align*}
\Omega_X^1\otimes \mathcal{Q}_0 \otimes K_X^*&\to \mathscr{H}om_{\mathcal{O}_X}\left( \mathcal{N}_{\mathcal{F}_0}^*, \frac{\Omega_X^1}{\mathcal{Q}_0\otimes K_X} \right) \to \mathscr{H}om_{\mathcal{O}_X}\left( \mathcal{N}_{\mathcal{F}_0}^*, T_{\mathcal{F}_0}^* \right) \\
i_{a_i}(dz_i^1\wedge \cdots \wedge dz_i^n)&\mapsto \left( w_{0i} \mapsto  \sum_{\alpha=1}^{n-1} i_{T_{0i}^1\wedge \cdots \wedge a_i^\alpha \wedge \cdots \wedge T_{0i}^{n-1}} \left(dz_i^1\wedge \cdots \wedge dz_i^n\right)    \right) \\
& \mapsto \left( w_{0i} \mapsto \left( T_{0i}^\alpha \mapsto i_{T_{0i}^\alpha}\left( \sum_{\alpha=1}^{n-1} i_{T_{0i}^1\wedge \cdots \wedge a_i^\alpha \wedge \cdots \wedge T_{0i}^{n-1}} \left(dz_i^1\wedge \cdots \wedge dz_i^n\right)      \right)     \right) \right)
\end{align*}
so that in order to show the commutativity we have to show
\begin{align*}
- i_{a_i^\alpha}(w_{0i})=  i_{T_{0i}^\alpha}\left( \sum_{\alpha=1}^{n-1} i_{T_{0i}^1\wedge \cdots \wedge a_i^\alpha \wedge \cdots \wedge T_{0i}^{n-1}} \left(dz_i^1\wedge \cdots \wedge dz_i^n\right)      \right)     
\end{align*}
Indeed,
\begin{align*}
&i_{T_{0i}^\beta}\left( \sum_{\alpha=1}^{n-1} i_{T_{0i}^1\wedge \cdots \wedge a_i^\alpha \wedge \cdots \wedge T_{0i}^{n-1}}(dz_i^1\wedge \cdots \wedge dz_i^n) \right) = \sum_{\alpha=1}^{n-1} i_{T_{0i}^1\wedge \cdots \wedge a_i^\alpha \wedge \cdots \wedge T_{0i}^{n-1} \wedge T_{0i}^\beta} (dz_i^1\wedge \cdots \wedge dz_i^n)\\
&= - i_{T_{0i}^1\wedge \cdots \wedge T_{0i}^{n-1} \wedge a_i^\beta}(dz_i^1\wedge \cdots \wedge dz_i^n) =  - i_{a_i^\beta} i_{T_{0i}^1\wedge \cdots \wedge T_{0i}^{n-1}}(dz_i^1\wedge \cdots \wedge dz_i^n) = - i_{a_i^\beta } (w_{0i} )
\end{align*}
Hence the above diagram commutes. Then we have the long exact sequence
\begin{align*}
0\to \mathbb{H}^0\left(X, \mathcal{E}_{\mathcal{Q}_0}^{\bullet}\right)\to   \mathbb{H}^0\left( X, \mathcal{E}_{\mathcal{Q}_0\otimes K_X^*}^\bullet  \right) \to 0 \to  \mathbb{H}^1\left(X, \mathcal{E}_{\mathcal{Q}_0}^{\bullet}\right) \hookrightarrow  \mathbb{H}^1\left( X, \mathcal{E}_{\mathcal{Q}_0\otimes K_X^*}^\bullet  \right)   \\
\to \mathbb{H}^0(X, \mathcal{P}^\bullet) \to \mathbb{H}^2 \left(X, \mathcal{E}_{\mathcal{Q}_0}^{\bullet}\right) \xrightarrow{\rho}   \mathbb{H}^2\left( X, \mathcal{E}_{\mathcal{Q}_0\otimes K_X^*}^\bullet  \right) \to \cdots
\end{align*}

We will show that the above map
\begin{align}\label{nt13}
\mathbb{H}^1\left(X , \mathcal{E}_{\mathcal{Q}_0}^\bullet \right) \hookrightarrow \mathbb{H}^1\left(X, \mathcal{E}_{\mathcal{Q}_0\otimes K_X^*}^\bullet \right)
\end{align}
induces a morphism of first-order deformations between $\textnormal{FDef}_{\left(X, T_{\mathcal{F}_0} \right)}$ and $\textnormal{FDef}_{\left( X, \mathcal{N}_{\mathcal{F}_0}^*\right)}$ and 
\begin{align} \label{nt12}
\rho:\mathbb{H}^2\left(X, \mathcal{E}_{\mathcal{Q}_0}^\bullet \right) \to \mathbb{H}^2\left( X, \mathcal{E}_{\mathcal{Q}_0\otimes K_X^*}^\bullet \right)
\end{align}
defines the obstruction map between $\textnormal{FDef}_{\left(X, T_{\mathcal{F}_0} \right)}$ and $\textnormal{FDef}_{\left( X, \mathcal{N}_{\mathcal{F}_0}^*\right)}$.

Let $\left(\mathcal{X}, T_\mathcal{F} \right)$ be a first-oder foliated deformation of $\left( X, T_{\mathcal{F}_0}\right)$ over $k[\epsilon]$ in terms of tangent sheaves. We keep the notations in the proof of Proposition \ref{ts11}. We choose a Stein open covering $\mathcal{U}=\left\{ U_i \right\}$ of $X$ with local coordinates $\left( z_i^1,..., z_i^n \right)$ on $U_i$ such that $\Gamma\left(U_i,  T_{\mathcal{F}_0} \right)$ is generated by $T_{0i}^1,..., T_{0i}^{n-1}\in \Gamma\left( U_i, T_X  \right)$ with $T_{0i}^\alpha= \sum_{\beta=1}^{n-1} r_{0ij}^{\alpha\beta} T_{0j}^\beta$ and $\left[ T_{0i}^\alpha, T_{0i}^\beta \right]= \sum_{\gamma=1}^{n-1} g_{0i\alpha\beta}^\gamma T_{0i}^\gamma$ and $\Gamma\left( U_i, T_{\mathcal{F}} \right)$ is generated by $T_{0i}^\alpha + \epsilon T_i'^\alpha$ and $\theta_{ij}$ corresponds to $\textnormal{Id} + \epsilon p_{ij}$ with $p_{ij}=\sum_{\gamma=1}^n p_{ij}^\gamma \frac{\partial}{\partial z_i^\gamma}\in \Gamma\left(U_{ij}, T_X \right)$ and
\begin{align*}
\theta_{ji}^*\left(T_{0i}^\alpha + \epsilon T_i'^\alpha \right)= T_{0i}^\alpha+ \epsilon T_i'^\alpha- \epsilon \left[ p_{ij}, T_{0i}^\alpha \right] = \sum_{\beta=1}^{n-1} \left(r_{0ij}^{\alpha\beta} + \epsilon b_{ij}^{\alpha\beta} \right)\left( T_{0j}^\beta + \epsilon  T_j'^\beta \right)
\end{align*}
We set $\det\left(r_{0ij}^{\alpha\beta} + \epsilon b_{ij}^{\alpha\beta} \right)= r_{0ij}+ \epsilon B_{ij}$. Then $\left(\left\{ T_i'^*=\sum_{\alpha=1}^{n-1}T_{0i}^1\wedge \cdots \wedge T_{0i}'^\alpha \wedge \cdots \wedge T_{0i}^{n-1} \right\}, \left\{ \left(p_{ij}, \frac{B_{ij}}{r_{0ij}} \right)\right\} \right)\in C^0\left( \mathcal{U}, \mathcal{R} \right) \bigoplus C^1\left( \mathcal{U}, \mathcal{E}_{\mathcal{Q}_0} \right)$ defines a $1$-cocycle of \v Cech resolution $\mathcal{E}_{\mathcal{Q}_0}^\bullet$ in terms of the covering $\mathcal{U}$.

On the other hand, let us describe the induced first-order foliated deformation $\left(\mathcal{X}, \mathcal{N}_\mathcal{F}^* \right)$ in terms of cotangent sheaves from $\left(\mathcal{X}, T_{\mathcal{F}} \right)$ over $k[\epsilon]$. We keep the content of the proof of Proposition \ref{si39}. $\Gamma\left( U_i, \mathcal{N}_\mathcal{F}^* \right)$ is generated by
\begin{align*}
i_{\left(T_{0i}^1+ \epsilon T_i'^1 \right) \wedge \cdots \wedge \left( T_{0i}^{n-1} + \epsilon T_i'^{n-1} \right)}\left( dz_i^1\wedge \cdots \wedge dz_i^n \right)= i_{T_{0i}^1\wedge \cdots \wedge T_{0i}^{n-1}}\left(  dz_i^1\wedge \cdots \wedge dz_i^n \right) + \epsilon i_{T_i'^*}\left( dz_i^1\wedge \cdots \wedge dz_i^n \right) := w_{0i} + \epsilon w_i'
\end{align*}
and
\begin{align*}
\theta_{ji}^*\left( w_{0i}+ \epsilon w_i' \right) &= w_{0i} + \epsilon w_i' -\epsilon \mathcal{L}_{p_{ij}}\left( w_{0i} \right) = \left(r_{0ij} + \epsilon B_{ij} \right)\left( K_{0ij} + \epsilon L_{ij}\right) \left( w_{0j} + \epsilon w_j' \right) \\
&=\left( r_{0ij} K_{0ij} + \epsilon\left(B_{ij}K_{0ij} + r_{0ij}L_{ij} \right) \right)\left( w_{0j} + \epsilon w_j' \right)
\end{align*}
where $L_{ij}$ is defined by the following property
\begin{align*}
\left( K_{0ij} + \epsilon L_{ij} \right)dz_j^1\wedge \cdots \wedge dz_j^n &= d\left( \left(\textnormal{Id}- \epsilon p_{ij} \right)\left( z_i^1\right) \right) \wedge \cdots \wedge d\left(\left(\textnormal{Id}-\epsilon p_{ij} \right)\left( z_i^n \right) \right) = d\left(z_i^1- \epsilon p_{ij}^1 \right)\wedge \cdots \wedge d\left(z_i^n - \epsilon p_{ij}^n \right) \\
&= \left( 1 - \epsilon \sum_{\alpha=1}^n \frac{\partial p_{ij}^\alpha}{\partial z_i^\alpha} \right)dz_i^1 \wedge \cdots \wedge dz_i^n = \left( 1 - \epsilon \sum_{\alpha=1}^n \frac{\partial p_{ij}^\alpha}{\partial z_i^\alpha} \right) K_{0ij} dz_j^1\wedge \cdots \wedge dz_j^n \\
\Longrightarrow &  \frac{ L_{ij} }{K_{0ij}} =  -  \sum_{\alpha=1}^n \frac{\partial p_{ij}^\alpha}{\partial z_i^\alpha}
\end{align*}
Then $\left( \left\{w_i' \right\} , \left\{ \left( p_{ij}, \frac{B_{ij} K_{0ij} + r_{0ij} L_{ij}}{r_{0ij} K_{0ij}} = \frac{B_{ij}}{r_{0ij}} + \frac{L_{ij}}{K_{0ij}}= \frac{B_{ij}}{r_{0ij}} - \sum_{\alpha=1}^n \frac{\partial p_{ij}^\alpha}{\partial z_i^\alpha} \right) \right\} \right)\in C^0\left( \mathcal{U}, \mathcal{S}^1 \right) \bigoplus C^1\left(\mathcal{U}, \mathcal{E}_{\mathcal{Q}_0\otimes K_X^*} \right)$ defines a $1$-cocycle of \v Cech resolution of $\mathcal{E}_{\mathcal{L}_0}^\bullet=\mathcal{E}_{Q_0\otimes K_X^*}^\bullet$ in terms of the covering $\mathcal{U}$ by the proof of Proposition \ref{si39}. Then it is clear that under the map from $(\ref{nt13})$
\begin{align*}
\mathbb{H}^1\left( X, \mathcal{E}_{\mathcal{Q}_0}^\bullet \right) &\hookrightarrow \mathbb{H}^1\left( X, \mathcal{E}_{\mathcal{Q}_0\otimes K_X^*}^\bullet  \right) \\
\left( \left\{ T_i'^* \right\}, \left\{ \left( p_{ij}, \frac{B_{ij}}{r_{0ij}} \right) \right\} \right)  &\mapsto \left( \left\{ w_i'\right\}, \left\{ \left( p_{ij}, \frac{B_{ij}}{r_{0ij}} - \sum_{\alpha=1}^n \frac{ \partial p_{ij}^\alpha}{\partial z_i^\alpha}  \right) \right\} \right)
\end{align*}

Next we describe the obstruction map between $\textnormal{FDef}_{\left(X, T_{\mathcal{F}_0} \right)}$ and $\textnormal{FDef}_{\left( X, \mathcal{N}_{\mathcal{F}_0}^* \right)}$. Let $\left( \mathcal{X}, T_\mathcal{F} \right)$ be a foliated deformation of $\left(\mathcal{X}, T_{\mathcal{F}} \right)$ over $A\in \bold{Art}$ and $e:0\to (t)\to \tilde{A} \to A \to 0$ be a small extension. We keep the notations of the proof of Proposition \ref{ts11}. Let $\left( \mathcal{X}, T_{\mathcal{F}} \right)$ be a foliated deformation of $\left(X, T_{\mathcal{F}_0} \right)$ over $A$. Let $\mathcal{U}=\left\{ U_i \right\}$ be a Stein open covering of $X$ with local coordinates $(z_i^1,..., z_i^n)$ on $U_i$ such that $\theta_i:\left( U_i\times \textnormal{Spec}(A), T_{\mathcal{F}_i} \right)\to \left(\mathcal{X}|_{U_i}, T_\mathcal{F}|_{U_i} \right)$ is a foliated isomorphism, and $\theta_{ij}:=\theta_j^{-1} \theta_i$, and $\Gamma\left( U_i, T_{\mathcal{F}_i}\right)$ is generated by $T_i^1,..., T_i^{n-1}\in \Gamma\left( U_i, T_X \right)\otimes A$ with ${\theta}_{ji}^* T_i^\alpha= \sum_{\beta=1}^{n-1} r_{ij}^{\alpha\beta} T_j^\beta$ for $r_{ij}^{\alpha\beta}\in \Gamma\left( U_{ij}, \mathcal{O}_X \right)\otimes A$ and $\left[ T_i^\alpha, T_i^\beta \right]= \sum_{\gamma=1}^{n-1} g_{i\alpha\beta}^\gamma T_i^\gamma$. Then we choose arbitrary automorphism $\left\{ \tilde{\theta}_{ij}\right\}$ which induces $\theta_{ij}$, and liftings $\tilde{T}_i^\alpha$ of $T_i^\alpha$,  liftings $\tilde{r}_{ij}^{\alpha\beta}$ of $r_{ij}^{\alpha\beta}$ and liftings $\tilde{g}_{i\alpha\beta}^\gamma$ of $g_{i\alpha\beta}^\gamma$. Then $\left(\left\{ P_i \right\}, \left\{ - W_{ij}^*\right\}, \left\{\left( d_{ijk}, - g_{ijk}^* \right) \right\} \right)\in C^0\left(\mathcal{U}, \mathscr{H}om_{\mathcal{O}_X}\left( \bigwedge^2 T_{\mathcal{F}_0}, \frac{T_X}{T_{\mathcal{F}_0}} \right) \right)\bigoplus C^1\left( \mathcal{U}, \mathcal{R} \right) \bigoplus C^2\left( \mathcal{U}, \mathcal{E}_{\mathcal{Q}_0} \right)$ defines a $2$-cocycle of $\mathcal{E}_{\mathcal{Q}_0}^\bullet$ where $P_i:\Gamma\left(U_i, \bigwedge^2 T_{\mathcal{F}_0} \right) \to \Gamma\left( U_i, \frac{T_X}{T_{\mathcal{F}_0}}\right)$ such that $\left[ \tilde{T}_i^\alpha, \tilde{T}_i^\beta\right] - \sum_{\gamma=1}^{n-1} \tilde{g}_{i\alpha\beta}^\gamma \tilde{T}_i^\gamma= t P_i^{\alpha\beta}$, and $W_{ij}^*=\sum_{\alpha=1}^{n-1} T_{0i}^1\wedge \cdots \wedge W_{ij}^\alpha \wedge \cdots \wedge T_{0i}^{n-1}\in \Gamma\left( U_i, \mathcal{R} \right)$ such that $\tilde{\theta}_{ji}^* \left(\tilde{T}_i^\alpha \right)- \sum_{\beta=1}^{n-1} \tilde{r}_{ij} \tilde{T}_j^\beta =  t W_{ij}^\alpha$, and $d_{ijk}\in \Gamma\left(U_{ijk}, T_X \right)$ such that $\textnormal{Id} +t d_{ijk}= \tilde{f}_{ij}\tilde{f}_{jk}\tilde{f}_{ki}$, and $g_{ijk}^*\in \Gamma\left( U_{ijk}, \mathcal{O}_X  \right)$ such that $1+ t g_{ijk}^*= \tilde{f}_{ij} \tilde{f}_{jk}\left( \tilde{r}_{ik} \right)\tilde{f}_{ij}\left( \tilde{r}_{kj} \right)\tilde{r}_{ji}$ where $\tilde{r}_{ij}= \det\left( \tilde{r}_{ij}^{\alpha\beta} \right)$.

Now let us describe the induced foliated deformation $\left( \mathcal{X}, \mathcal{N}_\mathcal{F}^* \right)$ over $A$ from $\left(\mathcal{X}, T_\mathcal{F} \right)$. We keep the notations in the proof of Proposition \ref{si39}. Then $\Gamma\left( U_i, \mathcal{N}_{\mathcal{F}_i}^* \right)$ is generated by $i_{T_i^1\wedge \cdots \wedge T_i^{n-1}}\left( dz_i^1\wedge \cdots \wedge dz_i^n \right)= \sum_{\alpha=1}^n A_i^\alpha dz_i^\alpha= w_i$, and
\begin{align*}
\theta_{ji}^*\left(w_i \right)&= \theta_{ji}^* \left(  i_{ T_i^1 \wedge \cdots \wedge T_i^{n-1}} \left( dz_i^1 \wedge \cdots \wedge dz_i^n \right) \right)= i_{\theta_{ji}^*\left( T_i^1 \right) \wedge \cdots \wedge \theta_{ji}^*\left( T_i^{n-1} \right)} \theta_{ji}^*\left( dz_i^1\wedge \cdots \wedge dz_i^n \right) \\
&= r_{ij}K_{ij} i_{T_j^1\wedge \cdots \wedge T_j^{n-1}}\left( dz_j^1 \wedge \cdots \wedge dz_j^n \right)= r_{ij} K_{ij} w_j = h_{ij} w_j,\,\,\,\,\,\,\,\textnormal{where}\,\,\, h_{ij}:= r_{ij} K_{ij}
\end{align*}
and $K_{ij}\in \Gamma\left(U_{ij}, \mathcal{O}_{\mathcal{X}}^* \right)$ is defined by the following property
\begin{align*}
\theta_{ji}^*\left( dz_i^1\wedge \cdots \wedge dz_i^n \right)= d\left(f_{ji}\left( z_i^1 \right) \right) \wedge \cdots \wedge d\left( f_{ji}\left(  z_i^n \right) \right)= K_{ij} dz_j^1\wedge \cdots \wedge dz_j^n
\end{align*}
and satisfy the integrability condition $w_i\wedge dw_i=0$, equivalently $i_{T_i^\alpha \wedge T_i^\beta}\left( dw_i \right)=0$. Then the liftings $\tilde{T}_i^\alpha$ corresponds to the lifting $\tilde{w}_i=i_{\tilde{T}_i^1\wedge \cdots \wedge \tilde{T}_i^{n-1}}\left( dz_i^1 \wedge \cdots \wedge dz_i^n \right)\in \Gamma\left( U_i,  \Omega_X^1 \right)\otimes \tilde{A}$ of $w_i$, and we note that since $\tilde{f}_{ij}$ is an isomorphism, there exists a lifting $\tilde{K}_{ij}\in \Gamma\left( U_{ij}, \mathcal{O}_X \right)\otimes \tilde{A}$ of $K_{ij}$ such that $\tilde{K}_{ij} dz_j^1 \wedge \cdots \wedge dz_j^n = d\left(\tilde{f}_{ji}\left(z_i^1\right) \right) \wedge \cdots \wedge d \left( \tilde{f}_{ji} \left( z_i^n \right) \right)$. Then we have
\begin{align*}
\tilde{\theta}_{ji}^*\left( \tilde{w}_i \right)&= \tilde{\theta}_{ji}^*\left(i_{\tilde{T}_i^1 \wedge \cdots \wedge \tilde{T}_i^{n-1}}\left( dz_i^1 \wedge \cdots \wedge  dz_i^n \right) \right) = i_{\tilde{\theta}_{ji}^*\left(\tilde{T}_i^1 \wedge \cdots \wedge \tilde{T}_i^{n-1} \right)}\left( d \left(\tilde{f}_{ji}\left(z_i^1\right)\right) \wedge \cdots \wedge d \left( \tilde{f}_{ji}\left( z_i^n \right) \right) \right) \\
&= \tilde{r}_{ij} i_{\tilde{T}_j^1 \wedge \cdots \wedge \tilde{T}_j^{n-1}}\left( d \left(\tilde{f}_{ji}\left(z_i^1\right)\right) \wedge \cdots \wedge d \left( \tilde{f}_{ji}\left( z_i^n \right) \right)\right) + t i_{\sum_{\alpha=1}^{n-1} T_{0i}^1 \wedge \cdots \wedge W_{ij}^\alpha \wedge \cdots \wedge T_{0i}^{n-1}}\left( dz_i^1 \wedge \cdots \wedge dz_i^n \right)\\
&= \tilde{r}_{ij} \tilde{K}_{ij}  \tilde{w}_j  + t i_{W_{ij}^*}\left( dz_i^1 \wedge \cdots \wedge dz_i^n \right)
\end{align*}

Then by combining the proof of Proposition \ref{si39} and the proof of Proposition \ref{si30}, we see that $\left( \left\{ \phi_i \right\} , \left\{ - R_{ij}^* \right\}, \left\{ \left( d_{ijk}, - h_{ijk}^* \right) \right\} \right) \in C^0\left( \mathcal{U}, \mathscr{H}om_{\mathcal{O}_X}\left( \mathcal{N}_{\mathcal{F}_0}^*, \bigwedge^2 T_{\mathcal{F}_0}^* \right) \right) \bigoplus C^1\left(\mathcal{U}, \Omega_X^1\otimes \mathcal{Q}_0\otimes K_X^* \right) \bigoplus C^2\left( \mathcal{U}, \mathcal{E}_{\mathcal{Q}_0\otimes K_X^*} \right)$ defines a $2$-cocycle in the \v Cech resolution of $\mathcal{E}_{\mathcal{Q}_0\otimes K_X^*}^\bullet$ from $(\ref{nt11})$ in terms of the covering $\mathcal{U}$, where $\phi_i:\Gamma\left( U_i, \mathcal{N}_{\mathcal{F}_0}^* \right)\to \Gamma\left( U_i, \bigwedge^2 T_{\mathcal{F}_0}^* \right), w_{0i} \mapsto \left(T_{0i}^\alpha \wedge T_{0i}^\beta \mapsto f_{i\alpha\beta} \right)$ such that $tf_{i\alpha\beta}= i_{\tilde{T}_i^\alpha \wedge \tilde{T}_i^\beta}\left( d \tilde{w}_i \right)$, and $R_{ij}^*= i_{W_{ij}^*}\left( dz_i^1 \wedge \cdots \wedge dz_i^n \right)$, and $d_{ijk}\in \Gamma\left( U_{ijk}, T_X \right)$ such that $1+ t d_{ijk}= \tilde{f}_{ij}\tilde{f}_{jk}\tilde{f}_{ki}$, and $h_{ijk}^*\in \Gamma\left( U_{ijk} , \mathcal{O}_X \right)$ such that $1+ t h_{ijk}^*= \tilde{f}_{ij}\tilde{f}_{jk}\left( \tilde{r}_{ik} \tilde{K}_{ik} \right) \cdot \tilde{f}_{ij}\left( \tilde{r}_{kj} \tilde{K}_{kj} \right) \cdot \tilde{r}_{ji}\tilde{K}_{ji}$. Then we have to show that under the map $\rho$ from $(\ref{nt12})$
\begin{align*}
\rho: \mathbb{H}^2\left( X, \mathcal{E}_{\mathcal{Q}_0}^\bullet \right) &\to  \mathbb{H}^2\left(  X  , \mathcal{E}_{\mathcal{Q}_0\otimes K_X^*}^\bullet \right) \\
\left( \left\{ P_i \right\}, \left\{ - W_{ij}^* \right\} , \left\{ \left( d_{ijk}, - g_{ijk}^* \right) \right\} \right) &\mapsto \left(\left\{ \phi\right\}, \left\{- R_{ij}^* \right\},  \left\{ \left( d_{ijk}, - h_{ijk}^* \right) \right\} \right)
\end{align*}

In fact, we note that
\begin{align*}
ti_{P_i^{\alpha\beta}}\left( w_{0i} \right)& = i_{\left[ \tilde{T}_i^\alpha, \tilde{T}_i^\beta \right]- \sum_{\eta=1}^{n-1} \tilde{g}_{i\alpha\beta}^\eta \tilde{T}_i^\eta }\left( \tilde{w}_i \right)= i_{\left[ \tilde{T}_i^\alpha, \tilde{T}_i^\beta \right]- \sum_{\eta=1}^{n-1} \tilde{g}_{i\alpha\beta}^\eta \tilde{T}_i^\eta }i_{ \tilde{T}_i^1 \wedge \cdots  \wedge \tilde{T}_i^{n-1}}\left( dz_i^1 \wedge \cdots \wedge dz_i^n \right) \\
&= i_{\left[ \tilde{T}_i^\alpha, \tilde{T}_i^\beta \right] } i_{\tilde{T}_i^1 \wedge \cdots \wedge \tilde{T}_i^{n-1}}\left( dz_i^1 \wedge \cdots \wedge dz_i^n \right) = -  i_{\tilde{T}_i^\alpha \wedge \tilde{T}_i^\beta} \left( d\tilde{w}_i    \right) = -  t f_{i\alpha\beta}
\end{align*}
This implies that $P_i$ is mapped to $\phi_i$, and it is clear that $- W_{ij}^*$ is mapped to $- R_{ij}^*$, and so it remains to show that $\left(d_{ijk}, - g_{ijk}^* \right)$ is mapped to $\left( d_{ijk}, - h_{ijk}^* \right)$. First we note that  since $\tilde{f}_{ij}\tilde{f}_{jk}\tilde{f}_{ki}= \textnormal{Id}+ t d_{ijk}$, by letting $d_{ijk}=\sum_{\alpha=1}^n d_{ijk}^\alpha \frac{\partial}{\partial z_i^\alpha}$, we have
\begin{align*}
&\tilde{K}_{ik} dz_k^1 \wedge \cdots \wedge dz_k^n = d\left(\tilde{f}_{ki}\left(z_i^1\right) \right) \wedge \cdots \wedge d \left( \tilde{f}_{ki} \left(z_i^n \right)\right) \\
&\Longrightarrow \tilde{f}_{jk}\left(\tilde{K}_{ik} \right) \tilde{K}_{kj} dz_j^1 \wedge \cdots \wedge dz_j^n = d\left( \tilde{f}_{jk} \tilde{f}_{kj}\left( z_i^1 \right)   \right) \wedge \cdots \wedge d \left( \tilde{f}_{jk} \tilde{f}_{ki}\left( z_i^n \right) \right) \\
&\Longrightarrow  \tilde{f}_{ij}\tilde{f}_{jk}\left( \tilde{K}_{ik} \right) \cdot \tilde{f}_{ij}\left( \tilde{K}_{kj} \right) \cdot \tilde{K}_{ji} dz_i^1 \wedge \cdots \wedge dz_i^n= d\left( \tilde{f}_{ij} \tilde{f}_{jk} \tilde{f}_{ki} \left( z_i^1 \right) \right) \wedge \cdots  \wedge d \left( \tilde{f}_{ij} \tilde{f}_{jk} \tilde{f}_{ki} \left( z_i^n \right) \right) \\
& \,\,\,\,\,\,\,\,\,\,\,\,\,\,= d\left( z_i^1 + t d_{ijk}\left(z_i^1\right)\right) \wedge \cdots \wedge d \left( z_i^n + t  d_{ijk} \left(z_i^n \right) \right) = \left( 1+ t \sum_{\alpha=1}^n \frac{\partial d_{ijk}^\alpha}{\partial z_i^\alpha} \right) dz_i^1 \wedge \cdots \wedge dz_i^n
\end{align*}
Hence we have
\begin{align} \label{nt14}
 \tilde{f}_{ij}\tilde{f}_{jk}\left( \tilde{K}_{ik} \right) \cdot \tilde{f}_{ij}\left( \tilde{K}_{kj} \right) \cdot \tilde{K}_{ji} = 1+ t \sum_{\alpha=1}^n \frac{\partial d_{ijk}^\alpha}{\partial z_i^\alpha}
\end{align}
Then we have from $(\ref{nt14})$
\begin{align*}
1+ t h_{ijk}^* &=\left(\tilde{f}_{ij} \tilde{f}_{jk} \right)\left(\tilde{r}_{ik} \tilde{K}_{ik} \right)\cdot \tilde{f}_{ij}\left( \tilde{r}_{kj}\tilde{K}_{kj} \right)\cdot \tilde{r}_{ji } \tilde{K}_{ji} \\
&=\left( \tilde{f}_{ij} \tilde{f}_{jk} \right)\left( \tilde{r}_{ik} \right)\cdot \tilde{f}_{ij}\left( \tilde{r}_{kj} \right)\cdot \tilde{r}_{ji} \cdot  \left( \tilde{f}_{ij} \tilde{f}_{jk} \right)\left( \tilde{K}_{ik} \right)\cdot \tilde{f}_{ij}\left( \tilde{K}_{kj} \right)\cdot \tilde{K}_{ji} \\
&=\left( 1 + t g_{ijk}^* \right)\cdot  \left( \tilde{f}_{ij} \tilde{f}_{jk} \right)\left( \tilde{K}_{ik} \right)\cdot \tilde{f}_{ij}\left( \tilde{K}_{kj} \right)\cdot \tilde{K}_{ji}  = \left( \tilde{f}_{ij} \tilde{f}_{jk} \right)\left( \tilde{K}_{ik} \right)\cdot \tilde{f}_{ij}\left( \tilde{K}_{kj} \right)\cdot \tilde{K}_{ji}  + t g_{ijk}^* \\
&= 1+ t \left( g_{ijk}^* + \sum_{\alpha=1}^n \frac{\partial d_{ijk}^\alpha}{\partial z_i^\alpha}  \right) \Longrightarrow - h_{ijk}^*= - g_{ijk}^* -  \sum_{\alpha=1}^n \frac{\partial d_{ijk}^\alpha}{\partial z_i^\alpha}
\end{align*} 
This proves that $\left(d_{ijk}, - g_{ijk}^* \right)$ is mapped to $\left( d_{ijk}, - h_{ijk}^* \right)$.

\subsection{Isomorphism between $\textnormal{FDef}_{\left(X, T_{\mathcal{F}_0}, \mathcal{N}_{\mathcal{F}_0}^* \right)}$ and $\textnormal{FDef}_{\left( X, T_{\mathcal{F}_0} \right)}$ with $T_{\mathcal{F}_0}$ locally free of rank $n-1$ }\

Let $\left(X, T_{\mathcal{F}_0} \right)$ be a foliated nonsingular variety with $T_{\mathcal{F}_0}$ locally free of rank $n-1$. We have shown that a foliated deformation of $\left(X, T_{\mathcal{F}_0}\right)$ of dimension $n-1$ also induces a foliated deformation of $\left(X, \mathcal{N}_{\mathcal{F}_0}^*\right)$ of codimension $1$. This implies that we have a morphism of deformation functors $\textnormal{FDef}_{\left(X, T_{\mathcal{F}_0} \right)} \to \textnormal{FDef}_{\left(X, T_{\mathcal{F}_0}, \mathcal{N}_{\mathcal{F}_0}^* \right)}$ which is the inverse of the natural morphism of deformation functors $\textnormal{FDef}_{\left(X, T_{\mathcal{F}_0}, \mathcal{N}_{\mathcal{F}_0}^* \right)}\to \textnormal{FDef}_{\left( X, T_{\mathcal{F}_0} \right)}$ by forgetting $\mathcal{N}_{\mathcal{F}_0}^*$, so that we have  
\begin{align}\label{nt24}
\textnormal{FDef}_{\left(X, T_{\mathcal{F}_0}, \mathcal{N}_{\mathcal{F}_0}^* \right)} \cong  \textnormal{FDef}_{\left( X, T_{\mathcal{F}_0} \right)}
\end{align}
In particular we have
\begin{align*}
\mathbb{H}^1\left( X, \mathcal{F}_0^\bullet \right)\cong \mathbb{H}^1\left(  X, T_{\mathcal{F}_0}^\bullet \right)
\end{align*}
Moreover we show that 
\begin{align*}
\mathbb{H}^2\left( X, \mathcal{F}_0^\bullet \right) \cong \mathbb{H}^2\left( X, T_{\mathcal{F}_0}^\bullet \right)
\end{align*}
In fact, we recall that we have an exact sequence of complexes of sheaves
{\tiny{\begin{center}
$\begin{CD}
@. \cdots @. \cdots @. \cdots \\
@. @AAA @AAA @AAA \\
0 @>>>\left( \mathcal{N}_{\mathcal{F}_0}^* \right)^*\otimes \tilde{\mathcal{S}}^3 \bigoplus \left(\mathcal{N}_{\mathcal{F}_0}^* \right)^* \otimes \bigwedge^2 T_{\mathcal{F}_0}^* @>>>\bigwedge^3 T_{\mathcal{F}_0}^*\otimes \frac{T_X}{T_{\mathcal{F}_0}} \bigoplus \left( \mathcal{N}_{\mathcal{F}_0}^* \right)^*\otimes \tilde{\mathcal{S}}^3 \bigoplus \left( \mathcal{N}_{\mathcal{F}_0}^* \right)^* \otimes \bigwedge^2 T_{\mathcal{F}_0}^* @>>>  \bigwedge^3 T_{\mathcal{F}_0}^* \otimes \frac{T_X}{T_{\mathcal{F}_0}} @>>> 0 \\
@.@AAA @AAA @AAA \\
0 @>>> \left( \mathcal{N}_{\mathcal{F}_0}^* \right)^*\otimes \tilde{\mathcal{S}}^2 \bigoplus \left(\mathcal{N}_{\mathcal{F}_0}^*\right)^*\otimes T_{\mathcal{F}_0}^*  @>>> \bigwedge^2 T_{\mathcal{F}_0}^*\otimes \frac{T_X}{T_{\mathcal{F}_0}} \bigoplus \left( \mathcal{N}_{\mathcal{F}_0}^* \right)^*\otimes \tilde{\mathcal{S}}^2\bigoplus \mathcal{N}_{\mathcal{F}_0}^*\otimes T_{\mathcal{F}_0}^* @>>> \bigwedge^2 T_{\mathcal{F}_0}^*\otimes \frac{T_X}{T_{\mathcal{F}_0}} @>>> 0 \\
@. @AAA @AAA @AAA \\
0@>>> \mathscr{H}om_{\mathcal{O}_X}\left( \mathcal{N}_{\mathcal{F}_0}^*,  \frac{\Omega_X^1}{\mathcal{N}_{\mathcal{F}_0}^*} \right)  @>>> \mathscr{H}om_{\mathcal{O}_X}\left(T_{\mathcal{F}_0}, \frac{T_X}{T_{\mathcal{F}_0}} \right)\bigoplus \mathscr{H}om_{\mathcal{O}_X}\left( \mathcal{N}_{\mathcal{F}_0}^*, \frac{\Omega_X^1}{\mathcal{N}_{\mathcal{F}_0}^*} \right) @>>>  \mathscr{H}om_{\mathcal{O}_X}\left( T_{\mathcal{F}_0}, \frac{T_X}{T_{\mathcal{F}_0}} \right)@>>> 0 \\
@. @AAA @AAA @AAA \\
0 @>>> 0 @>>> T_X @>>> T_X @>>> 0
\end{CD}$
\end{center}}}
Let us denote by $\mathcal{P}^\bullet$ the right complex of sheaves. Then we have the long exact sequence
\begin{align*}
0 \to \mathbb{H}^0\left(X, \mathcal{F}_0^\bullet \right) \to \mathbb{H}^0\left( X, T_{\mathcal{F}_0}^\bullet \right) \to \mathbb{H}^0\left(X, \mathcal{P}^\bullet \right)=0 \to \mathbb{H}^1\left(X, \mathcal{F}_0^\bullet \right) \to \mathbb{H}^1\left(X, T_{\mathcal{F}_0}^\bullet \right) \xrightarrow{\alpha} \mathbb{H}^1\left( X, \mathcal{P}^\bullet \right) \\
\to \mathbb{H}^2 \left( X, \mathcal{F}_0^\bullet \right)  \to \mathbb{H}^2\left(X, T_{\mathcal{F}_0}^\bullet \right) \to \mathbb{H}^2\left( X, \mathcal{P}^\bullet \right) \to \cdots
\end{align*}

Since $\mathbb{H}^1\left(X, \mathcal{F}_0^\bullet \right)\cong \mathbb{H}^1\left( X, T_{\mathcal{F}_0}^\bullet \right)$, $\alpha$ is the zero map. This implies that $\mathbb{H}^2\left( X, \mathcal{F}_0^\bullet \right) \to \mathbb{H}^2\left( X, T_{\mathcal{F}_0}^\bullet \right)$ is injective. It remains to show that it is surjective. In fact, first we note that we have $\mathbb{H}^2\left(X, T_{\mathcal{F}_0}^\bullet \right)\cong \mathbb{H}^2\left(X, \mathcal{E}_{\mathcal{Q}_0}^\bullet \right) \xrightarrow{\rho} \mathbb{H}^2\left(X, \mathcal{E}_{\mathcal{Q}_0\otimes K_X^*}^\bullet \right)\cong \mathbb{H}^2\left(X, \mathcal{N}_{\mathcal{F}_0}^{*\bullet} \right)$ from $(\ref{nt12})$. Let $\left( \left\{ a_i \right\}, \left\{ b_{ij}^*:=\sum_{\alpha=1}^{n-1} T_{0i}^1\wedge \cdots \wedge b_{ij}^\alpha \wedge \cdots \wedge T_{0i}^{n-1}\right\} , \left\{\left( c_{ijk}, e_{ijk} \right) \right\}\right)\in C^0\left( \mathcal{U}, \mathscr{H}om_{\mathcal{O}_X}\left( \bigwedge^2 T_{\mathcal{F}_0}, \frac{T_X}{T_{\mathcal{F}_0}} \right) \right) \bigoplus C^1\left(\mathcal{U}, \mathcal{R} \right)\bigoplus C^2\left(\mathcal{U}, \mathcal{E}_{\mathcal{Q}_0} \right)$ be a $2$-cocycle of $\mathcal{E}_{\mathcal{Q}_0}^\bullet$ in terms of a Stein open covering $\mathcal{U}=\left\{U_i \right\}$ of $X$ with local coordinates $\left(z_i^1,..., z_i^n \right)$ on $U_i$, where $a_i\left(T_{0i}^\alpha \wedge T_{0i}^\beta \right)= \overline{a_i^{\alpha\beta}}$ with $a_i^{\alpha\beta}\in \Gamma\left(U_i, T_X \right)$, and $c_{ijk}\in \Gamma\left( U_{ijk}, T_X\right)$ and $d_{ijk}\in \Gamma\left( U_{ijk}, \mathcal{O}_X \right)$. Then $\left(\left\{ a_i \right\}, \left\{ b_{ij}^*\right\}, \left\{\left( c_{ijk}, e_{ijk} \right) \right\} \right)$ is mapped to $\left(\left\{ a_i\right\}, \left\{ \overline{b_{ij}}\right\}, \left\{ c_{ijk}\right\} \right)\in C^0\left( \mathcal{U} , \mathscr{H}om_{\mathcal{O}_X}\left( \bigwedge^2 T_{\mathcal{F}_0}, \frac{T_X}{T_{\mathcal{F}_0}} \right) \right) \bigoplus C^1\left( \mathcal{U}, \mathscr{H}om_{\mathcal{O}_X}\left( T_{\mathcal{F}_0}, \frac{T_X}{T_{\mathcal{F}_0}} \right) \right) \bigoplus C^2\left( \mathcal{U}, T_X \right)$ under the isomorphism $\mathbb{H}^2\left(X, \mathcal{E}_{\mathcal{Q}_0}^\bullet \right)\xrightarrow{\cong} \mathbb{H}^2\left( X, T_{\mathcal{F}_0}^\bullet \right)$, where $\overline{b}_{ij}\left(T_{0i}^\alpha \right)= \overline{b_{ij}^\alpha}$. On the other hand, $\left( \left\{ a_i \right\}, \left\{ b_{ij}^* \right\}, \left\{\left( c_{ijk}, e_{ijk} \right) \right\} \right)$ is mapped to $\left( \left\{ - a_i'\right\}, \left\{\overline{b_{ij}'}  \right\} , \left\{ c_{ijk}\right\} \right)\in C^0\left(\mathcal{U}, \mathscr{H}om_{\mathcal{O}_X}\left(\mathcal{N}_{\mathcal{F}_0}^* , \bigwedge^2 T_{\mathcal{F}_0}^* \right) \right) \bigoplus C^1\left( \mathcal{U}, \mathscr{H}om_{\mathcal{O}_X}\left( \mathcal{N}_{\mathcal{F}_0}^*, \frac{\Omega_X^1}{\mathcal{N}_{\mathcal{F}_0}^*} \right) \right) \bigoplus C^2\left(\mathcal{U}, T_X \right)$ under $\mathbb{H}^2\left(X, \mathcal{E}_{\mathcal{Q}_0}^\bullet \right) \xrightarrow{\rho} \mathbb{H}^2\left( X, \mathcal{E}_{\mathcal{Q}_0\otimes K_X^*}^\bullet \right)\cong \mathbb{H}^2\left( X, \mathcal{N}_{\mathcal{F}_0}^{*\bullet} \right)$, where $a_i'\left(T_{0i}^\alpha \wedge T_{0i}^\beta \right)\left( w_{0i} \right)= i_{a_i^{\alpha\beta}}\left( w_{0i} \right)$, and $\overline{b_{ij}'}\left(w_{0i}\right)= \overline{i_{b_{ij}^*}\left( dz_i^1 \wedge \cdots \wedge dz_i^n \right) } $. Then we claim that
{\Small{\begin{align*}
&\left(\left\{ a_i \right\}, \left\{ -a_i' \right\} , 0,  \left\{ \overline{b_{ij}} \right\}, \left\{\overline{b_{ij}'} \right\}, \left\{c_{ijk} \right\} \right)\\
&\in C^0\left( \mathcal{U}, \mathscr{H}om_{\mathcal{O}_X}\left(\bigwedge^2  T_{\mathcal{F}_0}, \frac{T_X}{T_{\mathcal{F}_0}} \right) \right) \bigoplus C^0\left( \mathcal{U}, \mathscr{H}om_{\mathcal{O}_X}\left(\mathcal{N}_{\mathcal{F}_0}^* , \bigwedge^2 T_{\mathcal{F}_0}^* \right) \right) \bigoplus C^0\left( \mathcal{U}, \mathscr{H}om_{\mathcal{O}_X}\left( \mathcal{N}_{\mathcal{F}_0}^* , T_{\mathcal{F}_0}^* \right) \right)\\
&\,\,\,\,\, \bigoplus C^1\left( \mathcal{U}, \mathscr{H}om_{\mathcal{O}_X}\left( T_{\mathcal{F}_0}, \frac{T_X}{T_{\mathcal{F}_0}} \right) \right) \bigoplus C^1\left( \mathcal{U}, \mathscr{H}om_{\mathcal{O}_X}\left( \mathcal{N}_{\mathcal{F}_0}^*, \frac{\Omega_X^1}{\mathcal{N}_{\mathcal{F}_0}^*} \right) \right) \bigoplus C^2\left( \mathcal{U}, T_X  \right)
\end{align*}}}
defines a $2$-cocycle of $\mathcal{F}_0^\bullet$. Then this implies $\mathbb{H}^2\left(X, \mathcal{F}_0^\bullet \right)\to \mathbb{H}^2\left( X, T_{\mathcal{F}_0}^\bullet \right)$ is surjective. In fact, it is sufficient to show that
\begin{align} 
i_{a_i\left(T_{0i}^\alpha \wedge T_{0i}^\beta \right)}\left(w_{0i} \right) & - a_i'\left(T_{0i}^\alpha\wedge T_{0i}^\beta \right)\left( w_{0i} \right)=0 \label{nt20} \\
i_{b_{ij}\left(T_{0i}^\alpha \right)}\left( w_{0i} \right) & + i_{T_{0i}^\alpha}\left( b_{ij}'\left( w_{0i} \right) \right) = 0  \label{nt21}
\end{align}
In fact, $(\ref{nt20})$ follows from the definition of $a'$, and $(\ref{nt21})$ follows from
\begin{align*}
&i_{b_{ij}\left(T_{0i}^\alpha \right)}\left( w_{0i} \right) + i_{T_{0i}^\alpha}\left( b_{ij}'\left( w_{0i} \right) \right) =i_{b_{ij}^\alpha}\left( w_{0i} \right) + i_{T_{0i}^\alpha}\left( i_{b_{ij}^*}\left(dz_i^1\wedge \cdots \wedge dz_i^n \right)    \right)=i_{b_{ij}^\alpha}\left(w_{0i} \right) + i_{b_{ij}^*\wedge T_{0i}^\alpha}\left( dz_i^1\wedge \cdots \wedge dz_i^n \right)\\
&= i_{b_{ij}^\alpha}\left(w_{0i} \right) + i_{ T_{0i}^1\wedge \cdots \wedge b_{ij}^\alpha \wedge \cdots \wedge T_{0i}^{n-1}\wedge T_{0i}^\alpha}\left( dz_i^1 \wedge \cdots \wedge dz_i^n \right) = i_{b_{ij}^\alpha}\left(w_{0i} \right) - i_{ b_{ij}^\alpha}i_{T_{0i}^1 \wedge \cdots \wedge T_{0i}^{n-1}}\left( dz_i^1 \wedge \cdots \wedge dz_i^n \right) = 0 
\end{align*}

\section{ Relationship between foliated deformations of $\left(X, \mathcal{N}_{\mathcal{F}_0}^* \right)$ of codimension $n-1$ with $\mathcal{N}_{\mathcal{F}_0}^*$ locally free and foliated deformations of $\left( X, T_{\mathcal{F}_0} \right)$ of dimension $1$}

Let $\left(X, \mathcal{N}_{\mathcal{F}_0}^* \right)$ be a foliated nonsingular complex variety with $\dim_\mathbb{C} X=n$ and $\mathcal{N}_{\mathcal{F}_0}^*$ be locally free with rank $n-1$. Let $\mathcal{U}=\left\{ U_i \right\}$ be a Stein open covering of $X$ with local coordinates $\left( z_i^1,..., z_i^n \right)$ on $U_i$ such that $\Gamma\left( U_i, \mathcal{N}_{\mathcal{F}_0}^* \right)$ is generated by $w_{0i}^1,..., w_{0i}^{n-1}\in \Gamma\left( U_i, \Omega_X^1 \right)$ with $w_{0i}^\alpha= \sum_{\beta=1}^q h_{0ij}^{\alpha\beta} w_{0j}^\beta,\alpha=1,...,n-1$. We note that the integrability condition trivially holds. We have $w_{0i}^1 \wedge \cdots \wedge w_{0i}^{n-1}= h_{0ij} w_{0j}^1\wedge \cdots \wedge w_{0j}^{n-1}$, where $h_{0ij}= \det\left( h_{0ij}^{\alpha\beta} \right)$. Then $\mathcal{L}_0=\left( \mathcal{N}_{\mathcal{F}_0}^* \right)^*$ is defined by the transition functions $\left\{ h_{0ij} \right\}$ with $h_{0ij}\in \Gamma\left( U_{ij}, \mathcal{O}_X^* \right)$ and $w_0$ is defined locally on $U_i$ by  $w_0:=\left\{ w_{0i}^1\wedge \cdots \wedge w_{0i}^{n-1} \right\}\in H^0\left( X, \bigwedge^{n-1} \Omega_X^1 \otimes \mathcal{L}_0 \right)$ as in section \ref{nt5}.

We will describe the associated dimension $1$-foliation $T_{\mathcal{F}_0}$. Since $\bigwedge^{n-1} \Omega_X^1\times \Omega_X^1 \to \bigwedge^n \Omega_X^1= K_X$ is a perfect pairing, we have
\begin{align*}
\bigwedge^{n-1} \Omega_X^1 \cong \mathscr{H}om_{\mathcal{O}_X}\left(\Omega_X^1, K_X \right) \cong T_X\otimes K_X
\end{align*}
so that we have
\begin{align*}
\bigwedge^{n-1} \Omega_X^1\otimes \mathcal{L}_0 &\cong T_X\otimes K_X\otimes \mathcal{L}_0 \\
 w_{0i}^1\wedge \cdots \wedge w_{0i}^{n-1} &\mapsto \left( dz_i^\alpha \mapsto w_{0i}^1\wedge \cdots \wedge w_{0i}^{n-1}\wedge dz_i^\alpha         = a_i^\alpha dz_i^1\wedge \cdots \wedge dz_i^n \right) \mapsto  \sum_{\alpha=1}^n a_i^\alpha \frac{\partial}{\partial z_i^\alpha} := T_{0i}
\end{align*}

We note that
\begin{align*}
a_i^\alpha= i_{\frac{\partial}{\partial z_i^1}\wedge \cdots \wedge \frac{\partial}{\partial z_i^n}}\left(w_{0i}^1\wedge \cdots \wedge w_{0i}^{n-1} \wedge dz_i^\alpha \right) \Longrightarrow T_{0i}:=\sum_{\alpha=1}^n  i_{\frac{\partial}{\partial z_i^1}\wedge \cdots \wedge \frac{\partial}{\partial z_i^n}}\left(w_{0i}^1\wedge \cdots \wedge w_{0i}^{n-1} \wedge dz_i^\alpha \right) \frac{\partial}{\partial z_i^\alpha}
\end{align*}
and we also note that
\begin{align*}
i_{T_{0i}}\left( w_{0i}^\beta \right) &= \sum_{\alpha=1}^n i_{\frac{\partial}{\partial z_i^1}\wedge \cdots \wedge \frac{\partial}{\partial z_i^n}} \left( w_{0i}^1\wedge \cdots \wedge w_{0i}^{n-1} \wedge dz_i^\alpha\right)\cdot i_{\frac{\partial}{\partial z_i^\alpha}}\left( w_{0i}^\beta \right)\\
&= i_{\frac{\partial}{\partial z_i^1}\wedge \cdots \wedge \frac{\partial}{\partial z_i^n}} \left( w_{0i}^1\wedge \cdots \wedge w_{0i}^{n-1} \wedge \left(\sum_{\alpha=1}^n i_{\frac{\partial}{\partial z_i^\alpha}}\left( w_{0i}^\beta \right) dz_i^\alpha \right) \right) = i_{\frac{\partial}{\partial z_i^1}\wedge \cdots \wedge \frac{\partial}{\partial z_i^n}} \left( w_{0i}^1\wedge \cdots \wedge w_{0i}^{n-1}\wedge w_{0i}^\beta \right) = 0
\end{align*}
On the other hand, we have
\begin{align*}
\frac{\partial }{\partial z_i^1}\wedge \cdots \wedge \frac{\partial}{\partial z_i^n}= R_{0ij} \frac{\partial}{\partial z_j^1} \wedge \cdots \wedge \frac{\partial}{\partial z_j^n}\,\,\,\,\,\,\,\textnormal{for some}\,\,\, K_{0ij}\in \Gamma\left( U_{ij} , \mathcal{O}_X^*  \right)
\end{align*}
Then $K_X$ is defined by the transition functions $\left\{ R_{0ij}\right\} \in C^0\left( \mathcal{U}, K_X \right)$. Then
\begin{align*}
T_{0i}&= \sum_{\alpha=1}^n i_{\frac{\partial}{\partial z_i^1}\wedge \cdots \wedge \frac{\partial}{\partial z_i^n}} \left( w_{0i}^1\wedge \cdots \wedge w_{0i}^{n-1} \wedge dz_i^\alpha\right)\frac{\partial}{\partial z_i^\alpha}\\
&= R_{0ij}\cdot h_{0ij}\sum_{\alpha, \beta=1}^n i_{\frac{\partial}{\partial z_j^1}\wedge \cdots \wedge \frac{\partial}{\partial z_j^n}}\left( w_{0j}^1\wedge \cdots \wedge w_{0j}^{n-1}\wedge dz_j^\beta \right)\cdot \frac{\partial z_i^\alpha}{\partial z_j^\beta} \frac{\partial}{\partial z_i^\alpha} \\
& = R_{0ij} h_{0ij} \sum_{\beta=1}^n i_{\frac{\partial}{\partial z_j^1}\wedge \cdots \wedge \frac{\partial}{\partial z_j^n}}\left(w_{0j}^1\wedge \cdots \wedge w_{0j}^{n-1}\wedge dz_j^\beta \right) \frac{\partial}{\partial z_j^\beta} = R_{0ij} h_{0ij} T_{0j}
\end{align*}
This implies that
\begin{align*}
T_0=\left\{ T_{0i} \right\}\in H^0\left( X, T_X\otimes K_X\otimes \mathcal{L}_0 \right)=H^0\left(X, T_X \otimes \mathcal{Q}_0 \right),\,\,\,\,\,\,\, \mathcal{Q}_0 : = K_X \otimes \mathcal{L}_0=  T_{\mathcal{F}_0}^*
\end{align*}

\subsection{Foliated deformations of $\left(X, \mathcal{N}_{\mathcal{F}_0}^* \right)$ with $\mathcal{N}_{\mathcal{F}_0}^*$ locally free of rank $n-1$ in terms of cotangent sheaves and foliated deformations of $\left(X, T_{\mathcal{F}_0} \right)$}\

We keep the notations in the previous subsection. We will show that a foliated deformation $\left( \mathcal{X}, \mathcal{N}_\mathcal{F}^* \right)$ of $\left(X, \mathcal{N}_{\mathcal{F}_0}^* \right)$ with $\mathcal{N}_{\mathcal{F}_0}^*$ locally free of rank $n-1$ over $A\in \bold{Art}$ in terms of cotangent sheaves also defines a foliated deformation $\left( \mathcal{X}, T_\mathcal{F} \right)$ of $\left( X, T_{\mathcal{F}_0}\right)$ of dimension $1$ in terms of tangent sheaves, equivalently in terms of deformations of the section $T_0$. Let $\mathcal{U}=\left\{ U_i \right\}$ be a Stein open covering of $X$ with local coordinates $(z_i^1,,.., z_i^n)$ on $U_i$ such that for all $i$, we have an isomorphism
\begin{align*}
\theta_i:\left( U_i\times \textnormal{Spec}(A), \mathcal{N}_{\mathcal{F}_i}^* \right)\to \left( \mathcal{X}|_{U_i}, \mathcal{N}_{\mathcal{F}}^* |_{U_i} \right)
\end{align*}
and $\Gamma\left( U_i, \mathcal{N}_{\mathcal{F}_i}^* \right)$ is generated by $w_i^\alpha \in \Gamma\left(U_i, \Omega_X^1 \right)\otimes A$ which induces $w_{0i}^\alpha,\alpha=1,...,n-1$, and $\theta_{ij}=\theta_{j}^{-1}\theta_i$ is a foliated isomorphism with $\theta_{jk}\theta_{ij}=\theta_{ik}$ and
\begin{align*}
\theta_{ji}^*\left( w_i^\alpha \right) = \sum_{\beta=1}^{n-1} h_{ij}^{\alpha\beta} w_j^\beta,\,\,\,\,\,\,\,\,\,h_{ij}^{\alpha\beta} \in \Gamma\left( U_{ij}, \mathcal{O}_X \right)\otimes A,
\end{align*}
Let $\mathcal{L}=\left( \bigwedge^{n-1} \mathcal{N}_\mathcal{F}^* \right)^*$. Then $\{w_i^1\wedge \cdots \wedge w_i^{n-1}\} \in H^0\left(\mathcal{X}, \Omega_{\mathcal{X}/A}^1\otimes \mathcal{L} \right)$. We note that
\begin{align*}
\theta_{ji}^*\left(w_i^1\wedge \cdots \wedge w_i^{n-1}\right) = \theta_{ji}(w_i^1)\wedge \cdots \wedge \theta_{ji}^*(w_i^{n-1}) = h_{ij} w_j^1\wedge \cdots \wedge w_j^{n-1},\,\,\,\,\,\,\,h_{ij}=\det\left( h_{ij}^{\alpha\beta} \right)
\end{align*}

We note that since $\bigwedge^{n-1}\Omega_{\mathcal{X}/A}\times \Omega_{\mathcal{X}/A}^1\to \bigwedge^n \Omega_{\mathcal{X}/A}= K_{\mathcal{X}/A}$ is a perfect pairing, we have
\begin{align*}
\bigwedge^{n-1}\Omega_{\mathcal{X}/A}^1\cong \mathscr{H}om_{\mathcal{O}_X}\left( \Omega_{\mathcal{X}/A}^1,  K_{\mathcal{X}/A} \right)\cong T_{\mathcal{X}/A}\otimes K_{\mathcal{X}/A}
\end{align*}

Locally on $U_i$, we have an isomorphism
\begin{align*}
\Gamma\left( U_i, \bigwedge^{n-1} \Omega_X^1  \right)\otimes A  &\cong   \Gamma(U_i, T_X)\otimes A \\
w_i^1\wedge \cdots \wedge w_i^{n-1} &\mapsto \sum_{\alpha=1}^n A_i\frac{\partial}{\partial z_i^\alpha}=T_i,\,\,\,\,\,\,\,\, A_i^\alpha \in \Gamma(U_i, T_X)\otimes A
\end{align*}
where $A_i^\alpha$ is defined by $w_i^1\wedge \cdots \wedge w_i^{n-1} \wedge dz_i^\alpha= A_i^\alpha dz_i^1\wedge \cdots \wedge dz_i^n$, equivalently 
\begin{align*}
A_i^\alpha = i_{\frac{\partial}{\partial z_i^1}\wedge \cdots \wedge \frac{\partial}{\partial z_i^n}}\left( w_i^1\wedge \cdots \wedge w_i^{n-1} \wedge dz_i^\alpha \right) \Longrightarrow T_i=\sum_{\alpha=1}^n i_{\frac{\partial}{\partial z_i^1}\wedge \cdots \wedge \frac{\partial}{\partial z_i^n}}\left( w_i^1\wedge \cdots \wedge w_i^{n-1} \wedge dz_i^\alpha \right) \frac{\partial}{\partial z_i^\alpha}
\end{align*}
 Then we have
\begin{align*}
\left\{T_i=\sum_{\alpha=1}^n A_i^\alpha \frac{\partial}{\partial z_i^\alpha}\right\}\in H^0\left(\mathcal{X}, T_{\mathcal{X}/A}\otimes K_{\mathcal{X}/A}\otimes \mathcal{L}\right)
\end{align*}
which induces $\left\{\sum_{\alpha=1}^n a_i^\alpha \frac{\partial}{\partial z_i^\alpha} \right\}$. This is well-defined. In fact, $K_{\mathcal{X}/A}$ is described by $R_{ij}\in \Gamma\left( U_{ij}, \mathcal{O}_\mathcal{X}^* \right)$ in terms of $\frac{\partial}{\partial z_i^1}\wedge \cdots \wedge \frac{\partial}{\partial z_i^n}$ as
\begin{align*}
R_{ij}\frac{\partial}{\partial z_j^1}\wedge \cdots \wedge \frac{\partial}{\partial z_j^n}=  \theta_{ji}^*\left( \frac{\partial}{\partial z_1^1}  \right) \wedge \cdots \wedge \theta_{ji}^*\left( \frac{\partial}{\partial z_i^n} \right)
\end{align*}
Then
\begin{align*}
\theta_{ji}^*\left(\sum_{\alpha=1}^n A_i^\alpha \frac{\partial}{\partial z_i^\alpha} \right) &= \theta_{ji}^*\left(  \sum_{\alpha=1}^ni_{\frac{\partial}{\partial z_i^1}\wedge \cdots \wedge \frac{\partial}{\partial z_i^n}}\left( w_i^1\wedge \cdots \wedge w_i^{n-1}\wedge dz_i^\alpha \right)\frac{\partial}{\partial z_i^\alpha}  \right) \\
&= R_{ij} h_{ij} \sum_{\alpha=1}^n i_{ \frac{\partial}{\partial z_j^1 } \wedge \cdots \wedge \frac{\partial}{\partial z_j^n}}\left(w_j^1 \wedge \cdots \wedge w_j^{n-1} \wedge d\theta_{ji}^*(z_i^\alpha) \right)\cdot \theta_{ji}^*\left( \frac{\partial}{\partial z_i^\alpha} \right)\\
&= R_{ij} h_{ij} \sum_{\alpha=1}^n i_{ \frac{\partial}{\partial z_j^1 } \wedge \cdots \wedge \frac{\partial}{\partial z_j^n}}\left(w_j^1 \wedge \cdots \wedge w_j^{n-1} \wedge d\theta_{ji}^*(z_i^\alpha) \right)\cdot \left( \sum_{\beta=1}^n i_{\theta_{ji}^*\left( \frac{\partial}{\partial z_i^\alpha} \right)}\left(dz_j^\beta \right)\frac{\partial}{\partial z_j^\beta} \right)  \\
&= R_{ij} h_{ij} \sum_{\beta=1}^n i_{ \frac{\partial}{\partial z_j^1 } \wedge \cdots \wedge \frac{\partial}{\partial z_j^n}}\left(w_j^1 \wedge \cdots \wedge w_j^{n-1} \wedge \left( \sum_{\alpha=1}^n  i_{\theta_{ji}^*\left(\frac{\partial}{\partial z_i^\alpha} \right)}\left( dz_j^\beta \right) d\theta_{ji}^*(z_i^\alpha) \right) \right) \frac{\partial}{\partial z_j^\beta}\\
&= R_{ij} h_{ij} \sum_{\beta=1}^n i_{ \frac{\partial}{\partial z_j^1 } \wedge \cdots \wedge \frac{\partial}{\partial z_j^n}}\left(w_j^1 \wedge \cdots \wedge w_j^{n-1} \wedge \theta_{ji}^*\left( \sum_{\alpha=1}^n  i_{\frac{\partial}{\partial z_i^\alpha}}\left( d\theta_{ij}^*\left(z_j^\beta\right) \right) d z_i^\alpha \right) \right) \frac{\partial}{\partial z_j^\beta}\\
&= R_{ij} h_{ij} \sum_{\beta=1}^n i_{ \frac{\partial}{\partial z_j^1 } \wedge \cdots \wedge \frac{\partial}{\partial z_j^n}}\left(w_j^1 \wedge \cdots \wedge w_j^{n-1} \wedge \theta_{ji}^*\left(  d\theta_{ij}^*\left(z_j^\beta\right) \right)\right) \frac{\partial}{\partial z_j^\beta}\\
&= R_{ij} h_{ij} \sum_{\beta=1}^n i_{ \frac{\partial}{\partial z_j^1 } \wedge \cdots \wedge \frac{\partial}{\partial z_j^n}}\left(w_j^1 \wedge \cdots \wedge w_j^{n-1} \wedge d z_j^\beta \right) \frac{\partial}{\partial z_j^\beta} = R_{ij} h_{ij} \sum_{\beta=1}^n A_j^\beta \frac{\partial}{\partial z_j^\beta}
\end{align*}

On the other hand, let $T_\mathcal{F}$ be the locally free sheaf of $T_{\mathcal{X}/A}$ defined by $\left\{ T_i=\sum_{\alpha=1}^n A_i^\alpha \frac{\partial}{\partial z_i^\alpha} \right\}$, i.e. defined by the transition function $\left\{ \frac{1}{h_{ij} R_{ij}} \right\}$. Then we have the exact sequence $0\to T_\mathcal{F} \to T_{\mathcal{X}/A}\to \frac{T_{\mathcal{X}/A}}{T_{\mathcal{F}}}\to 0$. Since $T_\mathcal{F}$ and $T_{\mathcal{X}/A}$ are locally free, they are flat over $A$, and the exact sequence induces an injection $0\to T_{\mathcal{F}_0} \to T_X$ since $T_{\mathcal{F}_0}$ is the locally free sheaf defined by $\left\{T_{0i}= \sum_{\alpha=1}^n a_i^\alpha \frac{\partial}{\partial z_i^\alpha} \right\}$. This implies that $\frac{T_{\mathcal{X}/A}}{T_\mathcal{F}}$ is flat over $A$. It remains to show that $T_\mathcal{F}=\left(\frac{\Omega_{\mathcal{X}/A}^1}{\mathcal{N}_\mathcal{F}^*} \right)^*$. First we claim that
\begin{align} \label{nt6}
\mathcal{N}_\mathcal{F}^* = \left(\frac{T_{\mathcal{X}/A}}{T_{\mathcal{F}}} \right)^*.
\end{align}
 We have $\mathcal{N}_\mathcal{F}^* \subset  \left(\frac{T_{\mathcal{X}/A}}{T_{\mathcal{F}}} \right)^*$. On the other hand, We have the exact sequence $0\to \mathcal{N}_{\mathcal{F}}^* \to \Omega_{\mathcal{X}/A}^1 \to \frac{\Omega_{\mathcal{X}/A}^1}{\mathcal{N}_\mathcal{F}^*}\to 0$. Since $\frac{\Omega_{\mathcal{X}/A}^1}{\mathcal{N}_\mathcal{F}^*}$ is flat over $A$, and $\frac{\Omega_X^1}{\mathcal{N}_{\mathcal{F}_0}^*}$ is locally free on $U_i-S$, we see that $\frac{\Omega_{\mathcal{X}/A}^1}{\mathcal{N}_\mathcal{F}^*}$ is locally free on $U_i-S$.  From the exact sequence of locally free sheaves
\begin{align}\label{nt7}
0\to \mathcal{N}_{\mathcal{F}}^*|_{U_i-S} \to \Omega_{\mathcal{X}/A}^1|_{U_i-S} \to \frac{\Omega_{\mathcal{X}/A}^1}{\mathcal{N}_\mathcal{F}^*}|_{U_i-S} \to 0
\end{align}
we see that for $x\in U_i-S$,  there is a neighborhood $U_x\subset U_i-S$ such that $\Gamma\left(U_x, \frac{\Omega_{\mathcal{X}/A}^1}{\mathcal{N}_\mathcal{F}^*} \right)$ is generated by $\bar{W}_x$ such that $W_x\in \Gamma\left(U_x, \Omega_X^1 \right)\otimes A$ where $\bar{c}$ is the image of $c$ under $\Omega_{\mathcal{X}/A}^1 \to \frac{\Omega_{\mathcal{X}/A}^1}{\mathcal{N}_\mathcal{F}^*}$. Then $\Omega_{\mathcal{X}/A}^1|_{U_x}$ is generated by $w_i^1,..., w_i^n, W_x$. Let $Q\in \Gamma\left( U_i  , \Omega_{\mathcal{X}/A}^1 \right)= \Gamma(U_i, \Omega_X^1)\otimes A$ with $i_{\sum_{\alpha=1}^n A_i^\alpha \frac{\partial}{\partial z_i^\alpha}}(Q)=0$, i.e. $Q\in \left(\frac{\Omega_{\mathcal{X}/A}^1}{\mathcal{N}_\mathcal{F}^*} \right)^*$. We can write $Q= \sum_{\alpha=1}^{n-1} b_x^\alpha w_i^\alpha + a_x W_x$ on $U_x$. Then
\begin{align*}
0&=i_{\sum_{\alpha=1}^n A_i^\alpha \frac{\partial}{\partial z_i^\alpha}}\left(Q \right)= \sum_{\alpha=1}^n A_i^\alpha i_{\frac{\partial}{\partial z_i^\alpha}}\left( Q \right) = \sum_{\alpha=1}^n i_{\frac{\partial}{\partial z_i^1}\wedge \cdots \wedge \frac{\partial}{\partial z_i^n}}\left(w_i^1 \wedge \cdots \wedge w_i^{n-1}\wedge dz_i^\alpha \right) i_{\frac{\partial}{\partial z_i^\alpha}}(Q)\\
&= i_{\frac{\partial}{\partial z_i^1}\wedge \cdots \wedge \frac{\partial}{\partial z_i^n}}\left( w_i^1 \wedge \cdots \wedge w_i^{n-1} \wedge \left( \sum_{\alpha=1}^n i_{\frac{\partial}{\partial z_i^\alpha}}(Q) dz_i^\alpha \right)    \right) = i_{\frac{\partial}{\partial z_i^1}\wedge \cdots \wedge \frac{\partial}{\partial z_i^n}} \left( w_i^1 \wedge \cdots \wedge w_i^{n-1} \wedge Q \right)\\
&\Longrightarrow w_i^1\wedge \cdots \wedge w_i^{n-1} \wedge Q=0 \Longrightarrow a_x=0 \Longrightarrow Q\in \Gamma\left(U_x,  \mathcal{N}_{\mathcal{F}}^* \right) 
\end{align*}
Then $Q\in \Gamma\left(U_i-S, \mathcal{N}_\mathcal{F}^* \right)$, so that $Q\in \Gamma\left(U_i, \mathcal{N}_\mathcal{F}^* \right)$. This implies that $\mathcal{N}_\mathcal{F}^* = \left(\frac{T_{\mathcal{X}/A}}{T_{\mathcal{F}}} \right)^*$.

Next we will show that $T_\mathcal{F} = \left( \frac{\Omega_{\mathcal{X}/A}^1}{\mathcal{N}_\mathcal{F}^*} \right)^*$. We note that $T_\mathcal{F} \subset \left( \frac{\Omega_{\mathcal{X}/A}^1}{\mathcal{N}_{\mathcal{F}}^*} \right)^*$ by our construction. In fact,  
\begin{align*}
i_{\sum_{\alpha=1}^n A_i^\alpha \frac{\partial}{\partial z_i^\alpha}}\left( w_i^\beta \right)= \sum_{\alpha=1}^n i_{\frac{\partial}{\partial z_i^1}\wedge \cdots \wedge \frac{\partial}{\partial z_i^n}} \left( w_i^1 \wedge \cdots \wedge w_i^{n-1} \wedge dz_i^\alpha\right)i_{\frac{\partial}{\partial z_i^\alpha}}\left(w_i^\beta \right) = i_{\frac{\partial}{\partial z_i^1} \wedge \cdots \wedge \frac{\partial}{\partial z_i^n}}\left(w_i^1\wedge \cdots \wedge w_i^{n-1} \wedge w_i^\beta\right) = 0
\end{align*}

Conversely, assume that $V\in \Gamma\left(U_i, T_{\mathcal{X}/A} \right)= \Gamma(U_i, T_X)\otimes A$ such that $i_V(w_i^\alpha)=0, \alpha=1,...,n-1$, i.e. $V\in \left(\frac{\Omega_{\mathcal{X}/A}^1}{\mathcal{N}_{\mathcal{F}}^*} \right)^*\subset T_{\mathcal{X}/A}$. Then we will show that $V\in \Gamma\left( U_i, T_{\mathcal{F}} \right)$. In fact, since $\frac{T_{\mathcal{X}/A}}{T_{\mathcal{F}}}$ is flat over $A$ and $\frac{T_X}{T_{\mathcal{F}_0}}$ is locally free on $U_i-S$, we see that $\frac{T_{\mathcal{X}/A}}{T_{\mathcal{F}}}$ is locally free on $U_i-S$. Then from $(\ref{nt6})$ and $(\ref{nt7})$, we have the exact sequence
\begin{align*}
0 \to \left( \frac{\Omega_{\mathcal{X}/A}^1}{\mathcal{N}_\mathcal{F}^*} \right)^*|_{U_i-S} \to T_{\mathcal{X}/A} |_{U_i-S} \to  \left(\mathcal{N}_{\mathcal{F}}^*\right)^*|_{U_i-S} = \frac{T_{\mathcal{X}/A}}{T_{\mathcal{F}}}|_{U_i-S}\to 0
\end{align*}
This implies that $T_\mathcal{F}|_{U_i-S}= \left( \frac{\Omega_{\mathcal{X}/A}^1}{\mathcal{N}_\mathcal{F}^*}\right)^* |_{U_i-S}$. Then $V\in \Gamma\left( U_i-S, T_\mathcal{F} \right)$, so that $V\in \Gamma\left(U_i, T_\mathcal{F} \right)$. This implies that we have $ \left( \frac{\Omega_{\mathcal{X}/A}^1}{\mathcal{N}_\mathcal{F}^*} \right)^* \subset T_\mathcal{F}$. Hence $T_\mathcal{F} = \left( \frac{\Omega_{\mathcal{X}/A}^1}{\mathcal{N}_\mathcal{F}^*} \right)^*$.

Hence a foliated deformation $(\mathcal{X}, \mathcal{N}_\mathcal{F}^*)$ of $(X, \mathcal{N}_{\mathcal{F}_0}^*)$ of codimension $n-1$ also defines a foliated deformation $(\mathcal{X}, T_{\mathcal{F}})$ of $(X, T_{\mathcal{F}_0})$ of dimension $1$, so that we have a morphism of functors
\begin{align*}
\textnormal{FDef}_{\left(X, \mathcal{N}_{\mathcal{F}_0}^* \right)} \to \textnormal{FDef}_{\left(X, T_{\mathcal{F}_0}\right)}
\end{align*}

Now we describe the relationship between $\textnormal{FDef}_{\left(X, \mathcal{N}_{\mathcal{F}_0}^* \right)}$ and $\textnormal{FDef}_{\left(X, T_{\mathcal{F}_0}\right)}$. We shall describe a morphism of first-order deformations and obstructions between two functors. We recall section \ref{nt5}. We keep the notation in section \ref{nt5}. We note that foliated deformations of $\left(X, \mathcal{N}_{\mathcal{F}_0}^* \right)$ is controlled by the following complex of sheaves in terms of deformations of the section $w_0:=\left\{ w_{0i}^1\wedge \cdots \wedge w_{0i}^{n-1} \right\} \in H^0\left( X, \bigwedge^{n-1} \Omega_X^1 \otimes \mathcal{L}_0 \right)$
\begin{align}
\mathcal{E}_{\mathcal{L}_0}^\bullet : \mathcal{E}_{\mathcal{L}_0} \xrightarrow{\tilde{E}_0} \mathcal{S}^1 \to 0 \to 0 \to \cdots
\end{align}

On the other hand, we recall \ref{nt1}. We keep the notations in section \ref{nt1}. We note that foliated deformations of $\left( X , T_{\mathcal{F}_0} \right)$ is controlled by the following complex of sheaves in terms of deformations of the section $T_0= H^0\left( X, T_X \otimes \mathcal{Q}_0 \right)$, where $Q_0= K_X\otimes \mathcal{L}_0$
\begin{align*}
\mathcal{E}_{\mathcal{Q}_0}^\bullet: \mathcal{E}_{K_X\otimes \mathcal{L}_0}\xrightarrow{\tilde{D}_0} T_X\otimes K_X\otimes \mathcal{L}_0\to 0 \to 0 \to \cdots
\end{align*}
The relationship between $\textnormal{FDef}_{\left(X, \mathcal{N}_{\mathcal{F}_0}^* \right)}$ and $\textnormal{FDef}_{\left(X, T_{\mathcal{F}_0}\right)}$ is described by  the following exact sequence of complexes of sheaves by defining $\mathcal{E}_{\mathcal{L}_0}^\bullet \to \mathcal{E}_{K_X\otimes \mathcal{L}_0}^\bullet$
\begin{center}
$\begin{CD}
@. 0 @. 0 @. 0 \\
@. @AAA @AAA @AAA \\
0@>>> \mathcal{S}^1\subset \bigwedge^{n-1}\Omega_X^1\otimes \mathcal{L}_0 @>\pi>> T_X \otimes \mathcal{L}_0 \otimes K_X @>>> \mathcal{G} @>>> 0 \\
@. @A\tilde{E}_0AA @A\tilde{D}_0 AA  @AAA \\
0 @>>> \mathcal{E}_{\mathcal{L}_0} @>\eta >\cong > \mathcal{E}_{\mathcal{L}_0 \otimes K_X} @>>> 0 @>>> 0
\end{CD}$
\end{center}
We will show that the above diagram commutes and is exact. First we define $\pi:\mathcal{S}^1\to  T_X\otimes K_X \otimes \mathcal{L}_0$ is the restriction of the isomorphism $\bigwedge^{n-1} \Omega_X^1\otimes \mathcal{L}_0\cong T_X\otimes K_X \otimes \mathcal{L}_0$ to $\mathcal{S}^1$. Next we define an isomorphism 
\begin{align*}
\mathcal{E}_{\mathcal{L}_0 } \xrightarrow{\eta}  \mathcal{E}_{\mathcal{L}_0  \otimes K_X}
\end{align*}
which is locally defined in the following way:
\begin{align*}
\Gamma(U_i, \mathcal{O}_X )\bigoplus \Gamma(U_i, T_X ) &\xrightarrow{\eta}  \Gamma(U_i, \mathcal{O}_X) \bigoplus \Gamma(U_i,  T_X|_{U_i})\\
\left(a_i,b_i=\sum_{\alpha=1}^n b_i^\alpha \frac{\partial}{\partial z_i^\alpha}\right) &\mapsto \left( \sum_{\alpha=1}^n \frac{\partial b_i^\alpha}{\partial z_i^\alpha} + a_i, b_i     \right)
\end{align*}

We show that $\eta$ is well-defined. In other words, it is independent of local coordinates, i.e. the following diagram commutes
\begin{center}
$\begin{CD}
(a_i+[b_i, \log h_{0ij}], b_i)\in \mathcal{O}_X|_{U_j} \bigoplus T_X|_{U_j} @>\eta >> (a_i'+ \left[b_i, \log h_{0ij}\right] + \left[b_i, \log R_{0ij} \right] \in \mathcal{O}_X|_{U_j} \bigoplus T_X|_{U_j} \\
@AAA  @AAA \\
(a_i,b_i) \in \mathcal{O}_X|_{U_i} \bigoplus T_X|_{U_i} @>\eta>> (a_i', b_i) \in \mathcal{O}_X|_{U_i}\bigoplus T_X|_{U_i}
\end{CD}$
\end{center}
We set $b_i=\sum_{\gamma=1}^n b_i^\gamma \frac{\partial}{\partial z_i^\gamma}= \sum_{\gamma, \beta=1}^n b_i^\gamma \frac{\partial z_j^\beta}{\partial z_i^\gamma}\frac{\partial}{\partial z_j^\beta}$. Then we have to show that
\begin{align*}
a_i+ \sum_{\alpha=1}^n \frac{\partial b_i^\alpha}{\partial z_i^\alpha} + \left[b_i, \log h_{0ij} \right] +\left[b_i,  \log R_{0ij} \right] = a_i+ \left[b_i, \log h_{0ij} \right] +\sum_{\beta=1}^n \frac{\partial b_j^\beta}{\partial z_j^\beta} \iff   \sum_{\alpha=1}^n \frac{\partial b_i^\alpha}{\partial z_i^\alpha} +\left[b_i, \log R_{0ij} \right] =  \sum_{\beta=1}^n \frac{\partial b_j^\beta}{\partial z_j^\beta}
\end{align*}
which follows from $(\ref{nt8})$ since $R_{0ij}= \frac{1}{K_{0ij}}$. Hence $\eta$ is well-defined. We remark taht $\mathcal{E}_{\mathcal{L}_0}\to \mathcal{E}_{\mathcal{L}_0 \otimes K_X}$ is not a $\mathcal{O}_X$-module homomorphism.  Next we show that we have the following commutative diagram
\begin{center}
$\begin{CD}
\bigwedge^{n-1} \Omega_X^1 \otimes \mathcal{L}_0 @>\pi>\cong > T_X \otimes \mathcal{L}_0 \otimes K_X\\
@A\tilde{E}_0 AA @AA\tilde{D}_0 A \\
\mathcal{E}_{\mathcal{L}_0} @>\eta>\cong > \mathcal{E}_{\mathcal{L}_0\otimes K_X}
\end{CD}$
\end{center}
In fact, we have
\begin{align*}
&\tilde{E}_0(c_i, b_i)\wedge dz_i^\alpha= \mathcal{L}_{b_i}\left(w_{0i}^1\wedge \cdots \wedge w_{0i}^{n-1}\right)\wedge dz_i^\alpha + c_i w_{0i}^1\wedge \cdots \wedge w_{0i}^{n-1}\wedge dz_i^\alpha \\
&=\mathcal{L}_{b_i} \left(w_{0i}^1\wedge \cdots \wedge w_{0i}^{n-1}\wedge dz_i^\alpha \right) - w_{0i}^1\wedge \cdots \wedge w_{0i}^{n-1}\wedge \mathcal{L}_{b_i}\left(dz_i^\alpha\right) + c_ia_i^\alpha dz_i^1\wedge \cdots \wedge dz_i^n\\
&=\mathcal{L}_{b_i}\left(a_i^\alpha dz_i^1\wedge \cdots \wedge dz_i^n \right) - \sum_{\gamma=1}^n \frac{\partial b_i^\alpha}{\partial z_i^\gamma} a_i^\gamma dz_i^1\wedge \cdots \wedge dz_i^n + c_i a_i^\alpha dz_i^1\wedge \cdots \wedge dz_i^n\\
&= \sum_{\gamma=1}^n b_i^\gamma \frac{\partial a_i^\alpha}{\partial z_i^\gamma}dz_i^1\wedge \cdots \wedge dz_i^n + a_i^\alpha \sum_{\beta=1}^n \frac{\partial b_i^\beta}{\partial z_i^\beta} dz_i^1\wedge \cdots \wedge dz_i^n +\left(c_ia_i^\alpha - \sum_{\gamma=1}^n a_i^\gamma \frac{\partial b_i^\alpha}{\partial z_i^\gamma} \right) dz_i^1\wedge \cdots \wedge dz_i^n
\end{align*}
Then
\begin{align*}
\pi \tilde{E}_0\left(c_i, b_i\right)= \left[ b_i, T_{0i} \right]  +\left(c_i +\sum_{\beta=1}^n \frac{\partial b_i^\beta}{\partial z_i^\beta} \right) T_{0i}=\tilde{D}_0 \left(  c_i+ \sum_{\beta=1}^n \frac{\partial b_i^\beta}{\partial z_i^\beta}, b_i\right) = \tilde{D}_0 \eta \left(c_i, b_i\right) 
\end{align*}
Hence the diagram commutes. Then we have the long exact sequence
\begin{align*}
0\to \mathbb{H}^0\left( X, \mathcal{E}_{\mathcal{L}_0}^{\bullet} \right) \to \mathbb{H}^0 \left(X, \mathcal{E}_{\mathcal{L}_0\otimes K_X}^\bullet \right) \to 0 \to \mathbb{H}^1\left(X, \mathcal{E}_{\mathcal{L}_0 }^\bullet \right) \xrightarrow{\alpha} \mathbb{H}^1\left( X, \mathcal{E}_{\mathcal{L}_0 \otimes K_X}^\bullet \right) \to H^0\left(X, \mathcal{G}\right) \\
\to \mathbb{H}^2\left( X, \mathcal{E}_{\mathcal{L}_0}^\bullet \right) \xrightarrow{\beta} \mathbb{H}^2\left( X, \mathcal{E}_{\mathcal{L}_0 \otimes K_X}^\bullet \right) \to H^1\left(X, \mathcal{G} \right)\to \cdots
\end{align*}
We will show that $\alpha$ gives morphism of first-order deformations
\begin{align}\label{nt16}
\textnormal{FDef}_{\left(X,\mathcal{N}_{\mathcal{F}_0}^* \right)}\left(k[\epsilon]\right) \cong \mathbb{H}^1\left(X, \mathcal{E}_{\mathcal{L}_0}^\bullet \right) \hookrightarrow \textnormal{FDef}_{\left(X, T_{\mathcal{F}_0}\right)}\left(k[\epsilon]\right) \cong  \mathbb{H}^1\left( X, \mathcal{E}_{\mathcal{L}_0\otimes K_X}^\bullet \right)
\end{align}
and $\beta$ gives the obstruction map
\begin{align} \label{nt17}
\mathbb{H}^2\left(X,  \mathcal{E}_{\mathcal{L}_0}^\bullet  \right) \xrightarrow{\beta} \mathbb{H}^2\left(X, \mathcal{E}_{\mathcal{L}_0\otimes K_X}^\bullet \right)
\end{align}
between $\textnormal{FDef}_{\left( X , \mathcal{N}_{\mathcal{F}_0}^* \right)}$ and $\textnormal{FDef}_{\left( X, T_{\mathcal{F}_0} \right)}$.

Let $\left(\mathcal{X}, \mathcal{N}_\mathcal{F}^* \right)$ be a first-order foliated deformation of $\left( X, \mathcal{N}_{\mathcal{F}_0}^* \right)$ over $k[\epsilon]$ in terms of cotangent sheaves. We keep the notations in the proof of Proposition \ref{si39}. We choose a Stein open covering $\mathcal{U}=\left\{ U_i \right\}$ of $X$ with local coordinates $\left( z_i^1,..., z_i^n \right)$ on $U_i$ such that $\Gamma\left( U_i, \mathcal{N}_{\mathcal{F}_0}^* \right)$ is generated by $w_{0i}^1,..., w_{0i}^{n-1}\in \Gamma\left( U_i, \Omega_X^1 \right)$ with $w_{0i}^\alpha= \sum_{\beta=1}^{n-1} h_{0ij}^{\alpha\beta} w_{0j}^\beta$, and $\Gamma\left( U_i, \mathcal{N}_\mathcal{F}^*  \right)$ is generated by $w_{0i}^\alpha + \epsilon w_i'^\alpha,\alpha=1,...,n-1$ and $\theta_{ij}$ corresponds to $\textnormal{Id}+ \epsilon p_{ij}$ with $p_{ij} = \sum_{\gamma=1}^n p_{ij}^\gamma \frac{\partial}{\partial z_i^\gamma} \in \Gamma\left( U_i, T_X \right)$ and
\begin{align*}
\theta_{ji}^*\left( w_{0i}^\alpha + \epsilon w_i'^\alpha \right)  = w_{0i}^\alpha+ \epsilon w_i'^\alpha- \epsilon \mathcal{L}_{p_{ij}}\left( w_{0i}^\alpha \right) = \sum_{\beta=1}^{n-1}\left( h_{0ij}^{\alpha\beta} + \epsilon g_{ij}^{\alpha\beta} \right)\left( w_{0j}^\beta + \epsilon w_j'^\beta \right)
\end{align*}
We set $\det\left( h_{0ij}^{\alpha\beta} + \epsilon g_{ij}^{\alpha\beta} \right)= h_{0ij}+ \epsilon G_{ij}$. Then $\left( \left\{ w_i'^*=\sum_{\alpha=1}^{n-1} w_{0i}^1\wedge \cdots \wedge w_i'^\alpha \wedge \cdots \wedge w_{0i}^{n-1}\right\} , \left\{ \left( p_{ij}, \frac{G_{ij}}{h_{0ij}}\right) \right\} \right) \in C^0\left( \mathcal{U}, \mathcal{S}^1 \right) \bigoplus C^1\left( \mathcal{U}, \mathcal{E}_{\mathcal{L}_0} \right)$ defines a $1$-cocycle of \v Cech resolution of $\mathcal{E}_{\mathcal{L}_0}^\bullet$ in terms of the covering $\mathcal{U}$.

On the other hand, we describe the induced first-order foliated deformation $\left( \mathcal{X}, T_\mathcal{F} \right)$ in terms of tangent sheaves from $\left( \mathcal{X}, \mathcal{N}_\mathcal{F}^* \right)$. We keep the notations of the proof of Proposition \ref{ts11}. Then $\Gamma\left( U_i, T_\mathcal{F} \right)$ is generated by
\begin{align*}
&\sum_{\alpha=1}^n i_{\frac{\partial}{\partial z_i^1}\wedge \cdots \wedge \frac{\partial}{\partial z_i^n} }\left( \left( w_{0i}^1+ \epsilon w_i'^1 \right)\wedge \cdots \wedge \left( w_{0i}^{n-1} + \epsilon w_i'^{n-1} \right)  \wedge dz_i^\alpha    \right) \frac{\partial}{\partial z_i^\alpha}\\
&=\sum_{\alpha=1}^n i_{\frac{\partial}{\partial z_i^1}\wedge \cdots \wedge \frac{\partial}{\partial z_i^n} }\left( \left( w_{0i}^1\wedge \cdots \wedge w_{0i}^{n-1}\right)  \wedge dz_i^\alpha    \right) \frac{\partial}{\partial z_i^\alpha} + \epsilon  \sum_{\alpha=1}^n i_{\frac{\partial}{\partial z_i^1}\wedge \cdots \wedge \frac{\partial}{\partial z_i^n}}\left( w_i'^* \wedge dz_i^\alpha \right)\frac{\partial}{\partial z_i^\alpha}= T_{0i} + \epsilon T_i'
\end{align*}
and
\begin{align*}
\theta_{ji}^*\left( T_{0i}+ \epsilon T_i'\right) & = T_{0i} + \epsilon T_i' -\epsilon \left[ p_{ij}, T_{0i}\right] = \left(R_{0ij} + \epsilon M_{ij} \right)\left( h_{0 ij} + \epsilon G_{ij}   \right) \left( T_{0j} + \epsilon T_j' \right) \\
& = \left( R_{0ij} h_{0ij} +\epsilon\left( M_{ij} h_{0ij} + R_{ij} G_{ij} \right)\right)\left( T_{0j} + \epsilon T_j' \right)
\end{align*}
where $M_{ij}$ is defined by the following property
{\small{\begin{align}\label{nt15}
\left( R_{0ij} + \epsilon M_{ij} \right) \frac{\partial}{\partial z_j^1}\wedge \cdots \wedge \frac{\partial}{\partial z_j^n} =  \theta_{ji}^*\left( \frac{\partial}{\partial z_i^1} \right)  \wedge \cdots \wedge \theta_{ji}^*\left( \frac{\partial}{\partial z_i^n} \right) =  \left( 1+ \epsilon \sum_{\beta=1}^n \frac{\partial p_{ij}^\beta}{\partial z_i^\beta} \right) \frac{\partial}{\partial z_i^1}\wedge \cdots \wedge \frac{\partial}{\partial z_i^n}
\end{align}}}
which follows from 
\begin{align*}
\theta_{ji}^*\left( \frac{\partial}{\partial z_i^\alpha} \right) &=  \sum_{\beta=1}^n i_{\theta_{ji}^*\left( \frac{\partial}{\partial z_i^\alpha } \right)}\left( dz_i^\beta \right)\frac{\partial}{\partial z_i^\beta} = \sum_{\beta=1}^n \theta_{ji}^*\left( i_{\frac{\partial}{\partial z_i^\alpha}}\left(d \theta_{ij}^*\left(z_i^\beta \right) \right) \right)\frac{\partial}{\partial z_i^\beta}= \sum_{\beta=1}^n \theta_{ji}^*\left( i_{\frac{\partial}{\partial z_i^\alpha}}\left( dz_i^\beta + \epsilon d p_{ij}^\beta \right) \right)\frac{\partial}{\partial z_i^\beta} \\
&=\frac{\partial}{\partial z_i^\alpha} + \epsilon \sum_{\beta=1}^n \frac{\partial p_{ij}^\beta}{\partial z_i^\alpha} \frac{\partial}{\partial z_i^\beta}
\end{align*}
Then $(\ref{nt15})$ implies that $\frac{M_{ij}}{R_{0ij}}= \sum_{\beta=1}^n \frac{\partial p_{ij}^\beta}{\partial z_i^\beta}$. Then $\left( \left\{ T_i'\right\}, \left\{\left( p_{ij} , \frac{M_{ij}}{R_{0ij}} + \frac{G_{ij}}{ h_{0ij} } = \frac{G_{ij}}{h_{0ij}} + \sum_{\beta=1}^n \frac{\partial p_{ij}^\beta}{\partial z_i^\beta}  \right) \right\} \right) \in C^0\left(\mathcal{U}, \mathcal{R} \right) \bigoplus C^1\left( \mathcal{U}, \mathcal{E}_{\mathcal{L}_0\otimes K_X} \right)$ in terms of the covering $\mathcal{U}$ by the proof of Proposition \ref{ts11}. Then it is clear that under the map from $(\ref{nt16})$
\begin{align*}
\mathbb{H}^1\left( X, \mathcal{E}_{\mathcal{L}_0}^\bullet \right) &\hookrightarrow \mathbb{H}^1\left( X, \mathcal{E}_{\mathcal{L}_0\otimes K_X}^\bullet \right) \\
\left(\left\{ w_i'^*\right\}, \left\{ \left( p_{ij}, \frac{G_{ij}}{h_{0ij}} \right)\right\} \right) &\mapsto \left( \left\{ T_i' \right\}, \left\{ \left(  p_{ij}, \frac{G_{ij}}{h_{0ij}} + \sum_{\beta=1}^n \frac{\partial p_{ij}^\beta}{\partial z_i^\beta}   \right)   \right\} \right)
\end{align*}

Next we describe the obstruction map between $\textnormal{FDef}_{\left( X, \mathcal{N}_{\mathcal{F}_0}^* \right)}$ and $\textnormal{FDef}_{\left( X, T_{\mathcal{F}_0} \right)}$ for a small extension $0\to (t) \to \tilde{A} \to A \to 0$. Let $\left( \mathcal{X}, \mathcal{N}_\mathcal{F}^* \right)$ be a foliated deformation of $\left( X, \mathcal{N}_{\mathcal{F}_0}^* \right)$ over $A$. We keep the notations in the proof of Proposition \ref{si39}. Let $\mathcal{U}=\left\{ U_i \right\}$ be a Stein open covering of $X$ with local coordinates $\left( z_i^1 ,..., z_i^n \right)$ on $U_i$ such that $\theta_i :\left( U_i \times \textnormal{Spec}(A) , \mathcal{N}_{\mathcal{F}_i}^* \right) \to \left( \mathcal{X}|_{U_i}, \mathcal{N}_\mathcal{F}^*|_{U_i} \right)$ is a foliated isomorphism, and $\theta_{ij}:= \theta_j^{-1}\theta_i$, and $\Gamma\left( U_i, \mathcal{N}_{\mathcal{F}_i}^* \right)$ is generated by $w_i^1,..., w_i^{n-1}\in \Gamma\left( U_i, \Omega_X^1 \right)\otimes A$ with $\theta_{ji}^* w_i^\alpha= \sum_{\beta=1}^{n-1} h_{ij}^{\alpha\beta} w_j^\beta$ for $h_{ij}^{\alpha\beta}\in \Gamma\left( U_{ij}, \mathcal{O}_X \right)\otimes A$. Then we choose arbitrary automorphism $\left\{ \tilde{\theta}_{ij} \right\}$ which induces $\theta_{ij}$, and liftings $\tilde{w}_i^\alpha$ of $w_i^\alpha$, liftings $\tilde{h}_{ij}^{\alpha\beta}$ of $h_{ij}^{\alpha\beta}$. Then $\left( \left\{ - R_{ij}^* \right\} , \left\{ \left( d_{ijk}, - h_{ijk}^* \right) \right\} \right)\in C^1\left( \mathcal{U}, \mathcal{S}^1 \right) \bigoplus C^0\left( \mathcal{U}, \mathcal{E}_{\mathcal{L}_0} \right)$ defines a $2$-cocycle of $\mathcal{E}_{\mathcal{L}_0}^\bullet$, where $R_{ij}^*=\sum_{\alpha=1}^{n-1} w_{0i}^1\wedge \cdots \wedge R_{ij}^\alpha \wedge \cdots \wedge w_{0i}^{n-1} \in \Gamma\left( U_i, \mathcal{S}^1 \right)$ such that $\tilde{\theta}_{ji}^*\left( \tilde{w}_i^\alpha \right)  - \sum_{\beta=1}^{n-1} \tilde{h}_{ij}^{\alpha\beta} \tilde{w}_j^\beta = t R_{ij}^\alpha$, and $d_{ijk} \in \Gamma\left( U_{ijk} , T_X \right)$ such that $\textnormal{Id} +t d_{ijk}= \tilde{f}_{ij}\tilde{f}_{jk}\tilde{f}_{ki}$, and $h_{ijk}^*\in \Gamma\left( U_{ijk} , \mathcal{O}_X \right)$ such that $1+t h_{ijk}^* = \tilde{f}_{ij} \tilde{f}_{jk}\left( \tilde{h}_{ik}\right) \cdot \tilde{f}_{ij}\left( \tilde{h}_{kj} \right) \cdot \tilde{h}_{ji}$ where $\tilde{h}_{ij} = \det\left( \tilde{h}_{ij}^{\alpha\beta} \right)$.

Now let us describe the induced foliated deformation $\left( \mathcal{X}, T_\mathcal{F} \right)$ over $A$ from $\left( \mathcal{X}, \mathcal{N}_\mathcal{F}^* \right)$. We keep the notations in the proof of Proposition \ref{ts11}. Then $\Gamma\left( U_i, T_{\mathcal{F}_i} \right)$ is generated by $T_i=\sum_{\alpha=1}^n i_{\frac{\partial}{\partial z_i^1}\wedge \cdots \wedge \frac{\partial}{\partial z_i^n}}\left( \tilde{w}_i^1 \wedge \cdots \wedge \tilde{w}_i^{n-1} \wedge dz_i^\alpha  \right)\frac{\partial}{\partial z_i^\alpha}$, and
\begin{align*}
\theta_{ji}^*\left( T_i \right) & = \theta_{ji}^* \left(  \sum_{\alpha=1}^n  i_{\frac{\partial}{\partial z_i^1}\wedge \cdots \wedge \frac{\partial}{\partial z_i^n}}\left( w_i^1 \wedge \cdots \wedge w_i^{n-1} \wedge dz_i^\alpha  \right)\frac{\partial}{\partial z_i^\alpha} \right) \\
& = R_{ij} h_{ij} \sum_{\beta=1}^n i_{\frac{\partial}{\partial z_j^1} \wedge \cdots \wedge \frac{\partial}{\partial z_j^n} }\left( w_j^1 \wedge \cdots \wedge w_j^{n-1} \wedge dz_j^\beta \right) \frac{\partial}{\partial z_j^\beta} = R_{ij} h_{ij} T_j :=  r_{ij} T_j , \,\,\,\,\,\,\,\,\textnormal{where}\,\,\, r_{ij}= R_{ij} h_{ij}
\end{align*}
where $R_{ij} \in \Gamma\left( U_{ij} , \mathcal{O}_{\mathcal{X}}^* \right)$ is defined by the following property
\begin{align*}
\theta_{ji}^*\left( \frac{\partial}{\partial z_i^1}\wedge \cdots \wedge \frac{\partial}{\partial z_i^n} \right)= R_{ij} \frac{\partial}{\partial z_j^1}\wedge \cdots \wedge \frac{\partial}{\partial z_j^n}
\end{align*}
Then the liftings $\tilde{w}_i^\alpha$ corresponds to the lifting $\tilde{T}_i= i_{\frac{\partial}{\partial z_i^1}\wedge \cdots \wedge \frac{\partial}{\partial z_i^n}}\left( \tilde{w}_i^1 \wedge \cdots \wedge \tilde{w}_i^{n-1} \wedge dz_i^\alpha  \right)\frac{\partial}{\partial z_i^\alpha} \in \Gamma\left( U_i, T_X \right) \otimes \tilde{A}$ of $\tilde{T}_i$, and we note that since $\tilde{f}_{ij}$ is an isomorphism, there exists a lifting $\tilde{R}_{ij}\in \Gamma\left( U_{ij} , \mathcal{O}_X \right) \otimes \tilde{A}$ of $R_{ij}$ such that
\begin{align*}
\tilde{\theta}_{ji}^*\left( \frac{\partial}{\partial z_i^1}\wedge \cdots \wedge \frac{\partial}{\partial z_i^n} \right)= \tilde{R}_{ij} \frac{\partial}{\partial z_j^1}\wedge \cdots \wedge \frac{\partial}{\partial z_j^n}.
\end{align*}
Then we have
\begin{align*}
&\tilde{\theta}_{ji}^*\left( \tilde{T}_i \right)=\tilde{\theta}_{ji}^*\left(\sum_{\alpha=1}^n i_{\frac{\partial}{\partial z_i^1}\wedge \cdots \wedge \frac{\partial}{\partial z_i^n}}\left( \tilde{w}_i^1\wedge \cdots \wedge \tilde{w}_i^{n-1} \wedge dz_i^\alpha \right) \frac{\partial}{\partial z_i^\alpha} \right)\\
&= R_{ij} \sum_{\alpha=1}^n i_{\frac{\partial}{\partial z_j^1}\wedge \cdots \frac{\partial}{\partial z_j^n}} \left( \tilde{\theta}_{ji}^*\left(\tilde{w}_i^1\wedge \cdots \wedge \tilde{w}_i^{n-1} \right) \wedge \tilde{\theta}_{ji}^*\left( dz_i^\alpha \right)   \right)\tilde{\theta}_{ji}^*\left( \frac{\partial}{\partial z_i^\alpha} \right) \\
&= \tilde{R}_{ij}\tilde{h}_{ij} \sum_{\alpha=1}^n i_{\frac{\partial}{\partial z_j^1}\wedge \cdots \wedge \frac{\partial}{\partial z_j^n}}\left( \tilde{w}_j^1\wedge \cdots \wedge \tilde{w}_j^{n-1} \wedge d\tilde{\theta}_{ji}^* \left(z_i^\alpha\right) \right)\left( \sum_{\beta=1}^n i_{\tilde{\theta}_{ji}^*\left( \frac{\partial}{\partial z_i^\alpha} \right)}\left( dz_j^\beta \right)\frac{\partial}{\partial z_j^\beta} \right) + t \sum_{\alpha=1}^n i_{\frac{\partial}{\partial z_i^1 } \wedge \cdots  \wedge \frac{\partial}{\partial z_i^n}}\left( w_i'^*\wedge dz_i^\alpha \right)\frac{\partial}{\partial z_i^\alpha} \\
&=  \tilde{R}_{ij} \tilde{h}_{ij} \tilde{T}_j + t \sum_{\alpha=1}^n i_{\frac{\partial}{\partial z_i^1 } \wedge \cdots  \wedge \frac{\partial}{\partial z_i^n}}\left( w_i'^*\wedge dz_i^\alpha \right)\frac{\partial}{\partial z_i^\alpha}
\end{align*}
Then by the proof of Proposition \ref{ts11}, we see that $\left( \left\{ - W_{ij}^* \right\}, \left\{ \left( d_{ijk}, - g_{ijk}^*  \right) \right\}  \right) \in C^1\left( \mathcal{U}, \mathcal{R} \right) \bigoplus C^2\left( \mathcal{U}, \mathcal{E}_{\mathcal{L}_0\otimes K_X} \right)$ defines a $2$-cocycle in the \v Cech resolution of $\mathcal{E}_{\mathcal{L}_0\otimes K_X}^\bullet$, where $W_{ij}^*= \sum_{\alpha=1}^n i_{\frac{\partial}{\partial z_i^1 } \wedge \cdots  \wedge \frac{\partial}{\partial z_i^n}}\left( w_i'^*\wedge dz_i^\alpha \right)\frac{\partial}{\partial z_i^\alpha}$, and $d_{ijk} \in \Gamma\left(U_{ijk}, T_X \right)$ such that $\textnormal{Id} + t d_{ijk} = \tilde{f}_{ij} \tilde{f}_{jk} \tilde{f}_{ki}$, and $g_{ijk}^*\in \Gamma\left( U_{ijk} , \mathcal{O}_X \right)$ such that $1+ t g_{ijk}^*=\tilde{f}_{ij} \tilde{f}_{jk}\left( \tilde{h}_{ik} \tilde{R}_{ik} \right) \cdot \tilde{f}_{ij}\left( \tilde{h}_{kj} \tilde{R}_{kj}\right) \cdot \tilde{h}_{ji} \tilde{R}_{ji}$. Then we have to show that under the map $\beta$ from $(\ref{nt17})$
\begin{align*}
\beta : \mathbb{H}^2\left( X, \mathcal{E}_{\mathcal{L}_0}^\bullet \right) &\to \mathbb{H}^2\left( X, \mathcal{E}_{\mathcal{L}_0\otimes K_X}^\bullet \right) \\
\left( \left\{ - R_{ij}^* \right\}, \left\{\left( d_{ijk}, - h_{ijk}^* \right) \right\} \right) &\mapsto \left( \left\{ - W_{ij}^* \right\}, \left\{ \left( d_{ijk} ,  - g_{ijk}^* \right) \right\} \right)
\end{align*} 
In fact, it is clear that $- R_{ij}^*$ is mapped to $ - W_{ij}^*$, and so it remains to show that $\left( d_{ijk}, - h_{ijk}^* \right)$ is mapped to $\left( d_{ijk} , - g_{ijk}^* \right)$. First we note that
{\small{\begin{align} \label{nt19}
\tilde{\theta}_{ij}^*\tilde{\theta}_{jk}^*\tilde{\theta}_{ki}^*\left( \frac{\partial}{\partial z_i^\alpha} \right) &= \sum_{\beta=1}^n i_{ \tilde{\theta}_{ij}^*\tilde{\theta}_{jk}^*\tilde{\theta}_{ki}^* \left( \frac{\partial}{\partial z_i^\alpha} \right)} \left( dz_i^\beta \right)\frac{\partial}{\partial z_i^\beta} = \sum_{\beta=1}^n \tilde{\theta}_{ij}^*\tilde{\theta}_{jk}^*\tilde{\theta}_{ki}^* \left( i_{\frac{\partial}{\partial  z_i^\alpha}}\left( d \tilde{\theta}_{ik}^*\tilde{\theta}_{kj}^*\tilde{\theta}_{ji}^* \left(z_i^\beta\right) \right)   \right)\frac{\partial}{\partial z_i^\beta} \\
&= \sum_{\beta=1}^n \tilde{\theta}_{ij}^*\tilde{\theta}_{jk}^*\tilde{\theta}_{ki}^* \left( i_{\frac{\partial}{\partial z_i^\alpha}}\left( dz_i^\beta - td \left( d_{ijk}^\beta\right) \right) \right)\frac{\partial}{\partial z_i^\beta}  = \sum_{\beta=1}^n \tilde{\theta}_{ij}^*\tilde{\theta}_{jk}^* \tilde{\theta}_{ki}^*\left( i_{\frac{\partial}{\partial z_i^\alpha}}\left( dz_i^\beta \right) \right) \frac{\partial}{\partial z_i^\beta} -t \sum_{\beta=1}^n i_{\frac{\partial}{\partial z_i^\alpha}}\left( d\left( d_{ijk}^\beta \right) \right)\frac{\partial}{\partial z_i^\beta} \notag \\
&= \frac{\partial}{\partial z_i^\alpha} - t \sum_{\beta=1}^n i_{\frac{\partial}{\partial z_i^\alpha}}\left( d\left( d_{ijk}^\beta \right) \right) \frac{\partial}{\partial z_i^\beta} \notag
\end{align}}}

We also note that $\tilde{f}_{ij} \tilde{f}_{jk}\tilde{f}_{ki} = \textnormal{Id} + t d_{ijk}$, and by setting $d_{ijk}=\sum_{\gamma=1}^n d_{ijk}^\gamma \frac{\partial}{\partial z_i^\gamma}$, we have from $(\ref{nt19})$
\begin{align*}
&\tilde{R}_{ik}\frac{\partial}{\partial z_k^1}\wedge \cdots \wedge \frac{\partial}{\partial z_k^n}= \tilde{\theta}_{ki}^*\left( \frac{\partial}{\partial z_i^1}\wedge \cdots \wedge \frac{\partial}{\partial z_i^n} \right)  \\
&\Longrightarrow \tilde{f}_{jk}\left( \tilde{R}_{ik} \right) \tilde{R}_{kj}\frac{\partial}{\partial z_j^1}\wedge \cdots \wedge \frac{\partial}{\partial z_j^n} =  \tilde{\theta}_{jk}^* \tilde{\theta}_{ki}^*\left( \frac{\partial}{\partial z_i^1}\wedge \cdots \wedge \frac{\partial}{\partial z_i^n} \right) \\
&\Longrightarrow \tilde{f}_{ij} \tilde{f}_{jk}\left( \tilde{R}_{ik} \right)\cdot \tilde{f}_{ij}\left(\tilde{R}_{kj} \right)\cdot \tilde{R}_{ji}\frac{\partial}{\partial z_i^1}\wedge \cdots \wedge \frac{\partial}{\partial z_i^n} = \tilde{\theta}_{ij}^*\tilde{\theta}_{jk}^*\tilde{\theta}_{ki}^* \left( \frac{\partial}{\partial z_i^1}\wedge \cdots \wedge \frac{\partial}{\partial z_i^n} \right)\\
&\,\,\,\,\,\,\,\,\,\,\,\,\,= \left( 1 - t\sum_{\alpha=1}^n \frac{\partial d_{ijk}^\alpha}{\partial z_i^\alpha} \right) \frac{\partial}{\partial z_i^1}\wedge \cdots \wedge \frac{\partial}{\partial z_i^n}
\end{align*}
Hence we have
\begin{align} \label{nt18}
\tilde{f}_{ij} \tilde{f}_{jk}\left( \tilde{R}_{ik} \right)\cdot \tilde{f}_{ij}\left(\tilde{R}_{kj} \right)\cdot \tilde{R}_{ji} = 1- t \sum_{\alpha=1}^n \frac{\partial d_{ijk}^\alpha}{\partial z_i^\alpha}
\end{align}
Then we have from $(\ref{nt18})$
\begin{align*}
1+ t g_{ijk}^*&=\left( \tilde{f}_{ij} \tilde{f}_{jk} \right)\left(\tilde{h}_{ik} \tilde{R}_{ik} \right) \cdot \tilde{f}_{ij}\left( \tilde{h}_{kj} \tilde{R}_{kj} \right)\cdot \tilde{h}_{ji} \tilde{R}_{ji}\\
 &= \left( \tilde{f}_{ij} \tilde{f}_{jk} \right)\left(\tilde{h}_{ik} \right) \cdot \tilde{f}_{ij}\left( \tilde{h}_{kj}  \right)\cdot \tilde{h}_{ji} \cdot \left( \tilde{f}_{ij} \tilde{f}_{jk} \right)\left( \tilde{R}_{ik} \right) \cdot \tilde{f}_{ij}\left(\tilde{R}_{kj} \right)\cdot \tilde{R}_{ji} \\
 &= \left( 1+ t h_{ijk}^*\right) \left( \tilde{f}_{ij} \tilde{f}_{jk} \right)\left(\tilde{R}_{ik} \right) \cdot \tilde{f}_{ij}\left( \tilde{R}_{kj} \right)\cdot  \tilde{R}_{ji} =\left( \tilde{f}_{ij} \tilde{f}_{jk} \right)\left( \tilde{R}_{ik} \right) \cdot \tilde{f}_{ij}\left( \tilde{R}_{kj} \right)\cdot \tilde{R}_{ji} + t h_{ijk}^*\\
 &=1 + t\left(   h_{ijk}^* - \sum_{\alpha=1}^n \frac{\partial d_{ijk}^\gamma}{\partial z_i^\gamma} \right) \Longrightarrow - g_{ijk}^* = - h_{ijk}^* + \sum_{\alpha=1}^n \frac{\partial d_{ijk}^\alpha}{\partial z_i^\alpha}
\end{align*}
This proves that $\left( d_{ijk} , - h_{ijk}^* \right)$ is mapped to $\left(  d_{ijk} , - g_{ijk}^* \right)$.

\subsection{Isomorphism between $\textnormal{FDef}_{\left( X, T_{\mathcal{F}_0}, \mathcal{N}_{\mathcal{F}_0}^* \right)}$ and $\textnormal{FDef}_{\left( X, \mathcal{N}_{\mathcal{F}_0}^* \right)}$ with $\mathcal{N}_{\mathcal{F}_0}^*$ locally free of rank $n-1$}\

Let $\left( X, \mathcal{N}_{\mathcal{F}_0}^* \right)$ be a foliated nonsingular variety with $\mathcal{N}_{\mathcal{F}_0}^*$ locally free of rank $n-1$. We have shown that a foliated deformation of $\left( X, \mathcal{N}_{\mathcal{F}_0}^* \right)$ of codimension $n-1$ also induces a foliated deformation of $\left( X, T_{\mathcal{F}_0} \right)$ of dimension $1$. This implies that we have a morphism of deformation functors $\textnormal{FDef}_{\left(X, \mathcal{N}_{\mathcal{F}_0}^* \right)}\to \textnormal{FDef}_{\left( X, T_{\mathcal{F}_0}, \mathcal{N}_{\mathcal{F}_0}^* \right)}$ which is the inverse of the natural morphism of deformation functors $\textnormal{FDef}_{\left(X, T_{\mathcal{F}_0}, \mathcal{N}_{\mathcal{F}_0}^* \right)} \to \textnormal{FDef}_{\left( X, \mathcal{N}_{\mathcal{F}_0}^* \right)}$ by forgetting $T_{\mathcal{F}_0}$, so that we have
\begin{align}\label{nt25}
\textnormal{FDef}_{\left(X, T_{\mathcal{F}_0}, \mathcal{N}_{\mathcal{F}_0}^* \right)} \cong \textnormal{FDef}_{\left( X, \mathcal{N}_{\mathcal{F}_0}^* \right)}
\end{align}
In particular we have
\begin{align*}
\mathbb{H}^1\left(X, \mathcal{F}_0^\bullet \right) \cong \mathbb{H}^1\left( X, \mathcal{N}_{\mathcal{F}_0}^{*\bullet} \right)
\end{align*}
Moreover we show that
\begin{align*}
\mathbb{H}^2\left(X, \mathcal{F}_0^\bullet \right) \cong \mathbb{H}^2\left( X, \mathcal{N}_{\mathcal{F}_0}^{*\bullet} \right)
\end{align*}
In fact, we recall that we have an exact sequence of complexes of sheaves
{\tiny{\begin{center}
$\begin{CD}
@. 0 @. 0 @. 0 \\
@. @AAA @AAA @AAA \\
0 @>>> \mathscr{H}om_{\mathcal{O}_X}\left(T_{\mathcal{F}_0}, \left( \mathcal{N}_{\mathcal{F}_0}^* \right)^* \right) @>>> \mathscr{H}om_{\mathcal{O}_X}\left( T_{\mathcal{F}_0}, \left( \mathcal{N}_{\mathcal{F}_0}^* \right)^* \right) @>>> 0 @>>> 0 \\
@. @AAA @AAA @AAA \\
0 @>>> \mathscr{H}om_{\mathcal{O}_X}\left( T_{\mathcal{F}_0}, \frac{T_X}{T_{\mathcal{F}_0}} \right) @>>>\mathscr{H}om_{\mathcal{O}_X}\left( T_{\mathcal{F}_0}, \frac{T_X}{T_{\mathcal{F}_0}} \right) \bigoplus  \mathscr{H}om_{\mathcal{O}_X}\left( \mathcal{N}_{\mathcal{F}_0}^*, \frac{\Omega_X^1}{\mathcal{N}_{\mathcal{F}_0}^*} \right) @>>> \mathscr{H}om_{\mathcal{O}_X}\left( \mathcal{N}_{\mathcal{F}_0}^* , \frac{\Omega_X^1}{ \mathcal{N}_{\mathcal{F}_0}^*} \right) @>>> 0 \\
@. @AAA @AAA @AAA \\
0 @>>> 0 @>>> T_X @>>> T_X @>>> 0 
\end{CD}$
\end{center}}}
Let us denote the right complex of sheaves by $\mathcal{B}^\bullet$. Then we have the long exact sequence
\begin{align*}
0\to \mathbb{H}^0\left( X, \mathcal{F}_0^\bullet \right) \to \mathbb{H}^0\left( X, \mathcal{N}_{\mathcal{F}_0}^{*\bullet} \right) \to \mathbb{H}^0 \left(X, \mathcal{B}^\bullet \right)=0 \to \mathbb{H}^1\left( X, \mathcal{F}_0^\bullet \right) \to \mathbb{H}^1\left( X, \mathcal{N}_{\mathcal{F}_0}^{*\bullet} \right) \xrightarrow{\gamma} \mathbb{H}^1\left( X, \mathcal{B}^\bullet \right) \\
 \to \mathbb{H}^2\left( X, \mathcal{F}_0^\bullet \right) \to \mathbb{H}^2 \left(X, \mathcal{N}_{\mathcal{F}_0}^{*\bullet} \right) \to \mathbb{H}^2\left( X,  \mathcal{B}^\bullet \right) \to \cdots
\end{align*}
Since $\mathbb{H}^1\left(X, \mathcal{F}_0^\bullet \right) \cong \mathbb{H}^1\left( X, \mathcal{N}_{\mathcal{F}_0}^{*\bullet} \right)$, $\gamma$ is the zero map. This implies that $\mathbb{H}^2\left( X, \mathcal{F}_0^\bullet \right)\to \mathbb{H}^2\left( X, \mathcal{N}_{\mathcal{F}_0}^{*\bullet} \right) $ is injective. It remains to show that it is surjective. In fact, first we note that we have $\mathbb{H}^2\left( X, \mathcal{N}_{\mathcal{F}_0}^{*\bullet} \right) \cong \mathbb{H}^2\left( X, \mathcal{E}_{\mathcal{L}_0}^\bullet \right) \xrightarrow{\beta} \mathbb{H}^2\left(X, \mathcal{E}_{\mathcal{L}_0\otimes K_X}^\bullet \right)\cong \mathbb{H}^2\left( X, T_{\mathcal{F}_0}^\bullet \right)$ from $(\ref{nt17})$. Let $\left(  \left\{ b_{ij}^*:=\sum_{\alpha=1}^{n-1} w_{0i}^1\wedge \cdots \wedge b_{ij}^\alpha \wedge \cdots \wedge w_{0i}^{n-1} \right\} , \left\{\left( c_{ijk}, e_{ijk} \right) \right\} \right) \in C^1\left( \mathcal{U}, \mathcal{S}^1 \right)\bigoplus C^2\left(\mathcal{U}, \mathcal{E}_{\mathcal{L}_0}  \right)$ be a $2$-cocycle of $\mathcal{E}_{\mathcal{L}_0}^\bullet$ in terms of a Stein open covering $\mathcal{U}=\left\{ U_i \right\}$ of $X$ with local coordinates $\left( z_i^1,..., z_i^n \right)$ on $U_i$, where $c_{ijk}\in \Gamma\left(U_{ijk}, T_X \right)$ and $e_{ijk}\in \left( U_{ijk}, \mathcal{O}_X \right)$. Then $\left( \left\{ b_{ij}^*\right\}, \left\{\left( c_{ijk}, e_{ijk} \right) \right\}\right)$ is mapped to $\left( \left\{ \overline{b_{ij}}\right\}, \left\{ c_{ijk}\right\} \right)\in C^1\left(\mathcal{U}, \mathscr{H}om_{\mathcal{O}_X}\left( \mathcal{N}_{\mathcal{F}_0}^*,  \frac{\Omega_X^1}{\mathcal{N}_{\mathcal{F}_0}^*} \right) \right) \bigoplus C^0 \left( \mathcal{U}, T_X \right)$ under the isomorphism $\mathbb{H}^2\left( X, \mathcal{E}_{\mathcal{L}_0}^\bullet \right) \cong \mathbb{H}^2\left( X, \mathcal{N}_{\mathcal{F}_0}^{*\bullet} \right)$, where $\overline{b_{ij}}\left( w_{0i}^\alpha\right)= \overline{b_{ij}^\alpha}$. On the other hand, $\left( \left\{ b_{ij}^*\right\}, \left\{ \left( c_{ijk}, e_{ijk} \right)\right\} \right)$ is mapped to $\left( \left\{ \overline{b_{ij}'} \right\} \left\{ c_{ijk} \right\} \right)\in C^1\left( \mathcal{U}, \mathscr{H}om_{\mathcal{O}_X}\left( T_{\mathcal{F}_0}, \frac{T_X}{T_{\mathcal{F}_0}} \right) \right) \bigoplus C^2\left( \mathcal{U}, T_X \right)$ under $\mathbb{H}^2\left( X, \mathcal{E}_{\mathcal{L}_0}^\bullet \right)\xrightarrow{\beta} \mathbb{H}^2\left( X, \mathcal{E}_{\mathcal{L}_0\otimes K_X}^\bullet \right) \cong \mathbb{H}^2\left( X, T_{\mathcal{F}_0}^\bullet \right)$, where $\overline{b_{ij}'}\left( T_{0i} \right)=\overline{\sum_{\alpha=1}^n i_{\frac{\partial}{\partial z_i^1}\wedge \cdots \wedge \frac{\partial}{\partial z_i^n}}\left( b_{ij}^*\wedge dz_i^\alpha\right) \frac{\partial}{\partial z_i^\alpha}}$. Then we claim that

{\small{\begin{align*}
&\left(0, \left\{ \overline{b_{ij}'} \right\},  \left\{ \overline{b_{ij}} \right\} , \left\{ c_{ijk} \right\}  \right) \\&\in C^2\left( \mathcal{U}, \mathscr{H}om_{\mathcal{O}_X}\left(  \mathcal{N}_{\mathcal{F}_0}^* , T_{\mathcal{F}_0}^* \right)\right) \bigoplus  C^1\left( \mathcal{U}, \mathscr{H}om_{\mathcal{O}_X}\left( T_{\mathcal{F}_0}, \frac{T_X}{T_{\mathcal{F}_0}} \right) \right) \bigoplus  C^1\left( \mathcal{U}, \mathscr{H}om_{\mathcal{O}_X}\left( \mathcal{N}_{\mathcal{F}_0}^*, \frac{\Omega_X^1}{\mathcal{N}_{\mathcal{F}_0}^*} \right) \right) \bigoplus C^2\left( \mathcal{U}, T_X \right)
\end{align*}}}
defines a $2$-cocycle of $\mathcal{F}_0^\bullet$. Then this implies that $\mathbb{H}^2\left(X, \mathcal{F}_0^\bullet \right) \to \mathbb{H}^2\left( X, \mathcal{N}_{\mathcal{F}_0}^{*\bullet} \right)$ is surjective. In fact, it is sufficient to show that
\begin{align}\label{nt23}
i_{b_{ij}' \left( T_{0i}\right)}\left( w_{0i}^\beta \right) + i_{T_{0i}}\left( b_{ij}\left( w_{0i}^\beta \right) \right)=0
\end{align}
In fact, $(\ref{nt23})$ follows from 
\begin{align*}
&i_{b_{ij}' \left( T_{0i}\right)}\left( w_{0i}^\beta \right) + i_{T_{0i}}\left( b_{ij}\left( w_{0i}^\beta \right) \right)= \sum_{\alpha=1}^n i_{\frac{\partial}{\partial z_i^1}\wedge \cdots \wedge \frac{\partial}{\partial z_i^n}}\left(b_{ij}^*\wedge dz_i^\alpha  \right) i_{\frac{\partial}{\partial z_i^\alpha} }\left( w_{0i}^\beta \right) + i_{T_{0i}}\left( b_{ij}^\beta \right) \\
&= i_{\frac{\partial}{\partial z_i^1}\wedge \cdots \wedge \frac{\partial}{\partial z_i^n}}\left( b_{ij}^* \wedge w_{0i}^\beta \right) + i_{T_{0i}}\left( b_{ij}^\beta \right) = - i_{\frac{\partial}{\partial z_i^1}\wedge \cdots \wedge \frac{\partial}{\partial z_i^n}}\left( w_{0i}^1\wedge \cdots \wedge w_{0i}^{n-1} \wedge b_{ij}^\beta \right) + i_{T_{0i}}\left( b_{ij}^\beta \right) \\
&= - \sum_{\alpha=1}^n i_{\frac{\partial}{\partial z_i^1} \wedge \cdots \wedge \frac{\partial}{\partial z_i^n}}\left( w_{0i}^1 \wedge \cdots \wedge w_{0i}^{n-1} \wedge dz_i^\alpha \right) i_{\frac{\partial}{\partial z_i^\alpha}}\left( b_{ij}^\beta \right) + i_{T_{0i}}\left( b_{ij}^\beta \right) = 0.
\end{align*}

\begin{remark}
In particular, if $X$ is a compact complex surface, we have isomorphisms of deformation functors from $(\ref{nt24})$ and $(\ref{nt25})$
\begin{align*}
\textnormal{FDef}_{\left(X, T_{\mathcal{F}_0}, \mathcal{N}_{\mathcal{F}_0}^* \right)} \cong \textnormal{FDef}_{\left( X, T_{\mathcal{F}_0} \right)} \cong \textnormal{FDef}_{\left(X, \mathcal{N}_{\mathcal{F}_0}^* \right)}
\end{align*}
We also note that when $X$ is a compact complex surface, we have an isomorphism of complexes of sheaves
\begin{center}
$\begin{CD}
@. 0 @. 0 @. 0 \\
@. @AAA @AAA @AAA \\
0@>>> \Omega_X^1\otimes \mathcal{L}_0 @>\cong >> T_X \otimes \mathcal{L}_0\otimes K_X @>>> 0 @>>> 0 \\
@. @A\tilde{E}_0AA @A\tilde{D}_0 AA  @AAA \\
0 @>>> \mathcal{E}_{\mathcal{L}_0} @>\cong >> \mathcal{E}_{\mathcal{L}_0\otimes K_X} @>>> 0 @>>> 0
\end{CD}$
\end{center}
which directly shows that $\mathbb{H}^i\left( X, T_{\mathcal{F}_0}^\bullet \right)\cong \mathbb{H}^i\left( X, \mathcal{N}_{\mathcal{F}_0}^{*\bullet} \right), i\geq 0$.

\end{remark}

\bibliographystyle{amsalpha}
\bibliography{References-Rev9}

\end{document}